\documentclass{article}

\usepackage[latin1,utf8]{inputenc}
\usepackage[english]{babel}
\usepackage{amsmath, amsthm, amscd, amsfonts, amssymb, graphicx, color,mathrsfs,mathtools,enumerate, enumitem}
\usepackage{cancel,soul}
\usepackage{color}
\usepackage{amscd,amsxtra}
\usepackage{latexsym}
\usepackage{bm}
\usepackage{hyperref}
\usepackage[a4paper, top=2.7cm, bottom=3cm, left=2.7cm, right=2.7cm]{geometry}
\usepackage{authblk}

\allowdisplaybreaks

\newcommand{\mathd}{\mathrm{d}}

\newcommand{\tmtextit}[1]{\text{{\itshape{#1}}}}
\newenvironment{proof*}[1]{\noindent\textbf{#1\ }}{\hspace*{\fill}$\Box$\medskip}

\newtheorem{thm}{Theorem}[section]
\newtheorem{theorem}[thm]{Theorem}
\newtheorem{hyp}[thm]{Hypotheses}{\rm}
\newtheorem{lemm}[thm]{Lemma}
\newtheorem{cor}[thm]{Corollary}
\newtheorem{prop}[thm]{Proposition}
\newtheorem{defi}[thm]{Definition}
\newtheorem{rmk}[thm]{Remark}{\rm}

\newtheorem{corollary}[thm]{Corollary}
\newtheorem{definition}[thm]{Definition}
\newtheorem{lemma}[thm]{Lemma}
\newtheorem{proposition}[thm]{Proposition}
\newtheorem{remark}[thm]{Remark}

\newcommand{\C}{{\mathscr C}}
\newcommand{\R}{{\mathbb R}}
\newcommand{\N}{{\mathbb N}}
\newcommand{\E}{{\mathbb E}}
\newcommand{\eps}{\varepsilon}
\newcommand{\Id}{{\operatorname{I}}}
\newcommand{\abs}[1]{{\left|#1\right|}}
\newcommand{\norm}[1]{{\left\|#1\right\|}}
\newcommand{\scal}[2]{{\left\langle #1,#2\right\rangle}}

\newcommand{\eqsys}[1]{{\left\{\begin{array}{ll}#1\end{array}\right.}}

\newcommand{\plim}[1]{\underset{#1}{\pi\text{-}\operatorname{lim}\,}}
\newcommand{\Law}{\operatorname{Law}}
\newcommand{\free}{\operatorname{free}}
\newcommand{\inte}{\operatorname{int}}
\newcommand{\tot}{\operatorname{tot}}

\newcommand{\CM}{\operatorname{CM}}
\newcommand{\supp}{\operatorname{supp}}
\newcommand{\loc}{\operatorname{loc}}
\newcommand{\Gibbs}{\operatorname{Gibbs}}

\newcommand{\HH}{\mathcal H}
\newcommand{\Tr}{\operatorname{Tr}}
\newcommand{\dd}{\,\mathrm d}
\newcommand{\dis}{\operatorname{discrete}}
\newcommand{\Ising}{\operatorname{Ising}}

\begin{document}

\frenchspacing

\title{Integration by parts and DLR equations for SPDEs on the whole line}
\date{}

\author[1]{Davide A. Bignamini \thanks{\texttt{da.bignamini@uninsubria.it}}}
\author[2]{Francesco C. De Vecchi \thanks{\texttt{francescocarlo.devecchi@unipv.it}}}
\author[3]{Carlos Villanueva Mariz \thanks{\texttt{cavima2000@zedat.fu-berlin.de}}}

\affil[1]{ Dipartimento di Scienza e Alta Tecnologia, Università degli Studi dell'Insubria}
\affil[2]{Dipartimento di Matematica, Università degli Studi di Pavia}
\affil[3]{Institut f\"ur Mathematik, Freie Universit\"at Berlin}

\maketitle

\begin{abstract}
We study a stochastic reaction–diffusion equation on the whole real line, driven by space-time white noise, with a possibly nonlocal diffusion operator, such as a fractional Laplacian, and a coercive potential with polynomial growth. On the whole line, unlike on bounded intervals, its invariant measures cannot be described by an explicit density with respect to a Gaussian measure. We set the equation in suitable weighted spaces and show that invariant and infinitesimally invariant measures coincide, as do reversible measures and solutions to an integration-by-parts formula. We then show that the reversible measures are exactly the Gibbs measures defined through the Dobrushin–Lanford–Ruelle equations, without any moment condition, thanks to a priori estimates that also yield a Fernique-type theorem, uniform with respect to translations. In this way, we connect the Da Prato–Zabczyk approach to invariant measures of SPDEs with the DLR approach from statistical mechanics to Gibbs measures. Finally, we prove the existence of a reflection positive reversible measure and, in analogy with the one-dimensional long-range Ising model, we show that the reversible measure is unique, ergodic and translation invariant when the kernel of the nonlocal operator has a finite first moment, while it can be non-unique, for suitable double-well potentials, when the kernel has a power-law decay with infinite first moment.\\
\vspace{3pt}

 \end{abstract}

\tableofcontents

\section{Introduction}

\subsection{Motivations and setting}
The aim of this paper is to develop a systematic approach to the study of the invariant measures of stochastic partial differential equations (SPDEs) on the whole real line, and of their relation with Gibbs measures. Since this problem, in its full generality, combines difficulties of very different nature, here we focus on the regular case, namely on scalar random fields defined on $\R$. More precisely, we consider the SPDE
\begin{equation}\label{eq:LangevinSPDE1}
\partial_t\varphi_t=-(H^{\beta}+m^2)\varphi_t-V'(\varphi_t)+\sqrt{2}\dot{W}_t,
\end{equation}
where $\dot{W}$ is a space-time white noise, $m>0$, $H^{\beta}$ is a linear, possibly non-local, operator and $V$ is the potential of the nonlinear part of the drift. The operator $H^{\beta}$ is of the form
\[
H^{\beta}f(\xi)=-\beta_0\Delta f(\xi)+\frac{1}{2}\int_{\mathbb{R}}\big[2 f(\xi)-f(\xi+\xi')-f(\xi-\xi')\big]\bar{\beta}(\xi')\,\mathd \xi',
\]
where $\beta_0\geq0$ and $\bar\beta$ is the density of the L\'evy measure of a Brownian motion subordinated by a complete subordinator, such that $\bar\beta(\xi)\lesssim\abs{\xi}^{-1-2s_2}$ for $\abs{\xi}\geq1$, for some $s_2>0$, and, if $\beta_0=0$, $\bar\beta(\xi)$ is comparable with $\abs{\xi}^{-1-2s_1}$ near the origin, for some $s_1\in(\frac12,1)$ (see \ref{H0} for the precise assumptions). In particular, $-H^{\beta}$ is the generator of a symmetric L\'evy process, and this class of operators includes the Laplacian, the fractional Laplacians $(-\Delta)^s$ with $s\in(\frac12,1)$, their positive linear combinations and the operators associated with the quadratic pair interactions considered in \cite{OsadaSpohn1999}.\footnote{The class of operators $H^{\beta}$ is natural also from the point of view of Euclidean field theory. Indeed, the Gaussian measures with covariance $(H^{\beta}+m^2)^{-1}$ are reflection positive and, conversely, every reflection positive Gaussian measure on $\R$ with the natural invariance, regularity and mixing properties has a covariance of this form, apart from the conditions on the behavior of $\bar\beta$ at $0$ and at infinity and with $m\geq0$ (see Section \ref{appendix:reflectionpositivity}, and in particular Theorem \ref{theorem:rprepresentation2}).} 
The condition $s_1>\frac12$ guarantees that we are in the non-singular case, in which the stochastic convolution associated with \eqref{eq:LangevinSPDE1} takes values in spaces of continuous functions and no renormalization is needed (see Lemma \ref{lemma:stochconv}). The potential $V\in C^2(\R)$ satisfies, for some $q\geq1$, $C,\zeta_0>0$, $\zeta_1\geq0$ and $M_0,M_1\in\R$,
\begin{align*}
&\abs{V'(\xi)}\leq C(1+\abs{\xi}^q),\qquad \abs{V''(\xi)}\leq C(1+\abs{\xi}^{q-1}),\qquad V''(\xi)\geq M_0,\\
&V(\xi)\geq M_1,\qquad \xi V'(\xi)\geq\zeta_0\abs{\xi}^{q+1}-\zeta_1,\qquad \xi\in\R,
\end{align*}
(see \ref{H2}, \ref{H4} and \ref{H5}). The conditions in the first line guarantee the well-posedness of \eqref{eq:LangevinSPDE1} (see Proposition \ref{prop:wellposedness}), the lower bound on $V$ is needed to define the Gibbs measures on bounded intervals (see Section \ref{section:DLRequation}), and the last condition, which only concerns the behavior of $V$ at infinity, is used to obtain a priori estimates (see Section \ref{section:Fernique}). For instance, every polynomial of even degree with positive leading coefficient satisfies these assumptions (when $H^{\beta}=-\Delta$ the corresponding Gibbs measures are the so-called $P(\varphi)_1$ measures, see, e.g., \cite{Lorinczi2011}), as well as non-convex double-well potentials (see Remark \ref{rmk:H5}), which are relevant for phase transitions (see Section \ref{section:nonuniqueness}).\\

\noindent Equation \eqref{eq:LangevinSPDE1} is the Langevin dynamics associated with the (formal) action functional
\begin{equation}\label{eq:actionmain}
\mathcal{S}(\varphi)=\mathcal{S}_{\free}(\varphi)+\mathcal{S}_{\inte}(\varphi):=\frac{1}{2}\int_{\R}\varphi(\xi)[(H^{\beta}+m^2)\varphi(\xi)]\mathd \xi+\int_{\R}V(\varphi(\xi))\mathd \xi,
\end{equation}
namely it can be written as $\partial_t\varphi_t=-\delta_{\varphi}\mathcal{S}(\varphi_t)+\sqrt{2}\dot{W}_t$, and, as for finite dimensional gradient systems, its natural invariant measure is the Gibbs measure
\begin{equation}\label{eq:Gibbs1}
\mu_{\Gibbs}(\mathd \varphi)= \frac{\text{``$e^{-\mathcal{S}(\varphi)}$''}}{\mathcal{Z}}\mathcal{D}\varphi,
\end{equation}
where $\mathcal{D}\varphi$ is a non-existent Lebesgue measure on some space of functions. Formally, $\mu_{\Gibbs}$ is the unique invariant measure of \eqref{eq:LangevinSPDE1}, and the solution of \eqref{eq:LangevinSPDE1} with initial law $\mu_{\Gibbs}$ is a reversible process. The relation between Langevin (or, in the discrete case, Glauber) dynamics and Gibbs measures is a classical topic, with applications in statistical mechanics (see, e.g., \cite{Glauber1963,HolleyStroock1977} and \cite[Chapter 12]{Dap-Zab1996}) and in Euclidean quantum field theory, through the stochastic quantization program proposed in \cite{ParisiWu1981} (see, e.g., \cite{Iwata1987} for the $P(\varphi)_1$ field on $\R$, and \cite{DaPratoDebussche2003,Hairer2014RegularityStructures,MourratWeber2017,GubinelliHofmanova2021,AlbeverioKusuoka2020} for singular models).\\

\noindent Under the previous assumptions, equation \eqref{eq:LangevinSPDE1} can be studied in suitable spaces of functions on $\R$, which we discuss below, by means of the semigroup approach to SPDEs of Da Prato and Zabczyk (see \cite{Dap-Zab2014,Dap-Zab1996,Dap2004}). For every initial datum $x$ in such a space, \eqref{eq:LangevinSPDE1} has a unique (generalized) mild solution $\{\varphi_t\}_{t\geq0}$ (see Proposition \ref{prop:wellposedness}), and we can consider the associated Markov transition semigroup
\begin{equation}\label{eq:definitionsemigroup}
P(t)\Phi(x)=\mathbb{E}\left[\left.\Phi(\varphi_t)\right| \varphi_0=x\right],\qquad t\geq0,
\end{equation}
where $\Phi$ is a real bounded and Borel measurable function defined on a suitable functional space.
This gives a natural dynamical notion of invariant measure, namely a probability measure $\mu$ such that for every bounded continuous function $\Phi$ and $t\geq0$ we have
\begin{equation}\label{eq:def-inv-int}
\int P(t)\Phi\,\mathd\mu=\int\Phi\,\mathd\mu
\end{equation}
 and of reversible measure, namely a probability measure $\mu$ such that for every bounded continuous functions $\Phi,\Psi$ we have
\begin{equation}\label{eq:reversible}
\int \Psi(x) P(t)\Phi(x) \mu(\mathd x)=\int \Phi(x)P(t)\Psi(x) \mu(\mathd x),\qquad t\geq0.
\end{equation}
These notions have infinitesimal counterparts, expressed in terms of the Kolmogorov operator associated with \eqref{eq:LangevinSPDE1}; and also, for these gradient-type SPDEs, via the following integration by parts formula
\[
\int\scal{\delta_{\varphi}\Phi(\varphi)}{f}_{L^2}\mu(\mathd\varphi)=\int\Phi(\varphi)\left[\scal{\varphi}{(H^{\beta}+m^2)f}_{L^2}+\scal{V'\circ\varphi}{f}_{L^2}\right]\mu(\mathd\varphi),
\]
for every smooth cylinder function $\Phi$ and every $f\in C^{\infty}_c(\R)$ (see Section \ref{section:abstractreformulation}, and in particular Definitions \ref{definition:invariant} and \ref{definition:reversible}). Every reversible measure is invariant, as it can be seen by taking $\Psi=1$ in \eqref{eq:reversible}, while the converse is a delicate problem already for lattice systems (see, e.g., \cite{HolleyStroock1977}), which we do not address here.\\

\noindent However, the situation on the whole real line is very different from the one on a bounded interval (or on the torus). Indeed, under suitable conditions on $V$, equation \eqref{eq:LangevinSPDE1} (with suitable boundary conditions for $H^{\beta}$) has a unique invariant measure, which is reversible and has an explicit formulation given by
\begin{equation}\label{eq:Gibbs2}
 \mu_{\Gibbs}(\mathd \varphi):=\frac{e^{-\mathcal{S}_{\inte}(\varphi)}}{\mathcal{Z}_{\inte}}\mu_{\free}^{H^{\beta}+m^2}(\mathd \varphi),\qquad 
\end{equation}
where
\[
\mathcal{S}_{\inte}(\varphi):=\int_{\R}V(\varphi(\xi))\mathd \xi,\qquad \mathcal{Z}_{\inte}:=\int e^{-\mathcal{S}_{\inte}(\varphi)}\mu_{\free}^{H^{\beta}+m^2}(\mathd \varphi),
\]
 $\mu^{H^{\beta}+m^2}_{\free}$ is the centered Gaussian measure with covariance $(H^{\beta}+m^2)^{-1}$ (see \cite[Chapter 7 and Section 8.6]{Dap-Zab1996} and \cite{Zabczyk1989}, and \cite{DaPratoDebussche2003,TsatsoulisWeber2018,HairerMattingly2018} for analogous results in the singular setting of stochastic quantization). On $\R$ both these properties are lost. First, even if \eqref{eq:LangevinSPDE1} has a gradient structure, its invariant measures cannot be written in the form \eqref{eq:Gibbs2}, since $\mathcal{S}_{\inte}(\varphi)$ is not well defined for $\mu^{H^{\beta}+m^2}_{\free}$-almost every $\varphi$, and the Gibbs measures are expected to be singular with respect to $\mu^{H^{\beta}+m^2}_{\free}$ (for instance, since $\mu^{H^{\beta}+m^2}_{\free}$ is ergodic with respect to the translations of $\R$, every translation invariant and ergodic probability measure different from $\mu^{H^{\beta}+m^2}_{\free}$ is singular with respect to it). Second, the uniqueness of the invariant measure is not guaranteed. Indeed, \eqref{eq:LangevinSPDE1} cannot be studied in $L^2(\R)$, since the noise and the operator $H^{\beta}+m^2$ are invariant under translations; for every $t>0$ the stochastic convolution, namely the solution of \eqref{eq:LangevinSPDE1} with $V=0$ and zero initial datum, is a stationary Gaussian random field on $\R$ with positive constant variance, so that its norm in $L^2(\R)$ is almost surely infinite. Hence we are forced to work in spaces with an integrable weight, such as the spaces 
\begin{equation}\label{eq:weighted-introduction}
 L^p_{r,k}:=L^p(\R,\rho_{r,k}(\xi)\mathd\xi),\qquad \rho_{r,k}(\xi)=(1+k\abs{\xi}^2)^{-r/2},\qquad r\in(1,1+2s_2)
 \end{equation}
 and $k>0$ small enough (see Subsection \ref{section:weights} and \ref{H3}). In these spaces the classical strategies to prove uniqueness, based on the strong Feller property and on irreducibility, or on their asymptotic versions (see \cite[Chapter 7]{Dap-Zab1996}, \cite{Hairer2008,HairerMattingly2006}), cannot be applied in general. The choice of the space matters already for the linear equation: for $H^{\beta}=-\Delta$ and $V=0$, the invariant measure is unique in every space $L^2_{r,k}$, while there are infinitely many invariant measures in the $L^2$-spaces with weight $e^{-\kappa\abs{\xi}}$ and $\kappa^2>4m^2$ (see Remark \ref{rmk:H3} and \cite[Section 9.4]{Dap-Zab1996}). In the non-linear case, the dynamics associated with \eqref{eq:LangevinSPDE1} in $L^2_{r,k}$ is irreducible (see Proposition \ref{prop:irreducibility-c}), but in Section \ref{section:nonuniqueness} we present examples with several invariant measures, so that in general it is neither strong Feller nor asymptotically strong Feller.\\

\noindent Another approach to Gibbs measures on unbounded domains, which comes from the statistical mechanics of lattice systems, is based on the Dobrushin--Lanford--Ruelle (DLR) equations (see \cite{Dobrushin1968,LanfordRuelle1969}, \cite{Ruelle1970,LebowitzPresutti1976} for unbounded spins, and the monographs \cite{Georgii2011,FriedliVelenik2018}), and it has been extended to the continuum, in particular to Euclidean field theories and to Gibbs measures on path spaces (see, e.g., \cite{BetzLorinczi2003,CourregeRenouard1975,GuerraRosenSimon1975II,GuerraRosenSimon1976,Lorinczi2011,OsadaSpohn1999,RoyerYor1976}). The idea is not to define the Gibbs measure on $\R$ directly, but through its conditional probabilities. Even if \eqref{eq:Gibbs2} is meaningless on $\R$, the Gibbs measure on a bounded interval $\mathcal{I}$ with a prescribed exterior condition $\psi:\mathcal{I}^c\rightarrow\R$ is well defined and explicit, and a probability measure $\mu_{\Gibbs}$ on a space of functions on $\R$ is a Gibbs measure if its conditional probabilities given the values of the field outside the bounded intervals are given by these finite volume Gibbs measures. More precisely, for every open and bounded interval $\mathcal{I} \subset \R$ and (almost) every $\psi:\mathcal{I}^c \rightarrow \R$ we require
\begin{equation}\label{eq:GibbsDLR}
\mu_{\Gibbs}^{\mathcal{I},\psi}=\mu_{\Gibbs}(\mathd \varphi|\varphi|_{\mathcal{I}^c}=\psi)=\frac{1}{\mathcal{Z}^{\mathcal{I}}(\psi)}e^{-\int_{\mathcal{I}}V(\varphi(\xi))\mathd \xi}\mu_{\free}^{H^{\beta}+m^2}(\mathd \varphi|\varphi|_{\mathcal{I}^c}=\psi),
\end{equation}
where the last measure is the conditional probability of the Gaussian measure $\mu_{\free}^{H^{\beta}+m^2}$ given the values of the field outside $\mathcal{I}$ (see Section \ref{section:DLRequation}, and in particular Definition \ref{definition:DLR}); and where 
\[\mathcal{Z}^{\mathcal{I}}(\psi)=\int e^{-\int_{\mathcal{I}}V(\varphi(\xi))\mathd \xi}\mu_{\free}^{H^{\beta}+m^2}(\mathd \varphi|\varphi|_{\mathcal{I}^c}=\psi)\]
is a function of the boundary field $\psi$ having the role of normalization constant. In other words, the Gibbs measures on $\R$ are built from the Gibbs measures on the bounded intervals by means of conditional probabilities without considering the dynamics on the whole real line given by the SPDE \eqref{eq:LangevinSPDE1}. We stress that, when $H^{\beta}$ is non-local, the conditional measures in \eqref{eq:GibbsDLR} depend on the whole boundary condition $\psi$ and not only on its values at the endpoints of $\mathcal{I}$. The DLR equations \eqref{eq:GibbsDLR} are implicit: in general their solutions are not explicit and they are not necessarily unique, the non-uniqueness corresponding to the presence of phase transitions (see Section \ref{section:nonuniqueness}).\\

\noindent The main result of this paper builds a bridge between the dynamical approach of Da Prato and Zabczyk and the DLR approach from statistical mechanics. This goal is achieved in the main result of the paper, Theorem \ref{theorem:super-main}. Below, we present a concise and simplified version of the theorem, omitting the technical definitions required for its precise formulation.
\begin{theorem}\label{thm:main-intrdouzione}
    Let $\mu$ be a Borel probability measure on $L^{2}_{r,k}$ such that $\mu(L^{2q}_{r,k}\cap C(\R))=1$. Then the following conditions are equivalent:
    \begin{itemize}
        \item for every smooth cylinder function $\Phi$ and every $f\in C^{\infty}_c(\R)$ it holds
        \[
\int\scal{\delta_{\varphi}\Phi(\varphi)}{f}_{L^2}\mu(\mathd\varphi)=\int\Phi(\varphi)\left[\scal{\varphi}{(H^{\beta}+m^2)f}_{L^2}+\scal{V'\circ\varphi}{f}_{L^2}\right]\mu(\mathd\varphi),
\]
\item for every open and bounded interval $\mathcal{I} \subset \R$ and (almost) every $\psi:\mathcal{I}^c \rightarrow \R$ we require
\begin{equation*}
\mu_{\Gibbs}^{\mathcal{I},\psi}=\mu_{\Gibbs}(\mathd \varphi|\varphi|_{\mathcal{I}^c}=\psi)=\frac{1}{\mathcal{Z}^{\mathcal{I}}(\psi)}e^{-\int_{\mathcal{I}}V(\varphi(\xi))\mathd \xi}\mu_{\free}^{H^{\beta}+m^2}(\mathd \varphi|\varphi|_{\mathcal{I}^c}=\psi),\qquad
\end{equation*}
where the last measure is the conditional probability of the Gaussian measure $\mu_{\free}^{H^{\beta}+m^2}$ given the values of the field outside $\mathcal{I}$.
    \end{itemize}
\end{theorem}

\noindent Hence, among the probability measures concentrated on $L^{2q}_{r,k}\cap C(\R)$, the Gibbs measures on $\R$ in the DLR approach sense are exactly the reversible measures of the Langevin dynamics \eqref{eq:LangevinSPDE1}, and both the dynamical and the statistical mechanics tools can be used to study them. For lattice systems taking values in $\R^n$, the Gibbs measures can be characterized through their quasi-invariance with respect to shifts, which is an integrated form of the integration by parts formula (see \cite{AlbeverioKondratievRockner1997b, AlbeverioKondratievRockner1997a}), while in the continuum a characterization of the reversible measures as Gibbs measures was obtained in \cite{Funaki1991} for a local operator and a bounded non-linearity, under a moment condition. In the singular setting, the uniqueness of the invariant measure of the Langevin dynamics on the whole space has been recently proved in the high temperature regime: for the $\varphi^4_2$ model in the whole region of finite susceptibility \cite{BauerschmidtDagallierWeber2025}, by a continuum version of the Holley–Stroock strategy, which deduces the uniqueness from a log-Sobolev inequality uniform in the volume (proved in \cite{BauerschmidtDagallier2024}), and for the $\Phi^4_3$ model at small coupling \cite{DuchHairerYiZhao2025}, by showing that two solutions driven by the same noise converge to each other exponentially fast. These results concern all invariant measures, not only the reversible ones, but they rely on the absence of phase transitions and do not address the relation with the Gibbs measures defined through conditional probabilities (see \cite[Remark 1.6]{DuchHairerYiZhao2025}). Our results are complementary: in the non-singular one-dimensional case, we characterize the reversible measures as the solutions to the DLR equations for every potential $V$ satisfying our assumptions, also in the presence of phase transitions, while our uniqueness results concern only the reversible measures.\\

\noindent Finally, we present some other important related results contained in the paper. First, we prove the existence of a reversible measure, or equivalently of a solution to the DLR equations, as a limit of finite-volume approximations with Gaussian or Dirichlet boundary conditions; with Gaussian boundary conditions, the limit is also reflection positive (see Section \ref{section:existence}). The equivalence between reversible measures and solutions to the DLR equations also allows us to study their uniqueness, translation invariance, and ergodicity. If the reversible measure of \eqref{eq:LangevinSPDE1} is unique, then it is ergodic for the Langevin dynamics (see Subsection \ref{section:ergodicity}) and invariant with respect to translations; moreover, a translation-invariant reversible measure always exists, also when uniqueness fails (see Subsection \ref{section:translationinvariance}). Whether every invariant measure is reversible remains an open question in our setting, while for lattice systems it is known in some cases (see \cite{Holley1971,HolleyStroock1977,Fritz1982}). Finally, we show that, as for one-dimensional long-range Ising models (see \cite{Dyson1969,FrohlichSpencer1982}), the uniqueness of the solution to the DLR equations depends on the decay at infinity of the kernel $\bar\beta$ of the nonlocal part of $H^{\beta}$. These results extend the ones of \cite{OsadaSpohn1999} to general potentials and to the fractional case, in which the nonlocal part of $H^{\beta}$ is a fractional Laplacian (see Proposition \ref{proposition:regimes} in the last part of the introduction).

\subsection{Discussion of the ideas of the proofs}

\subsubsection{Equivalence of reversible measures and DLR solutions}

The proof of the main Theorem \ref{theorem:super-main} combines three ingredients: the dynamical approach to SPDEs of Da Prato and Zabczyk, based on infinite-dimensional analysis, the DLR approach of statistical mechanics, based on conditional measures, and some a priori estimates, which allow us to remove all the moment conditions.

\paragraph{Infinite-dimensional analysis and dynamical approach.} First of all, we study the well-posedness of \eqref{eq:LangevinSPDE1} in the mild sense, following the approach to SPDEs of Da Prato and Zabczyk (see \cite{Dap-Zab1996,Dap-Zab2002,Dap-Zab2014}). As already observed, even in the linear case $V\equiv0$, the mild solution of \eqref{eq:LangevinSPDE1} does not take values in $L^2(\R)$, so we are forced to study it in the weighted $L^p$-spaces defined by \eqref{eq:weighted-introduction}. A fundamental tool in the approach of Da Prato and Zabczyk is the semigroup generated by $-(H^\beta+m^2)$. Hence, in Section \ref{sect:analitic-Hbeta} we prove that, if the weight does not decay too fast at infinity, depending on the decay of the L\'evy density $\bar\beta$, the operator $-(H^\beta+m^2)$ generates a strongly continuous semigroup of negative type on $L^p_{r,k}$ (see Lemma \ref{lemma:semigroupLp}), and we prove some commutator estimates between $H^\beta$ and the weights (see Subsection \ref{section:commutator}). We underline that, to the best of our knowledge, the operators $H^\beta$ have not been studied in this generality in weighted spaces with integrable weights, so we believe that Section \ref{sect:analitic-Hbeta} is of independent interest. Once the well-posedness of the equation is established (see Proposition \ref{prop:wellposedness}), we prove that the associated dynamics is irreducible (see Proposition \ref{prop:irreducibility-c}). However, since the invariant measure can be non-unique (see Section \ref{section:nonuniqueness}), the transition semigroup \eqref{eq:definitionsemigroup} is in general neither strong Feller nor asymptotically strong Feller (see \cite{Hairer2008} and the remark after Proposition \ref{prop:irreducibility-c}), and the usual techniques to prove the uniqueness of the invariant measure cannot be applied.\\

\noindent In order to prove Theorem \ref{theorem:super-main}, we cannot directly use the notion of reversible measure \eqref{eq:reversible}, and so we use an infinitesimal counterpart of it, namely the following integration by parts formula (IPF) (see Definition \ref{definition:reversible}):
\[
\int\scal{h}{\delta_x\Phi(x)}_{L^2}\mu(\mathd x)=\int\left[\scal{(H^\beta+m^2)h}{x}_{L^2}+\scal{h}{V'\circ x}_{L^2}\right]\Phi(x)\mu(\mathd x),
\]
for every cylinder function $\Phi$, namely a smooth function of finitely many scalar products $\scal{x}{f_i}_{L^2}$ with $f_i\in C^\infty_c(\R)$, and for every $h\in C^\infty_c(\R)$. The equivalence between reversible measures and measures satisfying the IPF is related to the equivalence between invariant and infinitesimally invariant measures, namely solutions of the stationary Fokker--Planck--Kolmogorov equation tested against cylinder functions. In finite dimension these notions are equivalent under very general conditions, while in infinite dimension their relationship is more problematic (see \cite[Section 10.1]{BogachevKrylovRocknerShaposhnikov2015} and \cite{BogachevRockner1995}). The main problem is that it is not trivial to prove that there exists a ``good'' set of smooth functions on which the action of the second order Kolmogorov operator associated with \eqref{eq:LangevinSPDE1} is explicit, namely, in the language of evolution equations, a core for the generator of the transition semigroup. In finite dimension, this set is usually the space of smooth functions with compact support. Moreover, as in finite dimension, the semigroup \eqref{eq:definitionsemigroup} is in general not strongly continuous on spaces of continuous functions on $L^2_{r,k}$ endowed with the topology of uniform convergence. For this reason, we use a weaker notion of convergence, the $\pi$-convergence, namely the pointwise convergence of continuous functions which are uniformly bounded with respect to a weight depending on the polynomial growth of the drift. With respect to the $\pi$-convergence, we can define a domain for the second order Kolmogorov operator (see \cite{Pri1999}, and see also similar notions in \cite{Cer1994,Cer1995,Gol-Koc2001}). Hence, in Part \ref{part:abstract} we prove that, for a class of abstract dissipative SPDEs in Hilbert spaces, the cylinder functions form a $\pi$-core for the generator of the transition semigroup (see Theorem \ref{thm: coreN}, whose proof is based on \cite{Manca2008,Manca2009}), and we deduce the equivalence between invariant and infinitesimally invariant measures, and between reversible measures, infinitesimally reversible measures and measures satisfying the IPF (see Section \ref{sect: inv-infinve-abstract}). We underline that, to the best of our knowledge, these results are new for abstract dissipative SPDEs in Hilbert spaces, as well as for SPDEs on unbounded intervals. In particular, as far as the $\pi$-core is concerned, the results in the literature cover only SPDEs with Lipschitz continuous drift (see \cite{Manca2008,Manca2009}), reaction--diffusion SPDEs with polynomial nonlinearities on bounded domains (see \cite{Manca2009}) and the Burgers SPDE on a bounded interval (see \cite{Manca2010}). The counterparts of the abstract theorems in Part \ref{part:abstract} for the concrete equation \eqref{eq:LangevinSPDE1} are Theorems \ref{thm:core-concrete}, \ref{thm:main-inv-ii} and \ref{theorem:infinitesimalsymmetri}.\\

\noindent Finally, we prove that, if the reversible measure is unique, then it is ergodic (see Theorem \ref{thm:rev-ergo}). Similarly, since the translates of a reversible measure are reversible, the unique reversible measure is translation invariant and, in general, a translation invariant reversible measure can be obtained by averaging the translates of any reversible measure (see Subsection \ref{section:translationinvariance}).

\paragraph{DLR approach.} The DLR equations prescribe the conditional probabilities of a measure given the values of the field outside the bounded intervals (see \cite{Georgii2011,FriedliVelenik2018} for lattice systems). In infinite dimensions, the conditional probabilities are given by the disintegration of the measure (see Lemma \ref{lemma:disintegration} and \cite[Theorem 5.3.1]{AmbrosioGigliSavare}). In Section \ref{section:DLRequation}, following \cite{Rockner1985}, we compute the conditional probabilities of the Gaussian measure $\mu_{\free}$ with covariance $(H^\beta+m^2)^{-1}$ given the boundary condition $y$ outside a bounded interval $\mathcal I$. We prove that the conditional probability $\mu^{\mathcal I,y}_{\free}$ is the translation, by the harmonic extension $\varphi^y$ of $y$, of the centered Gaussian measure on $L^2(\mathcal I)$ whose covariance is the inverse of the realization of $H^\beta+m^2$ with homogeneous Dirichlet boundary conditions (see Proposition \ref{prop:gaussL2I} and Corollary \ref{corollary:Guassiantranslation}). Since the conditional probabilities corresponding to different boundary conditions are mutually singular, the DLR equations are written in terms of the densities with respect to $\mu^{\mathcal I,y}_{\free}$ (see Definition \ref{definition:DLR}), or equivalently in terms of the translated measures (see Corollary \ref{remark:DLRtranslated}). The harmonic extensions, namely the solutions of the exterior Dirichlet problem for $H^\beta+m^2$, are a fundamental tool in the DLR approach. Since in our setting the operator $H^\beta$ is in general nonlocal and the boundary conditions are only continuous functions in $L^2_{r,k}$, we had to prove again many classical results on harmonic extensions (see Section \ref{section:harmonicextension}, and \cite{Fukushima2011} for the theory of Dirichlet forms).\\

\noindent The equivalence between the IPF and the DLR equations (see Theorem \ref{theorem:main}) is proved in Section \ref{sect:int-DLR}. First, a measure satisfies the IPF on $\R$ if and only if its conditional probabilities satisfy the IPF on the bounded intervals with the corresponding boundary conditions. This localization holds only for the IPF associated with a reversible dynamics, and this is the reason for the restriction to reversible measures in Theorem \ref{theorem:main}. Then, we use that the SPDE on a bounded interval $\mathcal I$ with boundary condition $y$ has a unique invariant measure, which is the Gibbs measure prescribed by the DLR equations and which is characterized by the IPF on $\mathcal I$ (see Proposition \ref{prop:invmeasureI} and Corollary \ref{cor:ibpuniqueI}). Theorem \ref{theorem:main} is proved under the moment condition \eqref{eq:momenti-DLR}, which can be removed, under \ref{H5}, thanks to an a priori estimate for the solutions of the DLR equations (see Theorem \ref{thm:apriori} and the next paragraph).

\paragraph{A priori estimates.} The previous results require some moment conditions (see \eqref{eq:cond-momenti-inv} and \eqref{eq:momenti-DLR}). In Section \ref{section:Fernique}, under the coercivity assumption \ref{H5}, we remove them by means of some a priori estimates. First, we prove weighted energy estimates in $L^2_{r,k}$, $L^{q+1}_{r,k}$ and $H^{s_1}_{r,k}$, uniform with respect to the translations of the weight, for the solutions of \eqref{eq:LangevinSPDE1} with an additional predictable control (see Theorem \ref{theorem:apriori}); they rely on the commutator estimates of Subsection \ref{section:commutator}, which hold only if the weight does not decay too fast at infinity. These estimates imply that every invariant measure has finite moments (see Lemma \ref{lemma:invariantmoments}) and, by means of the Bou\'e--Dupuis formula (see \cite{BudhirajaDupuis2000} and Lemma \ref{lemma:BD}), a Fernique-type theorem (see Theorem \ref{theorem:Fernique}). The proof is based on the techniques developed to study the exponential integrability of quantum field theories by stochastic quantization (see, e.g., \cite{HairerSteele2022,GubinelliHofmanova2021,MourratWeber2017} and the related a priori estimates for (singular) stochastic quantization SPDEs in the full space \cite{ChandraMoinatWeber2023,ChevyrevGubinelli2025} for polynomial potentials and \cite{ChandraDeLimaFeltesWeber2026} for subquadratic potentials). The main novelty, which allows us to obtain the exponential integrability of the square of the norm, is the use of the Bou\'e--Dupuis formula as in \cite{DeVecchiJanakPriola2026}, where a similar method is applied to the Burgers SPDE on an interval (see also the applications of the Bou\'e--Dupuis formula to quantum field theory \cite{BarashkovGubinelli2020,BarashkovGubinelli2023,BarashkovGunaratnam2025,BarashkovGunaratnamHofstetter2023}). Our proof does not need the locality of $H^\beta$ and it does not distinguish between subquadratic and superquadratic potentials, since it does not use the ``coming down from infinity'' of the solutions as in \cite{ChandraMoinatWeber2023,ChevyrevGubinelli2025}; it also gives an alternative proof, in the regular case, of the optimal exponential integrability proved in \cite{HairerSteele2022}. \\

\noindent Second, we prove directly that every solution of the DLR equations satisfies the moment conditions (see Theorem \ref{thm:apriori}). Indeed, by the DLR equations, the conditional probabilities of a solution are given, up to a translation by the harmonic extension of the boundary condition, by the invariant measures of the dynamics on the bounded intervals, whose moments in the weighted $L^{2q}$ norm are bounded, uniformly with respect to the interval and to the boundary condition, in terms of the weighted $L^{2q}$ norm of the harmonic extension (see \eqref{eq:conditionalmoments} in the proof of Theorem \ref{thm:apriori}). Hence the fundamental estimate is the one on the harmonic extensions, which we obtain by averaging them with respect to the size of the intervals (see Lemma \ref{lemma:harmonicaverage}); this average is essential, since for a fixed interval the harmonic extension depends on the values of the boundary condition near the endpoints (see Remark \ref{rmk:harmonicaverage}). Theorem \ref{thm:apriori} is a continuum and nonlocal counterpart of the superstability estimates of \cite{Ruelle1970,LebowitzPresutti1976} (see Remark \ref{rmk:LebowitzPresutti}).\\

\noindent Combining the previous three ingredients, we obtain Theorem \ref{theorem:super-main}, which we believe to be the first systematic result in the literature that builds a bridge between these two approaches to the study of the invariant measures of \eqref{eq:LangevinSPDE1}. The proof of Theorem \ref{theorem:super-main} can be summarized in the following scheme
\[
\begin{array}{ccc}
			\text{\bf IPF on } \mathbb{R} & \Longleftrightarrow & \text{\bf IPF on }
			\mathcal{I} \text{ \bf with b.c. }y \\
			\phantom{a} & & \\
			\Big{\Downarrow} \pi\text{-core} \hspace{0.7em} \Big{\Uparrow}\text{It\^o formula} &  & \Big{\Downarrow} \pi \text{-core} \hspace{0.7em} \Big{\Uparrow} \text{It\^o formula}\\
			\phantom{a} & & \\
			P(t) \text{\bf-reversible} &  & P^{\mathcal{I}, y}(t) \text{\bf-reversible}\\
			\phantom{a} & & \\
			&  & \Big{\Downarrow}\text{uniqueness of i.m.} \hspace{0.7em}  \Big{\Uparrow}\text{Langevin eq.} \\
			\phantom{a} & &\\
			& & \text{\bf DLR equations}
\end{array},\]
where $P^{\mathcal{I},y}(t)$ is the transition semigroup associated with the SPDE on the bounded interval $\mathcal{I}$ with boundary condition $y$ (see Section \ref{sect:int-DLR}), and where the a priori estimates of Section \ref{section:Fernique} allow us to remove all the moment conditions.

\subsubsection{Applications to existence and uniqueness/non-uniqueness for DLR equations}

\paragraph{Existence and reflection positivity.} In Section \ref{section:existence} we prove the existence of solutions to the DLR equations, which are reversible measures, as limits of finite-volume approximations with Gaussian boundary conditions (see Theorem \ref{theorem:existence}) and with Dirichlet boundary conditions (see Theorem \ref{theorem:existence2}); the tightness follows from a priori estimates uniform with respect to the volume, as in \cite{AlbeverioKusuoka2020,GubinelliHofmanova2021} (see also \cite{de2022singular}), and the limits satisfy the IPF. In Section \ref{appendix:reflectionpositivity} we characterize the centered Gaussian measures which are invariant with respect to the Euclidean group and reflection positive (see \cite{OsterwalderSchrader1973,GlimmJaffe1987,Jaffe2017}): under natural regularity and mixing conditions, their covariances are exactly the Green functions of the operators $H^\beta+m^2$, with $-H^\beta$ generator of a subordinated Brownian motion with complete Bernstein symbol (see Theorem \ref{theorem:rprepresentation2}), which explains our choice of the class of operators \ref{H0}. As a consequence, the solution obtained with Gaussian boundary conditions is reflection positive (see Subsection \ref{section:ref.positive}). Moreover, if the solution of the DLR equations, or equivalently the reversible measure, is unique, then it is ergodic, reflection positive and translation invariant (see Corollary \ref{corollary:uniquetranslation}).

\paragraph{Uniqueness.} In Section \ref{section:uniquenessDLR} we prove that, if the L\'evy density $\bar\beta$ has a finite first moment at infinity (see \eqref{eq:firstmoment}), which means that the interaction between the two half-lines is finite, apart from the contribution of the singularity of $\bar\beta$ at the origin, then, for every potential $V$ satisfying our assumptions, the solution to the DLR equations is unique, and it is ergodic, reflection positive and translation invariant (see Theorem \ref{theorem:uniqueness} and Corollary \ref{corollary:uniqueness}). This is a continuum counterpart of the uniqueness of the Gibbs measure for one-dimensional lattice systems with summable interaction between half-lines (see \cite{Ruelle1968} and \cite[Section 8.3]{Georgii2011}). Since the conditional probabilities with different boundary conditions are mutually singular, the argument for lattice systems cannot be applied directly: we compare them after a Cameron--Martin translation localized near the boundary of the interval, and we control the resulting relative entropy on average with respect to the boundary conditions, by means of the uniform bounds of Section \ref{section:Fernique} and of an integrated estimate of the harmonic measure (see Lemma \ref{lemma:cutoff}). This gives the mutual absolute continuity of any two solutions (see Proposition \ref{prop:equivalenceDLR}), and the uniqueness follows by an argument on the tail $\sigma$-algebra. This result covers the Markovian case $\bar\beta\equiv0$, in which uniqueness was proved, for Gibbs measures relative to Brownian motion, by transfer operator methods (see \cite{CourregeRenouard1975,RoyerYor1976,BetzLorinczi2003,Lorinczi2011}), the uniqueness result of \cite[Theorem 5.4]{OsadaSpohn1999} in the case of quadratic interactions, without the domination assumptions on the potential required therein, and the case of the fractional Laplacians of order $s\in(\frac12,1)$, without Brownian component.

\paragraph{Non-uniqueness.} In Section \ref{section:nonuniqueness} we prove that, if $\beta_0>0$, $\bar\beta(\xi)$ decays as $|\xi|^{-1-2s_2}$ with $s_2\leq\frac12$ and $\bar\beta$ is not too singular at the origin (see \eqref{eq:conditionbarbeta}), then there are double-well potentials for which the DLR equations have at least two solutions (see Theorem \ref{theorem:nonuniqueness}). This is the continuum counterpart of the phase transition of the one-dimensional Ising model with ferromagnetic interaction $J_{i-j}\sim|i-j|^{-p}$, which occurs for $p\in(1,2]$ and large coupling constants (see \cite{Dyson1969,FrohlichSpencer1982}), while for $p>2$ the Gibbs measure is unique (see \cite{Ruelle1968}); in our setting the role of $J_{i-j}$ is played by $\bar\beta(\xi-\xi')$ and $p=1+2s_2$, so that the threshold $s_2=\frac12$ between Theorems \ref{theorem:uniqueness} and \ref{theorem:nonuniqueness} agrees with the one of the Ising model. The proof extends the one of \cite[Theorem 6.1]{OsadaSpohn1999}: we represent the conditional measures as Gibbs measures relative to the Brownian bridge (see Lemma \ref{lemma:extendeduniqueness1}), and, by means of the GKS and Wells inequalities (see \cite{Ginibre1970,Wells1977,BricmontLebowitzPfister1981,MadridSimonWells2022}), we compare the averages of the field on unit blocks with a long-range Ising model, whose coupling constant is in the phase transition regime thanks to the large deviation estimates of \cite{OttoWeberWestdickenberg2014} (see also \cite{BertiniButtaDiGesu2025}) for the double-well potential.

\paragraph{An explicit family of operators.} The previous results give a complete picture, with the exception of a single case, for the following family of operators, which includes the fractional Laplacians and the quadratic interactions considered in \cite{OsadaSpohn1999}, and for which the corresponding model can be seen as a continuum version of the one-dimensional long-range Ising model with interaction $J_{i-j}\sim|i-j|^{-1-2s}$. In the next proposition, we summarize the results in Section \ref{section:OScomparison}.

\begin{proposition}\label{proposition:regimes}
For $\beta_0,\beta_1\geq 0$ with $\beta_0+\beta_1>0$, $\lambda\geq 0$ and $s\in(0,1)$, let
\[
H^{\beta}f(\xi)=-\beta_0\Delta f(\xi)+\frac{\beta_1}{2}\int_{\R}\frac{2f(\xi)-f(\xi+\xi')-f(\xi-\xi')}{(\lambda+|\xi'|)^{1+2s}}\,\mathd \xi',
\]
where, if $\beta_0=0$, we assume $\lambda=0$ and $s\in\left(\frac12,1\right)$, and assume that \ref{H1} and \ref{H3} hold.
\begin{enumerate}
\item If $\beta_1=0$ or $s\in(\frac12,1)$, then, for every potential $V$ satisfying \ref{H2}, \ref{H4} and \ref{H5}, the DLR equations \eqref{eq:DLRequations} have a unique solution, which is also the unique reversible measure for \eqref{eq:LangevinSPDE1} among the Borel probability measures $\nu$ on $L^2_{r,k}$ such that $\nu(L^{2q}_{r,k}\cap C(\R))=1$, and it is ergodic, reflection positive and translation invariant.
\item If $\beta_0,\beta_1>0$ and either $\lambda>0$ and $s\in(0,\frac12]$, or $\lambda=0$ and $s\in(0,\frac12)$, then there are  potentials $V$ satisfying \ref{H2}, \ref{H4} and \ref{H5} for which the DLR equations \eqref{eq:DLRequations} have at least two solutions, namely \eqref{eq:LangevinSPDE1} has at least two reversible measures.
\end{enumerate}
\end{proposition}

\noindent In \cite{OsadaSpohn1999}, Gibbs measures relative to Brownian motion are considered, which, for the operators of Proposition \ref{proposition:regimes}, correspond to the case $\beta_0>0$ and $\lambda>0$ (see Theorem \ref{theorem:BB} and Remark \ref{remark:OSscaling} for the relation between the definition of \cite{OsadaSpohn1999} and Definition \ref{definition:DLR}). In this case, uniqueness is proved in \cite[Theorems 2.3 and 5.4]{OsadaSpohn1999} for $s>\frac12$, under a domination assumption on the potential, and non-uniqueness is proved in \cite[Theorems 2.4 and 6.1]{OsadaSpohn1999} for $\lambda=1$, $s\in(0,\frac12]$ and a quartic double-well potential. Proposition \ref{proposition:regimes} is more general: the uniqueness holds without domination assumptions on the potential, the non-uniqueness holds for general double-well potentials and, above all, Proposition \ref{proposition:regimes} covers the fractional case $\lambda=0$, both with $\beta_0>0$ and with $\beta_0=0$, which is not covered by the results of \cite{OsadaSpohn1999} (see Section \ref{section:OScomparison} for a detailed comparison). The only case which is not covered by Proposition \ref{proposition:regimes} is $\beta_0>0$, $\lambda=0$ and $s=\frac12$.

\section{Notation and basic definitions}\label{section:notation}

In this section we fix some notation and we recall some basic definitions which are used, but not introduced, in the rest of the paper.

\paragraph{Notation.}
\begin{itemize}
\item $\N:=\{1,2,\ldots\}$ and $\R_+:=[0,+\infty)$; for $a,b\in\R$, $a\wedge b:=\min\{a,b\}$ and $a\vee b:=\max\{a,b\}$.
\item $A\lesssim B$ means that $A\leq cB$ for a constant $c>0$ independent of the variables under consideration, and $A\lesssim_aB$ that $c$ depends on the parameter $a$; $A\simeq B$ means that $A\lesssim B$ and $B\lesssim A$ (analogously for $\simeq_a$), and $f(\xi)\sim g(\xi)$ as $\xi\rightarrow\xi_0$ means that $f(\xi)/g(\xi)\rightarrow1$ as $\xi\rightarrow\xi_0$.
\item For a set $B$, $B^c$ is its complement (in particular $\mathcal I^c=\R\setminus\mathcal I$) and $\mathbb I_B$ is its indicator function; if $B\subseteq\R$ is a Borel set, $\abs{B}$ is its Lebesgue measure. $\delta_x$ is the Dirac measure at $x$.
\item For measurable functions $f,g:\R\rightarrow\R$ such that $fg\in L^1(\R)$ we set $\scal{f}{g}_{L^2}:=\int_\R f(\xi)g(\xi)\,\mathd\xi$, also if $f$ or $g$ does not belong to $L^2$ (for instance, if $f\in L^2_{r,k}$ and $g\in C^\infty_c(\R)$).
\item For a Banach space $X$, $\mathcal L(X)$ is the space of bounded linear operators on $X$, endowed with the operator norm, $\Id$ (or $\Id_X$) is the identity and ${\rm Dom}(T)$ is the domain of an operator $T$; for a Hilbert space $H$, $\mathcal L_2(H)$ is the space of Hilbert--Schmidt operators on $H$ and ${\rm Tr}$ is the trace.
\item For a metric space $\mathcal X$, $\mathcal B(\mathcal X)$ is its Borel $\sigma$-algebra, $\mathcal P(\mathcal X)$ is the set of Borel probability measures on $\mathcal X$, and $B_b(\mathcal X)$ and $C_b(\mathcal X)$ are the spaces of bounded Borel and of bounded continuous functions $\Phi:\mathcal X\rightarrow\R$, endowed with the norm $\norm{\Phi}_\infty:=\sup_{x\in\mathcal X}\abs{\Phi(x)}$. For a Hilbert space $H$, $C^1_b(H)$ is the space of bounded and Fr\'echet differentiable functions $\Phi:H\rightarrow\R$ with bounded and continuous gradient $\nabla\Phi$, endowed with the norm $\norm{\Phi}_{C^1_b(H)}:=\norm{\Phi}_\infty+\sup_{x\in H}\norm{\nabla\Phi(x)}_H$.
\item For $B\subseteq\R$, $C_0(B)$ is the space of continuous functions on $B$ vanishing at the boundary of $B$ and at infinity, and $C^2_0(\R)$ is the space of the functions $f\in C^2(\R)$ such that $f,f',f''\in C_0(\R)$; $C^\infty_b(\R^n)$ is the space of smooth functions on $\R^n$ which are bounded together with all their derivatives.
\item For $\gamma\in(0,1)$, $C^{0,\gamma}(\R)$ is the space of bounded $\gamma$-H\"older continuous functions on $\R$, endowed with the norm
\[
\norm{f}_{C^{0,\gamma}(\R)}:=\sup_{\xi\in\R}\abs{f(\xi)}+\sup_{\xi\neq\xi'}\frac{\abs{f(\xi)-f(\xi')}}{\abs{\xi-\xi'}^{\gamma}},
\]
and, for a compact interval $[a,b]$, $C^\gamma([a,b])$ is the analogous space on $[a,b]$.
\item The Fourier transform of $f\in L^1(\R^n)$ is $\hat f(\xi):=\int_{\R^n}e^{-i\xi\cdot\eta}f(\eta)\,\mathd\eta$, $\xi\in\R^n$, and it is extended to tempered distributions by duality; with this normalization, $\norm{\hat f}^2_{L^2(\R^n)}=(2\pi)^n\norm{f}^2_{L^2(\R^n)}$ and $\widehat{f*g}=\hat f\hat g$. If $T$ is a linear operator such that $\widehat{Tf}=\psi\hat f$ for a function $\psi$, then $\psi$ is the (Fourier) symbol of $T$, and it is denoted by $\widehat T$.
\item For $s\geq0$, $H^s(\R)$ is the space of the functions $f\in L^2$ such that
\[
\norm{f}^2_{H^s(\R)}:=\frac{1}{2\pi}\int_\R(1+\abs{\xi}^2)^{s}\abs{\hat f(\xi)}^2\,\mathd\xi<+\infty,
\]
and, for $\varepsilon\in(0,1)$ and $p\geq1$, $W^{\varepsilon,p}(\R)$ is the fractional Sobolev space of the functions $f\in L^p$ such that
\[
\norm{f}^p_{W^{\varepsilon,p}(\R)}:=\norm{f}^p_{L^p}+\int_{\R^2}\frac{\abs{f(\xi)-f(\eta)}^p}{\abs{\xi-\eta}^{1+\varepsilon p}}\,\mathd\xi\,\mathd\eta<+\infty,
\]
with $W^{0,p}(\R):=L^p$ (see \cite{DiNezzaPalatucciValdinoci2012}).
\end{itemize}

\paragraph{Green functions.} If $D\subseteq\R^n$ is an open set and $T$ is a positive self-adjoint operator on $L^2(D)$, the Green function of $T$ is the (possibly distributional) integral kernel of $T^{-1}$, namely $G$ such that $T^{-1}f(\xi)=\int_DG(\xi,\eta)f(\eta)\,\mathd\eta$; for translation invariant operators on $\R^n$ we write $G(\xi-\eta)$ in place of $G(\xi,\eta)$.

\paragraph{L\'evy processes and complete Bernstein functions.} A real L\'evy process is a stochastic process starting from $0$, with independent and stationary increments and c\`adl\`ag trajectories (adding $\xi\in\R$ we obtain a L\'evy process starting from $\xi$), and it is symmetric if it has the same law as its opposite. A L\'evy measure on $\R$ is a measure $\nu$ on $\R\setminus\{0\}$ such that $\int(1\wedge\abs{z}^2)\nu(\mathd z)<+\infty$, and a L\'evy density is the density of a L\'evy measure with respect to the Lebesgue measure. A subordinator is a L\'evy process with nondecreasing trajectories, and its Laplace exponent is the function $\psi$ such that, for every $t\geq0$, the Laplace transform of the law of the subordinator at time $t$ is $q\mapsto e^{-t\psi(q)}$. A complete Bernstein function is a function $\psi:(0,+\infty)\rightarrow\R_+$ of the form
\[
\psi(q)=a+bq+\int_0^{+\infty}\big(1-e^{-tq}\big)\gamma(t)\,\mathd t,\qquad q>0,
\]
where $a,b\geq0$ and $\gamma$ is a completely monotonic function such that $\int_0^{+\infty}(t\wedge1)\gamma(t)\,\mathd t<+\infty$ (see the remark after \ref{H0}). A subordinator is complete if its Laplace exponent is a complete Bernstein function, and a subordinated Brownian motion is a Brownian motion evaluated at the random times given by an independent subordinator (see \cite{Applebaum2009,SchillingSongVondracek2012}).

\part{The non-local operator $H^{\beta}$}\label{part:H-beta}

One of the most interesting novelties of this paper is the consideration of SPDEs in which the differential part of the drift may also be nonlocal. Consequently, this section is devoted to the study of the properties of the following operator, which will form the basis for the results developed in the subsequent sections. We point out that, to the best of our knowledge, the operators $H^{\beta}$ defined in \ref{H0} below have not been studied in this generality in weighted $L^p$ spaces on $\R$ with integrable weights. For this reason, we believe that many of the results of this section are new and may be of independent interest.

\begin{enumerate}[start=0,label={{(H\arabic*})}]
  \item\label{H0}\emph{ \textbf{$H^{\beta}$ is a (possible) nonlocal operator} given by the following expression:
\begin{equation}\label{eq:Hbeta}
H^{\beta}f(\xi)=-\beta_0\Delta f(\xi)+\frac{1}{2}\int_{\mathbb{R}}\big[2 f(\xi)-f(\xi+\xi')-f(\xi-\xi')\big]\bar{\beta}(\xi')\,\mathd \xi',
\end{equation}
where $\beta_0\geq 0$ and $\bar{\beta}$ is the density of a L\'evy measure that is obtained by applying a complete subordinator to a one-dimensional Brownian motion, which is equivalent to asking that $\bar\beta$ is absolutely continuous with respect to the Lebesgue measure and admits the representation
\[
\bar{\beta}(\xi) = \int_0^{+\infty}\frac{e^{-\frac{|\xi|^2}{2t}}}{\sqrt{2 \pi t}}\gamma(t)\,\mathd t,
\]
where $\gamma(t)$ is a completely monotonic function such that
$\int_0^{+\infty} (t\wedge 1) \gamma(t) \,\mathd t <+\infty$.
We also ask that the function $\bar{\beta}$ satisfies the following conditions:}
\begin{itemize}
\item \emph{there exists $s_2 >0$ such that, for any $|\xi| \geq  1$,  \[  \bar{\beta}(\xi)  \lesssim \frac{1}{|\xi|^{1+2s_2}} ,\]}
\item \emph{if $\beta_0 =0$ there is $\frac{1}{2}< s_1 < 1$ such that
\begin{equation}\label{eq:comportamento-kernels_1}
0<\liminf_{|\xi| \rightarrow 0 } |\xi|^{1+2 s_1} \bar{\beta}(\xi) \leq   \limsup_{|\xi| \rightarrow 0 } |\xi|^{1+2 s_1} \bar{\beta}(\xi) <+\infty;
\end{equation}
if $\beta_0>0$ we set $s_1:=1$.}
\end{itemize}
\end{enumerate}

\begin{rmk}
We recall that a function $\gamma:(0,+\infty)\rightarrow\R$ is completely monotonic if $\gamma\in C^{\infty}((0,+\infty))$ and
\[
(-1)^n\frac{\mathd^n\gamma}{\mathd t^n}(t)\geq 0,\qquad t>0,\ n\in\N\cup\{0\},
\]
or equivalently, by Bernstein's theorem, if $\gamma$ is the Laplace transform of a positive measure on $[0,+\infty)$ (see \cite[Definition 1.3 and Theorem 1.4]{SchillingSongVondracek2012}). In \ref{H0}, $\gamma(t)\,\mathd t$ is the L\'evy measure of the complete subordinator, which explains the integrability condition on $\gamma$.
\end{rmk}

\begin{rmk}
There are many significant examples of operators $H^{\beta}$ that satisfy \ref{H0}:
\begin{itemize}
\item For $\beta_0,\beta_1\geq 0$ with $\beta_0+\beta_1>0$, $\lambda\geq 0$ and $s\in(0,1)$, the operator
\[
H^{\beta}f(\xi)=-\beta_0\Delta f(\xi)+\frac{\beta_1}{2}\int_{\R}\frac{2f(\xi)-f(\xi+\xi')-f(\xi-\xi')}{(\lambda+|\xi'|)^{1+2s}}\,\mathd \xi',
\]
where, if $\beta_0=0$, we assume $\lambda=0$ and $s\in\left(1/2,1\right)$. It satisfies \ref{H0} with $s_2=s$ (and $s_1=s$ if $\beta_0=0$), since $\bar{\beta}(\xi)=\beta_1(\lambda+|\xi|)^{-1-2s}$ admits the representation in \ref{H0} with $\gamma(t)=c\int_0^{+\infty}v^{1+2s}e^{-\lambda v-\frac{v^2t}{2}}\mathd v$, where $c$ is a suitable constant. For $\lambda=0$, the non-local part, up to a multiplicative constant, is the fractional Laplacian of index $s$, while for $\lambda>0$ it corresponds to the quadratic pair interaction considered in \cite{OsadaSpohn1999} for Gibbs measures relative to Brownian motion. We underline that if $\beta_0=0$ we need $\lambda=0$ and $s=s_1>1/2$ to work in the non-singular regime in which the stochastic convolution associated with \eqref{eq:LangevinSPDE1} takes values in continuous functions (see Lemma \ref{lemma:stochconv}). Moreover, if $\beta_0>0$ and $\lambda>0$, then $\bar{\beta}\in L^1(\R)$ and the nonlocal part of $H^{\beta}$ is a bounded operator on $L^p(\R)$ for every $p\in[1,+\infty]$, hence it has no regularizing effect.
\item The massive fractional Laplacian
\[
(-\Delta+c^{\frac{1}{s}})^{s}-c
\]
for any $s\in \left(1/2,1\right)$ and $c>0$.
\item The operator associated with the superposition of Ornstein--Uhlenbeck processes:
\[
H^{\beta}f=-\beta_0\Delta f+\frac{1}{2}\int_{\R} \big(2f(\cdot)-f(\cdot-\xi)-f(\cdot+\xi)\big) \left(\sum_{i=1}^{n-1}d_i^2 e^{-k_i |\xi|} \right)\mathd \xi,
\]
where $d_1,\ldots,d_{n-1},k_1,\ldots,k_{n-1}>0$. Another representation of this operator is given by
\[
H^{\beta}:=\left(\sum_{i=1}^n \frac{c_i^2}{-\Delta+m_i^2}\right)^{-1}-\left(\sum_{i=1}^n \frac{c_i^2}{m_i^2}\right)^{-1},
\]
where $c_1,\dots,c_n,m_1,\dots,m_n>0$ and $m_1 < k_1 < m_2 <\cdots < k_{n-1} < m_n$. A further representation of $H^{\beta}$ is given by
\[
H^{\beta}=\beta_0 \frac{\prod_{i=1}^n (-\Delta +m_i^2)}{\prod_{i=1}^{n-1} (-\Delta +k_i^2)}-\beta_0\frac{\prod_{i=1}^n m_i^2}{\prod_{i=1}^{n-1}k_i^2}.
\]
\end{itemize}
More generally, if we ignore the conditions on the behavior of $\bar{\beta}$ at $0$ and infinity, the operators of the form $H^{\beta}$ include all the operators on $\R$ associated with a Gaussian measure which are invariant with respect to translations and reflection positive (see the next Section \ref{appendix:reflectionpositivity} for a more detailed discussion).
\end{rmk}

\section{Analytic properties of $H^{\beta}$}\label{sect:analitic-Hbeta}

Throughout this section, we assume that \ref{H0} holds. Below, we collect some estimates for $\overline{\beta}$ that follow directly from \ref{H0} and will be used in the next subsections. First of all, we have the following integrability
\begin{equation}\label{eq:levynear}
\Theta_\beta:=\int_{\lvert z\rvert\leq  1}\lvert z\rvert^{2}\bar\beta(z)\mathd z<+\infty.
\end{equation}
There exists $C_\beta\geq1$ such that
\begin{equation}\label{eq:levygrowth}
\bar\beta(z)\leq   C_\beta\lvert z\rvert^{-1-2s_2}\qquad\text{if }\lvert z\rvert\geq 1,
\end{equation}
and, if in addition $\beta_0=0$, then
\begin{equation}\label{eq:levybilateral}
C_\beta^{-1}\lvert z\rvert^{-1-2s_1}\leq\bar\beta(z)\leq C_\beta\lvert z\rvert^{-1-2s_1}\qquad\text{if }\lvert z\rvert\leq1.
\end{equation}

In this section, we are going to study in some suitable weighted $L^p$ spaces the operator $H^{\beta}+m^2$, where $m\geq 0$ is an arbitrary constant. From here on we set $L^p:=L^p(\R)$, for every $p\geq 1$. Given a continuous and positive $\rho: \R\rightarrow\R$, $p\geq 1$ and $B\subseteq\R$, we introduce the (separable) weighted Banach spaces
\[
L^p_\rho(B)= \left\lbrace g \in L^1_{\loc}(B) : \int_{B} |g(\xi)|^p \rho(\xi)\mathd \xi <+\infty \right\rbrace, \qquad \norm{g}^p_{L^p_\rho(B)}:=\int_{B} |g(\xi)|^p \rho(\xi)\mathd \xi.
\] 
In the present work we consider weights of the form
\begin{equation}\label{eq:weight}
\rho_{r,k}(\xi)=(1+k|\xi|^2)^{-r/2},\qquad r>1,\,k\geq 0,
\end{equation}
and we set
\begin{align}
L^p_{r,k}(B)&:=L^p_{\rho_{r,k}}(B),\qquad r>1,\,k\geq 0,\ p\geq 1,\notag\\
L^p_{r,k}&:=L^p_{\rho_{r,k}},\qquad r>1,\,k\geq 0,\ p\geq 1,\label{eq:pes-wei}
\end{align}
we use the convention $L^p_{r,0}(B)=L^p(B)$.

\subsection{Operator and semigroup in weighted spaces}\label{section:weights}

We begin by collecting some properties of the weight given by \eqref{eq:weight} that will be used throughout the paper.

\begin{lemm}\label{lemma:weightproperties}
For every $r>1$ and $k>0$ and $\xi,\eta\in\R$, we have
\begin{align}
&\left\lvert\frac{\rho_{r,k}'(\xi)}{\rho_{r,k}(\xi)}\right\rvert\leq  \frac r2\,\sqrt k,\label{eq:rhoderivatives}\\[2pt]
&\left\lvert\frac{\rho_{r,k}''(\xi)}{\rho_{r,k}(\xi)}\right\rvert\leq  \frac{r(r+6)}{4}k ,\label{eq:rhoderivatives2}\\[2pt]
&\rho_{r,k}(\xi)\leq   4^{r/2}\rho_{r,k}(\eta)\qquad  \text{if } \lvert\eta\rvert\leq  2\lvert\xi\rvert,\label{eq:rhodoubling}\\[2pt]
&e^{-r/2}\leq  \frac{\rho_{r,k}(\xi)}{\rho_{r,k}(\eta)}\leq   e^{r/2}\qquad  \text{if } \lvert\xi-\eta\rvert\leq   k^{-1/2} .\label{eq:rhocomparable}
\end{align}
Moreover $\rho_{r,k}\in L^1(\R)$ and
\begin{equation}\label{eq:rhoL1}
\norm{\rho_{r,k}}_{L^1(\R)}=c_r\,k^{-1/2},
\qquad
c_r:=\int_\R\big(1+u^2\big)^{-r/2}\,\mathd u<+\infty ,
\end{equation}
and, for every $k'>0$,
\begin{equation}\label{eq:rhoequivk}
\min\left\{1,\frac{k'}{k}\right\}^{r/2}\leq  \frac{\rho_{r,k}(\xi)}{\rho_{r,k'}(\xi)}\leq  \max\left\{1,\frac{k'}{k}\right\}^{r/2},
\end{equation}
so that $L^p_{r,k}=L^p_{r,k'}$ up to an equivalent renorming, for every $p\geq 1$.
\end{lemm}

\begin{proof}
Let $r,k>0$. Set $\rho:=\rho_{r,k}$ for simplicity. Since $[\rho'/\rho](\xi)=-rk\xi(1+k\xi^2)^{-1}$ and $2\sqrt k\lvert\xi\rvert\leq  1+k\xi^2$ we deduce \eqref{eq:rhoderivatives}. Differentiating once more, we get
\[
\left[\frac{\rho''}{\rho}\right](\xi)=-\frac{rk}{1+k\xi^2}+\frac{r(r+2)k^2\xi^2}{(1+k\xi^2)^2},
\]
and $k\xi^2(1+k\xi^2)^{-2}\leq  \frac14$ yields \eqref{eq:rhoderivatives2}. Estimate \eqref{eq:rhodoubling} follows from
\[
\frac{\rho(\xi)}{\rho(\eta)}=\left(\frac{1+k\lvert\eta\rvert^2}{1+k\lvert\xi\rvert^2}\right)^{r/2}\leq  \left(\frac{1+4k\lvert\xi\rvert^2}{1+k\lvert\xi\rvert^2}\right)^{r/2}\leq  4^{r/2},
\]
while \eqref{eq:rhocomparable} follows by integrating \eqref{eq:rhoderivatives}, which gives $\lvert\log\rho(\xi)-\log\rho(\eta)\rvert\leq  \frac r2\sqrt k\,\lvert\xi-\eta\rvert$. Identity \eqref{eq:rhoL1} is the change of variables $u=\sqrt k\,\xi$, and \eqref{eq:rhoequivk} is immediate from
$\rho_{r,k}/\rho_{r,k'}=\big((1+k'\lvert\xi\rvert^2)/(1+k\lvert\xi\rvert^2)\big)^{r/2}$.
\end{proof}

We now prove an elementary consequence of the translation invariance of the family \eqref{eq:weight}.

\begin{lemm}\label{lemma:exponentchange}
Let $r>1$ and $r_0>0$. Then there exists $C_0=C_0(r,r_0)$ such that for every measurable $g:\R\rightarrow\R$, $k>0$ and $p\geq1$ we have
\begin{equation}\label{eq:exponentchange}
\norm{g}^p_{L^p_{r,k}}\leq C_0\sup_{\tau\in\R}\norm{g}^p_{L^p_{\rho^\tau_{r_0,k}}},
\end{equation}
where $\rho^\tau_{r,k}=\rho_{r,k}(\cdot-\tau)$.
\end{lemm}

\begin{proof}
If $\abs{\tau-\xi}\leq k^{-1/2}$ then $1+k\abs{\tau}^2\leq4(1+k\abs{\xi}^2)$, hence $\rho_{r,k}(\tau)\geq4^{-r/2}\rho_{r,k}(\xi)$ so that
\begin{align*}
\big(\rho_{r,k}*\rho_{r_0,k}\big)(\xi)&=\int_\R\rho_{r,k}(\tau)\rho_{r_0,k}(\xi-\tau)\mathd\tau
\\
&\geq4^{-r/2}\rho_{r,k}(\xi)\int_{\abs{z}\leq k^{-1/2}}\rho_{r_0,k}(z)\mathd z\geq c\,k^{-1/2}\rho_{r,k}(\xi),
\end{align*}
where $c=c(r,r_0)>0$, that is $\rho_{r,k}\leq c^{-1}k^{1/2}\big(\rho_{r_0,k}*\rho_{r,k}\big)$. Therefore, by the Fubini theorem,
\[
\int_\R\rho_{r,k}\abs{g}^p\,\mathd\xi
\leq c^{-1}k^{1/2}\int_\R\rho_{r,k}(\tau)\left(\int_\R\rho^\tau_{r_0,k}(\xi)\abs{g(\xi)}^p\mathd\xi\right)\mathd\tau
\leq c^{-1}k^{1/2}\norm{\rho_{r,k}}_{L^1(\R)}\sup_{\tau\in\R}\norm{g}^p_{L^p_{\rho^\tau_{r_0,k}}},
\]
and \eqref{eq:exponentchange} follows from \eqref{eq:rhoL1}.
\end{proof}

We are now in a position to prove the main result of this subsection, which will play a fundamental role in establishing several properties of the semigroup associated with $H^\beta+m^2$ where $m>0$ is an arbitrary constant.

Throughout this subsection and the next one, we denote by
\begin{equation}\label{eq:Hbeta-nl-only}
H^{\beta}_{nl}f(\xi):=\frac{1}{2}\int_{\R}\big[2f(\xi)-f(\xi+z)-f(\xi-z)\big]\bar\beta(z)\mathd z
\end{equation}
the nonlocal part of the operator $H^{\beta}$ introduced in \ref{H0}, so that $H^{\beta}=-\beta_0\Delta+H^{\beta}_{nl}$; the integral in \eqref{eq:Hbeta-nl-only} is understood as a principal value whenever necessary.

Note that $(H^{\beta}+m^2)h$ is not compactly supported even if $h$ is compactly supported, indeed outside the support of $h$ one has $[(H^{\beta}+m^2)h](\xi)=-\int_\R h(\eta)\bar\beta(\xi-\eta)\mathd\eta$, which decays only as fast as $\bar\beta$. The following estimate quantifies this decay and will be used repeatedly in the sequel.

\begin{lemm}\label{lemma:testdecay}
Assume \ref{H0} holds and let $1<r<1+2s_2$ and $k>0$. Then, for every $h\in\mathcal S(\R)$ there exists $C_h:=C_h(r,k,m,\bar\beta)>0$, such that
\begin{equation}\label{eq:testdecay}
\big\lvert\big[(H^{\beta}+m^2)h\big](\xi)\big\rvert\leq C_h\,\rho_{r,k}(\xi),\qquad \xi\in\R .
\end{equation}
In particular
\begin{equation}\label{eq:testdecay-L2}
\rho_{r,k}^{-1/2}\,(H^{\beta}+m^2)h\in L^2
\qquad\text{and}\qquad
\rho_{r,k}^{-1}\,(H^{\beta}+m^2)h\in L^2_{r,k}.
\end{equation}

\end{lemm}

\begin{proof}
Fix $\xi\in\R$. We note that $1+k\lvert\xi\rvert^2\leq\max\{1,k\}(1+\lvert\xi\rvert^2)$, so that $(1+\lvert\xi\rvert^2)^{-r/2}\leq\max\{1,k\}^{r/2}\rho_{r,k}(\xi)$; since $r<1+2s_2$, it is therefore enough to show that
\[
\big\lvert\big[(H^{\beta}+m^2)h\big](\xi)\big\rvert\leq K_h(1+\lvert\xi\rvert^2)^{(-1-2s_2)/2},\qquad\xi\in\R,
\]
for some $K_h>0$. For $N\in\N$ set
\[
C_N(h):=\sup_{\xi\in\R}\big(1+\lvert\xi\rvert^2\big)^{N/2}\big(\lvert h(\xi)\rvert+\lvert h''(\xi)\rvert\big),
\qquad
B_\beta:=\int_{\lvert z\rvert\geq1}\bar\beta(z)\mathd z ,
\]
which are finite since $h\in\mathcal S(\R)$ and $\bar\beta$ is a L\'evy density. Fix $N\in\N$ with $N\geq1+2s_2$. The local part is estimated directly by
\begin{equation}\label{eq:parte-loc}
\big\lvert-\beta_0\Delta h(\xi)+m^2h(\xi)\big\rvert\leq(\beta_0+m^2)\,C_N(h)\big(1+\lvert\xi\rvert^2\big)^{-N/2}.
\end{equation}

For the nonlocal part we split the integral in \eqref{eq:Hbeta-nl-only} into the regions $\{\lvert z\rvert\leq1\}$ and $\{\lvert z\rvert\geq1\}$. In the first one, a second-order Taylor expansion gives
\[
\big\lvert2h(\xi)-h(\xi+z)-h(\xi-z)\big\rvert\leq\lvert z\rvert^2\sup_{\lvert w-\xi\rvert\leq1}\lvert h''(w)\rvert
\leq 5^{N/2}\,C_N(h)\,\lvert z\rvert^2\big(1+\lvert\xi\rvert^2\big)^{-N/2},
\]
where we used that $\lvert w-\xi\rvert\leq1$ implies $1+\lvert w\rvert^2\geq\frac15\big(1+\lvert\xi\rvert^2\big)$. Hence, by \eqref{eq:levynear}, we have
\begin{equation}\label{eq:r1}
\frac12\int_{\lvert z\rvert\leq1}\big\lvert2h(\xi)-h(\xi+z)-h(\xi-z)\big\rvert\bar\beta(z)\mathd z
\leq\frac{5^{N/2}}{2}\,\Theta_\beta\,C_N(h)\big(1+\lvert\xi\rvert^2\big)^{-N/2} .
\end{equation}
For the second region we have
\begin{align*}
\frac12\int_{\lvert z\rvert>1}\big\lvert2h(\xi)-h(\xi+z)-h(\xi-z)\big\rvert\bar\beta(z)\mathd z
\leq \frac12\int_{\lvert z\rvert>1}\bigg[2\lvert h(\xi)\rvert+\lvert h(\xi+z)\rvert+\lvert h(\xi-z)\rvert\bigg]\bar\beta(z)\mathd z\\
\leq \lvert h(\xi)\rvert\int_{\lvert z\rvert\geq1}\bar\beta(z)\mathd z+\frac12\int_{\lvert z\rvert>1}\bigg[\lvert h(\xi+z)\rvert+\lvert h(\xi-z)\rvert\bigg]\bar\beta(z)\mathd z\\
\leq B_\beta\,C_N(h)\big(1+\lvert\xi\rvert^2\big)^{-N/2} +\frac12\int_{\lvert z\rvert>1}\bigg[\lvert h(\xi+z)\rvert+\lvert h(\xi-z)\rvert\bigg]\bar\beta(z)\mathd z.
\end{align*}
For the remaining two terms, by the change of variable $\eta:=\xi\pm z$, and, by the symmetry of $\bar\beta$, we obtain
\begin{align}\label{eq:r2-pre}
\frac12\int_{\lvert z\rvert>1}\big\lvert2h(\xi)-h(\xi+z)-h(\xi-z)\big\rvert\bar\beta(z)\mathd z
&\leq B_\beta C_N(h)\big(1+\lvert\xi\rvert^2\big)^{-N/2}\notag\\
&+\int_{\lvert\xi-\eta\rvert\geq1}\lvert h(\eta)\rvert\bar\beta(\xi-\eta)\mathd\eta.
\end{align}
To estimate the second term, we consider two cases: $\lvert\eta\rvert\leq\frac12\lvert\xi\rvert$ and $\lvert\eta\rvert\geq\frac12\lvert\xi\rvert$. In the first case $\lvert\xi-\eta\rvert\geq\lvert\xi\rvert-\lvert\eta\rvert\geq\frac12\lvert\xi\rvert$ and $\lvert\xi-\eta\rvert\geq1$, so that $\bar\beta(\xi-\eta)\leq C_\beta\max\{1,\tfrac12\lvert\xi\rvert\}^{-1-2s_2}$ by \eqref{eq:levygrowth}; since $\max\{1,\tfrac12\lvert\xi\rvert\}\geq\frac{1}{2\sqrt2}\big(1+\lvert\xi\rvert^2\big)^{1/2}$, we obtain
\[
\int_{\{\lvert\xi-\eta\rvert\geq1\}\cap\{\lvert\eta\rvert\leq\frac12\lvert\xi\rvert\}}\lvert h(\eta)\rvert\,\bar\beta(\xi-\eta)\mathd\eta
\leq\big(2\sqrt2\big)^{1+2s_2}C_\beta\norm{h}_{L^1(\R)}\big(1+\lvert\xi\rvert^2\big)^{-\frac{1+2s_2}{2}} .
\]
In the second case $1+\lvert\eta\rvert^2\geq\frac14\big(1+\lvert\xi\rvert^2\big)$, hence $\lvert h(\eta)\rvert\leq4^{N/2}C_N(h)\big(1+\lvert\xi\rvert^2\big)^{-N/2}$ and
\[
\int_{\{\lvert\xi-\eta\rvert\geq1\}\cap\{\lvert\eta\rvert\geq\frac12\lvert\xi\rvert\}}\lvert h(\eta)\rvert\bar\beta(\xi-\eta)\mathd\eta
\leq4^{N/2}B_\beta C_N(h)\big(1+\lvert\xi\rvert^2\big)^{-N/2} .
\]
Hence by \eqref{eq:r2-pre} we obtain
\begin{align}
    \frac12\int_{\lvert z\rvert>1}\big\lvert2h(\xi)-h(\xi+z)-h(\xi-z)\big\rvert\bar\beta(z)\mathd z
&\leq B_\beta C_N(h)\big(1+\lvert\xi\rvert^2\big)^{-N/2}\notag\\
&+\big(2\sqrt2\big)^{1+2s_2}C_\beta\norm{h}_{L^1(\R)}\big(1+\lvert\xi\rvert^2\big)^{-\frac{1+2s_2}{2}}\notag\\
&+4^{N/2}B_\beta C_N(h)\big(1+\lvert\xi\rvert^2\big)^{-N/2}.\label{eq:r2}
\end{align}
Since $N\geq1+2s_2$, we have $\big(1+\lvert\xi\rvert^2\big)^{-N/2}\leq\big(1+\lvert\xi\rvert^2\big)^{-\frac{1+2s_2}{2}}$, so that, collecting \eqref{eq:parte-loc}, \eqref{eq:r1} and \eqref{eq:r2}, we obtain
\[
\big\lvert\big[(H^{\beta}+m^2)h\big](\xi)\big\rvert\leq K_h\big(1+\lvert\xi\rvert^2\big)^{-\frac{1+2s_2}{2}},
\]
where
\[
K_h:=\left(\beta_0+m^2+\frac{5^{N/2}}{2}\Theta_\beta+\big(1+4^{N/2}\big)B_\beta\right)C_N(h)+\big(2\sqrt2\big)^{1+2s_2}C_\beta\norm{h}_{L^1(\R)} ,
\]
and \eqref{eq:testdecay} follows with $C_h:=\max\{1,k\}^{r/2}K_h$. Finally, by \eqref{eq:rhoL1} and \eqref{eq:testdecay}, we have
\[
\int_\R\big\lvert\rho_{r,k}^{-1/2}(H^{\beta}+m^2)h\big\rvert^2\mathd\xi
=\int_\R\big\lvert\rho_{r,k}^{-1}(H^{\beta}+m^2)h\big\rvert^2\rho_{r,k}\mathd\xi
\leq C_h^2\norm{\rho_{r,k}}_{L^1(\R)}<+\infty,
\]
which gives \eqref{eq:testdecay-L2}.
\end{proof}

\begin{rmk}\label{rmk:weighttranslation}
Since $H^{\beta}$ commutes with translations, for every $\tau\in\R$ the weight $\rho^{\tau}_{r,k}:=\rho_{r,k}(\cdot-\tau)$ satisfies $H^{\beta}\rho^{\tau}_{r,k}=\big(H^{\beta}\rho_{r,k}\big)(\cdot-\tau)$. Also the estimates \eqref{eq:rhoderivatives}--\eqref{eq:rhoL1} are invariant under translations, so that every result of this subsection and of Subsection \ref{section:commutator} holds for $\rho^{\tau}_{r,k}$ with the same constants (uniformly in $\tau\in\R$). 
\end{rmk}

\begin{prop}\label{prop:levy-weight-rho}
Assume \ref{H0} holds. For every $1<r<1+2s_2$, $k\in(0,1]$ and $0<\sigma<\min\{s_2,1\}$ there exists a constant $C_1=C_1(r,s_2,\sigma,C_\beta,\Theta_\beta)>0$ such that
\begin{equation}\label{eq:weightratio}
\sup_{\xi\in\R}\left\lvert\frac{H^{\beta}_{nl}\rho_{r,k}(\xi)}{\rho_{r,k}(\xi)}\right\rvert\leq   C_1k^{\sigma}.
\end{equation}
\end{prop}

\begin{proof}
Let $1<r<1+2s_2$, $k\in(0,1]$ and $\xi\in\R$. Set $\rho:=\rho_{r,k}$ and split the integral in \eqref{eq:Hbeta-nl-only} into the regions $\{\lvert z\rvert\leq   k^{-1/2}\}$ and $\{\lvert z\rvert\geq  k^{-1/2}\}$ respectively.\\
If $\lvert z\rvert\leq   k^{-1/2}$, then
\[
2\rho(\xi)-\rho(\xi+z)-\rho(\xi-z)=-\int_0^1(1-t)\big[\rho''(\xi+tz)+\rho''(\xi-tz)\big]z^2\,\mathd t,
\]
so that by \eqref{eq:rhoderivatives2} and \eqref{eq:rhocomparable} we obtain $\lvert2\rho(\xi)-\rho(\xi+z)-\rho(\xi-z)\rvert\lesssim_r k\rho(\xi)\lvert z\rvert^2$. Hence, by \eqref{eq:levynear} and \eqref{eq:levygrowth}, we have
\begin{align}
\left\lvert \frac{1}{2}\int_{\{\lvert z\rvert\leq   k^{-1/2}\}}\big[2\rho(\xi)-\rho(\xi+z)-\rho(\xi-z)\big]\bar\beta(z)\mathd z\right\rvert
&\leq
e^{r}k\rho(\xi)\int_{\lvert z\rvert\leq   k^{-1/2}}\lvert z\rvert^{2}\bar\beta(z)\mathd z\notag\\
&\leq
e^{r}k\rho(\xi)\left(\Theta_\beta+C_\beta\int_{1\leq  \lvert z\rvert\leq   k^{-1/2}}\lvert z\rvert^{1-2s_2}\mathd z\right)\notag\\
&\leq e^{r}\left(\Theta_\beta+\frac{C_\beta}{1-\sigma}\right)k^{\sigma}\rho(\xi),\label{eq:S<}
\end{align}
where $\Theta_\beta$ and $C_\beta$ are the constants of \eqref{eq:levynear} and \eqref{eq:levygrowth}, and in the last inequality we used $k\leq k^{\sigma}$ together with $\lvert z\rvert^{1-2s_2}\leq\lvert z\rvert^{1-2\sigma}$ for $\lvert z\rvert\geq1$, which gives
\[
k\int_{1\leq\lvert z\rvert\leq k^{-1/2}}\lvert z\rvert^{1-2s_2}\mathd z\leq 2k\int_1^{k^{-1/2}}t^{1-2\sigma}\mathd t\leq\frac{k^{\sigma}}{1-\sigma}.
\]
Now we study the case $\lvert z\rvert\geq k^{-1/2}$. Here we use $\lvert2\rho(\xi)-\rho(\xi+z)-\rho(\xi-z)\rvert\leq 2\rho(\xi)+\rho(\xi+z)+\rho(\xi-z)$, so that
\begin{align*}
\left\lvert\frac12\int_{\{\lvert z\rvert\geq k^{-1/2}\}}\big[2\rho(\xi)-\rho(\xi+z)-\rho(\xi-z)\big]\bar\beta(z)\mathd z\right\rvert
&\leq
\rho(\xi)\int_{\lvert z\rvert\geq k^{-1/2}}\bar\beta(z)\mathd z\\
&+\frac12\int_{\lvert z\rvert\geq k^{-1/2}}\rho(\xi+ z)\bar\beta(z)\mathd z \\
&+\frac12\int_{\lvert z\rvert\geq k^{-1/2}}\rho(\xi-z)\bar\beta(z)\mathd z .
\end{align*}
In the last two integrals of the right hand side we use the changes of variables $\eta:=\xi+z$ and $\eta:=\xi-z$ respectively; since $\bar\beta$ is even we obtain
\begin{align*}
\left\lvert\frac12\int_{\{\lvert z\rvert\geq k^{-1/2}\}}\big[2\rho(\xi)-\rho(\xi+z)-\rho(\xi-z)\big]\bar\beta(z)\mathd z\right\rvert
&\leq
\rho(\xi)\int_{\lvert z\rvert\geq k^{-1/2}}\bar\beta(z)\mathd z\\
&+\int_{\lvert\xi-\eta\rvert\geq k^{-1/2}}\rho(\eta)\,\bar\beta(\xi-\eta)\mathd\eta .
\end{align*}
We split the last integral over the regions $\lvert\eta\rvert\geq\frac12\lvert\xi\rvert$ and $\lvert\eta\rvert<\frac12\lvert\xi\rvert$ and we set
\begin{align*}
I_1(\xi)&:=\rho(\xi)\int_{\lvert z\rvert\geq k^{-1/2}}\bar\beta(z)\,\mathd z,\\
I_2(\xi)&:=\int_{\{\lvert\xi-\eta\rvert\geq k^{-1/2}\}\cap\{\lvert\eta\rvert\geq\frac12\lvert\xi\rvert\}}\rho(\eta)\,\bar\beta(\xi-\eta)\,\mathd\eta,\\
I_3(\xi)&:=\int_{\{\lvert\xi-\eta\rvert\geq k^{-1/2}\}\cap\{\lvert\eta\rvert<\frac12\lvert\xi\rvert\}}\rho(\eta)\,\bar\beta(\xi-\eta)\,\mathd\eta .
\end{align*}
We focus on $I_1$. By \eqref{eq:levygrowth}, since $k^{-1/2}\geq1$, we have
\begin{equation}\label{eq:I_1}
I_1(\xi)\leq C_\beta\rho(\xi)\int_{\lvert z\rvert\geq k^{-1/2}}\lvert z\rvert^{-1-2s_2}\mathd z=\frac{C_\beta}{s_2}k^{s_2}\rho(\xi).
\end{equation}
We focus on $I_2$. Since $\lvert\eta\rvert\geq\frac12\lvert\xi\rvert$, then $\rho(\eta)\leq 2^{r}\rho(\xi)$ by \eqref{eq:rhodoubling}, we get
\begin{equation}\label{eq:I_2}
I_2(\xi)\leq 2^{r}\frac{C_\beta}{s_2}k^{s_2}\rho(\xi).
\end{equation}
We focus on $I_3$. Since $\lvert\eta\rvert<\frac12\lvert\xi\rvert$, then $\lvert\xi-\eta\rvert\geq\lvert\xi\rvert-\lvert\eta\rvert\geq\frac12\lvert\xi\rvert$ and therefore $\lvert\xi-\eta\rvert\geq M(\xi):=\max\left\{k^{-1/2},\frac12\lvert\xi\rvert\right\}\geq1$. Since $t\mapsto t^{-1-2s_2}$ is decreasing, \eqref{eq:levygrowth} gives $\bar\beta(\xi-\eta)\leq C_\beta M(\xi)^{-1-2s_2}$ on $\{\lvert\xi-\eta\rvert\geq k^{-1/2}\}\cap\{\lvert\eta\rvert<\frac12\lvert\xi\rvert\}$, by \eqref{eq:rhoL1} we have
\[
I_3(\xi)\leq
C_\beta\,M(\xi)^{-1-2s_2}\int_\R\rho(\eta)\mathd\eta
=C_\beta c_r\,k^{-1/2}M(\xi)^{-1-2s_2}.
\]
We claim that
\begin{equation}\label{eq:I3claim}
k^{-1/2}M(\xi)^{-1-2s_2}\leq 2^{1+2s_2+\frac r2}k^{s_2}\,\rho(\xi),\qquad \xi\in\R .
\end{equation}
Indeed, if $\lvert\xi\rvert\leq k^{-1/2}$, then $M(\xi)=k^{-1/2}$, so that the left-hand side of \eqref{eq:I3claim} equals $k^{s_2}$, while $1+k\lvert\xi\rvert^2\leq2$ gives $\rho(\xi)\geq2^{-r/2}$ and hence \eqref{eq:I3claim} holds for $\lvert\xi\rvert\leq k^{-1/2}$. If instead $\lvert\xi\rvert\geq k^{-1/2}$, then $M(\xi)\geq\frac12\lvert\xi\rvert$ and $1+k\lvert\xi\rvert^2\leq2k\lvert\xi\rvert^2$, so that $\rho(\xi)^{-1}\leq2^{r/2}k^{r/2}\lvert\xi\rvert^{r}$ and
\[
\frac{k^{-1/2}M(\xi)^{-1-2s_2}}{\rho(\xi)}
\leq
2^{1+2s_2+\frac r2}\,k^{\frac{r-1}{2}}\,\lvert\xi\rvert^{\,r-1-2s_2}
\leq
2^{1+2s_2+\frac r2}\,k^{\frac{r-1}{2}}\,k^{-\frac{r-1-2s_2}{2}}
=2^{1+2s_2+\frac r2}\,k^{s_2},
\]
where we used $r-1-2s_2<0$ together with $\lvert\xi\rvert\geq k^{-1/2}$. This proves \eqref{eq:I3claim}, and therefore
\begin{equation}\label{eq:I_3}
I_3(\xi)\leq 2^{1+2s_2+\frac r2}C_\beta c_r\,k^{s_2}\rho(\xi).
\end{equation}
By \eqref{eq:I_1}, \eqref{eq:I_2} and \eqref{eq:I_3} we obtain
\begin{equation}\label{eq:S>}
\left\lvert\frac12\int_{\{\lvert z\rvert\geq k^{-1/2}\}}\big[2\rho(\xi)-\rho(\xi+z)-\rho(\xi-z)\big]\bar\beta(z)\mathd z\right\rvert
\leq C_\beta\left(\frac{1+2^{r}}{s_2}+2^{1+2s_2+\frac r2}c_r\right)k^{s_2}\rho(\xi).
\end{equation}
Combining \eqref{eq:S<} and \eqref{eq:S>} we obtain \eqref{eq:weightratio}.

\end{proof}

For the local part $\Delta$ of $H^{\beta}$ an explicit computation is available.

\begin{prop}\label{prop:laplacian-weight}
For every $r,k>0$ we have
\begin{equation}\label{eq:laplacianweight}
\frac{\Delta\rho_{r,k}(\xi)}{\rho_{r,k}(\xi)}=\frac{kr\big((r+1)k\lvert\xi\rvert^2-1\big)}{(1+k\lvert\xi\rvert^2)^2},
\qquad
\sup_{\xi\in\R}\frac{\Delta\rho_{r,k}(\xi)}{\rho_{r,k}(\xi)}=\frac{r(r+1)^2}{4(r+2)}\,k .
\end{equation}
\end{prop}
\begin{proof}
The identity follows from the expression of $\rho_{r,k}''/\rho_{r,k}$ in the proof of Lemma \ref{lemma:weightproperties}. Setting $u:=k\lvert\xi\rvert^2\geq0$, the function $u\mapsto((r+1)u-1)(1+u)^{-2}$ has the maximum at $u=\frac{r+3}{r+1}$, where it is equal to $\frac{(r+1)^2}{4(r+2)}$.
\end{proof}

Since $H^\beta\rho_{r,k}=-\beta_0\Delta\rho_{r,k}+H^\beta_{nl}\rho_{r,k}$, combining Propositions \ref{prop:levy-weight-rho} and \ref{prop:laplacian-weight} we obtain the estimate for the complete operator $H^\beta$.

\begin{cor}\label{cor:omega-rk}
Assume \ref{H0} holds. For every $1<r<1+2s_2$, $k\in(0,1]$ and $0<\sigma<\min\{s_2,1\}$ we have
\begin{equation}\label{eq:omega-rk-bound}
\omega_{r,k}:=\sup_{\xi\in\R}\left(-\frac{H^{\beta}\rho_{r,k}(\xi)}{\rho_{r,k}(\xi)}\right)\leq\beta_0\frac{r(r+1)^2}{4(r+2)}k+C_1k^{\sigma}<+\infty ,
\end{equation}
where $C_1$ is the constant given by Proposition \ref{prop:levy-weight-rho}. In particular, for every $1<r<1+2s_2$ and $k\in(0,1]$, we have
\begin{equation}\label{eq:excessive}
-H^{\beta}\rho_{r,k}(\xi)\leq  \omega_{r,k}\,\rho_{r,k}(\xi),\qquad\xi\in\R,
\end{equation}
where $\omega_{r,k}\rightarrow0$ as $k\rightarrow0^+$.
\end{cor}

\begin{rmk}\label{rmk:pure-laplacian-nobound}
The upper bound $r<1+2s_2$ in Proposition \ref{prop:levy-weight-rho} is due to the nonlocal part of $H^{\beta}$: if $\bar\beta\equiv0$, then $\omega_{r,k}=\beta_0\frac{r(r+1)^2}{4(r+2)}k$ for every $r>0$, with no upper restriction on $r$.
\end{rmk}

\begin{rmk}\label{rmk:exponentialweights}
If $\bar\beta$ decays exponentially at infinity, or vanishes identically, then the proof of Proposition \ref{prop:levy-weight-rho} applies similarly to exponential weights of the type
\[
\tilde\rho_{k}(\xi)=\exp\left(-(1+k|\xi|^2)^{1/2}\right).
\]
Hence Corollary \ref{cor:omega-rk} holds for $\tilde\rho_k$ as well. If instead $\bar\beta$ has a polynomial tail, then $H^{\beta}_{nl}\tilde\rho_k/\tilde\rho_k$ is unbounded, because the weight decays faster than the L\'evy tail; this is precisely the reason for the restriction $r<1+2s_2$.
\end{rmk}

We continue this subsection with the following Kato-type inequality.

\begin{lemm}\label{lemma:kato}
Assume \ref{H0} holds. For every $1<r<1+2s_2$, $k\in(0,1]$, $p\geq2$ and  $x\in C^\infty_c(\R)$, we have
\begin{equation}\label{eq:Kato-type}
\int_\R [H^{\beta}x](\xi){\rm sign}(x(\xi))\abs{x(\xi)}^{p-1}\rho_{r,k}(\xi)\mathd\xi
\geq \frac{1}{p}\int_\R \abs{x(\xi)}^{p}H^{\beta}\rho_{r,k}(\xi)\mathd\xi\geq -\frac{\omega_{r,k}}{p}\norm{x}^p_{L^p_{r,k}}.
\end{equation}
\end{lemm}

\begin{proof}
Fix $x\in C^\infty_c(\R)$. We set $\rho:=\rho_{r,k}$ and we study separately the local and the nonlocal part of $H^{\beta}$. 

For the local part, by the integration-by-parts formula we deduce
\[
-\beta_0\int_\R \Delta x\,{\rm sign}(x)\abs{x}^{p-1}\rho\,\mathd\xi
=\beta_0\int_\R \nabla x\cdot\nabla\!\left({\rm sign}(x)\abs{x}^{p-1}\rho\right)\mathd\xi .
\]
Noting that $\nabla\big({\rm sign}(x)\abs{x}^{p-1}\big)=(p-1)\abs{x}^{p-2}\nabla x$ (valid a.e., after a standard regularization of the sign function) we obtain
\[
-\beta_0\int_\R \Delta x\,{\rm sign}(x)\abs{x}^{p-1}\rho\,\mathd\xi
=\beta_0(p-1)\int_\R \abs{x}^{p-2}\abs{\nabla x}^2\rho\,\mathd\xi
+\frac{\beta_0}{p}\int_\R \nabla(\abs{x}^p)\cdot\nabla\rho\,\mathd\xi .
\]
The first term is nonnegative, and integrating by parts in the second one we get
\[
-\beta_0\int_\R \Delta x\,{\rm sign}(x)\abs{x}^{p-1}\rho\,\mathd\xi
\geq \frac1p\int_\R \abs{x}^{p}\,(-\beta_0\Delta\rho)\,\mathd\xi .
\]
For the nonlocal part, we exploit the following standard inequality
\begin{equation}\label{eq:convexity-pointwise}
{\rm sign}(a)\abs{a}^{p-1}\,(2a-b-c)\geq \frac1p\big(2\abs{a}^{p}-\abs{b}^{p}-\abs{c}^{p}\big),\qquad a,b,c\in\R,\ p\geq 2.
\end{equation}
Multiplying \eqref{eq:Hbeta-nl-only} by ${\rm sign}(x(\xi))\abs{x(\xi)}^{p-1}\rho(\xi)$, integrating, applying \eqref{eq:convexity-pointwise} with $a=x(\xi)$, $b=x(\xi+z)$ and $c=x(\xi- z)$, and using Fubini together with the change of variables $\xi\mapsto\xi\mp z$, we get
\[
\int_\R H^{\beta}_{nl}x\,{\rm sign}(x)\abs{x}^{p-1}\rho\,\mathd\xi
\geq \frac1p\int_\R \abs{x(\xi)}^{p}\,H^{\beta}_{nl}\rho(\xi)\,\mathd\xi .
\]
Adding the two estimates gives \eqref{eq:Kato-type}; the last assertion follows from \eqref{eq:excessive}.
\end{proof}

Finally, exploiting the Lumer--Phillips theorem and Lemma \ref{lemma:kato}, we can prove the following generation result.

\begin{lemm}\label{lemma:semigroupLp}
Assume \ref{H0} holds. For every $1<r<1+2s_2$, $m>0$ and $p\geq 2$, there exists $k_0:=k_0(r,m)$ such that for every $0<k<k_0$ the operator $-\left[H^\beta+m^2\right]:C^\infty_c(\R)\subseteq L^p_{r,k}\rightarrow L^p_{r,k}$ is closable and its closure $-\mathcal{A}_p:{\rm Dom}(\mathcal{A}_p)\subseteq L^p_{r,k}\rightarrow L^p_{r,k}$ is the infinitesimal generator of a strongly continuous semigroup $\{e^{-t\mathcal{A}_p}\}_{t\geq 0}$ and it holds
\begin{align}
&\norm{e^{-t\mathcal A_p}}_{\mathcal L(L^p_{r,k})}\leq e^{\left(\frac{\omega_{r,k}}{p}-m^2\right)t},\qquad t\geq0,\label{eq:Lp-semi}\\[3mm]
&\int_\R [\mathcal A_px](\xi){\rm sign}(x(\xi))\abs{x(\xi)}^{p-1}\rho_{r,k}(\xi)\,\mathd\xi\geq\left(m^2-\frac{\omega_{r,k}}{p}\right)\norm{x}_{L^p_{r,k}}^{p},\qquad x\in {\rm Dom}(\mathcal{A}_p)\label{eq:Ap-lower}
\end{align}
where $\omega_{r,k}$ is given by \eqref{eq:omega-rk-bound} and $m^2-\frac{\omega_{r,k}}{p}>0$. In particular, $C^\infty_c(\R)$ is a core for $-\mathcal{A}_p$, namely it is dense in ${\rm Dom}(\mathcal{A}_p)$ with respect to the graph norm.
\end{lemm}

\begin{proof}
For simplicity we set $\rho:=\rho_{r,k}$, $\norm{\cdot}_p:=\norm{\cdot}_{L^p_{r,k}}$. The idea is to apply the classical Lumer--Phillips \cite[Theorem 4.3]{Paz1983}, since, by standard arguments, the space $C^\infty_c(\R)$ is dense in $L^p_{r,k}$, we need only to prove the dissipativity and the range condition for $-(H^{\beta}+m^2)$.\\

By Corollary \ref{cor:omega-rk} we have $\omega_{r,k}\rightarrow0$ as $k\rightarrow0^+$, so we fix
$k_0=k_0(r,m)\in(0,1]$ such that $\omega_{r,k}<2m^2$ for every $k\in(0,k_0)$; in particular $\omega_{r,k}<pm^2$ for every $p\geq2$. By \eqref{eq:Kato-type} and the H\"older inequality
\begin{equation}\label{eq:resolvent-lower}
\norm{(\lambda+H^{\beta}+m^2)x}_p\geq\left(\lambda+m^2-\frac{\omega_{r,k}}{p}\right)\norm{x}_p ,\qquad x\in C^\infty_c(\R),\ \lambda>0,
\end{equation}
in particular, $-(H^{\beta}+m^2):C^\infty_c(\R)\subseteq L^p_{r,k}\rightarrow L^p_{r,k}$ is dissipative.

 Before proving the range condition, we need two preliminary observations. Fix $\chi\in C^\infty_c(\R)$ with $0\leq\chi\leq1$, $\chi\equiv1$ on
$[-1,1]$ and ${\rm supp}\,\chi\subseteq[-2,2]$, and set $\chi_n:=\chi(\cdot/n)$ for $n\in\N$. It is easy to see
that
\begin{equation}\label{eq:cutoffharmonic}
\lim_{n\rightarrow\infty}\norm{H^{\beta}\chi_n}_{L^\infty(\R)}=0 .
\end{equation}
Second, for every $v\in C^\infty(\R)$ with bounded derivatives of every order and every $x\in C^\infty_c(\R)$ with $x\geq0$ we have the following
Kato inequality
\begin{equation}\label{eq:Kato-pointwise}
\int_\R\abs{v}H^{\beta}x\mathd\xi\leq\int_\R x\,{\rm sign}(v)H^{\beta}v\mathd\xi ,
\end{equation}
where ${\rm sign}(0):=0$. For the nonlocal part this is immediate, since $\abs{v(\xi\pm z)}\geq{\rm sign}(v(\xi))v(\xi\pm z)$ gives
\[
\int_\R\abs{v}H^{\beta}_{nl}x\mathd\xi\leq\int_\R x\,{\rm sign}(v)H^{\beta}_{nl}v\mathd\xi,
\]
after multiplying by $x\geq0$, integrating and using the symmetry of $H^{\beta}_{nl}$; for the local part it is the classical Kato inequality.

Now we can prove the range condition. Fix $\lambda>0$, let $y\in L^{p'}_{r,k}$ with $p'=(p-1)^{-1}p$ be such that $\scal{(\lambda+H^{\beta}+m^2)x}{y\rho}_{L^2}=0$
for every $x\in C^\infty_c(\R)$, and set $u:=\rho y$. By the H\"older inequality and \eqref{eq:rhoL1}, we have
\[
\norm{u}_{L^1(\R)}=\int_\R\abs{y}\rho\mathd\xi\leq\norm{y}_{L^{p'}_{r,k}}\big(c_rk^{-1/2}\big)^{1/p}<+\infty.
\]
By the symmetry of $H^\beta$ in $L^2$, for every $x\in C^\infty_c(\R)$,  we have
\[
0=\scal{(\lambda+H^{\beta}+m^2)x}{\rho y}_{L^2}=\int_\R x(\lambda+H^{\beta}+m^2)u\mathd\xi,
\]
so that 
\[
\int_\R xH^{\beta}u\mathd\xi=-\int_\R x(\lambda+m^2)u\mathd\xi.
\]

Let $\eta_\varepsilon$ be a standard mollifier and $u_\varepsilon:=u*\eta_\varepsilon$. It is easy to see that
for every $x\in C^\infty_c(\R)$,  we have
\[
\int_\R xH^{\beta}u_\varepsilon\mathd\xi=-\int_\R x(\lambda+m^2)u_\varepsilon\mathd\xi,
\]
hence, applying \eqref{eq:Kato-pointwise} with $v=u_\varepsilon$ and $x=\chi_n$, we obtain
\[
(\lambda+m^2)\int_\R\abs{u_\varepsilon}\chi_n\mathd\xi
=-\int_\R \chi_n{\rm sign}(u_\varepsilon)\big(H^{\beta}u_\varepsilon\big)\mathd\xi
\leq-\int_\R\abs{u_\varepsilon}\,H^{\beta}\chi_n\mathd\xi
\leq\norm{u_\varepsilon}_{L^1(\R)}\norm{H^{\beta}\chi_n}_{L^\infty(\R)}.
\]
Letting $n\rightarrow\infty$, the left-hand side converges to
$(\lambda+m^2)\norm{u_\varepsilon}_{L^1(\R)}$ by monotone convergence, while the right-hand side
vanishes by \eqref{eq:cutoffharmonic}. Since $\lambda+m^2>0$, this gives $u_\varepsilon=0$ for every
$\varepsilon>0$, hence $u=0$ and $y=0$. By the Hahn--Banach theorem,
$[\lambda+H^{\beta}+m^2]\left(C^\infty_c(\R)\right)$ is dense in $L^p_{r,k}$ for every $\lambda>0$. By \cite[Theorem 4.5]{Paz1983} the operator $-(H^{\beta}+m^2)+m^2-\frac{\omega_{r,k}}{p}$ is closable and its closure $-\mathcal A_p+m^2-\frac{\omega_{r,k}}{p}$ is dissipative. In addition, the range of  $-\mathcal A_p+m^2-\frac{\omega_{r,k}}{p}$ is closed and so, by the density of $[\lambda+H^{\beta}+m^2]\left(C^\infty_c(\R)\right)$, coincides with $L^p_{r,k}$. Hence we can apply the Lumer--Phillips theorem \cite[Theorem 4.3]{Paz1983} to $-\mathcal A_p+m^2-\frac{\omega_{r,k}}{p}$ and so the proof is concluded (\eqref{eq:Ap-lower} follows from \eqref{eq:Kato-type} and standard density arguments).
\end{proof}

\begin{lemm}\label{lemma:core-adjoint}
Assume \ref{H0} holds, let $m>0$, $1<r<1+2s_2$ and $k\in(0,k_0)$, with $k_0=k_0(r,m)$ given by Lemma \ref{lemma:semigroupLp}, and let $\mathcal A^*$ be the adjoint of $\mathcal A:=\mathcal A_2$ in $L^2_{r,k}$. Then $C^\infty_c(\R)\subseteq{\rm Dom}(\mathcal A^*)$, $\mathcal A^*v=\rho_{r,k}^{-1}(H^\beta+m^2)(\rho_{r,k}v)$ for every $v\in C^\infty_c(\R)$, and $C^\infty_c(\R)$ is a core for $\mathcal A^*$.
\end{lemm}
\begin{proof}
We set $\rho:=\rho_{r,k}$. Let $v\in C^\infty_c(\R)$. By \eqref{eq:testdecay-L2}, $g_v:=\rho^{-1}(H^\beta+m^2)(\rho v)\in L^2_{r,k}$ and, by the symmetry of $H^\beta$ in $L^2$, $\scal{(H^\beta+m^2)u}{v}_{L^2_{r,k}}=\scal{u}{(H^\beta+m^2)(\rho v)}_{L^2}=\scal{u}{g_v}_{L^2_{r,k}}$ for every $u\in C^\infty_c(\R)$. Since $C^\infty_c(\R)$ is a core for $\mathcal A$, we get $v\in{\rm Dom}(\mathcal A^*)$ and $\mathcal A^*v=g_v$.

By \eqref{eq:Lp-semi} with $p=2$, $-\mathcal A^*$ generates the adjoint semigroup (see \cite[Chapter 1, Corollary 10.6]{Paz1983}), which satisfies the same estimate; hence, for $\lambda>0$, $(\lambda+\mathcal A^*)^{-1}$ is bounded from $L^2_{r,k}$ into ${\rm Dom}(\mathcal A^*)$ endowed with the graph norm, and it is enough to prove that $(\lambda+\mathcal A^*)(C^\infty_c(\R))$ is dense in $L^2_{r,k}$ for some $\lambda>0$ with $\lambda+m^2>\omega_{r,k}$. We argue as for the range condition in the proof of Lemma \ref{lemma:semigroupLp}. Let $w\in L^2_{r,k}$ be orthogonal to $(\lambda+\mathcal A^*)(C^\infty_c(\R))$. Since $\rho\,C^\infty_c(\R)=C^\infty_c(\R)$, we get
\[
\int_\R w\,(\lambda+H^{\beta}+m^2)x\,\mathd\xi=0,\qquad x\in C^\infty_c(\R),
\]
namely $w$ satisfies the same equation as $u=\rho y$ in the proof of Lemma \ref{lemma:semigroupLp}; the only difference is that $w\in L^1(\R;\rho\,\mathd\xi)$, by the H\"older inequality and \eqref{eq:rhoL1}, instead of $L^1(\R)$. Therefore we take $x=\rho$ in place of $x=\chi_n$: setting $w_\varepsilon:=w*\eta_\varepsilon$, as there we have $H^\beta w_\varepsilon=-(\lambda+m^2)w_\varepsilon$ and, by \eqref{eq:Kato-pointwise} and \eqref{eq:excessive},
\[
(\lambda+m^2)\int_\R\abs{w_\varepsilon}\rho\,\mathd\xi=-\int_\R\rho\,{\rm sign}(w_\varepsilon)H^{\beta}w_\varepsilon\,\mathd\xi\leq-\int_\R\abs{w_\varepsilon}H^{\beta}\rho\,\mathd\xi\leq\omega_{r,k}\int_\R\abs{w_\varepsilon}\rho\,\mathd\xi,
\]
so that $w_\varepsilon=0$ for every $\varepsilon>0$, and $w=0$. Here \eqref{eq:Kato-pointwise} is applied with $x=\rho\notin C^\infty_c(\R)$ and with $v=w_\varepsilon$, whose derivatives are not bounded: its proof works in the same way, since $w_\varepsilon,w_\varepsilon',w_\varepsilon''\in L^1(\R;\rho\,\mathd\xi)$ and $\abs{\rho'}+\abs{\rho''}+\abs{H^{\beta}_{nl}\rho}\leq C\rho$ by \eqref{eq:weightratio}, so that all the integrals involved are absolutely convergent (for the nonlocal part, after truncating $\bar\beta$ near $0$).
\end{proof}

\subsection{Commutator}\label{section:commutator}

Throughout this subsection we assume \ref{H0}, and we fix $1<r<1+2s_2$ and $k\in(0,1]$. 
We consider the space $\mathcal F:=\{f\in L^2:\mathcal E(f,f)<+\infty\}$, endowed with the norm $\norm{\cdot}_\mathcal{F}=\mathcal E(\cdot,\cdot)^{1/2}$ where $\mathcal{E}$ is the quadratic form associated with $H^{\beta}+m^2$ defined by
\begin{align}\label{eq:dirichletform}
&\mathcal E(f,g):=\mathcal E^{l}(f,g)+\mathcal E^{nl}(f,g)+m^2\scal{f}{g}_{L^2},\\
&\mathcal E^{l}(f,g):=\beta_0\int_\R\nabla f\cdot\nabla g \mathd\xi,\notag\\
&\mathcal E^{nl}(f,g):=\frac{1}{2}\int_{\R^2}\big(f(\xi)-f(\xi')\big)\big(g(\xi)-g(\xi')\big)\bar\beta(\xi-\xi') \mathd(\xi,\xi').\notag
\end{align}
 We underline that, for $\beta_0=0$ by \eqref{eq:comportamento-kernels_1} the norm of $\mathcal{F}$ is equivalent to $H^{s_1}$,  while for $\beta_0>0$ the local part dominates the norm of $\mathcal{F}$ so that is equivalent to $H^{1}$, hence, up to an equivalent renorming, we have
 \begin{equation}\label{eq:energy-is-Hs1}
\mathcal F=\begin{cases}
H^{s_1}(\R),\qquad &\beta_0=0\\
H^1(\R),\qquad &\beta_0>0.
\end{cases}
\end{equation}
  Since $s_1>\tfrac12$, by the fractional Morrey embedding \cite[Theorem 8.2]{DiNezzaPalatucciValdinoci2012} (applied with $n=1$, $p=2$, $s=s_1$),
\begin{equation}\label{eq:energy-embedding}
\mathcal F\subset C^{0,s_1-\frac12}(\R)\subset C^0(\R);
\end{equation}
in particular every $f\in\mathcal F$ admits a continuous representative. The quadratic form $\mathcal E^{l}+\mathcal E^{nl}$ is a regular symmetric Dirichlet form on $L^2$, whose associated Hunt process is a symmetric L\'evy process, hence we denote by $\{Y^\xi_t\}$ the process associated with $-H^{\beta}$ starting from $\xi\in\R$ (see \cite[Chapter 2]{Fukushima2011} and \ref{H0}) defined on some suitable probability space $(\Omega,\mathcal{G},\mathbb{P})$.

Moreover we note that $\rho_{r,k}^{1/2}=\rho_{r/2,k}$, so that all the estimates of Lemma \ref{lemma:weightproperties} apply to $\rho_{r,k}^{1/2}$ as well, with $r$ replaced by $r/2$. Accordingly, for $s\geq 0$ we set
\begin{equation}\label{eq:weightedsobolev}
\norm{f}_{H^{s}_{\rho_{r,k}}(\R)}:=\norm{\rho_{r/2,k}f}_{H^{s}(\R)} .
\end{equation}

We define the commutator
\begin{equation}\label{eq:commutatordef}
\mathcal C_{r,k}f:=H^{\beta}\big(\rho_{r/2,k}f\big)-\rho_{r/2,k}\,H^{\beta}f ,
\end{equation}
so that, for every $f,g\in C^\infty_c(\R)$,
\begin{equation}\label{eq:commutatorsplitting}
\scal{H^{\beta}f}{\rho_{r,k}\,g}_{L^2}=\scal{H^{\beta}\big(\rho_{r/2,k}f\big)}{\rho_{r/2,k}g}_{L^2}-\scal{\mathcal C_{r,k}f}{\rho_{r/2,k}g}_{L^2} .
\end{equation}
When $\beta_0=0$ the commutator has the explicit form
\begin{equation}\label{eq:commutatorexplicit}
\mathcal C_{r,k}f(\xi)=\int_\R f(\eta)\big(\rho_{r/2,k}(\xi)-\rho_{r/2,k}(\eta)\big)\bar\beta(\xi-\eta)\mathd\eta.
\end{equation}

Finally, for $\theta\in(0,1)$ we set
\begin{equation}\label{eq:gagliardo}
\mathcal G_\theta(g):=\frac12\int_\R\int_{\lvert\xi-\eta\rvert\leq1}\frac{\lvert g(\xi)-g(\eta)\rvert^2}{\lvert\xi-\eta\rvert^{1+2\theta}}\mathd\eta\mathd\xi .
\end{equation}

\begin{rmk}\label{rem:different_norms}
Splitting the Gagliardo seminorm of $\dot H^\theta(\R)$ at $\lvert\xi-\eta\rvert=1$ and estimating the outer region by the Minkowski inequality, we obtain
\begin{equation}\label{eq:Gtheta-Htheta}
\mathcal G_\theta(g)+\norm{g}^2_{L^2}\simeq\norm{g}^2_{H^{\theta}(\R)} .
\end{equation}
Recall from \eqref{eq:energy-is-Hs1} that $\scal{H^{\beta}g}{g}_{L^2}+\norm{g}^2_{L^2}\simeq\norm{g}^2_{H^{s_1}(\R)}$, or equivalently that the Fourier symbol of $H^{\beta}$ satisfies $1+\widehat{H^{\beta}}(\xi)\simeq(1+\abs{\xi}^2)^{s_1}$. Hence, if $\theta<s_1$, by interpolation for every $\varepsilon>0$ there exists $C_\varepsilon>0$ such that
\begin{equation}\label{eq:interpolationG}
\mathcal G_\theta(g)\leq\varepsilon\scal{H^{\beta}g}{g}_{L^2}+C_\varepsilon\norm{g}^2_{L^2} .
\end{equation}
\end{rmk}

\begin{prop}\label{prop:levy-commutator}
Assume \ref{H0} holds with $\beta_0=0$. For every $1<r<1+2s_2$, $k\in(0,1]$, $\theta\in(2s_1-1,s_1)$ and $0<\varsigma<\min\{1,2s_2\}$, there exists a constant $C_2=C_2(r,s_1,s_2,\theta,\varsigma,C_\beta,\Theta_\beta)>0$ such that
\begin{equation}\label{eq:commutatorestimate}
\norm{\mathcal C_{r,k}f}^2_{L^2}\leq C_2k^{\varsigma}\left(\norm{\rho_{r/2,k}f}^2_{L^2}+\mathcal G_\theta\big(\rho_{r/2,k}f\big)\right),
\qquad f\in C^\infty_c(\R).
\end{equation}
\end{prop}

\begin{proof}
Set $g:=\rho_{r/2,k}f$ and split $\mathcal C_{r,k}f=\mathcal C^{\mathrm n}+\mathcal C^{\mathrm f}$, where $\mathcal C^{\mathrm n}$ and $\mathcal C^{\mathrm f}$ denote the contributions to \eqref{eq:commutatorexplicit} of the regions $\{\lvert\xi-\eta\rvert\leq1\}$ and $\{\lvert\xi-\eta\rvert\geq1\}$; only the former requires the principal value. In what follows we repeatedly use Lemma \ref{lemma:weightproperties} with $r/2$ in place of $r$.

We focus first on the region $\{\lvert\xi-\eta\rvert\leq1\}$. Symmetrizing in $z\mapsto-z$ and adding and subtracting $f(\xi)$ we obtain the absolutely convergent expression
\[
\mathcal C^{\mathrm n}(\xi)=\frac12\int_{\lvert z\rvert\leq1}\left[f(\xi)\,\Lambda(\xi,z)+\sum_{\pm}\big(f(\xi\pm z)-f(\xi)\big)\big(\rho_{r/2,k}(\xi)-\rho_{r/2,k}(\xi\pm z)\big)\right]\bar\beta(z)\mathd z ,
\]
where $\Lambda(\xi,z):=2\rho_{r/2,k}(\xi)-\rho_{r/2,k}(\xi+z)-\rho_{r/2,k}(\xi-z)$. Since $\lvert z\rvert\leq1\leq k^{-1/2}$, by \eqref{eq:rhoderivatives}, \eqref{eq:rhoderivatives2} and \eqref{eq:rhocomparable} we have
\[
\lvert\Lambda(\xi,z)\rvert\lesssim_r k\,\rho_{r/2,k}(\xi)\lvert z\rvert^2,
\qquad
\big\lvert\rho_{r/2,k}(\xi)-\rho_{r/2,k}(\xi\pm z)\big\rvert\lesssim_r\sqrt k\,\rho_{r/2,k}(\xi)\lvert z\rvert ,
\]
while, writing $f=g/\rho_{r/2,k}$,
\[
\lvert f(\xi\pm z)-f(\xi)\rvert\lesssim_r\frac{\lvert g(\xi\pm z)-g(\xi)\rvert}{\rho_{r/2,k}(\xi)}+\frac{\lvert g(\xi)\rvert}{\rho_{r/2,k}(\xi)}\sqrt k\lvert z\rvert .
\]
Since $\bar\beta(z)\leq C_\beta\lvert z\rvert^{-1-2s_1}$ for $\lvert z\rvert\leq1$ by \eqref{eq:levybilateral}, and $\int_{\lvert z\rvert\leq1}\lvert z\rvert^{1-2s_1}\mathd z<+\infty$ because $s_1<1$, we get
\[
\lvert\mathcal C^{\mathrm n}(\xi)\rvert
\lesssim
k\lvert g(\xi)\rvert+\sqrt k\sum_{\pm}\int_{\lvert z\rvert\leq1}\frac{\lvert g(\xi\pm z)-g(\xi)\rvert}{\lvert z\rvert^{2s_1}}\mathd z .
\]
By the Cauchy--Schwarz inequality,
\[
\int_{\lvert z\rvert\leq1}\frac{\lvert g(\xi\pm z)-g(\xi)\rvert}{\lvert z\rvert^{2s_1}}\mathd z
\leq
\left(\int_{\lvert z\rvert\leq1}\frac{\lvert g(\xi\pm z)-g(\xi)\rvert^2}{\lvert z\rvert^{1+2\theta}}\mathd z\right)^{\!\frac12}
\left(\int_{\lvert z\rvert\leq1}\lvert z\rvert^{1+2\theta-4s_1}\mathd z\right)^{\!\frac12},
\]
and the second integral is finite precisely because $\theta>2s_1-1$. Squaring and integrating in $\xi$ we conclude that
\begin{equation}\label{eq:nearpart}
\norm{\mathcal C^{\mathrm n}}^2_{L^2}\lesssim k^2\norm{g}^2_{L^2}+k\,\mathcal G_\theta(g).
\end{equation}

Now we focus on the region $\{\lvert\xi-\eta\rvert\geq1\}$. Here the integral converges absolutely and
\[
\lvert\mathcal C^{\mathrm f}(\xi)\rvert\leq\int_{\lvert\xi-\eta\rvert\geq1}\frac{\big\lvert\rho_{r/2,k}(\xi)-\rho_{r/2,k}(\eta)\big\rvert}{\rho_{r/2,k}(\eta)}\lvert g(\eta)\rvert\bar\beta(\xi-\eta)\mathd\eta ,
\]
so we estimate the $L^2$ norm of the operator with kernel $\big\lvert\rho_{r/2,k}(\xi)-\rho_{r/2,k}(\eta)\big\rvert\rho_{r/2,k}(\eta)^{-1}\bar\beta(\xi-\eta)$ on the three regions
\[
F_1:=\left\{1\leq\lvert\xi-\eta\rvert\leq k^{-1/2}\right\},
\quad
F_2:=\left\{\lvert\xi-\eta\rvert\geq k^{-1/2},\ \lvert\eta\rvert\leq2\lvert\xi\rvert\right\},
\quad
F_3:=\left\{\lvert\xi-\eta\rvert\geq k^{-1/2},\ \lvert\eta\rvert\geq2\lvert\xi\rvert\right\}.
\]
Throughout we use that, by \eqref{eq:levygrowth} and $k^{-1/2}\geq1$,
\begin{equation}\label{eq:tailmass}
\int_{\lvert z\rvert\geq k^{-1/2}}\bar\beta(z)\mathd z\leq\frac{C_\beta}{s_2}k^{s_2}.
\end{equation}

On $F_1$ we have $\lvert\xi-\eta\rvert\leq k^{-1/2}$, so that \eqref{eq:rhoderivatives} and \eqref{eq:rhocomparable} give $\rho_{r/2,k}(\xi)\simeq_r\rho_{r/2,k}(\eta)$ and $\big\lvert\rho_{r/2,k}(\xi)-\rho_{r/2,k}(\eta)\big\rvert\lesssim_r\rho_{r/2,k}(\xi)\min\{1,\sqrt k\lvert\xi-\eta\rvert\}$; hence the kernel is bounded by the symmetric kernel $\mathbb I_{F_1}(\xi,\eta)\min\{1,\sqrt k\lvert\xi-\eta\rvert\}\lvert\xi-\eta\rvert^{-1-2s_2}$, whose two Schur constants coincide and satisfy
\[
\int_{\lvert z\rvert\geq1}\min\left\{1,\sqrt k\lvert z\rvert\right\}\lvert z\rvert^{-1-2s_2}\mathd z
\leq
\sqrt k\int_{1\leq\lvert z\rvert\leq k^{-1/2}}\lvert z\rvert^{-2s_2}\mathd z+\int_{\lvert z\rvert\geq k^{-1/2}}\lvert z\rvert^{-1-2s_2}\mathd z
\lesssim k^{\varsigma/2},
\]
so that Schur's test gives a contribution bounded by $k^{\varsigma}\norm{g}^2_{L^2}$.

On $F_2$ estimate \eqref{eq:rhodoubling} gives $\rho_{r/2,k}(\xi)\leq2^{r/2}\rho_{r/2,k}(\eta)$, so that the ratio in the kernel is bounded by a constant depending only on $r$ and the kernel is bounded by $\mathbb I_{\{\lvert\xi-\eta\rvert\geq k^{-1/2}\}}\bar\beta(\xi-\eta)$, whose Schur constants are both bounded by \eqref{eq:tailmass}; the corresponding contribution is bounded by $k^{2s_2}\norm{g}^2_{L^2}$.

On $F_3$ we have $\lvert\xi\rvert\leq\frac12\lvert\eta\rvert$, whence $\frac12\lvert\eta\rvert\leq\lvert\xi-\eta\rvert\leq\frac32\lvert\eta\rvert$ and, since $\lvert\xi-\eta\rvert\geq k^{-1/2}$, also $\lvert\eta\rvert\geq\frac23k^{-1/2}$. Using $\big\lvert\rho_{r/2,k}(\xi)-\rho_{r/2,k}(\eta)\big\rvert\leq\rho_{r/2,k}(\xi)+\rho_{r/2,k}(\eta)$, the contribution of $\rho_{r/2,k}(\eta)$ is estimated exactly as on $F_2$, while for the contribution of $\rho_{r/2,k}(\xi)$ the Cauchy--Schwarz inequality gives
\[
\left\lvert\mathcal C^{F_3}(\xi)\right\rvert^2
\lesssim
\rho_{r,k}(\xi)\left(\int_{\lvert\eta\rvert\geq2\lvert\xi\rvert}\lvert g\rvert^2\mathd\eta\right)J(\xi),
\qquad
J(\xi):=\int_{\lvert\eta\rvert\geq M(\xi)}\frac{\mathd\eta}{\rho_{r,k}(\eta)\lvert\eta\rvert^{2+4s_2}},
\]
where $M(\xi):=\max\left\{2\lvert\xi\rvert,\frac23k^{-1/2}\right\}$. Since $\lvert\eta\rvert\geq\frac23k^{-1/2}$ implies $1+k\lvert\eta\rvert^2\leq\frac{13}{4}k\lvert\eta\rvert^2$, we have $\rho_{r,k}(\eta)^{-1}\lesssim_r(k\lvert\eta\rvert^2)^{r/2}$ on the domain of integration, so that, using $r<1+2s_2<1+4s_2$,
\[
J(\xi)\lesssim_r k^{r/2}\int_{\lvert\eta\rvert\geq M(\xi)}\lvert\eta\rvert^{\,r-2-4s_2}\mathd\eta\lesssim k^{r/2}M(\xi)^{\,r-1-4s_2} .
\]
If $\lvert\xi\rvert\leq k^{-1/2}$, then $M(\xi)\simeq k^{-1/2}$ and $\rho_{r,k}(\xi)\leq1$, so that $\rho_{r,k}(\xi)J(\xi)\lesssim k^{\frac12+2s_2}$ and, integrating over $\{\lvert\xi\rvert\leq k^{-1/2}\}$, whose measure is $2k^{-1/2}$, we obtain a contribution bounded by $k^{2s_2}\norm{g}^2_{L^2}$. If $\lvert\xi\rvert\geq k^{-1/2}$, then $M(\xi)=2\lvert\xi\rvert$ and $\rho_{r,k}(\xi)\leq(k\lvert\xi\rvert^2)^{-r/2}$, whence $\rho_{r,k}(\xi)J(\xi)\lesssim\lvert\xi\rvert^{-1-4s_2}$ and, by the Fubini theorem,
\[
\int_{\lvert\xi\rvert\geq k^{-1/2}}\frac{1}{\lvert\xi\rvert^{1+4s_2}}\left(\int_{\lvert\eta\rvert\geq2\lvert\xi\rvert}\lvert g\rvert^2\mathd\eta\right)\mathd\xi
\leq
\int_\R\lvert g(\eta)\rvert^2\left(\int_{k^{-1/2}\leq\lvert\xi\rvert\leq\frac12\lvert\eta\rvert}\frac{\mathd\xi}{\lvert\xi\rvert^{1+4s_2}}\right)\mathd\eta
\lesssim k^{2s_2}\norm{g}^2_{L^2} .
\]
Collecting \eqref{eq:nearpart} and the estimates on $F_1$, $F_2$, $F_3$, and recalling that $\varsigma<\min\{1,2s_2\}$, we obtain \eqref{eq:commutatorestimate}.
\end{proof}

We now turn to the case $\beta_0>0$, in which \ref{H0} gives no information on $\bar\beta$ near the origin and $s_1=1$. The estimate below uses only \eqref{eq:levynear} and \eqref{eq:levygrowth}, and the seminorm $\mathcal G_\theta$ is replaced by the full gradient.

\begin{prop}\label{prop:commutator-full-operator}
Assume \ref{H0} holds with $\beta_0>0$. For every $1<r<1+2s_2$, $k\in(0,1]$ and $0<\varsigma<\min\{1,2s_2\}$, there exists a constant $C_3=C_3(r,s_2,\varsigma,\beta_0,C_\beta,\Theta_\beta)>0$ such that
\begin{equation}\label{eq:commutatorfull}
\norm{\mathcal C_{r,k}f}^2_{L^2}
\leq
C_3k^{\varsigma}\left(\norm{\rho_{r/2,k}f}^2_{L^2}+\norm{\nabla\big(\rho_{r/2,k}f\big)}^2_{L^2}\right),
\qquad f\in C^\infty_c(\R).
\end{equation}
Moreover, since $\beta_0\norm{\nabla g}^2_{L^2}\leq\scal{H^{\beta}g}{g}_{L^2}$, the right-hand side of \eqref{eq:commutatorfull} is bounded by $\beta_0^{-1}C_3k^{\varsigma}\big(\scal{H^{\beta}g}{g}_{L^2}+\beta_0\norm{g}^2_{L^2}\big)$ with $g=\rho_{r/2,k}f$, which is the analogue of the bound obtained from \eqref{eq:commutatorestimate} through \eqref{eq:interpolationG}.
\end{prop}

\begin{proof}
Set $g:=\rho_{r/2,k}f$ and split $\mathcal C_{r,k}f=-\beta_0\big(\Delta(\rho_{r/2,k}f)-\rho_{r/2,k}\Delta f\big)+\mathcal C^{nl}_{r,k}f$. The local contribution is
\[
-\beta_0\big(\Delta(\rho_{r/2,k}f)-\rho_{r/2,k}\Delta f\big)
=-\beta_0\left(2\frac{\rho_{r/2,k}'}{\rho_{r/2,k}}g'-2\left(\frac{\rho_{r/2,k}'}{\rho_{r/2,k}}\right)^{\!2}g+\frac{\rho_{r/2,k}''}{\rho_{r/2,k}}g\right),
\]
so that, by \eqref{eq:rhoderivatives} and \eqref{eq:rhoderivatives2} applied with $r/2$ in place of $r$,
\[
\norm{\beta_0\big(\Delta(\rho_{r/2,k}f)-\rho_{r/2,k}\Delta f\big)}^2_{L^2}\lesssim\beta_0^2k\left(\norm{g'}^2_{L^2}+k\norm{g}^2_{L^2}\right).
\]
For the nonlocal contribution we split as in the proof of Proposition \ref{prop:levy-commutator}. The estimates on $F_1$, $F_2$ and $F_3$ use only \eqref{eq:levygrowth} and Lemma \ref{lemma:weightproperties}, hence they apply verbatim and give a contribution bounded by $Ck^{\varsigma}\norm{g}^2_{L^2}$. For the near part, the same symmetrization as above gives
\[
\lvert\mathcal C^{\mathrm n}(\xi)\rvert
\lesssim
k\lvert g(\xi)\rvert\int_{\lvert z\rvert\leq1}\lvert z\rvert^2\bar\beta(z)\mathd z
+\sqrt k\sum_{\pm}\int_{\lvert z\rvert\leq1}\lvert z\rvert\lvert g(\xi\pm z)-g(\xi)\rvert\bar\beta(z)\mathd z ,
\]
and, writing $\lvert g(\xi\pm z)-g(\xi)\rvert\leq\lvert z\rvert\int_0^1\lvert g'(\xi\pm tz)\rvert\mathd t$, the Minkowski integral inequality yields
\[
\norm{\mathcal C^{\mathrm n}}_{L^2}
\lesssim
k\,\Theta_\beta\norm{g}_{L^2}+\sqrt k\,\Theta_\beta\norm{g'}_{L^2} ,
\]
with $\Theta_\beta$ as in \eqref{eq:levynear}. Collecting the three contributions and using $\varsigma\leq1$ we obtain \eqref{eq:commutatorfull}.
\end{proof}

We conclude the subsection with the coercivity estimate for the weighted energy, which is the form in which the commutator estimates are used in Section \ref{section:Fernique} and in Section \ref{section:existence}.

\begin{lemm}\label{lemma:coercivity-H3}
Assume \ref{H0} and let $m>0$, $1<r<1+2s_2$ and $k\in(0,k_0(r,m))$, with $k_0(r,m)$ given by Lemma \ref{lemma:semigroupLp}. Let $\mathcal E$ be the quadratic form associated with $H^\beta+m^2$ defined in \eqref{eq:dirichletform}.
Then there exist $C_4:=C_4(r,k,H^\beta,m),C_5:=C_5(r,k,H^\beta,m)>0$, such that for every $f,g \in H^{s_1}_{\rho_{r,k}}$ we have
\begin{equation}\label{eq:coercivity-H3}
\mathcal E(f,\rho_{r,k}f)\geq C_4\norm{f}^2_{H^{s_1}_{\rho_{r,k}}(\R)},\qquad
\abs{\mathcal E(f,\rho_{r,k}g)}\leq C_5\norm{f}_{H^{s_1}_{\rho_{r,k}}(\R)}\norm{g}_{H^{s_1}_{\rho_{r,k}}(\R)} .
\end{equation}
\end{lemm}

\begin{proof}
Let $f,g\in C^\infty_c(\R)$ and set $v:=\rho_{r/2,k}f$, so that $\norm{v}_{H^{s_1}(\R)}=\norm{f}_{H^{s_1}_{\rho_{r,k}}(\R)}$ by \eqref{eq:weightedsobolev}. Since $\rho_{r,k}g\in C^\infty(\R)$ is bounded with bounded derivatives, by the integration by parts formula and the symmetry of $\bar\beta$ we have
\begin{equation}\label{eq:E-vs-A}
\mathcal E(f,\rho_{r,k}g)=\scal{(H^\beta+m^2)f}{\rho_{r,k}g}_{L^2(\R)},\qquad\text{in particular}\qquad \mathcal E(f,\rho_{r,k}f)=\scal{\mathcal A_2f}{f}_{L^2_{r,k}} .
\end{equation}
By \eqref{eq:energy-is-Hs1} there exists $c_1>0$ such that $\mathcal E(v,v)\geq c_1\norm{v}^2_{H^{s_1}(\R)}$. By \eqref{eq:commutatorestimate} together with \eqref{eq:Gtheta-Htheta} and $\theta<s_1$ if $\beta_0=0$, by \eqref{eq:commutatorfull} if $\beta_0>0$, and since $k\leq1$, there exists $S>0$ such that $\norm{\mathcal C_{r,k}f}_{L^2(\R)}\leq S\norm{v}_{H^{s_1}(\R)}$. Hence, by \eqref{eq:commutatorsplitting}, \eqref{eq:E-vs-A} and the Young inequality,
\[
\mathcal E(f,\rho_{r,k}f)=\mathcal E(v,v)-\scal{\mathcal C_{r,k}f}{v}_{L^2(\R)}\geq\frac{c_1}{2}\norm{v}^2_{H^{s_1}(\R)}-K\norm{v}^2_{L^2(\R)},\qquad K:=\frac{S^2}{2c_1}.
\]
On the other hand, by \eqref{eq:E-vs-A} and \eqref{eq:Ap-lower} with $p=2$, we have $\mathcal E(f,\rho_{r,k}f)\geq c_2\norm{v}^2_{L^2(\R)}$, where $c_2:=m^2-\frac{\omega_{r,k}}{2}>0$ since $k<k_0(r,m)$. The convex combination of the two inequalities with weights $\frac{c_2}{c_2+K}$ and $\frac{K}{c_2+K}$ gives the first estimate in \eqref{eq:coercivity-H3} with $C_4:=\frac{c_1c_2}{2(c_2+K)}$. For the second one, by \eqref{eq:commutatorsplitting} and \eqref{eq:E-vs-A},
\[
\mathcal E(f,\rho_{r,k}g)=\mathcal E\big(\rho_{r/2,k}f,\rho_{r/2,k}g\big)-\scal{\mathcal C_{r,k}f}{\rho_{r/2,k}g}_{L^2(\R)},
\]
where the first term is bounded by $\norm{\rho_{r/2,k}f}_{\mathcal F}\norm{\rho_{r/2,k}g}_{\mathcal F}$ by the Cauchy--Schwarz inequality, and the second one by $S\norm{\rho_{r/2,k}f}_{H^{s_1}(\R)}\norm{\rho_{r/2,k}g}_{L^2(\R)}$. By \eqref{eq:energy-is-Hs1} and \eqref{eq:weightedsobolev} we get the second estimate in \eqref{eq:coercivity-H3}. By density, \eqref{eq:coercivity-H3} holds for every $f,g\in H^{s_1}_{\rho_{r,k}}(\R)$; in particular $C_4$ and $C_5$ do not depend on the supports of $f$ and $g$.
\end{proof}

\begin{rmk}\label{rmk:domainA2}
Let $1<r<1+2s_2$ and $k\in(0,k_0)$, with $k_0$ given by Lemma \ref{lemma:semigroupLp}. Then $H^{2s_1}_{\rho_{r,k}}\subseteq{\rm Dom}(\mathcal{A}_2)$. Indeed, for $\phi\in C^\infty_c(\R)$, by \eqref{eq:commutatordef} we have $\rho_{r/2,k}(H^{\beta}+m^2)\phi=(H^{\beta}+m^2)(\rho_{r/2,k}\phi)-\mathcal{C}_{r,k}\phi$, so that Remark \ref{rem:different_norms}, together with \eqref{eq:commutatorestimate} if $\beta_0=0$, \eqref{eq:commutatorfull} if $\beta_0>0$, and \eqref{eq:Gtheta-Htheta}, gives a constant $C>0$ such that $\norm{\mathcal{A}_2\phi}_{L^2_{r,k}}\leq C\norm{\phi}_{H^{2s_1}_{\rho_{r,k}}}$. Since $C^\infty_c(\R)$ is dense in $H^{2s_1}_{\rho_{r,k}}$ and $\mathcal{A}_2$ is closed, the claim follows.
\end{rmk}

\section{Harmonic extensions}\label{section:harmonicextension}

This section is devoted to the study of harmonic extensions of $H^\beta+m^2$. \\

\subsection{On the harmonic extensions}

We recall that $\mathcal{E}$ is the quadratic form associated to $H^\beta+m^2$ given by \eqref{eq:dirichletform} and $\mathcal{F}$ is the domain of $\mathcal{E}$ given by \eqref{eq:energy-is-Hs1}.

Let $\mathcal I\subset\R$ be a bounded open interval and let $\tau^{\mathcal I,\xi}=\inf\{t>0:Y_t^\xi\in\mathcal I^c\}$ be the first exit time from $\mathcal I$.

\begin{defi}\label{def:harmonicextension}
Assume condition \ref{H0} holds. Let $y:\mathcal I^c\to\R$ be measurable. Consider a function $\varphi:\R\rightarrow\R$ such that $\varphi|_{\mathcal I^c}=y$. We say that
\begin{enumerate}
\item $\varphi$ is an energy harmonic extension of $y$ if $\varphi\in\mathcal F$ and $\mathcal E(\varphi,h)=0$ for every $h\in\mathcal F$ with $\supp(h)\subset\mathcal I$;
\item $\varphi$ is the probabilistic harmonic extension of $y$ if
\begin{equation}\label{eq:harmonicextensionprob}
\varphi(\xi)=\mathbb E\left[e^{-m^2\tau^{\mathcal I,\xi}}\,y\big(Y_{\tau^{\mathcal I,\xi}}^\xi\big)\right],\qquad\xi\in\R,
\end{equation}
whenever the right-hand side is well defined (e.g. for $y$ bounded);
\item $\varphi$ is a weak harmonic extension of $y$ if $\varphi\in L^1_{\loc}(\R)$ is such that $\varphi\,(H^{\beta}+m^2)h\in L^1(\R)$ and $\int_\R\varphi\,(H^{\beta}+m^2)h\,\mathd\xi=0$ for every $h\in C^\infty_c(\mathcal I)$  (e.g. $\varphi\in L^1_{r,k}$ with $1<r<1+2s_2$ and $k>0$, see Lemma \ref{lemma:testdecay}).
\end{enumerate}
In the case where a function $\varphi$ satisfies all three notions, we call it a harmonic extension without further qualification.
\end{defi}

\begin{lemm}\label{lemma:harmonicextensionenergy}
Assume condition \ref{H0} holds. Let $y:\mathcal I^c\to\R$ be measurable such that there exists $\tilde y\in\mathcal F$ with $\tilde y_{|\mathcal{I}^c}=y$. Then, for $\varphi\in\mathcal F$ with $\varphi_{|\mathcal{I}^c}=y$, the three conditions of Definition \ref{def:harmonicextension} are equivalent, and there exists a unique $\varphi^y\in\mathcal F$ satisfying them. We call $\varphi^y$ harmonic extension of $y$.
\end{lemm}

\begin{proof}
The set $\mathcal F_{\mathcal I}:=\overline{\{h\in\mathcal F:\supp(h)\subset\mathcal I\}}^{\|\cdot\|_\mathcal{F}}$ is a closed subspace of the Hilbert space $(\mathcal F,\norm{\cdot}_{\mathcal{F}})$. Fix $\varphi\in\mathcal{F}$ with $\varphi_{|\mathcal{I}^c}=y$. By \eqref{eq:energy-embedding} every function in $\mathcal F$ is continuous and, since $s_1>\frac12$, a function of $\mathcal F$ vanishing on $\mathcal I^c$ belongs to $\mathcal F_{\mathcal I}$ (see \cite[Theorem 4.4.3]{Fukushima2011}). Therefore $w:=\tilde y-\varphi\in\mathcal F_{\mathcal I}$ and $\varphi$ is an energy harmonic extension if and only if 
\begin{equation}\label{eq:lax}
\mathcal E(w,h)=\mathcal E(\tilde y,h)\qquad \forall\, h\in\mathcal F_{\mathcal I}.
\end{equation}
Since $m>0$, $\mathcal E$ is coercive on $\mathcal F$, and by the Lax--Milgram theorem there exists a unique solution $w$ to \eqref{eq:lax} and so a unique  energy harmonic extension $\varphi^y$.\\

\noindent\textbf{1.$\Longleftrightarrow$3.}
Since $C^\infty_c(\mathcal I)$ is dense in $\{h\in\mathcal F:\supp(h)\subset\mathcal I\}$ with respect to $\norm{\cdot}_{\mathcal{F}}$ (see \cite[Chapter 2]{Fukushima2011}), condition in 1. of Definition \ref{def:harmonicextension} is equivalent to $\mathcal E(\varphi,h)=0$ for all $h\in C^\infty_c(\mathcal I)$. For $\varphi\in\mathcal F$ and $h\in C^\infty_c(\mathcal I)$ one has $\mathcal E(\varphi,h)=\int_\R\varphi(\xi) [(H^{\beta}+m^2)h](\xi) \mathd\xi$, so 1. and 3. of Definition \ref{def:harmonicextension} are equivalent.\\

\noindent\textbf{1.$\Longleftrightarrow$2.} By \eqref{eq:energy-embedding}, $\varphi$ is continuous and an energy harmonic extension in $\mathcal I$; by the probabilistic characterization of harmonic functions for regular Dirichlet forms (\cite[Theorem 4.3.2 and §4.4]{Fukushima2011}), $\varphi$ coincides with the function $\xi\mapsto\mathbb E[e^{-m^2\tau^{\mathcal I,\xi}}y(Y_{\tau^{\mathcal I,\xi}}^\xi)]$. Conversely, if $\varphi\in\mathcal F$ satisfies 2., then $\varphi$ coincides with this function, and so with $\varphi^y$, which satisfies 1. 
\end{proof}

Denote by $G^{\mathcal I}_{\free}$ the Green function of $H^{\beta}+m^2$ on $\mathcal I$ with zero exterior condition. Define the Poisson kernel
\begin{equation}\label{eq:poissonkernel}
P^{\mathcal I}_{\free}(\xi,B):=\mathbb E\left[e^{-m^2\tau^{\mathcal I,\xi}}\,\mathbb I_B\big(Y_{\tau^{\mathcal I,\xi}}^\xi\big)\right],\qquad B\in\mathcal B(\mathcal I^c).
\end{equation}

\begin{prop}[Ikeda--Watanabe formula]\label{prop:ikedawatanabe}
Assume condition \ref{H0} holds. Let $\mathcal I=(a,b)$. For every $\xi\in\mathcal I$ the Poisson kernel decomposes as
\begin{equation}\label{eq:IkedaWatanabeformula}
P^{\mathcal I}_{\free}(\xi,\mathd\xi')=f_a(\xi)\delta_a(\mathd\xi')+f_b(\xi)\delta_b(\mathd\xi')+\left(\int_{\mathcal I}G^{\mathcal I}_{\free}(\xi,\xi'')\bar\beta(\xi''-\xi')\mathd\xi''\right)\mathd\xi',
\end{equation}
where $f_a(\xi)=\mathbb E[e^{-m^2\tau^{\mathcal I,\xi}}\mathbb I_{\{a\}}(Y^\xi_{\tau^{\mathcal I,\xi}})]$ and $f_b(\xi)=\mathbb E[e^{-m^2\tau^{\mathcal I,\xi}}\mathbb I_{\{b\}}(Y^\xi_{\tau^{\mathcal I,\xi}})]$ are nonnegative and bounded. If $\beta_0=0$, then $f_a\equiv f_b\equiv 0$, and $P^{\mathcal I}_{\free}(\xi,\cdot)$ is absolutely continuous on $\mathcal I^c$ with density
\begin{equation}\label{eq:IkedaWatanabe-ac}
p^{\mathcal I}_{\free}(\xi,\xi')=\int_{\mathcal I}G^{\mathcal I}_{\free}(\xi,\xi'')\,\bar\beta(\xi''-\xi')\,\mathd\xi'',\qquad \xi\in\mathcal I,\ \xi'\in\mathcal I^c.
\end{equation}
Consequently, the harmonic extension $\varphi^y$ of $y\in C(\mathcal{I}^c;\R)$ admits the representation
\begin{equation}\label{eq:harmonic-representation}
\varphi^y(\xi)=y(\xi)\,\mathbb I_{\mathcal I^c}(\xi)+\left(f_a(\xi)y(a)+f_b(\xi)y(b)+\int_{\mathcal{I}^c}y(\xi')p^{\mathcal I}_{\free}(\xi,\xi')\mathd\xi'\right)\mathbb I_{\mathcal I}(\xi),\quad\xi\in\R,
\end{equation}
which for $\beta_0=0$ reduces to $\varphi^y(\xi)=y(\xi)\mathbb I_{\mathcal I^c}(\xi)+\big(\int_{\mathcal I^c}y(\xi')p^{\mathcal I}_{\free}(\xi,\xi')\,\mathd\xi'\big)\mathbb I_{\mathcal I}(\xi)$.
\end{prop}

\begin{proof}
The absolutely continuous part of \eqref{eq:IkedaWatanabeformula} comes from the standard Ikeda--Watanabe formula. Indeed by \cite[Lemma 4.5.5 and formula (4.5.16)]{Fukushima2011} for every $h\in C_c(\mathcal I^c)$ with $h(a)=h(b)=0$ we have
\begin{equation}\label{eq:IW-jump}
\mathbb E\left[e^{-m^2\tau^{\mathcal I,\xi}}h\big(Y_{\tau^{\mathcal I,\xi}}^\xi\big)\right]
=\int_{\mathcal I^c}\int_{\mathcal I} h(\xi')G^{\mathcal I}_{\free}(\xi,\xi'')\bar\beta(\xi''-\xi')\mathd\xi''  \mathd\xi'.
\end{equation}
It remains to account for the exit occurring \emph{without} a jump, i.e. on the event $\{Y^\xi_{\tau^{\mathcal I,\xi}-}=Y^\xi_{\tau^{\mathcal I,\xi}}\}$. By definition $P^{\mathcal I}_{\free}$ is a positive measure which integrated with respect to functions $h$ zero in $\partial\mathcal I=\{a,b\}$ is absolutely continuous and it has the expression \eqref{eq:IW-jump}. Thus the difference between the expression \eqref{eq:IW-jump} and  $P^{\mathcal I}_{\free}$ is a positive measure concentrated on $\{a,b\}$, which being a finite set is a combination of Dirac deltas.\\
If $\beta_0=0$, the process has no Gaussian part and it is symmetric. Hence, by \cite[Corollary 3.5 page 477]{Millar1973}, $Y^\xi$ can leave $\mathcal I$ through $\partial\mathcal I=\{a,b\}$ only by a jump landing exactly on $a$ or $b$, which is impossible because $\bar\beta\mathd z$ is absolutely continuous. Therefore $\mathbb P(Y^\xi_{\tau^{\mathcal I,\xi}}\in\{a,b\})=0$ and $f_a\equiv f_b\equiv0$.
\end{proof}

\begin{theorem}\label{thm:harmonic-continuity}
Assume condition \ref{H0} holds. Let $1<r<1+2s_2$ and $k>0$. For every $y\in C(\mathcal I^c,\R)\cap L^1_{r,k}(\mathcal I^c)$ the function $\varphi^y$ given by \eqref{eq:harmonicextensionprob} is continuous on $\R$ and is a weak harmonic extension of $y$ in the sense of Definition \ref{def:harmonicextension}.
Further for every $u\in\mathcal F\cap{\rm Dom}(H^{\beta})$ such that $\supp(u)\subseteq\overline{\mathcal I}$ and
$\sup_{\xi\in\mathcal I^c}\rho_{r,k}(\xi)^{-1}\lvert(H^{\beta}+m^2)u(\xi)\rvert<+\infty$, we have
\begin{equation}\label{eq:weakharmonicF}
\int_\R\varphi^y\,(H^{\beta}+m^2)u\,\mathd\xi=0.
\end{equation}
Moreover, the map $y\longmapsto\varphi^y$ is continuous from the metric space $C(\mathcal I^c,\R)\cap L^1_{r,k}(\mathcal I^c)$, endowed with the natural metric of $L^1_{r,k}$ together with the metric associated to uniform convergence on compact subsets of $\mathcal I^c$, into $L^\infty_{\loc}(\R)$.
\end{theorem}

\begin{proof}
Linearity is clear from \eqref{eq:harmonic-representation}; the estimates below also show that all the integrals appearing there are finite, so that $\varphi^y$ is well defined. We set $\mathcal N:= [a-1,a]\cup[b,b+1]$. Fix $K$ compact subset of $\R$.
On $K\cap\mathcal I^c$ we have $\varphi^y=y$, so it remains to bound $\varphi^y$ on $\mathcal I$. By \eqref{eq:harmonic-representation}, for $\xi\in\mathcal I$,
\begin{equation}\label{eq:cont-split}
\abs{\varphi^y(\xi)}\leq f_a(\xi)\abs{y(a)}+f_b(\xi)\abs{y(b)}+\int_{\mathcal N}\abs{y(\xi')}\,p^{\mathcal I}_{\free}(\xi,\xi')\,\mathd\xi'+\int_{\mathcal I^c\setminus\mathcal N}\abs{y(\xi')}\,p^{\mathcal I}_{\free}(\xi,\xi')\,\mathd\xi'.
\end{equation}
By \eqref{eq:poissonkernel} we have
\begin{equation}\label{eq:total-mass}
f_a(\xi)+f_b(\xi)+\int_{\mathcal I^c}p^{\mathcal I}_{\free}(\xi,\xi')\,\mathd\xi'=\mathbb E\left[e^{-m^2\tau^{\mathcal I,\xi}}\right]\leq1,\qquad\xi\in\mathcal I.
\end{equation}
\noindent Since $\mathcal N$ is compact and $y$ is continuous, by \eqref{eq:total-mass} we deduce
\[
\int_{\mathcal N}\abs{y(\xi')}\,p^{\mathcal I}_{\free}(\xi,\xi')\,\mathd\xi'\leq\norm{y}_{L^\infty(\mathcal N)}\int_{\mathcal I^c}p^{\mathcal I}_{\free}(\xi,\xi')\,\mathd\xi'\leq\norm{y}_{L^\infty(\mathcal N)},
\]
uniformly in $\xi\in\mathcal I$.\\
For the region $\mathcal I^c\setminus\mathcal N$, since $m^2>0$, we get
\begin{equation}\label{eq:green-mass}
\int_{\mathcal I}G^{\mathcal I}_{\free}(\xi,\xi'')\mathd\xi''=\mathbb E\left[\int_0^{\tau^{\mathcal I,\xi}}\!\!e^{-m^2 t}\,\mathd t\right]=\frac{1}{m^2}\,\mathbb E\left[1-e^{-m^2\tau^{\mathcal I,\xi}}\right]\leq\frac{1}{m^2}.
\end{equation}
In addition, recalling that condition \ref{H0} gives that $\bar\beta(\xi''-\xi')\lesssim|\xi''-\xi'|^{-1-2s_2}\lesssim (1+|\xi'|^2)^{-\frac{1+2s_2}{2}}$ uniformly for $\xi''\in\mathcal I$, by \eqref{eq:IkedaWatanabe-ac} we obtain
\[
p^{\mathcal I}_{\free}(\xi,\xi')\lesssim\big(1+|\xi'|^2\big)^{-\frac{1+2s_2}{2}}\int_{\mathcal I}G^{\mathcal I}_{\free}(\xi,\xi'')\mathd\xi''\leq\frac{1}{m^2}\big(1+|\xi'|^2\big)^{-\frac{1+2s_2}{2}},
\]
uniformly in $\xi\in\mathcal I$ and $\xi'\in\mathcal I^c\setminus\mathcal N$.
Since $1<r<1+2s_2$, there exists $C_0=C_0(r,k,s_2,\mathcal{I})$ such that
\[
\int_{\mathcal I^c\setminus\mathcal N}\abs{y(\xi')}\,p^{\mathcal I}_{\free}(\xi,\xi')\,\mathd\xi'\leq \frac{C_0}{m^2}\,\norm{y}_{L^1_{r,k}(\mathcal I^c)},
\]
uniformly in $\xi\in\mathcal I$. Hence by \eqref{eq:cont-split} we get
\[
\norm{\varphi^y}_{L^\infty(\mathcal I)}\lesssim\norm{y}_{L^\infty(\mathcal N)}+\frac{C_0}{m^2}\,\norm{y}_{L^1_{r,k}(\mathcal I^c)},
\]
and taking into account $\norm{\varphi^y}_{L^\infty(K\cap\mathcal I^c)}=\norm{y}_{L^\infty(K\cap\mathcal I^c)}$ we conclude
\begin{equation}\label{eq:harmonic-continuity-bound}
\norm{\varphi^y}_{L^\infty(K)}\leq C_K\,\big(\norm{y}_{L^1_{r,k}(\mathcal I^c)}+\norm{y}_{L^\infty(K\cap\mathcal I^c)}+\norm{y}_{L^\infty(\mathcal N)}\big).
\end{equation}
It remains to prove continuity and weak harmonicity. Let $\{y_n\}_{n\in\N}\subseteq C^\infty_c(\R)$ be such that $\norm{y-y_n}_{C(\mathcal N)}+\norm{y-y_n}_{L^1_{r,k}(\mathcal I^c)}\rightarrow0$ as $n\rightarrow+\infty$. Since $y_n\in\mathcal F$ for every $n\in\N$, by Lemma \ref{lemma:harmonicextensionenergy} the function $\varphi^{y_n}$ is the harmonic extension of $y_n$ and, by \eqref{eq:energy-embedding}, $\varphi^{y_n}\in C(\R)$. Applying \eqref{eq:harmonic-continuity-bound} to $y-y_n$ we deduce that $\varphi^{y_n}\rightarrow\varphi^y$ uniformly on $[a,b]$ as $n\rightarrow+\infty$, whence $\varphi^y$ is continuous on $\overline{\mathcal I}$; since $\varphi^y=y$ on $\mathcal I^c$ and $y$ is continuous, $\varphi^y\in C(\R)$.

Finally, $\varphi^y\in L^1_{r,k}(\R)$, since $\varphi^y=y$ on $\mathcal I^c$ and $\varphi^y$ is bounded on $\mathcal I$. Let $h\in C^\infty_c(\mathcal I)$. Since $\varphi^{y_n}-\varphi^y=y_n-y$ on $\mathcal I^c$, combining \eqref{eq:harmonic-continuity-bound} with \eqref{eq:testdecay} we get 
\[
\int_\R\big\lvert\varphi^{y_n}-\varphi^{y}\big\rvert\,\big\lvert(H^{\beta}+m^2)h\big\rvert\,\mathd\xi
\leq C_h\Big(\norm{\rho_{r,k}}_{L^1(\R)}\norm{\varphi^{y_n}-\varphi^{y}}_{L^\infty(\mathcal I)}+\norm{y_n-y}_{L^1_{r,k}(\mathcal I^c)}\Big)\longrightarrow0 .
\]
Since $\varphi^{y_n}$ is a weak harmonic extension of $y_n$, we pass to the limit in $\int_\R\varphi^{y_n}(H^{\beta}+m^2)h\mathd\xi=0$ and we conclude that $\varphi^y$ is a weak harmonic extension of $y$. Finally we prove \eqref{eq:weakharmonicF}. Let $u\in\mathcal F\cap{\rm Dom}(H^{\beta})$ be such that $\supp(u)\subseteq\overline{\mathcal I}$ and
$\sup_{\xi\in\mathcal I^c}\rho_{r,k}(\xi)^{-1}\lvert(H^{\beta}+m^2)u(\xi)\rvert<+\infty$. We note that $u$ belongs to $\mathcal F_{\mathcal I}=\overline{\{h\in\mathcal F:\supp(h)\subset\mathcal I\}}^{\|\cdot\|_\mathcal{F}}$. Since $y_n\in\mathcal F$, by Lemma \ref{lemma:harmonicextensionenergy} the function $\varphi^{y_n}$ is an energy harmonic extension of $y_n$, so that $\mathcal E(\varphi^{y_n},u)=0$; as $u\in{\rm Dom}(H^{\beta})$ and $\varphi^{y_n}\in\mathcal F$, hence
\[
\int_\R\varphi^{y_n}\,(H^{\beta}+m^2)u\,\mathd\xi=0 .
\]
Moreover $\varphi^{y_n}-\varphi^y=y_n-y$ on $\mathcal I^c$, hence, setting $\Lambda_u:=\sup_{\xi\in\mathcal I^c}\rho_{r,k}(\xi)^{-1}\lvert(H^{\beta}+m^2)u(\xi)\rvert$ and using that $(H^{\beta}+m^2)u\in L^2(\mathcal I)\subseteq L^1(\mathcal I)$, we obtain
\[
\int_\R\big\lvert\varphi^{y_n}-\varphi^{y}\big\rvert\,\big\lvert(H^{\beta}+m^2)u\big\rvert\,\mathd\xi
\leq\norm{(H^{\beta}+m^2)u}_{L^1(\mathcal I)}\norm{\varphi^{y_n}-\varphi^{y}}_{L^\infty(\mathcal I)}
+\Lambda_u\norm{y_n-y}_{L^1_{r,k}(\mathcal I^c)} ,
\]
and both terms tend to $0$ as $n\rightarrow+\infty$ by \eqref{eq:harmonic-continuity-bound}.
\end{proof}

\begin{rmk}\label{rmk:weakharmonicF}
Estimate \eqref{eq:weakharmonicF} is strictly stronger than point 3. of Definition \ref{def:harmonicextension}. Indeed passing from test functions in $C^\infty_c(\mathcal I)$ to test functions in $\mathcal F_{\mathcal I}$ by approximation would require convergence of $(H^{\beta}+m^2)u_n$ in $L^\infty_{\rho_{r,k}^{-1}}(\R)$ (for a positive and continuous function $\rho$ on $\R$, $L^\infty_\rho(\R)$, or simply $L^\infty_\rho$, is the space of the measurable functions $g$ such that $\norm{g}_{L^\infty_\rho}:={\rm ess\,sup}_{\xi\in\R}\abs{g(\xi)}\rho(\xi)<+\infty$), whereas only convergence in $\mathcal F$ is available. It is instead obtained by approximating $y$ rather than $u$, so that point 1. of Definition \ref{def:harmonicextension}, whose test functions already belong to $\mathcal F_{\mathcal I}$, applies to $\varphi^{y_n}\in\mathcal F$.
\end{rmk}

\begin{rmk}\label{rmk:harmonicextensionnotions}
Lemma \ref{lemma:harmonicextensionenergy} requires $y$ to be the restriction of a function of $\mathcal F$, which is a strong assumption, indeed for instance, if $y\equiv1$ then $\varphi^y\equiv1$ on $\mathcal I^c$, so that $\varphi^y\notin L^2$ and condition 1. of Definition \ref{def:harmonicextension} cannot hold. Nevertheless, by Theorem \ref{thm:harmonic-continuity}, for every $y\in C(\mathcal I^c,\R)\cap L^1_{r,k}(\mathcal I^c)$ with $1<r<1+2s_2$ the probabilistic harmonic extension of $y$ is well defined and is also a weak harmonic extension of $y$. Accordingly, for such $y$ we call $\varphi^y$ the harmonic extension of $y$ without any further qualification.
\end{rmk}

\begin{remark}\label{rmk:continuityharmonic}
With reasoning similar to the one in the proof of Theorem \ref{thm:harmonic-continuity} it is possible to prove that $\varphi^y$ belongs to $H^{s_1}_{\loc}(\mathcal I)$.
\end{remark}

\subsection{Some useful estimates about Harmonic extension}
In this subsection we collect some useful properties of the harmonic extension. First, we prove a lemma about the $L^p$-norm of the harmonic extension that will be fundamental for the a priori estimate of Theorem \ref{thm:apriori}.
\begin{lemm}\label{lemma:harmonicaverage}
Assume \ref{H0} and \ref{H3} hold. Assume $m>0$, $1<r<1+2s_2$ and $k\in (0,k_0)$ where $k_0:=k_0(r,m)$ is a constant given by Lemma \ref{lemma:semigroupLp} such that 
  \begin{equation}\label{eq:dissL1}
 m^2-\omega_{r,k}>0,\qquad \omega_{r,k}=\sup_{\xi\in\R}\left(-\,\frac{H^{\beta}\rho_{r,k}(\xi)}{\rho_{r,k}(\xi)}\right).
    \end{equation}
For every $p\geq1$ there exists a constant $\mathfrak{C}:=\mathfrak{C}(m,r,k)>0$ such that, for every $n\in\N$ and $x\in L^p_{r,k}\cap C(\R)$, we have
\begin{equation}\label{eq:harmonicaverage}
\int_n^{n+1}\int^T_{-T}\abs{\varphi^{x|_{\mathcal I_T^c}}(\xi)}^p\rho_{r,k}(\xi)\,\mathd\xi\,\mathd T\leq \mathfrak{C} \int_{\{\abs{\eta}\geq n\}}\abs{x(\eta)}^p\rho_{r,k}(\eta)\,\mathd\eta,
\end{equation}
where, for $T>0$, $\varphi^{x|_{\mathcal I_T^c}}$ is the harmonic extension of $x|_{\mathcal I_T^c}$ in $\mathcal I_T:=(-T,T)$. In particular, the left-hand side of \eqref{eq:harmonicaverage} vanishes as $n\rightarrow+\infty$.
\end{lemm}

\begin{rmk}\label{rmk:harmonicaverage}
The average with respect to $T\in[n,n+1]$ in \eqref{eq:harmonicaverage} is essential. Indeed, for a fixed $T>0$, the norm of $\varphi^{x|_{\mathcal I_T^c}}$ in $L^p_{r,k}(\mathcal I_T)$ cannot be estimated by the norm of $x$ in $L^p_{r,k}(\mathcal I_T^c)$, since it depends on the values of $x$ at the boundary of $\mathcal I_T$. This is evident in the case $H^\beta=-\Delta$, namely $\beta_0=1$ and $\bar\beta\equiv0$, where the harmonic extension is explicit:
\[
\varphi^{x|_{\mathcal I_T^c}}(\xi)=\frac{x(-T)\sinh\big(m(T-\xi)\big)+x(T)\sinh\big(m(T+\xi)\big)}{\sinh(2mT)},\qquad \xi\in\mathcal I_T .
\]
Hence $\varphi^{x|_{\mathcal I_T^c}}$ only depends on $x(-T)$ and $x(T)$, and it is easy to see that there exists $C>0$, independent of $T\geq1$ and $x$, such that
\[
\int_{\mathcal I_T}\abs{\varphi^{x|_{\mathcal I_T^c}}(\xi)}^p\rho_{r,k}(\xi)\,\mathd\xi\leq C\big(\abs{x(T)}^p\rho_{r,k}(T)+\abs{x(-T)}^p\rho_{r,k}(-T)\big).
\]
However, the left-hand side cannot be estimated in terms of the norm of $x$ in $L^p_{r,k}(\mathcal I_T^c)$, and a supremum norm of $x$ near $\partial\mathcal I_T$ would be needed. Indeed, let $T\geq1$ and, for $j\in\N$, let $x_j\in C(\R)$ be such that $x_j=0$ on $(-\infty,-T]$ and $x_j(\eta)=\max\{0,1-j(\eta-T)\}$ for $\eta\geq T$. Since $\rho_{r,k}$ is even and decreasing on $[0,+\infty)$, we have
\[
\norm{x_j}^p_{L^p_{r,k}(\mathcal I_T^c)}=\int_T^{T+\frac1j}\big(1-j(\eta-T)\big)^p\rho_{r,k}(\eta)\,\mathd\eta\leq\frac{\rho_{r,k}(T)}{(p+1)j}.
\]
On the other hand, since $x_j(-T)=0$ and $x_j(T)=1$, for every $\xi\in[T-1,T)$, using that $T+\xi\geq1$ and $T-\xi\leq1$, we have
\[
\varphi^{x_j|_{\mathcal I_T^c}}(\xi)=\frac{\sinh\big(m(T+\xi)\big)}{\sinh(2mT)}=e^{-m(T-\xi)}\,\frac{1-e^{-2m(T+\xi)}}{1-e^{-4mT}}\geq e^{-m}\big(1-e^{-2m}\big),
\]
so that
\begin{align*}
\int_{\mathcal I_T}\abs{\varphi^{x_j|_{\mathcal I_T^c}}(\xi)}^p\rho_{r,k}(\xi)\,\mathd\xi&\geq e^{-pm}\big(1-e^{-2m}\big)^p\int_{T-1}^T\rho_{r,k}(\xi)\,\mathd\xi\geq e^{-pm}\big(1-e^{-2m}\big)^p\rho_{r,k}(T)\\
&\geq(p+1)\,j\,e^{-pm}\big(1-e^{-2m}\big)^p\norm{x_j}^p_{L^p_{r,k}(\mathcal I_T^c)},
\end{align*}
where $(p+1)\,j\,e^{-pm}(1-e^{-2m})^p\rightarrow+\infty$ as $j\rightarrow+\infty$. Hence, the average with respect to $T$ removes the dependence on the boundary values, indeed, by the previous estimate, for every $n\geq1$ we have
\[
\int_n^{n+1}\int_{\mathcal I_T}\abs{\varphi^{x|_{\mathcal I_T^c}}(\xi)}^p\rho_{r,k}(\xi)\,\mathd\xi\,\mathd T\leq C\int_{\{n\leq\abs{\eta}\leq n+1\}}\abs{x(\eta)}^p\rho_{r,k}(\eta)\,\mathd\eta .
\]
In the general case the boundary values are replaced by the exit distribution $P^{\mathcal I_T}_{\free}(\xi,\cdot)$, which, for $\xi$ close to $\partial\mathcal I_T$, is concentrated near $\partial\mathcal I_T$ (and may have atoms at $\pm T$ if $\beta_0>0$, see \eqref{eq:IkedaWatanabeformula}).
\end{rmk}
\begin{proof}
Fix $n\in\N$ and $x\in L^p_{r,k}\cap C(\R)$, and set $\rho:=\rho_{r,k}$ and $\varphi_T:=\varphi^{x|_{\mathcal I_T^c}}$ for $T>0$. For every open interval $D\subseteq\R$, possibly unbounded, and $\xi\in\R$ we define $P^D_{\free}(\xi,\cdot)$ as in \eqref{eq:poissonkernel}, where we recall that $\{Y^\xi_t\}_{t\geq0}$ is the symmetric L\'evy process associated with $-H^\beta$ and starting from $\xi$ (see Section \ref{section:harmonicextension}), that $\tau^{D,\xi}=\inf\{t>0:Y_t^\xi\in D^c\}$ is its first exit time from $D$, and we use the convention $e^{-m^2\tau^{D,\xi}}=0$ on $\{\tau^{D,\xi}=+\infty\}$.\\

\noindent First, we prove the following preliminary estimate that will be essential in the proof. Let $\sigma\in\R$. Since $\rho^\sigma_{r,k}\in C^2_0(\R)$, by \ref{H3}, \eqref{eq:excessive} and Remark \ref{rmk:weighttranslation} we have
\[
(H^\beta+m^2)\rho^\sigma_{r,k}\geq(m^2-\omega_{r,k})\rho^\sigma_{r,k}\geq0.
\]
Hence, for every open interval $D$ and $\xi\in D$, by the Dynkin formula applied to $\{e^{-m^2s}\rho^\sigma_{r,k}(Y^\xi_s)\}_{s\geq0}$ and to the bounded stopping times $\tau^{D,\xi}\wedge t$, we have
\[
\mathbb E\Big[e^{-m^2(\tau^{D,\xi}\wedge t)}\rho^\sigma_{r,k}\big(Y^\xi_{\tau^{D,\xi}\wedge t}\big)\Big]\leq\rho^\sigma_{r,k}(\xi),\qquad t>0,
\]
letting $t\rightarrow+\infty$, by the Fatou lemma we get
\begin{equation}\label{eq:dynkinweight}
\int_{D^c}\rho^\sigma_{r,k}(\eta)P^D_{\free}(\xi,\mathd\eta)=\mathbb E\Big[e^{-m^2\tau^{D,\xi}}\rho^\sigma_{r,k}\big(Y^\xi_{\tau^{D,\xi}}\big)\Big]\leq\rho^\sigma_{r,k}(\xi),\qquad \xi\in D,\ \sigma\in\R .
\end{equation}
Moreover, we have
\begin{equation}\label{eq:weightconvolution}
\int_\R\rho_{r,k}(\zeta+u)\rho_{r,k}(\zeta)\,\mathd\zeta\leq S_1\rho_{r,k}(u),\qquad u\in\R,\qquad S_1:=2\cdot4^{r/2}c_rk^{-1/2}.
\end{equation}
Indeed, if $\abs{\zeta}\leq\frac12\abs{u}$, then $\abs{u}\leq2\abs{\zeta+u}$ and so $\rho_{r,k}(\zeta+u)\leq4^{r/2}\rho_{r,k}(u)$ by \eqref{eq:rhodoubling}, while if $\abs{\zeta}>\frac12\abs{u}$, then $\rho_{r,k}(\zeta)\leq4^{r/2}\rho_{r,k}(u)$ by \eqref{eq:rhodoubling}; hence \eqref{eq:weightconvolution} follows from \eqref{eq:rhoL1}.\\

\noindent Now we can estimate the weighted $L^p$ norm of the harmonic extension. By \eqref{eq:harmonicextensionprob}, \eqref{eq:poissonkernel}, and the H\"older inequality, we have
\begin{equation}\label{eq:harmonicaverage1}
\int_{\mathcal I_T}\abs{\varphi_T(\xi)}^p\rho(\xi)\,\mathd\xi\leq\int_{\mathcal I_T}\rho(\xi)\int_{[T,+\infty)}\abs{x(\eta)}^pP^{\mathcal I_T}_{\free}(\xi,\mathd\eta)\,\mathd\xi+\int_{\mathcal I_T}\rho(\xi)\int_{(-\infty,-T]}\abs{x(\eta)}^pP^{\mathcal I_T}_{\free}(\xi,\mathd\eta)\,\mathd\xi .
\end{equation}

We estimate the first term in the right-hand side of \eqref{eq:harmonicaverage1}. Let $T>0$, $\xi\in\mathcal I_T$ and let $B\subseteq[T,+\infty)$ be a Borel set. Since $\mathcal I_T\subseteq(-\infty,T)$, we have $\tau^{\mathcal I_T,\xi}\leq\tau^{(-\infty,T),\xi}$; on the other hand, on the event $\{Y^\xi_{\tau^{\mathcal I_T,\xi}}\in B\}$ we have $Y^\xi_{\tau^{\mathcal I_T,\xi}}\geq T$, so that $\tau^{(-\infty,T),\xi}\leq\tau^{\mathcal I_T,\xi}$. Hence on this event the two exit times and the two exit positions coincide, and
\[
P^{\mathcal I_T}_{\free}(\xi,B)\leq P^{(-\infty,T)}_{\free}(\xi,B)=P^{(-\infty,0)}_{\free}(\xi-T,B-T),
\]
where the last equality follows from the fact that the law of $\{Y^\xi_t-\xi\}_{t\geq0}$ does not depend on $\xi$. Therefore, setting $\xi=T+\zeta$ and $\eta=T+\delta$, and using \eqref{eq:rhocomparable} (recall that $k\leq1$), for every $T\in[n,n+1]$ we have
\[
\int_{\mathcal I_T}\rho(\xi)\int_{[T,+\infty)}\abs{x(\eta)}^pP^{\mathcal I_T}_{\free}(\xi,\mathd\eta)\,\mathd\xi\leq e^{r/2}\int_{-\infty}^0\rho(n+\zeta)\int_{[0,+\infty)}\abs{x(T+\delta)}^pP^{(-\infty,0)}_{\free}(\zeta,\mathd\delta)\,\mathd\zeta.
\]
Let $\delta\geq0$. If $s\in[n+\delta,n+\delta+1]$, then $0\leq s-n-\delta\leq1$, so that $\rho^{s-n}_{r,k}(\delta)\geq(1+k)^{-r/2}$; hence
\[
\int_n^{n+1}\abs{x(T+\delta)}^p\,\mathd T=\int_{n+\delta}^{n+\delta+1}\abs{x(s)}^p\,\mathd s\leq(1+k)^{r/2}\int_n^{+\infty}\abs{x(s)}^p\rho^{s-n}_{r,k}(\delta)\,\mathd s .
\]
By the Tonelli theorem, the last inequality, \eqref{eq:dynkinweight} with $D=(-\infty,0)$ and $\sigma=s-n$, and \eqref{eq:weightconvolution} with $u=s$, we obtain
\begin{align*}
\int_{\mathcal I_T}\rho(\xi)\int_{[T,+\infty)}&\abs{x(\eta)}^pP^{\mathcal I_T}_{\free}(\xi,\mathd\eta)\,\mathd\xi\leq e^{r/2}\int_{-\infty}^0\rho(n+\zeta)\int_{[0,+\infty)}\abs{x(T+\delta)}^pP^{(-\infty,0)}_{\free}(\zeta,\mathd\delta)\,\mathd\zeta\\
&\leq e^{r/2}(1+k)^{r/2}\int_n^{+\infty}\abs{x(s)}^p\int_{-\infty}^0\rho(n+\zeta)\int_{[0,+\infty)}\rho^{s-n}_{r,k}(\delta)P^{(-\infty,0)}_{\free}(\zeta,\mathd\delta)\,\mathd\zeta\,\mathd s\\
&\leq e^{r/2}(1+k)^{r/2}\int_n^{+\infty}\abs{x(s)}^p\int_\R\rho_{r,k}(n+\zeta)\rho_{r,k}(n+\zeta-s)\,\mathd\zeta\,\mathd s\\
&\leq \mathfrak{C}\int_n^{+\infty}\abs{x(s)}^p\rho(s)\,\mathd s,
\end{align*}
where $\mathfrak{C}:=e^{r/2}(1+k)^{r/2}S_1$ does not depend on $n$, $p$ and $x$.

Since $\{Y^\xi_t\}$ is a symmetric L\'evy process, we have $P^{\mathcal I_T}_{\free}(\xi,B)=P^{\mathcal I_T}_{\free}(-\xi,-B)$ for every $\xi\in\mathcal I_T$ and every Borel set $B\subseteq\mathcal I_T^c$; hence, since $\rho$ is even, the second term in the right-hand side of \eqref{eq:harmonicaverage1} coincides with the first one with $x$ replaced by $x(-\cdot)$, and so, we obtain
\begin{equation}\label{eq:harmonicaverage2}
\int_{\mathcal I_T}\abs{\varphi_T(\xi)}^p\rho(\xi)\,\mathd\xi\leq \mathfrak{C}\int_n^{+\infty}\abs{x(s)}^p\rho(s)\,\mathd s
\end{equation}
Finally, integrating \eqref{eq:harmonicaverage2} with respect to $T\in[n,n+1]$, we obtain \eqref{eq:harmonicaverage}. 
\end{proof}

The next lemma shows that a harmonic extension in a bounded open interval $\mathcal I$, even when its exterior data have no finite global energy, has finite energy after multiplication by a smooth cutoff with compact support in $\mathcal I$; this is what makes such products admissible Cameron--Martin directions in Section \ref{section:uniquenessDLR}. The bound contains only derivatives and differences of the cutoff; the mass term cancels by harmonicity.

\begin{lemm}\label{lemma:localization}
Assume \ref{H0} holds and let $1<r<1+2s_2$ and $k>0$. Let $\mathcal I\subset\R$ be a bounded open interval, let $\chi\in C^\infty_c(\mathcal I)$ and $y\in C(\mathcal I^c)\cap L^1_{r,k}(\mathcal I^c)$, and let $\varphi^y$ be the harmonic extension of $y$ in $\mathcal I$. Then $\chi\varphi^y\in\mathcal F_{\mathcal I}$, where $\mathcal F_{\mathcal I}$ is defined in the proof of Lemma \ref{lemma:harmonicextensionenergy}, and
\begin{equation}\label{eq:localization}
\mathcal E(\chi\varphi^y,\chi\varphi^y)\leq\beta_0\int_\R\abs{\chi'}^2\abs{\varphi^y}^2\mathd\xi+\frac12\int_{\R^2}\big(\chi(\xi)-\chi(\xi')\big)^2\abs{\varphi^y(\xi)}\abs{\varphi^y(\xi')}\bar\beta(\xi-\xi')\,\mathd(\xi,\xi').
\end{equation}
\end{lemm}

\begin{proof}
Let $\mathcal I=(a,b)$ and let $\{y_j\}_{j\in\N}\subseteq C^\infty_c(\R)$ be such that $y_j\rightarrow y$ in $L^1_{r,k}(\mathcal I^c)$ and uniformly on the compact subsets of $\mathcal I^c$. By Lemma \ref{lemma:harmonicextensionenergy} we have $\varphi^{y_j}\in\mathcal F$ and $\mathcal E(\varphi^{y_j},\chi^2\varphi^{y_j})=0$, since $\chi^2\varphi^{y_j}\in\mathcal F$ by \eqref{eq:energy-is-Hs1} and $\supp(\chi^2\varphi^{y_j})\subset\mathcal I$. The following product identities hold for smooth $u$ and extend to the form domain by approximation (the identity involving derivatives is used only when $\beta_0>0$):
\begin{align*}
&\abs{(\chi u)'}^2-u'(\chi^2u)'=\abs{\chi'}^2u^2,\\
&\big(\chi(\xi)u(\xi)-\chi(\xi')u(\xi')\big)^2-\big(u(\xi)-u(\xi')\big)\big(\chi^2(\xi)u(\xi)-\chi^2(\xi')u(\xi')\big)=\big(\chi(\xi)-\chi(\xi')\big)^2u(\xi)u(\xi'),
\end{align*}
hence, by \eqref{eq:dirichletform},
\begin{align}\label{eq:IMSj}
\mathcal E(\chi\varphi^{y_j},\chi\varphi^{y_j})&=\mathcal E(\chi\varphi^{y_j},\chi\varphi^{y_j})-\mathcal E(\varphi^{y_j},\chi^2\varphi^{y_j})\notag\\
&=\beta_0\int_\R\abs{\chi'}^2\abs{\varphi^{y_j}}^2\mathd\xi+\frac12\int_{\R^2}\big(\chi(\xi)-\chi(\xi')\big)^2\varphi^{y_j}(\xi)\varphi^{y_j}(\xi')\bar\beta(\xi-\xi')\,\mathd(\xi,\xi').
\end{align}
By Theorem \ref{thm:harmonic-continuity}, $\varphi^{y_j}\rightarrow\varphi^y$ uniformly on $[a-1,b+1]$, so that $\chi\varphi^{y_j}\rightarrow\chi\varphi^y$ in $L^2(\R)$, and we can pass to the limit in the right-hand side of \eqref{eq:IMSj}, where the integrand of the second term vanishes unless $\xi$ or $\xi'$ belongs to $\supp(\chi)$. Indeed, on $[a-1,b+1]^2$ this integrand is bounded by $c\min\{1,\abs{\xi-\xi'}^2\}\bar\beta(\xi-\xi')$, with $c>0$ independent of $j$, which is integrable on $[a-1,b+1]^2$ since $\bar\beta$ is a L\'evy density; if instead $\xi\in\supp(\chi)$ and $\xi'\notin[a-1,b+1]$, then, as in the proof of Theorem \ref{thm:harmonic-continuity}, $\bar\beta(\xi-\xi')\leq c\rho_{r,k}(\xi')$ with $c>0$ independent of $\xi$ and $\xi'$, and $\varphi^{y_j}(\xi')=y_j(\xi')\rightarrow y(\xi')$ in $L^1_{r,k}(\mathcal I^c)$ (and symmetrically if $\xi'\in\supp(\chi)$ and $\xi\notin[a-1,b+1]$). Hence the right-hand side of \eqref{eq:IMSj} converges to a quantity which is not greater than the right-hand side of \eqref{eq:localization}. Since $\{\chi\varphi^{y_j}\}_{j\in\N}$ is bounded in the Hilbert space $\mathcal F_{\mathcal I}$ and converges to $\chi\varphi^y$ in $L^2(\R)$, we conclude that $\chi\varphi^y\in\mathcal F_{\mathcal I}$ and that $\mathcal E(\chi\varphi^y,\chi\varphi^y)\leq\liminf_{j\rightarrow+\infty}\mathcal E(\chi\varphi^{y_j},\chi\varphi^{y_j})$, which gives \eqref{eq:localization}.
\end{proof}

In fact, we can prove something stronger: that the cut-off harmonic extension belongs to the domain of $H^\beta$. This will be used as well in Section \ref{section:uniquenessDLR}.

\begin{proposition}\label{proposition::cutoffdomain}
In the setting of Lemma \ref{lemma:localization}, we have also that $\chi\varphi^{y}\in \operatorname{Dom}(H^\beta)$.  
\end{proposition}
\begin{proof}
Since $\chi$ is supported in $\mathcal{I}$, we have by weak harmonicity that (at least at the level of distributions)
\begin{align*}
    \left(H^\beta+m^2\right)(\chi\varphi^y)=\left[H^\beta,\chi\right]\varphi^y,
\end{align*}
so we are done as long as we can prove that the commutator satisfies $\left[H^\beta,\chi\right]\varphi^y\in L^2(\mathbb{R})$. Consider $\widetilde{\chi}\in C_c^{\infty}(\mathcal{I})$ such that $\operatorname{supp} \chi\Subset \operatorname{supp} \widetilde{\chi}\Subset\mathcal{I}$ and $\widetilde{\chi}=1$ in a neighborhood of $\operatorname{supp}\chi$, so that
\begin{equation*}
    \left[H^\beta,\chi\right]\varphi^y=\left[H^\beta,\chi\right](\widetilde{\chi}\varphi^y)+\left[H^\beta,\chi\right]((1-\widetilde{\chi})\varphi^y).
\end{equation*}
The first term is equal to
\begin{align*}
    \left[H^\beta,\chi\right](\widetilde{\chi}\varphi^y)=&-\beta_0(\chi''(\widetilde{\chi}\varphi^y)+2\chi'(\widetilde{\chi}\varphi^y)')+H^\beta_{nl}\chi(\widetilde{\chi}\varphi^y)\\
    &+\int_{\mathbb{R}}\left(\chi(\cdot)-\chi(\xi')\right)\left(\widetilde{\chi}\varphi^y(\cdot)-\widetilde{\chi}\varphi^y(\xi')\right)\bar{\beta}(\cdot-\xi')\mathrm{d}\xi'.
\end{align*}
Since by Lemma \ref{lemma:localization} we know that $\widetilde{\chi}\varphi^y\subset\mathcal{F}_\mathcal{I}\subset\mathcal{F}$, we can conclude in a very similar way that $\left[H^\beta,\chi\right](\widetilde{\chi}\varphi^y)\in L^2(\mathbb{R})$. For the second term, notice that $\chi(1-\widetilde{\chi})=0$ implies
\begin{equation}
    \left[H^\beta,\chi\right]((1-\widetilde{\chi})\varphi^y)=\chi(\cdot)\int_{\mathbb{R}}((1-\tilde{\chi})\varphi^y)(\xi')\bar{\beta}(\cdot-\xi')\mathrm{d}\xi',
\end{equation}
and that by the choice of cut-offs there exists $C>0$ such that
\begin{equation*}
    \bar{\beta}(\xi-\xi')\leqslant C\rho_{r,k}(\xi')\quad\text{whenever }\xi\in\operatorname{supp}\chi\text{ and }\xi'\in \operatorname{supp}(1-\widetilde{\chi}). 
\end{equation*}
Then
\begin{equation*}
    \left|\left[H^\beta,\chi\right]((1-\widetilde{\chi})\varphi^y)(\xi)\right|\leqslant C\|(1-\widetilde{\chi})\varphi^y\|_{L^1_{r,k}}\chi(\xi),
\end{equation*}
so the second commutator is also in $L^2(\mathbb{R})$, as we wanted.
\end{proof}

\section{Reflection positive Gaussian measures and the operators $H^{\beta}$}\label{appendix:reflectionpositivity}

In this section we discuss the notion of reflection-positive probability measures and we characterize, in terms of the operators $H^{\beta}$, a natural class of reflection-positive Gaussian measures (see Theorem \ref{theorem:rprepresentation2}). This characterization will be used in Subsection \ref{section:ref.positive} to prove that, in the setting discussed in this paper, there is always a reflection-positive measure solving the DLR equations \eqref{eq:DLRequations}.\\

\noindent Our discussion is about reflection-positive measures in $\R^n$. Let $\theta: \mathcal{S}(\R^n) \rightarrow \mathcal{S}(\R^n)$ be the reflection map with respect to the hyperplane $\xi_1=0$ of  $\R^n$ applied to Schwartz functions, namely, if $f \in \mathcal{S}(\R^n)$, then 
\[\theta(f)(\xi_1,...,\xi_n)=f(-\xi_1,...,\xi_n).\]
By duality we can extend $\theta$ to the space of tempered distributions $\mathcal{S}'(\R^n)$, namely, if $\mathcal{D}\in \mathcal{S}'(\R^n)$ we have 
\[\langle \theta(\mathcal{D}), f \rangle =\langle \mathcal{D}, \theta(f) \rangle.\]

\begin{definition}\label{definition:reflectionpositivity}
Let $\mu$ be a probability measure on $\mathcal{S}'(\R^n)$. We say that $\mu$ is reflection positive if for any random distribution $X:(\Omega,\mathcal{F},\mathbb{P}) \rightarrow \mathcal{S}'(\R^n)$ with law $\Law(X)=\mu$, we have that for any cylinder function $\Phi$ such that $\operatorname{supp}(\Phi) \subset (0,+\infty)\times \R^{n-1}$ we have 
\begin{equation}\label{eq:rpdefinition}
\int\Phi(\varphi)\Phi(\theta(\varphi))\mu(\mathd \varphi)=\mathbb{E}\left[\Phi(X) \Phi(\theta(X)) \right] \geq  0.
\end{equation}
We also say that the random distribution $X$ is reflection positive. 
\end{definition}

\noindent Reflection positivity was introduced by Osterwalder and Schrader in \cite{OsterwalderSchrader1973} as the fundamental positivity property in Euclidean quantum field theory needed to reconstruct a physical Hilbert space and a positive Hamiltonian from Euclidean correlation functions. Reflection positivity later found important applications in statistical mechanics. It is the basis of powerful methods such as infrared bounds and chessboard estimates, which have been used to establish long-range order and phase transitions (see, e.g., \cite{FrohlichSimonSpencer1976,FrohlichEtAl1978}). More generally, reflection positivity connects probability theory, operator semigroups, harmonic analysis, and the representation theory of Lie groups (see, e.g., \cite{NeebOlafsson2018} for all these applications). See \cite{Jaffe2017} for a historical and more recent review on the applications of this property.\\

In what follows we consider measures $\mu$ on $\mathcal{S}'(\R^n)$ which are reflection positive, invariant (i.e. symmetric with respect to the action of the Euclidean group of $\R^n$) and Gaussian with mean $0$. If a measure is reflection positive and invariant we can use any reflection with respect to any hyperplane to define reflection positivity and $\theta$ does not play any special role (apart from being the reflection with respect to some hyperplane of $\R^n$).

\noindent Centered Gaussian measures with these properties are in one-to-one correspondence with bilinear forms $\Sigma:\mathcal{S}(\R^n) \times \mathcal{S}(\R^n) \rightarrow \R$, continuous with respect to the standard topology of the space of Schwartz functions $\mathcal{S}(\R^n)$. The form $\Sigma$ is the covariance of the measure $\mu$: if $X$ is a Gaussian random distribution with law $\mu$ then
\[\Sigma(f,g):=\mathbb{E}\left[\langle X , f \rangle \langle X , g \rangle  \right].\]
Since $\mu$ is invariant with respect to translations of $\R^n$, there is a tempered distribution $\bar\Sigma$ such that $\Sigma(f,g)=\scal{\bar\Sigma*f}{g}$. Furthermore, since $\Sigma$ is the covariance of a probability measure, $\bar\Sigma$ is a positive definite function, and thus, by Bochner theorem, there is a (tempered) positive measure $\nu$ on $\R^n$ such that 
\[\widehat{\bar\Sigma}(\mathd k)=\nu(\mathd k).\]
Finally since $\mu$ is invariant with respect to reflections and rotations also $\nu$ is invariant with respect to reflections and rotations.\\

\noindent In this situation we have a complete characterization of the measure $\mu$ (on $\mathcal{S}'(\R^n)$) through a complete characterization of the measure $\nu$ on $\R^n$. The characterization of this measure is very well understood (see, e.g., \cite{Klein1977,NeebOlafsson2015,NeebOlafsson2018,Simon1979FunctionalIntegration}). We prefer to report here the proof with many details, since they are not presented in full detail in the literature. Our main addition is the connection between these measures and the operators associated with subordinated Brownian motion.\\

\noindent Define the function $w_{N,n}:\R_+ \rightarrow \R_+$ where $n,N \in \mathbb{N}$ as 
\begin{equation}
w_{N,n}(s)=\mathbb{I}_{(1/2,+\infty)}(s)\frac{1}{(1+s)^N}+\mathbb{I}_{[0,1/2]}(s)\left\{\begin{array}{ll}
\frac{1}{\sqrt{s}} & \text{if }n=1\\
-\log(s)& \text{if }n=2\\
1 & \text{if }n>2
\end{array}\right.
\end{equation}

\begin{theorem}\label{theorem:rprepresentation1}
A centered Gaussian measure $\mu$ on $\mathcal{S}'(\R^n)$ is invariant with respect to the Euclidean group of $\R^n$ and reflection positive if and only if we have
\begin{equation}\label{eq:rprepresentation1}
\widehat{\bar\Sigma}(\mathd k)=\nu(\mathd k)=c_0 \delta_0(\mathd k)+ h(|k|^2) \mathd k ,\quad k \in \R^n
\end{equation}
where $h \geq  0$ is a Borel function growing at most polynomially at infinity. Furthermore, there is a unique positive measure $\sigma$ on $\R_+$ such that there is $N_0 \in \mathbb{N}$ for which we have
\begin{equation}\label{eq:integrabilitysigma}
\int_{0}^{+\infty}w_{N_0,n}(s)\sigma(\mathd s)<+\infty
\end{equation}
The measure $\sigma$ characterizes the function $h$ in the following way: for any $N \geq  N_0$ and any $q_0 >0$, $h$ can be written in a unique way as 
\begin{equation}\label{eq:formulahNq0}
h(q)=P^{N,q_0}(q)+F^{N,q_0}(q)=P^{N,q_0}(q)+\int_0^{+\infty}\frac{(-1)^{N-1}(q-q_0)^{N-1}}{(q+s)(q_0+s)^{N-1}}\sigma(\mathd s),
\end{equation}
where $P^{N,q_0}(q)$ is a polynomial. 
\end{theorem}

\noindent A sketch of the proof of Theorem \ref{theorem:rprepresentation1} can be found in Theorem 24.2 in \cite{Simon1979FunctionalIntegration} (see also  \cite{NeebOlafsson2015,NeebOlafsson2018}). We expand it a little here.\\ 
We introduce the notation $C^E_s$ to denote the Green function of the operator $-\Delta_{\R^n}+s$, when $s>0$. For $s=0$, we define $C_0^E(\xi)=(-\Delta_{\R^n})^{-1}(\xi)$ for $n>2$, and, conventionally, $C_0^E(\xi)=0$ for $n=1,2$. For any $\xi\not=0$ and $s\geq 0$, $C^E_s(\xi)$ is a well defined function that decays exponentially (for $s>0$) or polynomially (for $s=0$) at infinity.

\begin{lemma}\label{lemma:rprepresentation1}
Let $\sigma$ be a positive measure satisfying the condition \eqref{eq:integrabilitysigma}. Then the function $G^{\sigma}:\R^{n}\backslash\{0\}\rightarrow \R$ defined as 
\begin{equation}\label{eq:Gsigma}
G^{\sigma}(\xi):=\int_0^{+\infty}C_s^E(\xi)\sigma(\mathd s)
\end{equation}
is a well defined smooth function (on $\R^{n}\backslash\{0\}$) and, for any $N\geq N_0$ and $q_0>0$, we have that for any $\xi\not=0$
\begin{equation}\label{eq:Gsigmarepresentation}
G^{\sigma}(\xi)=\Big(\text{inverse Fourier transform of }F^{N,q_0}(|\cdot|^2)\Big)(\xi),
\end{equation}
where $F^{N,q_0}(|k|^2)=\int_0^{+\infty}\frac{(-1)^{N-1}(|k|^2-q_0)^{N-1}}{(|k|^2+s)(q_0+s)^{N-1}}\sigma(\mathd s)$ for $k\in\R^n\backslash\{0\}$.
\end{lemma}
\begin{proof}
By the heat-kernel representation of $C_s^E(\xi)$ (see, e.g.,  Chapter V, Section 3, formula (26) of \cite{Stein1970}) we have that 
\begin{equation}\label{eq:CEexpression}
C^E_s(\xi)=\frac{1}{(4\pi)^{n/2}}\int_0^{+\infty}\frac{e^{-\frac{|\xi|^2}{4t}-st}}{t^{\frac{n}{2}}}\mathd t.
\end{equation}
Using expression \eqref{eq:CEexpression} and Gamma and incomplete Gamma function properties we can prove that, for any compact set $K \subset \R^n \backslash\{0\}$ and any multi-index $\alpha \in \mathbb{N}_0^n$, we have
\[\sup_{\xi\in K}\left|\partial^{\alpha}_\xi C_s^E(\xi)\right|\lesssim_{K,\alpha} \mathbb{I}_{[0,1/2]}(s) w_{N,n}(s)+ \mathbb{I}_{(1/2,+\infty)}(s) e^{-c_K\sqrt{s}}, \]
where $c_K>0$ depends only on $K$ (for instance, $c_K=\frac{1}{2}\min_{\xi\in K}|\xi|$). This means that the integral \eqref{eq:Gsigma} reads
\begin{equation}\label{eq:Gsigma2}
G^{\sigma}(\xi)=\frac{1}{(4\pi)^{n/2}}\int_{\R_+^2}\frac{e^{-\frac{|\xi|^2}{4t}-st}}{t^{\frac{n}{2}}}\mathd t\sigma(\mathd s),
\end{equation}
which is a well defined integral when $\xi$ is in $\R^{n}\backslash\{0\}$. Furthermore, by the dominated convergence theorem, $G^{\sigma}$ is infinitely differentiable for $\xi \in \R^{n}\backslash\{0\}$. This means that the integral \eqref{eq:Gsigma} and all its derivatives with respect to $\xi$ converge uniformly when $\xi$ varies in a compact subset of $\R^{n}\backslash\{0\}$.\\
Finally the integral representation \eqref{eq:Gsigma2}, defines a continuous map $\sigma \longmapsto G^{\sigma} $ from the set of positive measures satisfying the condition \eqref{eq:integrabilitysigma} equipped with the metric 
\begin{equation}\label{eq:metricsigma}
d(\sigma,\sigma')=\int_0^{+\infty}w_{N_0,n}(s)|\sigma-\sigma'|(\mathd s),
\end{equation}
into the Fréchet space $C^{\infty}(\R^{n}\backslash\{0\})$ with the standard seminorms. We use this last observation to prove equation \eqref{eq:Gsigmarepresentation}.\\
Indeed, when $\sigma$ satisfies condition \eqref{eq:integrabilitysigma} with $N_0=1$, equation \eqref{eq:Gsigmarepresentation} trivially holds when $N=1$. For $N>1$ (and always $\sigma$ satisfying condition \eqref{eq:integrabilitysigma} with $N_0=1$), for any $q>0$, we have that 
\begin{equation}\label{eq:identityFNq0}
\int_0^{+\infty}\frac{1}{(q+s)}\sigma(\mathd s)=F^{N,q_0}(q)+\sum_{j=0}^{N-2}(q-q_0)^j\int_0^{+\infty}\frac{(-1)^j}{(q_0+s)^{j+1}}\sigma(\mathd s)= F^{N,q_0}(q)+Q^{N,q_0}(q)
\end{equation}
where $Q^{N,q_0}$ is some polynomial. Since $\sigma$ satisfies condition \eqref{eq:integrabilitysigma}, the function $F^{N,q_0}(|k|^2)$ can be extended to the whole $\R^n$ to be a locally integrable function. Furthermore, for $q\geq 1$ we have
\[|F^{N,q_0}(q)|\leq |q-q_0|^{N-1}\int_0^{+\infty}\frac{1}{(1+s)(q_0+s)^{N-1}}\sigma(\mathd s),\]
where the last integral is finite by \eqref{eq:integrabilitysigma}, so that $F^{N,q_0}$ grows at most polynomially at infinity. Thus, $F^{N,q_0}$ can be seen as a tempered distribution. For this reason, we can apply the inverse Fourier transform of both sides of equation \eqref{eq:identityFNq0}, and, if we use the fact that the inverse Fourier transform of $Q^{N,q_0}(|\cdot|^2)$ is a distribution supported in $0$, we obtain equality \eqref{eq:Gsigmarepresentation} in the case $\int_0^{+\infty}\frac{1}{(q+s)}\sigma(\mathd s)<+\infty$.\\
In the general case, we can always approximate any positive measure $\sigma$ satisfying \eqref{eq:integrabilitysigma} with $N_0 >1$, with a sequence of positive measures $\sigma_\ell$ satisfying the condition \eqref{eq:integrabilitysigma} with $N_0 =1$ and converging to $\sigma$ in the metric \eqref{eq:metricsigma}. The conclusion follows by the continuity of $\sigma \longmapsto G^{\sigma} $ and the representation \eqref{eq:Gsigmarepresentation} under this convergence, and the weak convergence of the approximations $F^{N,q_0,\ell}$, associated with the measures $\sigma_\ell$, as tempered distributions since, for any Schwartz function $f\in \mathcal{S}(\R^n)$, 
\[|\scal{F^{N,q_0,\ell}-F^{N,q_0}}{f}|\lesssim_f \int_0^{+\infty}w_{N,n}(s)|\sigma_\ell-\sigma|(\mathd s).\]
\end{proof}

\begin{proof}[Proof of Theorem \ref{theorem:rprepresentation1}]
Since $\mu$ is a reflection-positive Gaussian measure invariant with respect to the Euclidean group, it satisfies the Osterwalder--Schrader axioms apart from the clustering axiom (see, e.g., axiom (OS5) in Section II.3 of \cite{Simon1974}). Since the clustering property is needed for the uniqueness of the quantum vacuum in Osterwalder--Schrader reconstruction theorem (see, e.g., Remark 2 in Section II.4 of \cite{Simon1974}), this means that the distribution $\bar\Sigma$ is the two-point Schwinger function of a quantum field theory, with possibly multiple vacua, and, thus, the analytic continuation $W_2$ of $\bar\Sigma$ with respect to the first variable is a Wightman function of some (regular) quantum field theory (see, e.g., Chapter II, Theorem II.13 of \cite{Simon1974}). Thus, by K\"allén-Lehmann representation (Theorem IX.34 of \cite{ReedSimon1975II}), $W_2$ is of the form 
\[W_2(\xi)=c+\int_0^{+\infty}C_s^L(\xi)\sigma(\mathd s)\]
where $\sigma$ is some positive measure satisfying \eqref{eq:integrabilitysigma} and 
\[C_s^L(\xi)=\frac{1}{(2\pi)^{n-1}}\int_{\R^{n-1}}\frac{e^{i\xi_1\left(\sqrt{|p|^2+s}\right)+ip\cdot(\xi_2,...,\xi_n)^T}}{2\sqrt{|p|^2+s}}\mathd p\]
i.e. is the analytic continuation of the function $C^{E}_s(\xi)$ when the first variable is purely imaginary. This means, by doing the analytic continuation of $W_2$, we obtain, for any $\xi\in \R^n$ such that $\xi\not=0$, 
\begin{equation}\label{eq:barSigmanot0}
\bar\Sigma(\xi)=c+\int_0^{+\infty}C^E_s(\xi)\sigma(\mathd s)=c+G^{\sigma}(\xi).
\end{equation}
By the representation \eqref{eq:Gsigmarepresentation} of $G^{\sigma}(\xi)$, it is possible to extend $G^{\sigma}$ to a tempered distribution to the whole $\R^n$, by simply asking that, for any Schwartz function $f\in \mathcal{S}(\R^n)$, we have 
\[\scal{G^{\sigma}}{f}=\int_{\R^n}\int_{\R_+}\frac{(-1)^{N-1}(|k|^2-q_0)^{N-1}}{(|k|^2+s)(q_0+s)^{N-1}}\hat{f}(k)\sigma(\mathd s)\mathd k.\]
It is important to note that this extension of $G^{\sigma}$ is not canonical but depends, for example, on the choice of $q_0 \in (0,+\infty)$. The difference between two possible extensions of $G^{\sigma}$ is always of the form $\tilde{P}(-\Delta)\delta_0$, for some polynomial $\tilde{P}$. Our choice of $F^{N,q_0}$ fixes one of these possible extensions.\\
Since $\bar\Sigma-c$ and $G^{\sigma}$ coincide when tested against Schwartz functions supported in $\R^n\backslash\{0\}$, the difference $\bar\Sigma-c-G^{\sigma}$ is a tempered distribution supported in $\{0\}$. This means that there is a polynomial $P$ such that 
\[\bar\Sigma-c-G^{\sigma}=P(-\Delta)\delta_0\]
where we used that $\bar\Sigma-c-G^{\sigma}$ is invariant with respect to rotations and thus $P$ must be a function of $-\Delta$.\\
The converse of the previous statement follows from the fact that a Gaussian measure is reflection positive if and only if its covariance is a reflection positive kernel (see, e.g., \cite{Klein1977}), the fact that $C^E_s$ is a reflection positive kernel, the fact that positive linear combinations (or more generally positive integrals) of reflection positive kernels are reflection positive kernels (see, e.g., the proof of Proposition 8.8.2 of \cite{NeebOlafsson2018}, which gives this result when $\sigma$ is tame but that can be easily extended using the approximation exploited in the proof of Lemma \ref{lemma:rprepresentation1}), the fact that the covariance $\bar\Sigma=(2\pi)^{-n}c_0$ is reflection positive, and the fact that the covariances of the form $\bar\Sigma=P(-\Delta)\delta_0$ are reflection positive being supported on $\{0\}$.
\end{proof}

\noindent Clear and simple examples of Gaussian measures satisfying the equivalent conditions in Theorem \ref{theorem:rprepresentation1} are the derivatives of the Gaussian white noise where $\bar\Sigma=c(-\Delta)^j\delta$. In the case $j=0$ the previous $\bar\Sigma=c\delta$ is the covariance of the Gaussian white noise on $\mathcal{S}'(\R^n)$. We can measure which derivatives are present in $0$ (but also the growth of the function $h$ at infinity in general) in the following way: consider $f \in \mathcal{S}(\R^n)$ we write 
\[f_{\lambda}(\xi):=\lambda^{-n}f\left(\frac{\xi}{\lambda} \right).\]
If $\bar\Sigma=c(-\Delta)^j\delta$ then 
\[\mathbb{E}\left[|\scal{X}{f_{\lambda}}|^2\right]= O(\lambda^{-n-2j}),\quad \lambda \downarrow 0.\]
In particular if $\mu$ is the Gaussian white noise the asymptotic behavior of the covariance is $\lambda^{-n}$. This suggests the following definition.  

\begin{definition}\label{definition:ultraviolet1}
We say that the Gaussian measure $\mu$ has an ultraviolet behavior no worse than the Gaussian white noise if for any Schwartz function $f \in \mathcal{S}(\R^n)$ we have 
\begin{equation}\label{eq:ultraviolet}
\mathbb{E}\left[|\scal{X}{f_{\lambda}}|^2\right]=O(\lambda^{-n})
\end{equation}
as $\lambda \downarrow 0$.
\end{definition}

\begin{remark}
In the case where the measure $\mu$ is Gaussian, it is translation invariant, and the limit \eqref{eq:ultraviolet} is uniform when the function $f$ belongs to a suitable ball with respect to some seminorm of the space of Schwartz functions $\mathcal{S}(\R^n)$, the condition \eqref{eq:ultraviolet} implies that the measure $\mu$ is supported in the space $B^{-\frac{n}{2}-\varepsilon}_{\infty,\infty,\loc}(\R^n)$ (i.e. the Besov space of distributions with (local) infinite integrability and $-\frac{n}{2}-\varepsilon$ regularity), for any $\varepsilon>0$. More generally, by Lemma 10.2 of \cite{Hairer2014RegularityStructures}, if 
\begin{equation}\label{eq:ultraviolet2}
\mathbb{E}\left[|\scal{X}{f_{\lambda}}|^2\right]=O(\lambda^{-\alpha}),
\end{equation}
for $\alpha \geq  0$, then $\mu$ is concentrated in $B^{-\frac{\alpha}{2}-\varepsilon}_{\infty,\infty,\loc}(\R^n)$, for any $\varepsilon>0$.
\end{remark}

\noindent We say that the measure $\sigma$ in representation \eqref{eq:Gsigmarepresentation} (or the function $h$, or the covariance $\bar\Sigma$) is tame if condition \eqref{eq:integrabilitysigma} holds with $N_0=1$ (this name is taken from Appendix A in \cite{NeebOlafsson2015}, see this reference for other details).

\begin{proposition}\label{proposition:rprepresentation1}
Suppose that $\mu$ satisfies the hypotheses of Theorem \ref{theorem:rprepresentation1} then the following three conditions are equivalent:
\begin{itemize}
    \item the measure $\mu$ has ultraviolet behavior no worse than the Gaussian white noise;
    \item the function $h(q)$ in formula \eqref{eq:rprepresentation1} is bounded at infinity;
    \item the measure $\sigma$ in formula \eqref{eq:formulahNq0} is tame and $P^{1,q_0}$ is a constant.  
\end{itemize}
\end{proposition}
\begin{proof}
We prove that the last two conditions are equivalent; we then prove the equivalence of the first two. \\
Suppose that $h$ is bounded at infinity, and let $N_0$ be the minimum integer such that condition \eqref{eq:integrabilitysigma} holds. Suppose that $N_0>1$, then we have 
\begin{equation}\label{eq:FN0limit}
\lim_{q \rightarrow +\infty}\frac{|F^{N_0,q_0}|}{q^{N_0-2}}=\infty,\quad \lim_{q \rightarrow +\infty}\frac{|F^{N_0,q_0}|}{q^{N_0-1}}=0.
\end{equation}
The second limit follows from the fact that $\frac{1}{(q+s)(q_0+s)^{N_0-1}}\leq   \frac{1}{(1+s)(q_0+s)^{N_0-1}}\in L^1(\sigma(\mathd s))$, for $q\geq  1$, and that $\frac{1}{(q+s)(q_0+s)^{N_0-1}} \rightarrow 0$ as $q \rightarrow +\infty$. For the first limit, if $q\geq  2 q_0$, we have 
\begin{multline}
\frac{|F^{N_0,q_0}|}{q^{N_0-2}}\geq  \frac{1}{2^{N_0-1}}\int_0^{+\infty}\frac{q^{N_0-1}}{q^{N_0-2}(q+s)(q_0+ s)^{N_0-1}}\sigma(\mathd s)\\
=\frac{1}{2^{N_0-1}}\int_0^{+\infty}\frac{q}{(q+s)(q_0+ s)^{N_0-1}}\sigma(\mathd s)\rightarrow \frac{1}{2^{N_0-1}}\int_0^{+\infty}\frac{1}{(q_0+ s)^{N_0-1}}\sigma(\mathd s)=+\infty,
\end{multline}
as $q \rightarrow +\infty$. Since $h$ is bounded at infinity, by the limits \eqref{eq:FN0limit}, the polynomial $P^{N_0,q_0}$ is at most of degree $N_0-1$. But if $P^{N_0,q_0}$ is of degree exactly $N_0-1$, then $h(q) \sim P^{N_0,q_0} \rightarrow +\infty$. If $P^{N_0,q_0}$ is of degree at most $N_0-2$ then $|P^{N_0,q_0}(q)|=O(q^{N_0-2})$ and, by the first limit in \eqref{eq:FN0limit}, $|h(q)|\geq |F^{N_0,q_0}(q)|-|P^{N_0,q_0}(q)| \rightarrow +\infty$. This implies that $N_0=1$ and, since $F^{1,q_0}(q)\rightarrow 0$ as $q\rightarrow +\infty$ by the dominated convergence theorem, that $P^{N_0,q_0}$ is a constant. \\
Conversely, if $\sigma$ is tame and $P^{1,q_0}$ is a constant, then $h$ is bounded at infinity since $F^{1,q_0}(q) \rightarrow 0$, as $q\rightarrow +\infty$, by the dominated convergence theorem. \\

\noindent Suppose that $h$ is bounded at infinity, and consider $f\in \mathcal{S}(\R^n)$. We have 
\[\mathbb{E}\left[|\scal{X}{f_{\lambda}}|^2\right]=\scal{\bar\Sigma*f_{\lambda}}{f_{\lambda}}=(2\pi)^{-n}\Big(c_0 |\hat{f}(0)|^2+\lambda^{-n}\int_{\R^n} |\hat{f}(k)|^2 h\left(\frac{|k|^2}{\lambda^2}\right)\mathd k\Big).\]
By the first part of the proof, since $h$ is bounded at infinity, $\sigma$ is tame and $P^{1,q_0}$ is a constant; thus $h$ is decreasing. This means that when $\lambda <1$, we have $h(|k|^2/\lambda^2)\leq   h(|k|^2)$. Thus, by the dominated convergence theorem $\lambda^{n}\mathbb{E}\left[|\scal{X}{f_{\lambda}}|^2\right]$ is bounded when $\lambda \rightarrow 0$.\\
Conversely, suppose that $\mu$ has  not worse ultraviolet behavior than the Gaussian white noise, and consider a Schwartz function $f$ whose Fourier transform is supported in an annulus $a<|k|<b$. If $\sigma$ is not tame by the previous reasoning we have $h \rightarrow +\infty$ as $q \rightarrow +\infty$. This means that 
\[\lambda^{n}\mathbb{E}\left[|\scal{X}{f_{\lambda}}|^2\right]= (2\pi)^{-n}\int_{a\leq  |\cdot|\leq   b}|\hat f(k)|^2 h\left(\frac{|k|^2}{\lambda^2}\right)\mathd k\geq  \left|\inf_{q\geq  \frac{a^2}{\lambda^2}} h(q)\right| \|f\|_{L^2(\R^n)}^2 \longrightarrow +\infty \]
as $\lambda \rightarrow 0$. Now suppose that $\sigma$ is tame, but $P^{1,q_0}$ is non constant, then $\lim_{\lambda \rightarrow 0}\inf_{q\geq  |a/\lambda|^2} P^{1,q_0}(q)=+\infty$. Thus, by following the same reasoning as above 
\[\lim_{\lambda \rightarrow 0}\lambda^{n}\mathbb{E}\left[|\scal{X}{f_{\lambda}}|^2\right]\geq  \lim_{\lambda \rightarrow 0} \left(\inf_{q\geq  |a/\lambda|^2} P^{1,q_0}(q) \right) \|f\|_{L^2(\R^n)}^2 =+\infty.\]
This implies that $\sigma$ is tame and $P^{1,q_0}$ is a constant, and thus $h$ is bounded at infinity.    
\end{proof}

\begin{remark}\label{remark:betterwhitenoise}
Instead of asking that $\mu$ has ultraviolet behavior no worse than the Gaussian white noise, one could ask that $\mu$ has a better ultraviolet behavior than white noise, namely that for any $f \in \mathcal{S}(\R^n)$ we have 
\[\mathbb{E}\left[|\scal{X}{f_{\lambda}}|^2\right]=o(\lambda^{-n}).\]
Using the same techniques of the proof of Proposition \ref{proposition:rprepresentation1}, one can prove that this condition is equivalent to requiring that $h(q) \rightarrow 0$, as $q \rightarrow +\infty$, or to the condition that $\sigma$ is tame and that $P^{1,q_0}=0$.
\end{remark}

\noindent We introduce another interesting property for quantum field theory.

\begin{definition}\label{definition:mixing}
We say that a Gaussian measure with $0$ mean is mixing if for any compactly supported smooth functions $f,g\in C^{\infty}_c(\R^n)$ we have
\[\lim_{|a| \rightarrow +\infty}\mathbb{E}\left[\scal{X}{f}\scal{X}{\tau_a(g)}\right]=0\]
where $\tau_a(g)(\xi)=g(\xi-a)$.
\end{definition}

\begin{remark}
Definition \ref{definition:mixing} corresponds (in the case of Gaussian measures which are rotation invariant) to the clustering property (see, e.g., Section II.3 of \cite{Simon1974} or Section 6.1 of \cite{GlimmJaffe1987}) in Osterwalder--Schrader axioms, and corresponds to uniqueness of the vacuum of the Minkowskian quantum field.
\end{remark}

\begin{proposition}\label{proposition:rprepresentation2}
Suppose that the measure $\mu$ satisfies the hypotheses of Theorem \ref{theorem:rprepresentation1}, then $\mu$ is mixing if and only if $c_0=0$ in \eqref{eq:rprepresentation1}.
\end{proposition}
\begin{proof}
By \eqref{eq:rprepresentation1}, when tested against functions supported on $\R^n\backslash \{0\}$, $\bar\Sigma$ is a function and $\bar\Sigma(\xi)=(2\pi)^{-n}c_0+G^{\sigma}(\xi)$. For $a \in \R^n$ when $|a|$ is big enough the support of $f$ and $\tau_a(g)$ is disjoint and thus 
\[\mathbb{E}\left[\scal{X}{f}\scal{X}{\tau_a(g)}\right]=(2\pi)^{-n}c_0 \left(\int_{\R^n} f(\xi)\mathd \xi \right)\left(\int_{\R^n} g(\xi')\mathd \xi' \right)+\int_{\R^{2n}}f(\xi)G^{\sigma}(\xi-\xi'-a)g(\xi')\mathd \xi\mathd \xi'.\]
Since, by formula \eqref{eq:Gsigma2}, $G^{\sigma}(\xi-\xi'-a) \rightarrow 0$ as $|a|\rightarrow +\infty$ and $\xi,\xi'$ are in some compact sets, by the dominated convergence theorem, 
\[\mathbb{E}\left[\scal{X}{f}\scal{X}{\tau_a(g)}\right]\rightarrow (2\pi)^{-n}c_0 \left(\int_{\R^n} f(\xi)\mathd \xi \right)\left(\int_{\R^n} g(\xi')\mathd \xi' \right)\]
which gives the thesis.
\end{proof}

\begin{remark}
Although the conditions in Definition \ref{definition:ultraviolet1} and Definition \ref{definition:mixing} are not strictly necessary for the Osterwalder–Schrader reconstruction theorem (at least if one does not require the uniqueness of vacuum state), they are very often assumed in order to obtain better properties of the quantum field.\\
In the version of \cite{GlimmJaffe1987} (Chapter 6, Section 6.1, axiom OS1) of Osterwalder--Schrader axioms, a stronger version of regularity is required, namely, there is a constant $C>0$ and a $p \in[1,2]$, such that, for any $f\in \mathcal{S}(\R^n)$, one has
\[\left|\mathbb{E}\left[e^{\scal{X}{f}}\right]\right| \leq   e^{C(\|f\|_{L^1}+\|f\|_{L^{p}}^p)}.\footnote{In the case $p=2$ stronger property is required: the two point distribution $\bar{\Sigma}$ is a function of the difference variables, it belongs to the space of locally integrable
functions $L^1_{loc}(\R^n)$. In particular, it has only integrable singularities at coinciding
points.}\]
Using a scaling argument in the case of Gaussian measures considered in Theorem \ref{theorem:rprepresentation1}, this condition turns out to be equivalent to requiring  $c_0=0$ and $\sup_{q\in[0,+\infty)} |h(q)|<+\infty$ (which is actually stronger than the property in Definition \ref{definition:ultraviolet1}, since by Proposition \ref{proposition:rprepresentation1} $h$ could be unbounded around $0$).
\end{remark}

\noindent We are finally ready to prove the main result of the section connecting reflection-positive Gaussian measures with generators of subordinated Brownian motion.

\begin{theorem}\label{theorem:rprepresentation2}
Let $\mu$ be a nontrivial Gaussian measure on $\mathcal{S}'(\R^n)$. Then $\mu$ is invariant with respect to the Euclidean group, reflection positive, with an ultraviolet behavior no worse than the Gaussian white noise and mixing if and only if the covariance kernel $\bar\Sigma$ is the Green function of the operator $H^{\beta}+m^2$, where $m \geq  0$ and $H^{\beta}$ is the (opposite) generator of a subordinated Brownian motion on $\R^n$ with the following properties: its symbol is of the form $\widehat{H^{\beta}}(k)=\psi(|k|^2)$, where $\psi$ is a complete Bernstein function, and $1/(m^2+\psi(|k|^2))$ is integrable in a neighborhood of the origin.
\end{theorem}
\begin{proof}
The proof is essentially based on the properties of Stieltjes functions and complete
Bernstein functions (see, e.g., \cite{SchillingSongVondracek2012}). Suppose first that $\mu$ is invariant with respect to the Euclidean group,
reflection positive, has an ultraviolet behavior no worse than the Gaussian white
noise and is mixing. By Theorem \ref{theorem:rprepresentation1}, Proposition \ref{proposition:rprepresentation2} and Proposition \ref{proposition:rprepresentation1} we have
$\widehat{\bar\Sigma}(\mathd k)
=h(|k|^2)\mathd k$, where
\begin{equation}\label{eq:hStieltjesGreen}
h(q)
=
c+\int_0^{+\infty}\frac{1}{q+s}\sigma(\mathd s),
\qquad q>0.
\end{equation}
In particular $h$ is a non-zero
Stieltjes function (see, e.g., Definition 2.1 of \cite{SchillingSongVondracek2012}). By the one-to-one correspondence between Stieltjes functions and complete Bernstein functions (the two classes being reciprocals of one another, see, e.g., Theorem 7.3 of \cite{SchillingSongVondracek2012}), $1/h$ is a complete
Bernstein function. Since every complete Bernstein function is increasing, the
limit $m^2:=\lim_{q\rightarrow 0^+}1/h(q)\in[0,+\infty)$
exists. Define
\begin{equation}\label{eq:Phihq}
\psi(q):=\frac{1}{h(q)}-m^2,\qquad q>0,
\end{equation}
which is again a complete Bernstein function with $\psi(0+)=0$.
Consequently, $\psi$ is the Laplace exponent of a complete Lévy subordinator (see, e.g., Theorem 5.2 of \cite{SchillingSongVondracek2012}), and the operator $H^\beta:=\psi(-\Delta)$, namely the Fourier multiplier with symbol $\psi(|k|^2)$, is the opposite of the generator of the corresponding subordinated Brownian motion on $\R^n$ (see, e.g.,  Proposition 1.3.27 and Theorem 3.3.3 of \cite{Applebaum2009}). This is equivalent to saying that the symbol $\widehat{H^{\beta}}(k)=\psi(|k|^2)$, and thus the symbol of $(H^{\beta}+m^2)^{-1}$ is
\[(\widehat{H^{\beta}}(k)+m^2)^{-1}=(\psi(|k|^2)+m^2)^{-1}=h(|k|^2),\]
which is the density of $\widehat{\bar\Sigma}$. This means precisely that $\bar\Sigma$ is the Green function of
$H^\beta+m^2$. By the relation \eqref{eq:Phihq}, we obtain that the function $1/(m^2+\psi(|k|^2))=h(|k|^2)$ is integrable in a neighborhood of the origin since $h(|k|^2)\mathd k$ is a tempered measure.\\

\noindent Conversely, suppose that $\bar\Sigma$ is the Green function of
$H^\beta+m^2$, where $m\geq 0$ and $H^\beta$ is the opposite generator of
a complete subordinated Brownian motion. Let
$\psi$ denote the corresponding Laplace exponent (which by definition of complete subordinator is a complete Bernstein function), so that
\[
\widehat{H^\beta}(k)=\psi(|k|^2).
\]
Hence $m^2+\psi$ is a non-zero complete Bernstein function and therefore
\[
h(q):=\frac{1}{m^2+\psi(q)}
\]
is a Stieltjes function. Thus there exists a positive measure $\sigma$ on
$\R_+$ and a constant $c \geq  0$ such that
\[
h(q)=c+\int_0^{+\infty}\frac{1}{q+s}\sigma(\mathd s),
\qquad
\int_0^{+\infty}\frac{1}{1+s}\sigma(\mathd s)<+\infty.
\]
Since $h(|k|^2)$ is integrable in a neighborhood of the origin and $w_{1,n}(s)\lesssim\int_{|k|\leq 1}(|k|^2+s)^{-1}\mathd k$ for $s\in[0,1/2]$, by the Tonelli theorem we also have $\int_0^{1/2}w_{1,n}(s)\sigma(\mathd s)<+\infty$. In particular, $\sigma$ is tame and the polynomial term
$P^{1,q_0}$ is constant. Moreover,
\[
\widehat{\bar\Sigma}(\mathd k)
=
h(|k|^2)\mathd k,
\]
so that $c_0=0$ in \eqref{eq:rprepresentation1}. Theorem
\ref{theorem:rprepresentation1} then implies that $\mu$ is invariant with
respect to the Euclidean group and reflection positive. Therefore Proposition \ref{proposition:rprepresentation1} implies that $\mu$
has an ultraviolet behavior no worse than the Gaussian white noise. Finally,
since $c_0=0$, Proposition \ref{proposition:rprepresentation2} implies that
$\mu$ is mixing. This concludes the proof.
\end{proof}

\begin{remark}
If in Theorem \ref{theorem:rprepresentation2} we assume that $\mu$ has a better ultraviolet behavior than the Gaussian white noise (see Remark \ref{remark:betterwhitenoise}) then the symbol $m^2+\psi(|k|^2)$ of the operator $H^{\beta}+m^2$ must satisfy the additional property $\psi(q) \rightarrow +\infty$ as $q\rightarrow +\infty$. Conversely, if $H^{\beta}+m^2$ has all the properties required in Theorem \ref{theorem:rprepresentation2} and its symbol diverges at infinity then the Gaussian measure $\mu$ associated with its Green function has a better ultraviolet behavior than the Gaussian white noise.
\end{remark}

\begin{remark}
We stated Theorem \ref{theorem:rprepresentation2} in $\R^n$ since it has an independent interest. For $n=1$ the operators $H^{\beta}$ appearing in it are exactly the operators \eqref{eq:Hbeta} satisfying \ref{H0}, apart from the conditions on the behavior of $\bar\beta$ at $0$ and at infinity. Indeed, since in \ref{H0} the Brownian motion has generator $\frac{1}{2}\Delta$, the operator \eqref{eq:Hbeta} has symbol $\widehat{H^{\beta}}(k)=\psi(k^2)$ with
\[
\psi(q)=\beta_0q+\int_0^{+\infty}\big(1-e^{-\frac{tq}{2}}\big)\gamma(t)\,\mathd t,
\]
and, as $\beta_0\geq 0$ and $\gamma$ vary as in \ref{H0}, these are exactly the complete Bernstein functions $\psi$ such that $\psi(0+)=0$ (see, e.g., Chapter 6 of \cite{SchillingSongVondracek2012}).
\end{remark}

\part{Main results}\label{part:main}

\section{Abstract reformulation of the main equation}\label{section:abstractreformulation}
Consider the following SPDE
\begin{equation}
  \left[\partial_t + H^{\beta}+m^2 \right]\varphi(t,\xi) = -V'(\varphi(t,\xi))  + \sqrt{2}\,\dot{W}(t,\xi),\qquad \xi\in\R,\ t\geq 0,
  \label{eq:SPDEmain1psuedo}
\end{equation}
where $H^\beta$ is the operator defined by the assumption \ref{H0}, $m>0$, $V$ and $\dot{W}$ are defined by the following assumptions.

\begin{enumerate}[start=1,label={{(H\arabic*})}] 

  \item\label{H1}\emph{$\dot{W}$ is a space-time white noise defined on some suitable probability space $(\Omega,\mathcal{F},\mathbb{P})$, namely for every $t,s\geq 0$ and $\xi_1,\xi_2\in\R$ we have
  \begin{align*}
       &\mathbb{E}\left[\dot{W}(t,\xi_1)\right]=0,\\
       &\mathbb{E}\left[\dot{W}(t,\xi_1)\dot{W}(s,\xi_2)\right]=\delta(t-s)\delta(\xi_1-\xi_2),
  \end{align*}
 where $\delta$ is the Dirac delta in $0$;}

\item\label{H2}\emph{   $V:\R \rightarrow \R$ is a $C^2$ function and there exist $q\geq1$, $C>0$ and $M_0\in\R$ such that
\begin{align}
&|V'(\xi)|\leq C(1+|\xi|^q),\qquad \xi\in\R,\label{H2:pol}\\
&|V''(\xi)|\leq C(1+|\xi|^{q-1}),\phantom{aaa} \xi\in\R,\label{H2:der}\\
&V''(\xi)\geq M_0,\phantom{aaaaaaaaaa} \xi\in\R.\label{H2:diss}
\end{align}}
\end{enumerate}

The first part of this section is devoted to reformulating \eqref{eq:SPDEmain1psuedo} as an abstract stochastic evolution equation within the framework introduced in Section \ref{sect: Abstract-SPDE}. This reformulation will be fundamental for applying the results established in Section \ref{sect:core}. We also show that the associated dynamics is irreducible.\\
In the second part, we establish one of the main results of the present paper: the equivalence between invariant and infinitesimally invariant measures, and their relation to integration-by-parts formulas. Moreover, we show that uniqueness among reversible measures is sufficient to guarantee ergodicity. We point out that the proofs of the main theorems are given in Part \ref{part:abstract}, in an abstract setting typical of the Da Prato--Zabczyk approach. This choice does not reflect any lack of novelty or relevance of these results, we simply preferred not to interrupt the flow of the paper with sections set in abstract Hilbert and Banach spaces.

\subsection{Well-posedness of the stochastic evolution equation}
In this subsection, we introduce the functional framework in which we shall study \eqref{eq:SPDEmain1psuedo}, following the approach to SPDEs developed by Da Prato and Zabczyk.\\
Let $r>1$, $k>0$, $p\geq 2$, let $\rho_{r,k}$ be the weight defined by \eqref{eq:weight} and let $L^p_{r,k}=L^p_{\rho_{r,k}}(\R)$ be the weighted space defined by \eqref{eq:pes-wei}. The Da Prato-Zabczyk approach requires that the operator $-[H^\beta+m^2]$ generate a strongly continuous semigroup on $L^p_{r,k}$ for every $p\geq 2$. Furthermore, for the purposes of the theory developed in the present paper, it is useful to assume that the semigroup is of negative type. We therefore introduce the following additional assumption on the parameters $r$ and $k$.

\begin{enumerate}[start=3,label={{(H\arabic*})}]
\item\label{H3} \emph{Assume $m>0$, $1<r<1+2s_2$ and $k\in (0,k_0)$ where $k_0:=k_0(r,m)$ is a constant given by Lemma \ref{lemma:semigroupLp} such that for every $p\geq 2$
  \[
m^2-\frac{\omega_{r,k}}{p}\geq m^2-\omega_{r,k}>0,\qquad \omega_{r,k}:=\sup_{\xi\in\R}\left(-\,\frac{H^{\beta}\rho_{r,k}(\xi)}{\rho_{r,k}(\xi)}\right).
  \]}
\end{enumerate}

Assume \ref{H0} and \ref{H3} hold, then by Lemma \ref{lemma:semigroupLp} the operator $-[H^\beta+m^2]$ is closable in  $L^p_{r,k}$, for every $p\geq2$ and its closure $-\mathcal{A}_p$ is a strictly negative operator and it is the infinitesimal generator of a strongly continuous semigroup $\{e^{-t\mathcal{A}_p}\}_{t\geq 0}$, in particular \eqref{eq:Lp-semi} and \eqref{eq:Ap-lower} hold. Before proceeding, we note that for every $k_1,k_2\in(0,k_0]$, $1<r_1,r_2<1+2s_2$ and $p_1,p_2\geq 2$ we have
\begin{align*}
    &e^{-t\mathcal{A}_{p_1}}x=e^{-t\mathcal{A}_{p_2}}x,\qquad x\in L^{p_1}_{r_1,k_1}\cap L^{p_2}_{r_2,k_2},\\
    &\mathcal{A}_{p_1}x=\mathcal{A}_{p_2}x,\qquad x\in {\rm Dom}(\mathcal{A}_{p_1})\cap {\rm Dom}(\mathcal{A}_{p_2}).
\end{align*}
Hence, from now on, we omit the subscript $p$ and we simply use the notation $\mathcal{A}$ and the semigroup $\{e^{-t\mathcal{A}}\}_{t\geq 0}$.\\
\begin{rmk}\label{rmk:H3}
Hypothesis \ref{H3} plays a fundamental role, since it guarantees that at least the linear equation, namely \eqref{eq:SPDEmain1Abstract} with $V=0$, has a unique invariant measure in $L^2_{r,k}$. Indeed, by \eqref{eq:Lp-semi} with $p=2$ and \ref{H3}, the semigroup $\{e^{-t\mathcal A}\}_{t\geq0}$ is exponentially stable in $L^2_{r,k}$, so that $\delta_0$ is the only invariant measure of the deterministic equation $z'=-\mathcal Az$; hence, by Lemma \ref{lemma:stochconv} and \cite[Theorem 9.1.1(ii)]{Dap-Zab1996}, the linear equation has a unique invariant measure, namely the Gaussian measure $\mu_{\free}$ (see Remark \ref{rmk:conc-misura}). The choice of the space is essential, as shown in \cite[Section 9.4]{Dap-Zab1996} for $H^\beta=-\Delta$, namely $\beta_0=1$ and $\bar\beta\equiv0$. In this case, by \cite[Proposition 9.4.6(v)]{Dap-Zab1996}, the linear equation has exactly one invariant measure in $L^2_{r,k}$ for every $r>1$ and $k>0$, that is, in every space with polynomial weight. On the contrary, in the space $L^2_\rho$ with exponential weight $\rho(\xi)=e^{-\kappa\abs{\xi}}$, $\kappa>0$, by \cite[Proposition 9.4.6(iii)-(iv)]{Dap-Zab1996} the invariant measure is unique if $m^2>\kappa^2/2$, while there are infinitely many invariant measures if $m^2<\kappa^2/4$. Indeed, in the latter case the function $\bar x(\xi):=e^{m\xi}$ belongs to $L^2_\rho$ and satisfies $\bar x''=m^2\bar x$, so that $\delta_{\bar x}\neq\delta_0$ is an invariant measure of the deterministic equation $z'=\Delta z-m^2z$, and uniqueness fails by \cite[Theorem 9.1.1(ii)]{Dap-Zab1996}.
\end{rmk}

Equation \eqref{eq:SPDEmain1psuedo} can be rewritten as an abstract evolution equation in $L^2_{r,k}$ of the form
\begin{equation}
  \begin{cases}
      \mathd X(t) + \mathcal{A}X(t)\mathd t = \mathcal{B}(X(t))\mathd t+ \sqrt{2\rho_{r,k}}\mathd W_{r,k}(t),\qquad t> 0,\\
      X(0)=x,
  \end{cases}
  \label{eq:SPDEmain1Abstract}
\end{equation}
where
\begin{enumerate}
\item $\{W(t)\}_{t\geq 0}$ is an $L^2$-cylindrical Wiener process defined by
\[
W(t):=\sum_{j=1}^{+\infty}\beta_j(t)e_j,\qquad t\geq 0,
\]
$\{e_j\}_{j\in\N}$ is an orthonormal basis of $L^2$ and $\beta_1,\beta_2,\ldots$ are real independent Brownian motions defined on $(\Omega,\mathcal{F},\mathbb{P})$; and  $\{W_{r,k}(t)\}_{t\geq 0}$ is an $L_{r,k}^2(\R)$-cylindrical Wiener process defined by
\[
W_{r,k}(t):=\frac{W(t)}{\sqrt{\rho_{r,k}}}=\sum_{j=1}^{+\infty}\beta_j(t)\frac{e_j}{\sqrt{\rho_{r,k}}},\qquad t\geq 0,
\]
\item $\mathcal{B}:{\rm Dom}(\mathcal{B})\subseteq L^2_{r,k}\rightarrow L^2_{r,k}$ is the Nemytskii operator defined by
\[
\mathcal{B}(x)(\xi)=-V^{'}(x(\xi)),\qquad \xi\in\R.
\]
\end{enumerate}

In this section, we want to prove the well-posedness of \eqref{eq:SPDEmain1Abstract} in the mild sense using Theorem \ref{th:solmild}. Hence, first of all, we need to prove the regularity of the stochastic convolution associated with \eqref{eq:SPDEmain1Abstract},
\begin{equation}\label{eq:stochconv}
W_{\mathcal A}(t):=\sqrt{2}\int_0^t e^{-(t-s)\mathcal A}\,\mathd W(s),\qquad t\geq 0.
\end{equation}
Assume \ref{H3} holds, for $p\geq 2$ and $\eps\geq 0$, we denote by $W^{p,\eps}_{r,k}$ the weighted fractional Sobolev space
\[
W^{p,\eps}_{r,k}:=\left\{u\in L^p_{r,k}\ :\ \rho_{r,k}^{1/p}u\in W^{\eps,p}(\R)\right\},\qquad \norm{u}_{W^{p,\eps}_{r,k}}:=\norm{\rho_{r,k}^{1/p}u}_{W^{\eps,p}(\R)},
\]
where $W^{\eps,p}(\R)$ is the usual (unweighted) Sobolev space; in particular $W^{p,0}_{r,k}=L^p_{r,k}$. Since $\rho_{r,k}^{1/p}(\xi)=(1+k|\xi|^2)^{-r/(2p)}$ is smooth and satisfies
\begin{equation}\label{eq:weightclass}
\abs{D^\gamma \rho_{r,k}^{1/p}(\xi)}\leq c_\gamma\,\rho_{r,k}^{1/p}(\xi),\qquad
\rho_{r,k}^{1/p}(\xi)\leq c\,\rho_{r,k}^{1/p}(\xi')\big(1+\abs{\xi-\xi'}^2\big)^{\frac{r}{2p}},\qquad \xi,\xi'\in\R,
\end{equation}
the weight $\rho_{r,k}^{1/p}$ belongs to the class of admissible weights of \cite[Definition 6.1]{TriebelbookIII}; hence, by \cite[Theorem 6.5]{TriebelbookIII}, $W^{p,\eps}_{r,k}$ is a separable Banach space.

\begin{lemm}\label{lemma:stochconv}
Assume that conditions \ref{H0}, \ref{H1} and \ref{H3} hold. For every $p\geq 2$, $T>0$ and $\eps\in\big[0,s_1-\tfrac12\big)$ the stochastic convolution $\{W_{\mathcal A}(t)\}_{t\geq 0}$ defined in \eqref{eq:stochconv} admits a version with continuous trajectories in $W^{p,\eps}_{r,k}$ and, for every $n\in\N$, we have
\begin{equation}\label{eq:stochconvbound}
\mathbb{E}\left[\sup_{t\in[0,T]}\norm{W_{\mathcal A}(t)}_{W^{p,\eps}_{r,k}}^{n}\right]<+\infty.
\end{equation}
In particular, choosing $\eps\in\big(0,s_1-\tfrac12\big)$ and $p>1/\eps$, the trajectories of $\{W_{\mathcal A}(t)\}_{t\geq 0}$ take values in $C(\R)$ $\mathbb{P}$-a.s.
\end{lemm}

\begin{proof}
For simplicity, we set $\rho:=\rho_{r,k}$. Moreover, we recall that $\{e^{-t(\mathcal{A}-m^2)}\}_{t\geq 0}$ on $L^2$ corresponds to the strongly continuous convolution semigroup associated with a symmetric L\'evy process, whose transition kernel $g:[0,+\infty)\times\R\rightarrow [0,+\infty)$ is translation invariant and satisfies
\begin{equation}\label{eq:kernelL2}
\norm{g(\sigma,\cdot)}_{L^2}^2
=\frac{1}{2\pi}\int_\R e^{-2\sigma\,\widehat{H^\beta}(\xi)}\,\mathd\xi,\qquad \sigma>0,
\end{equation}
where $\widehat{H^\beta}$ denotes the Fourier symbol of $H^{\beta}$.

The continuity of the trajectories and the bound \eqref{eq:stochconvbound} follow from the factorization method \cite[Section 5.3]{Dap-Zab2014}. Fix $\theta\in\big(\eps,s_1-\tfrac12\big)$ and $\alpha\in\big(0,\tfrac12-\tfrac{1+2\theta}{4s_1}\big)$, which are nonempty ranges since $s_1>\tfrac12$; note that in particular $\alpha<\tfrac12-\tfrac{1}{4s_1}$. Consider the process
\[
Y(s):=\sqrt2\int_0^s (s-\sigma)^{-\alpha}e^{-(s-\sigma)\mathcal A}\,\mathd W(\sigma),\qquad s\in[0,T].
\]
By the standard factorization formula, for every $t>0$ we have
\begin{equation}\label{eq:factorization}
W_{\mathcal A}(t)=\frac{\sin(\pi\alpha)}{\pi}\int_0^t (t-s)^{\alpha-1}e^{-(t-s)\mathcal A}\,Y(s)\,\mathd s.
\end{equation}
For every $s\in[0,T]$ we set
\begin{equation}\label{eq:pointwisevariance}
v_\alpha(s):=2\int_0^s \tau^{-2\alpha}e^{-2m^2\tau}\norm{g(\tau,\cdot)}_{L^2}^2\,\mathd\tau
=\frac1\pi\int_0^s \frac{e^{-2m^2\tau}}{\tau^{2\alpha}}\int_\R e^{-2\tau\widehat{H^\beta}(\xi)}\mathd\xi\,\mathd\tau .
\end{equation}
By \eqref{eq:Hbeta} the Fourier symbol of $H^{\beta}$ is
\[
\widehat{H^\beta}(\xi)=\beta_0\xi^2+\int_\R\big(1-\cos(z\xi)\big)\bar\beta(z)\,\mathd z ,
\]
and we claim that there exists $c>0$ such that
\begin{equation}\label{eq:symbol-lower}
\widehat{H^\beta}(\xi)\geq c\,\abs{\xi}^{2s_1}\qquad\text{for }\abs{\xi}\geq1.
\end{equation}
Indeed, if $\beta_0>0$ then $s_1=1$ and \eqref{eq:symbol-lower} holds with $c=\beta_0$. If instead $\beta_0=0$, then for $\abs{\xi}\geq1$ we have $\abs{\xi}^{-1}\leq1$, so that \eqref{eq:levybilateral} and the change of variables $u=z\xi$ give
\[
\widehat{H^\beta}(\xi)\geq\int_{\abs{z}\leq\abs{\xi}^{-1}}\big(1-\cos(z\xi)\big)\bar\beta(z)\,\mathd z
\geq C_\beta^{-1}\int_{\abs{z}\leq\abs{\xi}^{-1}}\frac{1-\cos(z\xi)}{\abs{z}^{1+2s_1}}\,\mathd z
=c\,\abs{\xi}^{2s_1},
\]
with $c:=C_\beta^{-1}\int_{\abs{u}\leq1}(1-\cos u)\abs{u}^{-1-2s_1}\,\mathd u\in(0,+\infty)$, the integral being finite because $1-\cos u\simeq u^2/2$ near $0$ and $s_1<1$.

 Splitting the spatial integral in \eqref{eq:pointwisevariance} over $\{|\xi|\leq1\}$ and $\{|\xi|\geq1\}$, the first part is bounded by $2$, while for the second we use \eqref{eq:symbol-lower} and the change of variables $u=\tau^{1/(2s_1)}\xi$ to get, for $\tau\in(0,T]$,
\[
\int_{\abs{\xi}\geq1} e^{-2\tau\,\widehat{H^\beta}(\xi)}\,\mathd\xi
\leq \int_\R e^{-2c\tau\abs{\xi}^{2s_1}}\,\mathd\xi
= C_{s_1}\,\tau^{-\frac{1}{2s_1}} ,
\]
where $C_{s_1}=c^{-1/(2s_1)}\int_\R e^{-2\abs{u}^{2s_1}}\,\mathd u<+\infty$. Hence, since $2\alpha<1$ and $2\alpha+\tfrac{1}{2s_1}<1$, we get
\[
v_\alpha(s)\leq \frac1\pi\int_0^s\Big(2\tau^{-2\alpha}+C_{s_1}\,\tau^{-2\alpha-\frac{1}{2s_1}}\Big)\mathd\tau<+\infty,
\]
in particular $v_\alpha(s)$ is nondecreasing, so $\sup_{s\in[0,T]}v_\alpha(s)=v_\alpha(T)<+\infty$.

Analogously, for every $s\in[0,T]$ and $z\in\R$ we set
\begin{equation}\label{eq:incrementvariance}
\delta_\alpha(s,z):=\mathbb E\big[\abs{Y(s,\xi+z)-Y(s,\xi)}^2\big]
=\frac{2}{\pi}\int_0^s \frac{e^{-2m^2\tau}}{\tau^{2\alpha}}\!\int_\R e^{-2\tau\,\widehat{H^\beta}(\xi)}\big(1-\cos(z\xi)\big)\,\mathd\xi\,\mathd\tau,
\end{equation}
which does not depend on $\xi$ by the translation invariance of the kernel $g$. Using $1-\cos(a)\leq C_\theta\abs{a}^{2\theta}$, splitting again the spatial integral over $\{|\xi|\leq1\}$ and $\{|\xi|\geq1\}$ and using \eqref{eq:symbol-lower} with the change of variables $u=\tau^{1/(2s_1)}\xi$, we get
\[
\int_\R \abs{\xi}^{2\theta}e^{-2\tau\,\widehat{H^\beta}(\xi)}\,\mathd\xi\leq 2+C_{s_1,\theta}\,\tau^{-\frac{1+2\theta}{2s_1}},
\]
so that, since $2\alpha<1$ and $2\alpha+\tfrac{1+2\theta}{2s_1}<1$ by the choice of $\alpha$,
\begin{equation}\label{eq:incrementbound}
\delta_\alpha(s,z)\leq C\abs{z}^{2\theta}\int_0^s\Big(2\tau^{-2\alpha}+C_{s_1,\theta}\,\tau^{-2\alpha-\frac{1+2\theta}{2s_1}}\Big)\mathd\tau=:C_T\abs{z}^{2\theta},\qquad s\in[0,T].
\end{equation}
Moreover, it is easy to see that $\delta_\alpha(s,z)\leq4v_\alpha(s)$, so that, up to enlarging $C_T$, we have $\delta_\alpha(s,z)\leq C_T\min\{1,\abs{z}^{2\theta}\}$ for every $s\in[0,T]$ and $z\in\R$.\\
Since $\dot W$ is a space-time white noise and $v_\alpha(s)$ is finite for every $s\in[0,T]$, then for every $s\in[0,T]$ and $\xi\in\R$ the real random variable $Y(s,\xi)$ is Gaussian and
\[
\mathbb E\big[\,\abs{Y(s,\xi)}^2\,\big]=v_\alpha(s),
\]
in particular, by standard properties of the Gaussian measures, for every $p\geq 2$ we get
\begin{equation}\label{eq:Lpmoment}
\mathbb E\big[\norm{Y(s)}_{L^p_{r,k}}^p\big]
=\int_\R \mathbb E\big[\abs{Y(s,\xi)}^p\big]\rho(\xi)\,\mathd\xi
=c_p\,v_\alpha(s)^{p/2}\norm{\rho}_{L^1(\R)}<+\infty.
\end{equation}
We now estimate the Gagliardo seminorm of $\rho^{1/p}Y(s)$. Writing, for $\xi,\xi'\in\R$ and $z=\xi-\xi'$,
\[
\rho^{1/p}(\xi)Y(s,\xi)-\rho^{1/p}(\xi')Y(s,\xi')=\rho^{1/p}(\xi)\big(Y(s,\xi)-Y(s,\xi')\big)+\big(\rho^{1/p}(\xi)-\rho^{1/p}(\xi')\big)Y(s,\xi'),
\]
and using \eqref{eq:weightclass}, which gives $\abs{\rho^{1/p}(\xi)-\rho^{1/p}(\xi')}\leq C\min(1,\abs{z})\big(\rho^{1/p}(\xi)+\rho^{1/p}(\xi')\big)$, together with \eqref{eq:incrementbound} and $\theta\leq 1$, we obtain
\[
\mathbb E\big[\abs{\rho^{1/p}(\xi)Y(s,\xi)-\rho^{1/p}(\xi')Y(s,\xi')}^2\big]\leq C_T\big(\rho^{2/p}(\xi)+\rho^{2/p}(\xi')\big)\min\big(1,\abs{z}^{2\theta}\big).
\]
Since $Y(s,\xi)$ is Gaussian, $\mathbb E[\abs{\cdot}^p]=c_p(\mathbb E[\abs{\cdot}^2])^{p/2}$, hence by Tonelli
\[
\mathbb E\left[\int_{\R^2}\frac{\abs{\rho^{1/p}(\xi)Y(s,\xi)-\rho^{1/p}(\xi')Y(s,\xi')}^p}{\abs{\xi-\xi'}^{1+\eps p}}\,\mathd\xi\,\mathd\xi'\right]
\leq C_T\norm{\rho}_{L^1(\R)}\int_\R\frac{\min\big(1,\abs{z}^{p\theta}\big)}{\abs{z}^{1+\eps p}}\,\mathd z<+\infty,
\]
the last integral being finite because $\eps<\theta$ near $z=0$ and $\eps>0$ at infinity (for $\eps=0$ the seminorm term is absent). Together with \eqref{eq:Lpmoment} this gives
\begin{equation}\label{eq:Wpmoment}
\sup_{s\in[0,T]}\mathbb E\big[\norm{Y(s)}_{W^{p,\eps}_{r,k}}^p\big]<+\infty.
\end{equation}
Moreover, for every $p\geq 2$ and $q\geq 1$, there exists $C_{q,p}>0$ such that
\begin{equation}\label{eq:gaussianmoments}
\mathbb E\big[\norm{Y(s)}_{W^{p,\eps}_{r,k}}^q\big]\leq C_{q,p}\Big(\mathbb E\big[\norm{Y(s)}_{W^{p,\eps}_{r,k}}^p\big]\Big)^{q/p},
\end{equation}
so that, by \eqref{eq:Wpmoment}, $\sup_{s\in[0,T]}\mathbb E\big[\norm{Y(s)}_{W^{p,\eps}_{r,k}}^{q}\big]<+\infty$ for every $q\geq1$; in particular $Y\in L^{q}([0,T];W^{p,\eps}_{r,k})$ almost surely for every $q\geq1$.

On the other hand, $\{e^{-t\mathcal A}\}_{t\geq0}$ is a strongly continuous semigroup on $W^{p,\eps}_{r,k}$ with $\norm{e^{-t\mathcal A}}_{\mathcal L(W^{p,\eps}_{r,k})}\leq Ce^{(\omega_{r,k}p^{-1}-m^2) t}$: indeed $e^{-t\mathcal A}$ is a convolution operator, hence it commutes with the finite differences $\Delta_z u:=u(\cdot+z)-u$, so that by Lemma \ref{lemma:semigroupLp}
\[
\norm{\rho^{1/p}\Delta_z e^{-t\mathcal A}u}_{L^p}=\norm{e^{-t\mathcal A}\Delta_z u}_{L^p_{r,k}}\leq Ce^{(\omega_{r,k}p^{-1}-m^2)t}\norm{\Delta_z u}_{L^p_{r,k}},
\]
while the commutator between $\Delta_z$ and the multiplication by $\rho^{1/p}$ is bounded by $C\min(1,\abs{z})\norm{u}_{L^p_{r,k}}$ thanks to \eqref{eq:weightclass}; integrating against $\abs{z}^{-1-\eps p}$ and using $\eps<1$ gives the bound. The strong continuity follows from this bound on $[0,1]$ and from the convergence $e^{-t\mathcal A}u\rightarrow u$ in $W^{p,\eps}_{r,k}$ as $t\rightarrow0$ for $u\in C^\infty_c(\R)$, which is dense in $W^{p,\eps}_{r,k}$. Therefore the factorization operator
\[
(R_\alpha g)(t):=\int_0^t (t-s)^{\alpha-1}e^{-(t-s)\mathcal A}g(s)\,\mathd s
\]
maps $L^{q}([0,T];W^{p,\eps}_{r,k})$ continuously into $C([0,T];W^{p,\eps}_{r,k})$ if $q>\tfrac1\alpha$, see \cite[Proposition 5.9]{Dap-Zab2014}. Hence, by \eqref{eq:factorization} we conclude that $\{W_{\mathcal A}(t)\}_{t\geq 0}$ admits a version with continuous trajectories in $W^{p,\eps}_{r,k}$, and
\[
\mathbb E\Big[\sup_{t\in[0,T]}\norm{W_{\mathcal A}(t)}_{W^{p,\eps}_{r,k}}^{n}\Big]
\leq C\,\mathbb E\Big[\norm{Y}_{L^{q}([0,T];W^{p,\eps}_{r,k})}^{n}\Big]<+\infty
\qquad\text{for every }n\in\N.
\]
\end{proof}
\begin{rmk}\label{rmk:stochconvholder}
Let $\eps\in\big(0,s_1-\tfrac12\big)$ and $p>1/\eps$, and set $\gamma:=\eps-\tfrac1p>0$. Since $\rho_{r,k}^{1/p}=\rho_{r/p,k}$, by the Sobolev embedding $W^{\eps,p}(\R)\hookrightarrow C^{0,\gamma}(\R)$ we have
\[
W^{p,\eps}_{r,k}\hookrightarrow C^{\gamma}_{\rho_{r/p,k}}(\R):=\left\{u\in C(\R):\ \rho_{r/p,k}\,u\in C^{0,\gamma}(\R)\right\},
\qquad
\norm{u}_{C^{\gamma}_{\rho_{r/p,k}}}:=\norm{\rho_{r/p,k}\,u}_{C^{0,\gamma}(\R)},
\]
Consequently, by Lemma \ref{lemma:stochconv}, for every $T>0$ and $n\in\N$,
\begin{equation}\label{eq:stochconvholder}
\mathbb{E}\left[\sup_{t\in[0,T]}\norm{W_{\mathcal A}(t)}^{n}_{C^{\gamma}_{\rho_{r/p,k}}}\right]<+\infty .
\end{equation}
Since, for $0<b\leq a$, $\rho_{a,k}=\rho_{a-b,k}\rho_{b,k}$ and the multiplication by the bounded Lipschitz function $\rho_{a-b,k}$ is bounded on $C^{0,\gamma}(\R)$, we have $\norm{u}_{C^\gamma_{\rho_{a,k}}}\leq C\norm{u}_{C^\gamma_{\rho_{b,k}}}$; hence, choosing $p$ large enough, \eqref{eq:stochconvholder} holds with $\rho_{r/p,k}$ replaced by $\rho_{a,k}$, for every $a>0$ and $\gamma\in(0,s_1-\frac12)$.
Moreover, since $H^{\beta}$ commutes with translations, the law of $\{W_{\mathcal A}(t)\}_{t\geq 0}$ is invariant under translations, so that the left-hand side of \eqref{eq:stochconvholder} is unchanged if $\rho_{r/p,k}$ is replaced by $\rho^\tau_{r/p,k}=\rho_{r/p,k}(\cdot-\tau)$; in particular
\begin{equation}\label{eq:stochconvholderuniform}
\sup_{\tau\in\R}\ \mathbb{E}\left[\sup_{t\in[0,T]}\norm{W_{\mathcal A}(t)}^{n}_{C^{\gamma}_{\rho^\tau_{r/p,k}}}\right]<+\infty .
\end{equation}
The bound \eqref{eq:stochconvholder} is also uniform in time, namely
\[
\sup_{t\geq0}\mathbb{E}\left[\norm{W_{\mathcal A}(t)}^{n}_{C^{\gamma}_{\rho_{r/p,k}}}\right]<+\infty,
\]
since in the proof of Lemma \ref{lemma:stochconv} the variance $v_\alpha(s)$ is bounded uniformly in $s\geq0$ and, by \ref{H3}, the norm of $e^{-t\mathcal A}$ in $\mathcal L(W^{p,\eps}_{r,k})$ decays exponentially.
\end{rmk}

Now we can prove the well-posedness of \eqref{eq:SPDEmain1Abstract}.

\begin{prop}\label{prop:wellposedness}
    Assume \ref{H0}, \ref{H1}, \ref{H2} and \ref{H3} hold. Let $q\geq 1$ be the constant given by \eqref{H2:pol}. For every $x\in L^{2q}_{r,k}$ equation \eqref{eq:SPDEmain1Abstract} has a unique mild solution $\{X(t,x)\}_{t\geq 0}$ having continuous trajectories $\mathbb{P}$-a.s. in $L^{2q}_{r,k}$, namely
    \[
    X(t,x)=e^{-t\mathcal{A}}x+\int_0^te^{-(t-s)\mathcal{A}}\mathcal{B}(X(s,x))\mathd s+\int_0^te^{-(t-s)\mathcal{A}}\sqrt{2\rho_{r,k}}\mathd W_{r,k}(s),\qquad \mathbb{P}\text{-a.s.}
    \]
 Furthermore, for every $x\in L^{2}_{r,k}$ equation \eqref{eq:SPDEmain1Abstract} has a unique generalized mild solution $\{X(t,x)\}_{t\geq 0}$ having continuous trajectories $\mathbb{P}$-a.s. in $L^{2}_{r,k}$, namely for every sequence $\{x_n\}_{n\in\N}\subseteq L^{2q}_{r,k}$ such that $x_n\rightarrow x$ in $L^2_{r,k}$ we have
 \[
 \lim_{n\to+\infty}\mathbb{E}\left[\sup_{t\in [0,T]}\norm{X(t,x)-X(t,x_n)}^2_{L^2_{r,k}}\right]=0,\qquad T\geq 0,
 \]
 where $\{X(t,x_n)\}_{t\geq 0}$ is the unique mild solution to \eqref{eq:SPDEmain1Abstract} with initial datum $x_n$.
\end{prop}

\begin{proof}
The statement follows from Theorem \ref{th:solmild} once we check that the abstract Hypotheses \ref{Hyp:base} are satisfied with the identifications
\[
H=L^2_{r,k},\quad E=L^{2q}_{r,k},\quad R=\sqrt{2\rho_{r,k}}\,\Id,\quad A=-\mathcal A,\quad F(x)=\mathcal{B}(x)+M_0x,\quad G(x)=-M_0x,
\]
where $M_0\in\R$ is the constant given by \eqref{H2:diss}. Hypotheses \ref{Hyp:base} (1)-(4) follow immediately by the definitions of $H$, $E$, $R$ and $G$. Hypotheses \ref{Hyp:base} (5) follows by Lemma \ref{lemma:semigroupLp} and condition \ref{H3}, with $\varepsilon:=m^2-\frac{\omega_{r,k}}{2}>0$. Hypotheses \ref{Hyp:base} (7) follows from Lemma \ref{lemma:stochconv} applied with $p=2q^2$, since by \eqref{H2:pol} we have $\norm{F(x)}^{2q}_{E}\lesssim\norm{\rho_{r,k}}_{L^1(\R)}+\norm{x}^{2q^2}_{L^{2q^2}_{r,k}}$, and $W^{2q^2,\eps}_{r,k}\hookrightarrow L^{2q^2}_{r,k}\hookrightarrow L^{2q}_{r,k}$.\\
By \eqref{H2:diss} we have $V''(\xi)\geq M_0$ for every $\xi\in\R$, hence $\xi\mapsto V'(\xi)-M_0\xi$ is nondecreasing and the Nemytskii operator $F(x)=-\big(V'(x)-M_0x\big)$ is dissipative in $H$; standard results on Nemytskii operators of one-dimensional monotone maps give that $F$ is $m$-dissipative in $H$ and that its part $F_E$ is $m$-dissipative in $E$. Moreover, by the growth condition \eqref{H2:pol}, for every $x\in E$,
\[
\norm{F(x)}_H\leq \norm{V'(x)}_{L^2_{r,k}}+\abs{M_0}\norm{x}_{L^2_{r,k}}\leq C\big(1+\norm{x}_{L^{2q}_{r,k}}^{q}\big)+\abs{M_0}\norm{x}_{L^{2q}_{r,k}}\leq C'\big(1+\norm{x}_E^{q}\big),
\]
hence Hypotheses \ref{Hyp:base} (6) holds. The continuity of $F:E\rightarrow H$ follows from the continuity of $V'$ and \eqref{H2:pol}. Moreover, by \eqref{H2:der} and the H\"older inequality with exponents $\frac{q}{q-1}$ and $q$, for every $x,y,h\in E$ with $\norm{x}_E,\norm{y}_E\leq B$ we have
\[
\abs{\scal{F(x)-F(y)}{h}_H}\leq \int_\R\abs{x-y}\big(C(1+\abs{x}^{q-1}+\abs{y}^{q-1})+\abs{M_0}\big)\abs{h}\rho_{r,k}\mathd\xi\leq L_B\norm{x-y}_H\norm{h}_E,
\]
with $L_B:=C'(1+B^{q-1})$, namely \eqref{FdualEH} holds. Finally, the space $E=L^{2q}_{r,k}$ is reflexive and dense in $H$, so that Hypotheses \ref{Hyp:base}(1) holds.
\end{proof}

By the previous proposition, we can define the Markov transition semigroup associated to \eqref{eq:SPDEmain1Abstract} $\{P(t)\}_{t\geq 0}$ given by
\begin{equation}\label{eq:semiL2rk}
P(t)\Phi(x):=\mathbb{E}\left[\Phi(X(t,x))\right],\qquad \Phi\in B_b(L^2_{r,k}),\, x\in L^2_{r,k},\, t\geq 0,
\end{equation}
where $\{X(t,x)\}_{t\geq 0}$ is the unique generalized mild solution to \eqref{eq:SPDEmain1Abstract} with initial datum $x$.

\subsection{Irreducibility of the dynamics}

The next proposition shows that the dynamics associated with \eqref{eq:SPDEmain1Abstract} is irreducible.

\begin{prop}\label{prop:irreducibility-c}
Assume \ref{H0}, \ref{H1}, \ref{H2}, \ref{H3} hold and there exists $M_1\in\R$ such that $V(\xi)\geq M_1$ for every $\xi\in\R$. Then there exists $c_*\in\R$ such that $m^2c_*+V'(c_*)=0$ and, setting $x_*\equiv c_*$, for every $x\in L^2_{r,k}$ and $\eta>0$ there exists $t_0>0$ such that $\mathbb P\big(\norm{X(t,x)-x_*}_{L^2_{r,k}}<\eta\big)>0$ for every $t\geq t_0$.
\end{prop}
\begin{proof}
In the setting of Proposition \ref{prop:wellposedness}, namely $H=L^2_{r,k}$, $E=L^{2q}_{r,k}$, $A=-\mathcal A$, $F+G=\mathcal B$ and $R=\sqrt{2\rho_{r,k}}\Id$, Hypothesis \ref{hyp:irr} holds. Indeed, since there exists $M_1\in\R$ such that $V(\xi)\geq M_1$ for every $\xi\in\R$, then the continuous function $c\mapsto m^2c+V'(c)$ takes both signs: if it were negative for every $c$ large enough, integrating we would get $V(c)\rightarrow-\infty$ as $c\rightarrow+\infty$, and similarly at $-\infty$. Hence there exists $c_*\in\R$ with $m^2c_*+V'(c_*)=0$. The constant function $x_*\equiv c_*$ belongs to $E$ and to $H^{2s_1}_{\rho_{r,k}}$, since $\rho_{r,k}\in L^1(\R)$ and, by \eqref{eq:weightclass}, $\rho_{r/2,k}\in H^2(\R)$; hence $x_*\in{\rm Dom}(\mathcal A)$ by Remark \ref{rmk:domainA2}, and $\mathcal Ax_*=m^2c_*$ since $H^\beta$ vanishes on constants. Therefore $Ax_*+(F+G)(x_*)=-m^2c_*-V'(c_*)=0$. Let $E_0:=\{c_*+h:\ h\in C^\infty_c(\R)\}$, which is dense in $L^2_{r,k}$. For $x=c_*+h\in E_0$, the function $w(s):=e^{-s\mathcal A}h$ is the convolution of $h$ with the transition kernel of a L\'evy process times $e^{-m^2s}$, so that $\sup_{s\geq0}\big(\norm{w(s)}_{L^2(\R)}+\norm{w(s)}_{L^{2q}(\R)}\big)\leq\norm{h}_{L^2(\R)}+\norm{h}_{L^{2q}(\R)}$. By \eqref{H2:der},
\[
\abs{Ax_*+(F+G)(x_*+w(s))}=\abs{V'(c_*+w(s))-V'(c_*)}\leq C_*\abs{w(s)}\big(1+\abs{w(s)}^{q-1}\big),
\]
so that this function belongs to $L^\infty(0,T;L^2(\R))$ for every $T>0$, and $u(s):=-\big(Ax_*+(F+G)(x_*+w(s))\big)/\sqrt{2\rho_{r,k}}$ satisfies $\norm{u(s)}_{L^2_{r,k}}=2^{-1/2}\norm{Ax_*+(F+G)(x_*+w(s))}_{L^2(\R)}$, whence $u\in L^2(0,T;L^2_{r,k})$. Since Hypotheses \ref{Hyp:base} hold by the proof of Proposition \ref{prop:wellposedness}, the statement follows from Proposition \ref{prop:irreducibility}.
\end{proof}

\begin{rmk}
By Proposition \ref{prop:irreducibility-c}, under \ref{H0}, \ref{H1}, \ref{H2}, \ref{H3} and \ref{H4} the constant function $x_*$ is accessible for $\{P(t)\}_{t\geq0}$ in the sense of \cite[Definition 7.2]{Hairer2008}. However, in Section \ref{section:nonuniqueness} we present examples satisfying these assumptions for which the invariant measure is not unique. Hence, by \cite[Lemma 7.3, Corollary 7.8 and Theorem 11.7]{Hairer2008}, in general $\{P(t)\}_{t\geq0}$ is neither strong Feller nor asymptotically strong Feller. On the other hand, we expect that $\{P(t)\}_{t\geq0}$ is strong Feller along the directions of the noise, namely along $R(L^2_{r,k})=L^2(\R)$, where $R=\sqrt{2\rho_{r,k}}\Id$ is the diffusion coefficient (see the proof of Proposition \ref{prop:irreducibility-c}). Regularity results along the directions of a subspace for transition semigroups associated with SPDEs in Hilbert spaces are proved in \cite{Big-Fer2023PA,Big-Fer2023JDE,Big-Fer2024,CerraiLunardi2019,LunardiRockner2021}.
\end{rmk}

\subsection{The Kolmogorov operator associated with \eqref{eq:SPDEmain1psuedo}}

Since the results derived in the next subsection rely on the abstract theory developed in Section \ref{sect:core}, here we define the Kolmogorov operator associated with \eqref{eq:SPDEmain1Abstract} (or, equivalently, with \eqref{eq:SPDEmain1psuedo}). In this section, we fix $k>0$ and $r>1$ satisfying \ref{H3}. We define the space 
\begin{align}
\mathcal{F}C_b^{\infty}(L^2_{r,k}):=\bigg\{\Phi:L^2_{r,k}\rightarrow \R\ :\ \exists n\in\N,\ &\phi\in C^{\infty}_b(\R^n;\R),\ f_1,\ldots,f_n\in C^{\infty}_c(\R)\ :\\
& \ \Phi(x)=\phi\left(\scal{x}{f_1}_{L^2},\ldots, \scal{x}{f_n}_{L^2}\right),\ \forall x\in L^2_{r,k}\bigg\}\notag.
\end{align}
We note that for every $1\leq r_1\leq r_2$ and $k_1,k_2\geq 0$ and $\Phi\in \mathcal{F}C_b^{\infty}(L^2_{r_2,k_2})$ we have
\begin{equation}\label{eq:inc-cilindriche}
\Phi_{|L_{r_1,k_1}^2}\in \mathcal{F}C_b^{\infty}(L^2_{r_1,k_1}).
\end{equation}
For every $r,k\geq 0$ and $\Phi \in \mathcal{F}C_b^{\infty}(L^2_{r,k})$ we define 
\begin{align}
&\mathcal{N}_0\Phi(x) := \left( \sum_{i, j = 1}^n \partial_{y_i,
  y_j}  \phi\left(\scal{x}{f_1}_{L^2},\ldots, \scal{x}{f_n}_{L^2}\right)
  \scal{f_i}{f_j}_{L^2} \right) \nonumber\\
  &\phantom{aaaaaa}  - \sum_{i = 1}^n \left(\scal{(H^\beta + m^2)f_i}{x}_{L^2}+ \scal{f_i}{V'\circ x}_{L^2}\right) \partial_{y_i}\phi\left(\scal{x}{f_1}_{L^2},\ldots, \scal{x}{f_n}_{L^2}\right),\quad x\in L^{2q}_{r,k}\label{eq:operatorL}\\ 
&\delta_x \Phi (x) := \sum_{i = 1}^n \partial_{y_i} \phi\left(\scal{x}{f_1}_{L^2},\ldots, \scal{x}{f_n}_{L^2}\right)f_i,\qquad x\in L^2_{r,k}, 
\end{align}
where $\Phi(x)=\phi\left(\scal{x}{f_1}_{L^2},\ldots, \scal{x}{f_n}_{L^2}\right)$, $n\in\N$, $\phi\in C^{\infty}_b(\R^n;\R)$ and $ f_1,\ldots,f_n\in C^{\infty}_c(\R)$. We note that the right-hand side of \eqref{eq:operatorL} is well defined for every $x\in L^{2q}_{r,k}$. Indeed, $(H^{\beta}+m^2)f_i$ is not compactly supported, but by Lemma \ref{lemma:testdecay} we have $\abs{(H^{\beta}+m^2)f_i}\leq C_{f_i}\rho_{r,k}$, whence
\[
\abs{\scal{(H^{\beta}+m^2)f_i}{x}_{L^2}}\leq C_{f_i}\int_\R\abs{x}\rho_{r,k}\mathd\xi\leq C_{f_i}\norm{\rho_{r,k}}_{L^1(\R)}^{1/2}\norm{x}_{L^2_{r,k}},
\]
while $\scal{f_i}{V'\circ x}_{L^2}$ is finite by \eqref{H2:pol}. 

\begin{rmk}\label{rmk:idN0}
Let $r,k>0$ and $\Phi \in \mathcal{F}C_b^{\infty}(L^2_{r,k})$. We underline that $\Phi: L^2_{r,k}\rightarrow \R$ is Fréchet differentiable, in particular for every $h\in L^2_{r,k}$ we have
\begin{equation}\label{eq:derloc-derfun}
D\Phi(x)[h]=\scal{\nabla\Phi(x)}{h}_{L^2_{r,k}}=\scal{\delta_x\Phi(x)}{h\rho^{-1}_{r,k}}_{L^2_{r,k}},
\end{equation}
if in addition $h\in L^2$, then
\[
D\Phi(x)[h]=\scal{\delta_x\Phi(x)}{h}_{L^2}.
\]
Taking into account Section \ref{sect:core}, the infinite-dimensional second-order Kolmogorov operator associated to \eqref{eq:SPDEmain1Abstract} is the operator $N_0$ defined by \eqref{Op-OUP} with $H=L^2_{r,k}$, $E=L^{2q}_{r,k}$, $R=\sqrt{2\rho_{r,k}}\Id_{L^2_{r,k}}$, $A=-\mathcal{A}$ and $F+G=\mathcal{B}$, namely
\begin{align}
N_0\Phi(x):={\rm Tr}\left[\rho_{r,k}\nabla^2\Phi(x)\right]-\scal{x}{\mathcal{A}^*\nabla\Phi(x)}_{L^2_{r,k}}+\scal{\mathcal{B}(x)}{\nabla\Phi(x)}_{L^2_{r,k}},\quad x\in L^{2q}_{r,k},
\end{align}
    where $\Phi$ belongs to the space $\xi_{C^\infty_c(\R)}(L^2_{r,k})$ defined by 
\begin{align*}
\xi_{C^\infty_c(\R)}(L^2_{r,k}):={\rm Span}\{\mbox{real and imaginary parts of the functions } x\mapsto e^{i\scal{x}{h}_{L^2}}\,|\, h\in C^\infty_c(\R)\}.
\end{align*}
Noting that $e^{i\scal{x}{h}_{L^2}}=e^{i\scal{x}{h\rho^{-1}_{r,k}}_{L_{r,k}^2}}$ and that $h\rho^{-1}_{r,k}$ still belongs to $C^\infty_c(\R)$ we deduce
\[
\xi_{C^\infty_c(\R)}(L^2_{r,k})\subseteq \mathcal{F}C_b^{\infty}(L^2_{r,k}).
\]
Hence, taking $\Phi\in \xi_{C^\infty_c(\R)}(L^2_{r,k})$ such that $\Phi(x):=\phi\left(\scal{x}{f_1}_{L^2},\ldots, \scal{x}{f_n}_{L^2}\right)$, $n\in\N$, $\phi\in C^{\infty}_b(\R^n;\R)$ and $ f_1,\ldots,f_n\in C^{\infty}_c(\R)$, we have
\[
\Phi(x):=\phi\left(\scal{x}{f_1}_{L^2},\ldots, \scal{x}{f_n}_{L^2}\right)=\phi\left(\scal{x}{f_1\rho_{r,k}^{-1}}_{L_{r,k}^2},\ldots, \scal{x}{f_n\rho_{r,k}^{-1}}_{L_{r,k}^2}\right).
\]
Moreover, $N_0$ is defined in \eqref{Op-OUP} on $\xi_A(H)$ with with $\mathcal D=C^\infty_c(\R)$;; here $A^*=-\mathcal A^*$, where $\mathcal A^*$ denotes the adjoint of $\mathcal A$ in $L^2_{r,k}$, which by the symmetry of $H^{\beta}$ in $L^2$ is given by $\mathcal A^*v=\rho^{-1}_{r,k}(H^{\beta}+m^2)(\rho_{r,k}v)$. Since $\rho_{r,k}v\in C^\infty_c(\R)$ for every $v\in C^\infty_c(\R)$, the second part of \eqref{eq:testdecay-L2} gives $\mathcal A^*v\in L^2_{r,k}$, so that $C^\infty_c(\R)\subseteq{\rm Dom}(\mathcal A^*)$, and, by Lemma \ref{lemma:core-adjoint}, $C^\infty_c(\R)$ is a core for $\mathcal A^*$. Due to these observations, by \eqref{eq:operatorL} and \eqref{eq:derloc-derfun} we obtain
\[
N_0\Phi=\mathcal{N}_0\Phi.
\]
\end{rmk}
In view of the previous remark, Theorem \ref{thm: coreN} guarantees that  $\mathcal{F}C_b^{\infty}(L^2_{r,k})$ is a $\pi$-core for $\{P(t)\}_{t\geq 0}$ in a suitable space of continuous functions, see also Definition \ref{def:pi-gen-N-p}. In the next Theorem we simply rewrite Theorem \ref{thm: coreN} in this setting.

\begin{thm}\label{thm:core-concrete}
Assume \ref{H0}, \ref{H1}, \ref{H2} and \ref{H3} hold. Let $q\geq 1$ be the constant given by \eqref{H2:pol}. For every $\Phi\in C_b(L^{2q}_{r,k})$ such that 
\begin{align*}
    &\lim_{t\rightarrow 0^+}\dfrac{P(t)\Phi(x)-\Phi(x)}{t}=\Psi(x),\quad x\in L^{2q}_{r,k}\qquad {\rm with }\,\, \Psi\in C_b(L^{2q}_{r,k}),\\
    &\sup_{t\in [0,1]}\sup_{x\in L^{2q}_{r,k}}\dfrac{\abs{P(t)\Phi(x)-\Phi(x)}}{t(1+\norm{x}^q_{L^{2q}_{r,k}})}<+\infty,
\end{align*}
there exist $l\in\N$ and $\{\Phi_{n_1,\ldots,n_l}\}_{n_1\in\N,\ldots,n_l\in\N}\subseteq \mathcal{F}C^{\infty}_b(L^2_{r,k})$ such that
\begin{align*}
    &\lim_{n_1,\ldots,n_l\rightarrow +\infty}\Phi_{n_1,\ldots,n_l}(x)=\Phi(x),\quad x\in L^{2q}_{r,k},\\
    &\lim_{n_1,\ldots,n_l\rightarrow +\infty}\mathcal{N}_0\Phi_{n_1,\ldots,n_l}(x)=\Psi(x),\quad x\in L^{2q}_{r,k},\\
    &\sup_{n_1,\ldots,n_l\in\N}\sup_{x\in L^{2q}_{r,k}}\dfrac{\abs{\mathcal{N}_0\Phi_{n_1,\ldots,n_l}(x)}+\abs{\Phi_{n_1,\ldots,n_l}(x)}}{1+\norm{x}^q_{L^{2q}_{r,k}}}<+\infty.
\end{align*}
In this case, we set $\mathcal{N}\Phi=\Psi$.
\end{thm}

\begin{rmk}
    First of all, we note that the normalization $1+\norm{x}_{L^{2q}_{r,k}}^q$ is necessary, since $\mathcal{N}_0$ is a differential operator with an unbounded first-order term. Moreover, let $\mu$ be a probability measure on $L^{2q}_{r,k}$ with finite moment of order $q$. If the semigroup $\{P(t)\}_{t\geq 0}$ admits a unique extension to a strongly continuous semigroup on $L^2(L^{2q}_{r,k},\mu)$, then Theorem \ref{thm:core-concrete} implies that its infinitesimal generator on $L^2(L^{2q}_{r,k},\mu)$ coincides with the closure of $\mathcal{N}$ in $L^2(L^{2q}_{r,k},\mu)$. As in the finite-dimensional setting, the semigroup $\{P(t)\}_{t\geq 0}$ is not strongly continuous on $C_b(L^{2q}_{r,k})$. Indeed, the convergences established in Theorem \ref{thm:core-concrete} are, in general, not uniform, and therefore the notion of an infinitesimal generator is not available in this space. Nevertheless, in Section \ref{sect:core} we introduced a weaker notion of generator, which is sufficient for deriving the results presented in the following sections.
\end{rmk}

\subsection{Relation between invariant, infinitesimally invariant and reversible measures}

Here we state the main results of this section. We emphasize that these results are direct consequences of the abstract theory developed in Section \ref{sect: inv-infinve-abstract}. For this reason, all proofs are presented there in a more general framework, which applies to a broader class of SPDEs than \eqref{eq:SPDEmain1psuedo}. We also highlight that the key ingredient in the proofs is Theorem \ref{thm:core-concrete}.  First of all, we state the basic definitions.

\begin{defi}\label{definition:invariant}
Assume \ref{H0}, \ref{H1}, \ref{H2} and \ref{H3} hold. Let $q\geq 1$ be the constant given by \eqref{H2:pol}.
\begin{enumerate}
    \item We say that a probability measure $\mu$ is (semigroup) invariant for the SPDE \eqref{eq:SPDEmain1psuedo} on $L^2_{r,k}$ if for any $\Phi \in C_b(L^2_{r,k})$ and $t\geq 0$ we have
    \[
    \int_{L^2_{r,k}}P(t)\Phi(x)\mu(\mathd x)=\int_{L^2_{r,k}}\Phi(x)\mu(\mathd x).
    \]
    \item We say that a probability measure $\mu$ is infinitesimally invariant for the SPDE \eqref{eq:SPDEmain1psuedo} on $L^2_{r,k}$ if $\mu(L^{2q}_{r,k})=1$, 
    \[
    \int_{L^{2q}_{r,k}}\norm{x}^q_{L^{2q}_{r,k}}\mu(\mathd x)<+\infty
    \]
    and for every $\Phi\in\mathcal{F}C^{\infty}_b(L^2_{r,k})$ we have
    \[
    \int_{L^{2q}_{r,k}}[\mathcal{N}_0\Phi](x)\mu(\mathd x)=0.
    \]
\end{enumerate}
\end{defi}

The following result is a reformulation for \eqref{eq:SPDEmain1psuedo} of the abstract Theorem \ref{thm:Inv-inf}.

\begin{thm}\label{thm:main-inv-ii}
Assume \ref{H0}, \ref{H1}, \ref{H2} and \ref{H3} hold. Let $q\geq 1$ be the constant given by \eqref{H2:pol}. Let $\mu$ be a Borel probability measure on $L^2_{r,k}$ such that $\mu(L^{2q}_{r,k})=1$ and
    \begin{equation}\label{eq:cond-momenti-inv}
    \int_{L^{2q}_{r,k}}\norm{x}_{L^{2q}_{r,k}}^q\mu(\mathd x)<+\infty.
    \end{equation}
    Then the following conditions are equivalent.
  \begin{enumerate}
    \item $\mu$ is (semigroup) invariant for the SPDE \eqref{eq:SPDEmain1psuedo} on $L^2_{r,k}$.
    
    \item $\mu$ is infinitesimally invariant for the SPDE \eqref{eq:SPDEmain1psuedo} on $L^2_{r,k}$.
  \end{enumerate}
\end{thm}
\begin{rmk}
Thanks to the a priori estimates in Section \ref{section:Fernique}, we can remove condition \eqref{eq:cond-momenti-inv} (see Theorem \ref{theorem:super-main}).
\end{rmk}

Now we state the definitions about reversible measures.

\begin{defi}\label{definition:reversible}
Assume \ref{H0}, \ref{H1}, \ref{H2} and \ref{H3} hold. Let $q\geq 1$ be the constant given by \eqref{H2:pol}.
\begin{enumerate}
    \item We say that the measure $\mu$ is (semigroup) reversible (or equivalently the semigroup $\{P(t)\}$ is symmetric with respect to the measure $\mu$) for the SPDE \eqref{eq:SPDEmain1psuedo}  if for every $\Phi, \Psi \in C_b(L^{2}_{r,k})$ and $t>0$ we have
\[ 
\int_{L^{2}_{r,k}} \Psi(x)P(t)\Phi(x)\mu(\mathd x) = \int_{L^{2}_{r,k}} \Phi(x)P(t)\Psi(x)\mu(\mathd x). 
   \]
    \item  We say that the measure $\mu$ is infinitesimally reversible (or equivalently $\mathcal{N}_0$ is symmetric, or self-adjoint with respect to the measure $\mu$) if $\mu(L^{2q}_{r,k})=1$, 
    \[
    \int_{L^{2q}_{r,k}}\norm{x}^q_{L^{2q}_{r,k}}\mu(\mathd x)<+\infty
    \]
    and for every $\Phi,\Psi\in\mathcal{F}C^{\infty}_b(L^2_{r,k})$ we have
    \[ 
\int_{L^{2q}_{r,k}} \Psi(x)[\mathcal{N}_0\Phi](x)\mu(\mathd x) = \int_{L^{2q}_{r,k}} \Phi(x)[\mathcal{N}_0\Psi](x)\mu(\mathd x) 
   \]
 \item We say that $\mu$ satisfies the integration by parts formula associated with the SPDE \eqref{eq:SPDEmain1psuedo} (with the vector logarithmic derivative $-\mathcal{A}+\mathcal{B}$)
 if $\mu(L^{2q}_{r,k})=1$, 
    \[
    \int_{L^{2q}_{r,k}}\norm{x}^q_{L^{2q}_{r,k}}\mu(\mathd x)<+\infty,
    \]
 and for every $\Phi\in\mathcal{F}C^{\infty}_b(L^2_{r,k})$ and $h\in C^{\infty}_c(\R)$ we have
  \begin{equation}\label{eq:IPF}
  \int_{L^{2q}_{r,k}} \scal{h}{\delta_x \Phi (x)}_{L^2} \mu (\mathd x) = \int_{L^{2q}_{r,k}} 
    \left[\scal{(H^\beta +m^2)h}{x}_{L^2}+\scal{h}{V'\circ x}_{L^2}\right]\Phi (x) \mu (\mathd x) . 
    \end{equation}
\end{enumerate}
The right hand side of \eqref{eq:IPF} is well defined by the same reasons used for \eqref{eq:operatorL}. 
\end{defi}

The following result is a reformulation for the SPDE \eqref{eq:SPDEmain1psuedo} of the abstract Theorems \ref{thm: rev-infrev} and \ref{thm: rev-per-parti}.

\begin{thm}\label{theorem:infinitesimalsymmetri}
Assume \ref{H0}, \ref{H1}, \ref{H2} and \ref{H3} hold. Let $q\geq 1$ be the constant given by \eqref{H2:pol}. Let $\mu$ be a Borel probability measure on $L^2_{r,k}$ such that $\mu(L^{2q}_{r,k})=1$ and
    \[
    \int_{L^{2q}_{r,k}}\norm{x}_{L^{2q}_{r,k}}^{q}\mu(\mathd x)<+\infty.
    \]
    Then conditions 2. and 3. below are equivalent. If, in addition, $\int_{L^{2q}_{r,k}}\norm{x}_{L^{2q}_{r,k}}^{2q}\mu(\mathd x)<+\infty$, then conditions 1., 2. and 3. are equivalent.
\begin{enumerate}
    \item $\mu$ is (semigroup) reversible;
    
    \item $\mu$ is infinitesimally reversible;
    
    \item $\mu$ satisfies the integration by parts formula associated with the SPDE \eqref{eq:SPDEmain1psuedo}.
\end{enumerate}
\end{thm}
\begin{rmk}
We point out that, thanks to the a priori estimates in Section \ref{section:Fernique}, we can remove all the moment conditions in Theorem \ref{theorem:infinitesimalsymmetri}, see Theorem \ref{theorem:super-main}.
\end{rmk}

\subsection{Ergodicity of reversible invariant measure under uniqueness}\label{section:ergodicity}

In this section, we are going to prove that if there exists a unique reversible measure $\mu$ for \eqref{eq:SPDEmain1psuedo}, then $\mu$ is ergodic for the transition semigroup $\{P(t)\}_{t\geq 0}$ defined by \eqref{eq:semiL2rk}, in the following sense.

\begin{defi}\label{definition:ergodic}
Assume that \ref{H0}, \ref{H1}, \ref{H2} and \ref{H3} hold. Let $\mu$ be an invariant measure for \eqref{eq:SPDEmain1psuedo} on $L^2_{r,k}$ in the sense of Definition \ref{definition:invariant}. We say that $\mu$ is \emph{ergodic} for $\{P(t)\}_{t\geq 0}$ if every $A\in\mathcal{B}(L^2_{r,k})$ such that for every $t\geq 0$ 
\begin{equation}\label{eq:Ptinvariantset}
\left[P(t)\mathbb{I}_A\right](x)=\mathbb{I}_A(x),\qquad \mu\text{-a.a. }x\in L^2_{r,k},
\end{equation}
satisfies $\mu(A)\in\{0,1\}$.
\end{defi}

\begin{rmk}\label{rmk:ergodicequivalences}
Since $\mu$ is invariant, by Proposition \ref{prop:extensionL2} $\{P(t)\}_{t\geq 0}$ extends uniquely to a strongly continuous contraction semigroup on $L^2(L^2_{r,k},\mu)$, still denoted by $\{P(t)\}_{t\geq 0}$. By \cite[Theorem 3.2.4]{Dap-Zab1996}, $\mu$ is ergodic in the sense of Definition \ref{definition:ergodic} if and only if for every $\Phi\in L^2(L^2_{r,k},\mu)$ we have
\[
\lim_{T\rightarrow +\infty}\frac{1}{T}\int_0^TP(t)\Phi\,\mathd t=\int_{L^2_{r,k}}\Phi(x)\,\mu(\mathd x)\qquad \text{in }L^2(L^2_{r,k},\mu).
\]
\end{rmk}

\begin{thm}\label{thm:rev-ergo}
  \label{lemma:symmetricinvariantset} Assume that \ref{H0}, \ref{H1}, \ref{H2}, \ref{H3} and \ref{H4} hold. Let $q\geq 1$ be the constant given by \eqref{H2:pol}.  Let $\nu$ be a Borel probability measure on $L^2_{r,k}$ such that $\nu(L^{2q}_{r,k}\cap C(\R))=1$,
\begin{equation}\label{eq:momentergodic}
\int_{L^{2q}_{r,k}}\norm{x}_{L^{2q}_{r,k}}^{2q}\nu(\mathd x)<+\infty,
\end{equation}
and $\nu$ is a reversible measure for  $\{P(t)\}_{t\geq 0}$ in the sense of Definition \ref{definition:reversible}. If $\nu$ is the unique reversible measure for $\{P(t)\}_{t\geq 0}$ satisfying \eqref{eq:momentergodic}, then $\nu$ is ergodic.
\end{thm}

\begin{proof}
Since $\nu$ is reversible, it is invariant, and so, by Proposition \ref{prop:extensionL2}, $\{P(t)\}_{t\geq 0}$ extends to a strongly continuous semigroup of contractions on $L^2(L^2_{r,k},\nu)$ and, for every $t\geq 0$, we have
\begin{equation}\label{eq:reversibleL2}
\int_{L^2_{r,k}} \Psi(x)\left[P(t)\Phi\right](x)\,\nu(\mathd x)=\int_{L^2_{r,k}} \Phi(x)\left[P(t)\Psi\right](x)\,\nu(\mathd x),\qquad \Phi,\Psi\in L^2(L^2_{r,k},\nu).
\end{equation}
Assume by contradiction that $\nu$ is not ergodic, and let $A\in\mathcal{B}(L^2_{r,k})$ be such that $0<\nu(A)<1$ and \eqref{eq:Ptinvariantset} holds. Since $\{P(t)\}_{t\geq 0}$ is Markovian, by \eqref{eq:Ptinvariantset} for every $t\geq 0$ we have
\begin{equation}\label{eq:PtinvariantsetC}
\left[P(t)\mathbb{I}_{A^c}\right](x)=\mathbb{I}_{A^c}(x),\qquad \nu\text{-a.a. }x\in L^2_{r,k}.
\end{equation}
Hence, for every $\Psi\in B_b(L^2_{r,k})$ such that $0\leq \Psi\leq 1$, we have $0\leq P(t)(\mathbb{I}_A\Psi)\leq \mathbb{I}_A$ and $0\leq P(t)(\mathbb{I}_{A^c}\Psi)\leq \mathbb{I}_{A^c}$ $\nu$-a.e., so that
\begin{equation}\label{eq:multiplicative}
\mathbb{I}_A\,P(t)\Psi=\mathbb{I}_A\,P(t)(\mathbb{I}_A\Psi)+\mathbb{I}_A\,P(t)(\mathbb{I}_{A^c}\Psi)=P(t)(\mathbb{I}_A\Psi),\qquad \nu\text{-a.e.},
\end{equation}
and, by linearity and density, \eqref{eq:multiplicative} holds for every $\Psi\in L^2(L^2_{r,k},\nu)$ and every $t\geq 0$.\\
We now consider the Borel probability measure
\[
\nu^A(\mathd x):=\frac{1}{\nu(A)}\mathbb{I}_A(x)\,\nu(\mathd x),
\]
which is well defined since $\nu(A)>0$. Moreover, $\nu^A(L^{2q}_{r,k}\cap C(\R))=1$ and $\nu^A$ satisfies \eqref{eq:momentergodic}. For every $t\geq 0$ and $\Phi,\Psi\in C_b(L^2_{r,k})$, by \eqref{eq:reversibleL2} applied to $\Phi$ and $\mathbb{I}_A\Psi$ and by \eqref{eq:multiplicative}, we get
\begin{align*}
\int_{L^2_{r,k}}\Psi(x)\left[P(t)\Phi\right](x)\,\nu^A(\mathd x)&=\frac{1}{\nu(A)}\int_{L^2_{r,k}}\left[\mathbb{I}_A\Psi\right](x)\left[P(t)\Phi\right](x)\,\nu(\mathd x)\\
&=\frac{1}{\nu(A)}\int_{L^2_{r,k}}\Phi(x)\left[P(t)(\mathbb{I}_A\Psi)\right](x)\,\nu(\mathd x)\\
&=\int_{L^2_{r,k}}\Phi(x)\left[P(t)\Psi\right](x)\,\nu^A(\mathd x),
\end{align*}
namely $\nu^A$ is a reversible measure for $\{P(t)\}_{t\geq 0}$ satisfying \eqref{eq:momentergodic}. Since $\nu^A(A)=1$ and $\nu(A)<1$, we have $\nu^A\neq\nu$, which contradicts the uniqueness assumption. Hence $\nu$ is ergodic.
\end{proof}

\section{Gibbs measures and DLR equations}\label{section:DLRequation}

In this section, we introduce the DLR approach for SPDEs of the form \eqref{eq:SPDEmain1psuedo}.

\subsection{Conditional probability in infinite dimensions}\label{section:conditionalprobability}

In this subsection, we recall some abstract definitions and results that will be needed in the following sections. Throughout this subsection, $E$ denotes a general Banach space of Borel measurable functions $x:\mathbb{R}\to\mathbb{R}$, we denote by $\norm{\cdot}_E$ its norm. Let $\mathcal{I} \subset \R$ be an open interval and $y \in E$, we define the following subspace of $E$
\begin{equation}\label{eq:defEIy}
E^{\mathcal{I},y}:=\{x\in E\ :\ x|_{\mathcal{I}^c}=y|_{\mathcal{I}^c} \}.
\end{equation}
We call $y$, in the definition of $E^{\mathcal{I},y}$, \emph{the boundary condition}. 

Fixed the interval $\mathcal{I}$, it is convenient to consider the Banach space of all possible boundary conditions 
\[
E^{\mathcal{I}^c}:=\{y:\mathcal{I}^c\rightarrow \R\ :\ \exists \hat y\in E\ : \ \hat y|_{\mathcal{I}^c}=y\}.
\]
Hereafter, for any $y \in E^{\mathcal{I}^c}$, and whenever no ambiguity can arise, we use the same symbol $y$ to denote both the original element and any of its extensions. The space $E^{\mathcal{I}^c}$ can be equipped with the natural norm 
\[
\|y\|_{\mathcal{I}^c}:=\inf \{\|\hat{y}\|_E,\hat{y}\text{ is an extension of }y\}.
\]
\begin{rmk}
    We note that if $x\in E$ then $x\in E^{\mathcal{I},\,x|_{\mathcal{I}^c}}$, in particular
    \[
    \bigcup_{y \in E^{\mathcal{I}^c}}E^{\mathcal{I},y}=E.
    \]
\end{rmk}

\begin{rmk}
 If $E=L^p_{r,k}$, then $E^{\mathcal{I}^c}=L^p_{r,k}(\mathcal{I}^c)$ with the norm 
\[
\|y\|^p_{L^p_{r,k}(\mathcal{I}^c)}=\int_{\mathcal{I}^c}|y(\xi)|^p\rho_{r,k}(\xi)\mathd \xi.
\]
In a similar way if $E=C^{\alpha}_{\rho_{r,k}}$ for some $0\leq   \alpha <1$ then $E^{\mathcal{I}^c}=C^{\alpha}_{\rho_{r,k}}(\mathcal{I}^c)$ with equivalent norm 
\[\|y\|_{C^{\alpha}_{\rho_{r,k}}(\mathcal{I}^c)}=\|y \mathbb{I}_{\mathcal{I}^c}\|_{L^{\infty}_{\rho_{r,k}}}+\sup_{\xi\not= \xi',|\xi-\xi'|<1,\xi,\xi'\in \mathcal{I}^c}\rho_{r,k}(\xi)\frac{|y(\xi)-y(\xi')|}{|\xi-\xi'|^{\alpha}}.\]
More generally, if $E$ is a Besov space or a Triebel–Lizorkin space, the definition and the norm of $E^{\mathcal{I}^c}$ coincide with the definition and the norm of the corresponding Besov or Triebel–Lizorkin space on the (regular) domain $\mathcal{I}^c$ (see, e.g., Chapter 4 of \cite{TriebelbookIII}).
\end{rmk}
Due to the definition of the norm $\| \cdot\|_{\mathcal{I}^c}$ the linear map $R_{\mathcal{I}^c}:E\rightarrow E^{\mathcal{I}^c}$ defined as 
\begin{equation}\label{eq:Prestriction}
R_{\mathcal{I}^c}(x)=x|_{\mathcal{I}^c},
\end{equation}
is continuous (in the sense of the topology induced by the norms $\norm{\cdot}_E$ and $\norm{\cdot}_{E^{\mathcal{I}^c}}$, respectively). So, if $\mu$ is a probability measure on $E$, we define 
\begin{equation}\label{eq:legge-Ic}
\mu^{\mathcal{I}^c}:=R_{\mathcal{I}^c*}\mu=\mu\circ R^{-1}_{\mathcal{I}^c}
\end{equation}
namely $\mu^{\mathcal{I}^c}$ is the image measure of $\mu$ under the map $R_{\mathcal{I}^c}$, hence a probability measure on $E^{\mathcal{I}^c}$.

\begin{lemm}
  \label{lemma:disintegration}Assume that $E$ and $E^{\mathcal{I}^c}$
  are separable Banach spaces and let $\mathcal{I}$ be a bounded open interval of $\R$. Then there exists a family of probability measures
  $\{ \mu^{\mathcal{I},y} \}_{y \in E^{\mathcal{I}^c}} \subset \mathcal{P}(E)$ such that
  \begin{align}
  &\mu^{\mathcal{I},y} (E^{\mathcal{I}, y}) = 1,\qquad &\text{for } \mu^{\mathcal{I}^c}\text{-almost all } y \in E^{\mathcal{I}^c},\label{eq:prob-cond}\\
  &\int_E f (x) \mu (\mathd x) = \int_{E^{\mathcal{I}^c}} \left(
     \int_{E^{\mathcal{I}, y}} f (x) \mu^{\mathcal{I},y} (\mathd x) \right)
     \mu^{\mathcal{I}^c} (\mathd y)\qquad &f\in L^1(E;\mu;\R)\label{eq:formula-prob-cond}
     \end{align}
\end{lemm}
\begin{proof}
  The statement can be found in {\cite[Theorem 5.3.1]{AmbrosioGigliSavare}}.
\end{proof}

\begin{rmk}
   For $\mu^{\mathcal{I}^c}$-almost every $y\in E^{\mathcal{I}^c}$, the measure $\mu^{\mathcal{I},y}$ is a regular realization of the conditional probability
\[
\mu(\cdot\ |\ x|_{\mathcal{I}^c}=y|_{\mathcal{I}^c}).
\]
\end{rmk}

\begin{rmk}\label{remark:disintegrationsubset}
If $B\subset E$ is a Borel set with $\mu(B)=1$, then by Lemma \ref{lemma:disintegration} (applied to $f=\mathbb{I}_{E\setminus B}$) we have $\mu^{\mathcal{I},y}(B)=1$ for $\mu^{\mathcal{I}^c}$-almost every $y\in E^{\mathcal{I}^c}$.
\end{rmk}

\subsection{Conditional measure of Gaussian measures}
In this subsection, we focus on the Gaussian case. We denote by $X:\Omega \times \R \rightarrow \R$ the Gaussian random field such that for every $\varphi,\psi \in C^{\infty}_c$ we have
\begin{align*}
    \mathbb{E}\left[\scal{X}{\varphi}_{L^2}\right]&=0,\\
    \mathbb{E}\left[\scal{X}{\varphi}_{L^2}\scal{X}{\psi}_{L^2}\right]&=\scal{\varphi}{\psi}_{\CM}:=\scal{(H^{\beta}+m^2)^{-1}\varphi}{\psi}_{L^2}=\scal{\varphi}{(H^{\beta}+m^2)^{-1}\psi}_{L^2}\\
    &=\int_{\R} \varphi(\xi) (H^{\beta}+m^2)^{-1}(\psi)(\xi)\mathd \xi=\int_{\R^2}\varphi(\xi)G_{\free}(\xi-\xi')\psi(\xi')\mathd \xi \mathd \xi',
\end{align*}
where $G_{\free}$ is the Green function of $H^{\beta}+m^2$. We denote by $\mathcal{H}_{\CM}$ the Hilbert space given by the abstract closure of $C^\infty_c(\R)$ with respect to the norm $\norm{\cdot}^2_{\CM}:=\scal{\cdot}{\cdot}_{\CM}$.  For simplicity if $\varphi \in C^{\infty}_c(\R)$ we denote $X_\varphi:=\scal{X}{\varphi}_{L^2}=\int_{\R}\varphi(\xi) X(\xi) \mathd \xi$.\\

\begin{rmk}\label{rmk:conc-misura}
    Assume that \ref{H0}, \ref{H1} and \ref{H3} hold.  By standard arguments (see for instance \cite{Bogachev98}), the random field $X$ can be seen as a random variable $\tilde X:\Omega\rightarrow L^2_{r,k}$ whose law $\mu_{\free}$ is the centered Gaussian measure on $L^2_{r,k}$ with covariance operator $\mathcal{A}^{-1}M_{\rho_{r,k}}$, where $M_{\rho_{r,k}}$ denotes the multiplication by $\rho_{r,k}$. Taking into account Lemmas \ref{lemma:semigroupLp} and \ref{lemma:stochconv}, by \cite[Theorem 9.1.1(ii)]{Dap-Zab1996} the measure $\mu_{\free}$ is the unique invariant measure of the following SPDE
\begin{equation}
  \begin{cases}
      \mathd Z(t) + \mathcal{A}Z(t)\mathd t = \sqrt{2\rho_{r,k}}\mathd W_{r,k}(t),\qquad t> 0,\\
      Z(0)=x\in L_{r,k}^2.
  \end{cases}
  \label{eq:SPDElinearL2rk}
\end{equation}
    Finally, by  Lemma \ref{lemma:stochconv} we deduce $\mu_{\free}(C(\R))=1$.
\end{rmk}

For every open interval $\mathcal{I}$ of $\R$, following \cite{Rockner1985} we denote by
\begin{equation}\label{eq:HIc}
\sigma(\mathcal{I}^c):=\sigma\{X(\xi)\,|\,\xi \not\in \mathcal{I} \},\qquad \mathcal{H}_{\mathcal{I}^c}=\overline{\{f \in C^{\infty}_c(\R):\supp(f) \subset \mathcal{I}^c\}}^{\mathcal{H}_{\CM}}.
\end{equation}
We denote by $P_{\mathcal{I}^c}$ the orthogonal projection from $\mathcal{H}_{\CM}$ onto $\mathcal{H}_{\mathcal{I}^c}$ with respect to $\langle\cdot,\cdot\rangle_{\CM}$.
\begin{prop}\label{proposition:Rockner}
Assume \ref{H0}, \ref{H1} and \ref{H3} hold. For every $f\in C^\infty_c(\R)$ we have
\[\mathbb{E}\left[X_f|\sigma(\mathcal{I}^c)\right]=X_{P_{\mathcal{I}^c}(f)},\]
where, for $h\in\mathcal H_{\CM}$, $X_h$ denotes the $L^2(\Omega)$-limit of $X_{h_n}$ for any sequence $\{h_n\}_{n\in\N}\subseteq C^\infty_c(\R)$ converging to $h$ in $\mathcal H_{\CM}$.
\end{prop}
\begin{proof}
The statement of the proposition is equivalent to Lemma 6.6 of  \cite{Rockner1985}.
\end{proof}
 
 In the next lemma we provide a representation for $X_{P_{\mathcal{I}^c}(f)}$.

\begin{lemm}\label{lemma:conditional-mean}
Assume \ref{H0}, \ref{H1} and \ref{H3} hold and let $\mathcal I=(a,b)$ be a bounded open interval. Then for every $f\in C^\infty_c(\mathcal I)$ we have
\[
X_{P_{\mathcal{I}^c}(f)}=\scal{f}{\varphi^{R_{\mathcal I^c}X}}_{L^2}\qquad\mathbb P\text{-a.s.},
\]
where $\varphi^{R_{\mathcal I^c}X}$ is the harmonic extension of $R_{\mathcal I^c}X$.
\end{lemm}

\begin{proof}
Let $\mathcal G_{\mathcal I^c}:=\{X_h:\ h\in\mathcal H_{\mathcal I^c}\}$, which is a closed subspace of $L^2(\Omega)$ since $h\mapsto X_h$ is an isometry from $\mathcal H_{\CM}$ into $L^2(\Omega)$. Since $P_{\mathcal I^c}$ is the orthogonal projection onto $\mathcal H_{\mathcal I^c}$, $X_{P_{\mathcal I^c}(f)}$ is the orthogonal projection of $X_f$ onto $\mathcal G_{\mathcal I^c}$. Hence, setting $M_f:=\scal{f}{\varphi^{R_{\mathcal I^c}X}}_{L^2}$, it is enough to prove that $M_f\in\mathcal G_{\mathcal I^c}$ and that $\mathbb E[(X_f-M_f)X_h]=0$ for every $h\in C^\infty_c(\R)$ with $\supp(h)\subset\mathcal I^c$.

We recall that $\mathbb E[X(\xi)X(\xi')]=G_{\free}(\xi-\xi')$, where $G_{\free}$ is the inverse Fourier transform of $(m^2+\widehat{H^\beta})^{-1}$, which belongs to $L^1(\R)$ by \eqref{eq:symbol-lower} and $2s_1>1$; hence $G_{\free}$ is continuous and bounded, and $\xi\mapsto X(\xi)$ is continuous and bounded from $\R$ into $L^2(\Omega)$. By \eqref{eq:harmonic-representation} and the Fubini theorem
\[
M_f=\int_{\mathcal I^c}X(\xi')\,\kappa_f(\mathd\xi'),\qquad \kappa_f(B):=\int_{\mathcal I}f(\xi)P^{\mathcal I}_{\free}(\xi,B)\,\mathd\xi,\quad B\in\mathcal B(\mathcal I^c),
\]
where, by \eqref{eq:total-mass}, $\kappa_f$ is a finite signed measure with total variation at most $\norm{f}_{L^1}$.

If $\xi'$ belongs to the interior of $\mathcal I^c$, then $X(\xi')=\lim_{n}X_{h_n}$ in $L^2(\Omega)$, where $\{h_n\}_{n\in\N}$ is a sequence of mollifiers concentrating at $\xi'$ with $\supp(h_n)\subset\mathcal I^c$, so that $X(\xi')\in\mathcal G_{\mathcal I^c}$; by continuity, $X(a),X(b)\in\mathcal G_{\mathcal I^c}$ as well. Hence the Bochner integral $\int_{\mathcal I^c}X(\xi')\kappa_f(\mathd\xi')$ in $L^2(\Omega)$ belongs to $\mathcal G_{\mathcal I^c}$, and it coincides almost surely with $M_f$.

Let $h\in C^\infty_c(\R)$ with $\supp(h)\subset\mathcal I^c$ and $w_h:=(H^\beta+m^2)^{-1}h=G_{\free}*h$. Then $w_h\in\mathcal F\cap C_b(\R)$, $\mathbb E[X(\xi')X_h]=w_h(\xi')$ for every $\xi'\in\R$, and $\mathcal E(w_h,v)=\scal{h}{v}_{L^2}=0$ for every $v\in\mathcal F$ with $\supp(v)\subset\mathcal I$. Namely $w_h$ is an energy harmonic extension of $w_h|_{\mathcal I^c}$, so that $\varphi^{w_h|_{\mathcal I^c}}=w_h$ by Lemma \ref{lemma:harmonicextensionenergy}. Therefore, by the Fubini theorem and \eqref{eq:harmonic-representation},
\[
\mathbb E[M_fX_h]=\int_{\mathcal I}f(\xi)\int_{\mathcal I^c}w_h(\xi')P^{\mathcal I}_{\free}(\xi,\mathd\xi')\,\mathd\xi=\scal{f}{\varphi^{w_h|_{\mathcal I^c}}}_{L^2}=\scal{f}{w_h}_{L^2}=\mathbb E[X_fX_h],
\]
which concludes the proof.
\end{proof}

\begin{rmk}
    We note that the harmonic extension $\varphi^{R_{\mathcal I^c}X}$ of $R_{\mathcal I^c}X$ is well defined by Theorem \ref{thm:harmonic-continuity} since, by Remark \ref{rmk:conc-misura}, $R_{\mathcal I^c}X\in C(\mathcal I^c)\cap L^2_{r,k}(\mathcal I^c)$ almost surely.
\end{rmk}

Fixed an open bounded interval $\mathcal{I}$ of $\R$, by Lemma \ref{lemma:disintegration} (with $E=L^2_{r,k}$, $E^{\mathcal{I}^c}=L^2_{r,k}(\mathcal{I}^c)$ and $\mu=\mu_{\free}$) there exists $\{\mu^{\mathcal{I},y}_{\free} \}_{y \in L^2_{r,k}(\mathcal{I}^c)} \subset \mathcal{P}(L^2_{r,k})$ such that \eqref{eq:prob-cond} and \eqref{eq:formula-prob-cond} hold. Proposition \ref{proposition:Rockner} allows us to compute explicitly the law of the conditional probability of $\mu_{\free}$.

\begin{lemm}\label{lemma:disintegration-conditional}
Assume \ref{H0}, \ref{H1} and \ref{H3} hold. Let $\mathcal{I}$ be an open bounded interval of $\R$. Then, for  $\mu^{\mathcal{I}^c}_{\free}$-almost every $y\in L^2_{r,k}(\mathcal{I}^c)\cap C(\mathcal{I}^c)$ and for every $f,g\in C^\infty_c(\mathcal I)$, we have
\begin{align}
&\int_{E^{\mathcal{I},y}}\scal{x}{f}_{L^2}\mu^{\mathcal{I},y}_{\free}(\mathd x)=\scal{f}{\varphi^y}_{L^2},\label{eq:meanconditional}\\[4pt]
&\begin{aligned}
\int_{E^{\mathcal{I},y}}\scal{x}{f}_{L^2}\scal{x}{g}_{L^2}\mu^{\mathcal{I},y}_{\free}(\mathd x)
&-\int_{E^{\mathcal{I},y}}\scal{x}{f}_{L^2}\mu^{\mathcal{I},y}_{\free}(\mathd x)\int_{E^{\mathcal{I},y}}\scal{x}{g}_{L^2}\mu^{\mathcal{I},y}_{\free}(\mathd x)\\
&=\scal{Q_{\mathcal{I}^c}(f)}{Q_{\mathcal{I}^c}(g)}_{\CM},
\end{aligned}\label{eq:covariance-identification}
\end{align}
where $\varphi^y$ is the harmonic extension of $y$ given by Theorem \ref{thm:harmonic-continuity}, $Q_{\mathcal{I}^c}:=\Id-P_{\mathcal{I}^c}$, $E^{\mathcal{I},y}$ is given by \eqref{eq:defEIy} (with $E=L^2_{r,k}$) and $\mu_{\free}^{\mathcal{I},y}$ is given by Lemma \ref{lemma:disintegration}. In particular the covariance of $\mu^{\mathcal{I},y}_{\free}$ does not depend on $y$.
\end{lemm}

\begin{proof}
Fix $n\in\N$, $h_1,\dots,h_n\in C^{\infty}_c(\R)$ with $\supp(h_j)\subset\mathcal{I}^c$ for every $j\in\{1,\ldots,n\}$, and $\psi\in C_b(\R^n)$. We set $Y:=\psi(X_{h_1},\dots,X_{h_n})$, we note that $Y$ is bounded and $\sigma(\mathcal{I}^c)$-measurable. By \eqref{eq:formula-prob-cond}, for every $f\in C^\infty_c(\mathcal I)$ we obtain
\[
\mathbb{E}\left[Y X_f\right]
=\int_{L^2_{r,k}(\mathcal{I}^c)}\left(\int_{E^{\mathcal{I},y}}\psi\big(\scal{x}{h_1}_{L^2},\dots,\scal{x}{h_n}_{L^2}\big)\scal{x}{f}_{L^2}\mu^{\mathcal{I},y}_{\free}(\mathd x)\right)\mu^{\mathcal{I}^c}_{\free}(\mathd y).
\]
Since $\supp(h_j)\subset\mathcal{I}^c$ for every $j\in\{1,\ldots,n\}$, by \eqref{eq:prob-cond}, we have $\scal{x}{h_j}_{L^2}=\scal{y}{h_j}_{L^2}$ for $\mu^{\mathcal{I},y}_{\free}$-almost every $x$, hence
\begin{equation}\label{eq:test1}
\mathbb{E}\left[Y X_f\right]=\int_{L^2_{r,k}(\mathcal{I}^c)}\psi\big(\scal{y}{h_1}_{L^2},\dots,\scal{y}{h_n}_{L^2}\big)G_f(y)\mu^{\mathcal{I}^c}_{\free}(\mathd y),
\end{equation}
where $G_f(y):=\int_{E^{\mathcal{I},y}}\scal{x}{f}_{L^2}\mu^{\mathcal{I},y}_{\free}(\mathd x)$.
On the other hand, by Proposition \ref{proposition:Rockner} and Lemma \ref{lemma:conditional-mean},
\[
\mathbb{E}\left[YX_f\right]=\mathbb{E}\left[Y\mathbb{E}\left[X_f\ \big|\ \sigma(\mathcal{I}^c)\right]\right]=\mathbb{E}\left[YX_{P_{\mathcal{I}^c}(f)}\right]=\mathbb{E}\left[Y\scal{f}{\varphi^{R_{\mathcal I^c}X}}_{L^2}\right],
\]
so that, recalling that $\mu^{\mathcal I^c}_{\free}(C(\mathcal I^c))=1$ by Remark \ref{rmk:conc-misura}, we have
\begin{equation}\label{eq:test2}
\mathbb{E}\left[YX_{P_{\mathcal{I}^c}(f)}\right]=\int_{L^2_{r,k}(\mathcal{I}^c)}\psi\big(\scal{y}{h_1}_{L^2},\dots,\scal{y}{h_n}_{L^2}\big)H_f(y)\mu^{\mathcal{I}^c}_{\free}(\mathd y),
\end{equation}
where $H_f(y):=\scal{f}{\varphi^y}_{L^2}$, which is a Borel function of $y$ by Theorem \ref{thm:harmonic-continuity}. Comparing \eqref{eq:test1} and \eqref{eq:test2}, we get
\[
\int_{L^2_{r,k}(\mathcal{I}^c)}\psi\big(\scal{y}{h_1}_{L^2},\dots,\scal{y}{h_n}_{L^2}\big)\big(G_f(y)-H_f(y)\big)\,\mu^{\mathcal{I}^c}_{\free}(\mathd y)=0
\]
for every $n\in\N$, $h_1,\dots,h_n\in C^{\infty}_c(\R)$ with $\supp(h_j)\subset\mathcal{I}^c$ for every $j\in\{1,\ldots,n\}$, and $\psi\in C_b(\R^n)$. Since $L^2_{r,k}(\mathcal{I}^c)$ is separable, its Borel $\sigma$-algebra is generated by the maps $y\mapsto\scal{y}{h}_{L^2}$ with $h\in C^{\infty}_c(\mathcal{I}^c)$; hence by a monotone class argument we deduce $G_f=H_f$ $\mu^{\mathcal{I}^c}_{\free}$-almost everywhere, namely \eqref{eq:meanconditional} holds.\\

In order to prove \eqref{eq:covariance-identification} we first observe that, for every $f\in C^\infty_c(\R)$ we have
\[
X_f=X_{P_{\mathcal{I}^c}(f)}+X_{Q_{\mathcal{I}^c}(f)},\qquad Q_{\mathcal{I}^c}:=\Id-P_{\mathcal{I}^c},
\]
where $X_{P_{\mathcal{I}^c}(f)}$ is $\sigma(\mathcal{I}^c)$-measurable and $X_{Q_{\mathcal{I}^c}(f)}$ is independent of $\sigma(\mathcal{I}^c)$. 
Hence, by Proposition \ref{proposition:Rockner}, for every $f,g\in C^{\infty}_c(\R)$ we deduce
\begin{align}\label{eq:condcovariance}
\mathbb{E}\left[X_fX_g\ |\ \sigma(\mathcal{I}^c)\right]-\mathbb{E}\left[X_f\ |\ \sigma(\mathcal{I}^c)\right]\mathbb{E}\left[X_g\ |\ \sigma(\mathcal{I}^c)\right]&=\mathbb{E}\left[X_{Q_{\mathcal{I}^c}(f)}X_{Q_{\mathcal{I}^c}(g)}\ \big|\ \sigma(\mathcal{I}^c)\right]\notag\\
&=\mathbb{E}\left[X_{Q_{\mathcal{I}^c}(f)}X_{Q_{\mathcal{I}^c}(g)}\right]=\scal{Q_{\mathcal{I}^c}(f)}{Q_{\mathcal{I}^c}(g)}_{\CM}.
\end{align}
Arguing as for \eqref{eq:test1} with $\scal{x}{f}_{L^2}$ replaced by $\scal{x}{f}_{L^2}\scal{x}{g}_{L^2}$, we get
\[
\mathbb{E}\left[YX_fX_g\right]=\int_{L^2_{r,k}(\mathcal{I}^c)}\psi\big(\scal{y}{h_1}_{L^2},\dots,\scal{y}{h_n}_{L^2}\big)G_{f,g}(y)\mu^{\mathcal{I}^c}_{\free}(\mathd y),
\]
where $G_{f,g}(y):=\int_{E^{\mathcal{I},y}}\scal{x}{f}_{L^2}\scal{x}{g}_{L^2}\mu^{\mathcal{I},y}_{\free}(\mathd x)$.
While, by \eqref{eq:condcovariance} and the same computation used for \eqref{eq:test2} applied to $X_{P_{\mathcal{I}^c}(f)}X_{P_{\mathcal{I}^c}(g)}$, we get
\[
\mathbb{E}\left[YX_fX_g\right]=\int_{L^2_{r,k}(\mathcal{I}^c)}\psi\big(\scal{y}{h_1}_{L^2},\dots,\scal{y}{h_n}_{L^2}\big)\Big(H_f(y)H_g(y)+\scal{Q_{\mathcal{I}^c}(f)}{Q_{\mathcal{I}^c}(g)}_{\CM}\Big)\mu^{\mathcal{I}^c}_{\free}(\mathd y).
\]
Also here, by a monotone class argument 
\[
G_{f,g}(y)=H_f(y)H_g(y)+\scal{Q_{\mathcal{I}^c}(f)}{Q_{\mathcal{I}^c}(g)}_{\CM}\qquad\text{for }\mu^{\mathcal{I}^c}_{\free}\text{-almost every }y,
\]
which, combined with \eqref{eq:meanconditional}, gives \eqref{eq:covariance-identification}.
\end{proof}

\begin{prop}\label{proposition:conditionalGaussian}
Assume \ref{H0}, \ref{H1} and \ref{H3} hold. Then for every open bounded interval $\mathcal{I} \subset \R$ and for $\mu_{\free}^{\mathcal{I}^c}$-almost every $y \in L^2_{r,k}(\mathcal{I}^c)\cap C(\mathcal{I}^c)$, for every $f,g \in C^{\infty}_c(\R)$ with $\supp(f),\supp(g)\subseteq\mathcal{I}$  we have
\begin{align}
&\int_{E^{\mathcal{I},y}}\scal{x}{f}_{L^2}\,\mu^{\mathcal{I},y}_{\free}(\mathd x)=\scal{f}{\varphi^y}_{L^2},\label{eq:meanconditionalV}\\[4pt]
&\begin{aligned}
\int_{E^{\mathcal{I},y}}\scal{x}{f}_{L^2}\scal{x}{g}_{L^2}\,\mu^{\mathcal{I},y}_{\free}(\mathd x)
&-\int_{E^{\mathcal{I},y}}\scal{x}{f}_{L^2}\,\mu^{\mathcal{I},y}_{\free}(\mathd x)\int_{E^{\mathcal{I},y}}\scal{x}{g}_{L^2}\,\mu^{\mathcal{I},y}_{\free}(\mathd x)=\int_{\R} f(\xi) u(\xi) \mathd \xi
\end{aligned}\label{eq:covariance-identificationV}
\end{align}
where $E^{\mathcal{I},y}$ is given by \eqref{eq:defEIy} (with $E=L^2_{r,k}$), $\mu_{\free}^{\mathcal{I},y}$ is given by Lemma \ref{lemma:disintegration}, $\varphi^y$ is the harmonic extension of $y$ in the sense of Definition \ref{def:harmonicextension} and $u$  is the unique (weak) solution to the homogeneous Dirichlet problem
\begin{equation}\label{eq:Dirichletg}
\begin{cases}
    (H^{\beta}+m^2)u(\xi)=g(\xi)\qquad &\xi\in\mathcal{I},\\
    u(\xi)=0\qquad &\xi\in\mathcal{I}^c.
\end{cases}
\end{equation}
\end{prop}
\begin{proof}
By \eqref{eq:covariance-identification}, it is sufficient to prove that for every $f,g \in C^{\infty}_c(\R)$ with $\supp(f),\supp(g)\subseteq\mathcal{I}$ we have
\[
\scal{Q_{\mathcal{I}^c}(f)}{Q_{\mathcal{I}^c}(g)}_{\CM}=\int_{\R} f(\xi) u(\xi) \mathd \xi
\]
where $u=[(H^{\beta}+m^2)^{-1}Q_{\mathcal{I}^c}](g)$  is the unique (weak) solution to the Dirichlet problem \eqref{eq:Dirichletg}.
Taking into account that $Q_{\mathcal{I}^c}$ is an orthogonal projection on $\mathcal{H}_{\CM}$, if $h\in \mathcal{H}_{\mathcal{I}^c}$ we have $Q_{\mathcal{I}^c}h=0$ and 
\[
\int_{\R} h(\xi)u(\xi)\mathd \xi=\int_{\R} h(\xi)[(H^{\beta}+m^2)^{-1}Q_{\mathcal{I}^c}](g)(\xi) \mathd \xi=\scal{h}{Q_{\mathcal{I}^c}g}_{\CM}=\scal{Q_{\mathcal{I}^c}h}{Q_{\mathcal{I}^c}g}_{\CM}=0.
\]
By the definition \eqref{eq:HIc} of $\mathcal{H}_{\mathcal{I}^c}$, we get $u|_{\mathcal{I}^c}=0$. Furthermore, if $\bar h\in C^{\infty}_c(\mathcal{I})$,
\begin{multline}
\int_\R (H^{\beta}+m^2)\bar h(\xi) u(\xi)\mathd \xi=\int_\R \bar h(\xi) Q_{\mathcal{I}^c}(g)(\xi)\mathd \xi\\
=\int_\R  \bar h(\xi) g(\xi)\mathd \xi-\int_\R  \bar h(\xi) P_{\mathcal{I}^c}(g)(\xi)\mathd \xi=\int_\R  \bar h(\xi) g(\xi)\mathd \xi,
\end{multline}
where we used that $\supp(P_{\mathcal{I}^c}(g))\subset \overline{\mathcal{I}^c}$, so that $\bar hP_{\mathcal{I}^c}(g)=0$. Hence $u$ solves \eqref{eq:Dirichletg}. 

Finally, \eqref{eq:meanconditionalV} is \eqref{eq:meanconditional}.
\end{proof}

\begin{rmk}
By Remark \ref{rmk:conc-misura}, we have $\mu_{\free}(C(\R))=1$ and so that we can restrict to consider $y \in L^2_{r,k}(\mathcal{I}^c)\cap C(\mathcal{I}^c)$. This fact is fundamental to guarantee that $y$ admits a harmonic extension in the case where $\beta_0>0$ in condition \ref{H0}.
\end{rmk}
Let $\mathcal{I}$ be a bounded open interval of $\mathbb{R}$ and let $y\in L^2_{r,k}(\mathcal{I}^c)$. Taking into account that $E=L^2_{r,k}$ and recalling \eqref{eq:defEIy} (with $E=L^2_{r,k}$), the set $E^{\mathcal{I},y}$ is, in general, not a Banach space (indeed, it is not even a vector space), except in the case where $y\equiv 0$. Moreover, note that $E^{\mathcal{I},0}$ is isomorphic to $L^2(\mathcal{I})$. Therefore, $\mu_{\mathrm{free}}^{\mathcal{I},0}$ is a Gaussian measure on a Hilbert space, and standard infinite-dimensional analysis can be applied.

\begin{prop}\label{prop:gaussL2I}
Assume \ref{H0}, \ref{H1} and \ref{H3} hold.  Then for every open bounded interval $\mathcal{I} \subset \R$, the measure $\mu_{\free}^{\mathcal{I},0}$, namely the image of $\mu_{\free}^{\mathcal{I},y}$ under the map $x\mapsto x-\varphi^y$, which for $\mu_{\free}^{\mathcal{I}^c}$-almost every $y$ does not depend on $y$, is Gaussian on $E^{\mathcal{I},0}\cong L^2(\mathcal{I})$ with mean $0$ and covariance operator $\mathcal{A}_{\mathcal{I}}^{-1}$ where $\mathcal{A}_{\mathcal{I}}$ is the realization of $H^\beta+m^2$ in $L^2(\mathcal{I})$ with homogeneous Dirichlet boundary conditions. Moreover $\mu_{\free}^{\mathcal{I},0}$ is the unique invariant measure of the following SPDE
\begin{equation}
  \begin{cases}
      \mathd Z_{\mathcal{I}}(t) + \mathcal{A}_{\mathcal{I}}Z_{\mathcal{I}}(t)\mathd t= \sqrt{2}\mathd W_\mathcal{I}(t),\qquad t> 0,\\
      Z_{\mathcal{I}}(0)=x\in L^2(\mathcal{I}),
  \end{cases}
  \label{eq:SPDElinearL2(I)}
\end{equation}
where $\{W_{\mathcal{I}}(t)\}_{t\geq 0}$ is a $L^2(\mathcal{I})$-cylindrical Wiener process.
\end{prop}
\begin{proof}
We first check that $\mu^{\mathcal{I},0}_{\free}$ is a Gaussian measure on the Hilbert space $E^{\mathcal{I},0}\cong L^2(\mathcal{I})$. Since the topological dual of $L^2(\mathcal{I})$ is generated by the maps $x\mapsto\scal{x}{f}_{L^2}$ with $f\in C^{\infty}_c(\mathcal{I})$, it suffices to show that $\scal{X}{f}_{L^2}$ is a real Gaussian  random variable. Fix $f\in C^{\infty}_c(\mathcal{I})$ and $t\in\R$. Arguing as in the proof of Lemma \ref{lemma:disintegration-conditional}, with $\scal{x}{f}_{L^2}$ replaced by $e^{it\scal{x}{f}_{L^2}}$, we conclude
\[
\int_{E^{\mathcal{I},y}}e^{it\scal{x-\varphi^y}{f}_{L^2}}\,\mu^{\mathcal{I},y}_{\free}(\mathd x)
=e^{-\frac{t^2}{2}\scal{Q_{\mathcal{I}^c}(f)}{Q_{\mathcal{I}^c}(f)}_{\CM}}\qquad\text{for }\mu_{\free}^{\mathcal{I}^c}\text{-almost every }y.
\]
 The last statement follows from \cite[Section 6.2]{Dap-Zab1996}.
\end{proof}

Let $\mathcal{I}$ be a bounded open interval of $\mathbb{R}$ and let $y\in L^2_{r,k}(\mathcal{I}^c)$. As we just said $E^{\mathcal{I},y}$ is not even a vector space, however, in the Gaussian case, there is a natural bijection $\mathfrak{J}_y$ from $E^{\mathcal{I},y}$ to $E^{\mathcal{I},0}$ given by
\[
\mathfrak{J}_yf=f-\varphi^y
\]
where $\varphi^y$ is the harmonic extension of $y$ in the sense of Definition \ref{def:harmonicextension}. Due to this fact we obtain the following result, that will be fundamental to exploit Proposition \ref{prop:gaussL2I}.

\begin{cor}\label{corollary:Guassiantranslation}
Assume \ref{H0}, \ref{H1} and \ref{H3} hold.  Then for every open bounded interval $\mathcal{I} \subset \R$, for $\mu_{\free}^{\mathcal{I}^c}$-almost every $y \in L^2_{r,k}(\mathcal{I}^c)\cap C(\mathcal{I}^c)$ and $\Phi\in B_b(L^2_{r,k})$ we have
\begin{equation}\label{eq:translationconditionalGaussian}
  \int_{E^{\mathcal{I}, y}} \Phi (z) \mu_{\free}^{\mathcal{I}, y}
     (\mathd z) = \int_{E^{\mathcal{I}, 0}} \Phi (z + \varphi^y)
     \mu_{\free}^{\mathcal{I},0} (\mathd z),
\end{equation}
where $\varphi^y$ is the harmonic extension of $y$ in the sense of Definition \ref{def:harmonicextension} and $E^{\mathcal{I},y}$ is given by \eqref{eq:defEIy} (with $E=L^2_{r,k}$), $\mu_{\free}^{\mathcal{I},y}$ is given by Lemma \ref{lemma:disintegration}.\\
Moreover, the map $y \longmapsto\mu_{\free}^{\mathcal{I}, y}$ is continuous from $L^2_{r,k}(\mathcal{I}^c)\cap C(\mathcal{I}^c)$, endowed with the topology induced by the norm of $L^2_{r,k}(\mathcal{I}^c)$ together with the uniform convergence on compact subsets of $\mathcal{I}^c$, into $\mathcal{P}(L^2_{r,k})$, endowed with the topology of the weak convergence of measures. In particular \eqref{eq:translationconditionalGaussian} holds for every $y \in L^2_{r,k}(\mathcal{I}^c)\cap C(\mathcal{I}^c)$.
\end{cor}
\begin{proof}
The bijection $\mathfrak{J}_y$ and Propositions \ref{proposition:conditionalGaussian} and \ref{prop:gaussL2I} yield \eqref{eq:translationconditionalGaussian}.
For the second assertion, let $y_n\rightarrow y$ in $L^2_{r,k}(\mathcal{I}^c)\cap C(\mathcal{I}^c)$ with respect to the topology in the statement. Since $\varphi^{y_n}-\varphi^{y}=y_n-y$ on $\mathcal{I}^c$ and $\mathcal{I}$ is bounded, by Theorem \ref{thm:harmonic-continuity} we get
\[
\norm{\varphi^{y_n}-\varphi^{y}}^2_{L^2_{r,k}}
\leq\norm{\rho_{r,k}}_{L^1(\R)}\norm{\varphi^{y_n}-\varphi^{y}}^2_{L^\infty(\mathcal{I})}
+\norm{y_n-y}^2_{L^2_{r,k}(\mathcal{I}^c)}\longrightarrow0 ,
\]
so that $\varphi^{y_n}\rightarrow\varphi^{y}$ in $L^2_{r,k}$. Hence, for every $\Phi\in C_b(L^2_{r,k})$, by \eqref{eq:translationconditionalGaussian} and the dominated convergence theorem we obtain the continuity
\[
\int_{E^{\mathcal{I},y_n}}\Phi(z)\,\mu_{\free}^{\mathcal{I},y_n}(\mathd z)
=\int_{E^{\mathcal{I},0}}\Phi(z+\varphi^{y_n})\,\mu_{\free}^{\mathcal{I},0}(\mathd z)
\longrightarrow\int_{E^{\mathcal{I},0}}\Phi(z+\varphi^{y})\,\mu_{\free}^{\mathcal{I},0}(\mathd z).
\]
\end{proof}

\begin{rmk}
Let $y,y'\in L_{r,k}^2(\mathcal{I}^c)$. It is important to note that $\mu_{\free}^{\mathcal{I},y}$ and $\mu_{\free}^{\mathcal{I},y'}$ cannot be absolutely continuous with respect to each other when $y\not= y'$. Indeed, first, the two measures are supported in disjoint subsets, but also the translation vector transforming the first into the second, namely $\varphi^{y'}-\varphi^{y}=\varphi^{y'-y}$, cannot belong to the  space of the two measures since  $\varphi^{y'-y}|_{\partial \mathcal{I}}\not=0$. What one can prove is that, for any open interval $\mathcal{J}$ such that $\bar{\mathcal{J}} \subset \mathcal{I}$, the measures $\mu_{\free}^{\mathcal{I},y}$ and $\mu_{\free}^{\mathcal{I},y'}$, restricted to the set $\mathcal{J}$ are absolutely continuous with respect to each other. Indeed, by the regularity properties of harmonic extensions (see Remark \ref{rmk:continuityharmonic}) the function $\varphi^{y}|_{\mathcal{J}}-\varphi^{y'}|_{\mathcal{J}}$ belongs to $H^{s_1}_{\loc}$ which is equivalent (as normed space) to the Cameron-Martin spaces of $R_{\mathcal{J}*}\mu_{\free}^{\mathcal{I},y}$ and $R_{\mathcal{J}*}\mu_{\free}^{\mathcal{I},y'}$, see \eqref{eq:Prestriction}. Furthermore, one can also prove that $R_{\mathcal{J}*}\mu_{\free}^{\mathcal{I},y}$ and $R_{\mathcal{J}*}\mu_{\free}^{\mathcal{I},y'}$ are absolutely continuous with respect to $R_{\mathcal{J}*}\mu_{\free}$, where no boundary conditions are considered.
\end{rmk}

\subsection{DLR equations on $\R$ for \eqref{eq:SPDEmain1psuedo}}\label{sect:DLRequazioni}
Before stating the DLR equations for \eqref{eq:SPDEmain1psuedo}, we need an additional assumption to guarantee that the Gibbs measure is well-defined.
\begin{enumerate}[label={{(H4)}}]
\item\label{H4} \emph{There exists $M_1\in\R$ such that $V(\xi)\geq M_1$ for every $\xi\in\R$.}
\end{enumerate}
Indeed, under \ref{H4} we have $\exp(-\int_{\mathcal{I}}V(x(\xi))\mathd\xi)\leq \exp(-M_1|\mathcal{I}|)$ for every bounded interval $\mathcal{I}$, so that
\[
\int_{E^{\mathcal{I},y}} \exp\left(-\int_{\mathcal{I}}V(x(\xi))\mathd \xi\right) \mu^{\mathcal{I},y}_{\free}(\mathd x)<+\infty.
\]
\begin{rmk}
    Note that \ref{H4} is automatically satisfied when the constant $M_0$ in \eqref{H2:diss} is strictly positive.
\end{rmk}

\begin{rmk}\label{rmk:H4}
Hypothesis \ref{H4} excludes potentials of the type $V(\xi)=-\frac{\lambda}{2}\xi^2$ with $\lambda>0$. If $\lambda>m^2$, then \eqref{eq:SPDEmain1Abstract} becomes a linear equation associated with the operator $H^\beta+m^2-\lambda$, which is not positive, and it has no invariant measure. Moreover, the DLR equations \eqref{eq:DLRequations} are not even meaningful, since the normalization constant $\mathcal Z^{\mathcal I}(y)$ in Definition \ref{definition:DLR} is infinite whenever the smallest eigenvalue of the operator associated with the Dirichlet problem \eqref{eq:Dirichletg} is smaller than $\lambda$, and this happens for every bounded interval $\mathcal I$ large enough. If instead $\lambda<m^2$, the potential can be absorbed in the mass term, by replacing $m^2$ with $m^2-\lambda$ and $V$ with $0$.
\end{rmk}

\begin{defi}\label{definition:DLR}
Assume that \ref{H0},\ref{H1}, \ref{H2}, \ref{H3} and \ref{H4} hold. Given a Borel probability measure $\mu$ on $L^2_{r,k}$, we say that $\mu$ satisfies the Dobrushin-Lanford-Ruelle equations (or simply DLR equations) for \eqref{eq:SPDEmain1psuedo} if $\mu\left(L^{2q}_{r,k}\cap C(\R)\right)=1$ and, for every bounded open interval $\mathcal{I}$ of $\R$ and $\mu^{\mathcal I^c}$-almost every $y \in L^2_{r,k}(\mathcal{I}^c)\cap C(\mathcal{I}^c)$, we have
\begin{equation}
\frac{\mathd \mu^{\mathcal{I},y}}{\mathd \mu_{\free}^{\mathcal{I},y}}= \frac{1}{\mathcal{Z}^{\mathcal{I}}(y)}\exp\left(-\int_{\mathcal{I}}V(x(\xi))\mathd \xi\right),\label{eq:DLRequations}
\end{equation}
where
\[
\mathcal{Z}^{\mathcal{I}}(y)=\int_{E^{\mathcal{I},y}} \exp\left(-\int_{\mathcal{I}}V(x(\xi))\mathd \xi\right) \mu_{\free}^{\mathcal{I},y}(\mathd x),
\]
$E^{\mathcal{I},y}$ is given by \eqref{eq:defEIy} (with $E=L^2_{r,k}$), and $\mu^{\mathcal{I},y}$ and $\mu_{\free}^{\mathcal{I},y}$ are given by Lemma \ref{lemma:disintegration}.
\end{defi}

By Corollary \ref{corollary:Guassiantranslation} we deduce the following result.

\begin{cor}\label{remark:DLRtranslated}
    Assume that \ref{H0}, \ref{H1}, \ref{H2}, \ref{H3} and \ref{H4} hold.  Let $\mu$ be a Borel probability measure on $L^2_{r,k}$ such that $\mu(C(\R))=1$. Then  $\mu$ solves DLR equations for \eqref{eq:SPDEmain1psuedo} if and only if for every open bounded interval $\mathcal{I} \subset \R$, for $\mu^{\mathcal{I}^c}$-almost every $y \in L^2_{r,k}(\mathcal{I}^c)\cap C(\mathcal{I}^c)$ (see \eqref{eq:legge-Ic}) we have
    \begin{align}
    \int_{E^{\mathcal{I},y}}\Phi(x)\mu^{\mathcal{I},y}(\mathd x)&=\int_{E^{\mathcal{I},0}}\Phi(x+\varphi^y)\frac{\exp\left(-\int_{\mathcal{I}}V(x(\xi)+\varphi^y(\xi))\mathd \xi\right)}{\mathcal{Z}^{\mathcal{I}}(y)}\mu_{\free}^{\mathcal{I},0}(\mathd x),\label{eq:tras-mu}\\
\mathcal{Z}^{\mathcal{I}}(y)&=\int_{E^{\mathcal{I},0}} \exp\left(-\int_{\mathcal{I}}V(x(\xi)+\varphi^{y}(\xi))\mathd \xi\right) \mu^{\mathcal{I},0}_{\free}(\mathd x),\notag
 \end{align}
 where $\varphi^y$ is the harmonic extension of $y$ in the sense of Definition \ref{def:harmonicextension}, $E^{\mathcal{I},y}$ is given by \eqref{eq:defEIy} (with $E=L^2_{r,k}$) and $\mu^{\mathcal{I},y}$, $\mu_{\free}^{\mathcal{I},y}$ are given by Lemma \ref{lemma:disintegration}.\\
 Moreover, the map $y \longmapsto\mu^{\mathcal{I}, y}$ is continuous from $L^2_{r,k}(\mathcal{I}^c)\cap C(\mathcal{I}^c)$, endowed with the topology induced by the norm of $L^2_{r,k}(\mathcal{I}^c)$ together with the uniform convergence on compact subsets of $\mathcal{I}^c$, into $\mathcal{P}(L^2_{r,k})$, endowed with the topology of the weak convergence of measures and also the map $y\mapsto\mathcal{Z}^{\mathcal{I}}(y)$ is continuous with respect to the same topology. In particular \eqref{eq:tras-mu} holds for every $y \in L^2_{r,k}(\mathcal{I}^c)\cap C(\mathcal{I}^c)$.
\end{cor}

\section{Integration by parts formula and DLR equation}\label{sect:int-DLR}

In this section, we are going to prove the main result of the paper

\begin{thm}\label{theorem:main}
Assume that \ref{H0}, \ref{H1}, \ref{H2}, \ref{H3} and \ref{H4} hold. Let $q\geq 1$ be the constant given by \eqref{H2:pol}.  Let $\mu$ be a Borel probability measure on $L^2_{r,k}$ such that $\mu(L^{2q}_{r,k}\cap C(\R))=1$ and
\begin{equation}\label{eq:momenti-DLR}
\int_{L^{2q}_{r,k}}\norm{x}_{L^{2q}_{r,k}}^q\mu(\mathd x)<+\infty.
\end{equation}
Then the following conditions are equivalent.
\begin{enumerate}
  \item $\mu$ satisfies the integration by parts formula associated with the SPDE \eqref{eq:SPDEmain1psuedo} (with the vector logarithmic derivative $-\mathcal{A}+\mathcal{B}$), in the sense of point 3. of Definition \ref{definition:reversible};
  \item $\mu$ satisfies the DLR equations for \eqref{eq:SPDEmain1psuedo}, in the sense of Definition \ref{definition:DLR}.
\end{enumerate}
\end{thm}
\begin{rmk}
We point out that, thanks to the a priori estimates in Section \ref{section:Fernique}, condition \eqref{eq:momenti-DLR} can be removed (see Theorem \ref{theorem:super-main}).
\end{rmk}

Before proving the main theorem, we need several preliminary results concerning the behavior of the SPDE \eqref{eq:SPDEmain1Abstract} on $L^2(\mathcal{I})$, where $\mathcal{I}$ is a bounded open interval of $\R$.

\subsection{The SPDE on the bounded interval}
Assume that \ref{H0}, \ref{H1}, \ref{H2}, \ref{H3}  and \ref{H4} hold. Let $q\geq 1$ be the constant given by \eqref{H2:pol}.\\
In this subsection, we fix an open bounded interval $\mathcal{I}$ of $\R$ and a boundary condition $y \in L^2_{r,k}(\mathcal{I}^c)\cap C(\mathcal{I}^c)$. We consider the SPDE
\begin{equation}
  \begin{cases}
      \mathd X_{\mathcal{I}}(t) + \mathcal{A}_{\mathcal{I}}X_{\mathcal{I}}(t)\mathd t= -\nabla U(X_{\mathcal{I}}(t))\mathd t+\sqrt{2}\mathd W_\mathcal{I}(t),\qquad t> 0,\\
      X_{\mathcal{I}}(0)=x\in L^2(\mathcal{I}),
  \end{cases}
  \label{eq:SPDEL2(I)}
\end{equation}
where $\mathcal{A}_{\mathcal{I}}$ is the realization of $H^\beta+m^2$ in $L^2(\mathcal{I})$ with homogeneous Dirichlet boundary conditions (namely the operator associated with the problem \eqref{eq:Dirichletg}), $\{W_{\mathcal{I}}(t)\}_{t\geq 0}$ is a $L^2(\mathcal{I})$-cylindrical Wiener process, $U:L^2(\mathcal{I})\rightarrow \R$ is defined by
\begin{equation}
    U(x):=\begin{cases}
        \int_{\mathcal{I}}V(x(\xi)+\varphi^y(\xi))\mathd \xi,\qquad &x\in C(\overline{\mathcal{I}})\\
        +\infty\qquad &x\not\in C(\overline{\mathcal{I}})
    \end{cases}
\end{equation}
and $\varphi^y$ is the harmonic extension of $y$ in the sense of Definition \ref{def:harmonicextension} (see also Theorem \ref{thm:harmonic-continuity}). In this subsection, we are going to prove that \eqref{eq:SPDEL2(I)} is well-posed and admits a unique invariant measure.\\

First of all, we note that, by \cite{Kaleta-Kwa-Mal2016}, $\{e^{-t(\mathcal{A}_{\mathcal{I}}-m^2)}\}_{t\geq 0}$ is the transition semigroup of the L\'evy process associated with $-H^\beta$ killed upon leaving $\mathcal{I}$; in particular it is a symmetric sub-Markovian semigroup on $L^2(\mathcal{I})$, and so it extends to a strongly continuous contraction semigroup on $L^p(\mathcal{I})$ for every $p\in[1,+\infty)$ (see e.g. \cite{Fukushima2011}). Hence $-\mathcal{A}_{\mathcal{I}}$ generates a strongly continuous semigroup $\{e^{-t\mathcal{A}_{\mathcal{I}}}\}_{t\geq 0}$ on $L^p(\mathcal{I})$ with $\norm{e^{-t\mathcal{A}_{\mathcal{I}}}}_{\mathcal L(L^p(\mathcal{I}))}\leq e^{-m^2t}$ for every $p\geq 1$.

 We need to prove the regularity of the stochastic convolution associated with \eqref{eq:SPDEL2(I)},
\begin{equation}\label{eq:stochconvI}
W_{\mathcal A_{\mathcal{I}}}(t):=\sqrt{2}\int_0^t e^{-(t-s)\mathcal A_{\mathcal{I}}}\mathd W_{\mathcal{I}}(s),\qquad t\geq 0.
\end{equation}

\begin{lemm}\label{lemma:stochconvL2I}
Assume that conditions \ref{H0} and \ref{H1} hold. The stochastic convolution $\{W_{\mathcal A_{\mathcal{I}}}(t)\}_{t\geq 0}$ defined in \eqref{eq:stochconvI} admits a version with continuous trajectories in $C(\overline{\mathcal{I}})$ and, for every $T>0$ and $n\in\N$, we have
\begin{equation}\label{eq:stochconvboundI}
\mathbb{E}\left[\sup_{t\in[0,T]}\norm{W_{\mathcal A_{\mathcal{I}}}(t)}_{C(\overline{\mathcal{I}})}^{n}\right]<+\infty.
\end{equation}
\end{lemm}

\begin{proof}
We argue as in the proof of Lemma \ref{lemma:stochconv}. By \cite[Section 2.4]{Kaleta-Kwa-Mal2016}, $\{e^{-\tau\mathcal A_{\mathcal{I}}}\}_{\tau\geq0}$ is a strongly continuous contraction semigroup on $C_0(\mathcal{I})$. As in \cite[Section 3.5]{Kaleta-Kwa-Mal2016}, we denote by $g_{\mathcal{I}}:[0,+\infty)\times\mathcal{I}\times\mathcal{I}\rightarrow [0,+\infty)$ the transition kernel of $\{e^{-\tau(\mathcal{A}_{\mathcal{I}}-m^2)}\}_{\tau\geq0}$, and we recall that
\begin{equation}\label{eq:killeddomination}
0\leq g_{\mathcal{I}}(\tau,\xi,\xi')\leq g(\tau,\xi-\xi'),\qquad \tau>0,\qquad \xi,\xi'\in\mathcal{I},
\end{equation}
where $g$ is the kernel appearing in \eqref{eq:kernelL2}. Since $s_1>\tfrac12$, we can fix
\[
p>\frac{2}{2s_1-1},\qquad \alpha\in\left(\frac{1}{2s_1p},\ \frac12-\frac{1}{4s_1}\right),
\]
and set
\[
Y_{\mathcal{I}}(s):=\sqrt2\int_0^s (s-\sigma)^{-\alpha}e^{-(s-\sigma)\mathcal A_{\mathcal{I}}}\,\mathd W_{\mathcal{I}}(\sigma),\qquad s\in[0,T],
\]
By the It\^o isometry and \eqref{eq:killeddomination},
\[
\mathbb E\big[\abs{Y_{\mathcal{I}}(s,\xi)}^2\big]=2\int_0^s\tau^{-2\alpha}e^{-2m^2\tau}\norm{g_{\mathcal{I}}(\tau,\xi,\cdot)}^2_{L^2(\mathcal{I})}\,\mathd\tau
\leq 2\int_0^s\tau^{-2\alpha}e^{-2m^2\tau}\norm{g(\tau,\cdot)}^2_{L^2}\,\mathd\tau=v_\alpha(s),
\]
with $v_\alpha$ as in \eqref{eq:pointwisevariance}, which is finite since $\alpha<\tfrac12-\tfrac{1}{4s_1}$. Exactly as in \eqref{eq:Lpmoment} and \eqref{eq:gaussianmoments} we get, for every $l\geq 2$,
\begin{equation}\label{eq:LpmomentI}
\mathbb E\big[\norm{Y_{\mathcal{I}}(s)}^p_{L^p(\mathcal{I})}\big]=\int_{\mathcal{I}}\mathbb E\big[\abs{Y_{\mathcal{I}}(s,\xi)}^p\big]\mathd\xi\leq c_p\abs{\mathcal{I}}\,v_\alpha(s)^{p/2},\qquad
\sup_{s\in[0,T]}\mathbb E\big[\norm{Y_{\mathcal{I}}(s)}^{l}_{L^p(\mathcal{I})}\big]<+\infty,
\end{equation}
so that $Y_{\mathcal{I}}\in L^{l}([0,T];L^p(\mathcal{I}))$ almost surely.\\
Now we estimate $\norm{g(\tau,\cdot)}_{L^{p^*}(\R)}$, where $p^*$ denotes the conjugate exponent of $p$.  Splitting the integral over $\{\abs{\xi}\leq 1\}$ and $\{\abs{\xi}\geq 1\}$ and using \eqref{eq:symbol-lower},  for every $\tau\in(0,T]$, we get
\[
\norm{g(\tau,\cdot)}_{L^{\infty}(\R)}\leq\frac{1}{2\pi}\int_\R e^{-\tau\widehat{H^\beta}(\xi)}\,\mathd\xi
\leq\frac{1}{2\pi}\left(2+\int_\R e^{-c\tau\abs{\xi}^{2s_1}}\mathd\xi\right)\leq C_T\,\tau^{-\frac{1}{2s_1}},
\]
for some $C_T>0$. Since $\norm{g(\tau;\cdot)}_{L^{1}(\R)}\leq 1$ and recalling that $1-\tfrac{1}{p^*}=\tfrac1p$, we obtain
\begin{equation}\label{eq:ultracontr}
\norm{g(\tau,\cdot)}_{L^{p^*}(\R)}\leq\norm{g(\tau,\cdot)}_{L^{1}(\R)}^{\frac{1}{p^*}}\norm{g(\tau,\cdot)}_{L^{\infty}(\R)}^{1-\frac{1}{p^*}}\leq C_T\,\tau^{-\frac{1}{2s_1p}} .
\end{equation}
Hence, by \eqref{eq:killeddomination} and the H\"older inequality, for every $u\in L^p(\mathcal{I})$ and $\tau\in(0,T]$, we get
\[
\norm{e^{-\tau\mathcal A_{\mathcal{I}}}u}_{C(\overline{\mathcal{I}})}\leq\sup_{\xi\in\mathcal{I}}\norm{g_{\mathcal{I}}(\tau,\xi,\cdot)}_{L^{p^*}(\mathcal{I})}\norm{u}_{L^p(\mathcal{I})}
\leq\norm{g(\tau,\cdot)}_{L^{p^*}(\R)}\norm{u}_{L^p(\mathcal{I})}\leq C_T\tau^{-\frac{1}{2s_1p}}\norm{u}_{L^p(\mathcal{I})} .
\]
Recalling the factorization formula, for every $t>0$,
\begin{equation*}
W_{\mathcal A_{\mathcal{I}}}(t)=\frac{\sin(\pi\alpha)}{\pi}\int_0^t (t-s)^{\alpha-1}e^{-(t-s)\mathcal A_{\mathcal{I}}}Y_{\mathcal{I}}(s)\mathd s,
\end{equation*}
and so that
\[
\norm{W_{\mathcal A_{\mathcal{I}}}(t)}_{C(\overline{\mathcal{I}})}\leq  C_T\int_0^t(t-s)^{\alpha-1-\frac{1}{2s_1p}}\norm{Y_{\mathcal{I}}(s)}_{L^p(\mathcal{I})}\mathd s,
\]
Since $\alpha>\tfrac{1}{2s_1p}$ and choosing $l>\big(\alpha-\tfrac{1}{2s_1p}\big)^{-1}$, by the same arguments used in Lemma \ref{lemma:stochconv} the proof is concluded.
\end{proof}    

Taking into account Proposition \ref{prop:gaussL2I} and Lemma \ref{lemma:stochconvL2I} we have
\[
\mu_{\free}^{\mathcal{I},0}(C(\overline{\mathcal{I}}))=1.
\]
Moreover, by \ref{H4} we have $U(x)\geq M_1\abs{\mathcal{I}}$ for every $x\in L^2(\mathcal{I})$, and by \cite[Proposition 5.2]{Dap-Lun2014}, $U\in W^{1,2}(L^2(\mathcal{I}),\mu_{\free}^{\mathcal{I},0})$ and
\begin{equation}\label{eq:deriU}
\nabla U(x)(\xi)=V'(x(\xi)+\varphi^y(\xi)),\qquad x\in C(\overline{\mathcal{I}}).
\end{equation}
We can now prove the well-posedness of \eqref{eq:SPDEL2(I)}.

\begin{prop}\label{prop:wellposednessI}
    Assume that assumptions \ref{H0}, \ref{H1}, \ref{H2} and \ref{H3}  hold. Let $q\geq 1$ be the constant given by \eqref{H2:pol}. For every $x\in L^{2q}(\mathcal{I})$ equation \eqref{eq:SPDEL2(I)} has a unique mild solution $\{X(t,x)\}_{t\geq 0}$ having continuous trajectories $\mathbb{P}$-a.s. in $L^{2q}(\mathcal{I})$, namely
    \[
    X(t,x)=e^{-t\mathcal{A}_{\mathcal{I}}}x-\int_0^te^{-(t-s)\mathcal{A}_{\mathcal{I}}}\nabla U(X(s,x))\mathd s+\sqrt{2}\int_0^te^{-(t-s)\mathcal{A}_{\mathcal{I}}}\mathd W_{\mathcal{I}}(s),\qquad \mathbb{P}\text{-a.s.}
    \]
 Furthermore, for every $x\in L^{2}(\mathcal{I})$ equation \eqref{eq:SPDEL2(I)} has a unique generalized mild solution $\{X(t,x)\}_{t\geq 0}$ having continuous trajectories $\mathbb{P}$-a.s. in $ L^{2}(\mathcal{I})$, namely for every sequence $\{x_n\}_{n\in\N}\subseteq L^{2q}(\mathcal{I})$ such that $x_n\rightarrow x$ in $L^{2}(\mathcal{I})$ we have
 \[
 \lim_{n\to+\infty}\mathbb{E}\left[\sup_{t\in [0,T]}\norm{X(t,x)-X(t,x_n)}_{L^{2}(\mathcal{I})}\right]=0,\qquad T\geq 0,
 \]
 where $\{X(t,x_n)\}_{t\geq 0}$ is the unique mild solution to \eqref{eq:SPDEL2(I)} with initial datum $x_n$.
\end{prop}
\begin{proof}
The statement follows from Theorem \ref{th:solmild} exactly as in Proposition \ref{prop:wellposedness}, once Hypotheses \ref{Hyp:base} are checked with 
\[
H=L^2(\mathcal{I}),\quad E=L^{2q}(\mathcal{I}),\quad R=\sqrt{2}\,\Id_{L^2(\mathcal{I})},\quad A=-\mathcal{A}_{\mathcal{I}},\quad F(x)=-\nabla U(x)+M_0x,\quad G(x)=-M_0x,
\]
where $M_0$ is the constant in \eqref{H2:diss} and $\nabla U(x):=V'(x+\varphi^y)$ for every $x\in L^{2q}(\mathcal{I})$, consistently with \eqref{eq:deriU}. We only point out the differences with respect to Proposition \ref{prop:wellposedness}.\\
Hypotheses \ref{Hyp:base}(1) holds since $L^{2q}(\mathcal{I})$ is reflexive and dense in $L^2(\mathcal{I})$.\\
Hypotheses \ref{Hyp:base}(4) is checked as in Proposition \ref{prop:wellposedness}, with $\zeta=-M_0$.\\
Hypotheses \ref{Hyp:base}(5) holds with $\varepsilon=m^2$ by the contractivity of $\{e^{-t(\mathcal{A}_{\mathcal{I}}-m^2)}\}_{t\geq 0}$ on $L^2(\mathcal{I})$ and $L^{2q}(\mathcal{I})$ recalled at the beginning of this subsection.\\
Hypotheses \ref{Hyp:base}(6) follows from \eqref{H2:diss}, \eqref{H2:pol} and the fact that $\varphi^y$ is bounded on $\overline{\mathcal{I}}$ (see Theorem \ref{thm:harmonic-continuity}); in particular we note that for every $\xi\in\mathcal{I}$ the map $r\in\R\rightarrow V'(r+\varphi^y(\xi))-M_0r$ is nondecreasing.\\
Hypotheses \ref{Hyp:base}(7) follows from Lemma \ref{lemma:stochconvL2I}, since $C(\overline{\mathcal{I}})\subseteq L^{2q}(\mathcal{I})$ continuously, and from Hypotheses \ref{Hyp:base}(6).\\
Finally, Hypotheses \ref{Hyp:base}(\ref{Hyp:base8}) follows from \eqref{H2:der} and the H\"older inequality exactly as in Proposition \ref{prop:wellposedness}, since $\varphi^y$ is bounded on $\overline{\mathcal{I}}$.
\end{proof}

By Corollary \ref{cor:gradient-invariant} we obtain the following result.

\begin{prop}\label{prop:invmeasureI}
    Assume that \ref{H0}, \ref{H1}, \ref{H2}, \ref{H3}  and \ref{H4} hold. The SPDE \eqref{eq:SPDEL2(I)} has a unique invariant measure $\nu_{\mathcal{I},y}$ defined by
    \begin{equation}\label{eq:invmeasureI}
    \nu_{\mathcal{I},y}(\mathd x):=\frac{1}{\mathcal{Z}^{\mathcal{I}}(y)}\exp\left(-\int_{\mathcal{I}}V(x(\xi)+\varphi^y(\xi))\,\mathd\xi\right)\mu_{\free}^{\mathcal{I},0}(\mathd x),
\end{equation}
where $\mu_{\free}^{\mathcal{I},0}$ is the Gaussian measure with mean $0$ and covariance operator $\mathcal{A}_{\mathcal{I}}^{-1}$ given by Proposition \ref{prop:gaussL2I} and
\[
\mathcal{Z}^{\mathcal{I}}(y)=\int_{L^2(\mathcal{I})} \exp\left(-\int_{\mathcal{I}}V(x(\xi)+\varphi^{y}(\xi))\mathd \xi\right) \mu^{\mathcal{I},0}_{\free}(\mathd x)\in(0,+\infty).
\]
Moreover, for every $h\in {\rm Dom}(\mathcal{A}_{\mathcal{I}})$ and $\Phi\in\mathcal{F}C^{\infty}_b(L^2_{r,k})$,
\begin{equation}\label{eq:ibpI}
\int_{L^2(\mathcal{I})}\scal{h}{\delta_x\Phi(x)}_{L^2(\mathcal{I})}\nu_{\mathcal{I},y}(\mathd x)
=\int_{L^2(\mathcal{I})}\left[\scal{\mathcal{A}_{\mathcal{I}}h}{x}_{L^2(\mathcal{I})}+\scal{h}{V'(x+\varphi^y)}_{L^2(\mathcal{I})}\right]\Phi(x)\nu_{\mathcal{I},y}(\mathd x).
\end{equation}
In particular, if $h\in C^{\infty}_c(\mathcal{I})$, then $\mathcal{A}_{\mathcal{I}}h=(H^{\beta}+m^2)h$ on $\mathcal{I}$ and \eqref{eq:ibpI} reads
\begin{equation}\label{eq:ibpIbis}
\int_{L^2(\mathcal{I})}\scal{h}{\delta_x\Phi(x)}_{L^2}\nu_{\mathcal{I},y}(\mathd x)
=\int_{L^2(\mathcal{I})}\left[\scal{(H^{\beta}+m^2)h}{x}_{L^2}+\scal{h}{V'(x+\varphi^y)}_{L^2}\right]\Phi(x)\nu_{\mathcal{I},y}(\mathd x).
\end{equation}
\end{prop}

\begin{proof}
We apply Corollary \ref{cor:gradient-invariant} with $H=L^2(\mathcal I)$, $E=L^{2q}(\mathcal I)$, $R=\sqrt2\,\Id_{L^2(\mathcal I)}$ (so that $\gamma=\mu^{\mathcal I,0}_{\free}$ by Proposition \ref{prop:gaussL2I}), $A=-\mathcal A_{\mathcal I}$, which is self-adjoint, $\mathcal D={\rm Dom}(\mathcal A_{\mathcal I})$, which is contained in $E$ since its elements are continuous functions on $\R$ vanishing on $\mathcal I^c$ (see \cite[Theorem 4.4.3]{Fukushima2011} and \eqref{eq:energy-embedding}), and $F+G=-\nabla U$, where $U$ is bounded from below and belongs to $W^{1,2}(L^2(\mathcal I),\mu^{\mathcal I,0}_{\free})$. Moreover $\int\norm{x}^q_{L^{2q}(\mathcal I)}\,\nu_{\mathcal I,y}(\mathd x)<+\infty$, since, by \ref{H4}, $\nu_{\mathcal I,y}$ has a bounded density with respect to the Gaussian measure $\mu^{\mathcal I,0}_{\free}$, which is concentrated on $C(\overline{\mathcal I})$ and has finite moments of every order in $C(\overline{\mathcal I})$ by the Fernique theorem. Hence \eqref{eq:SPDEL2(I)} has a unique invariant measure given by \eqref{eq:invmeasureI}. The integration by parts formula \eqref{eq:ibpI} follows from Theorem \ref{thm:dap-grad-sist}, whose proof only uses that $\Phi\in C^1_b(L^2(\mathcal{I}))$ and therefore applies to every $\Phi\in\mathcal{F}C^{\infty}_b(L^2_{r,k})$.\\
Finally, if $h\in C^{\infty}_c(\mathcal{I})$, then $h\in{\rm Dom}(\mathcal{A}_{\mathcal{I}})$ and $\mathcal{A}_{\mathcal{I}}h=(H^{\beta}+m^2)h$ on $\mathcal{I}$; since every $x\in L^2(\mathcal{I})$ vanishes outside $\mathcal{I}$, we have $\scal{\mathcal{A}_{\mathcal{I}}h}{x}_{L^2(\mathcal{I})}=\scal{(H^{\beta}+m^2)h}{x}_{L^2}$, which gives \eqref{eq:ibpIbis}.
\end{proof}

We conclude this subsection with the following result, which is a consequence of Lemma \ref{lemma:sim-inv} and Theorem \ref{thm: rev-per-parti}.

\begin{cor}\label{cor:ibpuniqueI}
Assume that \ref{H0}, \ref{H1}, \ref{H2}, \ref{H3}  and \ref{H4} hold. Let $q\geq 1$ be the constant given by \eqref{H2:pol}. Let $\lambda$ be a Borel probability measure on $L^2(\mathcal{I})$ such that
\begin{align}
&\int_{L^2(\mathcal{I})}\norm{x}^q_{L^{2q}(\mathcal{I})}\lambda(\mathd x)<+\infty,\notag\\[3pt]
&\liminf_{\delta\rightarrow0^+}\frac1\delta\int_{L^2(\mathcal{I})}\int_{\partial_\delta\mathcal{I}}\abs{x(\xi)}\,\mathd\xi\,\lambda(\mathd x)=0,\qquad \partial_\delta\mathcal{I}:=(a,a+\delta)\cup(b-\delta,b),\label{eq:boundary-mean}
\end{align}
where $\mathcal{I}=(a,b)$, and, for every $h\in C^{\infty}_c(\mathcal{I})$ and every $\Phi\in\mathcal{F}C^{\infty}_b(L^2_{r,k})$,
\begin{equation}\label{eq:ibphp}
\int_{L^2(\mathcal{I})}\scal{h}{\delta_x\Phi(x)}_{L^2}\,\lambda(\mathd x)
=\int_{L^2(\mathcal{I})}\left[\scal{(H^{\beta}+m^2)h}{x}_{L^2}+\scal{h}{V'(x+\varphi^y)}_{L^2}\right]\Phi(x)\,\lambda(\mathd x).
\end{equation}
Then $\lambda=\nu_{\mathcal{I},y}$, where $\nu_{\mathcal{I},y}$ is the measure defined in \eqref{eq:invmeasureI}.
\end{cor}
\begin{proof}
We work in the setting of the proof of Proposition \ref{prop:invmeasureI}, with $\mathcal D={\rm Dom}(\mathcal A_{\mathcal I})$; all the scalar products are in $L^2(\mathcal I)$. By a standard approximation argument, \eqref{eq:ibphp} holds also for every $\Phi\in\xi_A(L^2(\mathcal I))$, with $\delta_x\Phi$ replaced by $\nabla\Phi$. Fix such a $\Phi$ and consider the Bochner integrals in $L^2(\mathcal I)$
\[
m_\Phi:=\int_{L^2(\mathcal I)}x\,\Phi(x)\,\lambda(\mathd x),\qquad w_\Phi:=\int_{L^2(\mathcal I)}\big[\nabla\Phi(x)-\Phi(x)V'(x+\varphi^y)\big]\lambda(\mathd x),
\]
which are well defined by the moment assumption, \eqref{H2:pol} and the boundedness of $\varphi^y$ on $\mathcal I$. Then \eqref{eq:ibphp} reads $\scal{\mathcal A_{\mathcal I}h}{m_\Phi}=\scal{h}{w_\Phi}$ for $h\in C^\infty_c(\mathcal I)$, so that, $\mathcal A_{\mathcal I}$ being self-adjoint, $u:=m_\Phi-\mathcal A_{\mathcal I}^{-1}w_\Phi$, extended by $0$ outside $\mathcal I$, satisfies $\int_\R u\,(H^\beta+m^2)h\,\mathd\xi=0$ for every $h\in C^\infty_c(\mathcal I)$. We claim that $u=0$. For $0<\varepsilon<\frac{b-a}4$, let $u_\varepsilon:=u*\eta_\varepsilon\in C^\infty_c(\R)$, with $\eta_\varepsilon:=\varepsilon^{-1}\eta(\cdot/\varepsilon)$ and $\eta\in C^\infty_c((-1,1))$ even, nonnegative, with $\int_\R\eta\,\mathd\xi=1$. Since $H^\beta$ commutes with convolutions, $(H^\beta+m^2)u_\varepsilon=0$ on $(a+\varepsilon,b-\varepsilon)$; moreover $u_\varepsilon=0$ outside $(a-\varepsilon,b+\varepsilon)$. If $\pm u_\varepsilon$ had a positive maximum at $\xi_\varepsilon\in(a+\varepsilon,b-\varepsilon)$, then, by \eqref{eq:Hbeta}, $0\leq\pm H^\beta u_\varepsilon(\xi_\varepsilon)=\mp m^2u_\varepsilon(\xi_\varepsilon)<0$. Hence $\abs{u_\varepsilon}$ attains its maximum in $J_\varepsilon:=[a-\varepsilon,a+\varepsilon]\cup[b-\varepsilon,b+\varepsilon]$ and, since $\int_B\abs{m_\Phi}\,\mathd\xi\leq\norm{\Phi}_\infty\int_{L^2(\mathcal I)}\int_B\abs{x}\,\mathd\xi\,\lambda(\mathd x)$ for every Borel set $B\subseteq\mathcal I$ (test $m_\Phi$ against ${\rm sign}(m_\Phi)\mathbb{I}_B$),
\[
\norm{u_\varepsilon}_\infty\leq\frac{\norm{\eta}_\infty}{\varepsilon}\int_{J_{2\varepsilon}}\abs{u}\,\mathd\xi\leq\frac{\norm{\eta}_\infty\norm{\Phi}_\infty}{\varepsilon}\int_{L^2(\mathcal I)}\int_{\partial_{2\varepsilon}\mathcal I}\abs{x}\,\mathd\xi\,\lambda(\mathd x)+8\norm{\eta}_\infty\sup_{J_{2\varepsilon}}\abs{\mathcal A_{\mathcal I}^{-1}w_\Phi}.
\]
Since $\mathcal A_{\mathcal I}^{-1}w_\Phi\in\mathcal D$ is continuous on $\R$ and vanishes at $a$ and $b$, by \eqref{eq:boundary-mean} the right-hand side vanishes along a sequence $\varepsilon_j\rightarrow0$; since $u_{\varepsilon_j}\rightarrow u$ in $L^2(\R)$, we get $u=0$. Hence $m_\Phi\in\mathcal D$ and $\mathcal A_{\mathcal I}m_\Phi=w_\Phi$, so that, for every $h\in\mathcal D$,
\[
\int_{L^2(\mathcal I)}\scal{\mathcal A_{\mathcal I}h}{x}\Phi(x)\lambda(\mathd x)=\scal{h}{\mathcal A_{\mathcal I}m_\Phi}=\scal{h}{w_\Phi}=\int_{L^2(\mathcal I)}\big[\scal{h}{\nabla\Phi(x)}-\scal{h}{V'(x+\varphi^y)}\Phi(x)\big]\lambda(\mathd x),
\]
which, since $R=\sqrt2\,\Id_{L^2(\mathcal I)}$ and $A^*=-\mathcal A_{\mathcal I}$, is \eqref{int-by-part} with $\mu=\lambda$. By Theorem \ref{thm: rev-per-parti}, $N_0$ is symmetric with respect to $\lambda$; by Lemma \ref{lemma:sim-inv}, $\lambda$ is invariant for the transition semigroup associated with \eqref{eq:SPDEL2(I)}, and by Proposition \ref{prop:invmeasureI} we conclude that $\lambda=\nu_{\mathcal{I},y}$.
\end{proof}

\subsection{Proof of Theorem \ref{theorem:main}}
Now we can prove Theorem \ref{theorem:main}.

\begin{proof}[Proof of Theorem \ref{theorem:main}: 2.$\Rightarrow$ 1.]$ $\\
Assume that $\mu$ satisfies the DLR equations for \eqref{eq:SPDEmain1psuedo}. Fix $\Phi\in\mathcal{F}C^{\infty}_b(L^2_{r,k})$ and $h\in C^{\infty}_c(\R)$, and let $\mathcal{I}$ be a bounded open interval of $\R$ such that $\supp(h)\subset\mathcal{I}$. By Corollary \ref{remark:DLRtranslated} and \eqref{eq:invmeasureI}, for every $y\in L^2_{r,k}(\mathcal{I}^c)\cap C(\mathcal{I}^c)$ we have
\begin{align}\label{eq:condvsinv}
\int_{E^{\mathcal{I},y}}\Phi(x)\mu^{\mathcal{I},y}(\mathd x)&=\int_{E^{\mathcal{I},0}}\Phi(x+\varphi^y)\frac{\exp\left(-\int_{\mathcal{I}}V(x(\xi)+\varphi^y(\xi))\mathd \xi\right)}{\mathcal{Z}^{\mathcal{I}}(y)}\mu_{\free}^{\mathcal{I},0}(\mathd x)\notag\\
&=\int_{L^2(\mathcal{I})}\Phi(x+\varphi^y)\nu_{\mathcal{I},y}(\mathd x),
\end{align}
where we identified $E^{\mathcal{I},0}$ with $L^2(\mathcal{I})$, $\mu^{\mathcal{I},y}$ is the conditional probability of $\mu$ given by Lemma \ref{lemma:disintegration} and $\nu_{\mathcal{I},y}$ is the unique invariant measure of \eqref{eq:SPDEL2(I)} given by Proposition \ref{prop:invmeasureI}.

Since $\varphi^y=y$ on $\mathcal{I}^c$ and $\varphi^y$ is bounded on the bounded interval $\mathcal{I}$, we have $\varphi^y\in L^2_{r,k}$; hence $\Phi(\cdot+\varphi^y)\in\mathcal{F}C^{\infty}_b(L^2_{r,k})$ and $\delta_x\big[\Phi(\cdot+\varphi^y)\big](x)=(\delta_x\Phi)(x+\varphi^y)$, so that we can apply \eqref{eq:ibpIbis} to the function $\Phi(\cdot+\varphi^y)$ and we obtain
\begin{align}\label{eq:intxpartiprov}
\int_{L^2(\mathcal{I})}&\scal{h}{(\delta_x\Phi)(x+\varphi^y)}_{L^2}\nu_{\mathcal{I},y}(\mathd x)\notag\\
&=\int_{L^2(\mathcal{I})}\left[\scal{(H^{\beta}+m^2)h}{x}_{L^2}+\scal{h}{V'(x+\varphi^y)}_{L^2}\right]\Phi(x+\varphi^y)\nu_{\mathcal{I},y}(\mathd x).
\end{align}
Moreover, $\varphi^y\in L^1_{r,k}(\R)$ by Theorem \ref{thm:harmonic-continuity}, so that by Lemma \ref{lemma:testdecay} we have $\varphi^y\,(H^{\beta}+m^2)h\in L^1(\R)$; since $\supp(h)\subset\mathcal{I}$ and $\varphi^y$ is a weak harmonic extension of $y$, point 3. of Definition \ref{def:harmonicextension} gives
\begin{equation}\label{eq:noresidue}
\scal{(H^{\beta}+m^2)h}{\varphi^y}_{L^2}=0,
\end{equation}
so that $\scal{(H^{\beta}+m^2)h}{x}_{L^2}=\scal{(H^{\beta}+m^2)h}{x+\varphi^y}_{L^2}$. Formula \eqref{eq:condvsinv} extends, by monotone convergence applied to positive and negative parts, to every Borel $\Phi$ such that $\int\abs{\Phi(x+\varphi^y)}\nu_{\mathcal{I},y}(\mathd x)<+\infty$, and both integrands in \eqref{eq:intxpartiprov} have this property, since $\nu_{\mathcal{I},y}$ has a bounded density with respect to $\mu^{\mathcal{I},0}_{\free}$, which has finite moments of every order in $C(\overline{\mathcal I})$. Hence, applying \eqref{eq:condvsinv} to both sides of \eqref{eq:intxpartiprov}, for every $y\in L^2_{r,k}(\mathcal{I}^c)\cap C(\mathcal{I}^c)$ 
\begin{equation}\label{eq:ibpcond}
\int_{E^{\mathcal{I},y}}\scal{h}{\delta_x\Phi(x)}_{L^2}\mu^{\mathcal{I},y}(\mathd x)
=\int_{E^{\mathcal{I},y}}\left[\scal{(H^{\beta}+m^2)h}{x}_{L^2}+\scal{h}{V'\circ x}_{L^2}\right]\Phi(x)\mu^{\mathcal{I},y}(\mathd x).
\end{equation}
Finally, integrating \eqref{eq:ibpcond} with respect to $\mu^{\mathcal{I}^c}(\mathd y)$, which is allowed since both sides of \eqref{eq:IPF} are $\mu$-integrable by the moment assumption of Theorem \ref{theorem:main}, and using \eqref{eq:formula-prob-cond}, we conclude
\[
\int_{L^{2q}_{r,k}}\scal{h}{\delta_x\Phi(x)}_{L^2}\mu(\mathd x)
=\int_{L^{2q}_{r,k}}\left[\scal{(H^{\beta}+m^2)h}{x}_{L^2}+\scal{h}{V'\circ x}_{L^2}\right]\Phi(x)\mu(\mathd x).
\]
\end{proof}

\begin{proof}[Proof of Theorem \ref{theorem:main}: 1.$\Rightarrow$ 2.]$ $\\
Assume that $\mu$ satisfies the integration by parts formula \eqref{eq:IPF} associated with \eqref{eq:SPDEmain1psuedo} and fix a bounded open interval $\mathcal{I}$ of $\R$.

We first recall two properties which hold for $\mu^{\mathcal{I}^c}$-almost every $y\in L^2_{r,k}(\mathcal{I}^c)$. Since $R^{-1}_{\mathcal{I}^c}\big(C(\mathcal{I}^c)\big)\supseteq C(\R)$ and $\mu(C(\R))=1$, by \eqref{eq:legge-Ic} we have
\begin{equation}\label{eq:ycontinuous}
\mu^{\mathcal{I}^c}\big(C(\mathcal{I}^c)\big)=1,
\end{equation}
so that, by Theorem \ref{thm:harmonic-continuity}, the harmonic extension $\varphi^y$ is well defined for $\mu^{\mathcal{I}^c}$-almost every $y$.
Moreover, applying \eqref{eq:formula-prob-cond} to $x\mapsto\norm{x}^q_{L^{2q}_{r,k}}$, we deduce
\begin{equation}\label{eq:momentcond}
\int_{E^{\mathcal{I},y}}\norm{x}^q_{L^{2q}_{r,k}}\,\mu^{\mathcal{I},y}(\mathd x)<+\infty
\qquad\text{for }\mu^{\mathcal{I}^c}\text{-almost every }y .
\end{equation}
Since $\mathcal{I}$ is bounded we have $c_{\mathcal{I}}:=\inf_{\mathcal{I}}\rho_{r,k}>0$, so that $\norm{x}_{L^{2q}(\mathcal{I})}\leq c_{\mathcal{I}}^{-1/(2q)}\norm{x}_{L^{2q}_{r,k}}$ for every $x\in L^{2q}_{r,k}$. 
Let $h\in C^{\infty}_c(\mathcal{I})$, let $\Phi\in\mathcal{F}C^{\infty}_b(L^2_{r,k})$ and let $\Psi\in\mathcal{F}C^{\infty}_b(L^2_{r,k})$ be of the form 
\begin{equation}\label{eq:psi}
 \Psi(x)=\psi(\scal{x}{g_1}_{L^2},\dots,\scal{x}{g_\ell}_{L^2})
\end{equation}
 where $l\in\N$, $g_1,\dots,g_\ell\in C^{\infty}_c(\mathcal{I}^c)$ and $\psi\in C_b^{\infty}(\R^\ell)$. Since $\supp(h)\subset\mathcal{I}$ we have $\scal{h}{\delta_x\Psi(x)}_{L^2}=0$, hence $\scal{h}{\delta_x(\Phi\Psi)(x)}_{L^2}=\Psi(x)\scal{h}{\delta_x\Phi(x)}_{L^2}$. Applying \eqref{eq:IPF} to $\Phi\Psi$, disintegrating $\mu$ by \eqref{eq:formula-prob-cond} and noting that $\Psi\equiv\Psi(y)$ on $E^{\mathcal{I},y}$, we get
\[
\int_{L^2_{r,k}(\mathcal{I}^c)}\Psi(y)\big[G_\Phi(y)-\tilde G_\Phi(y)\big]\mu^{\mathcal{I}^c}(\mathd y)=0,
\]
where
\begin{align}
G_\Phi(y)&:=\int_{E^{\mathcal{I},y}}\scal{h}{\delta_x\Phi(x)}_{L^2}\mu^{\mathcal{I},y}(\mathd x),\\
\tilde G_\Phi(y)&:=\int_{E^{\mathcal{I},y}}\left[\scal{(H^{\beta}+m^2)h}{x}_{L^2}+\scal{h}{V'\circ x}_{L^2}\right]\Phi(x)\mu^{\mathcal{I},y}(\mathd x).
\end{align}
Since the functions $\Psi$ of the form \eqref{eq:psi} generate $\mathcal{B}(L^2_{r,k}(\mathcal{I}^c))$, a monotone class argument gives $G_\Phi=\tilde G_\Phi$ $\mu^{\mathcal{I}^c}$-almost everywhere. Hence, for every $h\in C^\infty_c(\mathcal I)$ and $\Phi\in\mathcal{F}C^{\infty}_b(L^2_{r,k})$ there exists a $\mu^{\mathcal I^c}$-null set $N_{h,\Phi}$ such that for every $y\notin N_{h,\Phi}$ we have
\begin{equation}\label{eq:ibplocalizedv}
\int_{E^{\mathcal{I},y}}\scal{h}{\delta_x\Phi(x)}_{L^2}\mu^{\mathcal{I},y}(\mathd x)
=\int_{E^{\mathcal{I},y}}\left[\scal{(H^{\beta}+m^2)h}{x}_{L^2}+\scal{h}{V'\circ x}_{L^2}\right]\Phi(x)\mu^{\mathcal{I},y}(\mathd x),
\end{equation}
Since below we apply \eqref{eq:ibplocalizedv} to functions depending on $y$, we need a null set independent of $h$ and $\Phi$. It is enough to take $h$ and $\Phi$ in a countable family. We take $h$ in a countable subset of $C^\infty_c(\mathcal I)$, and $\Phi(x)=\phi(\scal{x}{f_1}_{L^2},\ldots,\scal{x}{f_n}_{L^2})$ with $f_1,\ldots,f_n$ in a countable subset of $C^\infty_c(\R)$ and $\phi$ in a countable subset of $C^\infty_b(\R^n)$. Let $N$ be the union of the corresponding sets $N_{h,\Phi}$  and of the null set where \eqref{eq:momentcond} fails; $N$ is $\mu^{\mathcal I^c}$-null. By standard approximation arguments for every $h\in C^\infty_c(\mathcal I)$, $\Phi\in\mathcal{F}C^{\infty}_b(L^2_{r,k})$ and $y\notin N$ we have
\begin{equation}\label{eq:ibplocalized}
\int_{E^{\mathcal{I},y}}\scal{h}{\delta_x\Phi(x)}_{L^2}\mu^{\mathcal{I},y}(\mathd x)
=\int_{E^{\mathcal{I},y}}\left[\scal{(H^{\beta}+m^2)h}{x}_{L^2}+\scal{h}{V'\circ x}_{L^2}\right]\Phi(x)\mu^{\mathcal{I},y}(\mathd x).
\end{equation}

Let $\varphi^y$ be the harmonic extension of $y$ given by Theorem \ref{thm:harmonic-continuity}, which, as observed above, belongs to $L^2_{r,k}$, and let
\[
\lambda_y:=\mu^{\mathcal{I},y}\circ \mathfrak{J}^{-1}_{y}
\]
where $\mathfrak{J}_y(x)=x-\varphi^y$. Hence for every $\Psi\in B_b(L^2_{r,k})$ we have
\begin{equation}\label{eq:immagine}
\int_{E^{\mathcal{I},y}}\Psi(x)\mu^{\mathcal{I},y}(\mathd x)=\int_{L^2(\mathcal{I})}\Psi(x+\varphi^y)\lambda_y(\mathd x).
\end{equation}
Applying \eqref{eq:ibplocalized} to $\Phi(\cdot-\varphi^y)\in\mathcal{F}C^{\infty}_b(L^2_{r,k})$, and using \eqref{eq:immagine} together with \eqref{eq:noresidue}, which is available by \eqref{eq:ycontinuous}, we obtain
\begin{align*}
\int_{L^2(\mathcal{I})}&\scal{h}{\delta_x\Phi(x)}_{L^2}\lambda_y(\mathd x)\\
&=\int_{L^2(\mathcal{I})}\left[\scal{(H^{\beta}+m^2)h}{x}_{L^2}+\scal{h}{V'(x+\varphi^y)}_{L^2}\right]\Phi(x)\lambda_y(\mathd x),
\end{align*}
namely $\lambda_y$ satisfies \eqref{eq:ibphp}. Moreover, since $\varphi^y$ is bounded on the bounded interval $\mathcal{I}$, by \eqref{eq:momentcond} and \eqref{eq:immagine} we get
\[
\int_{L^2(\mathcal{I})}\norm{x}^q_{L^{2q}(\mathcal{I})}\,\lambda_y(\mathd x)
=\int_{E^{\mathcal{I},y}}\norm{x-\varphi^y}^q_{L^{2q}(\mathcal{I})}\mu^{\mathcal{I},y}(\mathd x)<+\infty .
\]

We now prove that $\lambda_y$ satisfies \eqref{eq:boundary-mean} for $\mu^{\mathcal I^c}$-almost every $y$. Let $\mathcal I=(a,b)$. For $\mu^{\mathcal I^c}$-almost every $y$, $\mu^{\mathcal I,y}$-almost every $x$ is continuous on $\R$ and $x=\varphi^y$ on $\mathcal I^c$ (see Remark \ref{remark:disintegrationsubset} and \eqref{eq:prob-cond}), so that $x(a)=\varphi^y(a)$ and $x(b)=\varphi^y(b)$; hence, since $\varphi^y$ is continuous, by \eqref{eq:immagine} $\lambda_y$ satisfies \eqref{eq:boundary-mean} whenever
\begin{equation}\label{eq:Bcond}
\liminf_{\delta\rightarrow0^+}\frac1\delta\int_{E^{\mathcal I,y}}\int_0^\delta\big(\abs{x(a+s)-x(a)}+\abs{x(b-s)-x(b)}\big)\mathd s\,\mu^{\mathcal I,y}(\mathd x)=0.
\end{equation}
We first consider the intervals $(-T,T)$: for $x\in C(\R)$, $N\geq2$ and $0<\delta\leq1$, by the Fubini theorem
\begin{multline*}
\frac1\delta\int_1^N\!\!\int_0^\delta\big(\abs{x(-T+s)-x(-T)}+\abs{x(T-s)-x(T)}\big)\mathd s\,\mathd T\\
\leq\sup_{0<s\leq\delta}\big(\norm{x(\cdot+s)-x}_{L^1(-N,0)}+\norm{x(\cdot-s)-x}_{L^1(0,N)}\big),
\end{multline*}
which vanishes as $\delta\rightarrow0$ and is bounded by $C_N\norm{x}_{L^{2q}_{r,k}}$. Hence, by \eqref{eq:momenti-DLR}, the dominated convergence theorem, \eqref{eq:formula-prob-cond} and the Fatou lemma (in $T$ and then in $y$), for almost every $T>1$ \eqref{eq:Bcond} holds with $\mathcal I=(-T,T)$ for $\mu^{(-T,T)^c}$-almost every $y$, and so, by the argument below, \eqref{eq:tras-mu} holds with $\mathcal I=(-T,T)$. For a general $\mathcal I$, we fix such a $T$ with $\overline{\mathcal I}\subseteq(-T,T)$: by \eqref{eq:tras-mu} with $(-T,T)$, \ref{H4} and the Fernique theorem, $\sup_{\overline{\mathcal I}}\abs{x}$ is integrable with respect to $\mu^{(-T,T),y'}$ for $\mu^{(-T,T)^c}$-almost every $y'$, hence, by the tower property of conditional probabilities, with respect to $\mu^{\mathcal I,y}$ for $\mu^{\mathcal I^c}$-almost every $y$, and \eqref{eq:Bcond} follows by dominated convergence, the integrand being bounded by $4\sup_{\overline{\mathcal I}}\abs{x}$.

Hence, by Corollary \ref{cor:ibpuniqueI}, $\lambda_y=\nu_{\mathcal{I},y}$, where $\nu_{\mathcal{I},y}$ is the unique invariant measure of \eqref{eq:SPDEL2(I)} given by Proposition \ref{prop:invmeasureI}. Therefore, by \eqref{eq:immagine} and \eqref{eq:invmeasureI}, for $\mu^{\mathcal{I}^c}$-almost every $y$ we conclude
\begin{align*}
\int_{E^{\mathcal{I},y}}\Phi(x)\mu^{\mathcal{I},y}(\mathd x)&=\int_{L^2(\mathcal{I})}\Phi(x+\varphi^y)\nu_{\mathcal{I},y}(\mathd x)\\
&=\int_{E^{\mathcal{I},0}}\Phi(x+\varphi^y)\frac{\exp\left(-\int_{\mathcal{I}}V(x(\xi)+\varphi^y(\xi))\mathd \xi\right)}{\mathcal{Z}^{\mathcal{I}}(y)}\mu_{\free}^{\mathcal{I},0}(\mathd x),
\end{align*}
and so that, by Corollary \ref{remark:DLRtranslated}, $\mu$ verifies the DLR equations associated with \eqref{eq:SPDEmain1psuedo}.
\end{proof}
\part{Important related results}

\section{A Fernique-type theorem for solution to DLR equation}\label{section:Fernique}

Throughout this section we assume \ref{H0}, \ref{H1}, \ref{H2}, \ref{H3} and \ref{H4} hold. Let $q\geq1$ be the constant given by \eqref{H2:pol}.
Moreover, we require the following coercivity assumption on the potential, which is not a consequence of \ref{H2}.
\begin{enumerate}[label={{(H5)}}]
\item\label{H5}\emph{There exist $\zeta_0>0$ and $\zeta_1\geq0$ such that
\begin{equation}\label{H5:coerc}
\xi V'(\xi)\geq \zeta_0\abs{\xi}^{q+1}-\zeta_1,\qquad \xi\in\R .
\end{equation}}
\end{enumerate}

\begin{rmk}\label{rmk:H5}
Hypothesis \ref{H5} only concerns the behavior of $V$ at infinity, so it does not exclude potentials which are not convex on bounded sets, and in particular it does not exclude non-uniqueness. For instance, if $f(x)=(k^2-x^2)^2$, the potential $V^{\gamma_2}(y)=\gamma_2f(y)-\frac{m^2}{2}\abs{y}^2$ considered in Section \ref{section:nonuniqueness} satisfies \ref{H2} with $q=3$ and \ref{H5}, since
\[
\xi\,(V^{\gamma_2})'(\xi)=4\gamma_2\xi^4-\big(4\gamma_2k^2+m^2\big)\xi^2\geq2\gamma_2\xi^4-\frac{(4\gamma_2k^2+m^2)^2}{8\gamma_2},\qquad\xi\in\R .
\]
\end{rmk}

The next lemma collects some properties of $V$ that are consequences of \ref{H2} and \ref{H5}.

\begin{lemm}\label{lemma:potentialbounds}
Assume \ref{H2} and \ref{H5} hold. Then there exist $N_1>0$ and $N_2\geq0$ such that, for every $a,b\in\R$,
\begin{align}
aV'(a+b)&\geq N_1\abs{a}^{q+1}-N_2\big(1+\abs{b}^{q+1}\big),\label{eq:potential-coerc}\\
\abs{V'(a)}^2&\leq N_2\big(1+\abs{a}^{2q}\big).\label{eq:potential-growth}
\end{align}
\end{lemm}

\begin{proof}
Estimate \eqref{eq:potential-growth} is \eqref{H2:pol}. As for \eqref{eq:potential-coerc}, applying \eqref{H5:coerc} at the point $a+b$ and using \eqref{H2:pol} we get
\[
aV'(a+b)=(a+b)V'(a+b)-bV'(a+b)\geq \zeta_0\abs{a+b}^{q+1}-\zeta_1-C\abs{b}\big(1+\abs{a+b}^{q}\big).
\]
By the Young inequality there exists $C'>0$ such that $C\abs{b}\abs{a+b}^{q}\leq\frac{\zeta_0}{2}\abs{a+b}^{q+1}+C'\abs{b}^{q+1}$, whence
\[
aV'(a+b)\geq\frac{\zeta_0}{2}\abs{a+b}^{q+1}-C''\big(1+\abs{b}^{q+1}\big),
\]
and \eqref{eq:potential-coerc} follows from $\abs{a}^{q+1}\leq2^{q}\big(\abs{a+b}^{q+1}+\abs{b}^{q+1}\big)$.
\end{proof}

\subsection{A priori estimates for equation \eqref{eq:SPDEmain1psuedo}}\label{section:apriori}

Throughout this subsection we fix $1<r<1+2s_2$ and $k\in(0,k_0)$, according to \ref{H3}. 

Let $\{u(t)\}_{t\geq 0}$ be a predictable process such that $\mathbb{E}\big[\norm{u}_{L^2([0,T];L^2)}^2\big]<+\infty$ for every $T>0$. We consider the following equations in $L^2_{r,k}$
\begin{align}
&\begin{cases}
\mathd X^{u}(t)+\mathcal{A}X^{u}(t)\mathd t=\big[\mathcal{B}(X^{u}(t))+\sqrt{2}u(t)\big]\mathd t+\sqrt{2\rho_{r,k}}\mathd W_{r,k}(t),\qquad t>0,\\
X^{u}(0)=x,
\end{cases}\label{eq:systemZY1}\\[4pt]
&\begin{cases}
\mathd Z(t)+\mathcal{A}Z(t)\mathd t=\sqrt{2\rho_{r,k}}\mathd W_{r,k}(t),\qquad t>0,\\
Z(0)=z,
\end{cases}\label{eq:systemZY2}\\[4pt]
&\begin{cases}
\mathd Y^{u}(t)+\mathcal{A}Y^{u}(t)\mathd t=\big[\mathcal{B}(Y^{u}(t)+Z(t))+\sqrt{2}u(t)\big]\mathd t,\qquad t>0,\\
Y^{u}(0)=y,
\end{cases}\label{eq:systemZY3}
\end{align}
where $\mathcal{A}$, $\mathcal{B}$ and $\{W_{r,k}(t)\}_{t\geq 0}$ are the same coefficients as in \eqref{eq:SPDEmain1Abstract}. By Lemma \ref{lemma:stochconv}, for every $z\in L^2_{r,k}$, the equation \eqref{eq:systemZY2} has a unique mild solution $\{Z(t,z)\}_{t\geq 0}$ given by
\[
Z(t,z)=e^{-t\mathcal{A}}z+W_{\mathcal{A}}(t),\qquad t\geq 0,
\]
where $\{W_{\mathcal A}(t)\}_{t\geq 0}$ is the stochastic convolution process defined by \eqref{eq:stochconv}. Moreover, for every $x\in L^2_{r,k}$ equation \eqref{eq:systemZY1} admits a unique generalized mild solution $\{X^{u}(t,x)\}_{t\geq 0}$ in the sense of Proposition \ref{prop:wellposedness}. Indeed, \eqref{eq:systemZY1} is \eqref{eq:SPDEmain1Abstract} with the stochastic convolution $W_{\mathcal{A}}$ replaced by
\[
W^{u}_{\mathcal{A}}(t):=W_{\mathcal{A}}(t)+\sqrt2\int_0^te^{-(t-s)\mathcal{A}}u(s)\mathd s,\qquad t\geq0,
\]
and noting that, for every $p\geq 2$, the process $\{W^{u}_{\mathcal{A}}(t)\}_{t\geq 0}$ has continuous trajectories with values in $L^p_{r,k}$, the proof of Proposition \ref{prop:wellposedness} (see Theorem \ref{th:solmild}) also applies in this case. For every $y\in L^2_{r,k}$ also \eqref{eq:systemZY3} admits a unique generalized mild solution $\{Y^{u}(t,y)\}_{t\geq 0}$ and in particular if $x=y+z$ then
\[
Y^{u}(t,y)=X^{u}(t,x)-Z(t,z),\qquad t\geq0 .
\]
Following Remark \ref{rmk:weighttranslation}, for $\tau\in\R$ we write $\rho^\tau_{r,k}:=\rho_{r,k}(\cdot-\tau)$, and for $s\geq0$ we use the notation 
\[
H^{s}_{r,k}:=H^{s}_{\rho_{r,k}},
\] 
where $H^{s}_{\rho_{r,k}}$ is the weighted Sobolev space defined by \eqref{eq:weightedsobolev}. Recalling the weighted H\"older spaces introduced in Remark \ref{rmk:stochconvholder}, we shall use the following standard embedding, if $j\geq1$, $\gamma\geq0$, $\varepsilon>0$ and $r-j\varepsilon>1$, then
\begin{equation}\label{eq:weighted-holder-Lr}
\norm{f}_{L^j_{\rho^\tau_{r,k}}}^j\leq\norm{f}_{C^\gamma_{\rho^\tau_{\varepsilon,k}}}^j\int_\R\rho^\tau_{r,k}\big(\rho^\tau_{\varepsilon,k}\big)^{-j}\mathd\xi=N_3k^{-1/2}\norm{f}_{C^\gamma_{\rho^\tau_{\varepsilon,k}}}^j ,
\end{equation}
where we used that $\rho^\tau_{r,k}(\rho^\tau_{\varepsilon,k})^{-j}=\rho^\tau_{r-j\varepsilon,k}$ and \eqref{eq:rhoL1}, and  $N_3>0$. Moreover, since $\rho_{r,k}^{1/2}=\rho_{r/2,k}$, the Sobolev embedding $H^{s_1}(\R)\hookrightarrow C^{0,\gamma}(\R)$ gives, for every $\gamma\in\big(0,s_1-\frac12\big)$, a constant $N_4>0$, independent of $\tau$, such that
\begin{equation}\label{eq:sobolevweighted}
\norm{f}_{C^{\gamma}_{\rho^\tau_{r/2,k}}}\leq N_4\norm{f}_{H^{s_1}_{\rho^\tau_{r,k}}} .
\end{equation}

Before proving the a priori estimates, we need an auxiliary lemma ensuring that all the computations carried out below are justified.

\begin{lemm}\label{lemm:apriori}
Assume \ref{H0}, \ref{H1}, \ref{H2}, \ref{H3}, \ref{H4} and \ref{H5} hold. Then, for every $y\in H^{s_1}_{r,k}\cap L^{2q}_{r,k}$, $z\in L^{2q}_{r,k}$ and $T>0$, we have
\[
Y^{u}(\cdot,y)\in L^\infty\big([0,T];H^{s_1}_{r,k}\big)\cap L^2\big([0,T];H^{2s_1}_{r,k}\big),\qquad \mathbb{P}\text{-a.s.}
\]
\end{lemm}

\begin{proof}
For simplicity, we set $\rho:=\rho_{r,k}$ and $Y:=Y^{u}(\cdot,y)$. Let $y\in H^{s_1}_{r,k}\cap L^{2q}_{r,k}$, $z\in L^{2q}_{r,k}$ and $T>0$. We are going to work pathwise, so we fix a realization of $\{Y(t)\}_{t\geq0}$, $\{Z(t)\}_{t\geq0}$ and $\{u(t)\}_{t\geq 0}$.

Since  $y+z\in L^{2q}_{r,k}$, by Proposition \ref{prop:wellposedness} and estimate \eqref{eq:potential-growth} we have
\[
\norm{\mathcal{B}(Y(t)+Z(t))}^2_{L^2_{r,k}}\leq N_2\Big(\norm{\rho}_{L^1(\R)}+\norm{Y(t)+Z(t)}^{2q}_{L^{2q}_{r,k}}\Big),\qquad t\in[0,T].
\]
In addition $\norm{u(t)}_{L^2_{r,k}}\leq\norm{u(t)}_{L^2}$ since $\rho$ takes values in $(0,1]$. Hence, the function $g:=\mathcal{B}(Y+Z)+\sqrt2u$ belongs to $L^2([0,T];L^2_{r,k})$, and so $\rho^{1/2}g\in L^2([0,T];L^2)$. Hence $\{Y(t)\}_{t\in [0,T]}$ is the unique mild solution in $L^2_{r,k}$ of the linear problem
\begin{equation}\label{eq:apriori-linearY}
\partial_tY(t)+\mathcal{A}Y(t)=g(t),\qquad Y(0)=y.
\end{equation}
To prove the statement, we need a regularity estimate for  $H^{\beta}+m^2$ in $L^2$. We recall that the closure of $H^{\beta}+m^2$ in $L^2$ is self-adjoint, with domain contained in $H^{2s_1}(\R)$, and that $-(H^{\beta}+m^2)$ is dissipative; moreover, by \eqref{eq:energy-is-Hs1}, ${\rm Dom}\big((H^{\beta}+m^2)^{1/2}\big)=H^{s_1}(\R)$ with equivalent norms. Hence $-(H^{\beta}+m^2)$ generates an analytic semigroup of contractions on $L^2$ (see \cite[Corollary 4.8]{Eng-Nag2006}), so by the maximal $L^2$-regularity of analytic semigroups (see \cite{Des1964}) there exists a constant $S_{\rm mr}>0$ such that, for every $w_0\in H^{s_1}(\R)$ and $f\in L^2([0,T];L^2)$, the mild solution $w$ of $\partial_tw+(H^{\beta}+m^2)w=f$ with $w(0)=w_0$ belongs to $C([0,T];H^{s_1}(\R))\cap L^2([0,T];H^{2s_1}(\R))$ and
\begin{equation}\label{eq:apriori-maxreg}
\sup_{t\in[0,T]}\norm{w(t)}_{H^{s_1}(\R)}+\norm{w}_{L^2([0,T];H^{2s_1}(\R))}\leq S_{\rm mr}\Big(\norm{w_0}_{H^{s_1}(\R)}+\norm{f}_{L^2([0,T];L^2)}\Big).
\end{equation}
 By \eqref{eq:commutatorestimate} when $\beta_0=0$, by \eqref{eq:commutatorfull} when $\beta_0>0$, by \eqref{eq:Gtheta-Htheta} and by density, there exists a constant $S_{\mathcal K}>0$ such that
\begin{equation}\label{eq:apriori-commutator}
\norm{\mathcal{K}w}_{L^2}\leq S_{\mathcal K}\norm{w}_{H^{s_1}(\R)},\qquad w\in H^{s_1}(\R),
\end{equation}
where  $\mathcal{K}w:=\mathcal{C}_{r,k}\big(\rho^{-1/2}w\big)$ and $\mathcal{C}_{r,k}$ is the commutator given by \eqref{eq:commutatordef}. We consider the following equation
\begin{equation}\label{eq:aux}
\partial_tv(t)+(H^{\beta}+m^2)v(t)=\mathcal{K}v(t)+\rho^{1/2}g,\qquad v(0)=\rho^{1/2}y\in H^{s_1}(\R).
\end{equation}
Let $T_0\in(0,T]$ be such that $S_{\rm mr}S_{\mathcal K}T_0^{1/2}<1$. By \eqref{eq:apriori-maxreg}, \eqref{eq:apriori-commutator} and standard contraction arguments on the mild formulation, \eqref{eq:aux} has a unique mild solution $v$ in $C([0,T_0];H^{s_1}(\R))$, in particular $v\in C([0,T_0];H^{s_1}(\R))\cap L^2([0,T_0];H^{2s_1}(\R))$.

Set $\tilde Y:=\rho^{-1/2}v$, so that $\tilde Y\in L^2([0,T_0];H^{2s_1}_{r,k})$. By \eqref{eq:commutatordef} we have $\rho^{1/2}(H^{\beta}+m^2)\tilde Y=(H^{\beta}+m^2)v-\mathcal{K}v$, hence $\partial_t\tilde Y+(H^{\beta}+m^2)\tilde Y=g$ on $[0,T_0]$. Hence, by \eqref{eq:apriori-commutator} and Remark \ref{rmk:domainA2}, $\tilde Y$ is a strong solution to \eqref{eq:apriori-linearY} on $[0,T_0]$, and therefore also mild. By uniqueness $\tilde Y=Y$ on $[0,T_0]$. As $T_0$ depends only on $S_{\rm mr}$ and $S_{\mathcal K}$, iterating the argument on $[T_0,2T_0],[2T_0,3T_0],\ldots$ with initial data $Y(T_0),Y(2T_0),\ldots\in H^{s_1}_{r,k}$ gives the assertion on the whole of $[0,T]$.
\end{proof}

Now we can prove the a priori estimates that will be fundamental in the main results of the next sections.

\begin{thm}\label{theorem:apriori}
Assume \ref{H0}, \ref{H1}, \ref{H2}, \ref{H3}, \ref{H4} and \ref{H5} hold. Let $q\geq 1$ be the constant given by \eqref{H2:pol} and let
\[
\gamma\in\big(0,s_1-\frac12\big),\qquad 0<\varepsilon<\frac{r-1}{2q}.
\]
Then there exist $\ell_0>0$ and a constant $K_0>0$ such that, for every $\tau\in\R$, $\ell\in[0,\ell_0)$, $z\in C^{\gamma}_{\rho^\tau_{\varepsilon,k}}$ and $t\geq0$, the following estimates hold $\mathbb{P}$-a.s.
\begin{enumerate}[label=(\roman*)]
\item if $y\in L^2_{r,k}$, then
\begin{equation}\label{eq:aprioriL2}
\norm{Y^{u}(t,y)}^2_{L^2_{\rho^\tau_{r,k}}}+\int_0^te^{-\ell(t-s)}\Big(\norm{Y^{u}(s,y)}^2_{H^{s_1}_{\rho^\tau_{r,k}}}+\norm{Y^{u}(s,y)}^{q+1}_{L^{q+1}_{\rho^\tau_{r,k}}}\Big)\mathd s\leq K_0\Big(e^{-\ell t}\norm{y}^2_{L^2_{\rho^\tau_{r,k}}}+R_\ell(t)\Big),
\end{equation}
\item if $y\in L^{q+1}_{r,k}$, then
\begin{equation}\label{eq:aprioriLp}
\norm{Y^{u}(t,y)}^{q+1}_{L^{q+1}_{\rho^\tau_{r,k}}}+\int_0^te^{-\ell(t-s)}\norm{Y^{u}(s,y)}^{2q}_{L^{2q}_{\rho^\tau_{r,k}}}\mathd s\leq K_0\Big(e^{-\ell t}\norm{y}^{q+1}_{L^{q+1}_{\rho^\tau_{r,k}}}+R_\ell(t)\Big),
\end{equation}
\item if $y\in H^{s_1}_{r,k}\cap L^{q+1}_{r,k}$, then
\begin{equation}\label{eq:aprioriHs}
\norm{Y^{u}(t,y)}^2_{H^{s_1}_{\rho^\tau_{r,k}}}\leq K_0\Big(e^{-\ell t}\Big(\norm{y}^2_{H^{s_1}_{\rho^\tau_{r,k}}}+\norm{y}^{q+1}_{L^{q+1}_{\rho^\tau_{r,k}}}\Big)+R_\ell(t)\Big),
\end{equation}
\end{enumerate}
where $\rho^{\tau}_{r,k}:=\rho_{r,k}(\cdot-\tau)$ and
\begin{equation}\label{eq:R_l}
R_\ell(t):=\int_0^te^{-\ell(t-s)}\Big(1+\norm{Z(s,z)}^{2q}_{C^\gamma_{\rho^\tau_{\varepsilon,k}}}\Big)\mathd s+\norm{u}^2_{L^2([0,t];L^2)}.
\end{equation}
\end{thm}

\begin{proof}
For simplicity, we set $\rho:=\rho_{r,k}$, $Y:=Y^{u}(\cdot,y)$ and $Z:=Z(\cdot,z)$. By Remark \ref{rmk:weighttranslation}, it suffices to consider $\tau=0$. Moreover, by Proposition \ref{prop:wellposedness}, the Fatou lemma, and the weak lower semicontinuity of the norms, it suffices to prove (i), (ii) and (iii) for $y\in H^{s_1}_{r,k}\cap L^{2q}_{r,k}$. Let $y\in H^{s_1}_{r,k}\cap L^{2q}_{r,k}$ and $z\in L^{2q}_{r,k}$. We are going to work pathwise, so we fix a realization of  $Y:=Y^{u}(\cdot,y)$, $\{Z(t)\}_{t\geq0}$ and $\{u(t)\}_{t\geq0}$.

\noindent By Lemma \ref{lemma:coercivity-H3} and \eqref{eq:E-vs-A} (extended by density to $Y(t)\in H^{2s_1}_{r,k}$, see Remark \ref{rmk:domainA2}), setting $c:=C_4$ we have
\begin{equation}\label{eq:apriori-coercive}
\scal{\mathcal{A}Y}{Y}_{L^2_{r,k}}=\int_\R\big[(H^{\beta}+m^2)Y\big](\xi)Y(\xi)\rho(\xi)\mathd\xi=\mathcal E(Y,\rho Y)\geq c\norm{Y}^2_{H^{s_1}_{r,k}} .
\end{equation}
We set $\ell_0:=\min\big\{\frac c2,\ m^2-\frac{\omega_{r,k}}{2}\big\}$.

\noindent By \eqref{eq:weighted-holder-Lr}, for every $j\in\{q+1,2q\}$ we have 
\[
\norm{Z(t)}^{j}_{L^{j}_{r,k}}\leq N_{3}k^{-1/2}\big(1+\norm{Z(t)}^{2q}_{C^\gamma_{\rho_{\varepsilon,k}}}\big).
\]
Moreover, by Lemma \ref{lemm:apriori} and Remark \ref{rmk:domainA2}, $Y\in C([0,T];H^{s_1}_{r,k})\cap L^2([0,T];H^{2s_1}_{r,k})$ is a strong solution of \eqref{eq:systemZY3}, which justifies the computations below. Now we can prove the three points of the statement.\\

\noindent (i) Multiplying in $L^2_{r,k}$ both sides of \eqref{eq:systemZY3}  by $Y$, we get
\[
\frac12\frac{\mathd}{\mathd t}\norm{Y}^2_{L^2_{r,k}}+\scal{\mathcal{A}Y}{Y}_{L^2_{r,k}}-\scal{\mathcal{B}(X)}{Y}_{L^2_{r,k}}=\sqrt2\scal{u}{Y}_{L^2_{r,k}}.
\]
By \eqref{eq:potential-coerc} and \eqref{eq:apriori-coercive}, we have
\[
\frac12\frac{\mathd}{\mathd t}\norm{Y}^2_{L^2_{r,k}}\leq -c\norm{Y}^2_{H^{s_1}_{r,k}} -N_1\norm{Y}^{q+1}_{L^{q+1}_{r,k}}+N_2\big(\norm{\rho}_{L^1(\R)}+\norm{Z}^{q+1}_{L^{q+1}_{r,k}}\big)+\sqrt2\scal{u}{Y}_{L^2_{r,k}}
\]
Since $\rho$ takes values in $(0,1]$, we have $\sqrt2\abs{\scal{u}{Y}_{L^2_{r,k}}}\leq \frac{c}{2}\norm{Y}^2_{L^2_{r,k}}+\frac{1}{c}\norm{u}^2_{L^2}$ hence, by \eqref{eq:weighted-holder-Lr}, there exists $S_1>0$ such that
\[
\frac12\frac{\mathd}{\mathd t}\norm{Y}^2_{L^2_{r,k}}\leq -\frac{c}{2}\norm{Y}^2_{H^{s_1}_{r,k}} -N_1\norm{Y}^{q+1}_{L^{q+1}_{r,k}}+N_2\left(\norm{\rho}_{L^1(\R)}+S_1\Big(1+\norm{Z}^{2q}_{C^\gamma_{\rho_{\varepsilon,k}}}+\norm{u}^2_{L^2}\Big)\right),
\]
and so by the Gronwall inequality we obtain \eqref{eq:aprioriL2}.\\

\noindent (ii) Multiplying in $L^2_{r,k}$ both sides of \eqref{eq:systemZY3}  by $\abs{Y}^{q-1}Y$, we get 
\[
\frac{1}{q+1}\frac{\mathd}{\mathd t}\norm{Y}^{q+1}_{L^{q+1}_{r,k}}+\scal{\mathcal{A}Y}{\abs{Y}^{q-1}Y}_{L^2_{r,k}}-\scal{\mathcal{B}(X)}{\abs{Y}^{q-1}Y}_{L^2_{r,k}}=\sqrt2\scal{u}{\abs{Y}^{q-1}Y}_{L^2_{r,k}} .
\]
By \eqref{eq:Kato-type}, extended to $Y$ by density, we have 
\[
\scal{\mathcal{A}Y}{\abs{Y}^{q-1}Y}_{L^2_{r,k}}\geq\big(m^2-\frac{\omega_{r,k}}{q+1}\big)\norm{Y}^{q+1}_{L^{q+1}_{r,k}}.
\]
Multiplying \eqref{eq:potential-coerc} with $a=Y$ and $b=Z$ by $\abs{Y}^{q-1}$, and using the Young inequality with exponents $\frac{2q}{q-1}$ and $\frac{2q}{q+1}$ (for $q>1$), there exists a constant $S_2>0$ such that
\[
\scal{\mathcal{B}(X)}{\abs{Y}^{q-1}Y}_{L^2_{r,k}}\geq-\frac{N_1}{2}\norm{Y}^{2q}_{L^{2q}_{r,k}}-S_2\Big(\norm{\rho}_{L^1(\R)}+\norm{Z}^{2q}_{L^{2q}_{r,k}}\Big),
\]
while 
\[
\sqrt2\abs{\scal{u}{\abs{Y}^{q-1}Y}_{L^2_{r,k}}}\leq\frac{N_1}{4}\norm{Y}^{2q}_{L^{2q}_{r,k}}+\frac{2}{N_1}\norm{u}^2_{L^2}.
\]
Hence, taking into account also \eqref{eq:weighted-holder-Lr}, there exists a constant $S_3>0$ such that
\[
\frac{1}{q+1}\frac{\mathd}{\mathd t}\norm{Y}^{q+1}_{L^{q+1}_{r,k}}\leq - \big(m^2-\frac{\omega_{r,k}}{q+1}\big)\norm{Y}^{q+1}_{L^{q+1}_{r,k}}-\frac{N_1}{4}\norm{Y}^{2q}_{L^{2q}_{r,k}}+S_3\Big(1+\norm{Z}^{2q}_{C^\gamma_{\rho_{\varepsilon,k}}}+\norm{u}^2_{L^2}\Big),
\]
and so by the Gronwall inequality we obtain  \eqref{eq:aprioriLp}.\\

\noindent (iii) Set $v:=\rho_{r/2,k}Y$. By \eqref{eq:commutatordef}, multiplying \eqref{eq:systemZY3} by $\rho_{r/2,k}$ we get
\[
\partial_tv+\mathcal{A}v=\mathcal{C}_{r,k}Y+\rho_{r/2,k}\big(\mathcal{B}(X)+\sqrt2u\big),
\]
Multiplying in $L^2$ by $\mathcal{A}v$, using the Young inequality and recalling that $\rho$ takes values in $[0,1]$, we get
\[
\frac{\mathd}{\mathd t}\scal{\mathcal{A}v}{v}_{L^2}+\norm{\mathcal{A}v}^2_{L^2}\leq2\norm{\mathcal{C}_{r,k}Y}^2_{L^2}+4\norm{\mathcal{B}(X)}^2_{L^2_{r,k}}+8\norm{u}^2_{L^2} .
\]
We note that $\norm{\mathcal{A}v}^2_{L^2}\geq m^2\scal{\mathcal{A}v}{v}_{L^2}$, indeed, by the Cauchy-Schwarz and Young inequalities and since $\scal{\mathcal{A}v}{v}_{L^2}\geq m^2\norm{v}^2_{L^2}$, we have
\[
\scal{\mathcal{A}v}{v}_{L^2}\leq\frac{1}{2m^2}\norm{\mathcal{A}v}^2_{L^2}+\frac{m^2}{2}\norm{v}^2_{L^2}\leq\frac{1}{2m^2}\norm{\mathcal{A}v}^2_{L^2}+\frac12\scal{\mathcal{A}v}{v}_{L^2} .
\]
Whether by \eqref{eq:apriori-commutator} we have $\norm{\mathcal{C}_{r,k}Y}_{L^2}\leq S_2\norm{Y}_{H^{s_1}_{r,k}}$, and, by \eqref{eq:potential-growth} and \eqref{eq:weighted-holder-Lr}, we get
\[
\norm{\mathcal{B}(X)}^2_{L^2_{r,k}}\leq N_2\Big(\norm{\rho}_{L^1(\R)}+2^{2q-1}\norm{Y}^{2q}_{L^{2q}_{r,k}}+2^{2q-1}N_3k^{-1/2}\big(1+\norm{Z}^{2q}_{C^\gamma_{\rho_{\varepsilon,k}}}\big)\Big).
\]
Hence there exists a constant $S_4>0$, such that
\[
\frac{\mathd}{\mathd t}\scal{\mathcal{A}v}{v}_{L^2}\leq -m^2\scal{\mathcal{A}v}{v}_{L^2}+S_4\Big(\norm{Y}^2_{H^{s_1}_{r,k}}+\norm{Y}^{2q}_{L^{2q}_{r,k}}+1+\norm{Z}^{2q}_{C^\gamma_{\rho_{\varepsilon,k}}}+\norm{u}^2_{L^2}\Big),
\]
and so, applying first the Gronwall inequality, then using that, by \eqref{eq:energy-is-Hs1} and \eqref{eq:weightedsobolev}, $\scal{\mathcal{A}v}{v}_{L^2}$ is equivalent to $\norm{Y}^2_{H^{s_1}_{r,k}}$ (up to enlarging $S_4$), and finally \eqref{eq:aprioriL2} and \eqref{eq:aprioriLp}, we obtain \eqref{eq:aprioriHs}.
\end{proof}

\subsection{Moments of invariant measures}\label{section:consequences}

Let $\{P(t)\}_{t\geq0}$ be the transition semigroup associated to \eqref{eq:SPDEmain1Abstract} given by \eqref{eq:semiL2rk}. In this section we are going to prove that every invariant measure for $\{P(t)\}_{t\geq 0}$ in the sense of Definition \ref{definition:invariant} verifies the condition \eqref{eq:cond-momenti-inv}. Hence, this shows that assumption \eqref{eq:cond-momenti-inv} can be removed from Theorem \ref{thm:main-inv-ii}.

\begin{lemm}\label{lemma:Ttlimit}
Assume \ref{H0}, \ref{H1}, \ref{H2}, \ref{H3} and \ref{H4} hold. Let $\mu$ be an invariant measure for $\{P(t)\}_{t\geq 0}$ in the sense of Definition \ref{definition:invariant}. For every Borel function $\Phi:L^2_{r,k}\rightarrow [0,+\infty)$,  if there exist a Borel set $B\subseteq L^2_{r,k}$ with $\mu(B)=1$ and a constant $S\geq0$ such that
\[
\limsup_{t\rightarrow+\infty}\frac1t\int_0^tP(s)\Phi(x)\mathd s\leq S,\qquad x\in B,
\]
then $\int_{L^2_{r,k}}\Phi\mathd\mu\leq S$.
\end{lemm}

\begin{proof}
Since $P(s)(\min\{\Phi, n\})\leq P(s)\Phi$ for every $n\in\N$, by the monotone convergence theorem it suffices to consider $\Phi$ bounded. By the invariance of $\mu$ and the Fubini theorem, $\int \Phi\mathd\mu=\int\frac1t\int_0^tP(s)\Phi\mathd s\mathd\mu$ for every $t>0$. Since $0\leq\frac1t\int_0^tP(s)\Phi\mathd s\leq\norm{\Phi}_{\infty}$, by the Fatou lemma and $\mu(B)=1$ we obtain
\[
\int_{L^2_{r,k}}\Phi\mathd\mu\leq\int_B\limsup_{t\rightarrow+\infty}\frac1t\int_0^tP(s)\Phi\mathd s\,\mathd\mu\leq S. \qedhere
\]
\end{proof}

\begin{lemm}\label{lemma:invariantmoments}
Assume \ref{H0}, \ref{H1}, \ref{H2}, \ref{H3}, \ref{H4} and \ref{H5} hold. Let $q\geq 1$ be the constant given by \eqref{H2:pol} and let $\gamma\in\big(0,s_1-\frac12\big)$. Let $\mu$ be an invariant measure for $\{P(t)\}_{t\geq 0}$ in the sense of Definition \ref{definition:invariant}. Then there exists a constant $K_1>0$, independent of $\mu$, such that for every $\tau\in\R$, we have
\begin{equation}\label{eq:invariantmoments}
\sup_{\tau\in\R}\int_{L^2_{r,k}}\Big(\norm{x}^2_{C^\gamma_{\rho^\tau_{r/2,k}}}+\norm{x}^{q+1}_{L^{q+1}_{\rho^\tau_{r,k}}}+\norm{x}^{2q}_{L^{2q}_{\rho^\tau_{r,k}}}\Big)\mu(\mathd x)\leq K_1,
\end{equation}
where $\rho^{\tau}_{r,k}:=\rho_{r,k}(\cdot-\tau)$. In particular $\mu\big(C^\gamma_{\rho_{r/2,k}}\cap L^{2q}_{r,k}\big)=1$.
\end{lemm}

\begin{proof}
Fix $\tau\in\R$ and $\varepsilon\in\big(0,\frac{r-1}{2q}\big)$, and set $\rho:=\rho^\tau_{r,k}$. By Remark \ref{rmk:stochconvholder} and \eqref{eq:stochconvholderuniform}, there exists a constant $S_1>0$, independent of $t$ and of $\tau$, such that
\begin{equation}\label{eq:invariantmomentsZ}
\sup_{t\geq0}\mathbb{E}\Big[\norm{W_{\mathcal{A}}(t)}^{2}_{C^\gamma_{\rho^\tau_{r/2,k}}}+\norm{W_{\mathcal{A}}(t)}^{2q}_{C^\gamma_{\rho^\tau_{\varepsilon,k}}}\Big]\leq S_1.
\end{equation}
Set
\[
f_1(x):=\norm{x}^2_{C^\gamma_{\rho^\tau_{r/2,k}}}+\norm{x}^{q+1}_{L^{q+1}_{\rho}},\qquad x\in L^2_{r,k},
\]
with the convention $f_1(x)=+\infty$ if $x\notin C^\gamma_{\rho^\tau_{r/2,k}}\cap L^{q+1}_{\rho}$.

Let $x\in L^2_{r,k}$. Recalling that $X(s,x)=Y(s,x)+W_{\mathcal{A}}(s)$ for every $s\geq0$, by \eqref{eq:sobolevweighted} we have 
\[
f_1(X(s,x))\leq 2N_4^2\norm{Y(s)}^2_{H^{s_1}_{\rho}}+2^{q}\norm{Y(s)}^{q+1}_{L^{q+1}_{\rho}}+2\norm{W_{\mathcal{A}}(s)}^{2}_{C^\gamma_{\rho^\tau_{r/2,k}}}+2^{q}\norm{W_{\mathcal{A}}(s)}^{q+1}_{L^{q+1}_{\rho}} .
\]
Since $q+1\leq2q$ and $\varepsilon<\frac{r-1}{2q}$, we have $r-(q+1)\varepsilon>1$, so that by \eqref{eq:weighted-holder-Lr} with $j=q+1$ there exists $S_2>0$, independent of $\tau$, such that
\[
\norm{W_{\mathcal{A}}(s)}^{q+1}_{L^{q+1}_{\rho}}\leq N_{3}k^{-1/2}\norm{W_{\mathcal{A}}(s)}^{q+1}_{C^\gamma_{\rho^\tau_{\varepsilon,k}}}\leq S_2\Big(1+\norm{W_{\mathcal{A}}(s)}^{2q}_{C^\gamma_{\rho^\tau_{\varepsilon,k}}}\Big).
\]
Hence, by \eqref{eq:invariantmomentsZ}, there exists a constant $S_3>0$, independent of $x$, $s$ and $\tau$, such that
\begin{equation}\label{eq:invariantmomentsf1}
P(s)f_1(x)=\mathbb{E}\big[f_1(X(s,x))\big]\leq2N_4^2\mathbb{E}\Big[\norm{Y(s)}^2_{H^{s_1}_{\rho}}\Big]+2^{q}\mathbb{E}\Big[\norm{Y(s)}^{q+1}_{L^{q+1}_{\rho}}\Big]+S_3 .
\end{equation}

Taking the expectation in \eqref{eq:aprioriL2} with $\ell=0$ and $u=0$, and using \eqref{eq:invariantmomentsZ}, we obtain
\begin{equation}\label{eq:invariantmomentsf2}
\int_0^t\mathbb{E}\Big[\norm{Y(s)}^2_{H^{s_1}_{\rho}}+\norm{Y(s)}^{q+1}_{L^{q+1}_{\rho}}\Big]\mathd s\leq K_0\Big(\norm{x}^2_{L^2_{\rho}}+(1+S_1)t\Big),\qquad t>0 .
\end{equation}
Therefore, integrating \eqref{eq:invariantmomentsf1} over $[0,t]$ and dividing by $t$, by \eqref{eq:invariantmomentsf2} there exists a constant $S_4>0$, independent of $x$, $t$ and $\tau$, such that
\[
\frac1t\int_0^tP(s)f_1(x)\mathd s\leq S_4\Big(\frac1t\norm{x}^2_{L^2_{\rho}}+1\Big),\qquad t>0 .
\]
Letting $t\rightarrow+\infty$, the first term vanishes for every $x\in L^2_{r,k}$, so that by Lemma \ref{lemma:Ttlimit} with $B=L^2_{r,k}$ we have
\begin{equation}\label{eq:invariantmomentsstep1}
\int_{L^2_{r,k}}\Big(\norm{x}^2_{C^\gamma_{\rho^\tau_{r/2,k}}}+\norm{x}^{q+1}_{L^{q+1}_{\rho}}\Big)\mu(\mathd x)\leq S_4.
\end{equation}
In particular $\mu(L^{q+1}_{r,k})=1$. It remains to bound the $L^{2q}$-moment. Set $f_2(x):=\norm{x}^{2q}_{L^{2q}_\rho}$, with $f_2(x)=+\infty$ if $x\notin L^{2q}_\rho$, and let $x\in L^{q+1}_{r,k}$. Since $X(s,x)=Y(s)+W_{\mathcal A}(s)$ and $r-2q\varepsilon>1$, by \eqref{eq:weighted-holder-Lr} with $j=2q$ we get
\[
f_2(X(s,x))\leq2^{2q-1}\norm{Y(s)}^{2q}_{L^{2q}_\rho}+2^{2q-1}N_3k^{-1/2}\norm{W_{\mathcal A}(s)}^{2q}_{C^\gamma_{\rho^\tau_{\varepsilon,k}}} .
\]
Taking the expectation in \eqref{eq:aprioriLp} with $\ell=0$, $u=0$ and $z=0$, and using \eqref{eq:invariantmomentsZ}, we obtain
\[
\int_0^t\mathbb E\Big[\norm{Y(s)}^{2q}_{L^{2q}_\rho}\Big]\mathd s\leq K_0\Big(\norm{x}^{q+1}_{L^{q+1}_\rho}+(1+S_1)t\Big),\qquad t>0,
\]
so that there exists $S_5>0$, independent of $x$, $t$ and $\tau$, such that
$\frac1t\int_0^tP(s)f_2(x)\mathd s\leq S_5\big(\frac1t\norm{x}^{q+1}_{L^{q+1}_\rho}+1\big)$. By Lemma \ref{lemma:Ttlimit} with $B=L^{q+1}_{r,k}$ we get $\int f_2\,\mathd\mu\leq S_5$, which together with \eqref{eq:invariantmomentsstep1} gives \eqref{eq:invariantmoments} with $K_1:=S_4+S_5$.
\end{proof}

\subsection{Fernique-type theorem for invariant measures of \eqref{eq:SPDEmain1psuedo}}

Finally, in this subsection, we prove a Fernique-type theorem. To do so, we use the Boué-Dupuis formula.

For $t>0$ we denote by $\mathcal P_t$ the class of the admissible controls, that is the set of the predictable processes $u:[0,t]\times\Omega\rightarrow L^2$ such that $\norm{u}_{L^2([0,t];L^2)}\leq N$ almost surely, for some $N>0$. We fix a separable Hilbert space $U_1\supseteq L^2$ such that the embedding $L^2\hookrightarrow U_1$ is Hilbert-Schmidt, so that $\{W(s)\}_{s\geq0}$ is a $U_1$-valued Wiener process with continuous trajectories. By the strong well-posedness of the mild form of \eqref{eq:SPDEmain1Abstract} (see Proposition \ref{prop:wellposedness}), for every $t>0$ and $x\in L^{2q}_{r,k}$ there exists a Borel map
\[
\Theta_{t,x}:C([0,t];U_1)\longrightarrow C([0,t];L^{2q}_{r,k})
\]
such that, for every $L^2$-cylindrical Wiener process $\{\widetilde W(s)\}_{s\in[0,t]}$ on a filtered probability space, $\Theta_{t,x}(\widetilde W)$ is the mild solution to \eqref{eq:SPDEmain1Abstract} driven by $\widetilde W$; in particular $X(\cdot,x)=\Theta_{t,x}(W)$ $\mathbb P$-a.s.

\begin{lemm}[Boué-Dupuis formula]\label{lemma:BD}
Assume \ref{H0}, \ref{H1}, \ref{H2} and \ref{H3} hold. Let $t>0$, $x\in L^{2q}_{r,k}$ and let $G:C([0,t];L^{2q}_{r,k})\rightarrow[0,+\infty]$ be a Borel functional. Then
\begin{equation}\label{eq:BDSPDE}
\log\mathbb E\big[e^{G(X(\cdot,x))}\big]=\sup_{u\in\mathcal P_t}\mathbb E\left[G\big(X^{u}(\cdot,x)\big)-\frac12\norm{u}^2_{L^2([0,t];L^2)}\right],
\end{equation}
both sides being possibly equal to $+\infty$.
\end{lemm}

\begin{proof}
Let $u\in\mathcal P_t$. Since $\norm{u}_{L^2([0,t];L^2)}$ is bounded, the Novikov condition holds and, by the Girsanov theorem, $\widetilde W:=W+\int_0^{\cdot}u(s)\mathd s$ is an $L^2$-cylindrical Wiener process under a probability measure $\mathbb Q$ equivalent to $\mathbb P$. Replacing $W$ by $\widetilde W$, the mild formulation of \eqref{eq:SPDEmain1Abstract} is exactly the mild formulation of \eqref{eq:systemZY1}, so that $X^u(\cdot,x)=\Theta_{t,x}(\widetilde W)$ $\mathbb Q$-a.s., and hence
\begin{equation}\label{eq:shiftsolution}
\Theta_{t,x}\Big(W+\int_0^{\cdot}u(s)\mathd s\Big)=X^{u}(\cdot,x)\qquad\mathbb P\text{-a.s.},\qquad u\in\mathcal P_t .
\end{equation}
The Bou\'e-Dupuis formula for a cylindrical Wiener process on a separable Hilbert space (see \cite[Theorem 3.6]{BudhirajaDupuis2000}, where it is also shown that the supremum can be restricted to bounded controls) states that, for every bounded Borel functional $\Psi:C([0,t];U_1)\rightarrow\R$, we have
\begin{equation}\label{eq:BD1}
\log\mathbb E\big[e^{\Psi(W)}\big]=\sup_{u\in\mathcal P_t}\mathbb E\left[\Psi\Big(W+\int_0^{\cdot}u(s)\mathd s\Big)-\frac12\norm{u}^2_{L^2([0,t];L^2)}\right].
\end{equation}
If $G$ is bounded, applying \eqref{eq:BD1} to $\Psi:=G\circ\Theta_{t,x}$ and using \eqref{eq:shiftsolution} we obtain \eqref{eq:BDSPDE}.

In the general case, we set $G_N:=G\wedge N$ for $N\in\N$, so that \eqref{eq:BDSPDE} holds for $G_N$. As $N\rightarrow+\infty$, the left-hand side converges to $\log\mathbb E[e^{G(X(\cdot,x))}]$ and, for every fixed $u\in\mathcal P_t$, the right-hand side converges to $\mathbb E\big[G(X^{u}(\cdot,x))-\frac12\norm{u}^2_{L^2([0,t];L^2)}\big]$, both by the monotone convergence theorem. Since the suprema over $N$ and over $u$ commute, \eqref{eq:BDSPDE} follows.
\end{proof}

Now exploiting the Boué-Dupuis formula and the estimates in the previous subsections, we obtain the following lemma.
\begin{lemm}\label{lemma:Fernique}
Assume \ref{H0}, \ref{H1}, \ref{H2}, \ref{H3}, \ref{H4} and \ref{H5} hold. Let $q\geq 1$ be the constant given by \eqref{H2:pol}, let $\gamma\in\big(0,s_1-\frac12\big)$ and $\ell\in(0,\ell_0)$, with $\ell_0$ given by Theorem \ref{theorem:apriori}. Then there exist $\eta_1,\eta_2>0$ and constants $K_2,K_3>0$ such that, for every $\tau\in\R$, $x\in L^{2q}_{r,k}$ and $t\geq1$, we have
\begin{equation}\label{eq:Fernique1}
\log\mathbb{E}\left[\exp\Big(\eta_1\norm{X(t,x)}^2_{C^{\gamma}_{\rho^\tau_{r/2,k}}}+\eta_2\norm{X(t,x)}^{q+1}_{L^{q+1}_{\rho^\tau_{r,k}}}\Big)\right]\leq K_2e^{-\ell(t-1)}\Big(\norm{x}^2_{L^2_{\rho^\tau_{r,k}}}+\norm{x}^{q+1}_{L^{q+1}_{\rho^\tau_{r,k}}}\Big)+K_3 .
\end{equation}
\end{lemm}

\begin{proof}
Fix $\tau\in\R$ and $\varepsilon\in\big(0,\frac{r-1}{2q}\big)$, and set $\rho:=\rho^\tau_{r,k}$. By Remark \ref{rmk:stochconvholder} and \eqref{eq:stochconvholderuniform}, there exists a constant $S_1>0$, independent of $t$ and of $\tau$, such that
\begin{equation}\label{eq:FerniqueZ}
\sup_{t\geq0}\mathbb{E}\Big[\norm{W_{\mathcal{A}}(t)}^{2}_{C^\gamma_{\rho^\tau_{r/2,k}}}+\norm{W_{\mathcal{A}}(t)}^{2q}_{C^\gamma_{\rho^\tau_{\varepsilon,k}}}\Big]\leq S_1 .
\end{equation}
Let $u\in\mathcal P_t$ and let $Y^{u}:=Y^{u}(\cdot,x)$ be the solution of \eqref{eq:systemZY3} with $y=x$ and $z=0$, so that $X^{u}(s,x)=Y^{u}(s)+W_{\mathcal{A}}(s)$ and $Z(s,0)=W_{\mathcal{A}}(s)$ for every $s\in[0,t]$. Without loss of generality, assume $K_0\geq1$, where $K_0$ is the constant in Theorem \ref{theorem:apriori}.

We first estimate $\norm{Y^{u}(t)}_{H^{s_1}_{\rho}}$, for which \eqref{eq:aprioriHs} cannot be applied directly since $x\notin H^{s_1}_{r,k}$. For $n\in\N$ we set
\[
\sigma_n:=\inf\Big\{s\in[0,1]:\norm{Y^{u}(s)}_{H^{s_1}_{\rho}}\leq n\Big\},\qquad\inf\emptyset:=+\infty .
\]
It is easy to see that $\sigma_n$ is a stopping time. By \eqref{eq:aprioriL2} with $\ell=0$, the random variable
\[
\Lambda:=\int_0^1\norm{Y^{u}(s)}^2_{H^{s_1}_{\rho}}\mathd s
\]
is finite $\mathbb{P}$-a.s. If $\sigma_{n-1}>1$, then $\norm{Y^{u}(s)}_{H^{s_1}_{\rho}}>n-1$ for every $s\in[0,1]$, and so $(n-1)^2<\Lambda$. Hence $n^*:=\min\{n\in\N:\sigma_n\leq1\}$ is finite $\mathbb{P}$-a.s. and, since $n^2\leq2(n-1)^2+2$, it satisfies $(n^*)^2\leq2+2\Lambda$. In particular $Y^{u}(\sigma_{n^*})\in H^{s_1}_{r,k}$, and $Y^{u}(\sigma_{n^*})\in L^{q+1}_{r,k}$ by \eqref{eq:aprioriLp}.

Since the proof of Theorem \ref{theorem:apriori} is pathwise and \eqref{eq:systemZY3} is autonomous for fixed $Z$ and $u$, estimate \eqref{eq:aprioriHs} applies on $[\sigma_{n^*},t]$ with initial datum $Y^{u}(\sigma_{n^*})$ and so, noting that $e^{-\ell(t-\sigma_{n^*})}\leq e^{-\ell(t-1)}$, we get
\[
\norm{Y^{u}(t)}^2_{H^{s_1}_{\rho}}\leq K_0\Big(e^{-\ell(t-1)}\Big(2+2\Lambda+\sup_{s\in[0,1]}\norm{Y^{u}(s)}^{q+1}_{L^{q+1}_{\rho}}\Big)+R_\ell(t)\Big),
\]
where $R_\ell(t)$ is given by \eqref{eq:R_l}. Estimating $\Lambda$ and the supremum by \eqref{eq:aprioriL2} and \eqref{eq:aprioriLp} with $\ell=0$, we obtain
\begin{equation}\label{eq:Fernique2}
\norm{Y^{u}(t)}^2_{H^{s_1}_{\rho}}\leq2K_0^2e^{-\ell(t-1)}\Big(1+\norm{x}^2_{L^2_{\rho}}+\norm{x}^{q+1}_{L^{q+1}_{\rho}}+2R_0(1)\Big)+K_0R_\ell(t).
\end{equation}

By \eqref{eq:sobolevweighted} and \eqref{eq:Fernique2} we have
\[
\norm{X^{u}(t,x)}^2_{C^{\gamma}_{\rho^\tau_{r/2,k}}}\leq2N_4^2\norm{Y^{u}(t)}^2_{H^{s_1}_{\rho}}+2\norm{W_{\mathcal{A}}(t)}^2_{C^{\gamma}_{\rho^\tau_{r/2,k}}},
\]
while, by \eqref{eq:aprioriLp} with $y=x$ and $z=0$ and by \eqref{eq:weighted-holder-Lr} with $j=q+1$, we have
\[
\norm{X^{u}(t,x)}^{q+1}_{L^{q+1}_{\rho}}\leq2^{q}K_0\Big(e^{-\ell t}\norm{x}^{q+1}_{L^{q+1}_{\rho}}+R_\ell(t)\Big)+2^{q}N_{3}k^{-1/2}\Big(1+\norm{W_{\mathcal{A}}(t)}^{2q}_{C^\gamma_{\rho^\tau_{\varepsilon,k}}}\Big).
\]
Since $R_0(1)$ and $R_\ell(t)$ contain $\norm{u}^2_{L^2([0,t];L^2)}$ with coefficient at most $1$, in the two estimates above the coefficient of $\norm{u}^2_{L^2([0,t];L^2)}$ is at most
\[
S_2:=2\eta_1N_4^2\big(4K_0^2+K_0\big)+2^{q}\eta_2K_0 ,
\]
so that, choosing $\eta_1,\eta_2>0$ such that $S_2\leq\frac12$, this term is absorbed by $-\frac12\norm{u}^2_{L^2([0,t];L^2)}$. Moreover, by \eqref{eq:FerniqueZ} and $\ell>0$, we have $\mathbb{E}[R_0(1)]\leq1+S_1+\mathbb{E}\big[\norm{u}^2_{L^2([0,t];L^2)}\big]$ and $\mathbb{E}[R_\ell(t)]\leq\frac{1+S_1}{\ell}+\mathbb{E}\big[\norm{u}^2_{L^2([0,t];L^2)}\big]$, so that all the remaining terms are bounded either by $K_2e^{-\ell(t-1)}\big(\norm{x}^2_{L^2_{\rho}}+\norm{x}^{q+1}_{L^{q+1}_{\rho}}\big)$ or by a constant $K_3$, with $K_2$ and $K_3$ independent of $u$, $x$, $t$ and $\tau$.\\

\noindent Finally, the functional $G:C([0,t];L^{2q}_{r,k})\rightarrow[0,+\infty]$ defined by
\[
G(\zeta):=\eta_1\norm{\zeta(t)}^2_{C^{\gamma}_{\rho^\tau_{r/2,k}}}+\eta_2\norm{\zeta(t)}^{q+1}_{L^{q+1}_{\rho}}
\]
is Borel, since the evaluation at $t$ is continuous. Hence, by Lemma \ref{lemma:BD}, the left-hand side of \eqref{eq:Fernique1} equals
\[
\sup_{u\in\mathcal P_t}\mathbb{E}\left[\eta_1\norm{X^{u}(t,x)}^2_{C^{\gamma}_{\rho^\tau_{r/2,k}}}+\eta_2\norm{X^{u}(t,x)}^{q+1}_{L^{q+1}_{\rho}}-\frac12\norm{u}^2_{L^2([0,t];L^2)}\right],
\]
and the estimates above give \eqref{eq:Fernique1}.
\end{proof}

\noindent The following Fernique-type theorem is an immediate consequence of Lemmas \ref{lemma:Ttlimit} and \ref{lemma:Fernique}.

\begin{thm}\label{theorem:Fernique}
Assume \ref{H0}, \ref{H1}, \ref{H2}, \ref{H3}, \ref{H4} and \ref{H5} hold. Let $q\geq 1$ be the constant given by \eqref{H2:pol} and let $\gamma\in\big(0,s_1-\frac12\big)$. Let $\mu$ be an invariant measure for $\{P(t)\}_{t\geq 0}$ in the sense of Definition \ref{definition:invariant}. Then there exist $\eta_3,\eta_4>0$, independent of $\mu$, such that
\begin{equation}\label{eq:Fernique}
\sup_{\tau\in\R}\int_{L^2_{r,k}}\exp\Big(\eta_3\norm{x}_{C^{\gamma}_{\rho^\tau_{r/2,k}}}^2+\eta_4\norm{x}^{q+1}_{L^{q+1}_{\rho^\tau_{r,k}}}\Big)\mu(\mathd x)<+\infty .
\end{equation}
\end{thm}

\begin{proof}
Fix $\tau\in\R$ and let $\eta_1$, $\eta_2$, $K_2$ and $K_3$ be given by Lemma \ref{lemma:Fernique}. Set
\[
f(x):=\exp\Big(\eta_1\norm{x}^2_{C^{\gamma}_{\rho^\tau_{r/2,k}}}+\eta_2\norm{x}^{q+1}_{L^{q+1}_{\rho^\tau_{r,k}}}\Big),\qquad x\in L^2_{r,k},
\]
with the convention $f(x)=+\infty$ if $x\notin C^\gamma_{\rho^\tau_{r/2,k}}\cap L^{q+1}_{\rho}$. By \eqref{eq:Fernique1}, for every $x\in L^{2q}_{r,k}$ we have $\limsup_{t\rightarrow+\infty}P(t)f(x)\leq e^{K_3}$, and hence also $\limsup_{t\rightarrow+\infty}\frac1t\int_0^tP(s)f(x)\mathd s\leq e^{K_3}$. Since $\mu(L^{2q}_{r,k})=1$ by Lemmas \ref{lemma:Ttlimit} with $B=L^{2q}_{r,k}$ and \ref{lemma:invariantmoments} we have $\int f\mathd\mu\leq e^{K_3}$. As $K_3$ does not depend on $\tau$, we obtain \eqref{eq:Fernique}.
\end{proof}

\begin{rmk}\label{rmk:LSI}
Exponential integrability estimates as in Theorem \ref{theorem:Fernique} are closely related to logarithmic Sobolev inequalities and to the hypercontractivity, supercontractivity and ultracontractivity of the associated semigroups (see, e.g., \cite{BakryGentilLedoux2014,RocknerWang2003,Wang1997}). Hence it would be interesting to investigate whether a reversible measure $\mu$ of \eqref{eq:SPDEmain1psuedo} satisfies a logarithmic Sobolev inequality and whether $\{P(t)\}_{t\geq0}$ is ultracontractive in $L^2(\mu)$; results of this type are proved in \cite{Ang-Big-Fer2025} for transition semigroups associated with dissipative stochastic systems in Hilbert spaces (see also \cite{Ang-Big-Fer2023} for related Harnack-type inequalities).
\end{rmk}

\subsection{A priori estimates for DLR solutions}

Under the additional coercivity assumption \ref{H5} on the potential, we are able to prove that every solution of the DLR equations \eqref{eq:DLRequations} associated to \eqref{eq:SPDEmain1Abstract} always verifies the moment condition \eqref{eq:momenti-DLR}.

\begin{thm}\label{thm:apriori}
Assume \ref{H0}, \ref{H1}, \ref{H2}, \ref{H3}, \ref{H4} and \ref{H5} hold. Let $q\geq 1$ be the constant given by \eqref{H2:pol}. Then for every $\mu$ solution of the DLR equations \eqref{eq:DLRequations} associated to \eqref{eq:SPDEmain1Abstract} such that $\mu\big(L^{2q}_{r,k}\cap C(\R)\big)=1$ we have
\[\int_{L^{2q}_{r,k}}\norm{x}^{2q}_{L^{2q}_{r,k}}\mu(\mathd x)<+\infty.\]
\end{thm}

\begin{proof}
Let $\mu$ be as in the statement and, for $T>0$, let $\mathcal I_T:=(-T,T)$. We first show that there exists $S_0>0$ such that, for every $T>0$ and $y\in L^2_{r,k}(\mathcal I_T^c)\cap C(\mathcal I_T^c)$,
\begin{equation}\label{eq:conditionalmoments}
\int_{L^2(\mathcal I_T)}\norm{x+\varphi^y}^{2q}_{L^{2q}_{r,k}(\mathcal I_T)}\nu_{\mathcal I_T,y}(\mathd x)\leq S_0\Big(1+\norm{\varphi^y}^{2q}_{L^{2q}_{r,k}(\mathcal I_T)}\Big),
\end{equation}
where $\varphi^y$ is the harmonic extension of $y$ in $\mathcal I_T$ and $\nu_{\mathcal I_T,y}$ is the invariant measure of \eqref{eq:SPDEL2(I)} with $\mathcal I=\mathcal I_T$ given by Proposition \ref{prop:invmeasureI}. This is a standard energy estimate for the invariant measure of \eqref{eq:SPDEL2(I)}, and we only point out that all the constants are independent of $T$ and $y$. We set $\mathcal I:=\mathcal I_T$, $\rho:=\rho_{r,k}$ and $\varphi:=\varphi^y$, which is bounded on $\overline{\mathcal I}$ by Theorem \ref{thm:harmonic-continuity}. For $x\in L^{2q}(\mathcal I)$ we write the mild solution of \eqref{eq:SPDEL2(I)} given by Proposition \ref{prop:wellposednessI} as $X_{\mathcal I}(t,x)=Y_{\mathcal I}(t)+W_{\mathcal A_{\mathcal I}}(t)$, where $W_{\mathcal A_{\mathcal I}}$ is the stochastic convolution \eqref{eq:stochconvI}, so that, by \eqref{eq:deriU},
\[
\partial_tY_{\mathcal I}(t)+\mathcal A_{\mathcal I}Y_{\mathcal I}(t)=-V'\big(Y_{\mathcal I}(t)+W_{\mathcal A_{\mathcal I}}(t)+\varphi\big),\qquad Y_{\mathcal I}(0)=x .
\]
By \eqref{eq:killeddomination}, the variance of $W_{\mathcal A_{\mathcal I}}(s,\xi)$ is bounded by the quantity $v_0(s)$ given by \eqref{eq:pointwisevariance} with $\alpha=0$, which does not depend on $\mathcal I$ and is bounded uniformly in $s\geq0$ (see Remark \ref{rmk:stochconvholder}), so that, as in \eqref{eq:Lpmoment}, there exists $S_1>0$, independent of $s$ and $T$, such that $\mathbb E[\norm{W_{\mathcal A_{\mathcal I}}(s)}^{2q}_{L^{2q}_{\rho}(\mathcal I)}]\leq S_1$ for every $s\geq0$. Moreover, since $\scal{\mathcal A_{\mathcal I}w}{\abs{w}^{q-1}w\rho}_{L^2(\mathcal I)}=\mathcal E(w,\abs{w}^{q-1}w\rho)$ for every $w\in{\rm Dom}(\mathcal A_{\mathcal I})$, arguing as in the proof of Lemma \ref{lemma:kato} with $p=q+1$, and using in $\mathcal E^{nl}$ the inequality $(a-b)(\abs{a}^{q-1}a\rho_1-\abs{b}^{q-1}b\rho_2)\geq\frac{1}{q+1}(\abs{a}^{q+1}-\abs{b}^{q+1})(\rho_1-\rho_2)$, which follows from the convexity of $a\mapsto\abs{a}^{q+1}$, by \eqref{eq:excessive} and \ref{H3} we get $\scal{\mathcal A_{\mathcal I}w}{\abs{w}^{q-1}w\rho}_{L^2(\mathcal I)}\geq(m^2-\frac{\omega_{r,k}}{q+1})\norm{w}^{q+1}_{L^{q+1}_\rho(\mathcal I)}\geq0$. Finally, by \eqref{eq:potential-coerc} and the Young inequality, there exists $S_2>0$, depending only on $q$, $N_1$ and $N_2$, such that $\abs{a}^{q-1}aV'(a+b)\geq\frac{N_1}{2}\abs{a}^{2q}-S_2(1+\abs{b}^{2q})$ for every $a,b\in\R$. Hence, multiplying the equation for $Y_{\mathcal I}$ by $\abs{Y_{\mathcal I}}^{q-1}Y_{\mathcal I}\rho$ (the computation is justified by the maximal $L^2$-regularity of $\mathcal A_{\mathcal I}$, see \cite{Des1964}, and by approximation), we obtain
\[
\frac{1}{q+1}\frac{\mathd}{\mathd s}\norm{Y_{\mathcal I}(s)}^{q+1}_{L^{q+1}_{\rho}(\mathcal I)}+\frac{N_1}{2}\norm{Y_{\mathcal I}(s)}^{2q}_{L^{2q}_{\rho}(\mathcal I)}\leq S_2\Big(\norm{\rho}_{L^1(\R)}+\norm{W_{\mathcal A_{\mathcal I}}(s)+\varphi}^{2q}_{L^{2q}_{\rho}(\mathcal I)}\Big),\qquad s>0 .
\]
Since, by \ref{H4}, $\nu_{\mathcal I,y}$ has a bounded density with respect to the Gaussian measure $\mu^{\mathcal I,0}_{\free}$, which is concentrated on $C(\overline{\mathcal I})$, we have $\int\norm{x}^{q+1}_{L^{q+1}_\rho(\mathcal I)}\nu_{\mathcal I,y}(\mathd x)<+\infty$. Therefore, integrating over $[0,t]$, taking the expectation, integrating with respect to $\nu_{\mathcal I,y}$, using its invariance and $X_{\mathcal I}=Y_{\mathcal I}+W_{\mathcal A_{\mathcal I}}$, dividing by $t$ and letting $t\rightarrow+\infty$, we obtain \eqref{eq:conditionalmoments}.

Now, fix $N\in\N$ and set
\[
\Phi_N(x):=\min\Big\{N,\norm{x}^{2q}_{L^{2q}_{r,k}}\Big\},\qquad x\in L^2_{r,k},
\]
with the convention $\norm{x}_{L^{2q}_{r,k}}=+\infty$ if $x\notin L^{2q}_{r,k}$, and
\[
B_T(x):=\norm{\varphi^{x|_{\mathcal I_T^c}}}^{2q}_{L^{2q}_{r,k}(\mathcal I_T)},\qquad x\in L^2_{r,k}\cap C(\R),\ T>0,
\]
where $\varphi^{x|_{\mathcal I_T^c}}$ is the harmonic extension of $x|_{\mathcal I_T^c}$ in $\mathcal I_T$. Let $T>0$. By Corollary \ref{remark:DLRtranslated} and \eqref{eq:invmeasureI}, since $\norm{x}^{2q}_{L^{2q}_{r,k}}=\norm{x}^{2q}_{L^{2q}_{r,k}(\mathcal I_T)}+\norm{x}^{2q}_{L^{2q}_{r,k}(\mathcal I_T^c)}$, $\int\min\{N,g+c\}\mathd\nu\leq\min\{N,\int g\mathd\nu+c\}$ for every probability measure $\nu$, Borel function $g\geq0$ and $c\geq0$, and by \eqref{eq:conditionalmoments}, for every $y\in L^2_{r,k}(\mathcal I_T^c)\cap C(\mathcal I_T^c)$ we have
\[
\int_{E^{\mathcal I_T,y}}\Phi_N(x)\mu^{\mathcal I_T,y}(\mathd x)=\int_{L^2(\mathcal I_T)}\Phi_N(x+\varphi^y)\nu_{\mathcal I_T,y}(\mathd x)\leq\min\Big\{N,S_0\Big(1+\norm{\varphi^y}^{2q}_{L^{2q}_{r,k}(\mathcal I_T)}\Big)+\norm{y}^{2q}_{L^{2q}_{r,k}(\mathcal I_T^c)}\Big\}.
\]
Hence, by \eqref{eq:formula-prob-cond}, \eqref{eq:legge-Ic} and since $\mu\big(L^{2q}_{r,k}\cap C(\R)\big)=1$, we obtain
\begin{equation}\label{eq:DLRmomentsT}
\int_{L^{2q}_{r,k}}\Phi_N\,\mathd\mu\leq\int_{L^{2q}_{r,k}\cap C(\R)}\min\Big\{N,S_0\big(1+B_T(x)\big)+\norm{x}^{2q}_{L^{2q}_{r,k}(\mathcal I_T^c)}\Big\}\mu(\mathd x),\qquad T>0 .
\end{equation}
Let $n\in\N$. Integrating \eqref{eq:DLRmomentsT} with respect to $T\in[n,n+1]$ and using the Tonelli theorem, the inequality $\int_n^{n+1}\min\{N,g(T)\}\mathd T\leq\min\{N,\int_n^{n+1}g(T)\mathd T\}$ and the inclusion $\mathcal I_T^c\subseteq\{\abs{\eta}\geq n\}$ for $T\geq n$, we get
\[
\int_{L^{2q}_{r,k}}\Phi_N\,\mathd\mu\leq\int_{L^{2q}_{r,k}\cap C(\R)}\min\Big\{N,S_0\Big(1+\int_n^{n+1}B_T(x)\,\mathd T\Big)+\int_{\{\abs{\eta}\geq n\}}\abs{x(\eta)}^{2q}\rho_{r,k}(\eta)\,\mathd\eta\Big\}\mu(\mathd x).
\]
By Lemma \ref{lemma:harmonicaverage} with $p=2q$, for every $x\in L^{2q}_{r,k}\cap C(\R)$ both $\int_n^{n+1}B_T(x)\,\mathd T$ and the last integral vanish as $n\rightarrow+\infty$, and so, by the dominated convergence theorem, $\int\Phi_N\mathd\mu\leq S_0$. Letting $N\rightarrow+\infty$, by the monotone convergence theorem we obtain $\int_{L^{2q}_{r,k}}\norm{x}^{2q}_{L^{2q}_{r,k}}\mu(\mathd x)\leq S_0$, which concludes the proof.
\end{proof}

\begin{rmk}\label{rmk:LebowitzPresutti}
Theorem \ref{thm:apriori} can be seen as a continuum and nonlocal counterpart of \cite[Theorem 4.4]{LebowitzPresutti1976}. In \cite{LebowitzPresutti1976} unbounded spin systems on a lattice are considered, with summable interactions and a superstable energy, namely bounded from below, in every finite region, by a positive multiple of the sum of the squares of the spins minus a constant times the number of sites; it is proved that every tempered solution of the DLR equations satisfies Gaussian bounds on its finite dimensional marginals, with constants independent of the solution. As in Theorem \ref{thm:apriori}, a quadratic confinement is sufficient, which corresponds to $q=1$. The main differences are two. First, in \cite{LebowitzPresutti1976} temperedness is understood in the sense of \cite{Ruelle1970}, namely the solutions are concentrated on configurations whose quadratic averages over the cubes centered at the origin are bounded uniformly in the size of the cubes, while here we only assume $\mu\big(L^{2q}_{r,k}\cap C(\R)\big)=1$, which allows the averages of $\abs{x}^{2q}$ over the intervals centered at the origin to grow polynomially with the size of the intervals. Second, the proof in \cite{LebowitzPresutti1976} is based on the multiscale superstability estimates of \cite{Ruelle1970}, while here the interaction is Gaussian, so that the boundary condition enters the conditional measures only through the harmonic extension (see Corollary \ref{remark:DLRtranslated}), and the boundary effects can be controlled by averaging over the size of the intervals (see Lemma \ref{lemma:harmonicaverage}). We also point out that, by Theorem \ref{thm:apriori}, every solution of the DLR equations as in Theorem \ref{thm:apriori} satisfies the assumptions of Theorem \ref{theorem:FerniqueDLR}, and hence the exponential bounds \eqref{eq:Fernique}, which play the role of the Gaussian bounds of \cite{LebowitzPresutti1976}.\\

\noindent The method of \cite{LebowitzPresutti1976} goes back to \cite{Ruelle1970}, where superstability estimates were obtained for classical systems of particles. In \cite{Funaki1991} the reversible measures of a Ginzburg-Landau type continuum model with a local operator and a bounded nonlinearity are characterized as Gibbs states, whose definition includes a moment condition. Finally, in \cite{AlbeverioKondratievPasurekRockner2006} a priori estimates for Euclidean Gibbs states are obtained through integration by parts formulas and Lyapunov functionals, under a superquadratic growth assumption on the single-site potential.
\end{rmk}

\subsection{Main results without moments condition}

Under the additional coercivity assumption \ref{H5} on the potential, we can remove the moment conditions \eqref{eq:cond-momenti-inv} and \eqref{eq:momenti-DLR} from all the main results. Taking into account Definitions \ref{definition:invariant}, \ref{definition:reversible} and \ref{definition:DLR}, combining Theorems \ref{thm:main-inv-ii}, \ref{theorem:infinitesimalsymmetri}, \ref{theorem:main} and the a priori estimates Lemma \ref{lemma:invariantmoments} and Theorem \ref{thm:apriori} we obtain the following result.

\begin{thm}\label{theorem:super-main}
Assume \ref{H0}, \ref{H1}, \ref{H2}, \ref{H3}, \ref{H4} and \ref{H5} hold. Let $\mu$ be a Borel probability measure on $L^2_{r,k}$ such that $\mu(L^{2q}_{r,k}\cap C(\R))=1$. Then the following conditions are equivalent:
  \begin{enumerate}
    \item $\mu$ is (semigroup) invariant for the SPDE \eqref{eq:SPDEmain1psuedo} on $L^2_{r,k}$;
    
    \item $\mu$ is infinitesimally invariant for the SPDE \eqref{eq:SPDEmain1psuedo} on $L^2_{r,k}$.
  \end{enumerate}
  Moreover also the following conditions are equivalent:
\begin{enumerate}
    \item $\mu$ is (semigroup) reversible;
    
    \item $\mu$ is infinitesimally reversible;
    
    \item $\mu$ satisfies the integration by parts formula associated with the SPDE \eqref{eq:SPDEmain1psuedo};

    \item $\mu$ satisfies the DLR equations for \eqref{eq:SPDEmain1psuedo}.
\end{enumerate}
\end{thm}

\noindent The following Corollary \ref{corollary:DLRuniqueness} follows directly by Theorems \ref{theorem:super-main}, \ref{thm:rev-ergo} and Corollary \ref{corollary:uniquetranslation}.

\begin{cor}
  \label{corollary:DLRuniqueness}Assume that \ref{H0}, \ref{H1}, \ref{H2}, \ref{H3}, \ref{H4} and \ref{H5} hold. Let $\mu$ be a Borel probability measure on $L^2_{r,k}$ such that $\mu(L^{2q}_{r,k}\cap C(\R))=1$. If $\mu$ is the unique measure that satisfies the DLR equations for \eqref{eq:SPDEmain1psuedo}, in the sense of Definition \ref{definition:DLR}, then $\mu$ is translation-invariant and ergodic for the transition semigroup $\{P(t)\}_{t\geq 0}$ defined by \eqref{eq:semiL2rk}.
\end{cor}

\noindent We also have an important consequence of Theorem \ref{theorem:Fernique} which follows directly from Theorem \ref{theorem:super-main}.

\begin{thm}\label{theorem:FerniqueDLR}
Assume \ref{H0}, \ref{H1}, \ref{H2}, \ref{H3}, \ref{H4} and \ref{H5} hold.  Let $q\geq 1$ be the constant given by \eqref{H2:pol} and let $\gamma\in\big(0,s_1-\frac12\big)$. Let $\mu$ be a solution of the DLR equations \eqref{eq:DLRequations} associated to \eqref{eq:SPDEmain1Abstract} such that $\mu\big(L^{2q}_{r,k}\cap C(\R)\big)=1$. Then there exist $\eta_3,\eta_4>0$ such that \eqref{eq:Fernique} holds, namely
\begin{equation*}
\sup_{\tau\in\R}\int_{L^2_{r,k}}\exp\Big(\eta_3\norm{x}_{C^{\gamma}_{\rho^\tau_{r/2,k}}}^2+\eta_4\norm{x}^{q+1}_{L^{q+1}_{\rho^\tau_{r,k}}}\Big)\mu(\mathd x)<+\infty .
\end{equation*}
\end{thm}

\section{Existence of reversible measures/solutions to DLR equations}
\label{section:existence}

In this section we prove the existence of a reversible measure for \eqref{eq:SPDEmain1psuedo}, or equivalently of a solution to the DLR equations \eqref{eq:DLRequations}. We prove this result by proving that some measures which are suitable approximations of the solutions to the DLR equations are tight and that their limits solve the integration by parts formula \eqref{eq:IPF}. The proof of tightness is based on uniform a priori estimates of the solutions of the approximations of equation \eqref{eq:SPDEmain1psuedo} proposed in \cite{AlbeverioKusuoka2020,GubinelliHofmanova2021} (see also \cite{de2022singular}). Finally, in Subsection \ref{section:translationinvariance}, by averaging the translates of a reversible measure, we prove the existence of a translation invariant reversible measure.

Since, due to the general non-uniqueness of the invariant measures to the SPDE \eqref{eq:SPDEmain1psuedo} or to the DLR equation \eqref{eq:DLRequations} (see Section \ref{section:nonuniqueness}), different approximations could lead to different limits, we propose the construction in two cases: the measure obtained by Gaussian boundary conditions (corresponding to the free boundary conditions in \cite{GuerraRosenSimon1976}), and the case of (regular) Dirichlet boundary conditions (usually called fixed or configurations as boundary conditions in statistical mechanics, see, e.g., Chapter 3 of \cite{FriedliVelenik2018}). The reasons to focus on these two specific choices are: to prove existence of a reflection-positive measure solution to the DLR equation in Subsection \ref{section:ref.positive} (which uses Gaussian boundary conditions), and the proof of non-uniqueness of solutions to DLR equations in Section \ref{section:nonuniqueness} (which uses Dirichlet boundary conditions). Since the proofs with Gaussian and Dirichlet boundary conditions are very similar, we provide more details in the first case, and describe only the main differences in the second one. Other kinds of boundary conditions, such as periodic or Neumann boundary conditions (corresponding to the periodic or free boundary conditions of statistical mechanics, see, e.g., Chapter 3 of \cite{FriedliVelenik2018}), can be considered.

Throughout this section we assume \ref{H0}, \ref{H1}, \ref{H2}, \ref{H3}, \ref{H4} and \ref{H5} hold, and we let $q\geq1$ be the constant given by \eqref{H2:pol}. Moreover we fix
\[
1<r_0<r<1+2s_2,\qquad 0<\gamma<\gamma'<s_1-\frac12,\qquad 0<\varepsilon<\frac{r_0-1}{2q},
\]
and $k\in(0,k_0)$ with $k_0:=\min\{k_0(r_0,m);k_0(r,m)\}$ according to \ref{H3}, and for $\tau\in\R$ we write $\rho^\tau_{r,k}:=\rho_{r,k}(\cdot-\tau)$ as in Remark \ref{rmk:weighttranslation}. We recall that, by \ref{H3}, all the results of Subsection \ref{section:apriori} hold for $r_0$ as well as for $r$.

\subsection{Construction with Gaussian boundary conditions}
We recall that $\mu_{\free}$ is the Gaussian measure introduced in Subsection \ref{section:conditionalprobability}, that is the unique invariant measure of \eqref{eq:systemZY2}, and $\{\mu^{\mathcal I,y}_{\free}\}_{y\in L^2_{r,k}(\mathcal I^c)}$ are its conditional probabilities given by Proposition \ref{proposition:conditionalGaussian}, where $\mathcal{I}$ is an open bounded interval of $\R$. In this subsection we consider the family of measures
\begin{equation}\label{eq:existencemutheta}
\mu^\vartheta(\mathd x):=\frac{1}{\mathcal Z^\vartheta}\exp\left(-\int_{-\vartheta}^\vartheta V(x(\xi))\mathd\xi\right)\mu_{\free}(\mathd x),\qquad \vartheta>0,
\end{equation}
where, by \ref{H4}, the density is bounded by $e^{-2\vartheta M_1}$, so that $\mathcal Z^\vartheta\in(0,+\infty)$ and $\mu^\vartheta$ is well defined. Since $\mu^\vartheta$ is absolutely continuous with respect to $\mu_{\free}$ with bounded density, and since $\mu_{\free}\big(L^{2q}_{r,k}\cap C(\R)\big)=1$ and $\mu_{\free}$ has finite moments of every order in $L^{2q}_{r,k}$, we have
\begin{equation}\label{eq:existencemuthetamoments}
\mu^\vartheta\big(L^{2q}_{r,k}\cap C(\R)\big)=1,\qquad\int_{L^{2q}_{r,k}}\norm{x}^{q}_{L^{2q}_{r,k}}\mu^\vartheta(\mathd x)<+\infty .
\end{equation}

We call \eqref{eq:existencemutheta} Gaussian boundary conditions for the following reason. Let $\mathcal I_\vartheta:=(-\vartheta,\vartheta)$ and let $\{\mu^{\mathcal I_\vartheta,y}\}_{y}$ be the conditional probabilities of $\mu^\vartheta$ given $x|_{\mathcal I_\vartheta^c}=y$, in the sense of Subsection \ref{section:conditionalprobability}. Since the density in \eqref{eq:existencemutheta} depends on $x$ only through $x|_{\mathcal I_\vartheta}$, disintegrating $\mu_{\free}$ we obtain, for $\mu^\vartheta$-almost every $y$,
\begin{equation}\label{eq:existencecondtheta}
\frac{\mathd\mu^{\mathcal I_\vartheta,y}}{\mathd\mu_{\free}^{\mathcal I_\vartheta,y}}=\frac{1}{\mathcal Z^{\mathcal I_\vartheta}(y)}\exp\left(-\int_{\mathcal I_\vartheta}V(x(\xi))\mathd\xi\right),
\end{equation}
namely the conditional measures of $\mu^\vartheta$  on $\mathcal I_\vartheta$  are exactly the ones prescribed by the DLR equations \eqref{eq:DLRequations}, while the law of the boundary condition is
\[
\mu^\vartheta\circ R_{\mathcal I_\vartheta^c}^{-1}(\mathd y)=\frac{\mathcal Z^{\mathcal I_\vartheta}(y)}{\mathcal Z^\vartheta}\,\mu_{\free}\circ R_{\mathcal I_\vartheta^c}^{-1}(\mathd y).
\]
Thus $\mu^\vartheta$ is a solution of the DLR equations on $\mathcal I_\vartheta$, but not on larger intervals, where the interaction outside $\mathcal I_\vartheta$ is missing.

To describe $\mu^\vartheta$ on an arbitrary bounded interval, and to associate an SPDE with it, we set
\[
\chi_\vartheta:=\mathbb{I}_{[-\vartheta,\vartheta]},\qquad V_\vartheta(\xi,a):=\chi_\vartheta(\xi)V(a),\qquad \xi,a\in\R .
\]
The potential $V_\vartheta$ depends on $\xi$, but only through the factor $\chi_\vartheta$, and it satisfies \eqref{H2:diss}, \eqref{H2:pol} and \ref{H4} uniformly in $\xi$, with the constants $\min\{M_0,0\}$, $C$ and $\min\{M_1,0\}$ respectively; since the proofs of Part \ref{part:main} use no other property of $V$, all the results of Sections \ref{section:abstractreformulation}-\ref{sect:int-DLR} hold with $V$ replaced by $V_\vartheta$, with constants independent of $\vartheta$. We denote by
\[
\mathcal{B}_\vartheta(x):=-\chi_\vartheta\,(V'\circ x),\qquad x\in L^{2q}_{r,k},
\]
the associated Nemytskii operator, so that $\mathcal{B}_\vartheta(x)=\chi_\vartheta\mathcal{B}(x)$, and we consider the SPDE
\begin{equation}\label{eq:existenceSPDE}
\begin{cases}
\mathd X^{\vartheta}(t)+\mathcal{A}X^{\vartheta}(t)\mathd t=\mathcal{B}_\vartheta\big(X^{\vartheta}(t)\big)\mathd t+\sqrt{2\rho_{r,k}}\mathd W_{r,k}(t),\qquad t>0,\\
X^{\vartheta}(0)=x,
\end{cases}
\end{equation}
where $\mathcal{A}$ and $\{W_{r,k}(t)\}_{t\geq0}$ are the same as in \eqref{eq:SPDEmain1Abstract}. By Proposition \ref{prop:wellposedness} applied to $V_\vartheta$, equation \eqref{eq:existenceSPDE} has a unique generalized mild solution $\{X^{\vartheta}(t,x)\}_{t\geq0}$ for every $x\in L^{2}_{r,k}$, and we denote by $\{P^\vartheta(t)\}_{t\geq0}$ the associated transition semigroup on $L^2_{r,k}$.

The measure $\mu^\vartheta$ admits the following three equivalent descriptions.

\begin{enumerate}
\item \emph{DLR equations for $V_\vartheta$.} For every bounded open interval $\mathcal I$ of $\R$ and for $\mu^\vartheta$-almost every $y\in L^2_{r,k}(\mathcal I^c)\cap C(\mathcal I^c)$,
\begin{equation}\label{eq:existenceDLR}
\frac{\mathd\mu^{\vartheta,\mathcal I,y}}{\mathd\mu_{\free}^{\mathcal I,y}}=\frac{1}{\mathcal Z^{\vartheta,\mathcal I}(y)}\exp\left(-\int_{\mathcal I}\chi_\vartheta(\xi)V(x(\xi))\mathd\xi\right),
\end{equation}
which is proved exactly as \eqref{eq:existencecondtheta}, the density in \eqref{eq:existencemutheta} being the product of a function of $x|_{\mathcal I}$ and of a function of $x|_{\mathcal I^c}$.

\item \emph{Integration by parts formula.} By \eqref{eq:existencemuthetamoments}, the measure $\mu^\vartheta$ satisfies the assumptions of Theorem \ref{theorem:main} relative to $V_\vartheta$; hence, by point 1. and by that theorem,
\begin{equation}\label{eq:existenceIPF}
\int\scal{h}{\delta_x\Phi(x)}_{L^2}\mu^\vartheta(\mathd x)=\int\Big[\scal{(H^{\beta}+m^2)h}{x}_{L^2}+\scal{h}{\chi_\vartheta(V'\circ x)}_{L^2}\Big]\Phi(x)\mu^\vartheta(\mathd x),
\end{equation}
for every $\Phi\in\mathcal F C_b^\infty(L^2_{r,k})$ and $h\in C_c^\infty(\R)$.

\item \emph{Invariance.} By point 2. and by Theorem \ref{theorem:infinitesimalsymmetri} applied to $V_\vartheta$, the measure $\mu^\vartheta$ is infinitesimally reversible, and in particular infinitesimally invariant, for \eqref{eq:existenceSPDE}; hence, by \eqref{eq:existencemuthetamoments} and Theorem \ref{thm:main-inv-ii}, $\mu^\vartheta$ is an invariant measure for $\{P^\vartheta(t)\}_{t\geq0}$.
\end{enumerate}

In the rest of the subsection we prove that $\{\mu^\vartheta\}_{\vartheta>0}$ is tight and that every limit point of it solves the DLR equations \eqref{eq:DLRequations} on the whole of $\R$.

Let $y,z\in L^{2}_{r,k}$. Let $\{Z(t,z)\}_{t\geq 0}$ be the solution of \eqref{eq:systemZY2} and, for $x=y+z$, let
\[
Y^{\vartheta}(t,y):=X^{\vartheta}(t,x)-Z(t,z),\qquad t\geq0,
\]
which is the unique generalized mild solution of
\begin{equation}\label{eq:existenceY}
\begin{cases}
\mathd Y^{\vartheta}(t)+\mathcal{A}Y^{\vartheta}(t)\mathd t=\mathcal{B}_\vartheta\big(Y^{\vartheta}(t)+Z(t)\big)\mathd t,\qquad t>0,\\
Y^{\vartheta}(0)=y.
\end{cases}
\end{equation}

We prove an analogues of Theorem \ref{theorem:apriori} and  Lemma \ref{lemma:invariantmoments}, for $\{Y^{\vartheta}(t)\}_{t\geq 0}$ and $\mu^\vartheta$, respectively.
\begin{lemm}\label{lemma:existenceapriori}
Assume \ref{H0}, \ref{H1}, \ref{H2}, \ref{H3}, \ref{H4} and \ref{H5} hold. Let $q\geq 1$ be the constant given by \eqref{H2:pol} and let
\[
\gamma\in\big(0,s_1-\tfrac12\big),\qquad 0<\varepsilon<\frac{r_0-1}{2q}.
\]
Then there exist $\ell_0>0$ and a constant $K_4>0$ such that, for every $\vartheta>0$, $\tau\in\R$, $\ell\in[0,\ell_0)$, $z\in C^{\gamma}_{\rho^\tau_{\varepsilon,k}}$ and $t\geq0$, the following estimates hold $\mathbb{P}$-a.s.
\begin{enumerate}[label=(\roman*)]
\item if $y\in L^2_{r_0,k}$, then
{\footnotesize
\begin{equation}\label{eq:existenceapriori}
\norm{Y^{\vartheta}(t,y)}^2_{L^2_{\rho^\tau_{r_0,k}}}+\int_0^te^{-\ell(t-s)}\Big(\norm{Y^{\vartheta}(s,y)}^2_{H^{s_1}_{\rho^\tau_{r_0,k}}}+\norm{\chi_\vartheta Y^{\vartheta}(s,y)}^{q+1}_{L^{q+1}_{\rho^\tau_{r_0,k}}}\Big)\mathd s\leq K_4\Big(e^{-\ell t}\norm{y}^2_{L^2_{\rho^\tau_{r_0,k}}}+R_\ell(t)\Big),
\end{equation}}
\item if $y\in L^{q+1}_{r_0,k}$, then
\begin{equation}\label{eq:existenceaprioriLp}
\norm{Y^{\vartheta}(t,y)}^{q+1}_{L^{q+1}_{\rho^\tau_{r_0,k}}}+\int_0^te^{-\ell(t-s)}\norm{\chi_\vartheta Y^{\vartheta}(s,y)}^{2q}_{L^{2q}_{\rho^\tau_{r_0,k}}}\mathd s\leq K_4\Big(e^{-\ell t}\norm{y}^{q+1}_{L^{q+1}_{\rho^\tau_{r_0,k}}}+R_\ell(t)\Big),
\end{equation}
\end{enumerate}
where $\{Y^{\vartheta}(t,y)\}_{t\geq0}$ is the solution of \eqref{eq:existenceY}, $\rho^{\tau}_{r_0,k}:=\rho_{r_0,k}(\cdot-\tau)$ and $R_\ell(t)$ is given by \eqref{eq:R_l} with $u=0$.
\end{lemm}

\begin{proof}
The proof is the same as that of Theorem \ref{theorem:apriori}, points (i) and (ii), with $u=0$, with $r_0$ in place of $r$ and with $\mathcal{B}$ replaced by $\mathcal{B}_\vartheta$, so we only indicate the two changes. For simplicity, we set $\rho:=\rho^\tau_{r_0,k}$, $Y:=Y^{\vartheta}(\cdot,y)$, $Z:=Z(\cdot,z)$ and $X:=Y+Z$.

First, $\chi_\vartheta$ does not depend on $t$ and takes values in $\{0,1\}$, so that $\mathcal{B}_\vartheta$ satisfies \eqref{H2:diss} and \eqref{eq:potential-growth}; hence Lemma \ref{lemm:apriori} and Remark \ref{rmk:domainA2} apply to \eqref{eq:existenceY} as well, and this justifies the computations below. Second, the coercivity \eqref{eq:apriori-coercive} and the Kato-type inequality \eqref{eq:Kato-type} concern the operator $\mathcal{A}$ only, and are therefore unaffected by the truncation, while the potential term changes. Indeed, since $\chi_\vartheta^{q+1}=\chi_\vartheta$, multiplying \eqref{eq:potential-coerc} by $\chi_\vartheta$ and integrating we get
\[
\scal{\mathcal{B}_\vartheta(X)}{Y}_{L^2_{r_0,k}}\leq-N_1\norm{\chi_\vartheta Y}^{q+1}_{L^{q+1}_{\rho}}+N_2\Big(\norm{\rho_{r_0,k}}_{L^1(\R)}+\norm{Z}^{q+1}_{L^{q+1}_{\rho}}\Big),
\]
which replaces the corresponding estimate in point (i); analogously, multiplying \eqref{eq:potential-coerc} by $\chi_\vartheta\abs{Y}^{q-1}$ and using the Young inequality with exponents $\frac{2q}{q-1}$ and $\frac{2q}{q+1}$, there exists a constant $S_1>0$ such that
\[
\scal{\mathcal{B}_\vartheta(X)}{\abs{Y}^{q-1}Y}_{L^2_{r_0,k}}\leq-\frac{N_1}{2}\norm{\chi_\vartheta Y}^{2q}_{L^{2q}_{\rho}}+S_1\Big(\norm{\rho_{r_0,k}}_{L^1(\R)}+\norm{Z}^{2q}_{L^{2q}_{\rho}}\Big),
\]
which replaces the corresponding estimate in point (ii). The conclusion follows by the Gronwall inequality as in Theorem \ref{theorem:apriori}, taking into account \eqref{eq:weighted-holder-Lr}; all the constants involved are independent of $\vartheta$ and of $\tau$.
\end{proof}

\begin{rmk}
The presence of $\chi_\vartheta$ in the coercive terms is due to the truncation, indeed the nonlinear drift acts only on $[-\vartheta,\vartheta]$, and therefore \eqref{eq:existenceapriori} and \eqref{eq:existenceaprioriLp} control $Y^{\vartheta}$ in $L^{q+1}$ and in $L^{2q}$ only on that interval. This is the reason for the factor $\chi_\vartheta$ in the last term of \eqref{eq:existencemoments} below.
\end{rmk}

\begin{lemm}\label{lemma:existencemoments}
Assume \ref{H0}, \ref{H1}, \ref{H2}, \ref{H3}, \ref{H4} and \ref{H5} hold. Let $q\geq 1$ be the constant given by \eqref{H2:pol} and let $\gamma'\in\big(0,s_1-\frac12\big)$. Then there exists a constant $K_5>0$ such that, for every $\vartheta>0$,
\begin{equation}\label{eq:existencemoments}
\sup_{\tau\in\R}\int_{L^2_{r_0,k}}\Big(\norm{x}^2_{C^{\gamma'}_{\rho^\tau_{r_0/2,k}}}+\norm{x}^{q+1}_{L^{q+1}_{\rho^\tau_{r_0,k}}}+\norm{\chi_\vartheta x}^{2q}_{L^{2q}_{\rho^\tau_{r_0,k}}}\Big)\mu^\vartheta(\mathd x)\leq K_5 .
\end{equation}
\end{lemm}

\begin{proof}
Fix $\vartheta>0$, $\tau\in\R$ and $\varepsilon\in\big(0,\frac{r_0-1}{2q}\big)$, and set $\rho:=\rho^\tau_{r_0,k}$. The proof is the same as that of Lemma \ref{lemma:invariantmoments}, with $\{P(t)\}_{t\geq0}$ replaced by $\{P^\vartheta(t)\}_{t\geq0}$, with $\mu$ replaced by $\mu^\vartheta$, which is an invariant measure of $\{P^\vartheta(t)\}_{t\geq0}$ by point 3. above, and with Theorem \ref{theorem:apriori} replaced by Lemma \ref{lemma:existenceapriori}. We only indicate the changes.

Taking $z=0$ in \eqref{eq:existenceY}, we have $X^{\vartheta}(s,x)=Y^{\vartheta}(s,x)+W_{\mathcal{A}}(s)$ for every $s\geq0$, and \eqref{eq:invariantmomentsZ} holds with $\gamma'$ in place of $\gamma$ and with $r_0$ in place of $r$. Hence, arguing as in Lemma \ref{lemma:invariantmoments} with \eqref{eq:existenceapriori} in place of \eqref{eq:aprioriL2}, and applying Lemma \ref{lemma:Ttlimit} to $\{P^\vartheta(t)\}_{t\geq0}$ and to $\mu^\vartheta$ with $B=L^2_{r_0,k}$, we obtain a constant $S_2>0$, independent of $\vartheta$ and of $\tau$, such that
\[
\int_{L^2_{r_0,k}}\Big(\norm{x}^2_{C^{\gamma'}_{\rho^\tau_{r_0/2,k}}}+\norm{x}^{q+1}_{L^{q+1}_{\rho}}\Big)\mu^\vartheta(\mathd x)\leq S_2 .
\]
In particular $\mu^\vartheta\big(L^{q+1}_{r_0,k}\big)=1$. Repeating the argument with \eqref{eq:existenceaprioriLp} in place of \eqref{eq:existenceapriori} and with $B=L^{q+1}_{r_0,k}$, and using that
\[
\norm{\chi_\vartheta x}^{2q}_{L^{2q}_{\rho}}\leq2^{2q-1}\Big(\norm{\chi_\vartheta Y^{\vartheta}}^{2q}_{L^{2q}_{\rho}}+\norm{W_{\mathcal{A}}}^{2q}_{L^{2q}_{\rho}}\Big)
\]
together with \eqref{eq:weighted-holder-Lr}, we obtain a constant $S_3>0$, independent of $\vartheta$ and of $\tau$, bounding the last term in \eqref{eq:existencemoments}. This proves \eqref{eq:existencemoments} with $K_5:=S_2+S_3$.
\end{proof}

We are now going to prove the tightness of $\{\mu^\vartheta\}_{\vartheta>0}$ in the space
\[
E^\gamma_{r,k}:=L^2_{r,k}\cap C^\gamma_{\rho_{r/2,k}}(\R),\qquad\norm{x}_{E^\gamma_{r,k}}:=\norm{x}_{L^2_{r,k}}+\norm{x}_{C^\gamma_{\rho_{r/2,k}}} .
\]

\begin{lemm}\label{lemma:existencetightness}
Assume \ref{H0}, \ref{H1}, \ref{H2}, \ref{H3}, \ref{H4} and \ref{H5} hold. Let $q\geq 1$ be the constant given by \eqref{H2:pol}. Then the family $\{\mu^\vartheta\}_{\vartheta>0}$ is tight in $E^\gamma_{r,k}$. Moreover
\begin{equation}\label{eq:existenceuniformmoments}
\sup_{\vartheta>0}\int_{E^\gamma_{r,k}}\norm{x}^2_{E^\gamma_{r,k}}\mu^\vartheta(\mathd x)<+\infty,
\end{equation}
and, for every bounded open interval $\mathcal I$ of $\R$,
\begin{equation}\label{eq:existencelocalmoments}
\sup_{\vartheta:\,\mathcal I\subseteq[-\vartheta,\vartheta]}\int_{E^\gamma_{r,k}}\norm{x}^{2q}_{L^{2q}(\mathcal I)}\mu^\vartheta(\mathd x)<+\infty .
\end{equation}
\end{lemm}

\begin{proof}
For $R>0$ we set
\[
\mathcal K_R:=\Big\{x\in E^\gamma_{r,k}:\norm{x}_{C^{\gamma'}_{\rho_{r_0/2,k}}}+\norm{x}_{L^2_{r_0,k}}\leq R\Big\},
\]
and we claim that $\mathcal K_R$ is compact in $E^\gamma_{r,k}$. Indeed, let $\{x_n\}_{n\in\N}\subseteq\mathcal K_R$. Since $\gamma<\gamma'$, by the Ascoli-Arzel\`a theorem and a diagonal argument there exist a subsequence, still denoted by $\{x_n\}_{n\in\N}$, and $x\in\mathcal K_R$ such that $x_n\rightarrow x$ in $C^{\gamma}([-N,N])$ for every $N\in\N$. Since $r_0<r$, the quotients $\rho_{r/2,k}/\rho_{r_0/2,k}$ and $\rho_{r,k}/\rho_{r_0,k}$ vanish at infinity, so that
\[
\norm{x_n-x}_{C^{\gamma}_{\rho_{r/2,k}}}\leq\norm{x_n-x}_{C^{\gamma}([-N,N])}+2R\sup_{\abs{\xi}\geq N}\frac{\rho_{r/2,k}(\xi)}{\rho_{r_0/2,k}(\xi)},
\]
and analogously for $\norm{x_n-x}_{L^2_{r,k}}$; letting first $n\rightarrow+\infty$ and then $N\rightarrow+\infty$ we get $x_n\rightarrow x$ in $E^\gamma_{r,k}$, which proves the claim.

By the H\"older inequality with exponents $\frac{q+1}{2}$ and $\frac{q+1}{q-1}$ (if $q=1$ the following is an identity), and since $a^{2/(q+1)}\leq 1+a$ for every $a\geq0$, we have
\[
\norm{x}^2_{L^2_{r_0,k}}\leq\norm{\rho_{r_0,k}}_{L^1(\R)}^{\frac{q-1}{q+1}}\norm{x}^2_{L^{q+1}_{r_0,k}}\leq\norm{\rho_{r_0,k}}_{L^1(\R)}^{\frac{q-1}{q+1}}\Big(1+\norm{x}^{q+1}_{L^{q+1}_{r_0,k}}\Big),\qquad x\in L^{q+1}_{r_0,k}.
\]
Hence, since $(a+b)^2\leq2a^2+2b^2$, estimate \eqref{eq:existencemoments} with $\tau=0$ gives
\[
\sup_{\vartheta>0}\int_{L^2_{r_0,k}}\Big(\norm{x}_{C^{\gamma'}_{\rho_{r_0/2,k}}}+\norm{x}_{L^2_{r_0,k}}\Big)^2\mu^\vartheta(\mathd x)\leq 2K_5+2\norm{\rho_{r_0,k}}_{L^1(\R)}^{\frac{q-1}{q+1}}\big(1+K_5\big)<+\infty ,
\]
so that, by the Chebyshev inequality, for every $\eta>0$ there exists $R>0$, independent of $\vartheta$, with $\mu^\vartheta(\mathcal K_R)\geq1-\eta$; this proves the tightness. The same estimate gives \eqref{eq:existenceuniformmoments}, since $\norm{\cdot}_{E^\gamma_{r,k}}\leq S_4\big(\norm{\cdot}_{C^{\gamma'}_{\rho_{r_0/2,k}}}+\norm{\cdot}_{L^2_{r_0,k}}\big)$ for some $S_4>0$. Finally, if $\mathcal I\subseteq[-\vartheta,\vartheta]$, then $\chi_\vartheta=1$ on $\mathcal I$ and $\rho_{r_0,k}\geq\min_{\overline{\mathcal I}}\rho_{r_0,k}>0$ on $\mathcal I$, so that \eqref{eq:existencelocalmoments} follows from the last term in \eqref{eq:existencemoments}.
\end{proof}

\begin{thm}\label{theorem:existence}
Assume \ref{H0}, \ref{H1}, \ref{H2}, \ref{H3}, \ref{H4} and \ref{H5} hold. Let $q\geq 1$ be the constant given by \eqref{H2:pol}. Then there exists an unbounded and increasing sequence  $\{\vartheta_n\}_{n\in\N}\subseteq\R$ and a Borel probability measure $\mu$ such that $\{\mu^{\vartheta_n}\}_{n\in\N}$ converges weakly in $E^\gamma_{r,k}$ to $\mu$. Moreover $\mu$ satisfies the integration by parts formula \eqref{eq:IPF} associated with \eqref{eq:SPDEmain1psuedo}. Consequently $\mu$ is a reversible measure for \eqref{eq:SPDEmain1psuedo} and a solution of the DLR equations \eqref{eq:DLRequations}.
\end{thm}

\begin{proof}
The first part of the statement follows directly by Lemma \ref{lemma:existencetightness}. We now prove that
\begin{equation}\label{eq:existencelimitmoment}
\int_{L^{2q}_{r,k}}\norm{x}^{2q}_{L^{2q}_{r,k}}\mu(\mathd x)<+\infty .
\end{equation}
For every $N\in\N$ the map
\[
x\longmapsto\min\Big\{N,\int_{-N}^{N}\abs{x(\xi)}^{2q}\rho_{r,k}(\xi)\mathd\xi\Big\}
\]
is bounded and continuous on $E^\gamma_{r,k}$, so that, since $\rho_{r,k}\leq\rho_{r_0,k}$ and $\chi_{\vartheta_n}=1$ on $(-N,N)$ whenever $\vartheta_n\geq N$, estimate \eqref{eq:existencemoments} with $\tau=0$ gives
$\sup_{n:\,\vartheta_n\geq N}\int\int_{-N}^N\abs{x(\xi)}^{2q}\rho_{r,k}(\xi)\mathd\xi\,\mu^{\vartheta_n}(\mathd x)\leq K_5$, so that, by the weak convergence,
\[
\int_{E^\gamma_{r,k}}\min\Big\{N,\int_{-N}^{N}\abs{x(\xi)}^{2q}\rho_{r,k}(\xi)\mathd\xi\Big\}\mu(\mathd x)\leq\sup_{n}\int_{E^\gamma_{r,k}}\norm{x}^{2q}_{L^{2q}((-N,N))}\mu^{\vartheta_n}(\mathd x)\leq K_5 .
\]
Letting $N\rightarrow+\infty$ and using the monotone convergence theorem we obtain \eqref{eq:existencelimitmoment}. Since $E^\gamma_{r,k}\subseteq C(\R)$, this gives $\mu\big(L^{2q}_{r,k}\cap C(\R)\big)=1$ and, by the H\"older inequality, $\int\norm{x}^{q}_{L^{2q}_{r,k}}\mu(\mathd x)<+\infty$.

Let now $\Phi\in\mathcal F C_b^\infty(L^2_{r,k})$ and $h\in C_c^\infty(\R)$, and let $\mathcal I$ be a bounded open interval containing $\supp(h)$. For every $n$ such that $\mathcal I\subseteq[-\vartheta_n,\vartheta_n]$ we have $\chi_{\vartheta_n}h=h$, so that \eqref{eq:existenceIPF} reads
\begin{equation}\label{eq:existenceIPFn}
\int\scal{h}{\delta_x\Phi(x)}_{L^2}\mu^{\vartheta_n}(\mathd x)=\int\Big[\scal{(H^{\beta}+m^2)h}{x}_{L^2}+\scal{h}{V'\circ x}_{L^2}\Big]\Phi(x)\mu^{\vartheta_n}(\mathd x).
\end{equation}
The integrand on the left-hand side is bounded and continuous on $E^\gamma_{r,k}$. As for the right-hand side, by \eqref{eq:testdecay-L2} we have $\rho_{r,k}^{-1/2}(H^{\beta}+m^2)h\in L^2$, so that $x\mapsto\scal{(H^{\beta}+m^2)h}{x}_{L^2}$ is continuous on $E^\gamma_{r,k}$ and, by the Cauchy-Schwarz inequality, it is bounded in absolute value by $\norm{\rho_{r,k}^{-1/2}(H^{\beta}+m^2)h}_{L^2}\norm{x}_{L^2_{r,k}}$. Moreover, by \eqref{H2:pol} and the H\"older inequality, $x\mapsto\scal{h}{V'\circ x}_{L^2}$ is continuous on $E^\gamma_{r,k}$ and there exists a constant $S_6>0$ such that
\[
\abs{\scal{h}{V'\circ x}_{L^2}}\leq S_6\Big(1+\norm{x}^{q}_{L^{2q}(\mathcal I)}\Big).
\]
Hence, by \eqref{eq:existenceuniformmoments} and \eqref{eq:existencelocalmoments}, the integrands on the right-hand side of \eqref{eq:existenceIPFn} are bounded in $L^2(\mu^{\vartheta_n})$ uniformly in $n$, and therefore uniformly integrable; passing to the limit in \eqref{eq:existenceIPFn} we obtain \eqref{eq:IPF} for $\mu$.

Finally, by \eqref{eq:existencelimitmoment}, $\mu$ satisfies all the assumptions of Theorem \ref{theorem:main} and, including the moment condition of order $2q$, of Theorem \ref{theorem:infinitesimalsymmetri}. The former gives that $\mu$ solves \eqref{eq:DLRequations}, the latter that $\mu$ is reversible for \eqref{eq:SPDEmain1psuedo}.
\end{proof}

\subsection{Dirichlet boundary conditions}

Throughout this subsection we fix a boundary condition $y_0\in H^{s_1}_{r_0,k}\cap L^\infty$ and, for $\vartheta>0$, we set $\mathcal I_\vartheta:=(-\vartheta,\vartheta)$. We consider the family of measures
\begin{equation}\label{eq:existencemuthetay0}
\mu^{\vartheta,y_0}:=\frac{1}{\mathcal{Z}^{\mathcal I_\vartheta}(y_0|_{\mathcal I_\vartheta^c})}\exp\left(-\int_{-\vartheta}^\vartheta V(x(\xi))\mathd\xi\right)\mu_{\free}^{\mathcal I_\vartheta,y_0|_{\mathcal I_\vartheta^{c}}}=\mu^{\mathcal I_\vartheta,y_0|_{\mathcal I_\vartheta^c}},
\end{equation}
namely the solution of the DLR equations \eqref{eq:DLRequations} on $\mathcal I_\vartheta$ with the boundary condition $y_0$ outside $\mathcal I_\vartheta$; here $\mathcal{Z}^{\mathcal I_\vartheta}$ is the partition function of \eqref{eq:existencecondtheta} and $\mu_{\free}^{\mathcal I_\vartheta,y_0|_{\mathcal I_\vartheta^c}}$ is given by Proposition \ref{proposition:conditionalGaussian}. The assumption $y_0\in L^\infty(\R)$ covers in particular the constant boundary conditions, which belong to $H^{s_1}_{r_0,k}$ since $r_0>1$. It is convenient to consider the related measure
\[
\bar\mu^{\vartheta,y_0}:=\tau_{\varphi^{\vartheta,y_0}}\big(\mu^{\vartheta,y_0}\big)=\frac{1}{\mathcal{Z}^{\mathcal I_\vartheta}(y_0|_{\mathcal I_\vartheta^c})}\exp\left(-\int_{-\vartheta}^\vartheta V\big(x(\xi)+\varphi^{\vartheta,y_0}(\xi)\big)\mathd\xi\right)\mu_{\free}^{\mathcal I_\vartheta,0},
\]
where $\varphi^{\vartheta,y_0}$ is the harmonic extension of $y_0|_{\mathcal I_\vartheta^c}$ in $\mathcal I_\vartheta$, which is well defined by Theorem \ref{thm:harmonic-continuity} since $s_1>\frac12$ and $r_0>1$ imply $H^{s_1}_{r_0,k}\subseteq C(\R)\cap L^1_{\rho_{r_0,k}}(\R)$. We first prove that the influence of the boundary condition disappears as $\vartheta\rightarrow+\infty$.

\begin{lemm}\label{lemma:boundednessvarphitheta}
Assume \ref{H0}, \ref{H1} and \ref{H3} hold and let $y_0\in H^{s_1}_{r_0,k}\cap L^\infty$. Then there exists a constant $K_6>0$ such that
\begin{equation}\label{eq:boundednessvarphitheta}
\sup_{\vartheta>0}\norm{\varphi^{\vartheta,y_0}}_{H^{s_1}_{r_0,k}}\leq K_6\norm{y_0}_{H^{s_1}_{r_0,k}} .
\end{equation}
Moreover, for every $\eta\in(0,s_1)$, we have $\varphi^{\vartheta,y_0}\rightarrow0$ in $H^{s_1-\eta}_{\rho_{r_0+\eta,k}}$ as $\vartheta\rightarrow+\infty$.
\end{lemm}

\begin{proof}
In this proof we use the notation and the concepts of Section \ref{section:harmonicextension}, and we set $\rho:=\rho_{r_0,k}$ and $w^\vartheta:=y_0-\varphi^{\vartheta,y_0}$. By the discussion in the proof of Lemma \ref{lemma:harmonicextensionenergy} we have $w^\vartheta\in\mathcal{F}_{\mathcal{I}_\vartheta}$, and by \eqref{eq:lax}
\begin{equation}\label{eq:variational-wtheta}
\mathcal E(w^\vartheta,h)=\mathcal E(y_0,h),\qquad h\in\mathcal F_{\mathcal I_\vartheta}.
\end{equation}
Indeed, \eqref{eq:lax} is proved in Lemma \ref{lemma:harmonicextensionenergy} for $y_0\in H^{s_1}(\R)$, and it extends first to compactly supported smooth approximations of $y_0$ and then, by density, to $H^{s_1}_{r_0,k}$, using Theorem \ref{thm:harmonic-continuity} to identify the resulting harmonic extension. Here $\mathcal E$ is the quadratic form defined by \eqref{eq:dirichletform}, $\mathcal F$ is its domain and, as in the proof of Lemma \ref{lemma:harmonicextensionenergy}, $\mathcal F_{\mathcal I_\vartheta}$ is the closure in $\mathcal F$ of the functions of $\mathcal F$ with support contained in $\mathcal I_\vartheta$.

By Lemma \ref{lemma:coercivity-H3}, applied with $r_0$ in place of $r$ (recall that $k<k_0(r_0,m)$), we have
\begin{equation}\label{eq:weighted-energy-coercive}
\mathcal E(f,\rho f)\geq C_4\norm{f}^2_{H^{s_1}_{r_0,k}},\qquad\abs{\mathcal E(f,\rho g)}\leq C_5\norm{f}_{H^{s_1}_{r_0,k}}\norm{g}_{H^{s_1}_{r_0,k}},\qquad f,g\in H^{s_1}_{r_0,k},
\end{equation}
where $C_4,C_5>0$ do not depend on the supports of $f$ and $g$, and hence on $\vartheta$. Since the multiplication by $\rho$ maps $\mathcal F_{\mathcal I_\vartheta}$ into itself, we may take $h=\rho w^\vartheta$ in \eqref{eq:variational-wtheta}, so that \eqref{eq:weighted-energy-coercive} gives
\[
C_4\norm{w^\vartheta}^2_{H^{s_1}_{r_0,k}}\leq\mathcal E(w^\vartheta,\rho w^\vartheta)=\mathcal E(y_0,\rho w^\vartheta)\leq C_5\norm{y_0}_{H^{s_1}_{r_0,k}}\norm{w^\vartheta}_{H^{s_1}_{r_0,k}} .
\]
Hence $\sup_{\vartheta>0}\norm{w^\vartheta}_{H^{s_1}_{r_0,k}}\leq\frac{C_5}{C_4}\norm{y_0}_{H^{s_1}_{r_0,k}}$, and \eqref{eq:boundednessvarphitheta} follows from $\varphi^{\vartheta,y_0}=y_0-w^\vartheta$ with $K_6:=1+\frac{C_5}{C_4}$.

We now prove the weak convergence to zero. Let $\{\vartheta_n\}_{n\in\N}$ be an unbounded and increasing sequence. By \eqref{eq:boundednessvarphitheta}, after passing to a subsequence there exists $\varphi\in H^{s_1}_{r_0,k}$ such that $\varphi^{\vartheta_n,y_0}$ converges weakly to $\varphi$ in $H^{s_1}_{r_0,k}$. Let $h\in C^\infty_c(\R)$; for every $n$ large enough we have $\supp(h)\subseteq\mathcal I_{\vartheta_n}$, so that point 3. of Definition \ref{def:harmonicextension} gives
\begin{equation}\label{eq:varphithetany0}
\int_\R\varphi^{\vartheta_n,y_0}(\xi)\big[(H^{\beta}+m^2)h\big](\xi)\mathd\xi=0 .
\end{equation}
By \eqref{eq:testdecay-L2} we have $\rho^{-1/2}(H^{\beta}+m^2)h\in L^2$, so that
\[
u\longmapsto\int_\R u\,(H^{\beta}+m^2)h\,\mathd\xi=\scal{\rho^{1/2}u}{\rho^{-1/2}(H^{\beta}+m^2)h}_{L^2}
\]
is a continuous linear functional on $H^{s_1}_{r_0,k}$, and we may pass to the limit in \eqref{eq:varphithetany0}, obtaining that $\int_\R\varphi\,(H^{\beta}+m^2)h\,\mathd\xi=0$ for every $h\in C^\infty_c(\R)$. For $g\in C^\infty_c(\R)$ we can take $h=\rho_{r_0,k}g$, and, since both sides are continuous with respect to $\varphi\in H^{s_1}_{r_0,k}$ (the left-hand side by \eqref{eq:coercivity-H3} and the right-hand side by Lemma \ref{lemma:testdecay}, since $r_0<1+2s_2$), the identity $\mathcal E(f,\rho_{r_0,k}g)=\int_\R f\,(H^{\beta}+m^2)(\rho_{r_0,k}g)\,\mathd\xi$, valid for $f\in C^\infty_c(\R)$, extends to $f=\varphi$. Hence $\mathcal E(\varphi,\rho_{r_0,k}g)=0$ for every $g\in C^\infty_c(\R)$, and by density and \eqref{eq:coercivity-H3} we get $C_4\norm{\varphi}^2_{H^{s_1}_{r_0,k}}\leq\mathcal E(\varphi,\rho_{r_0,k}\varphi)=0$, namely $\varphi=0$. Since the same argument applies to every subsequence, we conclude that $\varphi^{\vartheta,y_0}$ converges weakly to $0$ in $H^{s_1}_{r_0,k}$, and the assertion follows since $H^{s_1}_{r_0,k}$ is compactly embedded in $H^{s_1-\eta}_{\rho_{r_0+\eta,k}}$.
\end{proof}

From now on we fix $\eta\in\big(0,s_1-\gamma-\frac12\big)$ such that $r_0+\eta<r$. Since $H^{s_1-\eta}_{\rho_{r_0+\eta,k}}$ is continuously embedded in $E^{\gamma}_{r,k}$, Lemma \ref{lemma:boundednessvarphitheta} gives $\varphi^{\vartheta,y_0}\rightarrow0$ in $E^{\gamma}_{r,k}$ as $\vartheta\rightarrow+\infty$. Hence $\{\mu^{\vartheta,y_0}\}_{\vartheta>0}$ is tight in $E^{\gamma}_{r,k}$ if and only if $\{\bar\mu^{\vartheta,y_0}\}_{\vartheta>0}$ is tight in $E^{\gamma}_{r,k}$, and, given an unbounded and increasing sequence $\{\vartheta_n\}_{n\in\N}$ and a Borel probability measure $\mu$, the sequence $\{\mu^{\vartheta_n,y_0}\}_{n\in\N}$ converges weakly to $\mu$ in $E^{\gamma}_{r,k}$ if and only if $\{\bar\mu^{\vartheta_n,y_0}\}_{n\in\N}$ converges weakly to $\mu$ in $E^{\gamma}_{r,k}$. We may therefore work with $\bar\mu^{\vartheta,y_0}$.

By Proposition \ref{prop:invmeasureI}, applied with $\mathcal I=\mathcal I_\vartheta$ and $y=y_0|_{\mathcal I_\vartheta^c}$, the measure $\bar\mu^{\vartheta,y_0}=\nu_{\mathcal I_\vartheta,y_0|_{\mathcal I_\vartheta^c}}$ is the unique invariant measure of the SPDE \eqref{eq:SPDEL2(I)} on $L^2(\mathcal I_\vartheta)$, namely of
\begin{equation}\label{eq:SPDEDiry0}
\begin{cases}
\mathd X^{\vartheta,y_0}(t)+\mathcal{A}_{\mathcal I_\vartheta}X^{\vartheta,y_0}(t)\mathd t=\mathcal{B}\big(X^{\vartheta,y_0}(t)+\varphi^{\vartheta,y_0}\big)\mathd t+\sqrt2\mathd W^{\vartheta}(t),\qquad t>0,\\
X^{\vartheta,y_0}(0)=x,
\end{cases}
\end{equation}
where $\mathcal{A}_{\mathcal I_\vartheta}$ is the realization of $H^{\beta}+m^2$ in $L^2(\mathcal I_\vartheta)$ with homogeneous Dirichlet boundary conditions and $\{W^{\vartheta}(t)\}_{t\geq0}$ is an $L^2(\mathcal I_\vartheta)$-cylindrical Wiener process. We denote by $\{P^{\vartheta,y_0}(t)\}_{t\geq0}$ the associated transition semigroup on $L^2(\mathcal I_\vartheta)$. Moreover, by \eqref{eq:ibpIbis}, the measure $\bar\mu^{\vartheta,y_0}$ solves the integration by parts formula
\begin{equation}\label{eq:IPFDiry0}
\int\scal{h}{\delta_x\Phi(x)}_{L^2}\bar\mu^{\vartheta,y_0}(\mathd x)=\int\Big[\scal{(H^{\beta}+m^2)h}{x}_{L^2}+\scal{h}{V'\circ\big(x+\varphi^{\vartheta,y_0}\big)}_{L^2}\Big]\Phi(x)\bar\mu^{\vartheta,y_0}(\mathd x),
\end{equation}
for every $\Phi\in\mathcal F C_b^\infty(L^2_{r,k})$ and $h\in C_c^\infty(\mathcal I_\vartheta)$.

Let $y,z\in L^2(\mathcal I_\vartheta)$. Let $\{Z^{\vartheta}(t,z)\}_{t\geq0}$ be the solution of \eqref{eq:SPDEL2(I)} with $V=0$, that is $Z^{\vartheta}(t,z)=e^{-t\mathcal{A}_{\mathcal I_\vartheta}}z+W_{\mathcal{A}_{\mathcal I_\vartheta}}(t)$, where $W_{\mathcal{A}_{\mathcal I_\vartheta}}$ is the stochastic convolution \eqref{eq:stochconvI}, and, for $x=y+z$, let
\[
Y^{\vartheta,y_0}(t,y):=X^{\vartheta,y_0}(t,x)-Z^{\vartheta}(t,z),\qquad t\geq0,
\]
which is the unique generalized mild solution of
\begin{equation}\label{eq:existenceYDir}
\begin{cases}
\mathd Y^{\vartheta,y_0}(t)+\mathcal{A}_{\mathcal I_\vartheta}Y^{\vartheta,y_0}(t)\mathd t=\mathcal{B}\big(Y^{\vartheta,y_0}(t)+Z^{\vartheta}(t)+\varphi^{\vartheta,y_0}\big)\mathd t,\qquad t>0,\\
Y^{\vartheta,y_0}(0)=y.
\end{cases}
\end{equation}
We prove the analogue of Lemma \ref{lemma:existenceapriori}, where the functions in $L^2(\mathcal I_\vartheta)$ are identified with the functions in $L^2$ vanishing outside $\mathcal I_\vartheta$.

\begin{lemm}\label{lemma:existenceaprioriDir}
Assume \ref{H0}, \ref{H1}, \ref{H2}, \ref{H3}, \ref{H4} and \ref{H5} hold. Let $q\geq 1$ be the constant given by \eqref{H2:pol} and let
\[
\gamma\in\big(0,s_1-\tfrac12\big),\qquad 0<\varepsilon<\frac{r_0-1}{2q}.
\]
Then there exist $\ell_0>0$ and a constant $K_7>0$ such that, for every $\vartheta>0$, $\tau\in\R$, $\ell\in[0,\ell_0)$, $z\in L^2(\mathcal I_\vartheta)$ and $t\geq0$, the following estimates hold $\mathbb{P}$-a.s.
\begin{enumerate}[label=(\roman*)]
\item if $y\in L^2(\mathcal I_\vartheta)$, then
{\footnotesize
\begin{equation}\label{eq:existenceaprioriDir}
\norm{Y^{\vartheta,y_0}(t,y)}^2_{L^2_{\rho^\tau_{r_0,k}}}+\int_0^te^{-\ell(t-s)}\Big(\norm{Y^{\vartheta,y_0}(s,y)}^2_{H^{s_1}_{\rho^\tau_{r_0,k}}}+\norm{Y^{\vartheta,y_0}(s,y)}^{q+1}_{L^{q+1}_{\rho^\tau_{r_0,k}}}\Big)\mathd s\leq K_7\Big(e^{-\ell t}\norm{y}^2_{L^2_{\rho^\tau_{r_0,k}}}+\widetilde R_\ell(t)\Big),
\end{equation}}
\item if $y\in L^{q+1}(\mathcal I_\vartheta)$, then
\begin{equation}\label{eq:existenceaprioriDirLp}
\norm{Y^{\vartheta,y_0}(t,y)}^{q+1}_{L^{q+1}_{\rho^\tau_{r_0,k}}}+\int_0^te^{-\ell(t-s)}\norm{Y^{\vartheta,y_0}(s,y)}^{2q}_{L^{2q}_{\rho^\tau_{r_0,k}}}\mathd s\leq K_7\Big(e^{-\ell t}\norm{y}^{q+1}_{L^{q+1}_{\rho^\tau_{r_0,k}}}+\widetilde R_\ell(t)\Big),
\end{equation}
\end{enumerate}
where $\{Y^{\vartheta,y_0}(t,y)\}_{t\geq0}$ is the solution of \eqref{eq:existenceYDir}, $\rho^{\tau}_{r_0,k}:=\rho_{r_0,k}(\cdot-\tau)$ and
\[
\widetilde R_\ell(t):=\int_0^te^{-\ell(t-s)}\Big(1+\norm{Z^{\vartheta}(s,z)}^{2q}_{C^\gamma_{\rho^\tau_{\varepsilon,k}}}+\norm{\varphi^{\vartheta,y_0}}^{2q}_{L^{2q}_{\rho^\tau_{r_0,k}}}\Big)\mathd s .
\]
\end{lemm}

\begin{proof}
The proof is the same as that of Lemma \ref{lemma:existenceapriori}, with two changes. For simplicity, we set $\rho:=\rho^\tau_{r_0,k}$, $Y:=Y^{\vartheta,y_0}(\cdot,y)$, $Z:=Z^{\vartheta}(\cdot,z)$ and $X:=Y+Z+\varphi^{\vartheta,y_0}$.

First, the computations are justified as in Lemma \ref{lemm:apriori}, since $\mathcal{A}_{\mathcal I_\vartheta}$ is self-adjoint and nonnegative in $L^2(\mathcal I_\vartheta)$, by maximal $L^2$-regularity, for $y\in\mathcal F_{\mathcal I_\vartheta}\cap L^{2q}(\mathcal I_\vartheta)$ we have $Y\in H^1(0,T;L^2(\mathcal I_\vartheta))\cap L^2(0,T;{\rm Dom}(\mathcal{A}_{\mathcal I_\vartheta}))$, and the general case follows by approximation. Since $Y(t)\in\mathcal F_{\mathcal I_\vartheta}$ and the multiplication by $\rho$ maps $\mathcal F_{\mathcal I_\vartheta}$ into itself, we have $\scal{\mathcal{A}_{\mathcal I_\vartheta}Y}{\rho Y}_{L^2(\mathcal I_\vartheta)}=\mathcal E(Y,\rho Y)$, and therefore Lemma \ref{lemma:coercivity-H3}, applied with $r_0$ in place of $r$ and with the translated weight $\rho=\rho^\tau_{r_0,k}$ (see Remark \ref{rmk:weighttranslation}), gives
\[
\scal{\mathcal{A}_{\mathcal I_\vartheta}Y}{\rho Y}_{L^2(\mathcal I_\vartheta)}=\mathcal E(Y,\rho Y)\geq C_4\norm{Y}^2_{H^{s_1}_{\rho}},
\]
while the proof of the Kato-type inequality \eqref{eq:Kato-type}, which only uses the convexity of $r\mapsto\abs{r}^{q+1}$, gives $\scal{\mathcal{A}_{\mathcal I_\vartheta}Y}{\rho\abs{Y}^{q-1}Y}_{L^2(\mathcal I_\vartheta)}\geq\big(m^2-\frac{\omega_{r_0,k}}{q+1}\big)\norm{Y}^{q+1}_{L^{q+1}_\rho}$; all the constants are independent of $\vartheta$ and $\tau$, since they depend neither on the supports of the functions involved nor on translations of the weight. Second, the truncation is now absent, while the boundary condition enters through \eqref{eq:potential-coerc} with $b=Z+\varphi^{\vartheta,y_0}$, which gives
\[
\scal{\mathcal{B}(X)}{Y}_{L^2_{\rho}}\leq-N_1\norm{Y}^{q+1}_{L^{q+1}_{\rho}}+2^{q}N_2\Big(\norm{\rho_{r_0,k}}_{L^1(\R)}+\norm{Z}^{q+1}_{L^{q+1}_{\rho}}+\norm{\varphi^{\vartheta,y_0}}^{q+1}_{L^{q+1}_{\rho}}\Big),
\]
and analogously, multiplying \eqref{eq:potential-coerc} by $\abs{Y}^{q-1}$ and using the Young inequality with exponents $\frac{2q}{q-1}$ and $\frac{2q}{q+1}$, there exists a constant $S_7>0$ such that
\[
\scal{\mathcal{B}(X)}{\abs{Y}^{q-1}Y}_{L^2_{\rho}}\leq-\frac{N_1}{2}\norm{Y}^{2q}_{L^{2q}_{\rho}}+S_7\Big(\norm{\rho_{r_0,k}}_{L^1(\R)}+\norm{Z}^{2q}_{L^{2q}_{\rho}}+\norm{\varphi^{\vartheta,y_0}}^{2q}_{L^{2q}_{\rho}}\Big).
\]
The last terms produce the last term of $\widetilde R_\ell$, since $\norm{\varphi^{\vartheta,y_0}}^{q+1}_{L^{q+1}_{\rho}}\leq\norm{\rho_{r_0,k}}_{L^1(\R)}+\norm{\varphi^{\vartheta,y_0}}^{2q}_{L^{2q}_{\rho}}$. The conclusion follows by the Gronwall inequality as in Lemma \ref{lemma:existenceapriori}, taking into account \eqref{eq:weighted-holder-Lr}; all the constants involved are independent of $\vartheta$ and of $\tau$.
\end{proof}

We need the following lemma about the stochastic convolution in this case.

\begin{lemm}\label{lemma:stochconvDirunif}
Assume \ref{H0}, \ref{H1} and \ref{H3} hold, and let $W_{\mathcal A_{\mathcal I_\vartheta}}$ be the stochastic convolution \eqref{eq:stochconvI} with $\mathcal I=\mathcal I_\vartheta$, extended by zero outside $\mathcal I_\vartheta$. Then for every $\gamma'\in\big(0,s_1-\frac12\big)$, $a>0$ and $n\in\N$ there exists $S=S(\gamma',a,n)>0$ such that
\begin{equation}\label{eq:stochconvDirunif}
\sup_{\vartheta>0}\,\sup_{t\geq0}\,\sup_{\tau\in\R}\,\mathbb E\Big[\norm{W_{\mathcal A_{\mathcal I_\vartheta}}(t)}^n_{C^{\gamma'}_{\rho^\tau_{a,k}}}\Big]\leq S .
\end{equation}
\end{lemm}

\begin{proof}
Fix $\vartheta>0$, $t\geq0$, and set $\mathcal I:=\mathcal I_\vartheta$ and $W_\vartheta:=W_{\mathcal A_{\mathcal I}}(t)$. We first compare covariances. The random variable $W_\vartheta$ is a centered Gaussian random variable on $L^2(\mathcal I)$ with covariance operator
\[
Q_t:=2\int_0^te^{-2s\mathcal A_{\mathcal I}}\mathd s=\mathcal A_{\mathcal I}^{-1}\big(\Id-e^{-2t\mathcal A_{\mathcal I}}\big)\leq\mathcal A_{\mathcal I}^{-1},
\]
the inequality being in the sense of nonnegative self-adjoint operators on $L^2(\mathcal I)$. By Proposition \ref{prop:gaussL2I}, $\mathcal A_{\mathcal I}^{-1}$ is the covariance of $\mu^{\mathcal I,0}_{\free}$. Moreover, let $X$ be the free field with law $\mu_{\free}$ and set $X_{\mathcal I}:=X-\varphi^{R_{\mathcal I^c}X}$. By Corollary \ref{corollary:Guassiantranslation}, the conditional law of $X_{\mathcal I}$ given $\sigma(\mathcal I^c)$ is $\mu^{\mathcal I,0}_{\free}$ for almost every boundary condition. Hence $X_{\mathcal I}$ has law $\mu^{\mathcal I,0}_{\free}$ and is independent of the $\sigma(\mathcal I^c)$-measurable field $\varphi^{R_{\mathcal I^c}X}$. Therefore, for every $f\in C^\infty_c(\R)$, since $W_\vartheta$ and $X_{\mathcal I}$ vanish outside $\mathcal I$,
\[
\mathbb E\big[\scal{W_\vartheta}{f}_{L^2}^2\big]=\scal{Q_tf\mathbb I_{\mathcal I}}{f\mathbb I_{\mathcal I}}_{L^2(\mathcal I)}\leq\mathbb E\big[\scal{X_{\mathcal I}}{f}_{L^2}^2\big]\leq\mathbb E\big[\scal{X_{\mathcal I}}{f}_{L^2}^2\big]+\mathbb E\big[\scal{\varphi^{R_{\mathcal I^c}X}}{f}_{L^2}^2\big]=\mathbb E\big[\scal{X}{f}_{L^2}^2\big].
\]
The fields $W_\vartheta$ (by Lemma \ref{lemma:stochconvL2I}, and since the killed semigroup acts on $C_0(\mathcal I)$) and $X$ have continuous versions. Taking $f=\eta_j(\cdot-\xi)-\eta_j(\cdot-\xi')$, or $f=\eta_j(\cdot-\xi)$, with $\{\eta_j\}_{j\in\N}$ a sequence of mollifiers, and letting $j\rightarrow+\infty$ (for centered Gaussian random variables almost sure convergence implies convergence of the variances), we obtain for every $\xi,\xi'\in\R$
\begin{equation}\label{eq:Dirvariancecomparison}
\mathbb E\big[\abs{W_\vartheta(\xi)}^2\big]\leq G_{\free}(0),\qquad
\mathbb E\big[\abs{W_\vartheta(\xi)-W_\vartheta(\xi')}^2\big]\leq\mathbb E\big[\abs{X(\xi)-X(\xi')}^2\big]=\frac1\pi\int_\R\frac{1-\cos\big((\xi-\xi')\eta\big)}{m^2+\widehat{H^\beta}(\eta)}\mathd\eta .
\end{equation}
Fix $\theta\in\big(\gamma',s_1-\frac12\big)$. Using $1-\cos(b)\leq C_\theta\min\{1,\abs{b}^{2\theta}\}$ and \eqref{eq:symbol-lower}, and since $2\theta-2s_1<-1$, the last integral is bounded by $C\min\{1,\abs{\xi-\xi'}^{2\theta}\}$, with $C$ independent of $\vartheta$, $t$, $\xi$ and $\xi'$.

The bounds \eqref{eq:Dirvariancecomparison} are exactly the two properties of $Y(s,\cdot)$ used in the proof of Lemma \ref{lemma:stochconv}: the pointwise variance bound, with $G_{\free}(0)$ in place of $v_\alpha(s)$, and the increment bound, with $C\min\{1,\abs{z}^{2\theta}\}$ in place of \eqref{eq:incrementbound}. That proof, from \eqref{eq:Lpmoment} to \eqref{eq:gaussianmoments}, uses no other property of $Y(s,\cdot)$. Fix $\varepsilon\in(\gamma',\theta)$ and $p\geq2$ such that $\varepsilon-\frac1p\geq\gamma'$ and $\frac{r_0}{p}\leq a$. Repeating it with the weight $\rho^\tau_{r_0,k}$, whose $L^1$ norm $c_{r_0}k^{-1/2}$ does not depend on $\tau$, we get, for every $n\in\N$,
\[
\sup_{\vartheta>0}\sup_{t\geq0}\sup_{\tau\in\R}\mathbb E\Big[\norm{(\rho^\tau_{r_0,k})^{1/p}\,W_\vartheta}^n_{W^{\varepsilon,p}(\R)}\Big]<+\infty .
\]
Finally, $(\rho^\tau_{r_0,k})^{1/p}=\rho^\tau_{r_0/p,k}$. By the embedding $W^{\varepsilon,p}(\R)\hookrightarrow C^{0,\varepsilon-1/p}(\R)\hookrightarrow C^{0,\gamma'}(\R)$ and by the argument of Remark \ref{rmk:stochconvholder} (multiplication by the bounded Lipschitz function $\rho^\tau_{a-r_0/p,k}$, whose Lipschitz constant does not depend on $\tau$), we obtain \eqref{eq:stochconvDirunif}.
\end{proof}

\begin{lemm}\label{lemma:existencemomentsDir}
Assume \ref{H0}, \ref{H1}, \ref{H2}, \ref{H3}, \ref{H4} and \ref{H5} hold. Let $q\geq 1$ be the constant given by \eqref{H2:pol}, let $\gamma'\in\big(0,s_1-\frac12\big)$ and let $y_0\in H^{s_1}_{r_0,k}\cap L^\infty$. Then there exists a constant $K_8>0$ such that, for every $\vartheta>0$,
\begin{equation}\label{eq:existencemomentsDir}
\sup_{\tau\in\R}\int_{L^2(\mathcal I_\vartheta)}\Big(\norm{x}^2_{C^{\gamma'}_{\rho^\tau_{r_0/2,k}}}+\norm{x}^{q+1}_{L^{q+1}_{\rho^\tau_{r_0,k}}}+\norm{x}^{2q}_{L^{2q}_{\rho^\tau_{r_0,k}}}\Big)\bar\mu^{\vartheta,y_0}(\mathd x)\leq K_8 .
\end{equation}
\end{lemm}

\begin{proof}
Fix $\vartheta>0$, $\tau\in\R$ and $\varepsilon\in\big(0,\frac{r_0-1}{2q}\big)$. The proof is the same as that of Lemma \ref{lemma:existencemoments}, with $\{P^{\vartheta}(t)\}_{t\geq0}$ replaced by $\{P^{\vartheta,y_0}(t)\}_{t\geq0}$, with $\mu^\vartheta$ replaced by $\bar\mu^{\vartheta,y_0}$, which is an invariant measure of $\{P^{\vartheta,y_0}(t)\}_{t\geq0}$ by Proposition \ref{prop:invmeasureI}, and with Lemma \ref{lemma:existenceapriori} replaced by Lemma \ref{lemma:existenceaprioriDir}. We only indicate the changes.
The last term of $\widetilde R_\ell$ is bounded uniformly in $\vartheta$ and $\tau$: indeed, by \eqref{eq:harmonicextensionprob} we have $\abs{\varphi^{\vartheta,y_0}(\xi)}\leq\norm{y_0}_{L^\infty(\R)}$ for every $\xi\in\R$, so that, by \eqref{eq:rhoL1},
\[
\sup_{\vartheta>0}\sup_{\tau\in\R}\norm{\varphi^{\vartheta,y_0}}^{2q}_{L^{2q}_{\rho^\tau_{r_0,k}}}\leq c_{r_0}k^{-1/2}\norm{y_0}^{2q}_{L^\infty(\R)}=:S_8 .
\]
Moreover, by Lemma \ref{lemma:stochconvDirunif}, applied with $(a,n)=(r_0/2,2)$ and $(a,n)=(\varepsilon,2q)$, the estimate \eqref{eq:invariantmomentsZ} holds for $W_{\mathcal{A}_{\mathcal I_\vartheta}}$ in place of $W_{\mathcal{A}}$, with $\gamma'$ in place of $\gamma$, with $r_0$ in place of $r$ and with a constant independent of $\vartheta$ and $\tau$. Hence, taking $z=0$ in \eqref{eq:existenceYDir} and arguing as in Lemma \ref{lemma:existencemoments}, first with \eqref{eq:existenceaprioriDir} and $B=L^2(\mathcal I_\vartheta)$ and then with \eqref{eq:existenceaprioriDirLp} and $B=L^{q+1}(\mathcal I_\vartheta)$, we obtain \eqref{eq:existencemomentsDir}.
\end{proof}

Since \eqref{eq:existencemomentsDir} is the analogues of \eqref{eq:existencemoments}, with $\bar\mu^{\vartheta,y_0}$ in place of $\mu^\vartheta$ and without the factor $\chi_\vartheta$, the proof of Lemma \ref{lemma:existencetightness} applies verbatim and gives the following result.

\begin{lemm}\label{lemma:existencetightnessDir}
Assume \ref{H0}, \ref{H1}, \ref{H2}, \ref{H3}, \ref{H4} and \ref{H5} hold. Let $q\geq 1$ be the constant given by \eqref{H2:pol} and let $y_0\in H^{s_1}_{r_0,k}\cap L^\infty$. Then the families $\{\bar\mu^{\vartheta,y_0}\}_{\vartheta>0}$ and $\{\mu^{\vartheta,y_0}\}_{\vartheta>0}$ are tight in $E^\gamma_{r,k}$. Moreover
\begin{equation}\label{eq:existenceuniformmomentsDir}
\sup_{\vartheta>0}\int_{E^\gamma_{r,k}}\norm{x}^2_{E^\gamma_{r,k}}\bar\mu^{\vartheta,y_0}(\mathd x)<+\infty,
\end{equation}
and, for every bounded open interval $\mathcal I$ of $\R$,
\begin{equation}\label{eq:existencelocalmomentsDir}
\sup_{\vartheta:\,\mathcal I\subseteq[-\vartheta,\vartheta]}\int_{E^\gamma_{r,k}}\norm{x}^{2q}_{L^{2q}(\mathcal I)}\bar\mu^{\vartheta,y_0}(\mathd x)<+\infty .
\end{equation}
\end{lemm}

\begin{thm}\label{theorem:existence2}
Assume \ref{H0}, \ref{H1}, \ref{H2}, \ref{H3}, \ref{H4} and \ref{H5} hold. Let $q\geq 1$ be the constant given by \eqref{H2:pol} and let $y_0\in H^{s_1}_{r_0,k}\cap L^\infty$. Then there exists an unbounded and increasing sequence $\{\vartheta_n\}_{n\in\N}\subseteq\R$ and a Borel probability measure $\mu$ such that $\{\mu^{\vartheta_n,y_0}\}_{n\in\N}$ converges weakly in $E^\gamma_{r,k}$ to $\mu$. Moreover $\mu$ satisfies the integration by parts formula \eqref{eq:IPF} associated with \eqref{eq:SPDEmain1psuedo}. Consequently $\mu$ is a reversible measure for \eqref{eq:SPDEmain1psuedo} and a solution of the DLR equations \eqref{eq:DLRequations}.
\end{thm}

\begin{proof}
The first part of the statement follows directly by Lemma \ref{lemma:existencetightnessDir}, and, as observed above, $\{\bar\mu^{\vartheta_n,y_0}\}_{n\in\N}$ converges weakly to $\mu$ in $E^\gamma_{r,k}$ as well. The proof of \eqref{eq:existencelimitmoment} is the same as in Theorem \ref{theorem:existence}: writing $x=\bar x+\varphi^{\vartheta_n,y_0}$ with $\bar x$ distributed according to $\bar\mu^{\vartheta_n,y_0}$, and using $\rho_{r,k}\leq\rho_{r_0,k}$, \eqref{eq:existencemomentsDir} with $\tau=0$ and $\abs{\varphi^{\vartheta_n,y_0}}\leq\norm{y_0}_{L^\infty(\R)}$, the bound $K_5$ is replaced by $2^{2q-1}\big(K_8+\norm{y_0}^{2q}_{L^\infty(\R)}\norm{\rho_{r,k}}_{L^1(\R)}\big)$.

Let now $\Phi\in\mathcal F C_b^\infty(L^2_{r,k})$ and $h\in C_c^\infty(\R)$, and let $\mathcal I$ be a bounded open interval containing $\supp(h)$. For every $n$ such that $\mathcal I\subseteq\mathcal I_{\vartheta_n}$ we have $h\in C^\infty_c(\mathcal I_{\vartheta_n})$, so that \eqref{eq:IPFDiry0} holds for $\bar\mu^{\vartheta_n,y_0}$. The integrands are estimated exactly as in the proof of Theorem \ref{theorem:existence}: the only difference is the term $\scal{h}{V'\circ(x+\varphi^{\vartheta_n,y_0})}_{L^2}$, for which \eqref{H2:pol} and the H\"older inequality give a constant $S_9>0$ such that
\[
\abs{\scal{h}{V'\circ\big(x+\varphi^{\vartheta_n,y_0}\big)}_{L^2}}\leq S_9\Big(1+\norm{x}^{q}_{L^{2q}(\mathcal I)}+\norm{\varphi^{\vartheta_n,y_0}}^{q}_{L^{2q}(\mathcal I)}\Big),
\]
which is bounded in $L^2(\bar\mu^{\vartheta_n,y_0})$ uniformly in $n$ by \eqref{eq:existencelocalmomentsDir} and by \eqref{eq:boundednessvarphitheta}. Moreover $\varphi^{\vartheta_n,y_0}\rightarrow0$ in $E^\gamma_{r,k}$, and hence in $L^{2q}(\mathcal I)$, by Lemma \ref{lemma:boundednessvarphitheta}; therefore, passing to the limit in \eqref{eq:IPFDiry0} we obtain \eqref{eq:IPF} for $\mu$.

Finally, by \eqref{eq:existencelimitmoment}, $\mu$ satisfies all the assumptions of Theorem \ref{theorem:main} and, including the moment condition of order $2q$, of Theorem \ref{theorem:infinitesimalsymmetri}. The former gives that $\mu$ solves \eqref{eq:DLRequations}, the latter that $\mu$ is reversible for \eqref{eq:SPDEmain1psuedo}.
\end{proof}

\subsection{Existence of a reflection-positive solution of the DLR equations}\label{section:ref.positive}

Thanks to the result of Section \ref{appendix:reflectionpositivity} we are able to prove the following theorem.

\begin{theorem}
Assume \ref{H0}, \ref{H1}, \ref{H2}, \ref{H3}, \ref{H4} and \ref{H5} hold. Then there exists a solution $\mu$ to the DLR equations \eqref{eq:DLRequations}, that is also a reversible measure of the SPDE \eqref{eq:LangevinSPDE1}, and is reflection positive (according to Definition \ref{definition:reflectionpositivity} in the case $n=1$).
\end{theorem}
\begin{proof}
Conditions \ref{H0} and \ref{H3} on $H^{\beta}$ and $m$ respectively, imply that $H^{\beta}+m^2$ satisfies the hypotheses of Theorem \ref{theorem:rprepresentation2}. This means that the Gaussian measure $\mu_{\free}$, with covariance $(H^{\beta}+m^2)^{-1}$, is a reflection-positive measure.\\
Since the local symmetric (with respect to reflections) perturbations of a reflection-positive measure are reflection positive (see, e.g., Theorem 10.4.3 of \cite{GlimmJaffe1987} or Section 8 of \cite{DuchDybalskiJahandideh2025}, see also \cite{JaffeJanssens2017}) we have that the measure $\mu^\vartheta$
\[\mu^{\vartheta}(\mathd x)=\frac{1}{\mathcal{Z}^\vartheta} \exp\left(-\int_{-\vartheta}^\vartheta V(x(\xi))\mathd \xi\right)\mu_{\free}(\mathd x),\]
defined in \eqref{eq:existencemutheta}, is reflection positive.\\
Indeed if condition \eqref{eq:rpdefinition} holds for any bounded cylinder function supported in $(0,+\infty)$, then it holds for every $L^2$-function that can be approximated by cylinder functions that are supported in $(0,+\infty)$. A function that can be approximated by cylinder functions supported in $(0,+\infty)$ is $G_\vartheta(x)=\exp\left(-\int_0^\vartheta V(x(\xi))\mathd \xi\right)$. Indeed there is a sequence $\Phi_n$ of bounded cylinder functions, supported in $(0,+\infty)$, converging to $G_\vartheta$ in probability with a uniform $L^{\infty}$ bound. Thus, the fact that $\mu^\vartheta$ is reflection positive follows from the observation that 
\[\exp\left(-\int_{-\vartheta}^\vartheta V(x(\xi))\mathd \xi \right)=G_\vartheta(x)G_\vartheta(\theta(x)).\]
By Theorem \ref{theorem:existence}, there exist an unbounded and increasing sequence $\{\vartheta_n\}_{n\in\N}$ and a solution $\mu$ to the DLR equations \eqref{eq:DLRequations}, which is also a reversible measure of the SPDE \eqref{eq:LangevinSPDE1}, such that $\{\mu^{\vartheta_n}\}_{n\in\N}$ converges weakly to $\mu$ in $E^\gamma_{r,k}$. Since reflection positivity is a property preserved by the weak limits of measures, the conclusion follows.
\end{proof}

\subsection{Translation invariance}\label{section:translationinvariance}

We conclude this section by studying the behavior of the reversible measures with respect to translations. For $a\in\R$ we denote by $\tau_a$ the translation $\tau_ax:=x(\cdot-a)$. Since $1+k|\xi|^2\leq2(1+k|\xi+a|^2)(1+k|a|^2)$, we have $\rho_{r,k}(\xi+a)\leq2^{r/2}(1+k|a|^2)^{r/2}\rho_{r,k}(\xi)$ for every $\xi\in\R$, so that, for every $p\geq1$,
\begin{equation}\label{eq:translationnorm}
\norm{\tau_ax}^p_{L^p_{r,k}}=\int_\R|x(\xi)|^p\rho_{r,k}(\xi+a)\mathd\xi\leq2^{r/2}(1+k|a|^2)^{r/2}\norm{x}^p_{L^p_{r,k}},\qquad x\in L^p_{r,k}.
\end{equation}

\begin{prop}\label{prop:translatedreversible}
Assume that \ref{H0}, \ref{H1}, \ref{H2} and \ref{H3} hold, and let $q\geq1$ be the constant given by \eqref{H2:pol}. Let $\mu$ be a reversible measure for \eqref{eq:SPDEmain1psuedo} such that $\mu(L^{2q}_{r,k})=1$ and $\int_{L^{2q}_{r,k}}\norm{x}^{2q}_{L^{2q}_{r,k}}\mu(\mathd x)<+\infty$. Then, for every $a\in\R$, the measure $\mu\circ\tau_a^{-1}$ is reversible and satisfies the same conditions.
\end{prop}
\begin{proof}
Let $a\in\R$ and $\mu_a:=\mu\circ\tau_a^{-1}$. By \eqref{eq:translationnorm}, $\mu_a(L^{2q}_{r,k})=1$ and $\int\norm{x}^{2q}_{L^{2q}_{r,k}}\mu_a(\mathd x)<+\infty$. Hence, by Theorem \ref{theorem:infinitesimalsymmetri}, $\mu$ satisfies \eqref{eq:IPF}, and it is enough to show that also $\mu_a$ satisfies \eqref{eq:IPF}. Let $\Phi\in\mathcal FC^\infty_b(L^2_{r,k})$ and $h\in C^\infty_c(\R)$. Then $\Phi\circ\tau_a\in\mathcal FC^\infty_b(L^2_{r,k})$ and $\scal{h}{(\delta_x\Phi)(\tau_ax)}_{L^2}=\scal{\tau_{-a}h}{\delta_x(\Phi\circ\tau_a)(x)}_{L^2}$ for every $x\in L^2_{r,k}$. Moreover, since $H^\beta$ commutes with translations and $V$ does not depend on $\xi$, for every $x\in L^{2q}_{r,k}$ we have
\[
\scal{(H^\beta+m^2)h}{\tau_ax}_{L^2}+\scal{h}{V'\circ\tau_ax}_{L^2}=\scal{(H^\beta+m^2)\tau_{-a}h}{x}_{L^2}+\scal{\tau_{-a}h}{V'\circ x}_{L^2}.
\]
Therefore, \eqref{eq:IPF} for $\mu_a$ is \eqref{eq:IPF} for $\mu$, with $\Phi\circ\tau_a$ and $\tau_{-a}h$ in place of $\Phi$ and $h$.
\end{proof}

\begin{cor}\label{corollary:uniquetranslation}
Assume that \ref{H0}, \ref{H1}, \ref{H2}, \ref{H3} and \ref{H4} hold, and let $\nu$ be as in Theorem \ref{thm:rev-ergo}. If $\nu$ is the unique reversible measure for $\{P(t)\}_{t\geq0}$ satisfying \eqref{eq:momentergodic}, then $\nu$ is ergodic and invariant with respect to translations, namely $\nu\circ\tau_a^{-1}=\nu$ for every $a\in\R$. If, in addition, \ref{H5} holds, then $\nu$ is also reflection positive (see Definition \ref{definition:reflectionpositivity}).
\end{cor}
\begin{proof}
By Proposition \ref{prop:translatedreversible}, for every $a\in\R$ the measure $\nu\circ\tau_a^{-1}$ is reversible, satisfies \eqref{eq:momentergodic} and is concentrated on $L^{2q}_{r,k}\cap C(\R)$, so that $\nu\circ\tau_a^{-1}=\nu$. The ergodicity follows from Theorem \ref{thm:rev-ergo}. Finally, if \ref{H5} holds, by the theorem in Subsection \ref{section:ref.positive} there exists a reflection positive solution $\mu$ of the DLR equations, which is reversible, is concentrated on $L^{2q}_{r,k}\cap C(\R)$ by Definition \ref{definition:DLR} and satisfies \eqref{eq:momentergodic} by Theorem \ref{thm:apriori}; hence $\mu=\nu$.
\end{proof}

\noindent If the reversible measure is not unique, the translates of a reversible measure may be different, but their averages give a translation invariant reversible measure. For a Borel probability measure $\mu$ on $L^2_{r,k}$ and $L>0$ we set
\[
\bar\mu_L:=\frac{1}{2L}\int_{-L}^{L}\mu\circ\tau_a^{-1}\,\mathd a,
\]
which is well defined since, by \eqref{eq:translationnorm} and the continuity of the translations in $L^2_{r,k}$, the map $a\mapsto\mu\circ\tau_a^{-1}$ is weakly continuous.

\begin{prop}\label{prop:translationinvariantexistence}
Assume that \ref{H0}, \ref{H1}, \ref{H2}, \ref{H3}, \ref{H4} and \ref{H5} hold, and let $q\geq1$ be the constant given by \eqref{H2:pol}. Let $\mu$ be a reversible measure for \eqref{eq:SPDEmain1psuedo} such that $\mu(L^{2q}_{r,k}\cap C(\R))=1$. Then the family $\{\bar\mu_L\}_{L>0}$ is tight in $L^2_{r,k}$, and every weak limit point $\bar\mu$ of $\bar\mu_L$ as $L\rightarrow+\infty$ is a reversible measure for \eqref{eq:SPDEmain1psuedo} such that $\bar\mu(L^{2q}_{r,k}\cap C(\R))=1$ and $\bar\mu\circ\tau_a^{-1}=\bar\mu$ for every $a\in\R$. In particular, by Theorem \ref{theorem:existence}, there exists a translation invariant reversible measure for \eqref{eq:SPDEmain1psuedo}.
\end{prop}
\begin{proof}
By Theorems \ref{theorem:super-main} and \ref{theorem:FerniqueDLR} we have
\[
\Lambda:=\sup_{\tau\in\R}\int_{L^2_{r,k}}\norm{x}^{2q}_{C^\gamma_{\rho^\tau_{r/2,k}}}\mu(\mathd x)<+\infty .
\]
Since $|x(\xi)|\leq\norm{x}_{C^\gamma_{\rho^\xi_{r/2,k}}}$ and, for every $N\in\N$, there is $C_N>0$ such that $\norm{x}_{C^{\gamma}([\tau-N,\tau+N])}\leq C_N\norm{x}_{C^\gamma_{\rho^\tau_{r/2,k}}}$ for every $\tau\in\R$, we get
\begin{equation}\label{eq:averagedbounds}
\sup_{\xi\in\R}\int|x(\xi)|^{2q}\nu(\mathd x)\leq\Lambda,\qquad\int\norm{x}^{2q}_{C^{\gamma}([-N,N])}\nu(\mathd x)\leq C_N^{2q}\Lambda,\qquad N\in\N,
\end{equation}
for $\nu=\mu$ and for $\nu=\bar\mu_L$, $L>0$. In particular $\int\norm{x}^{2q}_{L^{2q}_{r,k}}\mu(\mathd x)\leq\Lambda\norm{\rho_{r,k}}_{L^1(\R)}$, so that, by Proposition \ref{prop:translatedreversible}, $\mu\circ\tau_a^{-1}$ is reversible for every $a\in\R$, and so is $\bar\mu_L$, since the reversibility condition is linear with respect to the measure. By \eqref{eq:averagedbounds}, the Chebyshev inequality and the Arzel\`a-Ascoli theorem, and since
\[
\int\int_{\{|\xi|>N\}}|x(\xi)|^2\rho_{r,k}(\xi)\mathd\xi\,\bar\mu_L(\mathd x)\leq\Lambda^{\frac1q}\int_{\{|\xi|>N\}}\rho_{r,k}(\xi)\mathd\xi,\qquad L>0,\ N\in\N,
\]
the family $\{\bar\mu_L\}_{L>0}$ is tight in $L^2_{r,k}$.\\

\noindent Let $\bar\mu$ be the weak limit of $\{\bar\mu_{L_n}\}_{n\in\N}$, with $L_n\rightarrow+\infty$. By the lower semicontinuity on $L^2_{r,k}$ of the functionals $x\mapsto\int_{-N}^N|x(\xi)|^{2q}\rho_{r,k}(\xi)\mathd\xi$ and $x\mapsto\norm{x}_{C^{\gamma}([-N,N])}$ (which is set equal to $+\infty$ if $x$ has no $\gamma$-H\"older continuous version on $[-N,N]$), \eqref{eq:averagedbounds} gives $\int\norm{x}^{2q}_{L^{2q}_{r,k}}\bar\mu(\mathd x)\leq\Lambda\norm{\rho_{r,k}}_{L^1(\R)}$ and $\int\norm{x}^{2q}_{C^{\gamma}([-N,N])}\bar\mu(\mathd x)\leq C_N^{2q}\Lambda$ for every $N\in\N$, so that $\bar\mu(L^{2q}_{r,k}\cap C(\R))=1$. Moreover, by \eqref{solmild4} and the dominated convergence theorem, $P(t)\Phi\in C_b(L^2_{r,k})$ for every $\Phi\in C_b(L^2_{r,k})$ and $t\geq0$, so that we can pass to the limit in the reversibility condition for $\bar\mu_{L_n}$ (see Definition \ref{definition:reversible}), and $\bar\mu$ is reversible. Finally, let $b\in\R$. Since $\tau_b\circ\tau_a=\tau_{a+b}$, the total variation of $\bar\mu_L\circ\tau_b^{-1}-\bar\mu_L$ is not greater than $|b|/L$; since $\tau_b$ is continuous on $L^2_{r,k}$, letting $L=L_n\rightarrow+\infty$ we get $\bar\mu\circ\tau_b^{-1}=\bar\mu$.
\end{proof}

\begin{rmk}\label{rmk:latticestationary}
For lattice systems more is known. For stochastic Ising models, namely spin-flip dynamics of $\pm1$ spins on $\mathbb Z^d$, and for stochastic gradient systems, namely Langevin diffusions of unbounded real spins on $\mathbb Z^d$, with finite range interactions, every translation invariant stationary measure is a Gibbs measure (see \cite{Holley1971,Fritz1982}) and, if $d\leq2$, every stationary measure is a Gibbs measure (see \cite{HolleyStroock1977,Fritz1982}). For \eqref{eq:SPDEmain1psuedo}, the analogous question, namely whether every (translation invariant) invariant measure is reversible, and hence, by Theorem \ref{theorem:super-main}, a solution of the DLR equations, is open.
\end{rmk}

\part{Applications to uniqueness and non-uniqueness results}

\section{Relation with DLR equations defined via Brownian bridge}\label{section:OScomparison}

In this section we discuss the relation between the DLR equations introduced in Section \ref{section:DLRequation} and a similar approach using Brownian bridge proposed by \cite{OsadaSpohn1999} and then extended and further discussed in, e.g., \cite{BetzLorinczi2003,Gubinelli2006,GubinelliLorinczi2009,Lorinczi2011,LorincziMinlos2001}. At the end of the section we state the uniqueness and non-uniqueness results that we prove in Sections \ref{section:uniquenessDLR} and \ref{section:nonuniqueness}, and we compare them with the ones of \cite{OsadaSpohn1999}.\\

\noindent In order to make the comparison possible, in the first part of this section we consider an operator $H^{\beta}$, where $\beta_0>0$ and $\int_0^1\xi\bar\beta(\xi)\mathd \xi<+\infty$. In this case the operator $H^{\beta}$ can be reasonably treated as a ``regular'' perturbation of the Laplacian $-\beta_0 \Delta$. \\
First, we define a Brownian bridge $\mathcal W^{\mathcal I,(y^a,y^b)}$ on an interval $\mathcal{I}=(a,b)$, where $a<b$, and boundary values $y^a,y^b \in \R$. If we denote the heat kernel by 
\[
H_t(z)=(2\pi t)^{-1/2}\exp\left(-\frac{|z|^2}{2t}\right),
\]
we define the measure $\mathcal W^{\mathcal I,(y^a,y^b)}$ as the probability measure supported on $C(\mathcal{I},\R)$ such that 
its finite-dimensional distributions are
\begin{equation}\label{eq:Bbridge}
\int\Phi(x)\mathcal W^{\mathcal I,(y^a,y^b)}(\mathd x)
=\int_{\R^{r}}\bar\Phi(y^1,\ldots,y^r)
\frac{H_{\xi_1-a}(y^1-y^a)\cdots H_{b-\xi_r}(y^b-y^r)}
{H_{b-a}(y^b-y^a)}\mathd y^1\cdots\mathd y^r,
\end{equation}
for every $\Phi(x)=\bar\Phi(x(\xi_1),\ldots,x(\xi_r))$, where $a<\xi_1<\cdots<\xi_r<b$ and $\bar\Phi:\R^r\rightarrow\R$ is a bounded Borel function. We consider the Brownian bridge associated with the operator $-\beta_0 \Delta$: 
\begin{equation}\label{eq:betaBB}
\mathcal W_{\beta_0}^{\mathcal I,(z_1,z_2)}
:=
\left(x\mapsto \beta_0^{-1/2}x\right)_*
\mathcal W^{\mathcal I,(\sqrt{\beta_0}z_1,\sqrt{\beta_0}z_2)},
\end{equation}
where $\mathcal{W}^{\mathcal{I},(z_1,z_2)}$ has been defined in \eqref{eq:Bbridge}. Then we have the following theorem:

\begin{theorem}\label{theorem:BB}
Assume that \ref{H0}, \ref{H1}, \ref{H2}, \ref{H3} and \ref{H4} hold, with $\beta_0>0$ and $\int_0^1\xi\bar\beta(\xi)\mathd\xi<+\infty$, and let $\mu$ be a Borel probability measure on $L^2_{r,k}$ such that $\mu(L^{2q}_{r,k}\cap C(\R))=1$, where $q\geq1$ is the constant given by \eqref{H2:pol}. Then $\mu$ satisfies the DLR equations in the sense of Definition \ref{definition:DLR} if and only if, for every bounded open interval $\mathcal{I}=(a,b)$ and $\mu^{\mathcal I^c}$-almost every $y\in L^2_{r,k}(\mathcal I^c)\cap C(\mathcal I^c)$, we have
\begin{multline}\label{eq:DLRBB}
\frac{\mathd \mu^{\mathcal{I},y}}{\mathd \mathcal{W}^{\mathcal{I},(y(a),y(b))}_{\beta_0}}(x)=\frac{1}{\hat{\mathcal{Z}}^{\mathcal{I}}(y)}\exp\left(
-\int_{\mathcal{I}}V_{\tot}(x(\xi))\mathd \xi
-\frac{1}{4}\int_{\mathcal{I}^2}\bar{\beta}(\xi-\xi')|x(\xi)-x(\xi') |^2\mathd \xi\mathd \xi'
\right.\\
\left.
-\frac{1}{2}\int_{\mathcal{I}\times\mathcal{I}^c}\bar{\beta}(\xi-\xi')|x(\xi)-y(\xi') |^2\mathd \xi\mathd \xi'
\right),
\end{multline}
where $V_{\tot}(z):=\frac{m^2}{2}|z|^2+V(z)$, $z\in\R$,
\begin{multline*}
\hat{\mathcal{Z}}^{\mathcal{I}}(y)=\int \exp\left(
-\int_{\mathcal{I}}V_{\tot}(x(\xi))\mathd \xi
-\frac{1}{4}\int_{\mathcal{I}^2}\bar{\beta}(\xi-\xi')|x(\xi)-x(\xi') |^2\mathd \xi\mathd \xi'
\right.\\
\left.
-\frac{1}{2}\int_{\mathcal{I}\times\mathcal{I}^c}\bar{\beta}(\xi-\xi')|x(\xi)-y(\xi') |^2\mathd \xi\mathd \xi'
\right)
\mathcal{W}_{\beta_0}^{\mathcal{I},(y(a),y(b))}(\mathd x),
\end{multline*}
and, as in Lemma \ref{lemma:extendeduniqueness1} below, every continuous function $x$ on $\overline{\mathcal I}$ such that $x(a)=y(a)$ and $x(b)=y(b)$ is identified with the function equal to $x$ on $\mathcal I$ and to $y$ on $\mathcal I^c$.
\end{theorem}
In order to prove Theorem \ref{theorem:BB}, we start with an abstract lemma. We first introduce some notation and concepts. Let $\HH$ be a real separable Hilbert space, let $A$ be a densely defined self-adjoint operator on $\HH$ such that $A>0$ and $A^{-1}$ is trace class, and let $B$ be a nonnegative self-adjoint operator on $\HH$ such that 
\[ {\rm Dom}(A^{1/2})\subseteq{\rm Dom}(B^{1/2}),\quad \quad K:=(B^{1/2}A^{-1/2})^*(B^{1/2}A^{-1/2})\text{ is trace class}\] 
(if $B\in\mathcal L(\HH)$, then $K=A^{-1/2}BA^{-1/2}$ is trace class, since $A^{-1/2}$ is Hilbert-Schmidt). We set $\langle Bx,x\rangle_{\HH}:=\norm{B^{1/2}x}^2_{\HH}$ for $x\in{\rm Dom}(B^{1/2})$, and we denote by $A+B$ the self-adjoint operator associated with the closed quadratic form $x\mapsto\norm{A^{1/2}x}^2_{\HH}+\norm{B^{1/2}x}^2_{\HH}$, $x\in{\rm Dom}(A^{1/2})$. Let $\mu^A$ and $\mu^{A+B}$ be the Gaussian measures on $\HH$ with mean $0$ and covariance $A^{-1}$ and $(A+B)^{-1}$, respectively.

\begin{lemma}\label{lemma:Gaussianproof}
Using the notation and hypotheses described just above, we have that $\mu^A({\rm Dom}(B^{1/2}))=1$, $\mu^{A+B}$ and $\mu^A$ are equivalent, and
\[
    \frac{\dd\mu^{A+B}}{\dd\mu^A}(x)
    =
    \det\nolimits_F\!\left(
        I+K
    \right)^{1/2}
    \exp\!\left(-\frac12\langle Bx,x\rangle_{\HH}\right),
\]
where $\det_F$ denotes the Fredholm determinant, for $\mu^A$-almost every $x\in\HH$.
\end{lemma}
\begin{proof}
The operator $T:=B^{1/2}A^{-1/2}$ is bounded by the closed graph theorem, and $K=T^*T$. Let $\{v_j\}_{j\in\N}$ be an orthonormal basis of $\HH$ such that $Kv_j=\lambda_jv_j$, with $\lambda_j\geq0$ and $\sum_j\lambda_j=\Tr K<+\infty$. In particular, $\det_F(I+K)=\prod_j(1+\lambda_j)$ is well defined and strictly positive. Since $A^{-1/2}$ is Hilbert-Schmidt, a random variable $X$ with law $\mu^A$ can be written as $X=\sum_jg_jA^{-1/2}v_j$, where $\{g_j\}_{j\in\N}$ are independent standard Gaussian random variables and the series converges in $L^2(\Omega;\HH)$ and almost surely. The partial sums $X_n:=\sum_{j\leq n}g_jA^{-1/2}v_j$ belong to ${\rm Dom}(A^{1/2})$, and $B^{1/2}X_n=\sum_{j\leq n}g_jTv_j$ converges in $L^2(\Omega;\HH)$, since $\langle Tv_i,Tv_j\rangle_{\HH}=\lambda_j\delta_{ij}$. Since $B^{1/2}$ is closed, we get $X\in{\rm Dom}(B^{1/2})$ almost surely and
\[
\langle BX,X\rangle_{\HH}=\sum_j\lambda_jg_j^2,\qquad\langle h,X\rangle_{\HH}=\sum_j\langle A^{-1/2}h,v_j\rangle_{\HH}\,g_j,\qquad h\in\HH.
\]
Hence, by the independence of $\{g_j\}_{j\in\N}$,
\begin{align}
&\mathbb E\left[
    \exp\left(
        -\frac12\langle BX,X\rangle_{\HH}
        +i\langle h,X\rangle_{\HH}
    \right)
\right]                                                     \notag\\
&\quad =
\det\nolimits_F(I+K)^{-1/2}
\exp\left(
    -\frac12
    \left\langle
        h,
        A^{-1/2}(I+K)^{-1}A^{-1/2}h
    \right\rangle_{\HH}
\right),
\qquad h\in\HH.\label{eq:Gaussianrpoof1}
\end{align}
Taking $h=0$ in \eqref{eq:Gaussianrpoof1} gives
\[
\int_{\HH}\exp\left(-\frac12\langle Bx,x\rangle_{\HH}\right)\mu^A(\dd x)
=\det\nolimits_F(I+K)^{-1/2}.
\]
The characteristic functional of the probability measure defined by the density in the statement is therefore
\[
\exp\left(
    -\frac12
    \left\langle
        h,A^{-1/2}(I+K)^{-1}A^{-1/2}h
    \right\rangle_{\HH}
\right).
\]
Since
\[
(A+B)^{-1}=A^{-1/2}(I+K)^{-1}A^{-1/2},
\]
this is the characteristic functional of $\mu^{A+B}$. Indeed, for $h\in\HH$ the vector $u:=A^{-1/2}(I+K)^{-1}A^{-1/2}h$ belongs to ${\rm Dom}(A^{1/2})$ and, for every $v\in{\rm Dom}(A^{1/2})$, $\langle A^{1/2}u,A^{1/2}v\rangle_{\HH}+\langle B^{1/2}u,B^{1/2}v\rangle_{\HH}=\langle(I+K)(I+K)^{-1}A^{-1/2}h,A^{1/2}v\rangle_{\HH}=\langle h,v\rangle_{\HH}$.
\end{proof}

\begin{lemma}\label{lemma:extendeduniqueness1}
Assume that \ref{H0}, \ref{H1} and \ref{H3} hold, with $\beta_0>0$ and $\int_0^1\xi\bar\beta(\xi)\mathd\xi<+\infty$. Let $\mathcal I=(a,b)$ and $y\in L^2_{r,k}(\mathcal I^c)\cap C(\mathcal I^c)$, and let us identify every continuous function $x$ on $\overline{\mathcal I}$ such that $x(a)=y(a)$ and $x(b)=y(b)$ with the function equal to $x$ on $\mathcal I$ and to $y$ on $\mathcal I^c$. Then the integrals in the following exponent are finite for $\mathcal{W}_{\beta_0}^{\mathcal{I},(y(a),y(b))}$-almost every $x$, and we have
\begin{multline}\label{eq:otheruniqueness1}
\mu^{\mathcal{I},y}_{\free}(\mathd x)=\frac{1}{\bar{\mathcal{Z}}^{\mathcal{I}}(y)}
\exp\left(
-\int_{\mathcal{I}}V_{\free}(x(\xi))\mathd \xi
-\frac{1}{4}\int_{\mathcal{I}^2}\bar{\beta}(\xi-\xi')|x(\xi)-x(\xi') |^2\mathd \xi\mathd \xi'
\right.\\
\left.
-\frac{1}{2}\int_{\mathcal{I}\times\mathcal{I}^c}\bar{\beta}(\xi-\xi')|x(\xi)-y(\xi') |^2\mathd \xi\mathd \xi'
\right)
\mathcal{W}_{\beta_0}^{\mathcal{I},(y(a),y(b))}(\mathd x),
\end{multline}
where $V_{\free}(z)=m^2|z|^2/2$.
\end{lemma}
\begin{proof}
We denote by $H_{nl}$ the nonlocal part of $H^\beta$ and by $A$ the realization of $-\beta_0\Delta$ in $L^2(\mathcal I)$ with homogeneous Dirichlet boundary conditions, so that ${\rm Dom}(A^{1/2})=H^1_0(\mathcal I)$, $A^{-1}$ is trace class and the centered Gaussian measure with covariance $A^{-1}$ is $\mathcal W_{\beta_0}^{\mathcal I,(0,0)}$. In this proof, the functions in $L^2(\mathcal I)$ are extended by zero outside $\mathcal I$. For $x\in H^1_0(\mathcal I)$ we have $\mathcal E^{nl}(x,x)\leq\int_\R\bar\beta(z)\min\{|z|^2\norm{x'}^2_{L^2},4\norm{x}^2_{L^2}\}\mathd z<+\infty$, where $\mathcal E^{nl}$ is given by \eqref{eq:dirichletform}. Since, by the Fatou lemma, the quadratic form $x\mapsto\mathcal E^{nl}(x,x)+m^2\norm{x}^2_{L^2(\mathcal I)}$ is lower semicontinuous on $L^2(\mathcal I)$, its restriction to $H^1_0(\mathcal I)$ is closable; we denote by $B$ the nonnegative self-adjoint operator associated with its closure, and we note that, again by the Fatou lemma,
\begin{align*}
\langle Bx,x\rangle_{L^2(\mathcal I)}
&=\mathcal E^{nl}(x,x)+m^2\norm{x}^2_{L^2(\mathcal I)}\\
&=m^2\int_{\mathcal I}|x(\xi)|^2\mathd\xi
+\frac12\int_{\mathcal I^2}\bar\beta(\xi-\xi')|x(\xi)-x(\xi')|^2\mathd\xi\mathd\xi'
+\int_{\mathcal I\times\mathcal I^c}\bar\beta(\xi-\xi')|x(\xi)|^2\mathd\xi\mathd\xi'
\end{align*}
for every $x\in{\rm Dom}(B^{1/2})$. Since $\mathcal E=\mathcal E^l+\mathcal E^{nl}+m^2\langle\cdot,\cdot\rangle_{L^2}$ and, as $\beta_0>0$, $\mathcal F_{\mathcal I}=H^1_0(\mathcal I)$ (see \eqref{eq:energy-is-Hs1}), the operator $A+B$ defined in Lemma \ref{lemma:Gaussianproof} is the operator $\mathcal A_{\mathcal I}$ of Proposition \ref{prop:gaussL2I} (see \cite[Theorem 4.4.3]{Fukushima2011}). Moreover, let $\{e_j\}_{j\in\N}$ be the orthonormal basis of $L^2(\mathcal I)$ such that $Ae_j=\beta_0\pi^2j^2|\mathcal I|^{-2}e_j$. Since $\sum_j|\mathcal I|^2(\beta_0\pi^2j^2)^{-1}e_j(\xi)e_j(\xi')=\int x(\xi)x(\xi')\mathcal W^{\mathcal I,(0,0)}_{\beta_0}(\mathd x)=\beta_0^{-1}|\mathcal I|^{-1}(\xi\wedge\xi'-a)(b-\xi\vee\xi')$, we get
\begin{align*}
\Tr K&=\sum_j\norm{B^{1/2}A^{-1/2}e_j}^2_{L^2(\mathcal I)}=\int\left(\mathcal E^{nl}(x,x)+m^2\norm{x}^2_{L^2(\mathcal I)}\right)\mathcal W^{\mathcal I,(0,0)}_{\beta_0}(\mathd x)\\
&\leq\frac{m^2|\mathcal I|^2}{\beta_0}+\frac{1}{2\beta_0}\int_{\mathcal I^2}\bar\beta(\xi-\xi')|\xi-\xi'|\mathd\xi\mathd\xi'+\frac1{\beta_0}\int_{\mathcal I\times\mathcal I^c}\bar\beta(\xi-\xi')\min\{\xi-a,b-\xi\}\mathd\xi\mathd\xi',
\end{align*}
which is finite by \eqref{eq:levynear}, \eqref{eq:levygrowth} and since $\int_0^1\xi\bar\beta(\xi)\mathd\xi<+\infty$. Hence, by Proposition \ref{prop:gaussL2I} and Lemma \ref{lemma:Gaussianproof},
\[
\mu^{\mathcal I,0}_{\free}(\mathd x)
=\frac{1}{\bar{\mathcal Z}^{\mathcal I}(0)}
\exp\left(-\frac12\langle Bx,x\rangle_{L^2(\mathcal I)}\right)
\mathcal W_{\beta_0}^{\mathcal I,(0,0)}(\mathd x),
\]
which is \eqref{eq:otheruniqueness1} with $y=0$.\\

\noindent Let now $y$ be the restriction to $\mathcal I^c$ of a function of $C^\infty_c(\R)$, so that the harmonic extension $\varphi^y$ of $y$ belongs to $\mathcal F=H^1(\R)$ by Lemma \ref{lemma:harmonicextensionenergy}, and let $\psi^y$ be the function which is affine on $\overline{\mathcal I}$ and takes the values $y(a)$ and $y(b)$ at $a$ and $b$. By Proposition \ref{proposition::cutoffdomain}, $\varphi^y|_{\mathcal I}\in H^2(\mathcal I)$ and $-\beta_0(\varphi^y)''+H_{nl}\varphi^y+m^2\varphi^y=0$ in $\mathcal I$. In particular, since $\varphi^y$ is continuous on $\R$, $H^2(\mathcal I)\subseteq C^1(\overline{\mathcal I})$ and $\varphi^y$ coincides with a function of $C^\infty_c(\R)$ on $\mathcal I^c$, $\varphi^y$ is Lipschitz continuous on $\R$, and $d^y:=\varphi^y-\psi^y\in H^2(\mathcal I)\cap H^1_0(\mathcal I)={\rm Dom}(A)$, with $Ad^y=-\beta_0(\varphi^y)''$ in $\mathcal I$. By Corollary \ref{corollary:Guassiantranslation}, $\mu^{\mathcal I,y}_{\free}$ is the image of $\mu^{\mathcal I,0}_{\free}$ under the map $x\mapsto x+\varphi^y$, while $\mathcal W_{\beta_0}^{\mathcal I,(y(a),y(b))}$ is the image of $\mathcal W_{\beta_0}^{\mathcal I,(0,0)}$ under the map $x\mapsto x+\psi^y$. Since $d^y\in{\rm Dom}(A)$, the Cameron-Martin formula (see \cite[Corollary 2.4.3]{Bogachev98}) yields, up to a factor depending only on $y$,
\[
\mu^{\mathcal I,y}_{\free}(\mathd x)
=C(y)\exp\left(
-\frac12\langle B(x-\varphi^y),x-\varphi^y\rangle_{L^2(\mathcal I)}
+\langle A d^y,x\rangle_{L^2(\mathcal I)}
\right)\mathcal W_{\beta_0}^{\mathcal I,(y(a),y(b))}(\mathd x),
\]
where, for $\mathcal W_{\beta_0}^{\mathcal I,(y(a),y(b))}$-almost every $x$, $u:=x-\varphi^y=(x-\psi^y)-d^y$ belongs to ${\rm Dom}(B^{1/2})$ by Lemma \ref{lemma:Gaussianproof}, and it vanishes outside $\mathcal I$. Since $\varphi^y$ is bounded and Lipschitz continuous and $u$ is bounded, by the symmetry of $\bar\beta$ and the Fubini theorem we have $\mathcal E^{nl}(u,\varphi^y)=\int_{\mathcal I}uH_{nl}\varphi^y\mathd\xi$, and, writing $x=u+\varphi^y$ and $\langle Ad^y,u\rangle_{L^2(\mathcal I)}=-\int_{\mathcal I}u\,(H_{nl}\varphi^y+m^2\varphi^y)\mathd\xi$, we obtain
\[
-\frac12\langle Bu,u\rangle_{L^2(\mathcal I)}+\langle Ad^y,x\rangle_{L^2(\mathcal I)}
=-\frac12\mathcal E^{nl}(x,x)-\frac{m^2}2\norm{x}^2_{L^2(\mathcal I)}+C'(y),
\]
with $C'(y):=\frac12\mathcal E^{nl}(\varphi^y,\varphi^y)+\frac{m^2}{2}\norm{\varphi^y}^2_{L^2(\mathcal I)}+\langle Ad^y,\varphi^y\rangle_{L^2(\mathcal I)}$. Finally,
\begin{align*}
\frac12\mathcal E^{nl}(x,x)
={}&\frac14\int_{\mathcal I^2}\bar\beta(\xi-\xi')|x(\xi)-x(\xi')|^2\mathd\xi\mathd\xi'
+\frac12\int_{\mathcal I\times\mathcal I^c}
\bar\beta(\xi-\xi')|x(\xi)-y(\xi')|^2\mathd\xi\mathd\xi'\\
&+\frac14\int_{(\mathcal I^c)^2}\bar\beta(\xi-\xi')|y(\xi)-y(\xi')|^2\mathd\xi\mathd\xi',
\end{align*}
where the last term depends only on $y$. Absorbing the constants into the normalization gives \eqref{eq:otheruniqueness1}.\\

\noindent In the general case, let $\{y_n\}_{n\in\N}$ be restrictions to $\mathcal I^c$ of functions of $C^\infty_c(\R)$ such that $y_n(a)=y(a)$, $y_n(b)=y(b)$, and $y_n\rightarrow y$ in $L^2_{r,k}(\mathcal I^c)$ and uniformly on the compact subsets of $\mathcal I^c$. By Corollary \ref{corollary:Guassiantranslation}, $\mu^{\mathcal I,y_n}_{\free}$ converges weakly to $\mu^{\mathcal I,y}_{\free}$. On the other hand, the right-hand side of \eqref{eq:otheruniqueness1} depends on $y$ only through the last integral in the exponent, which is finite for every $x\in C(\overline{\mathcal I})$, since $\int_{\mathcal I}\int_{\mathcal I^c}\bar\beta(\xi-\xi')\mathd\xi'\mathd\xi<+\infty$ and, as in the proof of Theorem \ref{thm:harmonic-continuity}, $\int_{\mathcal I}\bar\beta(\xi-\xi')\mathd\xi\lesssim\rho_{r,k}(\xi')$ for $\xi'\in\mathcal I^c$ far from $\mathcal I$. For the same reasons, this integral with $y_n$ in place of $y$ converges, as $n\rightarrow+\infty$, to the one with $y$. Since the densities in \eqref{eq:otheruniqueness1} are bounded by $1$, up to the normalization constants, which converge to a positive limit, the right-hand side of \eqref{eq:otheruniqueness1} with $y_n$ converges in total variation to the one with $y$, and the proof is complete.
\end{proof}

\begin{proof}[Proof of Theorem \ref{theorem:BB}]
Let $\mathcal I=(a,b)$ be a bounded open interval and let $y\in L^2_{r,k}(\mathcal I^c)\cap C(\mathcal I^c)$. By Lemma \ref{lemma:extendeduniqueness1}, the probability measure $\mathcal Z^{\mathcal I}(y)^{-1}\exp\left(-\int_{\mathcal I}V(x(\xi))\mathd\xi\right)\mu^{\mathcal I,y}_{\free}(\mathd x)$, which appears in \eqref{eq:DLRequations}, has density with respect to $\mathcal W^{\mathcal I,(y(a),y(b))}_{\beta_0}$ given by the right-hand side of \eqref{eq:DLRBB}. Hence \eqref{eq:DLRequations} holds for $\mathcal I$ and $y$ if and only if \eqref{eq:DLRBB} holds, and the statement follows from Definition \ref{definition:DLR}.
\end{proof}

In the rest of this section we compare Definition \ref{definition:DLR} with the definition of Gibbs measures relative to Brownian motion of \cite{OsadaSpohn1999}, and the results of Sections \ref{section:uniquenessDLR} and \ref{section:nonuniqueness} with the ones of \cite{OsadaSpohn1999}. In \cite{OsadaSpohn1999} the Gibbs measures are probability measures on the space $C(\R)$ of the continuous functions on $\R$, endowed with the topology of the uniform convergence on compact sets, and they are defined by prescribing their conditional probabilities on the intervals $(-T,T)$, $T>0$, as Gibbs modifications of the Brownian bridge, through an external potential $\phi$ and an interaction potential $w$ (see \cite[Definition 2.2]{OsadaSpohn1999}).

\begin{remark}\label{remark:OSscaling}
Since the image of $\mathcal W_{\beta_0}^{\mathcal I,(z_1,z_2)}$ under the map $x\mapsto\sqrt{\beta_0}x$ is $\mathcal W^{\mathcal I,(\sqrt{\beta_0}z_1,\sqrt{\beta_0}z_2)}$, by Theorem \ref{theorem:BB} the conditional probabilities prescribed by the DLR equations are, up to this map, the ones prescribed in \cite[Definition 2.2]{OsadaSpohn1999}, with external potential and interaction potential given by
\begin{equation}\label{eq:OSpotentials}
\phi(z):=\frac{m^2}{2\beta_0}|z|^2+V\left(\frac{z}{\sqrt{\beta_0}}\right),\qquad w(\xi,z):=\frac{1}{2\beta_0}\bar\beta(\xi)|z|^2,\qquad \xi,z\in\R.
\end{equation}
\end{remark}

\begin{remark}\label{remark:OSdefinition}
In \cite{OsadaSpohn1999} more general interactions of the form $w(\xi,z)=\rho(\xi)v(z)$, with $v$ convex, and external potentials $\phi$ which are only locally integrable (and possibly equal to $+\infty$ outside an interval) are considered; on the other hand, in \cite{OsadaSpohn1999} the function $\rho$ is assumed to be bounded (see assumption (W3) therein), which in our setting means $\bar\beta(0)<+\infty$. Our Definition \ref{definition:DLR} does not require any condition on $\bar\beta$ near the origin besides \ref{H0}. If $\beta_0>0$ and $\int_0^1\xi\bar\beta(\xi)\mathd\xi=+\infty$, the expectation with respect to the Brownian bridge of the interaction term in \eqref{eq:otheruniqueness1} is infinite and Lemma \ref{lemma:extendeduniqueness1} does not apply, while if $\beta_0=0$ there is no Brownian component at all; in both cases the measures $\mu^{\mathcal I,y}_{\free}$, and hence Definition \ref{definition:DLR}, are still well defined.
\end{remark}

\noindent We conclude this section by stating, in the following explicit example, the uniqueness and non-uniqueness results that we prove in Sections \ref{section:uniquenessDLR} and \ref{section:nonuniqueness} below (see also Proposition \ref{proposition:regimes}), and by comparing them with the ones of \cite{OsadaSpohn1999}. For $\beta_0,\beta_1\geq 0$ with $\beta_0+\beta_1>0$, $\lambda\geq 0$ and $s\in(0,1)$, we consider
\[
H^{\beta}f(\xi)=-\beta_0\Delta f(\xi)+\frac{\beta_1}{2}\int_{\R}\frac{2f(\xi)-f(\xi+\xi')-f(\xi-\xi')}{(\lambda+|\xi'|)^{1+2s}}\,\mathd \xi',
\]
where, if $\beta_0=0$, we assume $\lambda=0$ and $s\in\left(\frac12,1\right)$. In this case $\bar\beta(\xi)=\beta_1(\lambda+|\xi|)^{-1-2s}$, and we assume $\beta_1>0$ (if $\beta_1=0$, then \eqref{eq:firstmoment} holds and the solution of the DLR equations is unique by Corollary \ref{corollary:uniqueness}). By Remark \ref{remark:OSscaling}, the setting of \cite{OsadaSpohn1999} corresponds to the case $\beta_0>0$ and $\lambda>0$, with $\phi$ and $w$ given by \eqref{eq:OSpotentials}.
\begin{itemize}
\item \emph{Uniqueness.} Condition \eqref{eq:firstmoment} holds if and only if $s\in(\frac12,1)$. In this case, by Theorem \ref{theorem:uniqueness} and Corollary \ref{corollary:uniqueness}, for every $\beta_0\geq0$, $\lambda\geq0$ and every potential $V$ satisfying \ref{H2}, \ref{H4} and \ref{H5}, the DLR equations have a unique solution. In \cite[Theorems 2.3 and 5.4]{OsadaSpohn1999} uniqueness is proved for $\beta_0>0$, $\lambda>0$ and $\int_0^{+\infty}\xi\rho(\xi)\mathd\xi<+\infty$, which is equivalent to $s>\frac12$, among the translation invariant Gibbs measures under which $x(0)$ has a finite moment of some order strictly greater than $2$, or among the $\mathbf d$-tame Gibbs measures (see \cite[Definition 5.3]{OsadaSpohn1999}), under the assumption that $\phi$ admits symmetric right and left dominators (see \cite[Definition 2.1]{OsadaSpohn1999}). For this family of operators, our result does not require the existence of dominators, it gives uniqueness among all the solutions of the DLR equations, and it also covers the fractional cases $\lambda=0$, both with $\beta_0>0$ and with $\beta_0=0$, which are not considered in \cite{OsadaSpohn1999} (on the other hand, in \cite{OsadaSpohn1999} potentials $\phi$ which do not satisfy \ref{H2} are also allowed, see Remark \ref{remark:OSdefinition}).
\item \emph{Non-uniqueness.} If $\beta_0>0$ and either $\lambda>0$ and $s\in(0,\frac12]$, or $\lambda=0$ and $s\in(0,\frac12)$, then \eqref{eq:conditionbarbeta} holds, and Theorem \ref{theorem:nonuniqueness} shows that, for every $f$ satisfying Definition \ref{definition:doublewell} and for suitable values of the parameters $\gamma_1$ and $\gamma_2$ introduced in Section \ref{section:nonuniqueness}, the DLR equations with operator $H^{\beta,\gamma_1}$ and potential $V^{\gamma_2}$ have at least two different solutions, which are reversible measures for \eqref{eq:LangevinSPDE1}. In \cite[Theorems 2.4 and 6.1]{OsadaSpohn1999} non-uniqueness is proved for $\phi(z)=\beta(z^4-z^2)$ and $w(\xi,z)=\frac{\alpha}{2}(1+|\xi|)^{-\gamma}|z|^2$, with $\alpha,\beta>0$ large enough and $\gamma\in(1,2]$, which corresponds to $\lambda=1$, $s=\frac{\gamma-1}2\in(0,\frac12]$ and to a quartic potential (in \cite{OsadaSpohn1999} the two Gibbs measures are also translation invariant). Hence, also in this case our result is more general, since it covers general double-well potentials and the fractional case $\lambda=0$, $s\in(0,\frac12)$.
\end{itemize}
These results can be summarized in the following table, where $\beta_1>0$, ``uniqueness'' means that the DLR equations have a unique solution for every potential $V$ satisfying \ref{H2}, \ref{H4} and \ref{H5}, and ``non-uniqueness'' means that there are potentials $V$ satisfying \ref{H2}, \ref{H4} and \ref{H5} for which the DLR equations have at least two solutions.
\begin{center}
\begin{tabular}{l|c|c|c}
 & $s\in(0,\frac12)$ & $s=\frac12$ & $s\in(\frac12,1)$\\
\hline
$\beta_0>0$, $\lambda>0$ & non-uniqueness & non-uniqueness & uniqueness\\
$\beta_0>0$, $\lambda=0$ & non-uniqueness & open & uniqueness\\
$\beta_0=0$, $\lambda=0$ & excluded by \ref{H0} & excluded by \ref{H0} & uniqueness
\end{tabular}
\end{center}
The first line corresponds to the setting of \cite{OsadaSpohn1999}, while the second and the third ones correspond to the fractional case $\lambda=0$, which is not covered by \cite{OsadaSpohn1999}. The only case which is not covered by our results is $\beta_0>0$, $\lambda=0$ and $s=\frac12$, where \eqref{eq:firstmoment} fails and $\int_0^1\xi\bar\beta(\xi)\mathd\xi=+\infty$, so that the interaction at short distances would require a renormalization.

\section{On uniqueness}\label{section:uniquenessDLR}

Throughout this section we assume that \ref{H0}, \ref{H1}, \ref{H2}, \ref{H3} and \ref{H4} hold, we let $q\geq1$ be the constant given by \eqref{H2:pol} and, for $T>0$, we let $\mathcal I_T:=(-T,T)$. We prove that, if the L\'evy density $\bar\beta$ has a finite first moment at infinity, namely
\begin{equation}\label{eq:firstmoment}
\int_{\{\abs{z}\geq1\}}\abs{z}\bar\beta(z)\,\mathd z<+\infty,
\end{equation}
then the DLR equations \eqref{eq:DLRequations} have at most one solution $\mu$ such that
\begin{equation}\label{eq:uniformmoments}
\sup_{\xi\in\R}\int_{L^2_{r,k}}\abs{x(\xi)}^{2q}\mu(\mathd x)<+\infty,
\end{equation}
and that, under \ref{H5}, every solution satisfies \eqref{eq:uniformmoments} (see Theorem \ref{theorem:uniqueness} and Corollary \ref{corollary:uniqueness}). By \ref{H0}, condition \eqref{eq:firstmoment} holds whenever $s_2>\frac12$; it is satisfied, for instance, by the Laplacian, by the operators associated with the superposition of Ornstein-Uhlenbeck processes and by the massive fractional Laplacians, whose L\'evy densities decay exponentially, and by the fractional Laplacians of order $s\in(\frac12,1)$ (see the examples after \ref{H0}).\\

\noindent Since
\[
\int_{-\infty}^{0}\int_{0}^{+\infty}\bar\beta(\xi-\xi')\,\mathd\xi'\,\mathd\xi=\int_0^{+\infty}z\bar\beta(z)\,\mathd z,
\]
condition \eqref{eq:firstmoment} means that the interaction between the two half-lines is finite, apart from the contribution of the singularity of $\bar\beta$ at the origin. Hence Theorem \ref{theorem:uniqueness} can be seen as a continuum counterpart of the classical uniqueness theorem for one-dimensional lattice systems with bounded spins and pair interactions $J$ such that $\sum_{j\geq1}j\abs{J(j)}<+\infty$ (see \cite{Ruelle1968} and \cite[Section 8.3]{Georgii2011}). As in the lattice case, this condition is sharp: if $\bar\beta(\xi)\sim\abs{\xi}^{-1-2s_2}$ as $\abs{\xi}\rightarrow+\infty$ with $s_2\leq\frac12$, then \eqref{eq:firstmoment} fails and the DLR equations can have more than one solution (see Section \ref{section:nonuniqueness}), in agreement with the phase transition of the one-dimensional Ising model with interaction $\abs{i-j}^{-1-2s_2}$ (see \cite{Dyson1969,FrohlichSpencer1982}).\\

\noindent For lattice systems with bounded spins, the proof of uniqueness relies on the fact that the conditional distributions in a bounded volume with two different boundary conditions have densities with respect to each other which are bounded uniformly with respect to the volume. In our setting this is false, since these conditional distributions are mutually singular (see the remark after Corollary \ref{corollary:Guassiantranslation}). We first translate one conditional measure so that it has the same exterior condition as the other, while leaving a smaller interval fixed. In centered Gaussian coordinates, the required correction is a cutoff harmonic extension belonging to the Cameron-Martin space (see Lemma \ref{lemma:localization} and Proposition \ref{proposition::cutoffdomain}). We control the resulting relative entropy on average with respect to the boundary conditions, by means of \eqref{eq:uniformmoments} and of an integrated estimate of the harmonic measure (see Lemma \ref{lemma:cutoff}). This gives the mutual absolute continuity of any two solutions satisfying \eqref{eq:uniformmoments} (see Proposition \ref{prop:equivalenceDLR}), and uniqueness follows by an argument on the tail $\sigma$-algebra (see \cite[Section 7.1]{Georgii2011} for its lattice counterpart).\\

\noindent More precisely, the proof consists of two parts:
\begin{enumerate}
\item We prove that any two solutions satisfying \eqref{eq:uniformmoments} are equivalent. The main estimate is a bound, uniform in $n\in \mathbb{N}$, on the average relative entropy between their conditional distributions restricted to the smaller interval $\mathcal I_n=(-n,n)$. The averaging is with respect to the two laws of the exterior conditions. Thus the moment assumption is used after averaging, and no bound uniform over all deterministic exterior conditions is needed.
\item We prove that this equivalence forces uniqueness. Conditioning a DLR solution on a tail event of probability strictly between zero and one produces two mutually singular DLR solutions that still satisfy \eqref{eq:uniformmoments}. Equivalence excludes such events. We then use the fact that the density of one DLR solution with respect to another is measurable with respect to the tail $\sigma$-algebra.
\end{enumerate}

\noindent The main results of this section are the following:

\begin{thm}\label{theorem:uniqueness}
Assume \ref{H0}, \ref{H1}, \ref{H2}, \ref{H3}, \ref{H4} and \eqref{eq:firstmoment} hold. Then there exists at most one solution of the DLR equations \eqref{eq:DLRequations} satisfying \eqref{eq:uniformmoments}.
\end{thm}

\begin{cor}\label{corollary:uniqueness}
Assume \ref{H0}, \ref{H1}, \ref{H2}, \ref{H3}, \ref{H4}, \ref{H5} and \eqref{eq:firstmoment} hold. Then the DLR equations \eqref{eq:DLRequations} have a unique solution $\mu$, which is also the unique reversible measure for \eqref{eq:SPDEmain1psuedo} among the Borel probability measures $\nu$ on $L^2_{r,k}$ such that $\nu(L^{2q}_{r,k}\cap C(\R))=1$. Moreover, $\mu$ is ergodic for $\{P(t)\}_{t\geq0}$, reflection positive and invariant with respect to translations.
\end{cor}

\begin{proof}
Let $\gamma\in(0,s_1-\frac12)$. For every $\tau\in\R$ and $x\in C^\gamma_{\rho^\tau_{r/2,k}}$ we have $\abs{x(\tau)}=\rho^\tau_{r/2,k}(\tau)\abs{x(\tau)}\leq\norm{x}_{C^\gamma_{\rho^\tau_{r/2,k}}}$, hence, by Theorem \ref{theorem:FerniqueDLR}, every solution of the DLR equations satisfies \eqref{eq:uniformmoments}. Therefore, by Theorems \ref{theorem:existence} and \ref{theorem:uniqueness}, the DLR equations have a unique solution $\mu$, which, by Theorem \ref{theorem:super-main}, is the unique reversible measure. By Theorem \ref{thm:apriori} and Corollary \ref{corollary:DLRuniqueness}, $\mu$ is ergodic and, by the theorem in Subsection \ref{section:ref.positive}, it is reflection positive. Finally, since $H^\beta$ commutes with translations and $V$ does not depend on $\xi$, by Corollary \ref{remark:DLRtranslated} the image of $\mu$ under the map $x\mapsto x(\cdot-a)$, with $a\in\R$, solves the DLR equations, and so it coincides with $\mu$.
\end{proof}

\subsection{Absolute continuity of solutions of the DLR equations}

In this subsection we are going to prove the following proposition.

\begin{prop}\label{prop:equivalenceDLR}
Assume \eqref{eq:firstmoment} holds and let $\mu_1$ and $\mu_2$ be two solutions of the DLR equations \eqref{eq:DLRequations} satisfying the moment condition \eqref{eq:uniformmoments}. Then $\mu_1$ and $\mu_2$ are equivalent.
\end{prop}

\noindent
Consider $\mu_1$ and $\mu_2$ satisfying the hypotheses of Proposition \ref{prop:equivalenceDLR}. Throughout this subsection, for each $n\in\N$, we write $\mathcal I:=\mathcal I_{n+1}$ and let
\[
\pi:=\mu_1^{\mathcal I^c_{n+1}}\otimes\mu_2^{\mathcal I^c_{n+1}},
\]
where $\mu_i^{\mathcal I^c_{n+1}}$ is the law of the restriction to $\mathcal I^c_{n+1}$ defined in \eqref{eq:legge-Ic}. Here the dependence of $\pi$ on $n$ is omitted from the notation. Thus, $\pi$ describes two independent exterior conditions, sampled from $\mu_1$ and $\mu_2$, respectively. By Corollary \ref{remark:DLRtranslated}, \eqref{eq:invmeasureI} and \eqref{eq:condvsinv}, for $\pi$-almost every $(y_1,y_2)$ and $i=1,2$,
\[
\mu_i^{\mathcal I_{n+1},y_i}=(T_{\varphi^{y_i}})_*\nu_{\mathcal I_{n+1},y_i},
\qquad T_u x:=x+u,
\]
where elements of $L^2(\mathcal I_{n+1})$ are extended by zero in $\mathcal I_{n+1}^c$. Here $y_i\in L^2_{r,k}(\mathcal I_{n+1}^c)\cap C(\mathcal I_{n+1}^c)$, and $\varphi^{y_i}$ is its harmonic extension.\\

\noindent For probability measures $P$ and $Q$ on the same measurable space, we recall that the relative entropy of $P$ with respect to $Q$ is
\begin{equation}\label{eq:entropy}
\operatorname{Ent}(P\Vert Q):=
\begin{cases}
\displaystyle\int\log\frac{\mathd P}{\mathd Q}\,\mathd P,
&P\ll Q,\\[1mm]
+\infty,&\text{otherwise}.
\end{cases}
\end{equation}
Since $\varphi^{y_1}\not = \varphi^{y_2}$ on the boundary of $\mathcal{I}_{n+1}$, when $y_1 \not=y_2$ on the boundary of $\mathcal{I}_{n+1}$, the measures $\mu_1^{\mathcal I,y_1}$ and $\mu_2^{\mathcal I,y_2}$ are, generically, supported on disjoint sets of full paths. Since we want to compare $\mu_1^{\mathcal I,y_1}$ and $\mu_2^{\mathcal I,y_2}$ using relative entropy, we compare their restrictions to $\mathcal I_n$, which is separated from $\mathcal I^c_{n+1}$ by a region of fixed width. Thus we denote by $\mu_i^{\mathcal I_{n+1},y_i}|_{\mathcal I_n}:=(R_{\mathcal I_n})_*\mu_i^{\mathcal I_{n+1},y_i}$ the image measure under the restriction map $R_{\mathcal I_n}$.

\begin{lemma}\label{lemma:supentropy}
Under the hypotheses of Proposition \ref{prop:equivalenceDLR} we have
\begin{align}\label{eq:uniqueness-conditional-entropy}
\sup_{n \in \mathbb{N}}\int\Big[&\operatorname{Ent}\big(
\mu_1^{\mathcal I_{n+1},y_1}|_{\mathcal I_n}\Vert
\mu_2^{\mathcal I_{n+1},y_2}|_{\mathcal I_n}\big)\notag\\
&+\operatorname{Ent}\big(
\mu_2^{\mathcal I_{n+1},y_2}|_{\mathcal I_n}\Vert
\mu_1^{\mathcal I_{n+1},y_1}|_{\mathcal I_n}\big)\Big]
\,\pi(\mathd y_1,\mathd y_2) < +\infty.
\end{align}
\end{lemma}

\noindent The proof of Lemma \ref{lemma:supentropy} requires several preliminary results. First of all we introduce some cutoff functions separating $\mathcal I_n$ from $\mathcal I_{n+1}$. We need to introduce such cutoffs since, for general continuous exterior conditions, the harmonic extension is only locally regular and its interior regularity need not extend to the boundary of the interval.\\

\noindent Let $\chi\in C^\infty(\R)$ be a nonincreasing function such that $\chi=1$ on $(-\infty,0]$ and $\chi=0$ on $[\frac12,+\infty)$, and for $n\in\N$ let $\chi_n(\xi):=\chi(\abs{\xi}-n)$, $\xi\in\R$. Then $\chi_n\in C^\infty_c(\mathcal I_{n+1})$, $0\leq\chi_n\leq1$ and $\chi_n=1$ on $[-n,n]$. Assumption \eqref{eq:firstmoment} is used only in the following lemma, which gives in \eqref{eq:cutoffharmonicmeasure} a bound, uniform with respect to $n$, on the integral, weighted by $\chi_n$, of the harmonic measure $P^{\mathcal I_{n+1}}_{\free}(\xi,\mathcal I_{n+1}^c)$, namely of the probability that the L\'evy process starting from $\xi$ leaves $\mathcal I_{n+1}$ before an independent exponential time of parameter $m^2$.

\begin{lemm}\label{lemma:cutoff}
Assume \ref{H0} and \eqref{eq:firstmoment} hold. Then
\begin{equation}\label{eq:cutoffL1}
\sup_{n\in\N}\norm{H^\beta\chi_n}_{L^1(\R)}<+\infty .
\end{equation}
Moreover, for every $n\in\N$ we have
\begin{align}
&\beta_0\int_\R\abs{\chi_n'(\xi)}^2\mathd\xi
+\frac12\int_{\R^2}\big(\chi_n(\xi)-\chi_n(\xi')\big)^2
\bar\beta(\xi-\xi')\,\mathd(\xi,\xi')
\leq\norm{H^\beta\chi_n}_{L^1(\R)},
\label{eq:cutoffenergy}\\
&m^2\int_{\mathcal I_{n+1}}\chi_n(\xi)
P^{\mathcal I_{n+1}}_{\free}(\xi,\mathcal I_{n+1}^c)\,\mathd\xi
\leq\norm{H^\beta\chi_n}_{L^1(\R)}.
\label{eq:cutoffharmonicmeasure}
\end{align}
\end{lemm}

\begin{remark}
In Lemma \ref{lemma:cutoff} the condition $\int_1^{+\infty}\xi \bar\beta(\xi)\mathd \xi <+\infty$ is used in an essential way. Indeed, it is possible to prove that $\sup_{n\in\N}\norm{H^\beta\chi_n}_{L^1(\R)}=+\infty$ when $\int_1^{+\infty}\xi \bar\beta(\xi)\mathd \xi=+\infty $.
\end{remark}

\begin{proof}
For $f\in C^\infty_c(\R)$ and $z\in\R$ we have
\[
2f(\xi)-f(\xi+z)-f(\xi-z)
=-\int_0^{z}\int_{-t}^{t}f''(\xi+s)\,\mathd s\,\mathd t
=-\int_0^z\big(f'(\xi+t)-f'(\xi-t)\big)\mathd t,
\qquad\xi\in\R,
\]
so that $\norm{2f-f(\cdot+z)-f(\cdot-z)}_{L^1(\R)}
\leq\min\big\{\abs{z}^2\norm{f''}_{L^1(\R)},
2\abs{z}\norm{f'}_{L^1(\R)}\big\}$.
Applying this estimate with $f=\chi_n$ in \eqref{eq:Hbeta}, and noting that $\norm{\chi_n'}_{L^1(\R)}=2\norm{\chi'}_{L^1(\R)}$ and $\norm{\chi_n''}_{L^1(\R)}=2\norm{\chi''}_{L^1(\R)}$, we obtain
\[
\norm{H^\beta\chi_n}_{L^1(\R)}
\leq(2\beta_0+\Theta_\beta)\norm{\chi''}_{L^1(\R)}
+2\norm{\chi'}_{L^1(\R)}
\int_{\{\abs{z}\geq1\}}\abs{z}\bar\beta(z)\,\mathd z,
\]
and \eqref{eq:cutoffL1} follows from \eqref{eq:levynear} and \eqref{eq:firstmoment}. By \eqref{eq:dirichletform}, the left-hand side of \eqref{eq:cutoffenergy} is equal to $\mathcal E(\chi_n,\chi_n)-m^2\norm{\chi_n}^2_{L^2}=\int_\R\chi_nH^\beta\chi_n\,\mathd\xi$ (see the proof of Lemma \ref{lemma:harmonicextensionenergy}), and \eqref{eq:cutoffenergy} follows since $0\leq\chi_n\leq1$. Finally, let $\varphi^1$ be the harmonic extension in $\mathcal I_{n+1}$ of the function equal to $1$ on $\mathcal I_{n+1}^c$. By \eqref{eq:harmonicextensionprob} and \eqref{eq:poissonkernel} we have $\varphi^1(\xi)=P^{\mathcal I_{n+1}}_{\free}(\xi,\mathcal I_{n+1}^c)\in[0,1]$ for $\xi\in\mathcal I_{n+1}$ and $\varphi^1=1$ on $\mathcal I_{n+1}^c$. Since the constant function $1$ belongs to $C(\mathcal I_{n+1}^c)\cap L^1_{r,k}(\mathcal I_{n+1}^c)$, by Theorem \ref{thm:harmonic-continuity} $\varphi^1$ is a weak harmonic extension, and choosing $h=\chi_n$ in point 3. of Definition \ref{def:harmonicextension} we obtain
\[
m^2\int_{\mathcal I_{n+1}}\chi_n\varphi^1\,\mathd\xi
=-\int_\R\varphi^1H^\beta\chi_n\,\mathd\xi
\leq\norm{H^\beta\chi_n}_{L^1(\R)},
\]
namely \eqref{eq:cutoffharmonicmeasure}.
\end{proof}

\noindent
Fix some $n\in \mathbb N$, and a pair of boundary conditions $(y_1,y_2)$ and set
\begin{equation}\label{eq:psi1}
\psi:=\varphi^{y_2}-\varphi^{y_1}=\varphi^{y_2-y_1},
\qquad \psi_{in}:=\chi_n\psi,
\qquad \psi_{out}:=(1-\chi_n)\psi.
\end{equation}
These functions satisfy
\begin{equation}\label{eq:psi2}
\psi=\psi_{in}+\psi_{out},\qquad
\psi_{out}=0\text{ on }\mathcal I_n,\qquad
\psi_{out}=y_2-y_1\text{ on }\mathcal I^c_{n+1},
\qquad \varphi^{y_1}+\psi_{out}=\varphi^{y_2}-\psi_{in}.
\end{equation}
Consequently, the translated measure
\[
\widetilde\mu^{\mathcal I_{n+1},y_1,y_2}
:=(T_{\psi_{out}})_*\mu_1^{\mathcal I_{n+1},y_1}
\]
has exterior condition $y_2$, and its restriction to $\mathcal I_n$ coincides with that of $\mu_1^{\mathcal I_{n+1},y_1}$. Thus it can be compared with $\mu_2^{\mathcal I_{n+1},y_2}$ without changing the inner marginal that interests us. Inside $\mathcal I_{n+1}$, the translation $\psi_{out}$ is supported in $\mathcal I_{n+1}\setminus\mathcal I_n$, whose Lebesgue measure is independent of $n$ (being two). In centered coordinates this translation is particularly simple. Indeed,
\[
\widetilde\mu^{\mathcal I_{n+1},y_1,y_2}
=(T_{\varphi^{y_2}})_*
\left[(T_{-\psi_{in}})_*\nu_{\mathcal I_{n+1},y_1}\right],
\qquad
\mu_2^{\mathcal I_{n+1},y_2}
=(T_{\varphi^{y_2}})_*\nu_{\mathcal I_{n+1},y_2}.
\]

\begin{lemma}\label{lemma:Sy1y2}
Under the hypotheses of Proposition \ref{prop:equivalenceDLR} we have
\begin{multline}\label{eq:definitionS1}
H_n(y_1,y_2):=
\operatorname{Ent}(\tilde{\mu}^{\mathcal I_{n+1},y_1,y_2}
\Vert \mu^{\mathcal I_{n+1},y_2}_2)
+\operatorname{Ent}(\mu^{\mathcal I_{n+1},y_2}_2
\Vert\tilde{\mu}^{\mathcal I_{n+1},y_1,y_2})\\
=\mathcal E(\psi_{in},\psi_{in})
+\int\!\int_{\mathcal I_{n+1}\setminus\mathcal I_n}
\Big[V\big(x(\xi)+\psi_{out}(\xi)\big)-V\big(x(\xi)\big)\Big]
\mathd\xi\,\mu_1^{\mathcal I_{n+1},y_1}(\mathd x)\\
+\int\!\int_{\mathcal I_{n+1}\setminus\mathcal I_n}
\Big[V\big(x(\xi)-\psi_{out}(\xi)\big)-V\big(x(\xi)\big)\Big]
\mathd\xi\,\mu_2^{\mathcal I_{n+1},y_2}(\mathd x)
+\int_{\mathcal I_{n+1}}\psi_{in}(\xi)
\left(b^{y_1,n}_1(\xi)-b^{y_2,n}_2(\xi)\right)\mathd\xi
\end{multline}
where the functions $b_i^{y_i,n}:\mathcal{I}_{n+1} \rightarrow \R$ are defined as
\begin{equation}\label{eq:bi}
b_i^{y_i,n}(\xi)=\int V'(x(\xi))
\mu^{\mathcal I_{n+1},y_i}_{i}(\mathd x),\quad i=1,2.
\end{equation}
\end{lemma}

\begin{proof}
We first prove that
\begin{equation}\label{eq:Hn1}
H_n(y_1,y_2)
=\int G^{y_1,y_2}(x)
[(T_{-\psi_{in}})_*\nu_{\mathcal I_{n+1},y_1}](\mathd x)
-\int G^{y_1,y_2}(x)\,
\nu_{\mathcal I_{n+1},y_2}(\mathd x),
\end{equation}
where $G^{y_1,y_2}$ is given by the expression \eqref{eq:uniqueness-logdensity} below.\\

\noindent In order to obtain expression \eqref{eq:Hn1}, we first prove that $\psi_{in}$ belongs to the domain of $\mathcal{A}_{\mathcal{I}_{n+1}}$, i.e. the realization of $H^{\beta}+m^2$ on $\mathcal{I}_{n+1}$ with homogeneous Dirichlet exterior conditions. Since $\rho_{r,k}\in L^1(\R)$, $y_2-y_1\in L^1_{r,k}(\mathcal I^c)\cap C(\mathcal I^c)$. Lemma \ref{lemma:localization} and Proposition \ref{proposition::cutoffdomain} applied to $\chi_n$ imply that the function $\psi_{in}$ (defined on the whole $\R$ and equal to $0$ on $\mathcal{I}_{n+1}^c$) satisfies
\[
\psi_{in}=\chi_n\psi\in\mathcal F_{\mathcal {I}_{n+1}}
\cap \operatorname{Dom}(H^\beta),
\]
where $\operatorname{Dom}(H^\beta)$ is the domain of the operator $H^{\beta}$ on the whole $\R$. This is the main reason for using the harmonic difference and the cutoff in the definition of the translated measure $\tilde{\mu}^{\mathcal I_{n+1},y_1,y_2}$.\\

\noindent The space $\mathcal F_{\mathcal I_{n+1}}$ is the form domain of $\mathcal A_{\mathcal I_{n+1}}$, and
\[
\scal{\mathcal A_{\mathcal I}u}{u}_{L^2(\mathcal I)}
=\mathcal E(u,u),\qquad u\in{\rm Dom}(\mathcal A_{\mathcal I});
\]
see \cite[Theorem 4.4.3]{Fukushima2011}. Since the covariance of $\mu^{\mathcal I_{n+1},0}_{\free}$ is $\mathcal A_{\mathcal I_{n+1}}^{-1}$, by Proposition \ref{prop:gaussL2I}, its Cameron-Martin space is precisely $\mathcal F_{\mathcal I_{n+1}}$, with squared norm $\mathcal E(u,u)$. One can furthermore check that
\[
\mathcal{F}_{\mathcal{I}_{n+1}}\cap\operatorname{Dom}(H^\beta)
\subset\operatorname{dom}(\mathcal{A}_{\mathcal{I}_{n+1}}),
\]
and deduce that the Gaussian linear functional associated with $\psi_{in}$ is represented by the continuous linear map $x\mapsto\scal{\mathcal A_{\mathcal I_{n+1}}\psi_{in}}{x}_{L^2(\mathcal I_{n+1})}$, which therefore makes sense pathwise and coincides with the Skorokhod integral of the vector $\psi_{in}$ in the Cameron-Martin space. The Cameron-Martin formula therefore gives
\begin{equation}\label{eq:CMformula}
\frac{\mathd(T_{-\psi_{in}})_*\mu^{\mathcal I_{n+1},0}_{\free}}
{\mathd\mu^{\mathcal I_{n+1},0}_{\free}}(x)
=\exp\left(
-\scal{\mathcal A_{\mathcal I_{n+1}}\psi_{in}}{x}_{L^2(\mathcal I_{n+1})}
-\frac12\mathcal E(\psi_{in},\psi_{in})\right);
\end{equation}
see \cite[Corollary 2.4.3]{Bogachev98}.
Combining this formula with \eqref{eq:invmeasureI}, we find that $(T_{-\psi_{in}})_*\nu_{\mathcal I,y_1}$ and $\nu_{\mathcal I,y_2}$ are equivalent and
\begin{align}\label{eq:uniqueness-logdensity}
G^{y_1,y_2}(x)
:=\log\frac{\mathd [(T_{-\psi_{in}})_*\nu_{\mathcal I_{n+1},y_1}]}
{\mathd\nu_{\mathcal I_{n+1},y_2}}(x)
={}&\log\frac{\mathcal Z^{\mathcal I_{n+1}}(y_2)}
{\mathcal Z^{\mathcal I_{n+1}}(y_1)}
-\scal{\mathcal A_{\mathcal I_{n+1}}\psi_{in}}{x}_{L^2(\mathcal I_{n+1})}
-\frac12\mathcal E(\psi_{in},\psi_{in})\notag\\
&-\int_{\mathcal I_{n+1}}\Big[
V\big(x(\xi)+\varphi^{y_2}(\xi)-\psi_{out}(\xi)\big)
-V\big(x(\xi)+\varphi^{y_2}(\xi)\big)
\Big]\mathd\xi.
\end{align}
Indeed, at a point $x$ the density of the translated measure $(T_{-\psi_{in}})_*\nu_{\mathcal I,y_1}$ contains the Gibbs weight of $\nu_{\mathcal I,y_1}$ evaluated at $x+\psi_{in}$, and
$x+\psi_{in}+\varphi^{y_1}=x+\varphi^{y_2}-\psi_{out}$.
For fixed $n,y_1,y_2$, the functions $\varphi^{y_i}$ and $\psi_{in},\psi_{out}$ are bounded on $\overline{\mathcal I_{n+1}}$. Moreover, by \ref{H4}, $\nu_{\mathcal I_{n+1},y_i}$ has a bounded density with respect to $\mu^{\mathcal I_{n+1},0}_{\free}$, which has finite moments of every order in $C(\overline{\mathcal I_{n+1}})$, as recalled in the proof of Proposition \ref{prop:invmeasureI}. Together with \eqref{H2:pol}, these facts imply that $G^{y_1,y_2}$ is integrable with respect to both $(T_{-\psi_{in}})_*\nu_{\mathcal I_{n+1},y_1}$ and $\nu_{\mathcal I_{n+1},y_2}$. Using the definition of relative entropy \eqref{eq:entropy} we obtain the expression \eqref{eq:Hn1}.\\

\noindent In order to prove the equality \eqref{eq:definitionS1} from equality \eqref{eq:Hn1} we use relations \eqref{eq:psi1} and \eqref{eq:psi2}, and the integration by parts formula \eqref{eq:ibpI}.\\

\noindent First we note that, by the definition of the translation map $T$, we get
\begin{multline}\label{eq:Gy1y2}
\int G^{y_1,y_2}(x)
[(T_{-\psi_{in}})_*\nu_{\mathcal I_{n+1},y_1}](\mathd x)\\
=\log\frac{\mathcal Z^{\mathcal I_{n+1}}(y_2)}
{\mathcal Z^{\mathcal I_{n+1}}(y_1)}
-\frac12\mathcal E(\psi_{in},\psi_{in})
-\int \scal{\mathcal A_{\mathcal I_{n+1}}\psi_{in}}
{x-\psi_{in}}_{L^2(\mathcal I_{n+1})}
\nu_{\mathcal{I}_{n+1},y_1}(\mathd x)\\
-\int\left(\int_{\mathcal I_{n+1}}\Big[
V\big(x(\xi)+\varphi^{y_2}(\xi)-\psi_{out}(\xi)-\psi_{in}(\xi)\big)
-V\big(x(\xi)+\varphi^{y_2}(\xi)-\psi_{in}(\xi)\big)
\Big]\mathd \xi\right)\nu_{\mathcal{I}_{n+1},y_1}(\mathd x).
\end{multline}
Since $\psi_{in}\in\operatorname{Dom}(\mathcal{A}_\mathcal{I})$ we have by linearity that
\[
\scal{\mathcal A_{\mathcal I_{n+1}}\psi_{in}}
{x-\psi_{in}}_{L^2(\mathcal I_{n+1})}
=\scal{\mathcal A_{\mathcal I_{n+1}}\psi_{in}}
{x}_{L^2(\mathcal I_{n+1})}
-\mathcal E(\psi_{in},\psi_{in}).
\]
Thus, by the previous equation, by using the relations \eqref{eq:psi1} and \eqref{eq:psi2}, by observing that $\psi_{out}=0$ on $\mathcal I_n$, and by exploiting the fact that $\mu^{\mathcal{I}_{n+1},y_1}_{1}=(T_{\varphi^{y_1}})_*(\nu_{\mathcal{I}_{n+1},y_1})$, we obtain
\begin{multline}\label{eq:Gy1y2bis}
\int G^{y_1,y_2}(x)
[(T_{-\psi_{in}})_*\nu_{\mathcal I_{n+1},y_1}](\mathd x)\\
=\log\frac{\mathcal Z^{\mathcal I_{n+1}}(y_2)}
{\mathcal Z^{\mathcal I_{n+1}}(y_1)}
+\frac12\mathcal E(\psi_{in},\psi_{in})
-\int \scal{\mathcal A_{\mathcal I_{n+1}}\psi_{in}}
{x}_{L^2(\mathcal I_{n+1})}
\nu_{\mathcal{I}_{n+1},y_1}(\mathd x)\\
+\int\!\int_{\mathcal I_{n+1}\setminus\mathcal I_n}
\Big[V\big(x(\xi)+\psi_{out}(\xi)\big)-V\big(x(\xi)\big)\Big]
\mathd\xi\,\mu_1^{\mathcal I_{n+1},y_1}(\mathd x).
\end{multline}
In a similar way we can prove that
\begin{multline}\label{eq:Gy1y2ter}
\int G^{y_1,y_2}(x)\,\nu_{\mathcal I_{n+1},y_2}(\mathd x)
=\log\frac{\mathcal Z^{\mathcal I_{n+1}}(y_2)}
{\mathcal Z^{\mathcal I_{n+1}}(y_1)}
-\frac12\mathcal E(\psi_{in},\psi_{in})
-\int \scal{\mathcal A_{\mathcal I_{n+1}}\psi_{in}}
{x}_{L^2(\mathcal I_{n+1})}
\nu_{\mathcal{I}_{n+1},y_2}(\mathd x)\\
-\int\!\int_{\mathcal I_{n+1}\setminus\mathcal I_n}
\Big[V\big(x(\xi)-\psi_{out}(\xi)\big)-V\big(x(\xi)\big)\Big]
\mathd\xi\,\mu_2^{\mathcal I_{n+1},y_2}(\mathd x).
\end{multline}
It remains to apply the integration by parts formula for the terms involving the integral of $\langle \mathcal{A}_{\mathcal{I}_{n+1}}\psi_{in},x\rangle_{L^2(\mathcal{I}_{n+1})}$. Indeed by using \eqref{eq:ibpI}, with $\Phi\equiv1$ in that formula, for $u\in{\rm Dom}(\mathcal A_{\mathcal I_{n+1}})$, we have
\begin{multline}\label{eq:integrationbypartsused}
\int\scal{\mathcal A_{\mathcal I_{n+1}}u}{x}_{L^2(\mathcal I_{n+1})}
\nu_{\mathcal I_{n+1},y}(\mathd x)
=-\int\scal{u}{V'(x+\varphi^y)}_{L^2(\mathcal I_{n+1})}
\nu_{\mathcal I_{n+1},y}(\mathd x)\\
=-\int\scal{u}{V'(x)}_{L^2(\mathcal I_{n+1})}
\mu^{\mathcal I_{n+1},y}(\mathd x)
\end{multline}
Using the definition of $b_i^{y_i,n}$ in \eqref{eq:bi}, and integration by parts for $u=\psi_{in}\in\operatorname{Dom}(\mathcal{A}_{\mathcal{I}_{n+1}})$ yields
\[
\int\scal{\mathcal A_{\mathcal I_{n+1}}\psi_{in}}{x}_{L^2(\mathcal I_{n+1})}
\nu_{\mathcal I_{n+1},y_i}(\mathd x)
=-\int_{\mathcal I_{n+1}}\psi_{in}(\xi)b_i^{y_i,n}(\xi)\mathd\xi.
\]
Combining the preceding identities we obtain the expression \eqref{eq:definitionS1}.
\end{proof}

\begin{proof}[Proof of Lemma \ref{lemma:supentropy}]
We first observe that, since $\tilde\mu^{\mathcal I_{n+1},y_1,y_2}|_{\mathcal I_n}=\mu^{\mathcal{I}_{n+1},y_1}_1|_{\mathcal{I}_{n}}$, by the monotonicity of relative entropy under restriction to a sub-$\sigma$-algebra, we have
\begin{align*}
\operatorname{Ent}(
\mu^{\mathcal{I}_{n+1},y_1}_1|_{\mathcal{I}_{n}}
\Vert\mu^{\mathcal{I}_{n+1},y_2}_2|_{\mathcal{I}_{n}})
&\leqslant\operatorname{Ent}(
\tilde\mu^{\mathcal{I}_{n+1},y_1,y_2}
\Vert\mu^{\mathcal{I}_{n+1},y_2}_2),\\
\operatorname{Ent}(
\mu^{\mathcal{I}_{n+1},y_2}_2|_{\mathcal{I}_{n}}
\Vert\mu^{\mathcal{I}_{n+1},y_1}_1|_{\mathcal{I}_{n}})
&\leqslant\operatorname{Ent}(
\mu^{\mathcal{I}_{n+1},y_2}_2
\Vert\tilde\mu^{\mathcal{I}_{n+1},y_1,y_2}).
\end{align*}
This implies that
\begin{equation}\label{eq:EntHn}
\operatorname{Ent}\big(
\mu_1^{\mathcal I_{n+1},y_1}|_{\mathcal I_n}\Vert
\mu_2^{\mathcal I_{n+1},y_2}|_{\mathcal I_n}\big)
+\operatorname{Ent}\big(
\mu_2^{\mathcal I_{n+1},y_2}|_{\mathcal I_n}\Vert
\mu_1^{\mathcal I_{n+1},y_1}|_{\mathcal I_n}\big)
\leqslant H_n(y_1,y_2).
\end{equation}
Thus, in order to prove Lemma \ref{lemma:supentropy}, it is sufficient to prove $\int H_n(y_1,y_2) \pi(\mathd y_1, \mathd y_2) \leqslant K < +\infty$, for some $K \in \R_+$ (independent of $n\in \mathbb{N}$).\\

\noindent Let $\Lambda$ be the maximum of the left-hand side of \eqref{eq:uniformmoments} with $\mu=\mu_1$ and $\mu=\mu_2$. Throughout the proof, $\lesssim$ denotes an inequality up to a positive constant depending only on $\Lambda$, $m$, $q$, the constant in \eqref{H2:pol} and the uniform bound in \eqref{eq:cutoffL1}. In particular, these constants are independent of $n$ and of the exterior conditions.\\

\noindent For $1\leq p\leq2q$, \eqref{eq:uniformmoments} implies
\[
\sup_{\xi\in\mathcal I^c_{n+1}}
\left(\int\abs{y_2(\xi)-y_1(\xi)}^p\mathd\pi\right)^{1/p}
\leq2\Lambda^{1/(2q)}.
\]
By the Poisson representation of the harmonic extension,
\[
\psi(\xi)=\int_{\mathcal I^c_{n+1}}
\big(y_2(\xi')-y_1(\xi')\big)
P^{\mathcal I}_{\free}(\xi,\mathd\xi'),
\qquad \xi\in\mathcal I,
\]
and thus we obtain
\begin{align}\label{eq:psimoments}
\left(\int\abs{\psi(\xi)}^p\mathd\pi\right)^{1/p}
&\leq2\Lambda^{1/(2q)}
P^{\mathcal I}_{\free}(\xi,\mathcal I^c_{n+1}),
&&\xi\in\mathcal I,\notag\\
\left(\int\abs{\psi(\xi)}^p\mathd\pi\right)^{1/p}
&\leq2\Lambda^{1/(2q)},&&\xi\in\R.
\end{align}
Here we used $P^{\mathcal I}_{\free}(\xi,\mathcal I^c_{n+1})\leq1$ and $\psi=y_2-y_1$ on $\mathcal I^c_{n+1}$. The first estimate retains the decay of the harmonic measure away from the exterior; this decay is needed for the term involving $b_1-b_2$.\\

\noindent We estimate the terms in \eqref{eq:definitionS1} separately. First, Lemma \ref{lemma:localization}, Cauchy-Schwarz and \eqref{eq:psimoments} yield
\begin{equation}\label{eq:ineEnt1}
\int\mathcal E(\psi_{in},\psi_{in})\mathd\pi
\leqslant 4\Lambda^{1/q}\left[
\beta_0\int_\R\abs{\chi_n'(\xi)}^2\mathd\xi\right.
\left.+\frac12\int_{\R^2}
\big(\chi_n(\xi)-\chi_n(\xi')\big)^2
\bar\beta(\xi-\xi')\mathd(\xi,\xi')\right]
\lesssim1,
\end{equation}
where the last inequality follows from \eqref{eq:cutoffenergy} and \eqref{eq:cutoffL1}. Harmonicity is essential here: it replaces a possible contribution from the whole interval by the cutoff energy, which is bounded independently of its length.\\

\noindent Second, conditional Jensen's inequality, the DLR disintegration and \eqref{H2:pol} give, for $i=1,2$ and $\xi\in\mathcal I$,
\begin{equation}\label{eq:ineEnt2}
\int\abs{b_i^{y_i,n}(\xi)}^2\mathd\pi
\leqslant\int\!\int\abs{V'(x(\xi))}^2
\mu_i^{\mathcal I_{n+1},y_i}(\mathd x)\,
\mu_i^{\mathcal I_{n+1}^c}(\mathd y_i)
=\int\abs{V'(x(\xi))}^2\mu_i(\mathd x)
\lesssim 1.
\end{equation}
Consequently, Cauchy-Schwarz on the boundary probability space and the first bound in \eqref{eq:psimoments} imply
\begin{equation}\label{eq:ineEnt3}
\int\!\int_{\mathcal I_{n+1}}
\abs{\psi_{in}(\xi)}
\abs{b_1^{y_1,n}(\xi)-b_2^{y_2,n}(\xi)}
\mathd\xi\,\mathd\pi(\mathd y_1, \mathd y_2)
\lesssim\int_{\mathcal I_{n+1}}
\chi_n(\xi)P^{\mathcal I}_{\free}(\xi,\mathcal I^c)\mathd\xi
\lesssim1,
\end{equation}
by \eqref{eq:cutoffharmonicmeasure}.\\

\noindent Finally, \eqref{H2:pol} and the fundamental theorem of calculus imply
\[
\abs{V(a+c)-V(a)}
\lesssim1+\abs{a}^{q+1}+\abs{c}^{q+1},
\qquad a,c\in\R.
\]
Since $q+1\leq2q$ and $\abs{\psi_{out}}\leq\abs{\psi}$, the DLR disintegration and \eqref{eq:psimoments} give
\begin{multline}\label{eq:ineEnt4}
\int\!\int\!\int_{\mathcal I_{n+1}\setminus\mathcal I_n}
\abs{V(x(\xi)+\psi_{out}(\xi))-V(x(\xi))}
\mathd\xi\,\mu_1^{\mathcal I_{n+1},y_1}(\mathd x)\,\mathd\pi\\
\lesssim\int_{\mathcal I_{n+1}\setminus\mathcal I_n}
\left[1+\int\abs{x(\xi)}^{q+1}\mu_1(\mathd x)
+\int\abs{\psi(\xi)}^{q+1}\mathd\pi\right]\mathd\xi
\lesssim 1.
\end{multline}
The last inequality uses that $\mathcal I_{n+1} \setminus\mathcal I_n$ has measure two. The same argument applies to the second potential-difference term in \eqref{eq:definitionS1}, with $-\psi_{out}$ and $\mu_2$.\\

\noindent Thus, putting together inequality \eqref{eq:EntHn}, the expression of $H_n$ in Lemma \ref{lemma:Sy1y2}, and inequalities \eqref{eq:ineEnt1}, \eqref{eq:ineEnt3} and \eqref{eq:ineEnt4}, we obtain the conclusion.
\end{proof}

\begin{proof}[Proof of Proposition \ref{prop:equivalenceDLR}]
In order to prove Proposition \ref{prop:equivalenceDLR}, it is sufficient to prove that the sum of the two relative entropies is finite, namely $\operatorname{Ent}(\mu_1\Vert\mu_2)+\operatorname{Ent}(\mu_2\Vert\mu_1)<+\infty$.\\
Since the $\sigma$-algebras generated by the restrictions to $\mathcal I_n$ increase to the full Borel $\sigma$-algebra on $L^2_{r,k}$, by continuity of relative entropy along increasing $\sigma$-algebras (see, e.g., Proposition (15.6) of \cite{Georgii2011}) we have
\begin{equation}\label{eq:entrp1}
\operatorname{Ent}(\mu_1\Vert\mu_2)
+\operatorname{Ent}(\mu_2\Vert\mu_1)
=\lim_{n \rightarrow +\infty}
\left(\operatorname{Ent}(
\mu_1|_{\mathcal{I}_{n}}\Vert\mu_2|_{\mathcal{I}_{n}})
+\operatorname{Ent}(
\mu_2|_{\mathcal{I}_{n}}\Vert\mu_1|_{\mathcal{I}_{n}})
\right)
\end{equation}
Furthermore, by using the fact that
\[
\int \mu^{\mathcal{I}_{n+1},y_i}_{i}|_{\mathcal{I}_{n}}\,
\mu_i^{\mathcal{I}_{n+1}^c}(\mathd y_i)
=\mu_{i}|_{\mathcal{I}_{n}},
\]
and the joint convexity of relative entropy (see, e.g., Lemma 2.4 of \cite{BudhirajaDupuis2019}), we have
\begin{multline}\label{eq:entrp2}
\operatorname{Ent}(
\mu_1|_{\mathcal{I}_{n}}\Vert\mu_2|_{\mathcal{I}_{n}})
+\operatorname{Ent}(
\mu_2|_{\mathcal{I}_{n}}\Vert\mu_1|_{\mathcal{I}_{n}})\\
\leqslant \int \left(
\operatorname{Ent}(
\mu_1^{\mathcal{I}_{n+1},y_1}|_{\mathcal{I}_{n}}
\Vert\mu_2^{\mathcal{I}_{n+1},y_2}|_{\mathcal{I}_{n}})
+\operatorname{Ent}(
\mu_2^{\mathcal{I}_{n+1},y_2}|_{\mathcal{I}_{n}}
\Vert\mu_1^{\mathcal{I}_{n+1},y_1}|_{\mathcal{I}_{n}})
\right)\pi(\mathd y_1, \mathd y_2).
\end{multline}
Thus using Lemma \ref{lemma:supentropy} together with equation \eqref{eq:entrp1} and inequality \eqref{eq:entrp2} the conclusion follows.
\end{proof}

\subsection{From absolute continuity to uniqueness}

Thanks to Proposition \ref{prop:equivalenceDLR}, the assertion of Theorem \ref{theorem:uniqueness} follows from the following result.

\begin{proposition}\label{proposition:equivalenceDLR2}
If any two solutions of the DLR equations satisfying the moment condition \eqref{eq:uniformmoments} are equivalent, then there is at most one such solution.
\end{proposition}

\noindent In order to prove Proposition \ref{proposition:equivalenceDLR2}, we need the following lemma. First, we recall that $R_{\mathcal{I}}$ denotes the map from $C(\R)$ into $C(\overline{\mathcal{I}})$ given by the restriction. We define the tail $\sigma$-algebra as
\[
\mathscr T:=\bigcap_{m\geq1}\sigma(R_{\mathcal I_m^c}).
\]
If $\mu_1$ and $\mu_2$ are two equivalent solutions of the DLR equations, set
\begin{equation}\label{eq:frd}
f=\frac{\mathd \mu_2}{\mathd \mu_1}.
\end{equation}

\begin{lemma}\label{lemma:tailsigma1}
Using the previous notation, if $\mu_1$ and $\mu_2$ are two equivalent probability measures solving the DLR equations, then there is a version of the Radon-Nikodym derivative $f$ which is measurable with respect to the tail $\sigma$-algebra $\mathscr T$.
\end{lemma}

\begin{proof}
Since $\mu_1$ and $\mu_2$ are equivalent solutions of the DLR equations, their exterior laws are equivalent and, for any bounded open interval $\mathcal{I}$, we have $\mu_1^{\mathcal I,y}=\mu_2^{\mathcal I,y}$ for $\mu_1^{\mathcal I^c}$-almost every $y$. Writing $f_n=\frac{\mathd \mu_2^{\mathcal I_n^c}}{\mathd \mu_1^{\mathcal I_n^c}}$, for every bounded Borel function $\Phi$ we have
\begin{multline}\label{eq:uniquenessequivalent1}
\int \Phi(x) f(x) \mu_1(\mathd x)
=\int \Phi(x) \mu_2(\mathd x)
=\int \left(\int\Phi(x)
\mu_2^{\mathcal{I}_n,y}(\mathd x)\right)
\mu_2^{\mathcal I^c_n}(\mathd y)\\
=\int \left(\int\Phi(x)
\mu_1^{\mathcal{I}_n,y}(\mathd x)\right)
f_n(y)\mu_1^{\mathcal I^c_n}(\mathd y)
=\int \Phi(x) f_n(R_{\mathcal{I}_n^c}(x))
\mu_1(\mathd x).
\end{multline}
Since equation \eqref{eq:uniquenessequivalent1} holds for any bounded Borel function $\Phi$ we have that $f(x)=f_n(R_{\mathcal I_n^c}(x))$ $\mu_1$-almost surely. This means that
\begin{equation}\label{eq:uniquenessequivalent2}
f=\limsup_{n \rightarrow +\infty}
f_n(R_{\mathcal I_n^c}(\cdot))
\end{equation}
$\mu_1$-almost surely. For every fixed $m\geq1$, the functions $f_n(R_{\mathcal I_n^c}(\cdot))$ with $n\geq m$ are measurable with respect to $\sigma(R_{\mathcal I_m^c})$. Hence the limsup in \eqref{eq:uniquenessequivalent2} is measurable with respect to every $\sigma(R_{\mathcal I_m^c})$, and therefore defines a version of $f$ measurable with respect to the tail $\sigma$-algebra $\mathscr T$.
\end{proof}

\begin{lemma}\label{lemma:tailsigma2}
Under the hypothesis of Proposition \ref{proposition:equivalenceDLR2}, let $\lambda$ be a solution of the DLR equations satisfying \eqref{eq:uniformmoments}. Then $\lambda$ is trivial on the tail $\sigma$-algebra $\mathscr T$.
\end{lemma}

\begin{proof}
Let $A \in \mathscr T$. We want to prove that $\lambda(A) \in \{0,1\}$. Indeed, suppose that there is $A \in \mathscr T$ such that $0<\lambda(A)<1$, and define
\[
\lambda_A:=\frac{\lambda(\cdot\cap A)}{\lambda(A)},
\qquad
\lambda_{A^c}:=\frac{\lambda(\cdot\cap A^c)}{\lambda(A^c)}.
\]
These measures inherit the support requirements of the DLR definition. They also satisfy the required moment bound, since
\[
\sup_{\xi\in\R}\int\abs{x(\xi)}^{2q}\lambda_A(\mathd x)
\leq\frac{1}{\lambda(A)}
\sup_{\xi\in\R}\int\abs{x(\xi)}^{2q}\lambda(\mathd x)
<+\infty,
\]
and similarly for $\lambda_{A^c}$.\\

\noindent Both $\lambda_A$ and $\lambda_{A^c}$ are solutions of the DLR equations. Indeed, given a bounded open interval $\mathcal I$, choose $m$ such that $\mathcal I\subseteq\mathcal I_m$. Since $A\in\sigma(R_{\mathcal I_m^c})\subseteq\sigma(R_{\mathcal I^c})$, there is a Borel set $A'$ of exterior conditions such that $A=R_{\mathcal I^c}^{-1}(A')$. For every bounded Borel function $\Phi$,
\[
\int\Phi\,\mathd\lambda_A
=\frac{1}{\lambda(A)}\int_{A'}
\left(\int\Phi\,\mathd\lambda^{\mathcal I,y}\right)
\lambda^{\mathcal I^c}(\mathd y).
\]
The exterior law of $\lambda_A$ is
$\mathbb I_{A'}(y)\lambda^{\mathcal I^c}(\mathd y)/\lambda(A)$,
so its conditional measures are still $\lambda^{\mathcal I,y}$ and satisfy \eqref{eq:DLRequations}. The same argument applies to $\lambda_{A^c}$.\\

\noindent By hypothesis $\lambda_A$ and $\lambda_{A^c}$ are equivalent. This is impossible, because one gives mass one to $A$ and the other gives mass zero to $A$. We conclude that every event in $\mathscr T$ has $\lambda$-probability zero or one.
\end{proof}

\begin{proof}[Proof of Proposition \ref{proposition:equivalenceDLR2}]
By hypothesis, if $\mu_1$ and $\mu_2$ are two solutions of the DLR equations satisfying the moment condition \eqref{eq:uniformmoments}, then $\mu_1$ and $\mu_2$ are equivalent, with Radon-Nikodym derivative $f$.\\
By Lemma \ref{lemma:tailsigma1} $f$ can be chosen to be measurable with respect to $\mathscr T$. By Lemma \ref{lemma:tailsigma2} $f$ must be $\mu_1$-almost surely a constant. Since $\int f(x)\mu_1(\mathd x)=1$, $f=1$ $\mu_1$-almost surely, and thus $\mu_1=\mu_2$.
\end{proof}

\section{On non-uniqueness}\label{section:nonuniqueness}

In this section, we discuss a family of potentials $V$ and operators $H^{\beta}$ for which it is possible to prove the existence of multiple solutions to the DLR equations \eqref{eq:DLRequations} satisfying the integrability condition \eqref{eq:momenti-DLR}. This proves the existence of multiple invariant measures for equation \eqref{eq:LangevinSPDE1} for some potentials and linear operators. The proof is essentially an extension of the argument of \cite[Theorem 6.1]{OsadaSpohn1999}; with respect to \cite{OsadaSpohn1999}, the L\'evy density $\bar\beta$ can be unbounded near the origin, as in the case of the fractional Laplacian of order $s\in(0,\frac12)$, and this requires a different treatment of the interaction at short distances (see Section \ref{section:OScomparison} for a detailed comparison).\\

\noindent First we discuss the hypotheses on the linear operator $H^{\beta}$. Fix $\beta_0>0$ and let $\bar\beta$ be a completely monotone function such that 
\begin{equation}\label{eq:conditionbarbeta}
\int_0^1\xi\bar\beta(\xi)\mathd\xi<+\infty,\quad \bar\beta(\xi) \sim \bar c_{\bar\beta}|\xi|^{-1-2s_2}\text{ as }|\xi|\rightarrow +\infty,\quad 0<s_2\leq  \frac{1}{2}.
\end{equation}
For any $\gamma_1>0$, define
\[
H^{\beta,\gamma_1}g
=-\beta_0\Delta g
+\gamma_1\int_{\R}\bigl(g(\cdot)-g(\xi)\bigr)
\bar\beta(\cdot-\xi)\mathd\xi.
\]
Since $\bar\beta$ is completely monotone and satisfies the conditions \eqref{eq:conditionbarbeta}, $H^{\beta,\gamma_1}$ satisfies
condition \ref{H0}. We stress that the first condition in \eqref{eq:conditionbarbeta} allows $\bar\beta$ to be unbounded near the origin, as for the explicit family of operators considered at the end of Section \ref{section:OScomparison}.\\

\noindent We now consider the conditions on the potential $V$, or better on $V_{\tot}$ (see Theorem \ref{theorem:BB}). First we introduce the following definition:

\begin{definition}\label{definition:doublewell}
We say that a function $f:\R\rightarrow\R_+$ is a (symmetric and
superquadratic) double-well-like potential if $f\in C^2(\R)$, $f$ is even,
there is $k>0$ such that
\[
f^{-1}(\{0\})=\{-k,k\},\qquad f'(x)=0\text{ for }x>0\text{ if and only if }x=k,
\qquad f''(\pm k)>0,
\]
and there are $\varepsilon,C,C'>0$ such that
\[
f(x)\geq  C|x|^{2+\varepsilon}-C'.
\]
We also assume, as required by \eqref{H2:pol}-\eqref{H2:diss} and \ref{H5}, that $f''$ is bounded from below and that there are $q>1$ and $C''>0$ such that
\[
|f'(x)|\leq C''(1+|x|^q),\qquad |f''(x)|\leq C''(1+|x|^{q-1}),\qquad xf'(x)\geq\frac{|x|^{q+1}}{C''}-C'',\qquad x\in\R.
\]
\end{definition}

\noindent A typical example satisfying Definition \ref{definition:doublewell} is
\[
f(x)=(k^2-x^2)^2.
\]
For any $\gamma_2>0$, define the total potential
\begin{equation}\label{eq:doublewell1}
V_{\tot}^{\gamma_2}(y)=\gamma_2 f(y).
\end{equation}
By Definition \ref{definition:doublewell}, for every $\gamma_2>0$ we can
write
\[
V_{\tot}^{\gamma_2}(y)=\frac{m^2}{2}|y|^2+V^{\gamma_2}(y),
\]
where $V^{\gamma_2}$ satisfies \eqref{H2:pol}-\eqref{H2:diss} and \ref{H5}, and is bounded
from below.

\begin{theorem}\label{theorem:nonuniqueness}
Consider $\beta_0,\bar\beta,H^{\beta,\gamma_1}$ as above, where $\bar\beta$ satisfies the conditions in \eqref{eq:conditionbarbeta}. Consider a function $f$ satisfying
Definition \ref{definition:doublewell}, and let $V_{\tot}^{\gamma_2}$ be
defined by \eqref{eq:doublewell1}. Then there are $\Gamma_1>0$ and a function
$\Gamma_2:\R_+\rightarrow\R_+$ such that, for every
$\gamma_1>\Gamma_1$ and $\gamma_2>\Gamma_2(\gamma_1)$, there are at least two
different measures $\mu_+$ and $\mu_-$ solving the DLR equations
\eqref{eq:DLRequations}. They satisfy
\[
\int \left(\int_{i}^{i+1} x(\xi)\mathd \xi\right)\mu_+(\mathd x)>0,
\qquad
\int  \left(\int_{i}^{i+1} x(\xi)\mathd \xi\right)\mu_-(\mathd x)<0,
\]
for every $i\in \N$.\\
Moreover, they satisfy the moment assumptions of Theorems
\ref{theorem:main} and \ref{theorem:infinitesimalsymmetri}; hence they are
reversible measures, and thus invariant measures of
\eqref{eq:LangevinSPDE1}.
\end{theorem}
\begin{proof}
The proof, which is based on Theorem 6.1 of \cite{OsadaSpohn1999}, is given in Subsection \ref{section:proofnonuniqueness}.
\end{proof}

The proof of Theorem \ref{theorem:nonuniqueness} is based on two main ingredients: the fact that the solution $\mu_+$ of the DLR equations with operator $H^{\beta,\gamma_1}$ and potential $V^{\gamma_2}$ dominates (in a sense detailed in Subsection \ref{section:GKSWells}) a long range Ising model which is known to have a phase transition (see, e.g., \cite{Dyson1969,FrohlichSpencer1982}); and the fact that we can always choose the constant of this Ising model in the multiple Gibbs measures regime thanks to the presence of two identical minima in the potential $V^{\gamma_2}_{\tot}$ (more concretely, this second result is based on a large-deviation principle for Gibbs measures in one dimension on an interval of increasing size of \cite{OttoWeberWestdickenberg2014}, see also \cite{BertiniButtaDiGesu2025}). Both ingredients rely on the representation, as in \cite{OsadaSpohn1999}, of the conditional measures $\mu^{\mathcal I,y}_{\free}$ as Gibbs measures relative to the Brownian bridge (see Theorem \ref{theorem:BB} and Lemma \ref{lemma:extendeduniqueness1}), which holds thanks to the first condition in \eqref{eq:conditionbarbeta}. When $\bar\beta$ is unbounded near the origin, the interaction between the points of the same block of the Ising approximation cannot be removed, and the large-deviation result of \cite{OttoWeberWestdickenberg2014} is transferred to the law of the blocks by an entropy argument (see Lemma \ref{lemma:convergencenugamma}).

\subsection{On GKS and Wells inequalities}\label{section:GKSWells}

We first introduce the principal notion of domination between probability measures we are going to use. 

\begin{definition}\label{definition:conecomparison}
Let $\mathcal X$ be a metric space and let
$\mathcal C\subset C(\mathcal X)$ be a family of continuous functions. We say
that a probability measure $\mu$ on $\mathcal X$ is dominated by a probability
measure $\nu$ with respect to $\mathcal C$ if
\[
\int_{\mathcal X}F(x)\mu(\mathd x)
\leq  
\int_{\mathcal X}F(x)\nu(\mathd x),
\qquad F\in\mathcal C,
\]
and we write $\mu\leq  _{\mathcal C}\nu$.
\end{definition}

\noindent In order to discuss meaningful results on the domination of measures, according to Definition \ref{definition:conecomparison}, we briefly recall the Griffiths-Kelly-Sherman
inequalities (see, e.g., Chapter VIII of \cite{Simon1974} or Chapter 4 of
\cite{GlimmJaffe1987}) and Wells' inequality (see, e.g.,
\cite{Wells1977,BricmontLebowitzPfister1981}), in the terminology of
\cite{MadridSimonWells2022} and the abstract formulation of
\cite{Ginibre1970}.\\

\noindent Let $\mathcal G(\mathcal C)$ be the smallest subset of $C(\mathcal X)$ that
contains $\mathcal C$ and the constant function $1$, and is closed under
products and positive linear combinations. We call
$\mathcal G(\mathcal C)$ the Ginibre cone generated by $\mathcal C$. We denote by $\mathcal{F}(\mathcal{C}) \subset \mathcal{G}(\mathcal{C})$ the subset of the Ginibre cone generated by the positive linear combinations of the constant function $1$, elements of $\mathcal{C}$, and products of elements of $\mathcal{C}$. We call $\mathcal{F}(\mathcal{C})$  the set of \emph{ferromagnetic potentials}.
Let $H \in \mathcal{F}(\mathcal{C})$  be a ferromagnetic potential. If $0<\int\exp(H)\mathd\mu<+\infty$, we write
\[
\mu^{\exp(H)}(\mathd x)
:=
\frac{\exp(H(x))}{\int_{\mathcal X}\exp(H(y))\mu(\mathd y)}\mu(\mathd x),
\]
and call it the tilted measure associated with $H$. Assume that $\mathcal C\subset L^p(\mu)\cap L^p(\nu)$ for every
$1\leq   p<+\infty$ and that $\exp(H) \in L^{p_0}(\mu) \cap L^{p_0}(\nu)$, for some $1< p_0 \leq   +\infty$. We say that $\mu$ is Ising dominated by $\nu$ with
respect to $\mathcal C$, and write
\[
\mu\leq  _{\Ising,\mathcal C}\nu,
\]
if, for every $H \in \mathcal{F}(\mathcal{C})$ such that the two tilted
measures $\mu^{\exp(H)},\nu^{\exp(H)}$ are well defined,
\[
\mu^{\exp(H)}
\leq  _{\mathcal G(\mathcal C)}
\nu^{\exp(H)}.
\]

\begin{remark}
Usually one has $\mathcal{X}=\R^n$ (or the Cartesian product of some compact subset of $\R$ or $\R^k$), and the cone $\mathcal{C}$ is generated by the standard coordinates $(y_1,...,y_n) \in \R^n$, 
\[\mathcal{C}=\{y_1,...,y_n\}.\]
The Ginibre cone $\mathcal{G}(\mathcal{C})$ is thus given by the expressions of the form 
\[\sum_h \sum_{i_1,...,i_h} c_{i_1,...,i_h} y_{i_1}\cdots y_{i_h} \]
where $c_{i_1,...,i_h} \geq  0$. The set $\mathcal{F}(\mathcal{C})$ is made of functions corresponding to the standard notion of ferromagnetic potential, namely expressions of the form 
\[\sum B_i y_i+ \sum J_{i,j} y_i y_j,\quad B_i, J_{i,j}\geq  0.\]

\end{remark}

\begin{definition}
We say that $(\mu,\mathcal X,\mathcal C)$ is a Ginibre system if
$\mathcal C\subset L^p(\mu)$ for every $1\leq   p<+\infty$ and, for every
$F_1,\ldots,F_n\in\mathcal C$,
\begin{align}
&\int_{\mathcal X}\prod_{i=1}^nF_i(x)\mu(\mathd x)\geq 0,
\label{eq:Ginibre}\\
&\int_{\mathcal X\times\mathcal X}
\prod_{i=1}^n\bigl(F_i(x)\pm F_i(y)\bigr)
\mu(\mathd x)\mu(\mathd y)\geq 0,
\label{eq:Ginibre2}
\end{align}
for all choices of the signs in \eqref{eq:Ginibre2}.
\end{definition}

\begin{remark}\label{remark:Ginibreproperties}
Ginibre systems were introduced and first studied in \cite{Ginibre1970}. \\
They are closed under tensor products, namely, if $(\mu_1,\mathcal{X}_1,\mathcal{C}_1)$ and $(\mu_2,\mathcal{X}_2,\mathcal{C}_2)$ are two Ginibre systems, then also $(\mu_1\otimes \mu_2,\mathcal{X}_1\times \mathcal{X}_2,\mathcal{C}_1 \otimes \mathcal{C}_2)$ is a Ginibre system.\\
Furthermore, Ginibre systems are closed with respect to the tilting of a measure: if $(\mu,\mathcal{X},\mathcal{C})$ is a Ginibre system and $H \in \mathcal{F}(\mathcal{C})$ (or even $H\in \mathcal{G}(\mathcal{C})$), satisfying the above hypotheses on $\exp(H)$, then $(\mu^{\exp(H)},\mathcal{X},\mathcal{C})$ is a new Ginibre system.
\end{remark}

\noindent For Ginibre systems the first and second Griffiths-Kelly-Sherman
inequalities (GKS) hold. We report here only GKSII. 

\begin{theorem}
Suppose that $(\mu,\mathcal X,\mathcal C)$ is a Ginibre system. Then, for
every $F,G\in\mathcal G(\mathcal C)$,
\[
\int_{\mathcal X}F(x)G(x)\mu(\mathd x)
\geq 
\left(\int_{\mathcal X}F(x)\mu(\mathd x)\right)
\left(\int_{\mathcal X}G(x)\mu(\mathd x)\right).
\]
\end{theorem}
\begin{proof}
The theorem is proved in 
\cite{Ginibre1970}; see also \cite{MadridSimonWells2022}.
\end{proof}

\begin{corollary}\label{corollary:monotonicityGKS}
Fix a Ginibre system $(\mu,\mathcal X,\mathcal C)$. Let $H$ be such that
$H\in\mathcal F(\mathcal C)$ and assume that, for any $t \in [0,1]$, $\exp(t H),H e^{tH}\in L^1(\mu)$.
Then
\[
\mu\leq  _{\mathcal G(\mathcal C)}\mu^{\exp(H)}.
\]
\end{corollary}
\begin{proof}
For $t\in[0,1]$, set $\mu_t:=\mu^{\exp(tH)}$. The expansion of
$\exp(tH(x)+tH(y))$ in \eqref{eq:Ginibre2} shows that
$(\mu_t,\mathcal X,\mathcal C)$ is again a Ginibre system. Thus, for every
$F\in\mathcal G(\mathcal C)$,
\[
\frac{\mathd}{\mathd t}\int_{\mathcal X}F\mathd\mu_t
=\operatorname{Cov}_{\mu_t}(F,H)\geq 0.
\]
Integration between $0$ and $1$ proves the claim.
\end{proof}

\begin{corollary}\label{corollary:limitmonotonicityGKS}
Fix a Ginibre system $(\mu,\mathcal X,\mathcal C)$. Let $H$ be such that
$\mu^{\exp(H)}$ is well defined and suppose that there are $H_n$ such that
$H_n\in\mathcal F(\mathcal C)$ and Corollary
\ref{corollary:monotonicityGKS} applies to every $H_n$. If, for every
$F\in\mathcal{G}(\mathcal{C})$,
\[
\exp(H_n)F\longrightarrow\exp(H)F
\quad\text{in }L^1(\mu),
\]
then
\[
\mu\leq  _{\mathcal G(\mathcal C)}\mu^{\exp(H)}.
\]
\end{corollary}
\begin{proof}
For every $F\in\mathcal G(\mathcal C)$, Corollary
\ref{corollary:monotonicityGKS} gives
\[
\int F\mathd\mu
\leq  
\frac{\int F\exp(H_n)\mathd\mu}{\int\exp(H_n)\mathd\mu}.
\]
The assumed $L^1$ convergence allows us to pass to the limit.
\end{proof}

\noindent Fix $\mathcal X_0=C([a,b],\R)$, with the topology of uniform convergence, and
let
\begin{equation}\label{eq:coneC0}
\mathcal C_0:=
\left\{x\longmapsto\int_a^b x(\xi)\kappa(\mathd\xi):
\kappa\text{ is a positive bounded Borel measure on }[a,b]\right\}.
\end{equation}
Let $\bar\beta,\bar\beta':[a,b]^2\rightarrow[0,+\infty)$ and $\gamma,\gamma':\R \rightarrow \R_+$ be measurable and
bounded, and set
\begin{align*}
B(x):=\int_{[a,b]^2}\bar\beta(t,s)x(s)x(t)\mathd s\mathd t,& \quad B'(x):=\int_{[a,b]^2}\bar\beta'(t,s)x(s)x(t)\mathd s\mathd t.\\
G(x):=\int_{[a,b]}\gamma(t)x(t)\mathd t&\quad G'(x):=\int_{[a,b]}\gamma'(t)x(t)\mathd t
\end{align*}

\begin{theorem}\label{theorem:kernelmonotonicityGKS}
Suppose that $(\mu,\mathcal X_0,\mathcal C_0)$ is a Ginibre system and that
there are $\varepsilon,\varepsilon'>0$ such that
\[
\int_{\mathcal X_0}
\exp\left(\varepsilon'\int_a^b|x(s)|^{2+\varepsilon}\mathd s\right)
\mu(\mathd x)<+\infty.
\]
If $\bar\beta\geq \bar\beta'$ and $\gamma \geq  \gamma'$ almost everywhere, then
\[
\frac{1}{\mathcal Z_{B'+G'}}\exp(B'+G')\mu
\leq  _{\mathcal G(\mathcal C_0)}
\frac{1}{\mathcal Z_{B+G}}\exp(B+G)\mu.
\]
\end{theorem}
\begin{proof}
We write explicitly the proof in the case $\gamma=\gamma'=0$. The general case is a simple extension. Set $D:=B-B'$ and $d:=\bar\beta-\bar\beta'$. Thus $d\geq 0$ almost
everywhere. Choose nonnegative rectangular step functions
$\bar\beta_n'$ and $d_n$, uniformly bounded in $n$, such that
\[
\bar\beta_n'\longrightarrow\bar\beta',
\qquad
d_n\longrightarrow d
\quad\text{almost everywhere on }[a,b]^2.
\]
We denote by $B_n'$ and $D_n$ the corresponding quadratic forms.
Since $\bar\beta'_n$ and $d_n$ are piecewise constant, $D_n,B'_n\in\mathcal F(\mathcal C_0)$. \\

\noindent The boundedness of the
kernels and H\"older's inequality give, for every $\delta>0$, a constant
$C_\delta$ such that
\begin{align*}
&|B'(x)|+|D(x)|+
\sup_n\bigl(|B_n'(x)|+|D_n(x)|\bigr)\\
&\qquad\leq  
C\left(\int_a^b|x(s)|\mathd s\right)^2
\leq   C_\delta+
\delta\int_a^b|x(s)|^{2+\varepsilon}\mathd s.
\end{align*}
By choosing $\delta$ smaller than the constant
$\varepsilon'$ in the assumption, the right-hand side provides an
integrable exponential majorant. Hence dominated convergence gives, for
every $F\in\mathcal G(\mathcal C_0)\cup\{1\}$,
\begin{align}
\exp(B_n')F&\longrightarrow\exp(B')F
&&\text{in }L^1(\mu),
\label{eq:Bnapproximation}\\
\exp(D_n)F&\longrightarrow\exp(D)F
&&\text{in }L^1(\mu^{\exp(B')}).
\label{eq:Dnapproximation}
\end{align}

\noindent For every $n$, by Remark \ref{remark:Ginibreproperties}, the tilted measure
$\nu_n:=\mu^{\exp(B_n')}$
is again a Ginibre system. By
\eqref{eq:Bnapproximation}, $\nu_n$ converges to
$\nu:=\mu^{\exp(B')}$ in all the moments appearing in those inequalities. Passing to the limit
shows that $(\nu,\mathcal X_0,\mathcal C_0)$ is a Ginibre system.\\
We may now apply Corollary \ref{corollary:limitmonotonicityGKS} to $\nu$,
with $H=D$ and $H_n=D_n$. Indeed,
$D_n\in\mathcal F(\mathcal C_0)$, and
\eqref{eq:Dnapproximation} is exactly the required convergence. We obtain
\[
\frac{1}{\mathcal Z_{B'}}\exp(B')\mu
=\nu
\leq  _{\mathcal G(\mathcal C_0)}
\nu^{\exp(D)}
=\frac{1}{\mathcal Z_B}\exp(B)\mu.
\]
\end{proof}

\noindent We now introduce a stronger notion than Ising dominated measures due to Wells \cite{Wells1977}.

\begin{definition}
For two probability measures $\mu$ and $\nu$ on $\mathcal X$, we write
\[
\mu\leq  _{\mathrm W,\mathcal C}\nu
\]
if, for every $F_1,\ldots,F_n\in\mathcal C$ and every choice of signs,
\begin{equation}\label{eq:Wellsdomination}
\int_{\mathcal X\times\mathcal X}
\prod_{i=1}^n\bigl(F_i(x)\pm F_i(y)\bigr)
\nu(\mathd x)\mu(\mathd y)\geq 0.
\end{equation}
Thus $\mu\leq  _{\mathrm W,\mathcal C}\nu$ means that $\nu$ Wells
dominates $\mu$.
\end{definition}

\begin{remark}\label{remark:Wellsproduct}
As for Ginibre systems (see Remark \ref{remark:Ginibreproperties}), Wells domination is stable under products: if $\mu_1,\nu_1$ and $\mu_2,\nu_2$ are pairs of probability measures on $\mathcal{X}_1$ and $\mathcal{X}_2$ respectively such that $\mu_1 \leq  _{\mathrm{W},\mathcal C_1} \nu_1$, and $\mu_2 \leq  _{\mathrm{W},\mathcal C_2} \nu_2$ then $\mu_1 \otimes \mu_2 \leq  _{\mathrm{W},\mathcal{C}_1 \otimes \mathcal{C}_2} \nu_1 \otimes \nu_2$. 
\end{remark}

\noindent For Wells dominated measures the following result holds.

\begin{theorem}[Wells]\label{theorem:Wells}
If $\mu\leq  _{\mathrm W,\mathcal C}\nu$, then
\[
\mu\leq  _{\Ising,\mathcal C}\nu.
\]
\end{theorem}
\begin{proof}
See \cite{BricmontLebowitzPfister1981,MadridSimonWells2022} for the complete
abstract formulation and its integrability extension.
\end{proof}

\noindent Wells condition is commonly used to compare generic symmetric measures with Bernoulli (symmetric) measures. We report here an easy sufficient criterion.

\begin{lemma}\label{lemma:Wellscriterion}
Let $\rho$ be an even atomless probability measure on $\R$ with all moments
finite, and fix $c \geq  \sqrt{2}$. If $a\geq0$ satisfies
\begin{equation}\label{eq:Wellschoicea}
\rho([0,a])\leq   (c^2-1) \rho([c a,+\infty)),
\end{equation}
then
\[
\frac12(\delta_{-a}+\delta_a)
\leq  _{\mathrm W,\{y\}}\rho.
\]
\end{lemma}
\begin{proof}
We verify \eqref{eq:Wellsdomination} directly. If there are $n$ plus signs
and $m$ minus signs, its left-hand side is
\begin{equation}\label{eq:Wellsdirectmoment}
\frac12\int_\R
\left[(x+a)^n(x-a)^m+(x-a)^n(x+a)^m\right]\rho(\mathd x).
\end{equation}
When $n+m$ is odd, the integrand is odd and the integral vanishes. Suppose
that $n+m$ is even. Put $r=\min\{n,m\}$ and $d=|n-m|$; then $d$ is even,
and the integrand in \eqref{eq:Wellsdirectmoment} equals
\[
(x^2-a^2)^r\,
\frac{(x+a)^d+(x-a)^d}{2}.
\]
The second factor is a polynomial in $x^2$ and $a^2$ with nonnegative
coefficients. It is therefore enough to show that
\begin{equation}\label{eq:Wellsauxmoment}
\int_\R x^{2\ell}(x^2-a^2)^r\rho(\mathd x)\geq 0,
\qquad \ell,r\in\mathbb N\cup\{0\}.
\end{equation}
This is immediate if $r$ is even. If $r$ is odd, the integrand is negative
only on $[-a,a]$, where  
\begin{equation}\label{eq:Wellsa}
\int_{-a}^{a}x^{2\ell}(x^2-a^2)^r\rho(\mathd x) \geq  -2a^{2\ell+2r} \rho([0,a])
\end{equation}
On the other hand, on $[c a,+\infty)$ and its reflection one has
\begin{equation}\label{eq:Wellsa2}
\int_{[c a,+\infty)} x^{2\ell}(x^2-a^2)^r\rho(\mathd x)\geq  c^{2 \ell} (c^2-1)^r a^{2\ell+2r}\rho([ca,+\infty)) \geq  (c^2-1) a^{2\ell+2r}\rho([ca,+\infty))
\end{equation}
and similarly for $\int_{(-\infty,-ca]}x^{2\ell}(x^2-a^2)^r\rho(\mathd x)$ by reflection symmetry. Putting together \eqref{eq:Wellsa} and \eqref{eq:Wellsa2} we obtain \eqref{eq:Wellsauxmoment}.
\end{proof}

\begin{remark}\label{corollary:deltaWells}
We can apply Lemma \ref{lemma:Wellscriterion} to the case where $a=0$, obtaining
\[
\delta_0\leq  _{\mathrm W,\{y\}}\rho,
\qquad
\delta_0\leq  _{\Ising,\{y\}}\rho.
\]
\end{remark}

\noindent We now return to $\mathcal X_0=C([a,c],\R)$ and to the cone $\mathcal C_0$ given by \eqref{eq:coneC0} with $[a,c]$ in place of $[a,b]$. Let $\mathcal I:=(a,c)$ and let $\mathcal W_{\beta_0}^{\mathcal I,(b,b)}$ be the
Brownian bridge defined in \eqref{eq:betaBB}, with $b\geq 0$. If
$a<t_1<\cdots<t_j<c$, let
\[
\mathcal W_{\beta_0}^{\mathcal I,(b,b),(t_1,\ldots,t_j)}
\]
denote the concatenation of the bridges from $b$ to $0$, from $0$ to $0$,
and from $0$ to $b$ obtained by pinning the path to zero at
$t_1,\ldots,t_j$.

\begin{lemma}\label{lemma:BBcomparison}
Let $V:\R\rightarrow\R$ be even and bounded from below. Then
\[
\frac{1}{\mathcal Z}
\exp\left(-\int_{\mathcal I}V(x(\xi))\mathd\xi\right)
\mathcal W_{\beta_0}^{\mathcal I,(b,b),(t_1,\ldots,t_j)}
\leq  _{\Ising,\mathcal C_0}
\frac{1}{\mathcal Z'}
\exp\left(-\int_{\mathcal I}V(x(\xi))\mathd\xi\right)
\mathcal W_{\beta_0}^{\mathcal I,(b,b)}.
\]
\end{lemma}
\begin{proof}
The argument is based on Trotter approximation, as in Appendix B of
\cite{GuerraRosenSimon1975II}. Let
\[
a=s_0<s_1<\cdots<s_n=c,
\qquad h_i:=s_{i+1}-s_i,
\]
be a sequence of partitions whose mesh tends to zero and which contains all
the points $t_1,\ldots,t_j$. When checking the inequality for fixed
cylinder functions, we also choose the partitions to contain all the times
on which those functions depend. Write $x_i=x(s_i)$ and
$x_0=x_n=b$. Up to constants independent of
$x_1,\ldots,x_{n-1}$, the discretized weighted Brownian density is
\begin{align*}
&\exp\left(
-\frac{\beta_0}{2}\sum_{i=0}^{n-1}
\frac{(x_{i+1}-x_i)^2}{h_i}
-\sum_{i=1}^{n-1}\ell_{i,n}V(x_i)
\right)\prod_{i=1}^{n-1}\mathd x_i\\
&\qquad=\exp\bigl(K_n(x)\bigr)
\prod_{i=1}^{n-1}\rho_{i,n}(\mathd x_i),
\end{align*}
where $\ell_{i,n}>0$ are the quadrature weights in the Trotter
approximation,
\[
\rho_{i,n}(\mathd z)
=\frac1{z_{i,n}}
\exp\left[-\frac{\beta_0}{2}
\left(\frac1{h_{i-1}}+\frac1{h_i}\right)z^2
-\ell_{i,n}V(z)\right]\mathd z,\]
\[K_n(x)=\beta_0\sum_{i=1}^{n-2}\frac{x_ix_{i+1}}{h_i}
+\frac{\beta_0b}{h_0}x_1
+\frac{\beta_0b}{h_{n-1}}x_{n-1}.
\]
Since
$V$ is even and bounded from below, every $\rho_{i,n}$ is a well-defined
even probability measure with all moments finite. If
\[
\mathcal C_n:=\{x\mapsto x_i:1\leq   i\leq   n-1\},
\]
then $K_n\in\mathcal F(\mathcal C_n)$: the nearest-neighbor coefficients
are positive and the two boundary-field coefficients are nonnegative because
$b\geq 0$.

\noindent For the bridge pinned at $t_1,\ldots,t_j$, the one-site measure
$\rho_{i,n}$ is replaced by $\delta_0$ whenever
$s_i\in\{t_1,\ldots,t_j\}$. At all other sites the two one-site measures
coincide. Remark \ref{corollary:deltaWells} gives
\[
\delta_0\leq  _{\mathrm W,\{x_i\}}\rho_{i,n}
\]
at a pinned site. Since Wells domination is stable under products by Remark \ref{remark:Wellsproduct}, after tilting
both product measures by the same ferromagnetic factor $\exp(K_n)$,
Theorem \ref{theorem:Wells} gives
\[
\mu_{n,V}^{(b,b),(t_1,\ldots,t_j)}
\leq  _{\Ising,\mathcal C_n}
\mu_{n,V}^{(b,b)},
\]
where the two sides are the normalized discretized measures.\\
The Trotter formula gives the convergence of the discretized measures to the two weighted Brownian bridges when integrated with respect to cylinder functions, namely polynomially bounded continuous functions depending on $x(a_1),...,x(a_h)$ for a finite number of times $a_1,...,a_h$. This convergence implies that 
\[
\frac{1}{\mathcal Z}
\exp\left(-\int_{\mathcal I}V(x(\xi))\mathd\xi\right)
\mathcal W_{\beta_0}^{\mathcal I,(b,b),(t_1,\ldots,t_j)}
\leq  _{\Ising,\mathcal C_{0,\operatorname{cyl}}}
\frac{1}{\mathcal Z'}
\exp\left(-\int_{\mathcal I}V(x(\xi))\mathd\xi\right)
\mathcal W_{\beta_0}^{\mathcal I,(b,b)},
\] 
where 
\[\mathcal C_{0,\operatorname{cyl}}:=\left\{\sum_{i=1}^h c_i x(a_i)\text{ such that }h\in \mathbb{N}, c_i>0,a_i \in [a,c]\right\}\]
Since $\mathcal C_{0,\operatorname{cyl}}$ is dense in $\mathcal C_{0}$ with respect to pointwise convergence dominated by the $L^p$ norms, we obtain the statement.
\end{proof}

\noindent The measure on the left-hand side of Lemma \ref{lemma:BBcomparison} is the
regular conditional law of the weighted bridge on the right-hand side given
\[
x(t_1)=\cdots=x(t_j)=0.
\]
Thus the content of the lemma is that conditioning a ferromagnetic bridge to
vanish at finitely many interior points decreases it in the Ginibre-Ising
order.

\begin{remark}\label{remark:BBinequality}
The same approximation gives, for $b\geq 0$,
\[
\frac{1}{\mathcal Z_0}
\exp\left(-\int_{\mathcal I}V(x(\xi))\mathd\xi\right)
\mathcal W_{\beta_0}^{\mathcal I,(0,0)}
\leq  _{\Ising,\mathcal C_0}
\frac{1}{\mathcal Z_b}
\exp\left(-\int_{\mathcal I}V(x(\xi))\mathd\xi\right)
\mathcal W_{\beta_0}^{\mathcal I,(b,b)}.
\]
Indeed, changing the endpoints from $0$ to $b$ adds a nonnegative boundary
field in every discretization, namely
\[
\frac{\beta_0b}{h_0}x_1
+\frac{\beta_0b}{h_{n-1}}x_{n-1}.
\]
\end{remark}

\begin{remark} Lemma \ref{lemma:BBcomparison} may be viewed as a one-dimensional version of the monotonicity and decoupling properties of Dirichlet boundary conditions discussed in Sections IV.3 and IV.4 of \cite {GuerraRosenSimon1975II} and Section 7.8 of \cite {GlimmJaffe1987}. Indeed, pinning the path to zero at an interior time introduces a Dirichlet boundary, and the Trotter formula is similar to the lattice approximation used in the cited works. A related comparison holds for the two-dimensional $P(\phi)_2$ measure (where the Wick renormalization is with respect to the full-space covariance operator, and $P$ is an even polynomial with positive leading coefficient) by imposing Dirichlet conditions on an interior curve (this is the so-called half-Dirichlet construction proposed in \cite{GuerraRosenSimon1975II,GuerraRosenSimon1976}). 
\end{remark}

\subsection{A large-deviation-type result}\label{section:largedeviation}

For $\gamma_2>0$, consider the probability measure
\begin{equation}\label{eq:kappagamma}
\kappa^{\gamma_2}(\mathd x)
:=\frac{1}{\mathcal Z^{\gamma_2}}
\exp\left(
-\gamma_2\int_0^1f(x(\xi))\mathd\xi
\right)
\mathcal W_{\beta_0}^{[0,1],(0,0)}(\mathd x),
\end{equation}
and define $g_{k,\max}$ to be the distance between a point $y \in \R$ and the set $\{\pm k\}$, namely
\[
g_{k,\max}(y):=\min\{|y-k|,|y+k|\}.
\]

\begin{lemma}\label{lemma:kappaconvergence}
Using the previous notation, for every $c>0$, we have
\begin{equation}\label{eq:kappaweakconvergence2}
\limsup_{\gamma_2\rightarrow+\infty}\frac{1}{\sqrt{\gamma_2}}\log
\kappa^{\gamma_2}
\left(
g_{k,\max}\left(\int_0^1x(\xi)\mathd\xi\right)\geq  c
\right)<0.
\end{equation}
Consequently, the laws $\Law_{\kappa^{\gamma_2}}\left(\int_0^1x(\xi)\mathd\xi\right)$ converge weakly to $\frac12\delta_{-k}+\frac12\delta_k$ as $\gamma_2\rightarrow+\infty$.

\end{lemma}

\noindent In order to prove Lemma \ref{lemma:kappaconvergence}, it is convenient to reformulate the previous convergence \eqref{eq:kappaweakconvergence2}  as a large-deviation problem for a tilted Gaussian measure. Put
\[
\varepsilon:=\gamma_2^{-1/2},
\qquad T_\varepsilon:=\varepsilon^{-1},
\qquad t:=\xi/\varepsilon.
\]
Under the map $u(t)=x(\varepsilon t)$, the measure
\eqref{eq:kappagamma} becomes
\begin{equation}\label{eq:kapparescaled}
\bar\kappa^{\varepsilon,[0,T_\varepsilon],(0,0)}(\mathd u)
=\frac{1}{\bar{\mathcal Z}^{\varepsilon}}
\exp\left(-\frac1\varepsilon
\int_0^{T_\varepsilon}f(u(t))\mathd t\right)
\mathcal W_{\beta_0/\varepsilon}^{[0,T_\varepsilon],(0,0)}(\mathd u).
\end{equation}
For an interval $[a,b]$ and endpoints $y_1,y_2$, we use the analogous
notation
$\bar\kappa^{\varepsilon,[a,b],(y_1,y_2)}$.\\
We say that $u$ has a $\delta$-crossing (in the interval $[a,b]$) if there are $a\leq   t_1<t_2 \leq   b$ such that
either
\[
u(t_1)=-k+\delta,\qquad u(t_2)=k-\delta,
\]
or the same relations hold with the signs reversed, and the path between
$t_1$ and $t_2$ remains between these two values.
We shall use the energy
\begin{equation}\label{eq:energyfunctional}
E^{[a,b]}(u):=
\begin{cases}
\displaystyle\frac{\beta_0}{2}\int_a^b|\partial_tu|^2\mathd t
+\int_a^bf(u(t))\mathd t,&u\in H^1([a,b]),\\[2mm]
+\infty,&\text{otherwise},
\end{cases}
\end{equation}
and write $E_{\inf}^{[a,b],(y_1,y_2)}$ for its infimum among paths with
endpoints $y_1,y_2$.\\

\noindent If we fix $a<b$ and $y_1,y_2 \in \R$, the behavior of the tilted Gaussian measure $\bar\kappa^{\varepsilon,[a,b],(y_1,y_2)}$ is very well studied and is related to the energy functional $E^{[a,b]}$ in \eqref{eq:energyfunctional}. For the measure $\bar\kappa^{\varepsilon,[0,T_{\varepsilon}],(0,0)}$ the problem is more complex since the length of the interval $[0,T_{\varepsilon}]$ of definition of $u$ is diverging as $\frac{1}{\varepsilon}$. For solving this problem we use the results of \cite{OttoWeberWestdickenberg2014}. In particular we have:

\begin{theorem}\label{theorem:OWW}
Fix $\delta\in(0,k)$. There are
$M_1,C,r_\delta>0$ such that, for every $M\geq  M_1$ and all sufficiently
small $\varepsilon$,
\begin{align}
&\sup_{t\in[0,T_\varepsilon]}
\bar\kappa^{\varepsilon,[0,T_\varepsilon],(0,0)}
(|u(t)|\geq  M)
\leq  \exp\left(-\frac{M}{C\varepsilon}\right),
\label{eq:OWWtail}\\
&\bar\kappa^{\varepsilon,[0,T_\varepsilon],(0,0)}
(\|u\|_{L^\infty}\leq   k-\delta)
\leq   CT_\varepsilon\exp\left(-\frac{r_\delta}{\varepsilon}\right),
\label{eq:OWWreach}\\
&\bar\kappa^{\varepsilon,[0,T_\varepsilon],(0,0)}
(u\text{ has a }\delta\text{-crossing})
\leq   CT_\varepsilon\exp\left(-\frac{r_\delta}{\varepsilon}\right).
\label{eq:OWWcrossing}
\end{align}
Moreover, if the interval length ranges in a fixed compact subset of
$(0,+\infty)$ and $|y_1|,|y_2|\leq   M$, then, for every measurable set
$\mathcal B$, every $\eta,\lambda>0$, and all sufficiently small
$\varepsilon$,
\begin{equation}\label{eq:OWWlocalLDP}
\bar\kappa^{\varepsilon,[a,b],(y_1,y_2)}(\mathcal B)
\leq  
\exp\left(-\frac{
\Delta E^{[a,b],(y_1,y_2)}(\mathcal B_\eta)-\lambda
}{\varepsilon}\right),
\end{equation}
where $\mathcal B_\eta$ is the $\eta$-neighborhood in the uniform norm and
\[
\Delta E^{[a,b],(y_1,y_2)}(\mathcal A)
:=\inf_{u\in\mathcal A}E^{[a,b]}(u)
-E_{\inf}^{[a,b],(y_1,y_2)}.
\]
\end{theorem}
\begin{proof}
The estimates follow from Proposition 3.4, Lemma 4.1, and the proof of
Theorem 1.5 of \cite{OttoWeberWestdickenberg2014}; see also Remark 1.8
therein concerning boundary conditions which do not force a transition.
\end{proof}

\noindent We need also the following elementary result.

\begin{lemma}\label{lemma:infeenergyuniform}
For every $M>0$, there is $K_{f,M}<+\infty$ such that
\[
\sup_{\substack{|y_1|,|y_2|\leq   M\\b-a\geq 1}}
E_{\inf}^{[a,b],(y_1,y_2)}\leq   K_{f,M}.
\]
\end{lemma}
\begin{proof}
By translation in time it is enough to work on $[0,T]$, $T=b-a\geq 1$. The result follows by considering the explicit path
\[
v(t)=
\begin{cases}
y_1+2t(k-y_1),&0\leq   t\leq  \frac12,\\
k,&\frac12\leq   t\leq   T-\frac12,\\
k+2(t-T+\frac12)(y_2-k),&T-\frac12\leq   t\leq   T.
\end{cases}
\]
using the fact that $f(\pm k)=0$.
\end{proof}

\begin{lemma}\label{lemma:convergenceintegralmodified}
Let
\[
D_{[a,b]}(u):=\frac1{b-a}\int_a^bg_{k,\max}(u(t))\mathd t.
\]
For every $c>0$,
\[
\limsup_{\varepsilon\rightarrow0}\,\varepsilon\log
\bar\kappa^{\varepsilon,[0,T_\varepsilon],(0,0)}
(D_{[0,T_\varepsilon]}(u)>c)<0.
\]
\end{lemma}
\begin{proof}
The hypotheses on $f$ imply that there is $c_f>0$ such that
\[
f(y)\geq  c_fg_{k,\max}(y)^2,
\qquad y\in\R.
\]
Indeed, near the wells this follows from $f''(\pm k)>0$; away from the
wells on a compact set it follows from strict positivity; and at infinity
it follows from the superquadratic lower bound. Thus the convex increasing
function $h(r):=c_fr^2$ satisfies $h(0)=0$, $h(r)>0$ for $r>0$, and
\[
f(y)\geq  h(g_{k,\max}(y)),\qquad y\in\R.
\]
Jensen's inequality therefore gives
\begin{equation}\label{eq:energybound}
E^{[a,b]}(u)\geq (b-a)h(D_{[a,b]}(u)).
\end{equation}

\noindent Fix $M\geq  M_1$ as in Theorem
\ref{theorem:OWW}. Choose $\ell\geq 1$ and
$\lambda>0$ so that
\[
\varrho:=\ell h(c/2)-K_{f,M}-\lambda>0.
\]
Let $N_\varepsilon:=\lfloor T_\varepsilon/\ell\rfloor$ and partition
$[0,T_\varepsilon]$ into $N_\varepsilon$ intervals $I_j$ of the common
length $\ell_\varepsilon:=T_\varepsilon/N_\varepsilon$. For small
$\varepsilon$, $\ell_\varepsilon\in[\ell,2\ell]$. If
$D_{[0,T_\varepsilon]}(u)>c$, then $D_{I_j}(u)>c$ for at least one $j$.
Since $g_{k,\max}$ is $1$-Lipschitz, the uniform $c/2$-neighborhood of
$\{D_{I_j}>c\}$ is contained in $\{D_{I_j}>c/2\}$. Hence
\eqref{eq:energybound}, Lemma \ref{lemma:infeenergyuniform}, and
\eqref{eq:OWWlocalLDP} show that, conditionally on endpoint values $y_1,y_2 \in \R$ such that $|y_1|,|y_2| \leq   M$,
\[
\bar\kappa^{\varepsilon,I_j,(y_1,y_2)}(D_{I_j}>c)
\leq   e^{-\varrho/\varepsilon}.
\]
Let $t_j$ denote the endpoints of the
partition. The Markov property gives, for any event $A_j$ depending only on
$u|_{I_j}$,
\[
\bar\kappa^{\varepsilon,[0,T_\varepsilon],(0,0)}(A_j)
=\int
\bar\kappa^{\varepsilon,I_j,(y_1,y_2)}(A_j)
\,\Law_{\bar{\kappa}^{\varepsilon,[0,T_{\varepsilon}],(0,0)}}(u(t_j),u(t_{j+1}))(\mathd y_1\mathd y_2).
\]
We split this integral according to
$\{|y_1|\vee|y_2|\leq   M\}$ and its complement. On the first part we
use the preceding local large-deviation estimate; on the second part we use
\eqref{eq:OWWtail}. Summing over all the intervals $I_j$ of the partition we obtain
\[
\bar\kappa^{\varepsilon,[0,T_\varepsilon],(0,0)}
(D_{[0,T_\varepsilon]}>c)
\leq  
2N_\varepsilon e^{-M/(C\varepsilon)}
+N_\varepsilon e^{-\varrho/\varepsilon}.
\]
Since $N_\varepsilon\leq   T_\varepsilon=\varepsilon^{-1}$, the
upper limit in the statement is not greater than $-\min\{M/C,\varrho\}<0$.
\end{proof}

\begin{proof}[Proof of Lemma \ref{lemma:kappaconvergence}]
Fix $\delta\in(0,k)$. If
\[
\mathcal T_\delta:=
\{\|u\|_{L^\infty}\leq   k-\delta\}
\cup\{u\text{ has a }\delta\text{-crossing}\},
\]
then, by Theorem \ref{theorem:OWW}, we have
\[
\bar\kappa^{\varepsilon,[0,T_\varepsilon],(0,0)}
(\mathcal T_\delta)
\leq   2C\varepsilon^{-1}e^{-r_\delta/\varepsilon}.
\]
On $\mathcal T_\delta^c$ the path can belong to one of two sets $H^+$ and
$H^-$: in
$H^+$ the path $u$ reaches $k-\delta$ and never reaches $-k+\delta$, and in $H^-$ the
opposite occurs. When $u \in H^{\pm}\subset \mathcal{T}^c_{\delta}$, there is $C_{k,\delta}<+\infty$ (not depending on $u$) such that
\[
|u(t)\mp k|\leq   C_{k,\delta}g_{k,\max}(u(t)),
\qquad u\in H^\pm.
\]
Consequently,
\[
\left|\frac1{T_\varepsilon}\int_0^{T_\varepsilon}u(t)\mathd t\mp k\right|
\leq   C_{k,\delta}D_{[0,T_\varepsilon]}(u),
\qquad u\in H^\pm.
\]
Thus, by the previous discussion we have
\begin{multline*}
\bar\kappa^{\varepsilon,[0,T_{\varepsilon}],(0,0)}\left(g_{k,\max}\left(\frac{1}{T_{\varepsilon}}\int_0^{T_{\varepsilon}}u(t)\mathd t\right)\geq c\right)\\
\leq   \bar\kappa^{\varepsilon,[0,T_{\varepsilon}],(0,0)}(\mathcal{T}_{\delta})+\bar\kappa^{\varepsilon,[0,T_{\varepsilon}],(0,0)}\left(D_{[0,T_{\varepsilon}]}(u)>\frac{c}{2C_{k,\delta}}\right),
\end{multline*}
so that, by Lemma \ref{lemma:convergenceintegralmodified}, the left-hand side is not greater than $e^{-r/\varepsilon}$ for some $r>0$ and every $\varepsilon$ small enough. Since $\varepsilon=\gamma_2^{-1/2}$ and
\[
\bar\kappa^{\varepsilon,[0,T_{\varepsilon}],(0,0)}\left(g_{k,\max}\left(\frac1{T_\varepsilon}\int_0^{T_\varepsilon}u(t)\mathd t\right)\geq c\right)
=\kappa^{\gamma_2}\left(g_{k,\max}\left(\int_0^1x(\xi)\mathd\xi\right)\geq c\right),
\]
equation \eqref{eq:kappaweakconvergence2} is proved.\\
The limit \eqref{eq:kappaweakconvergence2} implies that the laws $\Law_{\kappa^{\gamma_2}}\left(\int_0^1x(\xi)\mathd\xi\right)$ are tight and that all their limit points are supported on the set $\{\pm k\}$. Since the measures $\kappa^{\gamma_2}$ are symmetric with respect to the reflection $x \mapsto-x$ of the paths, we obtain the weak convergence to $\frac12\delta_{-k}+\frac12\delta_k$.
\end{proof}

\noindent In the proof of Theorem \ref{theorem:nonuniqueness} we also need the following estimates. We use the relative entropy $\operatorname{Ent}(P\Vert R)$ of a probability measure $P$ with respect to a probability measure $R$ on the same measurable space, introduced in the proof of Proposition \ref{prop:equivalenceDLR}, and we recall that, by the Jensen inequality, if $\operatorname{Ent}(P\Vert R)<+\infty$, then for every measurable function $\Phi$ bounded from below or from above we have
\begin{equation}\label{eq:Gibbsvariational}
\int\Phi\,\mathd P\leq\operatorname{Ent}(P\Vert R)+\log\int e^{\Phi}\mathd R.
\end{equation}

\begin{lemma}\label{lemma:kappaentropy}
There are $\Gamma_0,C_\kappa>0$ such that, for every $\gamma_2\geq\Gamma_0$,
\begin{equation}\label{eq:kappaentropy}
\sup_{\xi\in[0,1]}\int x(\xi)^2\kappa^{\gamma_2}(\mathd x)\leq C_\kappa,
\qquad
\operatorname{Ent}\left(\kappa^{\gamma_2}\big\Vert\mathcal W_{\beta_0}^{[0,1],(0,0)}\right)\leq C_\kappa\sqrt{\gamma_2}.
\end{equation}
\end{lemma}
\begin{proof}
Let $\varepsilon:=\gamma_2^{-1/2}$. Since $\kappa^{\gamma_2}(|x(\xi)|\geq M)=\bar\kappa^{\varepsilon,[0,T_\varepsilon],(0,0)}(|u(\xi/\varepsilon)|\geq M)$, by \eqref{eq:OWWtail}, for every $\xi\in[0,1]$ and every $\varepsilon$ small enough we have
\[
\int x(\xi)^2\kappa^{\gamma_2}(\mathd x)=\int_0^{+\infty}2M\kappa^{\gamma_2}(|x(\xi)|\geq M)\mathd M\leq M_1^2+\int_{M_1}^{+\infty}2Me^{-\frac{M}{C\varepsilon}}\mathd M\leq M_1^2+1.
\]
In order to prove the second estimate, we set $\mathcal W:=\mathcal W_{\beta_0}^{[0,1],(0,0)}$ and we denote by $A$ the operator introduced in the proof of Lemma \ref{lemma:extendeduniqueness1} with $\mathcal I=(0,1)$, so that $\mathcal W$ is the centered Gaussian measure with covariance $A^{-1}$ and the eigenvalues of $A$ are $\beta_0\pi^2j^2$, $j\in\N$. Since $f\geq0$, by \eqref{eq:kappagamma} and \eqref{eq:Gibbsvariational} with $R=\mathcal W$ and $\Phi(x)=-\gamma_2\int_0^1f(x(\xi))\mathd\xi$, for every probability measure $P$ we have
\begin{align*}
\operatorname{Ent}(\kappa^{\gamma_2}\Vert\mathcal W)&=-\gamma_2\int\int_0^1f(x(\xi))\mathd\xi\,\kappa^{\gamma_2}(\mathd x)-\log\mathcal Z^{\gamma_2}\\
&\leq-\log\mathcal Z^{\gamma_2}\leq\gamma_2\int\int_0^1f(x(\xi))\mathd\xi\,P(\mathd x)+\operatorname{Ent}(P\Vert\mathcal W).
\end{align*}
We choose as $P$ the law of $\psi+X$, where $\psi(\xi):=k\min\{1,\xi/\varepsilon,(1-\xi)/\varepsilon\}$ and $X$ is a centered Gaussian random variable with covariance $(A+\gamma_2)^{-1}$. By Lemma \ref{lemma:Gaussianproof}, applied with $B=\gamma_2 I$, and by the Cameron-Martin formula (see \cite[Corollary 2.4.3]{Bogachev98}), using that the law of $X$ is symmetric, we get
\begin{align*}
\operatorname{Ent}(P\Vert\mathcal W)&=\frac12\log\det\nolimits_F(I+\gamma_2A^{-1})-\frac{\gamma_2}{2}\mathbb E\left[\norm{X}^2_{L^2((0,1))}\right]+\frac{\beta_0}{2}\int_0^1|\psi'(\xi)|^2\mathd\xi\\
&\leq\frac12\sum_{j\in\N}\log\left(1+\frac{\gamma_2}{\beta_0\pi^2j^2}\right)+\frac{\beta_0k^2}{\varepsilon},
\end{align*}
and, since $\sum_{j\in\N}\log(1+a^2j^{-2})\leq\int_0^{+\infty}\log(1+a^2t^{-2})\mathd t=\pi a$ for $a>0$, we obtain $\operatorname{Ent}(P\Vert\mathcal W)\leq(\frac{1}{2\sqrt{\beta_0}}+\beta_0k^2)\sqrt{\gamma_2}$. Moreover, by the maximum principle, the Green function of $A+\gamma_2$ is bounded by the one of $-\beta_0\Delta+\gamma_2$ on $\R$, namely $(4\beta_0\gamma_2)^{-1/2}e^{-\sqrt{\gamma_2/\beta_0}|\xi-\xi'|}$, so that $\mathbb E[|X(\xi)|^2]\leq(4\beta_0\gamma_2)^{-1/2}$ for every $\xi\in[0,1]$. By Definition \ref{definition:doublewell}, since $f(k)=f'(k)=0$, there is $C_f>0$ such that $f(k+z)\leq C_f(|z|^2+|z|^{q+1})$ and $f(z)\leq C_f(1+|z|^{q+1})$ for every $z\in\R$. Therefore, since $X(\xi)$ is a centered Gaussian random variable and $\psi=k$ on $[\varepsilon,1-\varepsilon]$, there is $C'_f>0$, independent of $\gamma_2\geq1$, such that
\begin{align*}
\gamma_2\int\int_0^1f(x(\xi))\mathd\xi\,P(\mathd x)&\leq\gamma_2\int_\varepsilon^{1-\varepsilon}\mathbb E\left[f(k+X(\xi))\right]\mathd\xi+\gamma_2\int_{[0,\varepsilon]\cup[1-\varepsilon,1]}\mathbb E\left[f(\psi(\xi)+X(\xi))\right]\mathd\xi\\
&\leq C'_f\sqrt{\gamma_2},
\end{align*}
and the second estimate in \eqref{eq:kappaentropy} follows.
\end{proof}

\subsection{Proof of Theorem \ref{theorem:nonuniqueness}}\label{section:proofnonuniqueness}

For $N\in\N$, set
\[
\mathcal I_N=(-N,N),\qquad
\Lambda_N=\{-N,\ldots,N-1\},\qquad I_i=[i,i+1].
\]
On $\mathcal X^N=C([-N,N],\R)$, consider the cones
\begin{align*}
\mathcal C_c^N
&:=\left\{x\longmapsto\int_{-N}^Nx(\xi)\kappa(\mathd\xi):
\kappa\text{ is a positive bounded Borel measure}\right\},\\
\mathcal C_{\dis}^N
&:=\left\{x\longmapsto\sum_{i\in\Lambda_N}c_i
\int_i^{i+1}x(\xi)\mathd\xi:c_i\geq 0\right\}.
\end{align*}
Let $\bar{\mathcal X}_{\dis}^N=\R^{\Lambda_N}$ and let
$\bar{\mathcal C}_{\dis}^N$ be the cone generated by the coordinate maps.
Define
\[
\Pi_N:\mathcal X^N\longrightarrow\bar{\mathcal X}_{\dis}^N,
\qquad
(\Pi_Nx)_i=\int_i^{i+1}x(\xi)\mathd\xi.
\]
The map $\Pi_N$ is continuous and surjective, and
$\mathcal C_{\dis}^N$ is precisely the pullback under $\Pi_N$ of
$\bar{\mathcal C}_{\dis}^N$. Consequently,
\begin{equation}\label{eq:pushforwardcone}
\mu\leq  _{\mathcal C_{\dis}^N}\mu'
\quad\Longrightarrow\quad
\Pi_{N*}\mu\leq  _{\bar{\mathcal C}_{\dis}^N}\Pi_{N*}\mu'.
\end{equation}
For $n\neq0$, set
\[
J_n:=\inf_{\substack{t\in I_0\\s\in I_n}}\bar\beta(t-s),
\qquad J_0:=0.
\]
By translation invariance, this is the same coefficient for any two blocks at
distance $n$. For $i\in\mathbb Z$ and $\xi\in(i,i+1)$, set
\[
c_{\bar\beta}(\xi):=\frac12\int_{\R\setminus I_i}\bar\beta(\xi-\xi')\mathd\xi'.
\]
The function $c_{\bar\beta}$ is $1$-periodic and, by \eqref{eq:conditionbarbeta}, $\int_0^1c_{\bar\beta}(\xi)\mathd\xi=\int_0^{+\infty}\bar\beta(r)\min\{r,1\}\mathd r<+\infty$. Define
\begin{equation}\label{eq:nugamma1}
\nu^{\gamma_1,\gamma_2}
:=\Law_{\tilde\nu^{\gamma_1,\gamma_2}}
\left(\int_0^1x(\xi)\mathd\xi\right),
\end{equation}
where
\begin{multline}\label{eq:nugamma2}
\tilde\nu^{\gamma_1,\gamma_2}(\mathd x)
=\frac{1}{\tilde Z^{\gamma_1,\gamma_2}}
\exp\left(
-\int_0^1\big(\gamma_2f(x(\xi))+\gamma_1c_{\bar\beta}(\xi)x(\xi)^2\big)\mathd\xi
\right.\\
\left.
-\frac{\gamma_1}{4}\int_{[0,1]^2}\bar\beta(\xi-\xi')|x(\xi)-x(\xi')|^2\mathd\xi\mathd\xi'
\right)
\mathcal W_{\beta_0}^{[0,1],(0,0)}(\mathd x).
\end{multline}
The exponent in \eqref{eq:nugamma2} is finite for $\mathcal W_{\beta_0}^{[0,1],(0,0)}$-almost every $x$, since $\int|x(\xi)-x(\xi')|^2\mathcal W_{\beta_0}^{[0,1],(0,0)}(\mathd x)\leq\beta_0^{-1}|\xi-\xi'|$ for $\xi,\xi'\in[0,1]$. For $b\geq 0$, let
\begin{multline}\label{eq:discretemeasure}
\mu_{\dis}^{N,b,\gamma_1,\gamma_2}(\mathd M)
:=\frac1{Z_{\dis}^{N,b,\gamma_1,\gamma_2}}
\exp\left(
\frac{\gamma_1}{2}
\sum_{\substack{i,j\in\Lambda_N\\i\neq j}}J_{j-i}M_iM_j
\right.\\
\left.
+\gamma_1b\sum_{i\in\Lambda_N}\sum_{j\notin\Lambda_N}J_{j-i}M_i
\right)
\prod_{i\in\Lambda_N}\nu^{\gamma_1,\gamma_2}(\mathd M_i).
\end{multline}
Finally, let $\mu^{N,b}$ be the finite-volume DLR measure on $\mathcal I_N$
with constant exterior condition $b$, namely the measure $\mu^{\vartheta,y_0}$ given by \eqref{eq:existencemuthetay0} with $\vartheta=N$ and $y_0\equiv b$, for the parameters
$\gamma_1,\gamma_2$ fixed below.

\begin{lemma}\label{lemma:CNmu}
With the preceding notation,
\[
\mu_{\dis}^{N,b,\gamma_1,\gamma_2}
\leq  _{\bar{\mathcal C}_{\dis}^N}
\Pi_{N*}\mu^{N,b}.
\]
\end{lemma}
\begin{proof}
Since $\bar\beta$ may be unbounded near the origin, we first replace it by the bounded kernel $\bar\beta_L:=\min\{\bar\beta,L\}$, where $L\geq\bar\beta(1)$. We construct explicitly a continuous-path measure
$\tilde\mu^{N,b,\gamma_1,\gamma_2}_L$ such that
\begin{equation}\label{eq:twopropertiesmutilde}
\tilde\mu^{N,b,\gamma_1,\gamma_2}_L
\leq  _{\mathcal C_{\dis}^N}\mu^{N,b}_L,
\qquad
\Pi_{N*}\tilde\mu^{N,b,\gamma_1,\gamma_2}_L
=\mu_{\dis,L}^{N,b,\gamma_1,\gamma_2},
\end{equation}
where the measures $\mu^{N,b}_L$ and $\mu_{\dis,L}^{N,b,\gamma_1,\gamma_2}$, defined below, converge to $\mu^{N,b}$ and $\mu_{\dis}^{N,b,\gamma_1,\gamma_2}$ as $L\rightarrow+\infty$. The conclusion then follows from \eqref{eq:pushforwardcone} by letting $L\rightarrow+\infty$.\\
By Lemma \ref{lemma:extendeduniqueness1}, after combining the mass term and
$V^{\gamma_2}$ into $V_{\tot}^{\gamma_2}=\gamma_2f$, the finite-volume DLR
measure $\mu^{N,b}$ has density proportional to
\begin{multline}\label{eq:nonuniqueness1}
\exp\left(
-\gamma_2\int_{\mathcal I_N}f(x(\xi))\mathd\xi
-\frac{\gamma_1}{4}\int_{\mathcal I_N^2}
\bar\beta(\xi-\xi')|x(\xi)-x(\xi')|^2\mathd\xi\mathd\xi'
\right.\\
\left.
-\frac{\gamma_1}{2}\int_{\mathcal I_N\times\mathcal I_N^c}
\bar\beta(\xi-\xi')|x(\xi)-b|^2\mathd\xi\mathd\xi'
\right)
\mathcal W_{\beta_0}^{\mathcal I_N,(b,b)}(\mathd x),
\end{multline}
and we denote by $\mu^{N,b}_L$ the probability measure with density proportional to \eqref{eq:nonuniqueness1} with $\bar\beta$ replaced by $\bar\beta_L$. We expand the square in this expression. The total coefficient of the term involving the integral of $x(\xi)^2$ is given by the expression
\begin{align*}
&-\frac{\gamma_1}{2}
\int_{\mathcal I_N}x(\xi)^2
\left(\int_{\mathcal I_N}\bar\beta_L(\xi-\xi')\mathd\xi'
+\int_{\mathcal I_N^c}\bar\beta_L(\xi-\xi')\mathd\xi'\right)\mathd\xi=-\frac{\gamma_1}{2}
\norm{\bar\beta_L}_{L^1(\R)}
\int_{\mathcal I_N}x(\xi)^2\mathd\xi.
\end{align*}
Thus, after discarding constants depending only on $b$, $N$ and $L$, the exponent
in \eqref{eq:nonuniqueness1} with $\bar\beta_L$ in place of $\bar\beta$ is
\begin{multline}\label{eq:expandedcontinuuminteraction}
-\gamma_2\int_{\mathcal I_N}f(x(\xi))\mathd\xi
-\frac{\gamma_1}{2}\norm{\bar\beta_L}_{L^1(\R)}\int_{\mathcal I_N}x(\xi)^2\mathd\xi
+\frac{\gamma_1}{2}\int_{\mathcal I_N^2}
\bar\beta_L(\xi-\xi')x(\xi)x(\xi')\mathd\xi\mathd\xi'\\
+\gamma_1b\int_{\mathcal I_N\times\mathcal I_N^c}
\bar\beta_L(\xi-\xi')x(\xi)\mathd\xi\mathd\xi'.
\end{multline}

\noindent Since the local part $\gamma_2f(y)+\frac{\gamma_1}{2}\norm{\bar\beta_L}_{L^1(\R)}y^2$ of \eqref{eq:expandedcontinuuminteraction} is symmetric and the remaining terms in \eqref{eq:expandedcontinuuminteraction} can be approximated by objects in $\mathcal{C}_{c}^N$, by Lemma \ref{lemma:BBcomparison} and Remark \ref{remark:BBinequality}, we have
\begin{equation}\label{eq:muhatcomparison}
\widehat\mu^{N,b,\gamma_1,\gamma_2}_L
\leq  _{\Ising,\mathcal C_c^N}\mu^{N,b}_L,
\end{equation}
where $\widehat\mu^{N,b,\gamma_1,\gamma_2}_L$ is the probability measure whose density with respect to $\mathcal{W}^{\mathcal I_N,(0,0),(-N+1,\ldots,N-1)}_{\beta_0}$ is proportional to the exponential of \eqref{eq:expandedcontinuuminteraction}.\\

\noindent Define the blockwise kernel $\rho_L$ on $\R^2$ by
\[
\rho_L(t,s)=J_{j-i},
\qquad (t,s)\in I_i\times I_j,\quad i\neq j,
\]
and set $\rho_L(t,s)=\bar\beta_L(t-s)$ on $I_i\times I_i$. Since $\bar\beta$ is nonincreasing on $(0,+\infty)$ and $J_{j-i}\leq\bar\beta(1)\leq L$, we have
\[
0\leq  \rho_L(t,s)\leq  \bar\beta_L(t-s).
\]
Differently from \cite[Section 6]{OsadaSpohn1999}, we do not remove the interaction inside the blocks, since $\norm{\bar\beta_L}_{L^1(\R)}\rightarrow+\infty$ as $L\rightarrow+\infty$ when $\bar\beta$ is not integrable. By Remark \ref{remark:Ginibreproperties}, the measure with density proportional to $\exp(-\gamma_2\int_{\mathcal I_N}f(x(\xi))\mathd\xi-\frac{\gamma_1}{2}\norm{\bar\beta_L}_{L^1(\R)}\int_{\mathcal I_N}x(\xi)^2\mathd\xi)$ with respect to $\mathcal W_{\beta_0}^{\mathcal I_N,(0,0),(-N+1,\ldots,N-1)}$ is a Ginibre system, and we define
\begin{multline}\label{eq:nonuniqueness2}
\tilde\mu^{N,b,\gamma_1,\gamma_2}_L(\mathd x)
:=\frac1{\widetilde{\mathcal Z}_L}
\exp\left(
-\gamma_2\int_{\mathcal I_N}f(x(\xi))\mathd\xi
-\frac{\gamma_1}{2}\norm{\bar\beta_L}_{L^1(\R)}\int_{\mathcal I_N}x(\xi)^2\mathd\xi
\right.\\
\left.
+\frac{\gamma_1}{2}\int_{\mathcal I_N^2}
\rho_L(\xi,\xi')x(\xi)x(\xi')\mathd\xi\mathd\xi'
+\gamma_1b\int_{\mathcal I_N\times\mathcal I_N^c}
\rho_L(\xi,\xi')x(\xi)\mathd\xi\mathd\xi'
\right)\\
\times
\mathcal W_{\beta_0}^{\mathcal I_N,(0,0),(-N+1,\ldots,N-1)}(\mathd x),
\end{multline}
where $\rho_L(\xi,\xi'):=J_{j-i}$ also for $(\xi,\xi')\in I_i\times I_j$ with $i\in\Lambda_N$ and $j\notin\Lambda_N$.
Thus we can apply Theorem \ref{theorem:kernelmonotonicityGKS}
obtaining 
\[
\tilde\mu^{N,b,\gamma_1,\gamma_2}_L
\leq  _{\mathcal G(\mathcal C_c^N)}
\widehat\mu^{N,b,\gamma_1,\gamma_2}_L.
\]
Together with \eqref{eq:muhatcomparison} and
$\mathcal C_{\dis}^N\subset\mathcal C_c^N$, this proves the first relation
in \eqref{eq:twopropertiesmutilde}.\\

\noindent It remains to compute the pushforward. Since the bridge in
\eqref{eq:nonuniqueness2} is zero at every integer, its restrictions to the
unit blocks are independent $(0,0)$ bridges. Moreover, for $\xi\in I_i$ we have $\frac12\norm{\bar\beta_L}_{L^1(\R)}-\frac12\int_{I_i}\bar\beta_L(\xi-\xi')\mathd\xi'=c_{\bar\beta_L}(\xi)$, where $c_{\bar\beta_L}$ is defined as $c_{\bar\beta}$, with $\bar\beta_L$ in place of $\bar\beta$. Hence the terms in the exponent of \eqref{eq:nonuniqueness2} which involve only one block are
\[
-\sum_{i\in\Lambda_N}\left(\int_{I_i}\big(\gamma_2f(x(\xi))+\gamma_1c_{\bar\beta_L}(\xi)x(\xi)^2\big)\mathd\xi+\frac{\gamma_1}{4}\int_{I_i^2}\bar\beta_L(\xi-\xi')|x(\xi)-x(\xi')|^2\mathd\xi\mathd\xi'\right),
\]
which, since $c_{\bar\beta_L}$ is $1$-periodic, is a sum of identical one-block terms. Therefore
\begin{multline*}
\frac1{\hat{\mathcal Z}_L}
\exp\left(
-\sum_{i\in\Lambda_N}\left(\int_{I_i}\big(\gamma_2f(x(\xi))+\gamma_1c_{\bar\beta_L}(\xi)x(\xi)^2\big)\mathd\xi+\frac{\gamma_1}{4}\int_{I_i^2}\bar\beta_L(\xi-\xi')|x(\xi)-x(\xi')|^2\mathd\xi\mathd\xi'\right)
\right)\\
\times\mathcal W_{\beta_0}^{\mathcal I_N,(0,0),(-N+1,\ldots,N-1)}(\mathd x)=\bigotimes_{i\in\Lambda_N}
\tilde\nu^{\gamma_1,\gamma_2}_L(\mathd(x|_{I_i})),
\end{multline*}
where $\tilde\nu^{\gamma_1,\gamma_2}_L$ is defined by \eqref{eq:nugamma2} with $\bar\beta_L$ in place of $\bar\beta$.
If $M_i=\int_{I_i}x(\xi)\mathd\xi$, the remaining terms equal
\begin{align*}
\frac{\gamma_1}{2}\sum_{i\in\Lambda_N}\int_{I_i\times(\mathcal I_N\setminus I_i)}
\rho_L(\xi,\xi')x(\xi)x(\xi')\mathd\xi\mathd\xi'
&=\frac{\gamma_1}{2}
\sum_{\substack{i,j\in\Lambda_N\\i\neq j}}
J_{j-i}M_iM_j,\\
\gamma_1b\int_{\mathcal I_N\times\mathcal I_N^c}
\rho_L(\xi,\xi')x(\xi)\mathd\xi\mathd\xi'
&=\gamma_1b\sum_{i\in\Lambda_N}
\sum_{j\notin\Lambda_N}J_{j-i}M_i.
\end{align*}
Thus these terms depend only on $\Pi_Nx$, and pushing
\eqref{eq:nonuniqueness2} forward under $\Pi_N$ gives the measure $\mu_{\dis,L}^{N,b,\gamma_1,\gamma_2}$ defined by \eqref{eq:discretemeasure} with $\nu^{\gamma_1,\gamma_2}$ replaced by $\nu^{\gamma_1,\gamma_2}_L:=\Law_{\tilde\nu^{\gamma_1,\gamma_2}_L}(\int_0^1x(\xi)\mathd\xi)$. This proves the second relation in \eqref{eq:twopropertiesmutilde}.\\

\noindent Finally, we let $L\rightarrow+\infty$. Since $\bar\beta_L\uparrow\bar\beta$ and $c_{\bar\beta_L}\uparrow c_{\bar\beta}$, by the monotone convergence theorem the exponents of the densities of $\mu^{N,b}_L$ and $\tilde\nu_L^{\gamma_1,\gamma_2}$ with respect to $\mathcal W_{\beta_0}^{\mathcal I_N,(b,b)}$ and $\mathcal W_{\beta_0}^{[0,1],(0,0)}$, respectively, converge to the ones of $\mu^{N,b}$ and $\tilde\nu^{\gamma_1,\gamma_2}$, which are finite almost surely by Lemma \ref{lemma:extendeduniqueness1} and by the discussion after \eqref{eq:nugamma2}. Hence the normalization constants converge to positive limits and, since the densities are bounded by $\exp(-\gamma_2\int f(x(\xi))\mathd\xi)$, up to the normalization constants, and $f$ is superquadratic, by the dominated convergence theorem we obtain
\[
\lim_{L\rightarrow+\infty}\int F\circ \Pi_N\,\mathd\mu^{N,b}_L=\int F\circ \Pi_N\,\mathd\mu^{N,b},
\qquad
\lim_{L\rightarrow+\infty}\int F\,\mathd\mu_{\dis,L}^{N,b,\gamma_1,\gamma_2}=\int F\,\mathd\mu_{\dis}^{N,b,\gamma_1,\gamma_2},
\]
for every $F\in\bar{\mathcal C}_{\dis}^N$,
where, in the second limit, we also used that the exponent in \eqref{eq:discretemeasure} is a quadratic function of $M$ which does not depend on $L$. Therefore, by \eqref{eq:twopropertiesmutilde} and \eqref{eq:pushforwardcone}, letting $L\rightarrow+\infty$ we obtain the statement.
\end{proof}

\begin{lemma}\label{lemma:convergencenugamma}
For every fixed $\gamma_1>0$, $\nu^{\gamma_1,\gamma_2}$ converges weakly to $\frac12\delta_{-k}+\frac12\delta_k$ as $\gamma_2\rightarrow+\infty$.
\end{lemma}
\begin{proof}
Let $\mathcal W:=\mathcal W_{\beta_0}^{[0,1],(0,0)}$ and, for $\eta\in(0,1)$, let
\begin{align*}
G_1(x)&:=\gamma_1\int_0^1c_{\bar\beta}(\xi)x(\xi)^2\mathd\xi+\frac{\gamma_1}{4}\int_{\{(\xi,\xi')\in[0,1]^2:|\xi-\xi'|\geq\eta\}}\bar\beta(\xi-\xi')|x(\xi)-x(\xi')|^2\mathd\xi\mathd\xi',\\
G_2(x)&:=\frac{\gamma_1}{4}\int_{\{(\xi,\xi')\in[0,1]^2:|\xi-\xi'|<\eta\}}\bar\beta(\xi-\xi')|x(\xi)-x(\xi')|^2\mathd\xi\mathd\xi',
\end{align*}
so that, by \eqref{eq:kappagamma} and \eqref{eq:nugamma2}, $\tilde\nu^{\gamma_1,\gamma_2}$ has density $e^{-G_1-G_2}/\int e^{-G_1-G_2}\mathd\kappa^{\gamma_2}$ with respect to $\kappa^{\gamma_2}$. Since $G_1,G_2\geq0$, by the Jensen inequality we get
\[
\operatorname{Ent}(\tilde\nu^{\gamma_1,\gamma_2}\Vert\kappa^{\gamma_2})=-\int(G_1+G_2)\mathd\tilde\nu^{\gamma_1,\gamma_2}-\log\int e^{-G_1-G_2}\mathd\kappa^{\gamma_2}\leq\int(G_1+G_2)\mathd\kappa^{\gamma_2}.
\]
Since $\bar\beta$ is nonincreasing on $(0,+\infty)$, we have $G_1(x)\leq\gamma_1\int_0^1(c_{\bar\beta}(\xi)+\bar\beta(\eta))x(\xi)^2\mathd\xi$, so that, by Lemma \ref{lemma:kappaentropy}, $\int G_1\mathd\kappa^{\gamma_2}\leq\gamma_1C_\kappa(\int_0^1c_{\bar\beta}(\xi)\mathd\xi+\bar\beta(\eta))$ for $\gamma_2\geq\Gamma_0$. As far as $G_2$ is concerned, since $\int|x(\xi)-x(\xi')|^2\mathcal W(\mathd x)\leq\beta_0^{-1}|\xi-\xi'|$, we have
\[
\int G_2\mathd\mathcal W\leq\frac{\gamma_1}{2\beta_0}\int_0^\eta r\bar\beta(r)\mathd r=:\frac1{4\theta_\eta}.
\]
As in the proof of Lemma \ref{lemma:Gaussianproof}, applied with $A$ as in the proof of Lemma \ref{lemma:kappaentropy} and $B$ the nonnegative self-adjoint operator associated with the closure of the quadratic form $G_2$ on $H^1_0((0,1))$, there are $\lambda_j\geq0$, $j\in\N$, such that $\sum_j\lambda_j=\int G_2\mathd\mathcal W$ and $\int e^{\theta_\eta G_2}\mathd\mathcal W=\prod_j(1-2\theta_\eta\lambda_j)^{-1/2}$. Since $2\theta_\eta\lambda_j\leq\frac12$ and $-\log(1-z)\leq2z$ for $z\in[0,\frac12]$, we get $\log\int e^{\theta_\eta G_2}\mathd\mathcal W\leq2\theta_\eta\sum_j\lambda_j\leq\frac12$. Hence, by \eqref{eq:Gibbsvariational} with $P=\kappa^{\gamma_2}$, $R=\mathcal W$ and $\Phi=\theta_\eta G_2$, and by Lemma \ref{lemma:kappaentropy},
\[
\int G_2\mathd\kappa^{\gamma_2}\leq\frac{1}{\theta_\eta}\left(\operatorname{Ent}(\kappa^{\gamma_2}\Vert\mathcal W)+\frac12\right)\leq\frac{2\gamma_1}{\beta_0}\int_0^\eta r\bar\beta(r)\mathd r\left(C_\kappa\sqrt{\gamma_2}+\frac12\right),\qquad\gamma_2\geq\Gamma_0.
\]
Let now $\delta>0$. Since $\int_0^1r\bar\beta(r)\mathd r<+\infty$, we can choose $\eta$ such that $\frac{2\gamma_1C_\kappa}{\beta_0}\int_0^\eta r\bar\beta(r)\mathd r\leq\delta$, and we obtain a constant $C_\delta>0$, independent of $\gamma_2$, such that
\begin{equation}\label{eq:entropynukappa}
\operatorname{Ent}(\tilde\nu^{\gamma_1,\gamma_2}\Vert\kappa^{\gamma_2})\leq\delta\sqrt{\gamma_2}+C_\delta,\qquad\gamma_2\geq\Gamma_0.
\end{equation}
Let $c>0$ and let $E_c:=\{x:g_{k,\max}(\int_0^1x(\xi)\mathd\xi)\geq c\}$. By \eqref{eq:kappaweakconvergence2} there is $r_c>0$ such that $\kappa^{\gamma_2}(E_c)\leq e^{-r_c\sqrt{\gamma_2}}$ for every $\gamma_2$ large enough. Applying \eqref{eq:Gibbsvariational} with $P=\tilde\nu^{\gamma_1,\gamma_2}$, $R=\kappa^{\gamma_2}$ and $\Phi=t\mathbb I_{E_c}$, where $t:=\log(1+\kappa^{\gamma_2}(E_c)^{-1})\geq r_c\sqrt{\gamma_2}$, we get
\[
t\,\tilde\nu^{\gamma_1,\gamma_2}(E_c)\leq\operatorname{Ent}(\tilde\nu^{\gamma_1,\gamma_2}\Vert\kappa^{\gamma_2})+\log\left(1+(e^t-1)\kappa^{\gamma_2}(E_c)\right)=\operatorname{Ent}(\tilde\nu^{\gamma_1,\gamma_2}\Vert\kappa^{\gamma_2})+\log2,
\]
so that, by \eqref{eq:entropynukappa},
\[
\limsup_{\gamma_2\rightarrow+\infty}\nu^{\gamma_1,\gamma_2}\left(\{y\in\R:g_{k,\max}(y)\geq c\}\right)=\limsup_{\gamma_2\rightarrow+\infty}\tilde\nu^{\gamma_1,\gamma_2}(E_c)\leq\limsup_{\gamma_2\rightarrow+\infty}\frac{\delta\sqrt{\gamma_2}+C_\delta+\log2}{r_c\sqrt{\gamma_2}}=\frac{\delta}{r_c}.
\]
Since $\delta>0$ is arbitrary, the left-hand side is equal to $0$ for every $c>0$. Hence every limit point of $\{\nu^{\gamma_1,\gamma_2}\}$ is supported on $\{\pm k\}$, and, since $\nu^{\gamma_1,\gamma_2}$ is even, the statement follows.
\end{proof}

We can now introduce the long range Ising model which is dominated by $\mu^{N,b}$.

\begin{lemma}
We have
\begin{equation}\label{eq:disIsing}
a\int\sigma_i\,\Ising^{\Lambda_N,+}(\mathd\sigma)
\leq  
\int M_i\,\mu_{\dis}^{N,a,\gamma_1,\gamma_2}(\mathd M),
\qquad i\in\Lambda_N,
\end{equation}
where $\Ising^{\Lambda_N,+}$ is the measure of the Ising model having Hamiltonian
\begin{equation}\label{eq:Hamiltonian}
H_I^{\Lambda_N,+}(\sigma)
=-\frac{\gamma_1a^2}{2}
\sum_{\substack{i,j\in\Lambda_N\\i\neq j}}
J_{j-i}\sigma_i\sigma_j
-\gamma_1a^2\sum_{i\in\Lambda_N}
\sum_{j\notin\Lambda_N}J_{j-i}\sigma_i,
\qquad\sigma_i\in\{-1,1\},
\end{equation}
with $+1$ boundary conditions, and $a(\gamma_1,\gamma_2)>0$ is the only solution to the equation 
\begin{equation}\label{eq:defaWells}
F_{\gamma_1,\gamma_2}(a):=\nu^{\gamma_1,\gamma_2}([0,a])
-\nu^{\gamma_1,\gamma_2}([\sqrt2a,+\infty))=0.
\end{equation}
\end{lemma}
\begin{proof}
The measure $\nu^{\gamma_1,\gamma_2}$ is even. It is also atomless and gives
positive mass to every nonempty open interval. Indeed, under the Brownian
bridge the nonzero continuous linear functional
$x\mapsto\int_0^1x(\xi)\mathd\xi$ is a nondegenerate Gaussian random
variable, and the Gibbs density in \eqref{eq:nugamma2} is strictly positive.\\
From this we obtain that $F_{\gamma_1,\gamma_2}$  is continuous and strictly increasing. Moreover,
$F_{\gamma_1,\gamma_2}(0)=-1/2$ and
$F_{\gamma_1,\gamma_2}(a)\to1/2$ as $a\to+\infty$. This implies that equation \eqref{eq:defaWells} has a unique solution.  \\
Using the previous observation we can use Lemma \ref{lemma:Wellscriterion} obtaining
\begin{equation}\label{eq:singleWellscomparison}
\frac12\delta_{-a}+\frac12\delta_a
\leq  _{\mathrm W,\{y\}}
\nu^{\gamma_1,\gamma_2}.
\end{equation}
We use this inequality to show \eqref{eq:disIsing}. Take the product over
$i\in\Lambda_N$ of \eqref{eq:singleWellscomparison}, which is Wells dominated by the product $\otimes \nu^{\gamma_1,\gamma_2}$ due to Remark \ref{remark:Wellsproduct}. We then tilt both product measures by the same
ferromagnetic factor appearing in \eqref{eq:discretemeasure}  and apply
Theorem \ref{theorem:Wells}.\\
On the smaller product measure write
$M_i=a\sigma_i$, $\sigma_i\in\{-1,1\}$, and choose the exterior value
$b=a$. The resulting tilted measure is precisely the long-range Ising
measure with plus boundary conditions and Hamiltonian \eqref{eq:Hamiltonian}.
\end{proof}

\begin{proof}[Proof of Theorem \ref{theorem:nonuniqueness}]
We first prove that, given the conditions \eqref{eq:conditionbarbeta} on $\bar\beta$, the Ising model with Hamiltonian \eqref{eq:Hamiltonian} has a phase transition.\\
Observe that, for $t\in I_0$ and $s\in I_n$, one has $|t-s|=|n|+O(1)$. The asymptotic
assumption on $\bar\beta$ therefore implies
\[
J_n\sim\bar c_{\bar\beta}|n|^{-1-2s_2},
\]
and in particular $J_n\geq  C|n|^{-1-2s_2}$ for $n\neq0$, after
possibly decreasing $C>0$. Since $1<1+2s_2\leq  2$, the corresponding
Ising model has a phase transition; see
\cite{Dyson1969,FrohlichSpencer1982}. Thus there is $\alpha_0>0$, depending only on the sequence $(J_n)$, such
that, whenever $\gamma_1a^2>\alpha_0$,
\begin{equation}\label{eq:Isingpositive}
\lim_{N\rightarrow+\infty}
\int\sigma_i\,\Ising^{\Lambda_N,+}(\mathd\sigma)>0.
\end{equation}

\noindent Fix $0<a_0<k/\sqrt2$. By
Lemma \ref{lemma:convergencenugamma},
\[
\nu^{\gamma_1,\gamma_2}([0,a_0])\longrightarrow0,
\qquad
\nu^{\gamma_1,\gamma_2}([\sqrt2a_0,+\infty))\longrightarrow\frac12.
\]
Since the left-hand side minus the right-hand side of
\eqref{eq:defaWells} is $F_{\gamma_1,\gamma_2}(a)$, the last two limits say
that $F_{\gamma_1,\gamma_2}(a_0)<0$ for all sufficiently large $\gamma_2$.
The unique zero therefore satisfies $a(\gamma_1,\gamma_2)>a_0$. Since this
holds for every $a_0<k/\sqrt2$, we obtain
\[
\liminf_{\gamma_2\rightarrow+\infty}a(\gamma_1,\gamma_2)
\geq \frac{k}{\sqrt2}.
\]
Consequently, for every fixed $\gamma_1>0$, there is
$\Gamma_2(\gamma_1)$ such that
\[
a(\gamma_1,\gamma_2)\geq \frac{k}{2},
\qquad\gamma_2\geq \Gamma_2(\gamma_1).
\]
Choose, for instance,
\[
\Gamma_1:=\frac{8\alpha_0}{k^2}.
\]
Then $\gamma_1>\Gamma_1$ and
$\gamma_2\geq \Gamma_2(\gamma_1)$ imply
$\gamma_1a^2>2\alpha_0>\alpha_0$.\\

\noindent By \eqref{eq:Isingpositive}, \eqref{eq:disIsing} and Lemma \ref{lemma:CNmu}, applied to the coordinate function $M\mapsto M_i$, we have
\[
0<a\lim_{N\rightarrow+\infty}
\int\sigma_i\,\Ising^{\Lambda_N,+}(\mathd\sigma)\leq\liminf_{N\rightarrow+\infty}\int M_i\,\mu_{\dis}^{N,a,\gamma_1,\gamma_2}(\mathd M)
\leq\liminf_{N\rightarrow+\infty}
\int\left(\int_{i}^{i+1}x(\xi)\mathd\xi\right)\mu^{N,a}(\mathd x).
\]
Let $\mu_+$ be the limit of a converging subsequence of $\{\mu^{N,a}\}_{N\in\N}$, which exists and solves the DLR equations by Lemma \ref{lemma:existencetightnessDir} and Theorem \ref{theorem:existence2}. By Lemma \ref{lemma:existencetightnessDir} we can pass to the limit in the last integral, and we obtain
\[
\int\left(\int_{i}^{i+1}x(\xi)\mathd\xi\right)\mu_+(\mathd x)>0.
\]
\noindent Finally, the image of $\mu_+$ under $x\mapsto-x$ is a DLR measure $\mu_-$
and
\[
\int \left(\int_{i}^{i+1}x(\xi)\mathd \xi\right)\mu_-(\mathd x)<0.
\]
Thus $\mu_+\neq\mu_-$. The same uniform moment estimates pass to the limits, so that, by Theorem \ref{theorem:infinitesimalsymmetri}, both measures are
reversible and hence invariant for \eqref{eq:LangevinSPDE1}.
\end{proof}

\part{Abstract results on SPDEs}\label{part:abstract}  

\section{Preliminaries on Abstract SPDEs}\label{sect: Abstract-SPDE}
In this section we will consider the following SPDE
\begin{gather}\label{SPDE}
\eqsys{
\vspace{1pt}
\mathd X(t)=\big[AX(t)+F(X(t))+G(X(t))\big]\mathd t+R\mathd W(t), & t\in [0,T];\\
X(0)=x\in {H}.
}
\end{gather}
where the coefficients are defined by the following assumptions.

\begin{hyp}\label{Hyp:base} $ $
\begin{enumerate}
\item $H$ is a real separable Hilbert space and $E$ is a real separable reflexive Banach space continuously and densely embedded in $H$.
\vspace{3pt}
\item $\{W(t)\}_{t\geq 0}$ is an $H$-cylindrical Wiener process.
\vspace{3pt}
\item $R:H\rightarrow H$ is a linear bounded operator.
\vspace{3pt}
\item $G:H\rightarrow H$ and $G:E\rightarrow E$ is a twice Fréchet differentiable function with bounded and continuous derivatives and there exists $\zeta\in\R$ such that
\[
\scal{DG(x)h}{h}_H\leq\zeta\norm{h}^2_H,\qquad x,h\in H.
\]
and $G-\zeta\Id_E$ is dissipative in $E$.
\item $A:{\rm Dom}(A)\subseteq H\rightarrow H$ is a linear operator that generates a strongly continuous semigroup $\{e^{tA}\}_{t\geq 0}$ on $H$. $A_E:{\rm Dom}(A_E)\subseteq E\rightarrow E$ is a linear operator that generates a strongly continuous semigroup. Moreover, there exists $\varepsilon>0$, such that both $A+\varepsilon\Id_H$ and $A_E+\varepsilon\Id_E$ are dissipative.
\vspace{3pt}
\item $F:E\rightarrow H$ is a continuous function such that
\begin{enumerate}
\item $F$ is $m$-dissipative in ${H}$ and $F_{E}:{\rm Dom}(F_E)\subseteq E\rightarrow E$ is $m$-dissipative;
\item there exists $q\geq 1$ and $C>0$ such that
\[
\norm{F(x)}_H\leq C(1+\norm{x}_E^{q}),\quad x\in E.
\]
\end{enumerate}

\item The process $\{W_A(t)\}_{t\geq 0}$ defined by
\[
W_A(t):=\int_0^te^{(t-s)A}R\mathd W(s),
\]
is a continuous process in $E$ and for every $p\geq 1$ and $T\geq 0$ we have
\begin{equation}\label{Bound-convu}
\mathbb{E}\left[\sup_{t\in [0,T]}\norm{F\left(W_A(t)\right)}_E^p+\sup_{t\in [0,T]}\norm{W_A(t)}_E^q\right]<+\infty.
\end{equation}

\item\label{Hyp:base8} For every $B>0$ there exists $L_B>0$ such that 
\begin{equation}\label{FdualEH}
\abs{\scal{F(x)-F(y)}{h}_{H}}\leq L_B\norm{x-y}_H\norm{h}_E,\qquad  h\in E,
\end{equation}
for every $x,y\in E$ with $\norm{x}_E\leq B$ and $\norm{y}_E\leq B$.

\end{enumerate}
\end{hyp}
\begin{rmk}
    We recall that a (possibly nonlinear) function $\mathcal T:{\rm Dom}(\mathcal T)\subseteq E\rightarrow E$ is dissipative in $E$ if
    \[
    \norm{x-y-\lambda\big(\mathcal T(x)-\mathcal T(y)\big)}_E\geq\norm{x-y}_E,\qquad x,y\in{\rm Dom}(\mathcal T),\ \lambda>0,
    \]
    and that it is $m$-dissipative in $E$ if, in addition, the range of $\Id_E-\lambda\mathcal T$ is $E$ for some (equivalently, for every) $\lambda>0$; the analogous definitions in $H$ are obtained by replacing $E$ with $H$, and in this case dissipativity is equivalent to $\scal{\mathcal T(x)-\mathcal T(y)}{x-y}_H\leq0$ for every $x,y\in{\rm Dom}(\mathcal T)$ (see \cite[Appendix D]{Dap-Zab2014}). In Hypotheses \ref{Hyp:base}, $A_E$ and $F_E$ denote the parts of $A$ and $F$ in $E$, namely ${\rm Dom}(A_E):=\{x\in{\rm Dom}(A)\cap E:\ Ax\in E\}$, $A_Ex:=Ax$, ${\rm Dom}(F_E):=\{x\in E:\ F(x)\in E\}$ and $F_E(x):=F(x)$.
    The condition that there exists $\varepsilon>0$, such that both $A+\varepsilon\Id_H$ and $A_E+\varepsilon\Id_E$ are dissipative implies that $A$ and $A_E$ are dissipative. However, this is not a restrictive condition. Indeed, if $A$ does not satisfy the dissipativity condition, we can consider $A + \eta \Id$ for some $\eta > 0$ instead of $A$, and $G - \eta \Id$ instead of $G$. In this way, the SPDE remains unchanged, and the condition is satisfied.
\end{rmk}

We summarize the results presented in this section.
Subsection \ref{Ssect: well-posed-mild} addresses the well-posedness of equation \eqref{SPDE}.
In Subsection \ref{Ssect: app-F}, we introduce a useful regularizing sequence for the function $F$. Subsection \ref{Ssect: app-mild} leverages this sequence to construct an approximation of the mild solution to \eqref{SPDE}. Finally, in Subsection \ref{Ssect: derivata-mild}, we recall results concerning the regularity of the mild solution to the approximated equation \eqref{SPDE} with respect to the initial data.

\subsection{Existence and uniqueness of a generalized mild solution to \eqref{SPDE}}\label{Ssect: well-posed-mild}
In this subsection, we summarize some results about generalized mild solution to \eqref{SPDE}.

\begin{thm}[Theorems 3.16 and 3.17 of \cite{Big-Pri2026}]\label{th:solmild}
Assume that Hypotheses \ref{Hyp:base} hold true. For every $x\in E$, \eqref{SPDE} has unique $E$-valued continuous mild solution $\{X(t,x)\}_{t\geq 0}$. Moreover for every $t>0$ and $x,z\in E$ we have
{\small
\begin{align}
&\norm{X(t,x)}_E\leq e^{\zeta t}\norm{x}_E+3\int_0^te^{\zeta (t-s)}\left(\norm{[F+G](W_A(s))}_E+\max\{0,\zeta\}\norm{W_A(s)}_E\right)\mathd s+\norm{W_A(t)}_E,\;\;\mathbb{P}\text{-a.s.}\label{solmild1}	\\
&\norm{X(t,x)-X(t,z)}_E\leq  e^{\zeta t}\norm{x-z}_E,\;\;\mathbb{P}\text{-a.s.}\label{solmild2}
\end{align}
}
\noindent For every $x\in H$, \eqref{SPDE} has unique $H$-valued continuous generalized mild solution $\{X(t,x)\}_{t\geq 0}$, namely for every $\{x_n\}_{n\in\N}\subseteq E$ such that $x_n\rightarrow x$ in $H$ as $n\rightarrow+\infty$ we have
\begin{equation}
\lim_{n\rightarrow+\infty}\mathbb{E}\left[\sup_{t\in [0,T]}\norm{X(t,x_n)-X(t,x)}_H\right]=0,\label{eq:mild-gen}
\end{equation}
where $\{X(t,x_n)\}_{t\in [0,T]}$ is the unique mild solution to \eqref{SPDE} with initial datum $x_n$ given by the first part of this theorem. Moreover for every $t>0$ and $x,z\in E$ we have
{\small
\begin{align}
&\norm{X(t,x)}_{H}\leq e^{\zeta t}\norm{x}_{H}+3\int_0^te^{\zeta (t-s)}\left(\norm{[F+G](W_A(s))}_H+\max\{0,\zeta\}\norm{W_A(s)}_H\right)\mathd s+\norm{W_A(t)}_{H},\;\; \mathbb{P}\text{-a.s.}\label{solmild3}\\
&\norm{X(t,x)-X(t,z)}_{H}\leq  e^{\zeta t}\norm{x-z}_{H},\;\;\mathbb{P}\text{-a.s.}\label{solmild4}
\end{align}
}
\end{thm}

By condition \eqref{Bound-convu} and \eqref{solmild1}, we immediately deduce the following result.
\begin{prop}\label{prop:boundmild}
   Assume that Hypotheses \ref{Hyp:base} hold true and let $q\geq 1$ be the constant given by Hypotheses \ref{Hyp:base}. There exist $c_q\geq 2$ and $w_q>0$ such that for every  $t>0$ and $x\in E$ we have
    \begin{equation}\label{Bound-mild}
    \mathbb{E}\left[\norm{X(t,x)}^q_E\right]\leq 2^{-1}c_qe^{w_qt}(1+\norm{x}^q_E),
    \end{equation}
    where $\{X(t,x)\}_{t\geq 0}$ is the unique mild solution to \eqref{SPDE}.
\end{prop}

We denote by $\{P(t)\}_{t\geq 0}$ the transition semigroup associated to \eqref{SPDE} defined as 
\begin{equation}\label{semigruppo-transizione}
P(t)\Phi(x):=\mathbb{E}\left[\Phi(X(t,x))\right],\qquad \Phi\in C_b(H),\,t\geq 0,\,x\in H.
\end{equation}
It is possible to define a transition semigroup $\{P_E(t)\}_{t\geq 0}$ associated to \eqref{SPDE} even on $C_b(E)$. In particular
\[
P(t)\Phi(x)=P_E(t)\Phi(x),\qquad \Phi\in C_b(H),\,t\geq 0\,x\in E,
\]
hence, for simplicity, from here on we denote $\{P_E(t)\}_{t\geq 0}$ by $\{P(t)\}_{t\geq 0}$.

\subsection{The approximation for $F$}\label{Ssect: app-F}
Thanks to the $m$-dissipativity of $F:E\subseteq H\rightarrow H$ we can introduce the Yosida approximating (see \cite[Appendix A]{Cer2001}, \cite[Appendix D]{Dap-Zab2014}). Indeed for every $\delta>0$ and $x\in H$, there exists a unique $J_\delta(x)\in {\rm Dom}(F)$ such that
\begin{align}\label{eq:yosida}
J_\delta(x)-\delta F(J_\delta(x))=x.
\end{align}
For every $\delta>0$, we define $F_\delta:H\rightarrow H$ as
\begin{align}\label{def:yosida}
F_\delta(x):=F(J_\delta(x)),\qquad x\in H.
\end{align}

\begin{lemm}\label{lemma:yosida}
For every $\delta>0$ the function $F_\delta:H\rightarrow H$ is dissipative and Lipschitz continuous, and it satisfies the following inequalities
\begin{align}
&\lim_{\delta\rightarrow 0}\norm{F_{\delta}(x)-F(x)}_{H}=0,\qquad x\in E,\label{yosida1}\\
&\|F_{\delta}(x)\|_{H}\leq \|F(x)\|_{H}, \qquad x\in E,\label{yosida2}\\
&\norm{F_\delta(x)-F_\delta(y)}_{H}\leq 2\delta^{-1}\norm{x-y}_{H},\qquad x,y\in H.\label{yosida3}
\end{align}
\end{lemm}

Now we perform a further approximation to regularize $F_\delta$ on $H$. Let $H_0$ be a separable Hilbert space such that $H_0\hookrightarrow E\hookrightarrow H$ continuously, and let $\mathcal U:{\rm Dom}(\mathcal U)\subseteq H\rightarrow H$ be a negative self-adjoint operator such that
\begin{enumerate}[label=(U\arabic*)]
\item\label{U1} $\{e^{\beta\mathcal U}\}_{\beta\geq0}$ leaves $H_0$ invariant, its part $\mathcal U_0$ in $H_0$ is negative self-adjoint and $(-\mathcal U_0)^{-1/2}\in\mathcal L_2(H;H_0)$;
\item\label{U2} $M_E:=\sup_{\beta\in[0,1]}\norm{e^{\beta\mathcal U}}_{\mathcal L(E)}<+\infty$.
\end{enumerate} 
By \cite[Theorem 3.15]{Eng-Nag2006}  the operator $\mathcal{U}$ is the generator of a contractive analytic semigroup $\{e^{\beta\mathcal{U}}\}_{\beta\geq 0}$ on $H$ and by \cite[Proposition 2.1.1]{Lun1995} for every $\beta>0$ we have
\begin{equation}\label{controllability}
e^{\beta\mathcal{U}}(H)\subseteq Q_\beta^{1/2}(H).
\end{equation}
where 
\[
Q_\tau=\int_0^\tau e^{2s\mathcal{U}}\mathd s=\frac{1}{2}(-\mathcal{U})^{-1}(\Id-e^{2\tau\mathcal{U}}).
\]
For every $\delta,\tau>0$ we set
\begin{equation}\label{OU-approximating}
F_{\delta,\tau}(x):=\int_He^{\tau \mathcal{U}}F_{\delta}(e^{\tau \mathcal{U}}x+y)\mu_\tau(\mathd y),
\end{equation}
where $\mu_\tau$ is the Gaussian measure with mean $0$ and covariance operator $Q_\tau$. Since $Q_\tau=\frac12(-\mathcal U)^{-1}(\Id-e^{2\tau\mathcal U})\leq\frac12(-\mathcal U)^{-1}$, by \ref{U1} the measure $\mu_\tau$ is concentrated on $H_0$ and
$\int_{H_0}\norm{y}^2_{H_0}\mu_\tau(\mathd y)\leq\frac12\norm{(-\mathcal U_0)^{-1/2}}^2_{\mathcal L_2(H_0)}$ for every $\tau>0$; hence, by the standard properties of Gaussian measures and $H_0\hookrightarrow E$, for every $p\geq1$ we have
\begin{equation}\label{concentrazioneE}
\sup_{\tau\in (0,1)}\int_E\norm{y}_E^p\mu_\tau(\mathd y)<+\infty.
\end{equation}

\begin{prop}\label{prop:yosida-OU}
For every $\delta,\tau>0$ the function $F_{\delta,\tau}:H\rightarrow H
$ is twice Fréchet differentiable with continuous and bounded derivatives. For every $\delta>0$ we have
\begin{align}\label{yosida-OU1}
\lim_{\tau\rightarrow 0}\norm{F_{\delta,\tau}(x)-F_\delta(x)}_{H}=0,\qquad x\in H.
\end{align}
There exists a constant $M>0$ such that for every $\delta,\tau\in (0,1)$ we have
\begin{align}
&\norm{F_{\delta,\tau}(x)}_{H}\leq M(1+\norm{x}_{E}^q),\qquad x\in E,\label{yosida-OU2}\\
&\scal{DF_{\delta,\tau}(x)h}{h}_{H}\leq 0,\qquad x,h\in H,\label{yosida-OU3}
\end{align}
where $q\geq 1$ is the constant given by Hypotheses \ref{Hyp:base}(6a).
\end{prop}
\begin{proof}
Fix $\delta,\tau>0$. By \eqref{controllability} and the Cameron-Martin formula $F_{\delta,\tau}\in C^{1}(H;H)$, in particular, for every $h,v\in H$, we have
\begin{equation}\label{Der-cam}
\scal{DF_{\delta,\tau}(x)h}{v}_{H}=\int_{H}\scal{Q^{-1/2}_\tau e^{\tau\mathcal{U}}h}{Q^{-1/2}_\tau y}_{H}\scal{F_{\delta}(e^{\tau\mathcal{U}}x+y)}{v}_{H}\mu_\tau(\mathd y).
\end{equation}
Similarly, applying the Cameron-Martin formula again to $DF_{\delta,\tau}$, we obtain that $F_{\delta,\tau}$ belongs to $C^2(H;H)$.
By Lemma \ref{lemma:yosida} and \eqref{OU-approximating}, for every $x,h\in H$ we have 
\[
\scal{F_{\delta,\tau}(x+h)-F_{\delta,\tau}(x)}{h}_H=\int_H\scal{F_{\delta}(e^{\tau\mathcal{U}}(x+h)+y)-F_{\delta}(e^{\tau\mathcal{U}}x+y)}{e^{\tau\mathcal{U}}h}_H\mu_\tau(\mathd y)\leq 0,
\]
and so $F_{\delta,\tau}$ is dissipative. In particular, since $F_{\delta,\tau}$ is Fréchet differentiable then \eqref{yosida-OU3} holds. Since $F_{\delta}$ is Lipschitz continuous, by similar calculations we prove that $F_{\delta,\tau}$ is Lipschitz continuous and in particular it has bounded first order Fréchet derivative. Since $F_{\delta,\tau}\in C^2(H;H)$, it is sufficient to prove that $DF_{\delta,\tau}:H\rightarrow\mathcal{L}(H)$ is Lipschitz continuous to guarantee that also the second order derivative of $F_{\delta,\tau}$ is bounded. Given an orthonormal basis $\{e_k\}_{k\in\N}$ of $H$, by \eqref{Der-cam}, for every $z,w,h\in H$ we get
\begin{align*}
&\norm{DF_{\delta,\tau}(z)h-DF_{\delta,\tau}(w)h}^2_{H}=\sum_{k=1}^{+\infty}\scal{DF_{\delta,\tau}(z)h-DF_{\delta,\tau}(w)h}{e_k}_{H}^2\\
&\phantom{aaa}=\sum_{k=1}^{+\infty}\left(\int_{H}\scal{Q^{-1/2}_\tau e^{\tau\mathcal{U}}h}{Q^{-1/2}_\tau y}_{H}\scal{F_{\delta}(e^{\tau\mathcal{U}}z+y)-F_{\delta}(e^{\tau\mathcal{U}}w+y)}{e_k}_{H}\mu_{\tau}(\mathd y)\right)^2\\
&\phantom{aaa}\leq \int_{H}\scal{Q^{-1/2}_\tau e^{\tau\mathcal{U}}h}{Q^{-1/2}_\tau y}^2_{H}\sum_{k=1}^{+\infty}\scal{F_{\delta}(e^{\tau\mathcal{U}}z+y)-F_{\delta}(e^{\tau\mathcal{U}}w+y)}{e_k}_{H}^2\mu_{\tau}(\mathd y)\\
&\phantom{aaa}= \int_{H}\scal{Q^{-1/2}_\tau e^{\tau\mathcal{U}}h}{Q^{-1/2}_\tau y}_{H}^2\norm{F_{\delta}(e^{\tau\mathcal{U}}z+y)-F_{\delta,\tau}(e^{\tau\mathcal{U}}w+y)}^2_{H}\mu_{\tau}(\mathd y),
\end{align*}
and so $DF_{\delta,\tau}:H\rightarrow\mathcal{L}(H)$ is Lipschitz continuous.

We now prove \eqref{yosida-OU2}. Recalling that $\{e^{\tau\mathcal{U}}\}_{\tau\geq 0}$ is contractive, by Hypotheses \ref{Hyp:base}(6a), \eqref{U2}, \eqref{yosida2}, \eqref{OU-approximating} and \eqref{concentrazioneE}, for every $x\in E$ we have
\begin{align*}
\norm{F_{\delta,\tau}(x)}_H&\leq \int_E \norm{F_{\delta}(e^{\tau \mathcal{U}}x+y)}_H\mu_\tau(\mathd y)\leq \int_E \norm{F(e^{\tau \mathcal{U}}x+y)}_H\mu_\tau(\mathd y)\\
&\leq C\int_E \left(1+\norm{e^{\tau\mathcal{U}}x+y}^{q}_E\right)\mu_\tau(\mathd y)\leq C\left[1+2^{q-1}\left(M_E^q\norm{x}^q_E+\int_E\norm{y}_E^q\mu_\tau(\mathd y)\right)\right],
\end{align*}
and so that \eqref{yosida-OU2} holds.

We prove \eqref{yosida-OU1}. For every $\tau>0$ and $x\in H$ we have
\begin{align*}
F_{\delta,\tau}(x)-F_\delta(x)&=\int_{H} e^{\tau\mathcal{U}}\left(F_\delta(e^{\tau \mathcal{U}}x+y)-F_\delta(x)\right)\mu_\tau(\mathd y)\\
&+(e^{\tau\mathcal{U}}-\Id)F_\delta(x),
\end{align*}
by \eqref{yosida3} and since $\{e^{\tau\mathcal{U}}\}_{\tau\geq 0}$ is contractive there exists $C_\delta>0$ such that
\begin{align*}
\norm{F_{\delta,\tau}(x)-F_\delta(x)}_{H}&\leq C_\delta\left[\norm{(e^{\tau \mathcal{U}}-\Id)x}_{H} +\norm{(e^{\tau \mathcal{U}}-\Id)F_\delta(x)}_{H}+\int_{H} \norm{y}_{H}\mu_\tau(\mathd y)\right],
\end{align*}
since $\{e^{\tau\mathcal{U}}\}_{\tau\geq 0}$ is a strongly continuous semigroup  on $H$ and $\mu_\tau$ converges weakly to $\delta_0$ as $\tau\rightarrow 0$ then \eqref{yosida-OU1} holds true. 
\end{proof}

\subsection{Approximation of \eqref{SPDE}}\label{Ssect: app-mild}$ $\\
Now we exploit the previous approximation of $F$ to approximate $\eqref{SPDE}$. For every $\delta>0$ we denote by $F_\delta$ the Yosida approximations of $F$ given by \eqref{def:yosida} and we consider the SPDE
\begin{gather}\label{SPDE-delta}
\eqsys{
\mathd X_{\delta}(t)=\big(AX_{\delta}(t)+F_\delta(X_{\delta}(t))+G(X_{\delta}(t))\big)\mathd t+R\mathd W(t), & t>0;\\
X_\delta(0)=x\in H.
}
\end{gather}

\begin{lemm}\label{lemma:yosida-spde}
Assume that Hypotheses \ref{Hyp:base} hold true. For every $\delta>0$ and $x\in H$ the SPDE \eqref{SPDE-delta} has unique $H$-valued continuous mild solution $\{X_\delta(t,x)\}_{t\geq 0}$. Moreover for every $t>0$ it holds $\mathbb{P}\text{-a.s.}$
{\small
\begin{align}
&\norm{X_\delta(t,x)}_{H}\leq e^{\zeta t}\norm{x}_{H}+3\int_0^te^{\zeta (t-s)}\left(\norm{[F+G](W_A(s))}_H+\max\{0,\zeta\}\norm{W_A(s)}_H\right)\mathd s+\norm{W_A(t)}_{H},\quad x\in H,\label{yosida-spde1}\\
&\norm{X_\delta(t,x)}_E\leq e^{\zeta t}\norm{x}_E+3\int_0^te^{\zeta (t-s)}\left(\norm{[F+G](W_A(s))}_E+\max\{0,\zeta\}\norm{W_A(s)}_E\right)\mathd s+\norm{W_A(t)}_E,\quad x\in E\label{yosida-spde2}	\\
&\norm{X_\delta(t,x)-X_\delta(t,z)}_H\leq  e^{\zeta t}\norm{x-z}_H, \quad x,z\in H\label{yosida-spde3},\\
&\norm{X_\delta(t,x)-X_\delta(t,z)}_{E}\leq  e^{\zeta t}\norm{x-z}_{E},\quad x,z\in E\label{yosida-spde4}\\
&\lim_{\delta\rightarrow 0}\mathbb{E}\left[\norm{X_\delta(t,x)-X(t,x)}_{H}\right]=0,\quad x\in E.\label{yosida-spde5}
\end{align}
}
\end{lemm}
\begin{proof}
We note that
\[
X_\delta(t,x)=Y_\delta(t,x)+W_A(t),\qquad t\geq 0,\ x\in H, \mathbb{P}\text{-a.s.}
\]
where $\{Y_\delta(t,x)\}_{t\geq 0}$ is the unique mild solution to 
\begin{gather*}
\eqsys{
dY_{\delta}(t)=\big(AY_{\delta}(t)+F_\delta(Y_{\delta}(t)+W_A(t))+G(Y_{\delta}(t)+W_A(t))\big)\mathd t, & t>0;\\
Y_\delta(0)=x\in H.
}
\end{gather*}
Due to this observation, the statements follow by \cite[Proposition B.6]{Big-Pri2026}.
\end{proof}

For every $\delta,\tau>0$ let $F_{\delta,\tau}$ be the function defined by \eqref{OU-approximating} and we consider the SPDE
\begin{gather}\label{SPDE-delta-beta}
\eqsys{
\mathd X_{\delta,\tau}(t)=\big(AX_{\delta,\tau}(t)+F_{\delta,\tau}(X_{\delta,\tau}(t))+G(X_{\delta,\tau}(t))\big)\mathd t+R\mathd W(t), & t>0;\\
X_{\delta,\tau}(0)=x\in H.
}
\end{gather}

\begin{lemm}\label{lemma:yosida-beta-spde}
Assume that Hypotheses \ref{Hyp:base} hold true. For every $\delta,\tau >0$ and $x\in H$ the SPDE \eqref{SPDE-delta-beta} has a unique $H$-valued continuous mild solution $\{X_{\delta,\tau}(t,x)\}_{t\geq 0}$. Moreover for every $\tau,\delta>0$ and $t>0$ it holds $\mathbb{P}\text{-a.s.}$
{\small
\begin{align}
&\norm{X_{\delta,\tau}(t,x)}_{H}\leq e^{\zeta t}\norm{x}_{H}+3\int_0^te^{\zeta (t-s)}\left(\norm{[F+G](W_A(s))}_H+\max\{0,\zeta\}\norm{W_A(s)}_H\right)\mathd s+\norm{W_A(t)}_{H},\;\; x\in H,\label{yosida-beta-spde1}\\
&\norm{X_{\delta,\tau}(t,x)-X_{\delta,\tau}(t,z)}_H\leq  e^{\zeta t}\norm{x-z}_H, \quad x,z\in H\label{yosida-beta-spde3},\\
&\lim_{\tau\rightarrow 0}\mathbb{E}\left[\norm{X_{\delta,\tau}(t,x)-X_\delta(t,x)}_{H}\right]=0,\quad x\in H.\label{yosida-beta-spde5}
\end{align}
}
\end{lemm}
\begin{proof}
The first part of the proof follows by applying Lemma \ref{lemma:yosida-spde} with $F_\delta$ replaced by $F_{\delta,\tau}$. It remains to prove \eqref{yosida-beta-spde5}. Fix $x\in H$ and $\delta>0$. By Hypotheses \ref{Hyp:base} and the Lumer-Phillips theorem $\{e^{tA}\}_{t\geq 0}$ is a strongly continuous and contractive semigroup on $H$. By Proposition \ref{prop:yosida-OU} $F_{\delta,\tau}:H\rightarrow H$ is Lipschitz-continuous with Lipschitz constant independent of $\tau$. Hence there exists $C_\delta$ such that for every $\tau>0$ we have
\begin{align*}
\mathbb{E}\left[\norm{X_{\delta,\tau}(t,x)-X_{\delta}(t,x)}_{H}\right]& \leq \int_0^t\norm{e^{(t-s)A}}_{\mathcal{L}(H)}\mathbb{E}\left[\norm{F_{\delta,\tau}(X_{\delta,\tau}(s,x))-F_{\delta,\tau}(X_{\delta}(s,x))}_{H}\right]\mathd s\\
&+\int_0^t\norm{e^{(t-s)A}}_{\mathcal{L}(H)}\mathbb{E}\left[\norm{F_{\delta,\tau}(X_{\delta}(s,x))-F_{\delta}(X_{\delta}(s,x))}_{H}\right]\mathd s\\
&\leq C_\delta\int_0^t\mathbb{E}\left[\norm{X_{\delta,\tau}(s,x)-X_{\delta}(s,x)}_{H}\right]\mathd s\\
&+\int_0^t\mathbb{E}\left[\norm{F_{\delta,\tau}(X_{\delta}(s,x))-F_{\delta}(X_{\delta}(s,x))}_{H}\right]\mathd s.
\end{align*}
Hence by the Gronwall inequality there exists a constant $C'_\delta>0$ such that for every $\tau>0$ we have
\[
\mathbb{E}\left[\norm{X_{\delta,\tau}(t,x)-X_\delta(t,x)}_{H}\right]\leq C_\delta'\int_0^t\mathbb{E}\left[\norm{F_{\delta,\tau}(X_{\delta}(s,x))-F_{\delta}(X_{\delta}(s,x))}_{H}\right]\mathd s,
\]
hence \eqref{yosida-beta-spde5} follows from \eqref{yosida-OU1} and the dominated convergence theorem.
\end{proof}

By condition \eqref{Bound-convu} and \eqref{yosida-beta-spde1}, we immediately deduce the following result.
\begin{prop}
    Assume that Hypotheses \ref{Hyp:base} hold true. There exist $d_0\geq 2$ and $w_0>0$ such that for every $\delta,\tau >0$, $t>0$ and $x\in H$ we have
    \begin{equation}\label{Bound-mild-delta-beta}
    \mathbb{E}\left[\norm{X_{\delta,\tau}(t,x)}_H\right]\leq 2^{-1}d_0e^{w_0t}(1+\norm{x}_H),
     \end{equation}
    where $\{X_{\delta,\tau}(t,x)\}_{t\geq 0}$ is the unique mild solution to \eqref{SPDE-delta-beta}.
\end{prop}

For every $\delta,\tau >0$, we can define the transition semigroups associated to \eqref{SPDE-delta} and \eqref{SPDE-delta-beta} as
\begin{align}
&P_{\delta}(t)\Phi(x)=\mathbb{E}\left[\Phi(X_{\delta}(t,x))\right],\qquad t>0, x\in H, \Phi\in C_b(H),\label{app-semi-delta}\\
&P_{\delta,\tau}(t)\Phi(x)=\mathbb{E}\left[\Phi(X_{\delta,\tau}(t,x))\right],\qquad t>0, x\in H, \Phi\in C_b(H).\label{app-semi-delta-beta}
\end{align}
By \eqref{yosida-spde5} and \eqref{yosida-beta-spde5}, we deduce the following corollary.

\begin{cor}
For every $t\geq 0$, $x\in E$ and $\Phi\in C_b(H)$ we have
\begin{align}
&\lim_{\tau\rightarrow 0}\abs{P_{\delta,\tau}(t)\Phi(x)-P_\delta(t)\Phi(x)}=0,\qquad \delta>0,\label{convergenza-semi-delta-beta}\\
&\lim_{\delta\rightarrow 0}\abs{P_{\delta}(t)\Phi(x)-P(t)\Phi(x)}=0.\label{convergenza-semi-delta}
\end{align}
\end{cor}

\subsection{Estimate for derivative of the mild solution}\label{Ssect: derivata-mild}
We have the following result.
\begin{prop}\label{prop: reg-sol-mild}
Assume that Hypotheses \ref{Hyp:base} hold. For every $\delta,\tau >0$ and $t>0$, the map $x\in H\rightarrow X_{\delta,\tau}(t,x)\in H$ is $\mathbb{P}$-a.s. two times Gateaux differentiable with bounded derivatives. Moreover for every $\delta,\tau >0$ and $x,h\in H$ the process $\{DX_{\delta,\tau}(t,x)h\}_{t\geq 0}$ is the unique mild solution to 
\begin{gather}\label{SPDE-derivative}
\eqsys{
dY_{\delta,\tau}(t)=\big[A+DF_{\delta,\tau}(X_{\delta,\tau}(t,x))+DG(X_{\delta,\tau}(t,x))\big]Y_{\delta,\tau}(t)\mathd t, & t>0;\\
Y_{\delta,\tau}(0)=h,
}
\end{gather}
and
\begin{equation}\label{k-estimate}
\mathbb{E}\left[\norm{DX_{\delta,\tau}(t,x)h}_{H}\right]\leq e^{\zeta t}\norm{h}_{H}
\end{equation}
where $\zeta$ is the constant given by Hypotheses \ref{Hyp:base}.
\end{prop}
\begin{proof}
    By Proposition \ref{prop:yosida-OU}, for every $\delta,\tau>0$ the function $F_{\delta,\tau}$ is two times Fréchet differentiable with bounded and continuous derivatives, hence by \cite[Theorem 9.8 and 9.9]{Dap-Zab2014} or \cite[Theorem 3.6]{Dap2004} then the map $x\rightarrow X(t,x)$ is $\mathbb{P}$-a.s. two times Gateaux differentiable from $H$ into itself for every $t>0$. Estimate \eqref{k-estimate} follows by Hypotheses \ref{Hyp:base}, \eqref{yosida-OU3} and \cite[Theorem 3.6]{Dap2004}.
\end{proof}

By \eqref{app-semi-delta-beta}, \eqref{k-estimate} and the chain rule we obtain the following result.

\begin{cor}
Assume that Hypotheses \ref{Hyp:base} hold. For every $\delta,\tau>0$, $t>0$ and $\Phi\in C^1_b(H)$ we have
\begin{equation}\label{stima-grad-grad}
    \norm{\nabla P_{\delta,\tau}(t)\Phi(x)}_H\leq e^{\zeta t}\norm{\Phi}_{C^1_b(H)},\qquad x\in H.
\end{equation}
\end{cor}

\subsection{Strong Feller property}\label{Ssect: strong-Feller}

In this subsection we study the strong-Feller property of $\{P(t)\}_{t\geq 0}$ when $R=\alpha\Id_H$ for some $\alpha>0$..

\begin{prop}\label{prop:strongFeller}
Assume that Hypotheses \ref{Hyp:base} hold and that $R=\alpha\Id_H$ for some $\alpha>0$.
For every $t>0$, $\Phi\in B_b(H)$ and $x,z\in H$ we have
\begin{equation}\label{eq:strongFeller}
\abs{P(t)\Phi(x)-P(t)\Phi(z)}\leq K_t\norm{\Phi}_\infty\norm{x-z}_H,\qquad K_t:=\frac{e^{\max\{\zeta,0\}t}}{\alpha\sqrt t}.
\end{equation}
In particular $\{P(t)\}_{t\geq0}$ is strong Feller.
\end{prop}

\begin{proof}
Fix $t>0$ and $\delta,\tau>0$, and let $\{X_{\delta,\tau}(s,x)\}_{s\geq0}$ be the mild solution to \eqref{SPDE-delta-beta}. For $h\in H$ we set $\eta^h(s):=DX_{\delta,\tau}(s,x)h$, which by Proposition \ref{prop: reg-sol-mild} is the mild solution to \eqref{SPDE-derivative}. Since $A+\varepsilon\Id_H$ is dissipative, $\scal{DF_{\delta,\tau}(y)k}{k}_H\leq0$ by \eqref{yosida-OU3} and $\scal{DG(y)k}{k}_H\leq\zeta\norm{k}^2_H$ for every $y,k\in H$, the operator $A+DF_{\delta,\tau}(X_{\delta,\tau}(s,x))+DG(X_{\delta,\tau}(s,x))-\zeta\Id_H$ is dissipative for every $s\geq0$, so that
\begin{equation}\label{eq:pathwise-derivative}
\norm{\eta^h(s)}_H\leq e^{\zeta s}\norm{h}_H,\qquad s\geq0,\quad\mathbb P\text{-a.s.}
\end{equation}
Since $R=\alpha\Id_H$ is invertible and $F_{\delta,\tau}+G$ is twice Fr\'echet differentiable with bounded derivatives, the Bismut-Elworthy-Li formula (see e.g. \cite[Section 7.1]{Dap-Zab1996}) gives, for every $\Phi\in C_b(H)$,
\[
\scal{\nabla P_{\delta,\tau}(t)\Phi(x)}{h}_H=\frac{1}{\alpha t}\,\mathbb E\left[\Phi(X_{\delta,\tau}(t,x))\int_0^t\scal{\eta^h(s)}{\mathd W(s)}_H\right].
\]
Hence, by the Cauchy-Schwarz inequality, the It\^o isometry and \eqref{eq:pathwise-derivative},
\[
\abs{\scal{\nabla P_{\delta,\tau}(t)\Phi(x)}{h}_H}\leq\frac{\norm{\Phi}_\infty}{\alpha t}\left(\int_0^te^{2\zeta s}\mathd s\right)^{1/2}\norm{h}_H\leq K_t\norm{\Phi}_\infty\norm{h}_H,
\]
and therefore $\abs{P_{\delta,\tau}(t)\Phi(x)-P_{\delta,\tau}(t)\Phi(z)}\leq K_t\norm{\Phi}_\infty\norm{x-z}_H$ for every $x,z\in H$. Letting first $\tau\rightarrow0$ and then $\delta\rightarrow0$, by \eqref{convergenza-semi-delta-beta} and \eqref{convergenza-semi-delta} we obtain \eqref{eq:strongFeller} for every $\Phi\in C_b(H)$ and $x,z\in E$. By \eqref{solmild4} and the dominated convergence theorem $P(t)\Phi:H\rightarrow\R$ is continuous, hence, since $E$ is dense in $H$, \eqref{eq:strongFeller} holds for every $\Phi\in C_b(H)$ and $x,z\in H$. Finally, since Borel probability measures on the metric space $H$ are regular, the supremum of $\abs{\int_H\Phi\,\mathd(\nu_1-\nu_2)}$ over $\Phi\in C_b(H)$ with $\norm{\Phi}_\infty\leq1$ equals the total variation of $\nu_1-\nu_2$; applying this to the transition probabilities of $X(t,x)$ and $X(t,z)$ we obtain \eqref{eq:strongFeller} for every $\Phi\in B_b(H)$.
\end{proof}

\subsection{Accessibility and uniqueness of the invariant measure}\label{Ssect: irreducibility}

In this subsection we assume that Hypotheses \ref{Hyp:base} hold and the following additional condition.

\begin{hyp}\label{hyp:irr}
There exist $x_*\in{\rm Dom}(A)\cap E$ and $E_0\subseteq E$, dense in $H$, such that for every $x\in E_0$ and $T>0$ there exists $u\in L^2((0,T);H)$ such that $Ru(s)=-Ax_*-(F+G)\big(x_*+e^{sA}(x-x_*)\big)$ for a.e. $s\in[0,T]$.
\end{hyp}

\begin{rmk}\label{rmk:irr-alphaId}
If $R=\alpha\Id_H$ for some $\alpha>0$, then Hypothesis \ref{hyp:irr} holds with $x_*=0$ and $E_0=E$: indeed, for $x\in E$ we have $e^{\cdot A}x\in C([0,T];E)$ by \ref{Hyp:base}(5), so that, by \ref{Hyp:base}(4) and (6), $(F+G)(e^{\cdot A}x)\in C([0,T];H)$ and it suffices to take $u:=-\alpha^{-1}(F+G)(e^{\cdot A}x)$.
\end{rmk}

Following \cite[Definition 7.2]{Hairer2008}, we say that $x_0\in H$ is accessible for $\{P(t)\}_{t\geq0}$ if for every $x\in H$ and every open neighborhood $U$ of $x_0$ we have $\int_0^{+\infty}e^{-t}P(t)\mathbb I_U(x)\,\mathd t>0$.

\begin{prop}\label{prop:irreducibility}
Assume Hypotheses \ref{Hyp:base} and \ref{hyp:irr} hold. Then $x_*$ is accessible for $\{P(t)\}_{t\geq0}$. More precisely, for every $x\in H$ and $\eta>0$ there exists $t_0>0$ such that $\mathbb P\big(\norm{X(t,x)-x_*}_H<\eta\big)>0$ for every $t\geq t_0$.
\end{prop}

\begin{proof}
Fix $x\in H$ and $\eta>0$. We note that by Hypotheses \ref{Hyp:base}(5) and the Lumer-Phillips theorem we have $\norm{e^{tA}}_{\mathcal L(H)}\leq e^{-\varepsilon t}$ for every $t\geq0$. Hence, by \eqref{solmild4}, for every $t\geq 0$ and  $z\in E_0$ we have
\begin{align}\label{eq:triangolo1}
\norm{X(t,x)-x_*}_H&\leq\norm{X(t,x)-X(t,z)}_H+\norm{X(t,z)-\psi(t)}_H+\norm{e^{tA}(z-x_*)}_H\\
&\leq e^{\zeta t}\norm{x-z}_H+e^{-\varepsilon t}\norm{z-x_*}_H+\norm{X(t,z)-\psi(t)}_H,\qquad \mathbb{P}\text{-a.s.}
\end{align}
where $\psi(t):=x_*+e^{tA}(z-x_*)$ and $x_*$ is given by Hypotheses \ref{hyp:irr}. We  fix $t_0>0$ such that $e^{-\varepsilon t_0}\left(1+\norm{x-x_*}_H\right)\leq \eta/3$.\\
Since $E_0$ is dense in $H$, for every $t\geq t_0$ there exists $z_t\in E_0$, depending on $t$, such that $\norm{x-z_t}_H< \min\{1;e^{-\max\{0;\zeta\} t}\eta/3\}$. Noting that $\norm{z_t-x_*}_H\leq \norm{z_t-x}_H+\norm{x-x_*}_H<1+\norm{x-x_*}_H$, from \eqref{eq:triangolo1} we deduce
\begin{equation}\label{eq:triangolo3}
\norm{X(t,x)-x_*}_H< \frac{2\eta}{3}+\norm{X(t,z_t)-\psi(t)}_H, \qquad \mathbb{P}\text{-a.s.}
\end{equation}
Fix $t\geq t_0$ and $z_t\in E_0$ such that \eqref{eq:triangolo3} holds, to conclude the proof it is sufficient to prove that
\begin{equation}\label{eq:irrid1}
\mathbb{P}\left(\norm{X(t,z_t)-\psi(t)}_H< \frac{\eta}{3}\right)>0 .
\end{equation}

Let $u\in L^2((0,t);H)$ be given by Hypothesis \ref{hyp:irr} with $x=z_t$ and $T=t$. Since $x_*\in{\rm Dom}(A)$ we have $x_*-e^{sA}x_*=-\int_0^se^{(s-r)A}Ax_*\,\mathd r$, so that
\[
\psi(s)=e^{sA}z_t+\int_0^se^{(s-r)A}\big[(F+G)(\psi(r))+Ru(r)\big]\mathd r,\qquad s\in[0,t],
\]
and $\psi\in C([0,t];E)$ by Hypotheses \ref{Hyp:base}(5). The law $\gamma$ of $\{W_A(s)\}_{s\in[0,t]}$ is a centered Gaussian measure on the separable Banach space $C([0,t];E)$, whose Cameron-Martin space is $\{\int_0^\cdot e^{(\cdot-r)A}Rv(r)\mathd r:\ v\in L^2((0,t);H)\}$ (see \cite[Theorem 5.2(iii), Corollary B.5 and Proposition 2.10]{Dap-Zab2014}); in particular $g_u:=\int_0^\cdot e^{(\cdot-r)A}Ru(r)\mathd r$ belongs to it, and $g_u\in C([0,t];E)$.

We have $X(s,z_t)=Y(s)+W_A(s)$, where $Y(s)=e^{sA}z_t+\int_0^se^{(s-r)A}(F+G)(Y(r)+W_A(r))\mathd r$. Setting $d:=Y-\psi+g_u$ and $b(r):=(F+G)\big(\psi(r)+W_A(r)-g_u(r)\big)-(F+G)(\psi(r))$, $d$ is the mild solution of
\[
d'=Ad+(F+G)(Y+W_A)-(F+G)(Y-d+W_A)+b,\qquad d(0)=0 .
\]
Since $F+G-\zeta\Id_H$ is dissipative in $H$, by standard a priori estimates arguments we get $\norm{d(t)}_H\leq\int_0^te^{\zeta(t-r)}\norm{b(r)}_H\mathd r$, and, since $X(t,z_t)-\psi(t)=d(t)+W_A(t)-g_u(t)$,
\begin{equation}\label{eq:irr-comparison}
\norm{X(t,z_t)-\psi(t)}_H\leq t\,e^{\max\{\zeta,0\}t}\sup_{r\in[0,t]}\norm{b(r)}_H+c_E\norm{W_A-g_u}_{C([0,t];E)},
\end{equation}
where $c_E$ is the norm of the embedding $E\hookrightarrow H$.

The set $K:=\psi([0,t])$ is compact in $E$. Since $F:E\rightarrow H$ is continuous and $G$ is Lipschitz continuous in $H$, a standard compactness argument gives, for every $\delta>0$, a $\kappa_\delta>0$ such that $\norm{(F+G)(p+v)-(F+G)(p)}_H<\delta$ for every $p\in K$ and $v\in E$ with $\norm{v}_E<\kappa_\delta$. Choosing $\delta>0$ with $t\,e^{\max\{\zeta,0\}t}\delta<\eta/6$ and $\kappa:=\min\{\kappa_\delta,\eta/(6c_E)\}$, by \eqref{eq:irr-comparison} we get $\norm{X(t,z_t)-\psi(t)}_H<\eta/3$ on the event $\{\norm{W_A-g_u}_{C([0,t];E)}<\kappa\}$. 
Since $g_u$ belongs to the Cameron-Martin space of $\mu_{W_A}$ and, by \cite[Corollary B.5]{Dap-Zab2014}, its Cameron-Martin norm is bounded by $\norm{u}_{L^2((0,t);H)}$, by \cite[Lemma 12.7]{Dap-Zab2014} we have
\[
\mathbb P\big(\norm{W_A-g_u}_{C([0,t];E)}<\kappa\big)\geq e^{-\frac12\norm{u}^2_{L^2((0,t);H)}}\,\mathbb P\big(\norm{W_A}_{C([0,t];E)}<\kappa\big)>0,
\]
which gives \eqref{eq:irrid1}. By \eqref{eq:triangolo3} we conclude $\mathbb P\big(\norm{X(t,x)-x_*}_H<\eta\big)>0$ for every $t\geq t_0$.

Finally, let $U$ be an open neighborhood of $x_*$ and let $\eta>0$ be such that the ball of center $x_*$ and radius $\eta$ is contained in $U$. Then $P(t)\mathbb I_U(x)\geq\mathbb P\big(\norm{X(t,x)-x_*}_H<\eta\big)>0$ for every $t\geq t_0$, so that $\int_0^{+\infty}e^{-t}P(t)\mathbb I_U(x)\,\mathd t>0$, namely $x_*$ is accessible.
\end{proof}

\begin{cor}\label{cor:uniqueness-invariant}
Assume Hypothesis \ref{hyp:irr} holds with $R=\alpha\Id_H$ for some $\alpha>0$. Then $\{P(t)\}_{t\geq0}$ has at most one invariant measure.
\end{cor}

\begin{proof}
By Proposition \ref{prop:strongFeller}, $\{P(t)\}_{t\geq0}$ is strong Feller at every point of $H$, and by Proposition \ref{prop:irreducibility} the point $x_*$ is accessible. The statement follows from \cite[Corollary 7.8]{Hairer2008} (see also Lemma 7.3 and Proposition 7.7 therein).
\end{proof}

\section{The core}\label{sect:core}
We consider the following second-order differential operators
\begin{align}
&L_0\Phi(x):=\frac{1}{2}{\rm Tr}\left[RR^*\nabla^2\Phi(x)\right]+\scal{x}{A^*\nabla\Phi(x)}_H,\label{Op-OU}\\
 &   N_0\Phi(x):=L_0\Phi(x)+\scal{F(x)+G(x)}{\nabla\Phi(x)}_H\label{Op-OUP},
\end{align}
where $x\in E$, $\Phi\in \xi_A(H)$,
\begin{align}\label{xia}
\xi_A(H):={\rm Span}\{\mbox{real and imaginary parts of the functions } x\mapsto e^{i\scal{x}{h}}\,|\, h\in \mathcal{D}\},
\end{align}
and $\mathcal{D}\subseteq{\rm Dom}(A^*)\cap E$ is a linear subspace which is a core for $A^*$, namely for every $h\in{\rm Dom}(A^*)$ there exists $\{h_n\}_{n\in\N}\subseteq\mathcal D$ such that $h_n\rightarrow h$ and $A^*h_n\rightarrow A^*h$ in $H$ (in particular $\mathcal D$ is dense in $H$, and $\xi_A(H)$ is an algebra containing the constant functions). First of all, we note that the coefficients of the first order of the operators $L_0$ and $N_0$ are unbounded; to be precise, they exhibit polynomial growth. Consequently, we will study the operators $L_0$ and $N_0$ on the following weighted spaces of continuous functions.
\begin{rmk}\label{rmk:cilvsxiai}
We emphasize that all the results established in this section remain valid if $\xi_A(H)$ is replaced by the following space of smooth cylindrical functions: \begin{align*} \mathcal{F}C_b^{\infty}(H):=\bigg\{\Phi:H\rightarrow \R\ :\ \exists n\in\N,\ &\phi\in C^{\infty}_b(\R^n;\R),\ h_1,\ldots,h_n\in \mathcal{D}\ :\\ & \Phi(x)=\phi\left(\scal{x}{h_1}_{H},\ldots, \scal{x}{h_n}_{H}\right),\ \forall x\in H\bigg\}. \end{align*}
\end{rmk}

\begin{defi}\label{p-space}
Let $\mathcal{X}$ be a separable Banach space and let $p\in\N$.  We denote by $\C_{p}(\mathcal{X})$ the Banach space defined as 
\[
\C_{p}(\mathcal{X}):=\Bigg\lbrace \Phi:\mathcal{X}\rightarrow \R\ : \  x\mapsto\frac{\Phi(x)}{1+\norm{x}^p_{\mathcal{X}}}\text{ belongs to } C_b(\mathcal{X}) \Bigg\rbrace.
\] 
endowed with the norm
\[
\norm{\Phi}_{\C_p(\mathcal{X})}:=\sup_{x\in \mathcal{X}}\left(\frac{\vert \Phi(x)\vert}{1+\norm{x}^p_{\mathcal{X}}}\right),\qquad \Phi\in \C_{p}(\mathcal{X}).
\]
\end{defi}

In this section, we study the relation between the transition semigroup $\{P(t)\}_{t\geq 0}$ and $N_0$ on the spaces of continuous functions defined above.
However, in general, for $m\in\N$ the semigroup $\{P(t)\}_{t\geq 0}$ is not strongly continuous on $\C_m(H)$ or $\C_m(E)$. Various approaches have been introduced in the literature to handle this issue, including the notion of weakly continuous semigroups \cite{Cer1994,Cer1995}, the concept of $\pi$-semigroups \cite{Pri1999}, and methods based on the mixed topology \cite{Gol-Koc2001}. In many cases, these different approaches turn out to be equivalent. In this paper, we adopt the one based on the $\pi$-semigroup, so we begin to define the notion of $\pi$-convergence.

\begin{defi}\label{pi-conv}
Let $\mathcal{X}$ be a separable Banach space and let $p\in\N$.
\begin{enumerate}
\item Given a sequence $\{\Phi_{n}\}_{n\in\N}\subseteq \C_p(\mathcal{X})$ we say that $\Phi_{n}$ $\pi$-converge to $\Phi$ in $\C_p(\mathcal{X})$ if
\begin{align*}
    &\lim_{n\rightarrow+\infty}\Phi_{n}(x)=\Phi(x),\qquad \forall x\in\mathcal{X},\\
    &\sup_{n\in\N}\norm{\Phi_n}_{\C_p(\mathcal{X})}<+\infty.
\end{align*}
In this case, we use the notation
\[
\plim{n\to +\infty} \Phi_n=\Phi, \qquad \text{in } \C_p(\mathcal{X})
\]
\item Given a $k$-sequence $\{\Phi_{n_1,\ldots,n_k}\}_{n_1\in\N, \ldots, n_k\in\N}\subseteq \C_p(\mathcal{X})$ we say that $\Phi_{n_1,\ldots,n_k}$ $\pi$-converges to $\Phi$ in $\C_p(\mathcal{X})$ if for every $i\in\{1,\ldots,k-1\}$ there exists an $i$-sequence $\{\Phi_{n_1,\ldots,n_i}\}_{n_1,\ldots,n_i\in\N}\subseteq\C_p(\mathcal{X})$ such that, for all $n_1,\ldots,n_i\in\N$,
\[
\plim{n_{i+1}\to+\infty}\Phi_{n_1,\ldots,n_{i+1}}=\Phi_{n_1,\ldots,n_i}\quad\text{in }\C_p(\mathcal{X}),\qquad\text{and}\qquad \plim{n_1\to+\infty}\Phi_{n_1}=\Phi\quad\text{in }\C_p(\mathcal{X}).
\]
In this case, we use the notation
\[
\plim{n_1 \to+\infty,\ldots,n_k\rightarrow+\infty}\Phi_{n_1,\ldots,n_k}=\Phi, \qquad \text{in } \C_p(\mathcal{X}).
\]
Since the limits are pointwise, they are taken iteratively and cannot be replaced by a diagonal sequence; the uniform bound is required only with respect to the last index, for fixed values of the previous ones (see \cite[Definition 2.1]{Manca2009}).

\item Similarly, given a family $\{\Phi_{t}\}_{t\geq 0}\subseteq \C_p(\mathcal{X})$ and $\Phi\in \C_p(\mathcal{X})$. We use the notation
\[
\plim{t\to 0} \Phi_t=\Phi, \qquad \text{in } \C_p(\mathcal{X})
\]
to denote that
\begin{align*}
    &\lim_{t\rightarrow 0}\Phi_{t}(x)=\Phi(x),\qquad \forall x\in\mathcal{X},\\
    &\sup_{t\in (0,1)}\norm{\Phi_t}_{\C_p(\mathcal{X})}<+\infty.
\end{align*}
\end{enumerate}
\end{defi}

\begin{prop}\label{densità-xia}
    Let $\mathcal{X}$ be a separable Banach space continuously embedded in $H$ and let $p\in\N$. There exists a constant $D_p$ such that
    \begin{equation}\label{imme!}
    \norm{\Phi}_{\C_p(\mathcal{X})}\leq D_p\norm{\Phi}_{\C_1(H)},\qquad \Phi\in\C_1(H),
    \end{equation}
    in particular, if $\Phi\in\C_1(H)$ and $\{\Phi_n\}_{n\in\N}\subseteq\C_1(H)$ then
    \begin{equation}\label{imme}
    \plim{n\rightarrow+\infty}\Phi_n=\Phi,\quad \text{ in } \C_1(H),\quad\Rightarrow\quad \plim{n\rightarrow+\infty}{\Phi_n}_{|\mathcal{X}}=\Phi_{|\mathcal{X}},\quad \text{ in } \C_1(\mathcal{X}).
    \end{equation}
        Moreover, for every $p\in\N\cup\{0\}$ and $\Phi\in \C_{p}(\mathcal{X})$, where $\C_0(\mathcal X)=C_b(\mathcal X)$, there exist $l\in\N$ and $\{\Phi_{n_1,\ldots,n_l}\}\subseteq \xi_A(H)$ such that
    \begin{equation}\label{densità}
    \plim{n_1,\ldots,n_l\rightarrow +\infty}\Phi_{n_1,\ldots,n_l}=\Phi,\qquad \text{ in } \C_p(\mathcal{X}).
    \end{equation}  
\end{prop}
\begin{proof}
Since $p\geq 1$ and $\mathcal{X}$ is continuously embedded in $H$, there exists a constant $D_p>0$ such that for every $x\in \mathcal{X}$ we have
\begin{equation}\label{stimaEH}
\frac{1}{1+\norm{x}_\mathcal{X}^p}\leq \frac{D_p}{1+\norm{x}_H},
\end{equation}
and so, taking into account Definition \ref{p-space}, \eqref{imme!} holds. Taking into account Definitions \ref{p-space} and \ref{pi-conv}, estimate \eqref{imme!} yields \eqref{imme}.

Now we prove that \eqref{densità} holds. Fix $\Phi\in \C_p(\mathcal{X})$. We consider the sequence $\{\Phi_{k}\}_{k\in\N}\subseteq C_b(\mathcal{X})$ defined by
    \[
    \Phi_k(x)=\dfrac{\Phi(x)}{1+k^{-1}\norm{x}_\mathcal{X}^p},\qquad
    \]
    It is easy to see that
    \begin{equation}\label{D1}
    \plim{k\rightarrow +\infty}\Phi_k=\Phi,\qquad \text{ in } \C_p(\mathcal{X}).
    \end{equation}
    Fixed $k\in\N$, since $\mathcal{X}$ is continuously embedded in $H$, by \cite[Proposition 6.1.5]{Cer2001}, there exists a sequence $\{\Phi_{h,k}\}_{h\in\N}\subseteq C_b(H)$ such that
     \begin{equation}\label{D2}
    \plim{h\rightarrow +\infty}\Phi_{h,k}=\Phi_k,\qquad \text{ in } C_b(\mathcal{X}).
    \end{equation}
    Finally, fixed $k,h\in\N$, by \cite[Lemma 2.4]{Dap-Tub2001} there exists a $3$-sequence  $\{\Phi_{n_1,n_2,n_3,h,k}\}_{n_1,n_2,_3\in\N}\subseteq \xi_A(H)$ such that
     \begin{equation}\label{D3}
    \plim{n_1,n_2,n_3\rightarrow +\infty}\Phi_{n_1,n_2,n_3,h,k}=\Phi_{h,k},\qquad \text{ in } C_b(H).
    \end{equation}
    Here \cite[Lemma 2.4]{Dap-Tub2001} is stated for exponential functions with directions in ${\rm Dom}(A^*)$; since $\mathcal D$ is dense in $H$ and $\abs{e^{i\scal{x}{g_n}}-e^{i\scal{x}{g}}}\leq\norm{x}_H\norm{g_n-g}_H$, adding one more index we can take directions in $\mathcal D$. Combining \eqref{D1}, \eqref{D2} and \eqref{D3} we obtain \eqref{densità}.
\end{proof}

The main goal of this section is to prove that $\xi_A(H)$ is a $\pi$-core for the semigroup $\{P(t)\}_{t\geq 0}$ on $\C_q(E)$, where $q\geq1$ is the constant given by Hypotheses \ref{Hyp:base}. Before proceeding, we need to show that formula \eqref{semigruppo-transizione} indeed defines a semigroup on $\C_q(E)$.

\begin{prop}\label{semigruppoP}
Assume Hypotheses \ref{Hyp:base} hold. Let $q\geq1$ be the constant given by Hypotheses \ref{Hyp:base}. Then the formula \eqref{semigruppo-transizione} defines a semigroup $\{P(t)\}_{t\geq 0}$ on $\C_{q}(E)$ such that
\begin{equation}\label{Bound-semi}
\norm{P(t)}_{\mathcal{L}(\C_{q}(E))}\leq c_qe^{w_qt}
\end{equation}
where $q$ is the constant given by Hypotheses \ref{Hyp:base}, and $c_q\geq 2$ and $w_q>0$ are the same constants in \eqref{Bound-mild}.
\end{prop}
\begin{proof}
   Let us only show that $P(t)$ is a bounded linear operator on $\C_q(E)$, for every $t\geq 0$. For every $\Phi\in \C_q(E)$, $t\geq 0$ and $x\in E$ we have 
    \begin{align*}
        \abs{P(t)\Phi(x)}\leq \mathbb{E}\left[\dfrac{\abs{\Phi(X(t,x))}(1+\norm{X(t,x)}_E^q)}{1+\norm{X(t,x)}_E^q} \right]\leq \norm{\Phi}_{\C_q(E)}\left(1+\mathbb{E}\left[\norm{X(t,x)}_E^q\right]\right),
    \end{align*}
hence, by Hypotheses \eqref{Bound-mild} we have
\[
\abs{P(t)\Phi(x)}\leq c_qe^{w_qt}(1+\norm{x}_E^q)\norm{\Phi}_{\C_q(E)},
\]
and so that we conclude
\[
\norm{P(t)\Phi}_{\C_q(E)}\leq c_qe^{w_qt}\norm{\Phi}_{\C_q(E)}.
\]
\end{proof}
Now we can define the notion of $\pi$-generator of $\{P(t)\}_{t\geq 0}$ in $\C_q(E)$.
\begin{defi}\label{def:pi-gen-N-p}
    Assume Hypotheses \ref{Hyp:base} hold. Let $q\geq1$ be the constant given by Hypotheses \ref{Hyp:base}. We denote by $N_q$ the $\pi$-generator of $\{P(t)\}_{t\geq 0}$ in $\C_q(E)$, namely the operator $N_q:{\rm Dom}(N_q)\subseteq \C_{q}(E)\rightarrow \C_{q}(E)$ defined by
    \begin{align*}
        &{\rm Dom}(N_q):=\Bigg\lbrace\Phi\in \C_{q}(E) : \exists \Psi\in \C_{q}(E) : \plim{t\rightarrow 0}\dfrac{P(t)\Phi-\Phi}{t}= \Psi \text{ in } \C_q(E)\Bigg\rbrace\\
        &N_q\Phi(x):=\lim_{t\rightarrow 0}\dfrac{P(t)\Phi(x)-\Phi(x)}{t},\qquad \Phi\in {\rm Dom}(N_q),\ x\in E.
    \end{align*}
\end{defi}
We recall the following useful properties of the $\pi$-semigroup
\begin{prop}\label{prop:commutazione}
    Assume Hypotheses \ref{Hyp:base} hold. Let $q\geq1$ be the constant given by Hypotheses \ref{Hyp:base}. For every $\Phi\in {\rm Dom}(N_q)$ and $r\geq 0$ we have $P(r)\Phi\in{\rm Dom}(N_q)$ and
    \[
    \left[N_qP(r)\Phi\right](x)=P(r)\left[N_q\Phi\right](x),\qquad \ x\in E.
    \]
\end{prop}
\begin{proof}
    Fix $r\geq 0$ and $\Phi\in{\rm Dom}(N_q)$. For every $x\in E$, we have
\begin{align*}
    \lim_{t\rightarrow 0}\dfrac{P(t)[P(r)\Phi](x)-P(r)\Phi(x)}{t}=\lim_{t\rightarrow 0}P(r)\left[\dfrac{P(t)\Phi-\Phi}{t}\right](x),
\end{align*}
by Definitions \ref{pi-conv} and \ref{def:pi-gen-N-p} and the dominated convergence theorem we deduce
\begin{align*}
    \lim_{t\rightarrow 0}\dfrac{P(t)[P(r)\Phi](x)-P(r)\Phi(x)}{t}=P(r)N_q\Phi(x).
\end{align*}
Moreover, by \eqref{Bound-semi} we have
\begin{align*}
    \sup_{t\in (0,1)}\norm{\dfrac{P(t)P(r)\Phi-P(r)\Phi}{t}}_{\C_q(E)}\leq c_qe^{w_qr}\sup_{t\in (0,1)}\norm{\dfrac{P(t)\Phi-\Phi}{t}}_{\C_q(E)},
\end{align*}
and so
\begin{align*}
    \plim{t\rightarrow 0}\dfrac{P(t)P(r)\Phi-P(r)\Phi}{t}=P(r)N_q\Phi,\qquad \text{ in } \C_q(E).
\end{align*}
Hence, taking into account Definition \ref{def:pi-gen-N-p}, we conclude $P(r)\Phi\in{\rm Dom}(N_q)$ and
\begin{align*}
\left[N_qP(r)\Phi\right](x):=\lim_{t\rightarrow 0}\dfrac{P(t)[P(r)\Phi](x)-P(r)\Phi(x)}{t}=P(r)\left[N_q\Phi\right](x),\qquad \ x\in E.
\end{align*}

\end{proof}

Now we can state the main theorem of the section.
\begin{thm}\label{thm: coreN}
Assume that Hypotheses \ref{Hyp:base} hold true. Let $q\geq1$ be the constant given by Hypotheses \ref{Hyp:base}. $N_q$ is an extension in $\C_{q}(E)$ of the operator $N_0$ defined in \eqref{Op-OUP}. Moreover $\xi_A(H)$ is a $\pi$-core for $N_q$, namely for every $\Phi\in {\rm Dom}(N_q)$ there exist $k\in\N$ and $\{\Phi_{n_1,\ldots,n_k}\}_{n_1\in\N,\ldots,n_k\in\N}\subseteq\xi_A(H)$ such that
\begin{align}
\plim{n_1\rightarrow+\infty,\ldots,n_k\rightarrow+\infty}\Phi_{n_1,\ldots,n_k}&=\Phi,\qquad \text{in } \C_q(E), \label{coreN1}\\
\plim{n_1\rightarrow+\infty,\ldots,n_k\rightarrow+\infty}N_q\Phi_{n_1,\ldots,n_k}&=N_q\Phi,\qquad \text{in }\C_q(E).\label{coreN2}
\end{align}
\end{thm}

Before proving the previous theorem, we need to recall some preliminary results.
\subsection{Ornstein-Uhlenbeck operator}

We begin to recall some properties of the Ornstein-Uhlenbeck semigroup on $\C_1(H)$.

\begin{prop}[Proposition 2.2 of \cite{Manca2009}]\label{prop: semi-OU}
Assume that Hypotheses \ref{Hyp:base} hold true and that $F=G=0$. Then the formula \eqref{semigruppo-transizione} defines a semigroup $\{T(t)\}_{t\geq 0}$ on $\C_{1}(H)$.
\end{prop}

\begin{defi}
    We denote by $L_1$ the $\pi$-generator of $\{T(t)\}_{t\geq 0}$, namely the operator $L_1:{\rm Dom}(L_1)\subseteq \C_{1}(H)\rightarrow \C_{1}(H)$ defined by
    \begin{align*}
        &{\rm Dom}(L_1):=\Bigg\lbrace\Phi\in \C_{1}(H) : \exists \Psi\in \C_{1}(H) : \plim{t\rightarrow 0}\dfrac{T(t)\Phi-\Phi}{t}= \Psi \text{ in } \C_1(H)\Bigg\rbrace\\
        &L_1\Phi(x):=\lim_{t\rightarrow 0}\dfrac{T(t)\Phi(x)-\Phi(x)}{t},\qquad \Phi\in {\rm Dom}(L_1),\ x\in H.
    \end{align*}
\end{defi}

\begin{prop}[Proposition 4.3 of \cite{Manca2009}]\label{prop: Manca-core}
Assume that Hypotheses \ref{Hyp:base} hold true. $L_1$ is an extension in $\C_{1}(H)$ of the operator $L_0$ defined by \eqref{Op-OU}.
Moreover $\xi_A(H)$ is a $\pi$-core for $L_1$, namely for every $\Phi\in {\rm Dom}(L_1)$ there exists $k\in\N$ and $\{\Phi_{n_1,\ldots,n_k}\}_{n_1\in\N,\ldots,n_k\in\N}\subseteq\xi_A(H)$ such that
\begin{align}
\plim{n_1\rightarrow+\infty,\ldots,n_k\rightarrow+\infty}\Phi_{n_1,\ldots,n_k}&=\Phi,\qquad \text{ in } \C_1(H)\label{Mancacore1}\\
\plim{n_1\rightarrow+\infty,\ldots,n_k\rightarrow+\infty} L_0\Phi_{n_1,\ldots,n_k}&=L_1\Phi,\qquad \text{ in } \C_1(H)\label{Mancacore2}.
\end{align}
If in addition $\Phi\in C^1_b(H)$ then for every $h\in H$ we have
\begin{align}
\plim{n_1\rightarrow+\infty,\ldots,n_k\rightarrow+\infty}\scal{\nabla\Phi_{n_1,\ldots,n_k}}{h}_H=\scal{\nabla\Phi}{h}_H,\qquad \text{ in } C_b(H)\label{Mancacore3}.
\end{align}
\end{prop}
\begin{proof}
In \cite[Proposition 4.3]{Manca2009} the statement is proved for the space $\mathcal E_A(H)$ defined as $\xi_A(H)$ with $\mathcal D$ replaced by ${\rm Dom}(A^*)$. Hence it is enough to approximate every element of $\mathcal E_A(H)$ by a sequence in $\xi_A(H)$, adding one more index to the multi-sequences. By linearity, let $\Psi(x)=e^{i\scal{x}{k}_H}$ with $k\in{\rm Dom}(A^*)$, let $\{k_n\}_{n\in\N}\subseteq\mathcal D$ be such that $k_n\rightarrow k$ and $A^*k_n\rightarrow A^*k$ in $H$, and set $\Psi_n(x):=e^{i\scal{x}{k_n}_H}$. Then $\Psi_n\rightarrow\Psi$ and $\scal{\nabla\Psi_n}{h}_H=i\scal{k_n}{h}_H\Psi_n\rightarrow\scal{\nabla\Psi}{h}_H$ pointwise and boundedly, for every $h\in H$, and
\[
L_0\Psi_n(x)=\Big[-\frac12\norm{R^*k_n}^2_H+i\scal{x}{A^*k_n}_H\Big]\Psi_n(x)\longrightarrow L_0\Psi(x),\qquad \abs{L_0\Psi_n(x)}\leq C(1+\norm{x}_H),
\]
with $C$ independent of $n$, namely $L_0\Psi_n$ $\pi$-converges to $L_0\Psi$ in $\C_1(H)$.
\end{proof}

\begin{rmk}\label{rmk:Manca-hyp}
Let $\Phi\in{\rm Dom}(L_1)\cap C^1_b(H)$ and let $\{\Phi_{n_1,\ldots,n_k}\}\subseteq\xi_A(H)$ be given by Proposition \ref{prop: Manca-core}. By \eqref{Mancacore3} and the uniform boundedness principle, $\sup_{x\in H}\norm{\nabla\Phi_{n_1,\ldots,n_k}(x)}_H$ is bounded with respect to the last index. Hence, if $\mathcal T:E\rightarrow H$ is a continuous function such that $\norm{\mathcal T(x)}_H\leq C(1+\norm{x}^q_E)$ for every $x\in E$, then
\[
\plim{n_1\rightarrow+\infty,\ldots,n_k\rightarrow+\infty}\scal{\mathcal T}{\nabla\Phi_{n_1,\ldots,n_k}}_H=\scal{\mathcal T}{\nabla\Phi}_H,\qquad \text{ in } \C_q(E).
\]
\end{rmk}

\subsection{The transition semigroup in the Lipschitz continuous case}

In this subsection, we recall the results in the case where $F$ is Lipschitz continuous in $H$ contained in \cite{Manca2009}.

\begin{prop}[Proposition 2.2 of \cite{Manca2009}]\label{prop: semi-OUP}
Assume Hypotheses \ref{Hyp:base} hold. Then the formula \eqref{app-semi-delta-beta} defines a semigroup $\{P_{\delta,\tau}(t)\}_{t\geq 0}$ on $\C_{1}(H)$ such that
\begin{equation}\label{Bound-semi-delta-beta}
\norm{P_{\delta,\tau}(t)}_{\mathcal{L}(\C_{1}(H))}\leq d_0e^{w_0t}
\end{equation}
where $d_0\geq 2$ and $w_0>0$ are the same constants in \eqref{Bound-mild-delta-beta}.
\end{prop}

\begin{defi}
    For every $\delta,\tau>0$. We denote by $N_{\delta,\tau}$ the $\pi$-generator of $\{P_{\delta,\tau}(t)\}_{t\geq 0}$ in $\C_1(H)$, namely the operator $N_{\delta,\tau}:{\rm Dom}(N_{\delta,\tau})\subseteq \C_{1}(H)\rightarrow \C_{1}(H)$ defined by
    \begin{align*}
        &{\rm Dom}(N_{\delta,\tau}):=\Bigg\lbrace\Phi\in \C_{1}(H) : \exists \Psi\in \C_{1}(H) : \plim{t\rightarrow 0}\dfrac{P_{\delta,\tau}(t)\Phi-\Phi}{t}= \Psi \text{ in } \C_1(H)\Bigg\rbrace\\
        &N_{\delta,\tau}\Phi(x):=\lim_{t\rightarrow 0}\dfrac{P_{\delta,\tau}(t)\Phi(x)-\Phi(x)}{t},\qquad \Phi\in {\rm Dom}(N_{\delta,\tau}),\ x\in H.
    \end{align*}
\end{defi}

\begin{prop}[Proposition 2.4 of \cite{Manca2009}]\label{prop: Manca-core2}
Assume that Hypotheses \ref{Hyp:base} hold true. For every $\delta,\tau>0$, let $\lambda>w_0$, $\Psi\in \C_{1}(H)$ and let $\Phi:H\rightarrow\R$ be the function defined by
\[
\Phi(x):=\int_0^{+\infty}e^{-\lambda t}P_{\delta,\tau}(t)  \Psi(x)\mathd t.
\]
Then $\Phi$ belongs to ${\rm Dom}(N_{\delta,\tau})$ and verifies
\[
\lambda\Phi-N_{\delta,\tau}\Phi=\Psi.
\]
\end{prop}

\begin{prop}[Theorem 4.1 of \cite{Manca2009}]\label{prop: Manca-core3}
Assume that Hypotheses \ref{Hyp:base} hold true. For every $\delta,\tau>0$ we have ${\rm Dom}(L_1)\cap C^1_b(H)={\rm Dom}(N_{\delta,\tau})\cap C^1_b(H)$ and 
\[
N_{\delta,\tau}\Phi=L_1\Phi+\scal{F_{\delta,\tau}+G}{\nabla\Phi}_H,\quad \Phi\in {\rm Dom}(L_1)\cap C^1_b(H).
\]
\end{prop}

\subsection{Proof of Theorem \ref{thm: coreN}}

\begin{prop}\label{semigruppo-proprieta}
Assume that Hypotheses \ref{Hyp:base} hold true.  Let $q\geq 1$ be the constant given by Hypotheses \ref{Hyp:base} and let $\lambda>w_q$ where $w_q$ is the constant appearing in \eqref{Bound-semi}. Let $\Psi\in \C_{q}(E)$ and let $\Phi:E\rightarrow\R$ be the function defined by
\begin{equation}\label{Res-semi}
\Phi(x):=\int_0^{+\infty}e^{-\lambda t}P(t)\Psi(x)\mathd t.
\end{equation}
Then $\Phi$ belongs to ${\rm Dom}(N_q)$ and verifies
\[
\lambda\Phi-N_q\Phi=\Psi.
\]
\end{prop}
\begin{proof}
    Fix $\Psi\in \C_{q}(E)$, $\lambda>w_q$ and $\Phi$ defined by \eqref{Res-semi}. For every $h\geq 0$ and $x\in E$, we have
    \begin{align*}
        P(h)\Phi(x)-\Phi=(e^{\lambda h}-1)\int_h^{\infty}e^{-\lambda t}P(t)\Psi(x)\mathd t-\int_0^he^{-\lambda t}P(t)\Psi(x)\mathd t.
    \end{align*}
    By the integral mean value theorem, we obtain
    \[
    \plim{h\rightarrow 0}\dfrac{P(h)\Phi-\Phi}{h}=\lambda \Phi-\Psi,\qquad \text{ in } \C_q(E),
    \]
    and so that, by Definition \ref{def:pi-gen-N-p}, we obtain $\Phi\in {\rm Dom}(N_q)$ and
    \[
    N_q\Phi(x)=\lambda\Phi(x)-\Psi(x).
    \]
\end{proof}

\begin{lemm}
Assume that Hypotheses \ref{Hyp:base} hold true. For every $\Phi\in\xi_A(H)$, $t\geq 0$ and $x\in E$ we have
\begin{align}\label{sviluppo-taylor}
P(t)\Phi(x)=\E[\Phi(X(t,x))]&=\Phi(x)+\E\left[\int^t_0N_0\Phi(X(s,x))\mathd s\right]\notag\\
&=\Phi(x)+\int^t_0 P(s)\left[N_0\Phi\right](x)\mathd s.
\end{align}
\end{lemm}

\begin{proof}
    Fix $\Phi\in\xi_A(H)$ and $x\in E$. Let $\{X_\delta(t,x)\}_{t\geq 0}$ be the unique mild solution to \eqref{SPDE-delta} and let $\{P_\delta(t)\}_{t\geq 0}$ be the semigroup defined by \eqref{app-semi-delta}. By \cite[Theorem 3.19]{Dap2004}, for every $\delta>0$, $t\geq 0$ we have 
\begin{align*}
P_\delta(t)\Phi(x)=\E[\Phi(X_\delta(t,x))]&=\Phi(x)+\E\left[\int^t_0N_{0,\delta}\Phi(X_\delta(s,x))\mathd s\right]\notag\\
&=\Phi(x)+\int^t_0 P_\delta(s)\left[N_{0,\delta}\Phi\right](x)\mathd s,
\end{align*}
where
{\small
\[
N_{0,\delta}\Phi(X_\delta(t,x))=\frac{1}{2}{\rm Tr}\left[RR^*\nabla^2\Phi(X_\delta(t,x))\right]+\scal{X_\delta(t,x)}{A^*\nabla\Phi(X_\delta(t,x))}+\scal{F_\delta(X_\delta(t,x))+G(X_\delta(t,x))}{\nabla\Phi(X_\delta(t,x))}.
\]
}
We let $\delta\rightarrow0$. By \eqref{yosida-spde5}, $X_\delta(s,x)\rightarrow X(s,x)$ in probability in $H$ for every $s\in[0,t]$. Let $\Gamma$ be the supremum over $s\in[0,t]$ of the right-hand side of \eqref{yosida-spde2}; by \eqref{yosida-spde2} and \eqref{solmild1} we have $\norm{X_\delta(s,x)}_E+\norm{X(s,x)}_E\leq 2\Gamma$ for every $\delta>0$ and $s\in[0,t]$, and $\Gamma$ has finite moments of every order by \eqref{Bound-convu}. The terms of $N_{0,\delta}\Phi(X_\delta(s,x))$ which do not involve $F_\delta$ converge in probability to the corresponding terms of $N_0\Phi(X(s,x))$, since they are continuous functions on $H$. As far as the remaining term is concerned, we recall that $F_\delta=F\circ J_\delta$ and, by \eqref{eq:yosida} and \eqref{yosida2}, $\norm{J_\delta(y)-y}_H=\delta\norm{F_\delta(y)}_H\leq\delta\norm{F(y)}_H$ for every $y\in E$; moreover, since $F_E$ is $m$-dissipative in $E$, the restriction of $J_\delta$ to $E$ is the resolvent of $F_E$, so that there exists $C_0>0$ with $\norm{J_\delta(y)}_E\leq\norm{y}_E+C_0$ for every $y\in E$ and $\delta\in(0,1]$. Since $\mathcal D\subseteq E$, we have $C'_\Phi:=\sup_{y\in H}\norm{\nabla\Phi(y)}_E<+\infty$, and $\nabla\Phi:H\rightarrow H$ is Lipschitz continuous. Hence, on the event $\{2\Gamma+C_0\leq B\}$, by \eqref{FdualEH} we get
\[
\abs{\scal{F_\delta(X_\delta(s,x))-F(X(s,x))}{\nabla\Phi(X_\delta(s,x))}_H}\leq L_BC'_\Phi\Big(\delta\norm{F(X_\delta(s,x))}_H+\norm{X_\delta(s,x)-X(s,x)}_H\Big),
\]
so that $\scal{F_\delta(X_\delta(s,x))}{\nabla\Phi(X_\delta(s,x))}_H\rightarrow\scal{F(X(s,x))}{\nabla\Phi(X(s,x))}_H$ in probability. Finally, by \eqref{yosida2} and Hypotheses \ref{Hyp:base}(6),
\[
\abs{N_{0,\delta}\Phi(X_\delta(s,x))}\leq C_\Phi(1+\Gamma^q),\qquad \delta\in(0,1],\ s\in[0,t],
\]
with $\E[\Gamma^q]<+\infty$. Hence, letting $\delta\rightarrow0$, the Vitali convergence theorem, the dominated convergence theorem in $s$ and \eqref{convergenza-semi-delta} yield \eqref{sviluppo-taylor}.
\end{proof}

\begin{lemm}\label{lemma:passo2}
    Assume that Hypotheses \ref{Hyp:base} hold true. Let $q\geq 1$ be the constant in Hypotheses \ref{Hyp:base}. If $\Phi\in {\rm Dom}(L_1)\cap C_b^1(H)$ then $\Phi\in {\rm Dom}(N_q)$ and
    \begin{equation}\label{p-conv}
        \plim{t\rightarrow 0}\dfrac{P(t)\Phi-\Phi}{t}= L_1\Phi+\scal{F+G}{\nabla\Phi}_H,\qquad  \text{ in } \C_q(E)
    \end{equation}
    In particular, $\xi_A(H)\subseteq {\rm Dom}(N_q)$ and $N_q$ is an extension in $\C_{q}(E)$ of the operator $N_0$.
\end{lemm}
\begin{proof}
Fix $\Phi\in {\rm Dom}(L_1)\cap C_b^1(H)$.  We have to prove that
   \begin{equation}\label{Tesi}
        \plim{t\rightarrow 0}\dfrac{P(t)\Phi-\Phi}{t}= L_1\Phi+\scal{F+G}{\nabla\Phi}_H,\qquad \text{ in } \C_q(E).
   \end{equation}
By Proposition \ref{prop: Manca-core} there exist $k\in\N$ and a family $\{\Phi_{n_1,\ldots,n_k}\}_{n_1,\ldots,n_k\in\N}\subseteq\xi_A(H)$ such that \eqref{Mancacore1}, \eqref{Mancacore2} and \eqref{Mancacore3} hold. 
By \eqref{Op-OUP}, \eqref{imme!}, \eqref{Mancacore1}, \eqref{Mancacore2} and Remark \ref{rmk:Manca-hyp} with $\mathcal T=F+G$, we have
\begin{align}
\plim{n_1,\ldots n_k\rightarrow +\infty}N_0\Phi_{n_1,\ldots,n_k}&=\plim{n_1,\ldots n_k\rightarrow +\infty}L_0\Phi_{n_1,\ldots,n_k}+\scal{F+G}{\nabla\Phi_{n_1,\ldots,n_k}}_H\notag\\
&=L_1\Phi+\scal{F+G}{\nabla\Phi}_H, \qquad \text{ in } \C_q(E).\label{S1}
\end{align}
By \eqref{sviluppo-taylor}, for every $x\in E$ we have
\begin{equation*}
P(t)\Phi_{n_1,\ldots,n_k}(x)-\Phi_{n_1,\ldots,n_k}(x)=\int^t_0P(s)\left[N_0\Phi_{n_1,\ldots,n_k}\right](x)\mathd s
\end{equation*}
hence by Hypotheses \ref{Hyp:base}, Proposition \ref{prop: Manca-core}, \eqref{S1} and the dominated convergence theorem we have
\begin{align}
\dfrac{P(t)\Phi(x)-\Phi(x)}{t}&=\lim_{n_1,\ldots n_k\rightarrow +\infty}\dfrac{P(t)\Phi_{n_1,\ldots,n_k}(x)-\Phi_{n_1,\ldots,n_k}(x)}{t}\notag\\
&=\lim_{n_1,\ldots n_k\rightarrow +\infty}\dfrac{\int^t_0P(s)\left[L_0\Phi_{n_1,\ldots,n_k}+\scal{F+G}{\nabla\Phi_{n_1,\ldots,n_k}}_H\right](x)\mathd s}{t},\notag\\
&=\dfrac{\int^t_0P(s)\left[L_1\Phi+\scal{F+G}{\nabla\Phi}_H\right](x)\mathd s}{t}.\label{S2}
\end{align}

By definition $L_1\Phi\in \C_{1}(H)$, so by \eqref{imme!} we have $L_1\Phi\in \C_{q}(E)$. Moreover by Hypotheses \ref{Hyp:base} we have $\scal{F+G}{\nabla\Phi}_H\in \C_{q}(E)$. Hence, by Proposition \ref{semigruppoP}, \eqref{S2} and the mean value theorem, for every $t\in(0,1)$, we obtain
{\small
\begin{align*}
&\norm{\dfrac{P(t)\Phi-\Phi}{t}}_{\C_{q}(E)}\leq \dfrac{\int^t_0\norm{P(s)\left[L_1\Phi+\scal{F+G}{\nabla\Phi}_H\right]}_{\C_{q}(E)}\mathd s}{t}\leq c_qe^{w_q}\norm{L_1\Phi+\scal{F+G}{\nabla\Phi}_H}_{\C_{q}(E)},\\
 &\lim_{t\rightarrow 0}\dfrac{P(t)\Phi(x)-\Phi(x)}{t}= L_1\Phi(x)+\scal{F(x)+G(x)}{\nabla\Phi(x)}_H,\qquad x\in E.
\end{align*}
} and so that \eqref{Tesi} holds.

\end{proof}

\begin{lemm}\label{lemma:passo3}
  Assume that Hypotheses \ref{Hyp:base} hold true. Let $q\geq 1$ be the constant in Hypotheses \ref{Hyp:base}. For every $\Phi\in {\rm Dom}(N_q)$ there exist $k\in\N$ and $\{\Phi_{n_1,\ldots,n_k}\}_{n_1\in\N,\ldots,n_k\in\N}\subseteq  {\rm Dom}(L_1)\cap C_b^1(H)$ such that
\begin{align*}
&\plim{n_1\rightarrow+\infty,\ldots,n_k\rightarrow+\infty}\Phi_{n_1,\ldots,n_k}=\Phi,\qquad \textit{ in } \C_q(E),\\
&\plim{n_1\rightarrow+\infty,\ldots,n_k\rightarrow+\infty}L_1\Phi_{n_1,\ldots,n_k}+\scal{F(\cdot)+G(\cdot)}{\nabla\Phi_{n_1,\ldots,n_k}}_H=N_q\Phi,\qquad  \textit{ in } \C_q(E),
\end{align*}
\end{lemm}
\begin{proof}
Fix $\Phi\in {\rm Dom}(N_q)$. $\Psi=\lambda\Phi-N_q\Phi\in \C_{q}(E)$ and so that, by \eqref{densità}, there exist $m\in\N$ and $\{\Psi_{n_1,\ldots,n_m}\}\subseteq \xi_A(H)$ such that
\begin{equation}\label{conv-fff}
\plim{n_1,\ldots,n_m\rightarrow +\infty}\Psi_{n_1,\ldots,n_m}= \Psi=\lambda\Phi-N_q\Phi,\qquad \text{ in } \C_q(E),
\end{equation}
where $\lambda>\max\{w_0,w_q,\zeta\}$ and $w_q$, $w_0$ and $\zeta$ are the constants appearing in \eqref{Bound-semi}, \eqref{Bound-semi-delta-beta} and \eqref{stima-grad-grad}, respectively. For every $n_1,\ldots n_m\in\N$ and $\delta,\tau>0$, by Proposition \ref{prop: Manca-core2} we have $\Phi_{n_1,\ldots,n_m,\delta,\tau}\in{\rm Dom}(N_{\delta,\tau})$, and, since $\Psi_{n_1,\ldots,n_m}\in C^1_b(H)$ and $F_{\delta,\tau}+G$ is Lipschitz continuous, by \cite[Proposition 4.5]{Manca2009} we have $\Phi_{n_1,\ldots,n_m,\delta,\tau}\in C^1_b(H)$. Hence, by Proposition \ref{prop: Manca-core3}, $\Phi_{n_1,\ldots,n_m,\delta,\tau}\in{\rm Dom}(L_1)\cap C^1_b(H)$ and 
\[
\Phi_{n_1,\ldots,n_m,\delta,\tau}(x):=\int_0^{+\infty}e^{-\lambda t}P_{\delta,\tau}(t)\Psi_{n_1,\ldots,n_m}(x)\mathd t.
\]
 We first note that, by \eqref{convergenza-semi-delta-beta}, \eqref{convergenza-semi-delta}, \eqref{Bound-semi}, Proposition \ref{semigruppo-proprieta}, \eqref{conv-fff} and the dominated convergence theorem we obtain 
\begin{align}
\plim{n_1,\ldots,n_m\rightarrow +\infty}\plim{\delta\rightarrow 0}\plim{\tau\rightarrow 0}\Phi_{n_1,\ldots,n_m,\delta,\tau}&=\plim{n_1,\ldots,n_m\rightarrow +\infty}\int_0^{+\infty}e^{-\lambda t}P(t)\Psi_{n_1,\ldots,n_m}\mathd t\notag\\
&=\int_0^{+\infty}e^{-\lambda t}P(t)\Psi\mathd t=\Phi,\qquad \text{ in } \C_q(E).\label{ccc1}
\end{align}
For every $n_1,\ldots n_m\in\N$ and $\delta,\tau>0$ by Propositions \ref{prop: Manca-core2} and \ref{prop: Manca-core3} we have
\begin{equation*}
    \lambda\Phi_{n_1,\ldots,n_m,\delta,\tau}(x)-L_1\Phi_{n_1,\ldots,n_m,\delta,\tau}(x)=\Psi_{n_1,\ldots,n_m}(x)+\scal{F_{\delta,\tau}(x)+G(x)}{\nabla\Phi_{n_1,\ldots,n_m,\delta,\tau}(x)}_H,\qquad x\in E.
    \end{equation*}
By Lemma \ref{lemma:passo2} we have
\begin{equation*}
    \lambda\Phi_{n_1,\ldots,n_m,\delta,\tau}(x)-N_q\Phi_{n_1,\ldots,n_m,\delta,\tau}(x)=\Psi_{n_1,\ldots,n_m}(x)+\scal{F_{\delta,\tau}(x)-F(x)}{\nabla\Phi_{n_1,\ldots,n_m,\delta,\tau}(x)}_H,\qquad x\in E.
    \end{equation*}
By \eqref{stima-grad-grad}, for every $n_1,\ldots n_m\in\N$ and $\delta,\tau>0$ we have
\[
\sup_{x\in H}\norm{\nabla \Phi_{n_1,\ldots,n_m,\delta,\tau}(x)}_H\leq (\lambda-\zeta)^{-1}\norm{\Psi_{n_1,\ldots,n_m}}_{C^1_b(H)},
\]
hence, taking into account \eqref{imme!}, by \eqref{yosida1}, \eqref{yosida-OU1} and \eqref{conv-fff} we obtain
  \begin{align*}
\plim{n_1,\ldots,n_m\rightarrow +\infty}\lim_{\delta\rightarrow 0}\lim_{\tau\rightarrow 0}\Psi_{n_1,\ldots,n_m}(x)&+\scal{F_{\delta,\tau}(x)-F(x)}{\nabla\Phi_{n_1,\ldots,n_m,\delta,\tau}(x)}_H\\
&=\plim{n_1,\ldots,n_m\rightarrow +\infty}\Psi_{n_1,\ldots,n_m}=\Psi,\qquad \textit{ in } \C_q(E)
    \end{align*}
    and so that, recalling that $\Psi=\lambda\Phi-N_q\Phi$, we obtain
    \begin{equation}\label{ccc2}
\plim{n_1,\ldots,n_m\rightarrow +\infty}\lim_{\delta\rightarrow 0}\lim_{\tau\rightarrow 0}\lambda\Phi_{n_1,\ldots,n_m,\delta,\tau}-N_q\Phi_{n_1,\ldots,n_m,\delta,\tau}= \Psi=\lambda\Phi-N_q\Phi,\qquad \text{ in } \C_q(E).
    \end{equation}
Hence, by \eqref{ccc1} and \eqref{ccc2}, the statement holds true.
\end{proof}

Now we can prove the main result of this section.

\begin{proof}[Proof of Theorem \ref{thm: coreN}]
Let $\Phi\in {\rm Dom}(N_q)$. By Lemma \ref{lemma:passo3} there exist $v\in\N$ and $\{\Phi_{n_1,\ldots,n_v}\}_{n_1\in\N,\ldots,n_v\in\N}\subseteq  {\rm Dom}(L_1)\cap C_b^1(H)$ such that
\begin{align}
&\plim{n_1\rightarrow+\infty,\ldots,n_v\rightarrow+\infty}\Phi_{n_1,\ldots,n_v}=\Phi,\qquad\text{ in } \C_q(E),\label{C1}\\
&\plim{n_1\rightarrow+\infty,\ldots,n_v\rightarrow+\infty}L_1\Phi_{n_1,\ldots,n_v}+\scal{F(\cdot)+G(\cdot)}{\nabla\Phi_{n_1,\ldots,n_v}}_H=N_q\Phi,\qquad\text{ in } \C_q(E)\label{C2}.
\end{align}
Fixed $n_1,\ldots,n_v\in\N$, by Proposition \ref{prop: Manca-core} there exist $w\in\N$ and\\
$\{\Phi_{n_1,\ldots,n_v,n_{v+1},\ldots,n_{v+w}}\}_{n_{v+1},\ldots,n_{v+w}\in\N}\subseteq  \xi_A(H)$ such that
\begin{align}
&\plim{n_{v+1}\rightarrow+\infty,\ldots,n_{v+w}\rightarrow+\infty}\Phi_{n_1,\ldots,n_v,n_{v+1},\ldots,n_{v+w}}=\Phi_{n_1,\ldots,n_v},\qquad\text{ in } \C_1(H)\label{C3}\\
&\plim{n_{v+1}\rightarrow+\infty,\ldots,n_{v+w}\rightarrow+\infty}\scal{\nabla\Phi_{n_1,\ldots,n_v,n_{v+1},\ldots,n_{v+w}}}{h}_H=\scal{\nabla\Phi_{n_1,\ldots,n_v}}{h}_H,\qquad\text{ in } \C_1(H)\quad h\in H,\label{C4}\\
&\plim{n_{v+1}\rightarrow+\infty,\ldots,n_{v+w}\rightarrow+\infty}L_0\Phi_{n_1,\ldots,n_v,n_{v+1},\ldots,n_{v+w}}=L_1\Phi_{n_1,\ldots,n_v},\qquad\text{ in } \C_1(H).\label{C5}
\end{align}
Set $k=v+w$. Taking into account \eqref{imme!}, by \eqref{C1} and $\eqref{C3}$ we obtain \eqref{coreN1}, further by \eqref{C2}, \eqref{C4} and \eqref{C5} we deduce \eqref{coreN2}.
\end{proof}

\section{Relation between invariant measures and infinitesimal invariant measures}\label{sect: inv-infinve-abstract}

In this section, exploiting the $\pi$-core of $N_0$ established in Theorem \ref{thm: coreN}, we investigate the relationship between invariant measures and infinitesimally invariant measures for the transition semigroup $\{P(t)\}_{t\geq 0}$ defined in \eqref{semigruppo-transizione}. First of all, we recall a standard result concerning the behavior of the semigroup $\{P(t)\}_{t\geq 0}$ on $L^p$ spaces with respect to an invariant measure.

\begin{prop}\label{prop:extensionL2}
    Assume that Hypotheses \ref{Hyp:base} hold true and let $\mu$ be a Borel probability measure on $H$. If $\mu$ is an invariant measure for $\{P(t)\}_{t\geq 0}$, namely for every $\Phi\in C_b(H)$ we have
        \begin{equation}\label{inv}
        \int_HP(t)\Phi(x)\mu(\mathd x)=\int_H\Phi(x)\mu(\mathd x),
        \end{equation}
then for every $t\geq 0$ and $p\geq 1$ the operator $P(t)$ defined by \eqref{semigruppo-transizione} extends uniquely to a linear contraction of $L^p(H,\mu)$, still denoted by $P(t)$, and $\{P(t)\}_{t\geq 0}$ is a strongly continuous contraction semigroup on $L^p(H,\mu)$. If, in addition, $\int_H\norm{x}^{pq}_E\mu(\mathd x)<+\infty$, where $q$ is the constant given by Hypotheses \ref{Hyp:base}, then the infinitesimal generator $N_p$ of $\{P(t)\}_{t\geq0}$ in $L^p(H;\mu)$ is the closure of $N_0:\xi_A(H)\subseteq L^p(H;\mu)\rightarrow L^p(H;\mu)$.\\
If, in addition,
        \begin{equation}\label{eq:symmetryCb}
        \int_H\Psi(x)P(t)\Phi(x)\mu(\mathd x)=\int_H\Phi(x)P(t)\Psi(x)\mu(\mathd x),\qquad \Phi,\Psi\in C_b(H),\ t\geq 0,
        \end{equation}
        then $P(t)$ is self-adjoint in $L^2(H,\mu)$ for every $t\geq 0$. 
\end{prop}
\begin{proof}
It is sufficient to note that, for every $p\geq 1$, $t\geq 0$ and $\Phi\in C_b(H)$ by \eqref{inv}, we obtain
\[
\int_H\abs{P(t)\Phi(x)}^p\mu(\mathd x)\leq \int_H P(t)(\abs{\Phi}^p)(x)\mu(\mathd x)=\int_H\abs{\Phi(x)}^p\mu(\mathd x),
\]
and so that
\[
\norm{P(t)\Phi}_{L^p(H;\mu)}\leq \norm{\Phi}_{L^p(H;\mu)}.
\]
Hence, since $C_b(H)$ is dense in $L^p(H;\mu)$ for every $p\geq 1$ and $\xi_A(H)$ is $\pi$-dense in $C_b(H)$ (see Proposition \ref{densità-xia}), taking into account Theorem \ref{thm: coreN} the proof is complete.
\end{proof}

Now we can prove the main result of this section.

\begin{thm}\label{thm:Inv-inf}
    Assume that Hypotheses \ref{Hyp:base} hold true. Let $\mu$ be a Borel probability measure such that
    \begin{equation}\label{cond-mom}
       \int_H\norm{x}^q_E\mu(\mathd x)<+\infty, 
    \end{equation}
    where $q$ is the constant given by Hypotheses \ref{Hyp:base}. The following propositions are equivalent.
    \begin{enumerate}
        \item\textbf{(Invariant)} $\mu$ is an invariant measure for $\{P(t)\}_{t\geq 0}$.
        \item \textbf{(Infinitesimal invariant)} $\mu$ is an infinitesimal invariant measure for $\{P(t)\}_{t\geq 0}$, namely for every $\Phi\in \xi_A(H)$ we have
        \begin{equation}\label{inf-inv}
        \int_HN_0\Phi(x)\mu(\mathd x)=0
        \end{equation}
    \end{enumerate}
\end{thm}

\begin{proof} $ $\\
    \textbf{1. $\Rightarrow$ 2.}\\
    Since $\xi_A(H)\subseteq C_b(H)$, by \eqref{densità} and the dominated convergence theorem, formula \eqref{inv} holds for every $\Phi\in \mathscr{C}_p(E)$. Moreover, we recall that, by the definitions  \eqref{Op-OUP} and \eqref{xia} of $N_0$ and $\xi_A(H)$, for every $\Phi\in\xi_A(H)$ we have $N_0\Phi\in \mathscr{C}_p(E)$. Hence, by \eqref{sviluppo-taylor}, for every $\Phi\in\xi_A(H)$ we have
    \begin{align*}
        0&=\int_H\left(P(t)\Phi(x)-\Phi(x)\right)\mu(\mathd x)=\int_H\left(\int_0^tP(s)\left[N_0\Phi\right](x)\mathd s\right)\mu(\mathd x)\\
        &=\int_0^t\left(\int_HP(s)\left[N_0\Phi\right](x)\mu(\mathd x)\right)\mathd s=\int_0^t\left(\int_HN_0\Phi(x)\mu(\mathd x)\right)\mathd s=t\int_HN_0\Phi(x)\mu(\mathd x),\\
    \end{align*}
and so that \eqref{inf-inv} holds.\\
\textbf{2. $\Rightarrow$ 1.}\\
First of all, we note that, by Theorem \ref{thm: coreN} and \eqref{inf-inv} we have
\begin{equation}\label{inf-inv-p}
        \int_HN_q\Phi(x)\mu(\mathd x)=0,\qquad \Phi\in{\rm Dom}(N_q),
\end{equation}
where $N_q$ is the $\pi$-generator of $\{P(t)\}_{t\geq 0}$ given by Definition \ref{def:pi-gen-N-p}. Fix $\Phi\in \xi_A(H)$ and $t>0$. Since $\Phi_{|E}\in \mathscr{C}_p(E)$, by Proposition \ref{prop:commutazione}, \eqref{sviluppo-taylor} and \eqref{p-conv} we deduce
\begin{align*}
        \int_H\left(P(t)\Phi(x)-\Phi(x)\right)\mu(\mathd x)&=\int_H\left(\int_0^tP(s)\left[N_0\Phi\right](x)\mathd s\right)\mu(\mathd x)\\
        &=\int_H\left(\int_0^tN_q\left[P(s)\Phi\right](x)\mathd s\right)\mu(\mathd x)\\
        &=\int_0^t\left(\int_HN_q\left[P(s)\Phi\right](x)\mu(\mathd x)\right)\mathd s=0,
    \end{align*}
    and so that, since $\xi_A(H)$ is $\pi$-dense in $C_b(H)$ (see Proposition \ref{densità-xia}), we obtain \eqref{inv}.
\end{proof}

\begin{remark}$ $
    \begin{enumerate}
        \item Condition \eqref{cond-mom} guarantees that $\mu(E)=1$ and \eqref{inf-inv} makes sense, since the function $F$  appearing in the first order part of $N_0$ (see \eqref{Op-OUP}) has a polynomial growth.

        \item Theorem \ref{thm: coreN} is necessary only in the proof of \textbf{2. $\Rightarrow$ 1.} to guarantee \eqref{inf-inv-p} holds.
    \end{enumerate}
\end{remark}

We now investigate the symmetry properties of the transition semigroup.

\begin{lemma}\label{lemma:sim-inv}
    Assume that Hypotheses \ref{Hyp:base} hold true. Let $\mu$ be a Borel probability measure such that
    \begin{equation*}
       \int_H\norm{x}^q_E\mu(\mathd x)<+\infty, 
    \end{equation*}
    where $q$ is the constant given by Hypotheses \ref{Hyp:base}. If $N_0$ is symmetric with respect to $\mu$, namely
        \begin{equation}\label{sim-N_0}
            \int_H\Psi(x)N_0\Phi(x)\mu(\mathd x)=\int_H\Phi(x)N_0\Psi(x)\mu(\mathd x),\qquad \Phi,\Psi\in \xi_A(H),
        \end{equation}
        then $\mu$ is an invariant measure of $\{P(t)\}_{t\geq 0}$. 
\end{lemma}
\begin{proof}
    Choosing $\Psi\equiv 1 $ in \eqref{sim-N_0} we obtain \eqref{inf-inv}, and so that, by Theorem \ref{thm:Inv-inf}, also \eqref{inv} holds.
\end{proof}

\begin{thm}\label{thm: rev-infrev}
    Assume that Hypotheses \ref{Hyp:base} hold true. Let $\mu$ be a Borel probability measure such that
    \begin{equation*}
       \int_H\norm{x}^{2q}_E\mu(\mathd x)<+\infty, 
    \end{equation*}
    where $q$ is the constant given by Hypotheses \ref{Hyp:base}. The following propositions are equivalent
    \begin{enumerate}
        \item $\{P(t)\}_{t\geq 0}$ is symmetric with respect to $\mu$, namely \eqref{eq:symmetryCb} holds.
        
        \item $N_0$ is symmetric with respect to $\mu$.
    \end{enumerate}
\end{thm}
\begin{proof} $ $\\
    \textbf{1. $\Rightarrow$ 2.}\\
Fix $ \Phi,\Psi\in \xi_A(H)$. By Definition \ref{def:pi-gen-N-p} and Theorem \ref{thm: coreN} we deduce
\begin{align*}
    \int_H\Psi(x)N_0\Phi(x)\mu(\mathd x)&=\int_H\Psi(x)\left(\lim_{t\rightarrow 0}\dfrac{P(t)\Phi(x)-\Phi(x)}{t}\right)\mu(\mathd x)\\
    &=\lim_{t\rightarrow 0}t^{-1}\left[\int_H\Psi(x)P(t)\Phi(x)\mu(\mathd x)-\int_H\Psi(x)\Phi(x)\mu(\mathd x)\right],
\end{align*}
     hence, applying \eqref{eq:symmetryCb} and taking into account Definition \ref{def:pi-gen-N-p}, we conclude
    \begin{align*}
    \int_H\Psi(x)N_0\Phi(x)\mu(\mathd x)&=\lim_{t\rightarrow 0}t^{-1}\left[\int_H\Phi(x)P(t)\Psi(x)\mu(\mathd x)-\int_H\Psi(x)\Phi(x)\mu(\mathd x)\right]\\
    &=\int_H\Phi(x)\left(\lim_{t\rightarrow 0}\dfrac{P(t)\Psi(x)-\Psi(x)}{t}\right)\mu(\mathd x)=\int_H\Phi(x)N_0\Psi(x)\mu(\mathd x).
\end{align*} 
\textbf{2. $\Rightarrow$ 1.}\\
By Lemma \ref{lemma:sim-inv} the measure $\mu$ is invariant for the semigroup $\{P(t)\}_{t\geq 0}$. By Proposition \ref{prop:extensionL2},  the semigroup $\{P(t)\}_{t\geq 0}$ is uniquely extendable to a strongly continuous semigroup in $L^2(H,\mu)$, still denoted by $\{P(t)\}_{t\geq 0}$ and we denote its infinitesimal generator by $N$. Moreover, by Proposition \ref{prop:extensionL2} with $p=2$, the operator $N$ is the closure of $N_0$ on $\xi_A(H)$ in $L^2(H;\mu)$; hence, by \eqref{sim-N_0}, $N$ is symmetric, and since it generates a contraction semigroup it is self-adjoint.
\end{proof}

Finally, we investigate the relationship between the integration-by-parts formulas and the symmetry of $N_0$.

\begin{thm}\label{thm: rev-per-parti}
    Assume that Hypotheses \ref{Hyp:base} hold true and that $R$ is self-adjoint. Let $\mu$ be a Borel probability measure such that
    \begin{equation*}
       \int_H\norm{x}^q_E\mu(\mathd x)<+\infty, 
    \end{equation*}
    where $q$ is the constant given by Hypotheses \ref{Hyp:base}. The following propositions are equivalent
    \begin{enumerate}
        \item $N_0$ is symmetric with respect to $\mu$.
        \item $\mu$ verifies, for every $h\in {\rm Dom}(A^*)$ and $\Phi\in\xi_A(H)$, the following integration by part formula 
        \begin{equation}\label{int-by-part}
    \int_H\scal{R\nabla\Phi(x)}{Rh}\mu(\mathd x)=-2\int_H\left[\scal{x}{A^*h}+\scal{F(x)+G(x)}{h}\right]\Phi(x)\mu(\mathd x).
        \end{equation}

    \end{enumerate}
    Moreover, in 2. it is enough to require \eqref{int-by-part} for every $h\in\mathcal D$ and $\Phi\in\xi_A(H)$.
\end{thm}
\begin{proof} $ $\\
    \textbf{1. $\Rightarrow$ 2.}\\
First of all, we note that  for every $\Phi,\Psi\in\xi_A(H)$ we have
\[
N_0(\Phi\Psi)=\Phi N_0\Psi+\Psi N_0\Phi+\scal{R\nabla\Phi}{R\nabla\Psi}_H,
\]
integrating both sides on $H$ with respect the measure $\mu$, by Lemma \ref{lemma:sim-inv} we obtain
\[
0=\int_HN_0(\Phi\Psi)(x)\mu(\mathd x)=\int_H\left[\Phi(x)N_0\Psi(x)+\Psi(x)N_0\Phi(x)+\scal{R\nabla\Phi(x)}{R\nabla\Psi(x)}_H\right]\mu(\mathd x),
\]
using \eqref{sim-N_0} in the previous formula we conclude
\begin{equation}\label{quadratica}
\int_H\Phi(x)N_0\Psi(x)\mu(\mathd x)=-\frac{1}{2}\int_H\scal{R\nabla\Phi(x)}{R\nabla\Psi(x)}_H\mu(\mathd x).
\end{equation}
We can extend \eqref{quadratica} to every $\Psi\in{\rm Dom}(L_1)\cap C^1_b(H)$, with $N_0\Psi$ replaced by $N_q\Psi=L_1\Psi+\scal{F+G}{\nabla\Psi}_H$ (see Lemma \ref{lemma:passo2}). Indeed, if $\{\Psi_{n_1,\ldots,n_k}\}\subseteq\xi_A(H)$ is given by Proposition \ref{prop: Manca-core}, then, by \eqref{Mancacore1}, \eqref{Mancacore2} and Remark \ref{rmk:Manca-hyp} with $\mathcal T=F+G$ and with $\mathcal T=R^*R\nabla\Phi$, $N_0\Psi_{n_1,\ldots,n_k}$ $\pi$-converges to $N_q\Psi$ in $\C_q(E)$ and $\scal{R\nabla\Phi}{R\nabla\Psi_{n_1,\ldots,n_k}}_H$ $\pi$-converges to $\scal{R\nabla\Phi}{R\nabla\Psi}_H$ in $\C_q(E)$, and it suffices to apply iteratively the dominated convergence theorem, since $\Phi$ is bounded and $\int_H\norm{x}^q_E\mu(\mathd x)<+\infty$. Now fix $h\in{\rm Dom}(A^*)$ and, for $\varepsilon\in(0,1]$, set $\Psi_\varepsilon(x):=\varepsilon^{-1}\sin(\varepsilon\scal{x}{h}_H)$. By \cite[Proposition 4.3]{Manca2009}, $\Psi_\varepsilon\in{\rm Dom}(L_1)\cap C^1_b(H)$, with $\nabla\Psi_\varepsilon(x)=\cos(\varepsilon\scal{x}{h}_H)h$ and
\[
N_q\Psi_\varepsilon(x)=-\frac\varepsilon2\norm{R^*h}^2_H\sin(\varepsilon\scal{x}{h}_H)+\cos(\varepsilon\scal{x}{h}_H)\Big[\scal{x}{A^*h}_H+\scal{F(x)+G(x)}{h}_H\Big].
\]
Since $\abs{N_q\Psi_\varepsilon(x)}\leq C_h(1+\norm{x}^q_E)$ uniformly with respect to $\varepsilon\in(0,1]$, applying \eqref{quadratica} to $\Psi_\varepsilon$ and letting $\varepsilon\rightarrow0$ by the dominated convergence theorem, we obtain \eqref{int-by-part}, recalling that $R$ is self-adjoint.
 $ $\\
    \textbf{2. $\Rightarrow$ 1.}\\
We extend $N_0$ and both sides of \eqref{int-by-part} to test functions with complex values, and so $\scal{\cdot}{\cdot}$ denotes the complex-bilinear extension of the scalar product of $H$. Since $\mathcal D$ is a linear subspace, the complex linear span of $\{e^{i\scal{\cdot}{h}}:\ h\in\mathcal D\}$ is an algebra containing $\xi_A(H)$ (indeed $\cos\scal{x}{h}=\frac12(e^{i\scal{x}{h}}+e^{i\scal{x}{-h}})$ and $\sin\scal{x}{h}=\frac1{2i}(e^{i\scal{x}{h}}-e^{i\scal{x}{-h}})$), and both sides of \eqref{sim-N_0} are bilinear in $(\Phi,\Psi)$; hence it is enough to prove \eqref{sim-N_0} for $\Phi(x)=e^{i\scal{x}{h}_H}$ and $\Psi(x)=e^{i\scal{x}{k}_H}$ with $h,k\in\mathcal D$. Applying \eqref{int-by-part} to the product $\Phi\Psi$ and multiplying both sides by $i$, for every $h\in\mathcal D$ we obtain
\begin{align}
        i\int_H\Psi(x)\scal{R\nabla\Phi(x)}{Rh}\mu(\mathd x)=&-i\int_H\Phi(x)\scal{R\nabla\Psi(x)}{Rh}\mu(\mathd x)\notag\\
        &-2i\int_H\left[\scal{x}{A^*h}+\scal{F(x)+G(x)}{h}\right]\Phi(x)\Psi(x)\mu(\mathd x).\label{int-by-partV}
        \end{align}

 Let $\Phi(x):=e^{i\scal{x}{h}_H}$ and $\Psi(x):=e^{i\scal{x}{k}_H}$ with with $h,k\in\mathcal D$. By \eqref{Op-OUP} we obtain
\begin{align}\label{N0-xiA}
N_0\Phi(x)=\left[-\frac{1}{2}\norm{Rh}_H^2+i\scal{x}{A^*h}_H+i\scal{G(x)+F(x)}{h}\right]\Phi(x)
\end{align}
and by \eqref{int-by-partV} 
\begin{align}
-\int_H\norm{Rh}^2_H\Psi(x)\Phi(x)\mu(\mathd x)=&-\int_H\scal{R\nabla\Psi(x)}{R\nabla\Phi(x)}_H\mu(\mathd x)\notag\\
&-2i\int_H\left[\scal{x}{A^*h}+\scal{F(x)+G(x)}{h}\right]\Phi(x)\Psi(x)\mu(\mathd x),\label{int-by-part-xia}
\end{align}
finally, combining \eqref{N0-xiA} and \eqref{int-by-part-xia}, we conclude 
\begin{equation}
    \int_H\Psi(x)N_0\Phi(x)\mu(\mathd x)=-\frac{1}{2}\int_H\scal{R\nabla\Phi(x)}{R\nabla\Psi(x)}\mu(\mathd x),
\end{equation}
and so that \eqref{sim-N_0} holds.
\end{proof}

\subsection{Application to the gradient system}

Throughout this subsection we assume that Hypotheses \ref{Hyp:base} hold with $R=\alpha\Id_H$ for some $\alpha>0$ and that $A$ is self-adjoint. By Hypotheses \ref{Hyp:base}(5) the operator $A$ is strictly negative, so that the transition semigroup $\{T(t)\}_{t\geq 0}$ associated with \eqref{SPDE} with $F=G=0$ has the unique invariant measure $\gamma=\mathcal{N}(0,Q)$, where $Q:=\frac{\alpha^2}{2}(-A)^{-1}$ (see \cite{Dap-Zab1996}). By \cite[Proposition 9.2.4]{Dap-Zab2002} the operator
\[
\nabla:C_b^1(H)\subseteq L^2(H,\gamma)\rightarrow L^2(H,\gamma;H)
\]
is closable, and we denote by $W^{1,2}(H;\gamma)$ the domain of its closure.

\begin{thm}\label{thm:dap-grad-sist}
Assume that there exists $U\in W^{1,2}(H;\gamma)$, bounded from below, such that
\begin{equation}\label{eq:gradient-drift}
F+G=-\frac{\alpha^2}{2}\nabla U\qquad\gamma\text{-a.e.},
\end{equation}
and define the Borel probability measure
\begin{equation}\label{Gibbs}
\nu(\mathd x):=\Big[\int_He^{-U}\mathd\gamma\Big]^{-1}e^{-U(x)}\gamma(\mathd x).
\end{equation}
Then for every $\Phi\in\xi_A(H)$ and $h\in{\rm Dom}(A)$ the integration by parts formula \eqref{int-by-part} holds with $\mu=\nu$. If in addition $\int_H\norm{x}^q_E\,\nu(\mathd x)<+\infty$, then $N_0$ is symmetric with respect to $\nu$.
\end{thm}

\begin{proof}
Let $\Phi\in\xi_A(H)$ and $h\in{\rm Dom}(A)$. Since $U$ is bounded from below, $e^{-U}\in W^{1,2}(H;\gamma)\cap L^\infty(H,\gamma)$ with $\nabla e^{-U}=-e^{-U}\nabla U$, hence $\Phi e^{-U}\in W^{1,2}(H;\gamma)$. Moreover $\alpha^2h=Q\big(2(-A)h\big)\in Q(H)$ and $Q^{-1}(\alpha^2h)=-2Ah$, so that the Gaussian integration by parts formula (see e.g. \cite{Bogachev98}) gives
\[
\int_H\scal{\nabla(\Phi e^{-U})(x)}{\alpha^2h}_H\,\gamma(\mathd x)=-2\int_H\scal{x}{Ah}_H\Phi(x)e^{-U(x)}\,\gamma(\mathd x).
\]
Since $\scal{\nabla(\Phi e^{-U})}{\alpha^2h}_H=\big(\scal{R\nabla\Phi}{Rh}_H-\Phi\scal{\alpha^2\nabla U}{h}_H\big)e^{-U}$ and, by \eqref{eq:gradient-drift}, $\alpha^2\nabla U=-2(F+G)$, dividing by $\int_He^{-U}\mathd\gamma$ we obtain
\[
\int_H\scal{R\nabla\Phi(x)}{Rh}_H\,\nu(\mathd x)=-2\int_H\Big[\scal{x}{Ah}_H+\scal{F(x)+G(x)}{h}_H\Big]\Phi(x)\,\nu(\mathd x),
\]
which is \eqref{int-by-part} since $A^*=A$. The symmetry of $N_0$ follows from Theorem \ref{thm: rev-per-parti}, whose assumptions are satisfied since $R=\alpha\Id_H$ is self-adjoint.
\end{proof}

\begin{cor}\label{cor:gradient-invariant}
Under the assumptions of Theorem \ref{thm:dap-grad-sist}, assume moreover that $\int_H\norm{x}^q_E\,\nu(\mathd x)<+\infty$. Then $\nu$ is the unique invariant measure of $\{P(t)\}_{t\geq0}$.
\end{cor}
\begin{proof}
By Theorem \ref{thm:dap-grad-sist} and Lemma \ref{lemma:sim-inv}, $\nu$ is invariant. Since $R=\alpha\Id_H$, Hypothesis \ref{hyp:irr} holds by Remark \ref{rmk:irr-alphaId}, so that uniqueness follows from Corollary \ref{cor:uniqueness-invariant}.
\end{proof}

\section*{Declaration of generative AI use}
\addcontentsline{toc}{section}{Declaration of generative AI use}
We began work on this paper more than two years ago and developed its overall structure and most of its results before first using AI in August 2026. At that stage, we used an AI tool to review a nearly final draft of the manuscript. The tool identified some gaps and suggested ways to improve certain results (in particular the uniqueness proof in Section \ref{section:uniquenessDLR}). Over the following three months, we assessed these suggestions and incorporated them into the manuscript, using AI tools during this final stage to assist us in refining the presentation and improving those results. We reviewed and edited the AI-generated output and take full responsibility for the content of the article.

\section*{Acknowledgments}
\addcontentsline{toc}{section}{Acknowledgments}
The authors D. A. Bignamini and F. C. De Vecchi are members of GNAMPA-INdAM. The authors are very grateful to Carlo Orrieri for the inspiring discussions, which were fundamental to the development of this paper.\\
\noindent Carlos Villanueva Mariz is grateful for funding by the Deutsche Forschungsgemeinschaft (DFG, German Research
Foundation) through EXC 2046, Berlin Mathematical School, and through IRTG 2544 “Stochastic
Analysis in Interaction” - Project ID 410208580.

\bibliographystyle{plain}
\bibliography{Gibbsbibliography}

\begin{thebibliography}{99}

\bibitem{DiNezzaPalatucciValdinoci2012}
E.~Di Nezza, G.~Palatucci, E.~Valdinoci,
\emph{Hitchhiker's guide to the fractional Sobolev spaces},
Bull. Sci. Math. \textbf{136} (2012), no.~5, 521--573.

\bibitem{BBGM2012} Bakry, D., Bolley, F., Gentil, I., and Maheux, P. (2012). Weighted Nash inequalities. Revista Matematica Iberoamericana, 28(3), 879-906.

\bibitem{Big-Def2025}
Bignamini D. A., De Fazio P.
\newblock{\it Stochastic and deterministic reaction-diffusion equations},
eprint {\bf arXiv}:2510.10842v1 (2025).

\bibitem{Big2021} 
Bignamini D.A.
\newblock{\it $L^2$-theory for transitions semigroups associated to dissipative systems}, Stoch. Partial Differ. Equ. Anal. Comput., 2023.

\bibitem{Cer1994}  
Cerrai, S.
\newblock{\small \it A Hille--Yosida theorem for weakly continuous semigroups}, Semigroup Forum {\bf 49} (1994), 349--367.

\bibitem{Cer1995}  
Cerrai, S.
\newblock{\small \it Weakly continuous semigroups in the space of functions with polynomial growth}, Dynam. Syst. Appl. {\bf 4} (1995), 351--371.



\bibitem{Cer2001}  
Cerrai, S.
\newblock{ \it Second order {PDE}'s in finite and infinite dimension}, Lecture Notes in Mathematics, vol. 1762, Springer-Verlag, Berlin, 2001.

\bibitem{CCE2025} Chaintron, Louis-Pierre, Giovanni Conforti, and Katharina Eichinger. "Propagation of weak log-concavity along generalised heat flows via Hamilton-Jacobi equations." arXiv preprint arXiv:2508.07931 (2025).


\bibitem{Dap2004} 
Da Prato, G. 
\newblock{ \it Kolmogorov equations for stochastic PDEs}, Advanced Courses in Mathematics. CRM Barcelona, Birkh\"{a}user Verlag, Basel, 2004.

\bibitem{Dap-Tub2001} 
Da Prato, G., Tubaro, L.
\newblock{\it Some results about dissipativity of Kolmogorov operators}, Czechoslovak Mathematical Journal, {\bf 51} (2001), 685--699.

\bibitem{Dap-Zab2014}
Da Prato, G., Zabczyk, J.
\newblock{ \it Stochastic equations in infinite dimensions}, Encyclopedia of Mathematics and its Applications, {\bf 152}, Cambridge University Press, Cambridge, 2014.

\bibitem{Dap-Zab1996}
Da Prato, G., Zabczyk, J.,
\newblock{ \it Ergodicity for Infinite Dimensional Systems}, Cambridge University Press, Cambridge, 1996.


\bibitem{Dap-Zab2002}
Da Prato, G., Zabczyk, J.,
\newblock{ \it Second order partial differential equations in {H}ilbert spaces}, London Mathematical Society Lecture Note Series, vol. 293, Cambridge University Press, Cambridge, 2002.

\bibitem{Eng-Nag2006} 
Engel, K.-J., Nagel, R.
\newblock{ \it A short course on operator semigroups}, Universitext, Springer, New York, 2006.

\bibitem{Gol-Koc2001} 
Goldys, B., Kocan, M.
\newblock{ \it Diffusion Semigroups in Spaces of Continuous Functions with Mixed Topology}, J. Differential Equations {\bf 173} (2001), 17--39.


\bibitem{Hai2008}
M. Hairer
\newblock{ \it Ergodic theory for Stochastic PDEs}, www.hairer.org/notes/Imperial.pdf


\bibitem{Lun1995} 
Lunardi, A.
\newblock{ \it Analytic semigroups and optimal regularity in parabolic problems}, Modern Birkh\"{a}user Classics, Birkh\"{a}user/Springer Basel AG, Basel, 1995.

\bibitem{Manca2009} 
Manca, L.
\newblock{\it Fokker-Planck equation for Kolmogorov operators with unbounded coefficients}, Stoch. Anal. Appl. \textbf{27}, 747--769, 2009.

\bibitem{MT1993III} Meyn, Sean P., and Richard L. Tweedie. "Stability of Markovian processes III: Foster–Lyapunov criteria for continuous-time processes." Advances in Applied Probability 25.3 (1993): 518-548.

\bibitem{Gol-Gat1997}
D. Gatarek, B. Goldys, 
\newblock{\it On Invariant Measures for Diffusions on Banach Spaces.} Potential Analysis 7, 533–553 (1997). https://doi.org/10.1023/A:1008663614438

\bibitem{Kaleta-Kwa-Mal2016}
K.~Kaleta, M.~Kwa\'snicki, J.~Ma\l ecki,
\emph{Asymptotic estimate of eigenvalues of pseudo-differential operators in an interval},
J. Math. Anal. Appl. \textbf{439} (2016), no.~2, 896--924.


\bibitem{Paz1983}
A. Pazy
\newblock{ \it Semigroups of Linear Operators and Applications to Partial Differential Equations}, Springer New York, NY, 1983.


\bibitem{Pri1999}
Priola E., 
\newblock{\it On a class of Markov type semigroups in spaces of uniformly continuous and bounded functions}, Studia Math., 136, 271--295, 1999.




 \bibitem{AmbrosioGigliSavare}Luigi Ambrosio, Nicola Gigli, and  Giuseppe
  Savar{\'e}. {\newblock}\tmtextit{Gradient flows in metric spaces and in the
  space of probability measures}. {\newblock}Lectures in Mathematics ETH
  Z{\"u}rich. Birkh{\"a}user Verlag, Basel, Second  edition, 2008.{\newblock}
  
  \bibitem{BogachevDifferentiablebook}Vladimir~I.~Bogachev.
  {\newblock}\tmtextit{Differentiable measures and the Malliavin calculus}, 
  volume  164  of \tmtextit{Mathematical Surveys and Monographs}.
  {\newblock}American Mathematical Society, Providence, RI, 2010.{\newblock}
  
  \bibitem{BogachevRockner1995}Vladimir~I.~Bogachev  and  Michael
  R{\"o}ckner. {\newblock}Regularity of invariant measures on finite- and
  infinite-dimensional spaces and applications. {\newblock}\tmtextit{J. Funct.
  Anal.}, 133(1):168--223, 1995.{\newblock}
  
  \bibitem{Cerraibook}Sandra Cerrai. {\newblock}\tmtextit{Second order
  PDE's in finite and infinite dimension},  volume  1762  of \tmtextit{Lecture
  Notes in Mathematics}. {\newblock}Springer-Verlag, Berlin, 2001.
  {\newblock}A probabilistic approach.{\newblock}
  
  \bibitem{FukushimaChen2012}Zhen-Qing Chen  and  Masatoshi Fukushima.
  {\newblock}\tmtextit{Symmetric Markov processes, time change, and boundary
  theory},  volume~35  of \tmtextit{London Mathematical Society Monographs
  Series}. {\newblock}Princeton University Press, Princeton, NJ,
  2012.{\newblock}
  
  \bibitem{CruzNaibo2016}David Cruz-Uribe  and  Virginia Naibo.
  {\newblock}Kato-Ponce inequalities on weighted and variable Lebesgue spaces.
  {\newblock}\tmtextit{Differential Integral Equations}, 29(9-10):801--836,
  2016.{\newblock}
  
  \bibitem{DaPratoKolmogorov}Giuseppe Da Prato.
  {\newblock}\tmtextit{Kolmogorov equations for stochastic PDEs}.
  {\newblock}Advanced Courses in Mathematics. CRM Barcelona. Birkh{\"a}user
  Verlag, Basel, 2004.{\newblock}
  
  \bibitem{DaPratoZabczykred}Giuseppe Da Prato  and  Jerzy Zabczyk.
  {\newblock}\tmtextit{Stochastic equations in infinite dimensions}, 
  volume~44  of \tmtextit{Encyclopedia of Mathematics and its Applications}.
  {\newblock}Cambridge University Press, Cambridge, 1992.{\newblock}
  
  \bibitem{de2022singular}Francesco~C.~De Vecchi, Massimiliano Gubinelli,
  and  Mattia Turra. {\newblock}A singular integration by parts formula for
  the exponential Euclidean QFT on the plane. {\newblock}\tmtextit{ArXiv
  preprint arXiv:2212.05584}, 2022.{\newblock}
  
  \bibitem{dyda2017eigenvalues}Bart{\l}omiej Dyda, Alexey Kuznetsov, and 
  Mateusz Kwa{\'s}nicki. {\newblock}Eigenvalues of the fractional laplace
  operator in the unit ball. {\newblock}\tmtextit{Journal of the London
  Mathematical Society}, 95(2):500--518, 2017.{\newblock}
  
  \bibitem{Fukushima2011}Masatoshi Fukushima, Yoichi Oshima, and 
  Masayoshi Takeda. {\newblock}\tmtextit{Dirichlet forms and symmetric Markov
  processes},  volume~19  of \tmtextit{De Gruyter Studies in Mathematics}.
  {\newblock}Walter de Gruyter \& Co., Berlin, extended  edition,
  2011.{\newblock}
  
  \bibitem{Lorinczi2011}J{\'o}zsef L{\H o}rinczi, Fumio Hiroshima, and 
  Volker Betz. {\newblock}\tmtextit{Feynman-Kac-type theorems and Gibbs
  measures on path space},  volume~34  of \tmtextit{De Gruyter Studies in
  Mathematics}. {\newblock}Walter de Gruyter \& Co., Berlin, 2011.
  {\newblock}With applications to rigorous quantum field theory.{\newblock}
  
  \bibitem{Rockner1985}Michael R{\"o}ckner. {\newblock}Generalized Markov
  fields and Dirichlet forms. {\newblock}\tmtextit{Acta Appl. Math.},
  3(3):285--311, 1985.{\newblock}

  \bibitem{TriebelbookIII} Triebel, H. (2006). Theory of function spaces III. Basel: Birkhauser Basel.

 \bibitem{Jamison1974} B. Jamison, ``Reciprocal processes,'' \emph{Z. Wahrscheinlichkeitstheorie verw. Gebiete} 30 (1974), 65--86.

 \bibitem{LeonardRoellyZambrini2014} C. L\'eonard, S. R{\oe}lly, and J.-C. Zambrini,
 ``Reciprocal processes: A measure-theoretical point of view,'' \emph{Probability Surveys} 11 (2014), 237--269.

 \bibitem{Agmon1982} Shmuel Agmon. \newblock \emph{Lectures on Exponential Decay of Solutions of Second-Order Elliptic Equations: Bounds on Eigenfunctions of $N$-Body Schr\"odinger Operators}. \newblock Mathematical Notes 29, Princeton University Press, 1982.

 \bibitem{CordobaCordoba2004}
C\'ordoba, A., C\'ordoba, D.
\newblock{\it A maximum principle applied to quasi-geostrophic equations},
Commun. Math. Phys. {\bf 249} (2004), 511--528.

\bibitem{GilbargTrudingerbook}
D.~Gilbarg and N.~S.~Trudinger,
\emph{Elliptic Partial Differential Equations of Second Order},
Classics in Mathematics, Springer, Berlin, 2001, reprint of the 1998 edition.


\bibitem{Georgii2011}
H.-O.~Georgii,
\newblock \emph{Gibbs Measures and Phase Transitions},
\newblock 2nd extended ed., De Gruyter Studies in Mathematics, vol.~9,
Walter de Gruyter, Berlin, 2011.

\bibitem{Ruelle1968}
D.~Ruelle,
\newblock Statistical mechanics of a one-dimensional lattice gas,
\newblock \emph{Commun. Math. Phys.} \textbf{9} (1968), 267--278.

\bibitem{CourregeRenouard1975}
Ph.~Courr\`ege and P.~Renouard,
\newblock Oscillateur anharmonique. Mesures quasi-invariantes sur
$C(\mathbb R,\mathbb R)$ et th\'eorie quantique des champs en dimension
$d=1$,
\newblock \emph{Ast\'erisque} \textbf{22--23} (1975), 3--245.

\bibitem{RoyerYor1976}
G.~Royer and M.~Yor,
\newblock Repr\'esentation int\'egrale de certaines mesures quasi-invariantes
sur $C(\mathbb R)$; mesures extr\'emales et propri\'et\'e de Markov,
\newblock \emph{Ann. Inst. Fourier (Grenoble)} \textbf{26} (1976), no.~2,
7--24.

\bibitem{BetzLorinczi2003}
V.~Betz and J.~L{\H o}rinczi,
\newblock Uniqueness of Gibbs measures relative to Brownian motion,
\newblock \emph{Ann. Inst. H. Poincar\'e Probab. Statist.} \textbf{39}
(2003), no.~5, 877--889.

\bibitem{Baker1961}
G.~A.~Baker, Jr.,
\newblock One-dimensional order-disorder model which approaches a
second-order phase transition,
\newblock \emph{Phys. Rev.} \textbf{122} (1961), 1477--1484.

\bibitem{KacThompson1966}
M.~Kac and C.~J.~Thompson,
\newblock On the mathematical mechanism of phase transition,
\newblock \emph{Proc. Natl. Acad. Sci. USA} \textbf{55} (1966), 676--683.

\bibitem{HilgertMayer2002}
J.~Hilgert and D.~Mayer,
\newblock Transfer operators and dynamical zeta functions for a class of
lattice spin models,
\newblock \emph{Commun. Math. Phys.} \textbf{232} (2002), 19--58.






\end{thebibliography}

\end{document}